\documentclass{surv-l}

\usepackage{setspace}
\usepackage{macros}
\standardpackages\standardrenews\standardmakes\standardlabeling
\colorcommentstrue

\makeatletter
\def\l@figure{\@tocline{0}{3pt plus2pt}{0pt}{3em}{}}
\makeatother

\makeindex

\makeatletter
\usepackage{xstring}
\renewcommand{\tocappendix}[3]{%
  \indentlabel{\IfStrEq{#3}{Bibliography}{}{#1}\@ifnotempty{#2}{ #2.\quad}}#3}
\makeatother

\def\ellipseA{[rotate around={39.01185162422447:(-4.15,2.2600000000000002)}] (-4.15,2.2600000000000002) ellipse (1.4310425486624263cm and 1.0070664208890356cm)}
\def\ellipseB{[rotate around={31.06516488549863:(-2.7499999999999996,3.64)}] (-2.7499999999999996,3.64) ellipse (1.4391779605649575cm and 1.0641114613497555cm)}
\def\ellipseC{[rotate around={-31.930682103717796:(1.8900000000000008,3.85)}] (1.8900000000000008,3.85) ellipse (1.3500580747303617cm and 1.0778018394605968cm)}
\def\ellipseD{[rotate around={-42.089162173832236:(2.9,2.5199999999999996)}] (2.9,2.5199999999999996) ellipse (1.4375510221432186cm and 1.1698516748994354cm)}
\def\ellipseE{[rotate around={0.6030911943805329:(-1.6300000000000001,-0.8499999999999999)}] (-1.6300000000000001,-0.8499999999999999) ellipse (1.3858345482283274cm and 1.0089288354800905cm)}
\def\ellipseF{[rotate around={6.137255949261987:(0.06999999999999999,-1.1000000000000003)}] (0.06999999999999999,-1.1000000000000003) ellipse (1.4251601572520578cm and 1.0752587938811333cm)}

\newcommand{\dexp}{\alpha}
\newcommand{\lexp}{\dexp_*}
\newcommand{\uexp}{\dexp^*}
\newcommand{\thicken}{N}
\newcommand{\schottkycl}{\cl}

\let\O\undefined
\DeclareMathOperator{\O}{O}
\DeclareMathOperator{\Aut}{Aut}
\DeclareMathOperator{\Out}{Out}
\renewcommand{\prec}{\preceq}
\newcommand{\LambdaG}{\Lambda}

\newcommand{\ROSS}{ROSSONCT}
\newcommand{\ROSSs}{ROSSONCTs}
\newcommand{\CROSS}{ROSSONCT}
\newcommand{\CROSSs}{ROSSONCTs}
\newcommand{\bp}{p}
\newcommand{\bq}{q}
\newcommand{\dox}[1]{\|#1\|}
\newcommand{\doxvar}[1]{\dist(\zero,#1)}
\newcommand{\dogo}[1]{\|#1\|}
\newcommand{\dogovar}[1]{\dist(\zero,#1(\zero))}

\newcommand{\point}[2]{p(#1,#2)}
\newcommand{\set}[2]{A(#1,#2)}
\newcommand{\bij}[2]{\psi_{#1,#2}}

\newcommand{\SigmaVWA}{\Sigma_{\text{div}}}

\theoremlevel{section}

\numberwithin{equation}{section}
\numberwithin{figure}{section}

\theoremstyle{theorem}

\newtheorem{problem}{Problem}[]

\theoremstyle{definition}
\maketheorems{standingassumption}{Standing Assumption}
\maketheorems{edelsteinexample}{Edelstein-type Example}

\draftfalse

\newcommand{\thispaper}{this monograph}

\newcommand{\limitset}{\mathcal J}
\newcommand{\refl}{\mathrm{ref \hspace{0 in} l}}
\DeclareMathOperator{\interior}{Int}
\newcommand{\notzero}{w}

\begin{document}

\frontmatter

\title[Geometry and dynamics in hyperbolic metric spaces]{
\uppercase{
\large Geometry and dynamics in Gromov hyperbolic metric spaces \\ {\small{\textnormal{With an emphasis on non-proper settings}}}
}
}


\author{Tushar Das}
\address{University of Wisconsin -- La Crosse, Department of Mathematics \& Statistics, 1725 State Street, La Crosse, WI 54601, USA}
\email{tdas\at uwlax.edu}
\urladdr{\url{https://sites.google.com/a/uwlax.edu/tdas/}}

\author{David Simmons}
\address{University of York, Department of Mathematics, Heslington, York YO10 5DD, UK}
\email{David.Simmons\at york.ac.uk}
\urladdr{\url{https://sites.google.com/site/davidsimmonsmath/}}

\author{Mariusz Urba\'nski}
\address{University of North Texas, Department of Mathematics, 1155 Union Circle \#311430, Denton, TX 76203-5017, USA}
\email{urbanski\at unt.edu}
\urladdr{\url{http://www.urbanskimath.com/}}




\begin{Abstract}
Our monograph presents the foundations of the theory of groups and semigroups acting isometrically on Gromov hyperbolic metric spaces. We make it a point to avoid any assumption of properness/compactness, keeping in mind the motivating example of $\mathbb H^\infty$, the infinite-dimensional rank-one symmetric space of noncompact type over the reals. We have not skipped over parts that might be thought of as ``trivial'' extensions of the finite-dimensional/proper theory, as our intuition has turned out to be wrong often enough regarding these matters that we feel it is worth writing everything down explicitly and with an eye to detail. Moreover, we feel that it is a methodological advantage to have the entire theory presented from scratch, in order to provide a basic reference for the theory.

Though our work unifies and extends a long list of results obtained by many authors, Part \ref{partpreliminaries} of this monograph may be treated as mostly expository. The remainder of this work, some of whose highlights are described in brief below, contains several new methods, examples, and theorems. In Part \ref{partbishopjones}, we introduce a modification of the Poincar\'e exponent, an invariant of a group which provides more information than the usual Poincar\'e exponent, which we then use to vastly generalize the Bishop--Jones theorem relating the Hausdorff dimension of the radial limit set to the Poincar\'e exponent of the underlying semigroup. We construct examples which illustrate the surprising connection between Hausdorff dimension and various notions of discreteness which show up in non-proper settings. Part \ref{partexamples} of the monograph provides a number of examples of groups acting on $\mathbb H^\infty$ which exhibit a wide range of phenomena not to be found in the finite-dimensional theory. Such examples often demonstrate the optimality of our theorems.

In Part \ref{partahlforsthurston}, we construct Patterson--Sullivan measures for groups of divergence type without any compactness assumption on either the boundary or the limit set. This is carried out by first constructing such measures on the Samuel--Smirnov compactification of the bordification of the underlying hyperbolic space, and then showing that the measures are supported on the (non-compactified) bordification. 
We end with a study of quasiconformal measures of geometrically finite groups in terms of \emph{doubling} and \emph{exact dimensionality}. Our analysis characterizes exact dimensionality in terms of Diophantine approximation on the boundary. We demonstrate that though all doubling Patterson--Sullivan measures are exact dimensional, there exist Patterson--Sullivan measures that are exact dimensional but not doubling, as well as ones that are neither doubling nor exact dimensional.
\end{Abstract}
\maketitle

\cleardoublepage
\thispagestyle{empty}
\vspace*{13.5pc}
\begin{center}
Dedicated to the memory of our friend\\ 
{\bf \Large Bernd O. Stratmann}\\
Mathematiker\\
17th July 1957 -- 8th August 2015 
\end{center}
\cleardoublepage

\setcounter{page}{7}

\tableofcontents

\listoffigures

\chapter*{Prologue}

\bigskip

\begin{quotation}
\emph{\dots Cela suffit pour faire comprendre que dans les cinq m\'emoires des Acta mathematica que j'ai consacr\'es \`a l'\'etude des transcendantes fuchsiennes et klein\'eennes, je n'ai fait qu'effleurer un sujet tr\`es vaste, qui fournira sans doute aux g\'eom\`etres l'occasion de nombreuses et importantes d\'ecouvertes.\Footnote{This is enough to make it apparent that in these five memoirs in Acta Mathematica which I have dedicated to the study of Fuschian and Kleinian transcendants, I have only skimmed the surface of a very broad subject, which will no doubt provide geometers with the opportunity for many important discoveries.}}\flushright{ -- H. Poincar\'e, Acta Mathematica, {\bf 5}, 1884, p. 278.}
\end{quotation}

\smallskip

The theory of discrete subgroups of real hyperbolic space has a long history. It was inaugurated by Poincar\'e, who developed the two-dimensional (Fuchsian) and three-dimensional (Kleinian) cases of this theory in a series of articles published between 1881 and 1884 that included numerous notes submitted to the C. R. Acad. Sci. Paris,
a paper at Klein's request in Math. Annalen, and five memoirs commissioned by Mittag-Leffler for his then freshly-minted Acta Mathematica. One must also mention the complementary work of the German school that came before Poincar\'e and continued well after he had moved on to other areas, viz. that of Klein, Schottky, Schwarz, and Fricke. See \cite[Chapter 3]{Gray2} for a brief exposition of this fascinating history, and \cite{Gray1, Saint-Gervais} for more in-depth presentations of the mathematics involved.

We note that in finite dimensions, the theory of \emph{higher-dimensional Kleinian groups}, i.e., discrete isometry groups of the hyperbolic $d$-space $\H^d$ for $d \geq 4$, is markedly different from that in $\H^3$ and $\H^2$. For example, the Teichm\"uller theory used by the Ahlfors--Bers school (viz. Marden, Maskit, J\o rgensen, Sullivan, Thurston, etc.) to study three-dimensional Kleinian groups has no generalization to higher dimensions. Moreover, the recent resolution of the Ahlfors measure conjecture \cite{Agol,CalegariGabai} has more to do with three-dimensional topology than with analysis and dynamics. Indeed, the conjecture remains open in higher dimensions \cite[p. 526, last paragraph]{Kapovich}. Throughout the twentieth century, there are several instances of theorems proven for three-dimensional Kleinian groups whose proofs extended easily to $n$ dimensions (e.g. \cite{BeardonMaskit,Mostow}), but it seems that the theory of higher-dimensional Kleinian groups was not really considered a subject in its own right until around the 1990s. For more information on the theory of higher-dimensional Kleinian groups, see the survey article \cite{Kapovich}, which describes the state of the art up to the last decade, emphasizing connections with homological algebra.

But why stop at finite $n$? Dennis Sullivan, in his IH\'ES \emph{Seminar on Conformal and Hyperbolic Geometry} \cite{Sullivan_seminar} that ran during the late 1970s and early '80s, indicated a possibility of developing the theory of discrete groups acting by hyperbolic isometries on the open unit ball of a separable infinite-dimensional Hilbert space.\Footnote{This was the earliest instance of such a proposal that we could find in the literature, although (as pointed out to us by P. de la Harpe) infinite-dimensional hyperbolic spaces without groups acting on them had been discussed earlier \cite[\627]{Michal}, \cite{Michal2, Harpe_thesis}. It would be of interest to know whether such an idea may have been discussed prior to that.} Later in the early '90s, Misha Gromov observed the paucity of results regarding such actions in his seminal lectures \emph{Asymptotic Invariants of Infinite Groups} \cite{Gromov4} where he encouraged their investigation in memorable terms:
``The spaces like this [infinite-dimensional symmetric spaces] \ldots ~look as cute and sexy to me as their finite dimensional siblings but they have been for years shamefully neglected by geometers and algebraists alike''.

Gromov's lament had not fallen to deaf ears, and the geometry and representation theory of infinite-dimensional hyperbolic space $\H^\infty$ and its isometry group have been studied in the last decade by a handful of mathematicians, see e.g. \cite{BIM, DelzantPy, MonodPy}. However, infinite-dimensional hyperbolic geometry has come into prominence most spectacularly through the recent resolution of a long-standing conjecture in algebraic geometry due to Enriques from the late nineteenth century. Cantat and Lamy \cite{CantatLamy} proved that the Cremona group (i.e. the group of birational transformations of the complex projective plane) has uncountably many non-isomorphic normal subgroups, thus disproving Enriques' conjecture. Key to their enterprise is the fact, due to Manin \cite{Manin}, that the Cremona group admits a faithful isometric action on a non-separable infinite-dimensional hyperbolic space, now known as the Picard--Manin space.

Our project was motivated by a desire to answer Gromov's plea by exposing a coherent general theory of groups acting isometrically on the infinite-dimensional hyperbolic space $\H^\infty$. In the process we came to realize that a more natural domain for our inquiries was the much larger setting of semigroups acting on Gromov hyperbolic metric spaces -- that way we could simultaneously answer our own questions about $\H^\infty$ and construct a theoretical framework for those who are interested in more exotic spaces such as the curve graph, arc graph, and arc complex \cite{HPW, MasurSchleimer, HilionHorbez} and the free splitting and free factor complexes \cite{HandelMosher, BestvinaFeighn, KapovichRafi, HilionHorbez}. These examples are particularly interesting as they extend the well-known dictionary \cite[p.375]{Bestvina_survey} between mapping class groups and the groups $\Out(\F_N)$. In another direction, a dictionary is emerging between mapping class groups and Cremona groups, see \cite{BlancCantat, DHKK}. We speculate that developing the Patterson--Sullivan theory in these three areas would be fruitful and may lead to new connections and analogies that have not surfaced till now.


In a similar spirit, we believe there is a longer story for which this monograph lays the foundations. In general, infinite-dimensional space is a wellspring of outlandish examples and the wide range of new phenomena we have started to uncover has no analogue in finite dimensions. The geometry and analysis of such groups should pique the interests of specialists in probability, geometric group theory, and metric geometry. More speculatively, our work should interact with the ongoing and still nascent study of geometry, topology, and dynamics in a variety of infinite-dimensional spaces and groups, especially in scenarios with sufficient negative curvature. Here are three concrete settings that would be interesting to consider: the universal Teichm\"uller space, the group of volume-preserving diffeomorphisms of $\R^3$ or a $3$-torus, and the space of K\"ahler metrics/potentials on a closed complex manifold in a fixed cohomology class equipped with the Mabuchi--Semmes--Donaldson metric. We have been developing a few such themes. The study of thermodynamics (equilibrium states and Gibbs measures) on the boundaries of Gromov hyperbolic spaces will be investigated in future work \cite{DSU2}.  We speculate that the study of stochastic processes (random walks and Brownian motion) in such settings would be fruitful. Furthermore, it would be of interest to develop the theory of discrete isometric actions and limit sets in infinite-dimensional spaces of higher rank.

\bigskip

{\bf Acknowledgements.} This monograph is dedicated to our colleague Bernd O. Stratmann, who passed away on the 8th of August, 2015. Various discussions with Bernd provided inspiration for this project and we remain grateful for his friendship. The authors thank D. P. Sullivan, D. Mumford, B. Farb, P. Pansu, F. Ledrappier, A. Wilkinson, K. Biswas, E. Breuillard, A. Karlsson, I. Assani, M. Lapidus, R. Guo, Z. Huang, I. Gekhtman, G. Tiozzo, P. Py, M. Davis, M. Roychowdury, M. Hochman, J. Tao, P. de la Harpe, T. Barthelme, J. P. Conze, and Y. Guivarc'h for their interest and encouragement, as well as for invitations to speak about our work at various venues. We are grateful to S. J. Patterson, J. Elstrodt, and \'E. Ghys for enlightening historical discussions on various themes relating to the history of Fuchsian and Kleinian groups and their study through the twentieth century, and to D. P. Sullivan and D. Mumford for suggesting work on diffeomorphism groups and the universal Teichm\"uller space. We thank X. Xie for pointing us to a solution to one of the problems in our Appendix \ref{appendixopenproblems}. The research of the first-named author was supported in part by 2014-2015 and 2016-2017 Faculty Research Grants from the University of Wisconsin--La Crosse. The research of the second-named author was supported in part by the EPSRC Programme Grant EP/J018260/1. The research of the third-named author was supported in part by the NSF grant DMS-1361677.

\chapter{Introduction and Overview}

The purpose of this monograph is to present the theory of groups and semigroups acting isometrically on Gromov hyperbolic metric spaces in full detail as we understand it, with special emphasis on the case of infinite-dimensional algebraic hyperbolic spaces $X = \H_\F^\infty$, where $\F$ denotes a division algebra. We have not skipped over the parts which some would call ``trivial'' extensions of the finite-dimensional/proper theory, for two main reasons: first, intuition has turned out to be wrong often enough regarding these matters that we feel it is worth writing everything down explicitly; second, we feel it is better methodologically to present the entire theory from scratch, in order to provide a basic reference for the theory, since no such reference exists currently (the closest, \cite{BridsonHaefliger}, has a fairly different emphasis). Thus Part \ref{partpreliminaries} of this monograph should be treated as mostly expository, while Parts \ref{partbishopjones}-\ref{partahlforsthurston} contain a range of new material. For experts who want a geodesic path to significant theorems, we list here five such results that we prove in this monograph: Theorems \ref{theorembishopjonesregular} and \ref{theoremglobalmeasureintro} provide generalizations of the Bishop--Jones theorem \cite[Theorem 1]{BishopJones} and the Global Measure Formula \cite[Theorem 2]{StratmannVelani}, respectively, to Gromov hyperbolic metric spaces. Theorem \ref{theoremahlforsthurstongeneral} guarantees the existence of a $\delta$-quasiconformal measure for groups of divergence type, even if the space they are acting on is not proper. Theorem \ref{theoremhlogseriesintro} provides a sufficient condition for the exact dimensionality of the Patterson-Sullivan measure of a geometrically finite group, and Theorem  \ref{theoremequivalentVWAintro} relates the exact dimensionality to Diophantine properties of the measure. However, the reader should be aware that a sharp focus on just these results, without care for their motivation or the larger context in which they are situated, will necessarily preclude access to the interesting and uncharted landscapes that our work has begun to uncover. The remainder of this chapter provides an overview of these landscapes.

\begin{convention}
\label{conventionimplied}
The symbols $\lesssim$, $\gtrsim$, and $\asymp$ will denote coarse asymptotics; a subscript of $\plus$ indicates that the asymptotic is additive, and a subscript of $\times$ indicates that it is multiplicative. For example, $A\lesssim_{\times,K} B$ means that there exists a constant $C > 0$ (the \emph{implied constant}\label{pageimpliedconstant}), depending only on $K$, such that $A\leq C B$. Moreover, $A\lesssim_{\plus,\times}B$ means that there exist constants $C_1,C_2 > 0$ so that $A\leq C_1 B + C_2$. In general, dependence of the implied constant(s) on universal objects such as the metric space $X$, the group $G$, and the distinguished point $\zero\in X$ (cf. Notation \ref{standingassumptions}) will be omitted from the notation.
\end{convention}

\begin{convention}
\label{conventiontendston}
The notation $x_n\tendsto n x$ means that $x_n\to x$ as $n\to \infty$, while the notation $x_n\tendsto{n,\plus} x$ means that
\[
x \asymp_\plus \limsup_{n\to\infty}x_n \asymp_\plus \liminf_{n\to\infty} x_n,
\]
and similarly for $x_n\tendsto{n,\times} x$.
\end{convention}

\begin{convention}
The symbol $\triangleleft$ is used to indicate the end of a nested proof.
\end{convention}

\begin{convention}
We use the Iverson bracket notation: 
$$[\text{statement}] = \begin{cases} 1 & \text{statement true} \\ 0 & \text{statement false} \end{cases}$$
\end{convention}

\begin{convention}
Given a distinguished point $\zero\in X$, we write
\begin{align*}
\dox x = \dist(\zero,x) \text{ and } \dogo g = \dox{g(\zero)}.
\end{align*}
\end{convention}

\section{Preliminaries}
\subsection{Algebraic hyperbolic spaces}

Although we are mostly interested in this monograph in the \emph{real} infinite-dimensional hyperbolic space $\H_\R^\infty$, the complex and quaternionic hyperbolic spaces $\H_\C^\infty$ and $\H_\Q^\infty$ are also interesting. In finite dimensions, these spaces constitute (modulo the Cayley hyperbolic plane\Footnote{We omit all discussion of the Cayley hyperbolic plane $\amsbb H_{\amsbb O}^2$, as the algebra involved is too exotic for our taste; cf. Remark \ref{remarkH2O}.}) the \emph{rank one symmetric spaces of noncompact type}. In the infinite-dimensional case we retain this terminology by analogy; cf. Remark \ref{remarkROSSONCTterminology}. For brevity we will refer to a rank one symmetric space of noncompact type as an {\it algebraic hyperbolic space}.

There are several equivalent ways to define algebraic hyperbolic spaces; these are known as ``models'' of hyperbolic geometry. We consider here the hyperboloid model, ball model (Klein's, not Poincar\'e's), and upper half-space model (which only applies to algebraic hyperbolic spaces defined over the reals, which we will call \emph{real hyperbolic spaces}), which we denote by $\H_\F^\alpha$, $\B_\F^\alpha$, and $\E^\alpha$, respectively. Here $\F$ denotes the base field (either $\R$, $\C$, or $\Q$), and $\alpha$ denotes a cardinal number. We omit the base field when it is $\R$, and denote the exponent by $\infty$ when it is $\#(\N)$, so that $\H^\infty = \H_\R^{\#(\N)}$ is the unique separable infinite-dimensional real hyperbolic space.

The main theorem of Chapter \ref{sectionROSSONCTs} is Theorem \ref{theoremisometries}, which states that any isometry of an algebraic hyperbolic space must be an ``algebraic'' isometry. The finite-dimensional case is given as an exercise in Bridson--Haefliger \cite[Exercise II.10.21]{BridsonHaefliger}. We also describe the relation between totally geodesic subsets of algebraic hyperbolic spaces and fixed point sets of isometries (Theorem \ref{theoremtotallygeodesic}), a relation which will be used throughout the paper.

\begin{remark}
Key to the study of finite-dimensional algebraic hyperbolic spaces is the theory of quasiconformal mappings (e.g., as in Mostow and Pansu's rigidity theorems \cite{Mostow, Pansu1}). Unfortunately, it appears to be quite difficult to generalize this theory to infinite dimensions. For example, it is an open question \cite[p.1335]{Heinonen2} whether every quasiconformal homeomorphism of Hilbert space is also quasisymmetric.
\end{remark}

\subsection{Gromov hyperbolic metric spaces}
\label{subsubsectionlist}

Historically, the first motivation for the theory of negatively curved metric spaces came from differential geometry and the study of negatively curved Riemannian manifolds. The idea was to describe the most important consequences of negative curvature in terms of the metric structure of the manifold. This approach was pioneered by Aleksandrov \cite{Aleksandrov}, who discovered for each $\kappa\in\R$ an inequality regarding triangles in a metric space with the property that a Riemannian manifold satisfies this inequality if and only if its sectional curvature is bounded above by $\kappa$, and popularized by Gromov, who called Aleksandrov's inequality the ``CAT($\kappa$) inequality'' as an abbreviation for ``comparison inequality of Alexandrov--Toponogov'' \cite[p.106]{Gromov3}.\Footnote{It appears that Bridson and Haefliger may be responsible for promulgating the idea that the C in CAT refers to E. Cartan \cite[p.159]{BridsonHaefliger}. We were unable to find such an indication in \cite{Gromov3}, although Cartan is referenced in connection with some theorems regarding CAT($\kappa$) spaces (as are Riemann and Hadamard).} A metric space is called CAT($\kappa$) if the distance between any two points on a geodesic triangle is smaller than the corresponding distance on the ``comparison triangle'' in a model space of constant curvature $\kappa$; see Definition \ref{definitionCAT}.

The second motivation came from geometric group theory, in particular the study of groups acting on manifolds of negative curvature. For example, Dehn proved that the word problem is solvable for finitely generated Fuchsian groups \cite{Dehn}, and this was generalized by Cannon to groups acting cocompactly on manifolds of negative curvature \cite{Cannon}. Gromov attempted to give a geometric characterization of these groups in terms of their Cayley graphs; he tried many definitions (cf. \cite[\66.4]{Gromov1}, \cite[\64]{Gromov2}) before converging to what is now known as Gromov hyperbolicity in 1987 \cite[1.1, p.89]{Gromov3}, a notion which has influenced much research. A metric space is said to be \emph{Gromov hyperbolic} if it satisfies a certain inequality that we call \emph{Gromov's inequality}; see Definition \ref{definitiongromovhyperbolic}. A finitely generated group is then said to be \emph{word-hyperbolic} if its Cayley graph is Gromov hyperbolic.

The big advantage of Gromov hyperbolicity is its generality. We give some idea of its scope by providing the following nested list of metric spaces which have been proven to be Gromov hyperbolic:

\begin{itemize}
\item CAT(-1) spaces (Definition \ref{definitionCAT})
\begin{itemize}
\item Riemannian manifolds (both finite- and infinite-dimensional) with sectional curvature $\leq -1$
\begin{itemize}
\item Algebraic hyperbolic spaces (Definition \ref{definitionROSSONCT})
\begin{itemize}
\item Picard--Manin spaces of projective surfaces defined over algebraically closed fields \cite{Manin}, cf. \cite[\63.1]{Cantat_Cremona}
\end{itemize}
\end{itemize}
\item $\R$-trees (Definition \ref{definitionRtree})
\begin{itemize}
\item Simplicial trees
\begin{itemize}
\item Unweighted simplicial trees
\end{itemize}
\end{itemize}
\end{itemize}
\item Cayley metrics (Example \ref{examplecayleygraph}) on word-hyperbolic groups
\item Green metrics on word-hyperbolic groups \cite[Corollary 1.2]{BHM}
\item Quasihyperbolic metrics of uniform domains in Banach spaces \cite[Theorem 2.12]{Vaisala2}
\item Arc graphs and curve graphs \cite{HPW} and arc complexes \cite{MasurSchleimer, HilionHorbez} of finitely punctured oriented surfaces
\item Free splitting complexes \cite{HandelMosher, HilionHorbez} and free factor complexes \cite{BestvinaFeighn, KapovichRafi, HilionHorbez}
\end{itemize}

\begin{remark}
Many of the above examples admit natural isometric group actions:
\begin{itemize}
\item The Cremona group acts isometrically on the Picard--Manin space \cite{Manin}, cf. \cite[Theorem 3.3]{Cantat_Cremona}.
\item The mapping class group of a finitely punctured oriented surface acts isometrically on its arc graph, curve graph, and arc complex.
\item The outer automorphism group $\Out(\F_N)$ of the free group on $N$ generators acts isometrically on the free splitting complex $\FF\SS(\F_N)$ and the free factor complex $\FF\FF(\F_N)$.
\end{itemize}
\end{remark}

\begin{remark}
\label{remarkproper}
Most of the above examples are examples of \emph{non-proper} hyperbolic metric spaces. Recall that a metric space is said to be \emph{proper} if its distance function $x\mapsto \dox x = \doxvar x$ is proper, or equivalently if closed balls are compact. Though much of the existing literature on CAT(-1) and hyperbolic metric spaces assumes that the spaces in question are proper, it is often not obvious whether this assumption is really essential. However, since results about proper metric spaces do not apply to infinite-dimensional algebraic hyperbolic spaces, we avoid the assumption of properness.
\end{remark}

\begin{remark}
\label{remarkgeodesic}
One of the above examples, namely, Green metrics on word-hyperbolic groups, is a natural class of \emph{non-geodesic} hyperbolic metric spaces.\Footnote{Quasihyperbolic metrics on uniform domains in Banach spaces can also fail to be geodesic, but they are \emph{almost geodesic} which is almost as good. See e.g. \cite{Vaisala} for a study of almost geodesic hyperbolic metric spaces.} However, Bonk and Schramm proved that all non-geodesic hyperbolic metric spaces can be isometrically embedded into geodesic hyperbolic metric spaces \cite[Theorem 4.1]{BonkSchramm}, and the equivariance of their construction was proven by Blach\`ere, Ha\"issinsky, and Mathieu \cite[Corollary A.10]{BHM}. Thus, one could view the assumption of geodesicity to be harmless, since most theorems regarding geodesic hyperbolic metric spaces can be pulled back to non-geodesic hyperbolic metric spaces. However, for the most part we also avoid the assumption of geodesicity, mostly for methodological reasons rather than because we are considering any particular non-geodesic hyperbolic metric space. Specifically, we felt that Gromov's definition of hyperbolicity in metric spaces is a ``deep'' definition whose consequences should be explored independently of such considerations as geodesicity. We do make the assumption of geodesicity in Chapter \ref{sectionGF}, where it seems necessary in order to prove the main theorems. (The assumption of geodesicity in Chapter \ref{sectionGF} can for the most part be replaced by the weaker assumption of almost geodesicity \cite[p.271]{BonkSchramm}, but we felt that such a presentation would be more technical and less intuitive.)
\end{remark}

We now introduce a list of standing assumptions and notations. They apply to all chapters except for Chapters \ref{sectionROSSONCTs}, \ref{sectiongeometry1}, and \ref{sectiondiscreteness} (see also \6\ref{standingassumptions2}).

\begin{notation}
\label{standingassumptions}
Throughout the introduction,
\begin{itemize}
\item $X$ is a Gromov hyperbolic metric space (cf. Definition \ref{definitiongromovhyperbolic}),
\item $\dist$ denotes the distance function of $X$,
\item $\del X$ denotes the Gromov boundary of $X$, and $\bord X$ denotes the bordification $\bord X = X\cup \del X$ (cf. Definition \ref{definitiongromovboundary}),
\item $\Dist$ denotes a visual metric on $\del X$ with respect to a parameter $b > 1$ and a distinguished point $\zero\in X$ (cf. Proposition \ref{propositionDist}). By definition, a visual metric satisfies the asymptotic
\begin{equation}
\label{distanceasymptoticv1}
\Dist_{b,\zero}(\xi,\eta) \asymp_\times b^{-\lb\xi|\eta\rb_\zero},
\end{equation}
where $\lb\cdot|\cdot\rb$ denotes the Gromov product (cf. \eqref{gromovproduct}).
\item $\Isom(X)$ denotes the isometry group of $X$. Also, $G\leq\Isom(X)$ will mean that $G$ is a subgroup of $\Isom(X)$, while $G\prec\Isom(X)$ will mean that $G$ is a subsemigroup of $\Isom(X)$.
\end{itemize}
\end{notation}

A prime example to have in mind is the special case where $X$ is an infinite-dimensional algebraic hyperbolic space, in which case the Gromov boundary $\del X$ can be identified with the natural boundary of $X$ (Proposition \ref{propositionboundariesequivalent}), and we can set $b = e$ and get equality in \eqref{distanceasymptoticv1} (Observation \ref{observationDist}).

Another important example of a hyperbolic metric space that we will keep in our minds is the case of $\R$-trees alluded to above. $\R$-trees are a generalization of simplicial trees, which in turn are a generalization of unweighted simplicial trees, also known as ``$\Z$-trees'' or just ``trees''. $\R$-trees are worth studying in the context of hyperbolic metric spaces for two reasons: first of all, they are ``prototype spaces'' in the sense that any finite set in a hyperbolic metric space can be roughly isometrically embedded into an $\R$-tree, with a roughness constant depending only on the cardinality of the set \cite[pp.33-38]{GhysHarpe1}; second of all, $\R$-trees can be equivariantly embedded into infinite-dimensional real hyperbolic space $\H^\infty$ (Theorem \ref{theoremBIM}), meaning that any example of a group acting on an $\R$-tree can be used to construct an example of the same group acting on $\H^\infty$. $\R$-trees are also much simpler to understand than general hyperbolic metric spaces: for any finite set of points, one can draw out a list of all possible diagrams, and then the set of distances must be determined from one of these diagrams (cf. e.g., Figure \ref{figurequadruple}).

Besides introducing $\R$-trees, CAT(-1) spaces, and hyperbolic metric spaces, the following things are done in Chapter \ref{sectiongeometry1}: construction of the Gromov boundary $\del X$ and analysis of its basic topological properties (Section \ref{subsectionboundary}), proof that the Gromov boundary of an algebraic hyperbolic space is equal to its natural boundary (Proposition \ref{propositionboundariesequivalent}), and the construction of various metrics and metametrics on the boundary of $X$ (Section \ref{subsectionmetametrics}). None of this is new, although the idea of a metametric (due to V\"ais\"al\"a \cite[\64]{Vaisala}) is not very well known.

In Chapter \ref{sectiongeometry2}, we go more into detail regarding the geometry of hyperbolic metric spaces. We prove the geometric mean value theorem for hyperbolic metric spaces (Section \ref{subsectionderivatives}), the existence of geodesic rays connecting two points in the boundary of a CAT(-1) space (Proposition \ref{propositiongeodesicconvergence}), and various geometrical theorems regarding the sets
\[
\Shad_z(x,\sigma) := \{\xi\in\del X : \lb x|\xi\rb_z \leq \sigma\},
\]
which we call ``shadows'' due to their similarity to the famous shadows of Sullivan \cite[Fig. 2]{Sullivan_density_at_infinity} on the boundary of $\H^d$ (Section \ref{subsectionshadows}). We remark that most proofs of the existence of geodesics between points on the boundary of complete CAT(-1) spaces, e.g. \cite[Proposition II.9.32]{BridsonHaefliger}, assume properness and make use of it in a crucial way, whereas we make no such assumption in Proposition \ref{propositiongeodesicconvergence}. Finally, in Section \ref{subsectionpolar} we introduce ``generalized polar coordinates'' in a hyperbolic metric space. These polar coordinates tell us that the action of a loxodromic isometry (see Definition \ref{definitionclassification1}) on a hyperbolic metric space is roughly the same as the map $\xx\mapsto \lambda\xx$ in the upper half-plane $\E^2$.




\subsection{Discreteness}
\label{subsubsectiondiscreteness}

The first step towards extending the theory of Kleinian groups to infinite dimensions (or more generally to hyperbolic metric spaces) is to define the appropriate class of groups to consider. This is less trivial than might be expected. Recalling that a $d$-dimensional Kleinian group is defined to be a discrete subgroup of $\Isom(\H^d)$, we would want to define an infinite-dimensional Kleinian group to be a discrete subgroup of $\Isom(\H^\infty)$. But what does it mean for a subgroup of $\Isom(\H^\infty)$ to be discrete? In finite dimensions, the most natural definition is to call a subgroup discrete if it is discrete relative to the natural topology on $\Isom(\H^d)$; this definition works well since $\Isom(\H^d)$ is a Lie group. But in infinite dimensions and especially in more exotic spaces, many applications require stronger hypotheses (e.g., Theorem \ref{theorembishopjonesregular}, Chapter \ref{sectionGF}). In Chapter \ref{sectiondiscreteness}, we discuss several potential definitions of discreteness, which are inequivalent in general but agree in the case of finite-dimensional space $X = \H^d$ (Proposition \ref{propositionfinitedimequivalent}):

\makeatletter
\def\rep@title{Definitions \ref{definitiondiscreteness} and \ref{definitionparametricdiscreteness}}
\begin{rep@definition}
Fix $G\leq\Isom(X)$.
\begin{itemize}
\item $G$ is called \emph{strongly discrete (SD)} if for every bounded set $B \subset X$, we have
\[
\#\{g\in G: g(B) \cap B \neq \emptyset\} < \infty.
\]
\item $G$ is called \emph{moderately discrete (MD)} if for every $x\in X$, there exists an open set $U$ containing $x$ such that
\[
\#\{g\in G: g(U) \cap U \neq \emptyset\} < \infty.
\]
\item $G$ is called \emph{weakly discrete (WD)} if for every $x\in X$, there exists an open set $U$ containing $x$ such that
\[
g(U) \cap U \neq \emptyset \Rightarrow g(x) = x.
\]
\item $G$ is called \emph{COT-discrete (COTD)} if it is discrete as a subset of $\Isom(X)$ when $\Isom(X)$ is given the compact-open topology (COT).
\item If $X$ is an algebraic hyperbolic space, then $G$ is called \emph{UOT-discrete (UOTD)} if it is discrete as a subset of $\Isom(X)$ when $\Isom(X)$ is given the uniform operator topology (UOT; cf. Section \ref{subsectiontopologies}).
\end{itemize}
\end{rep@definition}\makeatother

As our naming suggests, the condition of strong discreteness is stronger than the condition of moderate discreteness, which is in turn stronger than the condition of weak discreteness (Proposition \ref{propositionSDMDWD}). Moreover, any moderately discrete group is COT-discrete, and any weakly discrete subgroup of $\Isom(\H^\infty)$ is COT-discrete (Proposition \ref{propositionparametricdiscreteness}). These relations and more are summarized in Table \ref{figurediscreteness} on p. \pageref{figurediscreteness}.

Out of all these definitions, strong discreteness should perhaps be thought of as the best generalization of discreteness to infinite dimensions. Thus, we propose that the phrase ``infinite-dimensional Kleinian group'' should mean ``strongly discrete subgroup of $\Isom(\H^\infty)$''. However, in this monograph we will be interested in the consequences of all the different notions of discreteness, as well as the interactions between them.

\begin{remark}
Strongly discrete groups are known in the literature as \emph{metrically proper}, and moderately discrete groups are known as \emph{wandering}. However, we prefer our terminology since it more clearly shows the relationship between the different notions of discreteness.
\end{remark}

\bigskip
\subsection{The classification of semigroups}

After clarifying the different types of discreteness which can occur in infinite dimensions, we turn to the question of classification. This question makes sense both for individual isometries and for entire semigroups.\Footnote{In Chapters \ref{sectionclassification}-\ref{sectionschottky}, we work in the setting of semigroups rather than groups. Like dropping the assumption of geodesicity (cf. Remark \ref{remarkgeodesic}), this is done partly in order to broaden our class of examples and partly for methodological reasons -- we want to show exactly where the assumption of being closed under inverses is being used. It should be also noted that semigroups sometimes show up naturally when one is studying groups; cf. Proposition \ref{propositionnonelementaryequivalent}(B).} Historically, the study of classification began in the 1870s when Klein proved a theorem classifying isometries of $\H^2$ and attached the words ``elliptic'', ``parabolic'', and ``hyperbolic'' to these classifications. Elliptic isometries are those which have at least one fixed point in the interior, while parabolic isometries have exactly one fixed point, which is a neutral fixed point on the boundary, and hyperbolic isometries have two fixed points on the boundary, one of which is attracting and one of which is repelling. Later, the word ``loxodromic'' was used to refer to isometries in $\H^3$ which have two fixed points on the boundary but which are geometrically ``screw motions'' rather than simple translations. In what follows we use the word ``loxodromic'' to refer to all isometries of $\H^n$ (or more generally a hyperbolic metric space) with two fixed points on the boundary -- this is analogous to calling a circle an ellipse. Our real reason for using the word ``loxodromic'' in this instance, rather than ``hyperbolic'', is to avoid confusion with the many other meanings of the word ``hyperbolic'' that have entered usage in various scenarios.

To extend this classification from individual isometries to groups, we call a group ``elliptic'' if its orbits are bounded, ``parabolic'' if it has a unique neutral global fixed point on the boundary, and ``loxodromic'' if it contains at least one loxodromic isometry. The main theorem of Chapter \ref{sectionclassification} (viz. Theorem \ref{classificationofsemigroups}) is that every subsemigroup of $\Isom(X)$ is either elliptic, parabolic, or loxodromic.

Classification of groups has appeared in the literature in various contexts, from Eberlein and O'Neill's results regarding visiblility manifolds \cite{EberleinOneill}, through Gromov's remarks about groups acting on strictly convex spaces \cite[\63.5]{Gromov1} and word-hyperbolic groups \cite[\63.1]{Gromov3}, to the more general results of Hamann \cite[Theorem 2.7]{Hamann}, Osin \cite[\63]{Osin}, and Caprace, de Cornulier, Monod, and Tessera \cite[\63.A]{CCMT} regarding geodesic hyperbolic metric spaces.\Footnote{We remark that the results of \cite[\63.A]{CCMT} can be generalized to non-geodesic hyperbolic metric spaces by using the Bonk--Schramm embedding theorem \cite[Theorem 4.1]{BonkSchramm} (see also \cite[Corollary A.10]{BHM}).} Many of these theorems have similar statements to ours (\cite{Hamann} and \cite{CCMT} seem to be the closest), but we have not kept track of this carefully, since our proof appears to be sufficiently different to warrant independent interest anyway.

After proving Theorem \ref{classificationofsemigroups}, we discuss further aspects of the classification of groups, such as the further classification of loxodromic groups given in \6\ref{subsubsectionloxodromic}: a loxodromic group is called ``lineal'', ``focal'', or ``of general type'' according to whether it has two, one, or zero global fixed points, respectively. (This terminology was introduced in \cite{CCMT}.) The ``focal'' case is especially interesting, as it represents a class of nonelementary groups which have global fixed points.\Footnote{Some sources (e.g. \cite[\65.5]{Ratcliffe}) define nonelementarity in a way such that global fixed points are automatically ruled out, but this is not true of our definition (Definition \ref{definitionelementary}).} We show that certain classes of discrete groups cannot be focal (Proposition \ref{propositionnonfocal}), which explains why such groups do not appear in the theory of Kleinian groups. On the other hand, we show that in infinite dimensions, focal groups can have interesting limit sets even though they satisfy only a weak form of discreteness; cf. Remark \ref{remarkfocalfractal}.

\subsection{Limit sets}

An important invariant of a Kleinian group $G$ is its \emph{limit set} $\Lambda = \Lambda_G$, the set of all accumulation points of the orbit of any point in the interior. By putting an appropriate topology on the bordification of our hyperbolic metric space $X$ (\6\ref{subsubsectiontopologyclX}), we can generalize this definition to an arbitrary subsemigroup of $\Isom(X)$. Many results generalize relatively straightforwardly\Footnote{As is the case for many of our results, the classical proofs use compactness in a crucial way -- so here ``straightforwardly'' means that the statements of the theorems themselves do not require modification.} to this new context, such as the minimality of the limit set (Proposition \ref{propositionminimal}) and the connection between classification and the cardinality of the limit set (Proposition \ref{propositioncardinalitylimitset}). In particular, we call a semigroup \emph{elementary} if its limit set is finite.

In general, the convex hull of the limit set may need to be replaced by a quasiconvex hull (cf. Definition \ref{definitionconvexhull}), since in certain cases the convex hull does not accurately reflect the geometry of the group. Indeed, Ancona \cite[Corollary C]{Ancona} and Borbely \cite[Theorem 1]{Borbely} independently constructed examples of CAT(-1) three-manifolds $X$ for which there exists a point $\xi\in\del X$ such that the convex hull of any neighborhood of $\xi$ is equal to $\bord X$. Although in a non-proper setting the limit set may no longer be compact, compactness of the limit set is a reasonable geometric condition that is satisfied for many examples of subgroups of $\Isom(\H^\infty)$ (e.g. Examples \ref{exampleinfinitepoincareBIM}, \ref{exampleAutTBIM}). We call this condition \emph{compact type} (Definition \ref{definitioncompacttype}).

\bigskip
\section{The Bishop--Jones theorem and its generalization}

The term \emph{Poincar\'e series} classically referred to a variety of averaging procedures, initiated by Poincar\'e in his aforementioned Acta memoirs, with a view towards uniformization of Riemann surfaces via the construction of automorphic forms. Given a Fuchsian group $\Gamma$ and a rational function $H:\what\C\to\what\C$ with no poles on $\del\B^2$, Poincar\'e proved that for every $m\geq 2$ the series
\[
\sum_{\gamma \in \Gamma} H(\gamma(z))(\gamma'(z))^m
\]
(defined for $z$ outside the limit set of $\Gamma$) converges uniformly to an automorphic form of dimension $m$; see \cite[p.218]{Saint-Gervais}. Poincar\'e called these series ``$\theta$-fuchsian series of order $m$'', but the name ``Poincar\'e series'' was later used to refer to such objects.\Footnote{The modern definition of Poincar\'e series (cf. Definition \ref{definitionpoincareexponent}) is phrased in terms of hyperbolic geometry rather than complex analysis, but it agrees with the special case of Poincar\'e's original definition which occurs when $H\equiv 1$ and $z = 0$, with the caveat that $\gamma'(z)^m$ should be replaced by $|\gamma'(z)|^m$.} The question of for which $m < 2$ the Poincar\'e series still converges was investigated by Schottky, Burnside, Fricke, and Ritter; cf. \cite[pp.37-38]{AdelmannGerbracht}.

In what would initially appear to be an unrelated development, mathematicians began to study the ``thickness'' of the limit set of a Fuchsian group: in 1941 Myrberg \cite{Myrberg} showed that the limit set $\Lambda$ of a nonelementary Fuchsian group has positive logarithmic capacity; this was improved by Beardon \cite{Beardon2} who showed that $\Lambda$ has positive Hausdorff dimension, thus deducing Myrberg's result as a corollary (since positive Hausdorff dimension implies positive logarithmic capacity for compact subsets of $\R^2$ \cite{Taylor}). The connection between this question and the Poincar\'e series was first observed by Akaza, who showed that if $G$ is a Schottky group for which the Poincar\'e series converges in dimension $s$, then the Hausdorff $s$-dimensional measure of $\Lambda$ is zero \cite[Corollary of Theorem A]{Akaza2}. Beardon then extended Akaza's result to finitely generated Fuchsian groups \cite[Theorem 5]{Beardon3}, as well as defining the \emph{exponent of convergence} (or \emph{Poincar\'e exponent}) $\delta = \ \delta_G$ of a Fuchsian or Kleinian group to be the infimum of $s$ for which the Poincar\'e series converges in dimension $s$ (cf. Definition \ref{definitionpoincareexponent} and \cite{Beardon1}). The reverse direction was then proven by Patterson \cite{Patterson2} using a certain measure on $\Lambda$ to produce the lower bound, which we will say more about below in \6\ref{subsectionconformal}. Patterson's results were then generalized by Sullivan \cite{Sullivan_density_at_infinity} to the setting of geometrically finite Kleinian groups. The necessity of the geometrical finiteness assumption was demonstrated by Patterson \cite{Patterson3}, who showed that there exist Kleinian groups of the first kind (i.e. with limit set equal to $\del\H^d$) with arbitrarily small Poincar\'e exponent \cite{Patterson3} (see also \cite{Hopf_statistik} or \cite[Example 8]{Starkov} for an earlier example of the same phenomenon).

Generalizing these theorems beyond the geometrically finite case requires the introduction of the \emph{radial} and \emph{uniformly radial} limit sets. In what follows, we will denote these sets by $\Lr$ and $\Lur$, respectively. Note that the radial and uniformly radial limit sets as well as the Poincar\'e exponent can all (with some care) be defined for general hyperbolic metric spaces; see Definitions \ref{definitionradialconvergence}, \ref{definitionlimitset}, and \ref{definitionpoincareexponent}. The radial limit set was introduced by Hedlund in 1936 in his analysis of transitivity of horocycles \cite[Theorem 2.4]{Hedlund}.

After some intermediate results \cite{FernandezMelian,Stratmann1}, Bishop and Jones \cite[Theorem 1]{BishopJones} generalized Patterson and Sullivan by proving that if $G$ is a nonelementary Kleinian group, then $\HD(\Lr) = \HD(\Lur) = \delta$.\Footnote{Although Bishop and Jones' theorem only states that $\HD(\Lr) = \delta$, they remark that their proof actually shows that $\HD(\Lur) = \delta$ \cite[p.4]{BishopJones}.} Further generalization was made by Paulin \cite{Paulin}, who proved the equation $\HD(\Lr) = \delta$ in the case where $G\leq\Isom(X)$, and $X$ is either a word-hyperbolic group, a CAT(-1) manifold, or a locally finite unweighted simplicial tree which admits a discrete cocompact action. We may now state the first major theorem of this monograph, which generalizes all the aforementioned results:


\begin{theorem}
\label{theorembishopjonesregular}
Let $G\leq\Isom(X)$ be a nonelementary group. Suppose either that
\begin{itemize}
\item[(1)] $G$ is strongly discrete,
\item[(2)] $X$ is a CAT(-1) space and $G$ is moderately discrete,
\item[(3)] $X$ is an algebraic hyperbolic space and $G$ is weakly discrete, or that
\item[(4)] $X$ is an algebraic hyperbolic space and $G$ acts irreducibly (cf. Section \ref{subsectionirreducibly}) and is COT-discrete.
\end{itemize}
Then there exists $\sigma > 0$ such that
\begin{equation}
\label{bishopjonesregular}
\HD(\Lr) = \HD(\Lur) = \HD(\Lur\cap\Lrsigma) = \delta
\end{equation}
(cf. Definitions \ref{definitionradialconvergence} and \ref{definitionlimitset} for the definition of $\Lrsigma$); moreover, for every $0 < s < \delta$ there exist $\tau > 0$ and an Ahlfors $s$-regular\Footnote{Recall that a measure $\mu$ on a metric space $Z$ is called \emph{Ahlfors $s$-regular} if for all $z\in Z$ and $0 < r \leq 1$, we have that $\mu(B(z,r)) \asymp_\times r^s$. The topological support of an Ahlfors $s$-regular measure is called an Ahlfors $s$-regular set. \label{footnoteahlforsregular}} set $\limitset_s\subset\Lurtau\cap\Lrsigma$.
\end{theorem}

For the proof of Theorem \ref{theorembishopjonesregular}, see the comments below Theorem \ref{theorembishopjonesmodified}.

\begin{remark*}
We note that weaker versions of Theorem \ref{theorembishopjonesregular} already appeared in \cite{DSU0} and \cite{FSU4}, each of which has a two-author intersection with the present paper. In particular, case (1) of Theorem \ref{theorembishopjonesregular} appeared in \cite{FSU4} and the proofs of Theorem \ref{theorembishopjonesregular} and \cite[Theorem 5.9]{FSU4} contain a number of redundancies. This was due to the fact that we worked on two projects which, despite having fundamentally different objectives, both required essentially the same argument to produce ``large, nice'' subsets of the limit set: in the present monograph, this argument forms the core of the proof of our generalization of the Bishop--Jones theorem, while in \cite{FSU4}, the main use of the argument is in proving the full dimension of the set of badly approximable points, in two different senses of the phrase ``badly approximable" (approximation by the orbits of distinguished points, vs. approximation by rational vectors in an ambient Euclidean space). There are also similarities between the proof of Theorem \ref{theorembishopjonesregular} and the proof of the weaker version found in \cite[Theorem 8.13]{DSU0}, although in this case the presentation is significantly different. However, we remark that the main Bishop--Jones theorem of this monograph, Theorem \ref{theorembishopjonesmodified}, is significantly more powerful than both \cite[Theorem 5.9]{FSU4} and \cite[Theorem 8.13]{DSU0}.
\end{remark*}

\begin{remark*}
The ``moreover'' clause is new even in the case which Bishop and Jones considered, demonstrating that the limit set $\Lur$ can be approximated by subsets which are particularly well distributed from a geometric point of view. It does not follow from their theorem since a set could have large Hausdorff dimension without having any closed Ahlfors regular subsets of positive dimension (much less full dimension); in fact it follows from the work of Kleinbock and Weiss \cite{KleinbockWeiss1} that the set of well approximable numbers forms such a set.\Footnote{It could be objected that this set is not closed and therefore should not constitute a counterexample. However, since it has full measure, it has closed subsets of arbitrarily large measure (which in particular still have dimension 1).} In \cite{FSU4}, a slight strengthening of this clause was used to deduce the full dimension of badly approximable vectors in the radial limit set of a Kleinian group \cite[Theorem 9.3]{FSU4}.
\end{remark*}

\begin{remark*}
It is possible for a group satisfying one of the hypotheses of Theorem \ref{theorembishopjonesregular} to also satisfy $\delta = \infty$ (Examples \ref{exampleinfinitepoincare}-\ref{examplenotstronglydiscreteBIM} and \ref{examplepoincareextensiontwo}-\ref{exampleMDWD});\Footnote{For the parabolic examples, take a Schottky product (Definition \ref{definitionschottkyproduct}) with a lineal group (Definition \ref{definitionCCMT}) to get a nonelementary group, as suggested at the beginning of Chapter \ref{sectionexamples}.} note that Theorem \ref{theorembishopjonesregular} still holds in this case.
\end{remark*}

\begin{remark*}
A natural question is whether \eqref{bishopjonesmodified} can be improved by showing that there exists some $\sigma > 0$ for which $\HD(\Lursigma) = \delta$ (cf. Definitions \ref{definitionradialconvergence} and \ref{definitionlimitset} for the definition of $\Lursigma$). The answer is negative. For a counterexample, take $X = \H^2$ and $G = \SL_2(\Z)\leq\Isom(X)$; then for all $\sigma > 0$ there exists $\epsilon > 0$ such that $\Lursigma\subset\BA(\epsilon)$, where $\BA(\epsilon)$ denotes the set of all real numbers with Lagrange constant at most $1/\epsilon$. (This follows e.g. from making the correspondence in \cite[Observation 1.15 and Proposition 1.21]{FSU4} explicit.) It is well-known (see e.g. \cite{Kurzweil} for a more precise result) that $\HD(\BA(\epsilon)) < 1$ for all $\epsilon > 0$, demonstrating that $\HD(\Lursigma) < 1 = \delta$.
\end{remark*}

\begin{remark*}
Although Theorem \ref{theorembishopjonesregular} computes the Hausdorff dimension of the radial and uniformly radial limit sets, there are many other subsets of the limit set whose Hausdorff dimension it does not compute, such as the horospherical limit set (cf. Definitions \ref{definitionhorosphericalconvergence} and \ref{definitionlimitset}) and the ``linear escape'' sets $(\Lambda_\alpha)_{\alpha\in (0,1)}$ \cite{Lundh}. We plan on discussing these issues at length in \cite{DSU2}.
\end{remark*}

Finally, let us also remark that the hypotheses (1) - (4) cannot be weakened in any of the obvious ways:
\begin{proposition}
\label{propositionbishopjones}
We may have $\HD(\Lr) < \delta$ even if:
\begin{itemize}
\item[(1)] $G$ is moderately discrete (even properly discontinuous) (Example \ref{exampleAutTparttwo}).
\item[(2)] $X$ is a proper CAT(-1) space and $G$ is weakly discrete (Example \ref{exampleAutT}).
\item[(3)] $X = \H^\infty$ and $G$ is COT-discrete (Example \ref{examplepoincareextension}).
\item[(4)] $X = \H^\infty$ and $G$ is irreducible and UOT-discrete (Example \ref{exampleAutTBIM}).
\item[(5)] $X = \H^2$ (Example \ref{examplenondiscrete}).
\end{itemize}
In each case the counterexample group $G$ is of general type (see Definition \ref{definitionCCMT}) and in particular is nonelementary.
\end{proposition}

\subsection{The modified Poincar\'e exponent}
The examples of Proposition \ref{propositionbishopjones} illustrate that the Poincar\'e exponent does not always accurately calculate the Hausdorff dimension of the radial and uniformly radial limit sets. In Chapter \ref{sectionmodified} we introduce a modified version of the Poincar\'e exponent which succeeds at accurately calculating $\HD(\Lr)$ and $\HD(\Lur)$ for all nonelementary groups $G$. (When $G$ is an elementary group, $\HD(\Lr) = \HD(\Lur) = 0$, so there is no need for a sophisticated calculation in this case.) Some motivation for the following definition is given in \6\ref{subsectionmodified}.

\begin{repdefinition}{definitionmodifiedexponent}
Let $G$ be a subsemigroup of $\Isom(X)$.
\begin{itemize}
\item For each set $S\subset X$ and $s\geq 0$, let
\begin{align*}
\Sigma_s(S) &= \sum_{x\in S} b^{-s\dox x}\\
\Delta(S) &= \{s\geq 0: \Sigma_s(S) = \infty\}\\
\delta(S) &= \sup\Delta(S).
\end{align*}
\item The \emph{modified Poincar\'e set} of $G$ is the set
\begin{repequation}{modifiedpoincaredef}
\w\Delta_G = \bigcap_{\rho > 0} \bigcap_{S_\rho} \Delta(S_\rho),
\end{repequation}
where the second intersection is taken over all maximal $\rho$-separated sets $S_\rho \subset G(\zero)$.
\item The number $\w\delta_G = \sup\w\Delta_G$ is called the \emph{modified Poincar\'e exponent} of $G$. If $\w\delta_G\in\w\Delta_G$, we say that $G$ is of \emph{generalized divergence type},\Footnote{We use the adjective ``generalized'' rather than ``modified'' because all groups of convergence/divergence type are also of generalized convergence/divergence type; see Corollary \ref{corollarygeneralizedtypes} below.} while if $\w\delta_G\in\Rplus\butnot\w\Delta_G$, we say that $G$ is of \emph{generalized convergence type}. Note that if $\w\delta_G = \infty$, then $G$ is neither of generalized convergence type nor of generalized divergence type.
\end{itemize}
\end{repdefinition}


We may now state the most powerful version of our Bishop--Jones theorem:

\begin{theorem}[Proven in Chapter \ref{sectionbishopjones}]
\label{theorembishopjonesmodified}
Let $G\prec\Isom(X)$ be a nonelementary semigroup. There exists $\sigma > 0$ such that
\begin{equation}
\label{bishopjonesmodified}
\HD(\Lr) = \HD(\Lur) = \HD(\Lur\cap\Lrsigma) = \w\delta.
\end{equation}
Moreover, for every $0 < s < \w\delta$ there exist $\tau > 0$ and an Ahlfors $s$-regular set $\limitset_s\subset\Lurtau\cap\Lrsigma$.
\end{theorem}

Theorem \ref{theorembishopjonesregular} can be deduced as a corollary of Theorem \ref{theorembishopjonesmodified}; specifically, Propositions \ref{propositionbasicmodified}(ii) and \ref{propositionpoincareregular} show that any group satisfying the hypotheses of Theorem \ref{theorembishopjonesregular} satisfies $\delta = \w\delta$, and hence for such a group \eqref{bishopjonesmodified} implies \eqref{bishopjonesregular}. On the other hand, Proposition \ref{propositionbishopjones} shows that Theorem \ref{theorembishopjonesmodified} applies in many cases where Theorem \ref{theorembishopjonesregular} does not.

We call a group \emph{Poincar\'e regular} if its Poincar\'e exponent $\delta$ and modified Poincar\'e exponent $\w\delta$ are equal. In this language, Proposition \ref{propositionpoincareregular}/Theorem \ref{theorembishopjonesregular} describes sufficient conditions for a group to be Poincar\'e regular, and Proposition \ref{propositionbishopjones} provides a list of examples of groups which are Poincar\'e irregular.

Though Theorem \ref{theorembishopjonesmodified} requires $G$ to be nonelementary, the following corollary does not:

\begin{corollary}
\label{corollaryLrLur}
Fix $G\prec\Isom(X)$. Then for some $\sigma > 0$,
\begin{equation}
\label{LrLur}
\HD(\Lr) = \HD(\Lur) = \HD(\Lur\cap\Lrsigma).
\end{equation}
\end{corollary}
\begin{proof}
If $G$ is nonelementary, then \eqref{LrLur} follows from \eqref{bishopjonesmodified}. On the other hand, if $G$ is elementary, then all three terms of \eqref{LrLur} are equal to zero.
\end{proof}

\bigskip
\section{Examples}
\label{subsectionexamples}

A theory of groups acting on infinite-dimensional space would not be complete without some good ways to construct examples. Techniques used in the finite-dimensional setting, such as arithmetic construction of lattices and Dehn surgery, do not work in infinite dimensions. (The impossibility of constructing lattices in $\Isom(\H^\infty)$ as a direct limit of arithmetic lattices in $\Isom(\H^d)$ is due to known lower bounds on the covolumes of such lattices which blow up as the dimension goes to infinity; see Proposition \ref{propositionarithmeticvolume} below.) Nevertheless, there is a wide variety of groups acting on $\H^\infty$, including many examples of actions which have no analogue in finite dimensions.

\subsection{Schottky products}
The most basic tool for constructing groups or semigroups on hyperbolic metric spaces is the theory of Schottky products. This theory was created by Schottky in 1877 when he considered the Fuchsian group generated by a finite collection of loxodromic isometries $g_i$ described by a disjoint collection of balls $B_{i}^+$ and $B_{i}^-$ with the property that $g_i(\H^2\butnot B_{i}^-) = B_{i}^+$. It was extended further in 1883 by Klein's Ping-Pong Lemma, and used effectively by Patterson \cite{Patterson3} to construct a ``pathological'' example of a Kleinian group of the first kind with arbitrarily small Poincar\'e exponent.

We consider here a quite general formulation of Schottky products: a collection of subsemigroups of $\Isom(X)$ is said to be in \emph{Schottky position} if open sets can be found satisfying the hypotheses of the Ping-Pong lemma whose closure is not equal to $X$ (cf. Definition \ref{definitionschottkyproduct}). This condition is sufficient to guarantee that the product of groups in Schottky position (called a \emph{Schottky product}) is always COT-discrete, but stronger hypotheses are necessary in order to prove stronger forms of discreteness. There is a tension here between hypotheses which are strong enough to prove useful theorems and hypotheses which are weak enough to admit interesting examples. For the purposes of this monograph we make a fairly strong assumption (the \emph{strong separation condition}, Definition \ref{definitionstrongseparation}), one which rules out infinitely generated Schottky groups whose generating regions have an accumulation point (for example, infinitely generated Schottky subgroups of $\Isom(\H^d)$). However, we plan on considering weaker hypotheses in future work \cite{DSU2}.

One theorem of significance in Chapter \ref{sectionschottky} is Theorem \ref{theoremschottkylimitset}, which relates the limit set of a Schottky product to the limit set of its factors together with the image of a Cantor set $\del\Gamma$ under a certain symbolic coding map $\pi:\del\Gamma\to\del X$. As a consequence, we deduce that the properties of compact type and geometrical finiteness are both preserved under finite strongly separated Schottky products (Corollary \ref{corollaryschottkycompacttype} and Proposition \ref{propositionSproductGF}, respectively). A result analogous to Theorem \ref{theoremschottkylimitset} in the setting of infinite alphabet conformal iterated function systems can be found in \cite[Lemma 2.1]{MauldinUrbanski1}.

In \6\ref{subsectionschottkyexamples}, we discuss some (relatively) explicit constructions of Schottky groups, showing that Schottky products are fairly ubiquitous - for example, any two groups which act properly discontinuously at some point of $\del X$ may be rearranged to be in Schottky position, assuming that $X$ is sufficiently symmetric (Proposition \ref{propositionschottkyproduct}).

\subsection{Parabolic groups}

A major point of departure where the theory of subgroups of $\Isom(\H^\infty)$ becomes significantly different from the finite-dimensional theory is in the study of parabolic groups. As a first example, knowing that a group admits a discrete parabolic action on $\Isom(X)$ places strong restrictions on the algebraic properties of the group if $X = \H_\F^d$, but not if $X = \H_\F^\infty$. Concretely, discrete parabolic subgroups of $\Isom(\H_\F^d)$ are always virtually nilpotent (virtually abelian if $\F = \R$), but any group with the Haagerup property admits a parabolic strongly discrete action on $\H^\infty$ (indeed, this is a reformulation of one of the equivalent definitions of the Haagerup property; cf. \cite[p.1, (4)]{CCJJV}). Examples of groups with the Haagerup property include all amenable groups and free groups. Moreover, strongly discrete parabolic subgroups of $\Isom(\H^\infty)$ need not be finitely generated; cf. Example \ref{exampleQparabolic}.

Moving to infinite dimensions changes not only the algebraic but also the geometric properties of parabolic groups. For example, the cyclic group generated by a parabolic isometry may fail to be discrete in any reasonable sense (Example \ref{exampleedelstein}), or it may be discrete in some senses but not others (Example \ref{examplevalette}). The Poincar\'e exponent of a parabolic subgroup of $\Isom(\H_\F^d)$ is always a half-integer \cite[Proof of Lemma 3.5]{CorletteIozzi}, but the situation is much more complicated in infinite dimensions. We prove a general lower bound on the Poincar\'e exponent of a parabolic subgroup of $\Isom(X)$ for any hyperbolic metric space $X$, depending only on the algebraic structure of the group (Theorem \ref{theoremparaboliclowerbound}); in particular, the Poincar\'e exponent of a parabolic action of $\Z^k$ on a hyperbolic metric space is always at least $k/2$. Of course, it is well-known that all parabolic actions of $\Z^k$ on $\H^d$ achieve equality. By contrast, we show that for every $\delta > k/2$ there exists a parabolic action of $\Z^k$ on $\H^\infty$ whose Poincar\'e exponent is equal to $\delta$ (Theorem \ref{theoremnilpotentembedding}).

\subsection{Geometrically finite and convex-cobounded groups}

It has been known for a long time that every finitely generated Fuchsian group has a finite-sided convex fundamental domain (e.g. \cite[Theorem 4.6.1]{Katok_book}). This result does not generalize beyond two dimensions (e.g. \cite{Bers, Jorgensen}), but subgroups of $\Isom(\H^3)$ with finite-sided fundamental domains came to be known as \emph{geometrically finite} groups. Several equivalent definitions of geometrical finiteness in the three-dimensional setting became known, for example Beardon and Maskit's condition that the limit set is the union of the radial limit set $\Lr$ with the set $\Lbp$ of bounded parabolic points \cite{BeardonMaskit}, but the situation in higher dimensions was somewhat murky until Bowditch \cite{Bowditch_geometrical_finiteness} wrote a paper which described which equivalences remain true in higher dimensions, and which do not. The condition of a finite-sided convex fundamental domain is no longer equivalent to any other conditions in higher dimensions (e.g. \cite{Apanasov}), so a higher-dimensional Kleinian group is said to be \emph{geometrically finite} if it satisfies any of Bowditch's five equivalent conditions (GF1)-(GF5).

In infinite dimensions, conditions (GF3)-(GF5) are no longer useful (cf. Remark \ref{remarkbowditchcomparision}), but appropriate generalizations of conditions (GF1) (convex core is equal to a compact set minus a finite number of cusp regions) and (GF2) (the Beardon--Maskit formula $\Lambda = \Lr\cup\Lbp$) are still equivalent for groups of compact type. In fact, (GF1) is equivalent to (GF2) + compact type (Theorem \ref{theoremGFcompact}). We define a group to be \emph{geometrically finite} if it satisfies the appropriate analogue of (GF1) (Definition \ref{definitionGF}). A large class of examples of geometrically finite subgroups of $\Isom(\H^\infty)$ is furnished by combining the techniques of Chapters \ref{sectionschottky} and \ref{sectionparabolic}; specifically, the strongly separated Schottky product of any finite collection of parabolic groups and/or cyclic loxodromic groups is geometrically finite (Corollary \ref{corollaryschottkyGF}).

It remains to answer the question of what can be proven about geometrically finite groups. This is a quite broad question, and in \thispaper\ we content ourselves with proving two theorems. The first theorem, Theorem \ref{theoremGF}, is a generalization of the Milnor--Schwarz lemma \cite[Proposition I.8.19]{BridsonHaefliger} (see also Theorem \ref{theoremCCB}), and describes both the algebra and geometry of a geometrically finite group $G$: firstly, $G$ is generated by a finite subset $F\subset G$ together with a finite collection of parabolic subgroups $G_\xi$ (which are not necessarily finitely generated, e.g. Example \ref{exampleQparabolic}), and secondly, the orbit map $g\mapsto g(\zero)$ is a quasi-isometric embedding from $(G,\dist_G)$ into $X$, where $\dist_G$ is a certain weighted Cayley metric (cf. Example \ref{examplecayleygraph} and \eqref{l0def}) on $G$ whose generating set is $F\cup \bigcup_\xi G_\xi$. As a consequence (Corollary \ref{corollaryGF}), we see that if the groups $G_\xi$, $\xi\in\Lbp$, are all finitely generated, then $G$ is finitely generated, and if these groups have finite Poincar\'e exponent, then $G$ has finite Poincar\'e exponent.

\subsection{Counterexamples}

A significant class of subgroups of $\Isom(\H^\infty)$ that has no finite-dimensional analogue is provided by the \emph{Burger--Iozzi--Monod (BIM) representation theorem} \cite[Theorem 1.1]{BIM}, which states that any unweighted simplicial tree can be equivariantly and quasi-isometrically embedded into an infinite-dimensional real hyperbolic space, with a precise relation between distances in the domain and distances in the range. We call the embeddings provided by their theorem \emph{BIM embeddings}, and the corresponding homomorphisms provided by the equivariance we call \emph{BIM representations}. We generalize the BIM embedding theorem to the case where $X$ is a separable $\R$-tree rather than an unweighted simplicial tree (Theorem \ref{theoremBIM}).

If we have an example of an $\R$-tree $X$ and a subgroup $\Gamma\leq\Isom(X)$ with a certain property, then the image of $\Gamma$ under a BIM representation generally has the same property (Remark \ref{remarkBIM}). Thus, the BIM embedding theorem allows us to translate counterexamples in $\R$-trees into counterexamples in $\H^\infty$. For example, if $\Gamma$ is the free group on two elements acting on its Cayley graph, then the image of $\Gamma$ under a BIM representation provides a counterexample both to an infinite-dimensional analogue of Margulis's lemma (cf. Example \ref{examplemargulis}) and to an infinite-dimensional analogue of I. Kim's theorem regarding length spectra of finite-dimensional algebraic hyperbolic spaces (cf. Remark \ref{remarklengthspectrum}).

Most of the other examples in Chapter \ref{sectionexamples} are concerned with our various notions of discreteness (cf. \6\ref{subsubsectiondiscreteness} above), the notion of Poincar\'e regularity (i.e. whether or not $\delta = \w\delta$), and the relations between them. Specifically, we show that the only relations are the relations which were proven in Chapter \ref{sectiondiscreteness} and Proposition \ref{propositionpoincareregular}, as summarized in Table \ref{figurediscreteness}, p.\pageref{figurediscreteness}. Perhaps the most interesting of the counterexamples we give is Example \ref{exampleAutTBIM}, which is the image under a BIM representation of (a countable dense subgroup of) the automorphism group $\Gamma$ of the $4$-regular unweighted simplicial tree. This example is notable because discreteness properties are not preserved under taking the BIM representation: specifically, $\Gamma$ is weakly discrete but its image under the BIM representation is not. It is also interesting to try to visualize this image geometrically (cf. Figure \ref{figureexampleAutTBIM}).

\subsection{$\R$-trees and their isometry groups}

Motivated by the BIM representation theorem, we discuss some ways of constructing $\R$-trees which admit natural isometric actions. Our first method is the cone construction, in which one starts with an ultrametric space $(Z,\Dist)$ and builds an $\R$-tree $X$ as a ``cone'' over $Z$. This construction first appeared in a paper of F. Choucroun \cite{Choucroun}, although it is similar to several other known cone constructions: \cite[1.8.A.(b)]{Gromov3}, \cite{TrotsenkoVaisala}, \cite[\67]{BonkSchramm}. $\R$-trees constructed by the cone method tend to admit natural parabolic actions, and in Theorem \ref{theoremRtreeorbitalcounting} we provide a necessary and sufficient condition for a function to be the orbital counting function of some parabolic group acting on an $\R$-tree.

Our second method is to staple $\R$-trees together to form a new $\R$-tree. We give sufficient conditions on a graph $(V,E)$, a collection of $\R$-trees $(X_v)_{v\in V}$, and a collection of sets $\set vw\subset X_v$ and bijections $\bij vw:\set vw\to\set wv$ ($(v,w)\in E$) such that stapling the trees $(X_v)_{v\in V}$ along the isometries $(\bij vw)_{(v,w)\in E}$ yields an $\R$-tree (Theorem \ref{theoremstapledunion}). In \6\ref{subsectionstaplingexamples}, we give three examples of the stapling construction, including looking at the cone construction as a special case of the stapling construction. The stapling construction is somewhat similar to a construction of G. Levitt \cite{Levitt}.

\bigskip
\section{Patterson--Sullivan theory}
\label{subsectionconformal}

The connection between the Poincar\'e exponent $\delta$ of a Kleinian group and the geometry of its limit set is not limited to Hausdorff dimension considerations such as those in the Bishop--Jones theorem. As we mentioned before, Patterson and Sullivan's proofs of the equality $\HD(\Lambda) = \delta$ for geometrically finite groups rely on the construction of a certain measure on $\Lambda$, the \emph{Patterson--Sullivan measure}, whose Hausdorff dimension is also equal to $\delta$. In addition to connecting the Poincar\'e exponent and Hausdorff dimension, the Patterson--Sullivan measure also relates to the spectral theory of the Laplacian (e.g. \cite[Theorem 3.1]{Patterson2}, \cite[Proposition 28]{Sullivan_density_at_infinity}) and the geodesic flow on the quotient manifold \cite{Kaimanovich2}. An important property of Patterson--Sullivan measures is \emph{conformality}. Given $s > 0$, a measure $\mu$ on $\del\B^d$ is said to be \emph{$s$-conformal} with respect to a discrete group $G\leq\Isom(\B^d)$ if
\begin{equation}
\label{conformal}
\mu(g(A)) = \int_A |g'(\xi)|^s\;\dee\mu(\xi) \all g\in G \all A\subset \del\B^d.
\end{equation}
The Patterson--Sullivan theorem on the existence of conformal measures may now be stated as follows: For every Kleinian group $G$, there exists a $\delta$-conformal measure on $\Lambda$, where $\delta$ is the Poincar\'e exponent of $G$ and $\Lambda$ is the limit set of $G$.

When dealing with ``coarse'' spaces such as arbitrary hyperbolic metric spaces, it is unreasonable to expect equality in \eqref{conformal}. Thus, a measure $\mu$ on $\del X$ is said to be \emph{$s$-quasiconformal} with respect to a group $G\leq\Isom(X)$ if
\[
\mu(g(A)) \asymp_\times \int_A \overline g'(\xi)^s\;\dee\mu(\xi) \all g\in G \all A\subset \del X.
\]
Here $\overline g'(\xi)$ denotes the upper metric derivative of $g$ at $\xi$; cf. \6\ref{subsubsectionderivativemaps}. We remark that if $X$ is a CAT(-1) space and $G$ is countable, then every quasiconformal measure is coarsely asymptotic to a conformal measure (Proposition \ref{propositionQCtoC}).

In Chapter \ref{sectionconformal}, we describe the theory of conformal and quasiconformal measures in hyperbolic metric spaces. The main theorem is the existence of $\w\delta$-conformal measures for groups of compact type (Theorem \ref{theorempattersonsullivangeneral}). An important special case of this theorem has been proven by Coornaert \cite[Th\'eor\`eme 5.4]{Coornaert} (see also \cite[\61]{BurgerMozes}, \cite[Lemme 2.1.1]{Roblin2}): the case where $X$ is proper and geodesic and $G$ satisfies $\delta < \infty$. The main improvement from Coornaert's theorem to ours is the ability to construct quasiconformal measures for Poincar\'e irregular ($\w\delta < \delta = \infty$) groups; this improvement requires an argument using the class of uniformly continuous functions on $\bord X$.

The big assumption of Theorem \ref{theorempattersonsullivangeneral} is the assumption of compact type. All proofs of the Patterson--Sullivan theorem seem to involve taking a weak-* limit of a sequence of measures in $X$ and then proving that the limit measure is (quasi)conformal, but how can we take a weak-* limit if the limit set is not compact? In fact, Theorem \ref{theorempattersonsullivangeneral} becomes false if you remove the assumption of compact type. In Proposition \ref{propositioncounterexample}, we construct a group acting on an $\R$-tree and satisfying $\delta < \infty$ which admits no $\delta$-conformal measure on its limit set, and then use the BIM embedding theorem (Theorem \ref{theoremBIM}) to get an example in $\H^\infty$.

Surprisingly, it turns out that if we replace the hypothesis of compact type with the hypothesis of \emph{divergence type}, then the theorem becomes true again. Specifically, we have the following:

\begin{theorem}[Proven in Chapter \ref{sectionahlforsthurston}]
\label{theoremahlforsthurstongeneral}
Let $G\leq\Isom(X)$ be a nonelementary group of generalized divergence type (see Definition \ref{definitionmodifiedexponent}). Then there exists a $\w\delta$-quasiconformal measure $\mu$ for $G$ supported on $\Lambda$, where $\w\delta$ is the modified Poincar\'e exponent of $G$. It is unique up to a multiplicative constant in the sense that if $\mu_1,\mu_2$ are two such measures then $\mu_1\asymp_\times \mu_2$ (cf. Remark \ref{remarkmuasymp}). In addition, $\mu$ is ergodic and gives full measure to the radial limit set of $G$.
\end{theorem}


To motivate Theorem \ref{theoremahlforsthurstongeneral}, we recall the connection between the divergence type condition and Patterson--Sullivan theory in finite dimensions. Although the Patterson--Sullivan theorem guarantees the existence of a $\delta$-conformal measure, it does not guarantee its uniqueness. Indeed, the $\delta$-conformal measure is often not unique; see e.g. \cite{AFT}. However, it turns out that the hypothesis of divergence type is enough to guarantee uniqueness. In fact, the condition of divergence type turns out to be quite important in the theory of conformal measures:

\begin{theorem}[Hopf--Tsuji--Sullivan theorem, {\cite[Theorem 8.3.5]{Nicholls}}]
\label{theoremahlforsthurston}
Fix $d\geq 2$, let $G\leq\Isom(\H^d)$ be a discrete group, and let $\delta$ be the Poincar\'e exponent of $G$. Then for any $\delta$-conformal measure $\mu\in\MM(\Lambda)$, the following are equivalent:
\begin{itemize}
\item[(A)] $G$ is of divergence type.
\item[(B)] $\mu$ gives full measure to the radial limit set $\Lr(G)$.
\item[(C)] $G$ acts ergodically on $(\LambdaG,\mu)\times (\LambdaG,\mu)$.
\end{itemize}
In particular, if $G$ is of divergence type, then every $\delta$-conformal measure is ergodic, so there is exactly one (ergodic) $\delta$-conformal probability measure.
\end{theorem}
We remark that our sentence ``In particular \ldots'' stated in theorem above was not included in \cite[Theorem 8.3.5]{Nicholls} but it is well-known and follows easily from the equivalence of (A) and (C).

\begin{remark}
Theorem \ref{theoremahlforsthurston} has a long history. The equivalence (B) \iff (C) was first proven by E. Hopf in the case $\delta = d - 1$\Footnote{In this paragraph, when we say that someone proves the case $\delta = d - 1$, we mean that they considered the case where $\mu$ is Hausdorff $(d - 1)$-dimensional measure on $S^{d - 1}$.} \cite{Hopf_fuchsian, Hopf_statistik} (1936, 1939). The equivalence (A) \iff (B) was proven by Z. Y\^uj\^ob\^o in the case $\delta = d - 1 = 1$ \cite{Yujobo} (1949), following an incorrect proof by M. Tsuji \cite{Tsuji_multiply_connected} (1944).\Footnote{See \cite[p.484]{Sullivan_ergodic_theory_at_infinity} for some further historical remarks on the case $\delta = d - 1 = 1$.} Sullivan proved (A) \iff (C) in the case $\delta = d - 1$ \cite[Theorem II]{Sullivan_ergodic_theory_at_infinity}, then generalized this equivalence to the case $\delta > (d - 1)/2$ \cite[Theorem 32]{Sullivan_density_at_infinity}. He also proved (B) \iff (C) in full generality \cite[Theorem 21]{Sullivan_density_at_infinity}. Next, W. P. Thurston gave a simpler proof of (A) \implies (B)\Footnote{By this point, it was considered obvious that (B) \implies (A).} in the case $\delta = d - 1$ \cite[Theorem 4 of Section VII]{Ahlfors}. P. J. Nicholls finished the proof by showing (A) \iff (B) in full generality \cite[Theorems 8.2.2 and 8.2.3]{Nicholls}. Later S. Hong re-proved (A) \implies (B) in full generality twice in two independent papers \cite{Hong1,Hong2}, apparently unaware of any previous results. Another proof of (A) \implies (B) in full generality, which was conceptually similar to Thurston's proof, was given by P. Tukia \cite[Theorem 3A]{Tukia}. Further generalization was made by C. Yue \cite{Yue} to negatively curved manifolds, and by T. Roblin \cite[Th\'eor\`eme 1.7]{Roblin1} to proper CAT(-1) spaces.
\end{remark}

Having stated the Hopf--Tsuji--Sullivan theorem, we can now describe why Theorem \ref{theoremahlforsthurstongeneral} is true, first on an intuitive level and then giving a sketch of the real proof. On an intuitive level, the fact that divergence type implies both ``existence and uniqueness'' of the $\delta$-conformal measure in finite dimensions indicates that perhaps the compactness assumption is not needed -- the sequence of measures used to construct the Patterson--Sullivan measure converges already, so it should not be necessary to use compactness to take a convergent subsequence.

The real proof involves taking the Samuel--Smirnov compactification of $\bord X$, considered as a metric space with respect to a visual metric (cf. \6\ref{subsectionvisual}). The Samuel--Smirnov compactification of a metric space (cf. \cite[\67]{NaimpallyWarrack}) is conceptually similar to the more familiar Stone--\v Cech compactification, except that only uniformly continuous functions on the metric space extend to continuous functions on the compactification, not all continuous functions. If we used the Stone--\v Cech compactification rather than the Samuel--Smirnov compactification, then our proof would only apply to groups with finite Poincar\'e exponent; cf. Remark \ref{remarkstonecech} and Remark \ref{remarksamuelsmirnov}.


\begin{proof}[Sketch of the proof of Theorem \ref{theoremahlforsthurstongeneral}]
We denote the Samuel--Smirnov compactification of $\bord X$ by $\what X$. By a nonstandard analogue of Theorem \ref{theorempattersonsullivangeneral} (viz. Lemma \ref{lemmapattersonsullivannonstandard}), there exists a $\w\delta$-quasiconformal measure $\what\mu$ on $\what{\del X}$. By a generalization of Theorem \ref{theoremahlforsthurston} (viz. Proposition \ref{propositionahlforsthurston}), $\what\mu$ gives full measure to the radial limit set $\what\Lr$. But a simple computation (Lemma \ref{lemmaradialstandard}) shows that $\what\Lr = \Lr$, demonstrating that $\what\mu\in\MM(\Lambda)$.
\end{proof}

\subsection{Quasiconformal measures of geometrically finite groups}
Let us consider a geometrically finite group $G\leq\Isom(X)$ with Poincar\'e exponent $\delta < \infty$, and let $\mu$ be a $\delta$-quasiconformal measure on $\Lambda$. 
Such a measure exists since geometrically finite groups are of compact type (Theorem \ref{theoremGFcompact} and Theorem \ref{theorempattersonsullivangeneral}), and is unique as long as $G$ is of divergence type (Corollary \ref{corollaryahlforsthurston}). When $X = \H^d$, the geometry of $\mu$ is described by the Global Measure Formula \cite[Theorem on p.271]{Sullivan_entropy}, \cite[Theorem 2]{StratmannVelani}: the measure of a ball $B(\eta,e^{-t})$ is coarsely asymptotic to $e^{-\delta t}$ times a factor depending on the location of the point $\eta_t := \geo\zero\eta_t$ in the quotient manifold $\H^d/G$. Here $\geo\zero\eta_t$ is the unique point on the geodesic connecting $\zero$ and $\eta$ with distance $t$ from $\zero$; cf. Notations \ref{notationgeopqt}, \ref{notationgeopqt2}.

In a general hyperbolic metric space $X$ (indeed, already for $X = \H^\infty$), one cannot get a precise asymptotic for $\mu(B(\eta,e^{-t}))$, due to the fact that the measure $\mu$ may fail to be doubling (Example \ref{exampledoublingnondoubling}). Instead, our version of the global measure formula gives both an upper bound and a lower bound for $\mu(B(\eta,e^{-t}))$. Specifically, we define a function $m:\Lambda\times\Rplus\to(0,\infty)$ (for details see \eqref{metat}) and then show:
\begin{theorem}[Global measure formula, Theorem \ref{theoremglobalmeasure}; proven in Section \ref{subsectionglobalmeasureproof}]
\label{theoremglobalmeasureintro}
For all $\eta\in\Lambda$ and $t > 0$,
\begin{equation}
\label{globalmeasureformulaintro}
m(\eta,t + \sigma) \lesssim_\times \mu(B(\eta,e^{-t})) \lesssim_\times m(\eta,t - \sigma),
\end{equation}
where $\sigma > 0$ is independent of $\eta$ and $t$.
\end{theorem}

It is natural to ask for which groups \eqref{globalmeasureformulaintro} can be improved to an exact asymptotic, i.e. for which groups $\mu$ is doubling. We address this question in Section \ref{subsectiondoubling}, proving a general result (Proposition \ref{propositiondoubling2}), a special case of which is that if $X$ is a finite-dimensional algebraic hyperbolic space, then $\mu$ is doubling (Example \ref{exampledoubling}). Nevertheless, there are large classes of examples of groups $G\leq\Isom(\H^\infty)$ for which $\mu$ is not doubling (Example \ref{exampledoublingnondoubling}), illustrating once more the wide difference between $\H^\infty$ and its finite-dimensional counterparts.

It is also natural to ask about the implications of the Global Measure Formula for the dimension theory of the measure $\mu$. For example, when $X = \H^d$, the Global Measure Formula was used to show that $\HD(\mu) = \delta$ \cite[Proposition 4.10]{StratmannVelani}. In our case we have: 
\begin{theorem}[Cf. Theorem \ref{theoremhlogseries}]
\label{theoremhlogseriesintro}
If for all $\bp\in P$, the series
\begin{equation}
\label{hlogseriesintro}
\sum_{h\in G_\bp} e^{-\delta\dogo h} \dogo h
\end{equation}
converges, then $\mu$ is exact dimensional (cf. Definition \ref{definitionexactdimensional}) of dimension $\delta$. In particular, 
\[
\HD(\mu) = \PD(\mu) = \delta ~.
\]
\end{theorem}
The hypothesis that \eqref{hlogseriesintro} converges is a very non-restrictive hypothesis. For example, it is satisfied whenever $\delta > \delta_\bp$ for all $\bp\in P$ (Corollary \ref{corollarydeltagtrdeltap}). Combining with Proposition \ref{propositiondeltagtrdeltap} shows that any counterexample must satisfy
\[
\sum_{h\in G_\bp} e^{-\delta\dogo h} < \infty = \sum_{h\in G_\bp} e^{-\delta\dogo h} \dogo h
\]
for some $\bp\in P$, creating a very narrow window for the orbital counting function $\NN_\bp$ (cf. Notation \ref{notationmetat}) to lie in. Nevertheless, we show that there exist counterexamples (Example \ref{examplenotexactdim}) for which the series \eqref{hlogseriesintro} diverges. After making some simplifying assumptions, we are able to prove (Theorem \ref{theoremhlogseriesconverse}) that the Patterson--Sullivan measures of groups for which \eqref{hlogseriesintro} diverges cannot be exact dimensional, and in fact satisfy $\HD(\mu) = 0$.

There is a relation between exact dimensionality of the Patterson--Sullivan measure and the theory of Diophantine approximation on the boundary of $\del X$, as described in \cite{FSU4}. Specifically, if $\VWA_\xi$ denotes the set of points which are very well approximable with respect to a distinguished point $\xi$ (cf. \6\ref{subsubsectiondiophantineapproximation}), then we have the following:

\begin{theorem}[Cf. Theorem \ref{theoremequivalentVWA}]
\label{theoremequivalentVWAintro}
The following are equivalent:
\begin{itemize}
\item[(A)] For all $p \in P$, $\mu(\VWA_\bp) = 0$.
\item[(B)] $\mu$ is exact dimensional.
\item[(C)] $\HD(\mu) = \delta$.
\item[(D)] For all $\xi \in \Lambda$, $\mu(\VWA_\xi) = 0$.
\end{itemize}
\end{theorem}

In particular, combining with Theorem \ref{theoremhlogseriesintro} demonstrates that the equation 
\[
\mu(\VWA_\xi) = 0
\] 
holds for a large class of geometrically finite groups $G$ and for all $\xi\in\Lambda$. This improves the results of \cite[\61.5.3]{FSU4}.

\section{Appendices}

We conclude this monograph with two appendices. Appendix \ref{appendixopenproblems} contains a list of open problems, and Appendix \ref{appendixindex} an index of defined terms.



\mainmatter
\part{Preliminaries}
\label{partpreliminaries}
This part will be divided as follows: In Chapter \ref{sectionROSSONCTs} we define the class of algebraic hyperbolic spaces, which are often called rank one symmetric spaces of noncompact type. In Chapters \ref{sectiongeometry1} and \ref{sectiongeometry2}, we define the class of hyperbolic metric spaces and study their geometry. In Chapter \ref{sectiondiscreteness}, we explore different notions of discreteness for groups of isometries of a metric space. In Chapter \ref{sectionclassification} we prove two classification theorems, one for isometries (Theorem \ref{classificationofisometries}) and one for semigroups of isometries (Theorem \ref{classificationofsemigroups}). Finally, in Chapter \ref{sectionlimitsets} we define and study the limit set of a semigroup of isometries.


\chapter{Algebraic hyperbolic spaces}
\label{sectionROSSONCTs}

In this chapter we introduce our main objects of interest, algebraic hyperbolic spaces in finite and infinite dimensions. References for the theory of finite-dimensional algebraic hyperbolic spaces, which are often called \emph{rank one symmetric spaces of noncompact type}, include \cite{BridsonHaefliger,CFKP,MackayTyson}. Infinite-dimensional algebraic hyperbolic spaces, as well as some non-hyperbolic infinite-dimensional symmetric spaces, have been discussed in \cite{Duchesne}.

\bigskip
\section{The definition}
\label{subsectionROSSONCT}

Finite-dimensional rank one symmetric spaces of noncompact type come in four flavors, corresponding to the classical division algebras $\R$, $\C$, $\Q$ (quaternions), and $\mathbb O$ (octonions).\Footnote{We denote the quaternions by $\Q$ in order to avoid confusion with the rank one symmetric space of noncompact type (defined over $\Q$) itself, which we will denote by $\H$. Be aware that $\Q$ should not be confused with the set of rational numbers.} The first three division algebras have corresponding rank one symmetric spaces of noncompact type of arbitrary dimension, but there is only one rank one symmetric space of noncompact type corresponding to the octonions; it occurs in dimension two (which corresponds to real dimension 16). Consequently, the octonion rank one symmetric space of noncompact type (known as the \emph{Cayley hyperbolic plane}\Footnote{Not to be confused with the \emph{Cayley plane}, a different mathematical object.}) does not have an infinite-dimensional analogue, while the other three classes do admit infinite-dimensional analogues.

The rank one symmetric spaces of noncompact type corresponding to $\R$ have constant negative curvature. However, those corresponding to the other division algebras have variable negative curvature \cite[Lemmas 2.3, 2.7, 2.11]{Quint} (see also \cite[Corollary of Proposition 4]{Heintze}).

\begin{remark}
\label{remarkH2O}
In this monograph we will use the term ``algebraic hyperbolic spaces'' to refer to all rank one symmetric spaces of noncompact type \emph{except} the Cayley hyperbolic plane $\H_{\mathbb O}^2$, in order to avoid dealing with the complicated algebra of the octonions.\Footnote{The complications come from the fact that the octonions are not associative, thus making it somewhat unclear what it means to say that $\amsbb O^3$ is a vector space ``over'' the octonions, since in general $(\xx a) b \neq \xx (a b)$.} However, we feel confident that all the theorems regarding algebraic hyperbolic spaces in this monograph can be generalized to the Cayley hyperbolic plane (possibly after modifying the statements slightly). We leave this task to an algebraist.

For the reader interested in learning more about the Cayley hyperbolic plane, see \cite[pp.136-139]{Mostow}, \cite{SpringerVeldkamp}, or \cite{Allcock}; see also \cite{Baez} for an excellent introduction to octonions in general.
\end{remark}

Fix $\F\in\{\R,\C,\Q\}$ and an index set $J$, and let us construct an algebraic hyperbolic space (i.e. a rank one symmetric space of noncompact type) over the field $\F$ in dimension $\#(J)$. We remark that usually we will let $J = \Namer = \{1,2,\ldots\}$, but occasionally $J$ may be an uncountable set. Let
\[
\HH = \HH_\F^J := \left\{ \xx = (x_i)_{i\in J}\in\F^J \left\vert\;\; \sum_{i\in J} |x_i|^2 < \infty
\right.\right\},
\]
and for $\xx\in\HH$ let
\[
\|\xx\| := \left( \sum_{i\in J} |x_i|^2 \right)^{1/2}.
\]
We will think of $\HH$ as a right $\F$-module, so scalars will always act on the right.\Footnote{The advantage of this convention is that it allows operators to act on the left.} Note that
\[
\|\xx a\| = |a|\cdot \|\xx\| \all \xx\in \HH \all a\in\F.
\]
A \emph{sesquilinear form on $\HH$} is an $\R$-bilinear map $B(\cdot,\cdot):\HH\times\HH\to \F$ satisfying
\[
B(\xx a,\yy) = \wbar a B( \xx,\yy) \text{ and } B( \xx,\yy a) = B( \xx,\yy) a.\Footnote{In the case $\amsbb F = \amsbb C$, this disagrees with the usual convention; we follow here the convention of \cite[\63.3.1]{MackayTyson}.}
\]
Here and from now on $\wbar a$ denotes the conjugate of a complex or quaternionic number $a\in\F$; if $\F = \R$, then $\wbar a = a$.

A sesquilinear form is said to be \emph{skew-symmetric} if $B(\yy,\xx) = \wbar{B(\xx,\yy)}$. For example, the map
\[
B_\EE(\xx,\yy) := \sum_{i\in J} \wbar{x_i} y_i
\]
is a skew-symmetric sesquilinear form. Note that
\[
\EE(\xx) := B_\EE(\xx,\xx) = \|\xx\|^2.
\]

\section{The hyperboloid model}
\label{subsectionhyperboloidmodel}
Assume that $0\notin J$, and let
\[
\LL = \LL_\F^{J\cup\{0\}} = \HH_\F^{J\cup\{0\}} = \left\{\xx = (x_i)_{i\in J\cup\{0\}}\in\F^{J\cup\{0\}} \left\vert\sum_{i\in J\cup\{0\}} |x_i|^2 < \infty\right.\right\}.
\]
Consider the skew-symmetric sesquilinear form $B_\QQ:\LL\times\LL\to\F$ defined by
\[
B_\QQ(\xx,\yy) := -\wbar x_0 y_0 + \sum_{i\in J} \wbar x_i y_i
\]
and its associated quadratic form
\begin{equation}
\label{Qdef}
\QQ(\xx) := B_\QQ(\xx,\xx) = -|x_0|^2 + \sum_{i\in J} |x_i|^2.
\end{equation}
We observe that the form $\QQ$ is not positive definite, since $\QQ(\ee_0) = -1$.

\begin{remark}
If $\F = \R$, then the form $\QQ$ is called a \emph{Lorentzian} quadratic form, and the pair $(\LL,\QQ)$ is called a \emph{Minkowski space}.
\end{remark}

Let $\proj(\LL)$ denote the \emph{projectivization} of $\LL$, i.e. the quotient of $\LL\butnot\{\0\}$ under the equivalence relation $\xx\sim \xx a$ ($\xx\in\LL\butnot\{\0\}$, $a\in \F\butnot\{0\}$). Let
\[
\H = \H_\F^J := \{ [\xx]\in\proj(\LL_\F^{J\cup\{0\}}) : \QQ(\xx) < 0 \},
\]
and consider the map $\dist_\H: \H\times\H\to \Rplus$ defined by the equation
\begin{equation}
\label{distanceinL}
\cosh\dist_\H([\xx],[\yy]) = \frac{|B_\QQ(\xx,\yy)|}{\sqrt{|\QQ(\xx)|\cdot|\QQ(\yy)}|}, \;\;\; [\xx],[\yy]\in \H.
\end{equation}
Note that the map $\dist_\H$ is well-defined because the right hand side is invariant under multiplying $\xx$ and $\yy$ by scalars.

\begin{proposition}
\label{propositionROSSONCT}
The map $\dist_\H$ is a metric on $\H$ that is compatible with the natural topology, when viewed as a subspace of the quotient space $\proj(\LL)$. Moreover, for any two distinct points $[\xx],[\yy]\in\HH$ there exists a unique isometric embedding $\gamma:\R\to \H$ such that $\gamma(0) = [\xx]$ and $\gamma\circ\dist_\H([\xx],[\yy]) = [\yy]$.
\end{proposition}
\begin{remark}
The second sentence is the \emph{unique geodesic extension} property of $\H$. It holds more generally for Riemannian manifolds (cf. Remark \ref{remarkriemannian} below), but is an important distinguishing feature in the larger class of uniquely geodesic metric spaces.
\end{remark}
\begin{proof}[Proof of Proposition \ref{propositionROSSONCT}]
The key to the proof is the following lemma, which may also be deduced from the infinite-dimensional analogue of Sylvester's law of inertia \cite[Lemma 3]{Maddocks}.
\begin{lemma}
\label{lemmasylvester}
Fix $\zz\in\LL$ with $\QQ(\zz) < 0$, and let $\zz^\perp = \{\ww : B_\QQ(\zz,\ww) = 0\}$. Then $\QQ\given\zz^\perp$ is positive definite.
\end{lemma}
\begin{subproof}[Proof of Lemma \ref{lemmasylvester}]
By contradiction, suppose $\QQ(\yy) \leq 0$ for some $\yy\in\zz^\perp$. There exist $a,b\in\F$, not both zero, such that $y_0 a + z_0 b = 0$. But then
\[
0 < \QQ(\yy a + \zz b) = |a|^2\QQ(\yy) + |b|^2\QQ(\zz) \leq 0,
\]
which provides a contradiction.
\end{subproof}

Now fix $[\xx],[\yy],[\zz]\in\H$, and let $\xx,\yy,\zz\in\LL\butnot\{\0\}$ be representatives which satisfy
\[
B_\QQ(\xx,\zz) = B_\QQ(\yy,\zz) = \QQ(\zz) = -1.
\]
Then
\begin{align*}
\cosh\dist_\H([\xx],[\zz]) &= \frac{1}{\sqrt{1 - \QQ(\xx - \zz)}} &
\cosh\dist_\H([\yy],[\zz]) &= \frac{1}{\sqrt{1 - \QQ(\yy - \zz)}}\\
\sinh\dist_\H([\xx],[\zz]) &= \frac{\sqrt{\QQ(\xx - \zz)}}{\sqrt{1 - \QQ(\xx - \zz)}} &
\sinh\dist_\H([\yy],[\zz]) &= \frac{\sqrt{\QQ(\yy - \zz)}}{\sqrt{1 - \QQ(\yy - \zz)}}\cdot
\end{align*}
By the addition law for hyperbolic cosine we have
\[
\cosh(\dist_\H([\xx],[\zz]) + \dist_\H([\yy],[\zz])) = \frac{1 + \sqrt{\QQ(\xx - \zz)}\sqrt{\QQ(\yy - \zz)}}{\sqrt{1 - \QQ(\xx - \zz)}\sqrt{1 - \QQ(\yy - \zz)}} \cdot
\]
On the other hand, we have
\[
\cosh\dist_\H([\xx],[\yy]) = \frac{1}{\sqrt{1 - \QQ(\xx - \zz)}}\frac{1}{\sqrt{1 - \QQ(\yy - \zz)}}\left|-1 + B_\QQ(\xx - \zz,\yy - \zz)\right|.
\]
Since $\xx - \zz,\yy - \zz\in\zz^\perp$, the Cauchy--Schwartz inequality together with Lemma \ref{lemmasylvester} gives
\[
\left|-1 + B_\QQ(\xx - \zz,\yy - \zz)\right| \leq 1 + \sqrt{\QQ(\xx - \zz)}\sqrt{\QQ(\yy - \zz)},
\]
with equality if and only if $\xx - \zz$ and $\yy - \zz$ are proportional with a negative real constant of proportionality. This demonstrates the triangle inequality.

To show that $\dist_\H$ is compatible with the natural topology, it suffices to show that if $U$ is a neighborhood in the natural topology of a point $[\xx]\in\H$, then there exists $\epsilon > 0$ such that $B([\xx],\epsilon) \subset U$. Indeed, fix a representative $\xx\in [\xx]$; then there exists $\delta > 0$ such that $\|\yy - \xx\|\leq \delta$ implies $[\yy]\in U$. Now, given $[\yy]\in B([\xx],\epsilon)$, choose a representative $\yy\in [\yy]$ such that $\zz := \yy - \xx$ satisfies $B_\QQ(\xx,\zz) = 0$; this is possible since any representative $\yy\in[\yy]$ satisfies $B_\QQ(\xx,\yy)\neq 0$ by Lemma \ref{lemmasylvester}. Then
\[
\cosh\dist_\H([\xx],[\yy]) = \frac{|\QQ(\xx)|}{\sqrt{|\QQ(\xx)| \cdot |\QQ(\xx) + \QQ(\zz)|}} = \sqrt{\left|\frac{\QQ(\xx)}{\QQ(\xx) + \QQ(\zz)}\right|}\cdot
\]
So if $\dist_\H([\xx],[\yy]) \leq \epsilon$, then $\QQ(\zz)\leq \QQ(\xx)[1 - 1/\cosh(\epsilon)]$. By Lemma \ref{lemmasylvester}, there exists $C > 0$ such that $\|\zz\|^2 \leq C\QQ(\xx)[1 - 1/\cosh(\epsilon)]$. In particular, we may choose $\epsilon$ so that $C\QQ(\xx)[1 - 1/\cosh(\epsilon)] \leq \delta$, which completes the proof that $\dist_\H$ is compatible with the natural topology.

Now suppose that $\gamma:\R\to \H$ is an isometric embedding, and let $[\zz] = \gamma(0)$. Choose a representative $\zz\in\LL\butnot\{\0\}$ such that $\QQ(\zz) = -1$, and for each $t\in\R\butnot\{0\}$ choose a representative $\xx_t\in \LL\butnot\{\0\}$ such that $B_\QQ(\xx_t,\zz) = -1$. The preceding argument shows that for $t_1 < 0 < t_2$, $\xx_{t_1} - \zz$ and $\xx_{t_2} - \zz$ are proportional with a negative constant of proportionality. Together with \eqref{distanceinL}, this implies that
\begin{equation}
\label{geodesic}
\xx_t = \zz + \tanh(t)\ww
\end{equation}
for some $\ww\in \zz^\perp$ with $\QQ(\ww) = 1$. Conversely, direct calculation shows that the equation \eqref{geodesic} defines an isometric embedding $\gamma_{\zz,\ww}:\R\to\H$ via the formula $\gamma_{\zz,\ww}(t) = [\xx_t]$.
\end{proof}

\begin{definition}
\label{definitionROSSONCT}
An \emph{algebraic hyperbolic space} is a pair $(\H_\F^J,\dist_\H)$, where $\F\in\{\R,\C,\Q\}$ and $J$ is a nonempty set such that $0\notin J$. 
\end{definition}

\begin{remark}
\label{remarkROSSONCTterminology}
In finite dimensions, the class of algebraic hyperbolic spaces is identical (modulo the Cayley hyperbolic plane, cf. Remark \ref{remarkH2O}) to the class of rank one symmetric spaces of noncompact type. This follows from the classification theorem for finite-dimensional symmetric spaces, see e.g. \cite[p.518]{Helgason}.\Footnote{In the notation of \cite{Helgason}, the spaces $\amsbb H_{\amsbb R}^p$, $\amsbb H_{\amsbb C}^p$, $\amsbb H_{\amsbb Q}^p$, and $\amsbb H_{\amsbb O}^2$ are written as $\SO(p,1)/\SO(p)$, $\SU(p,1)/\SU(p)$, $\Sp(p,1)/\Sp(p)$, and $(\mathfrak f_{4(-20)},\so(9))$, respectively.} It is not clear whether an analogous theorem holds in infinite dimensions (but see \cite{Duchesne3} for some results in this direction).
\end{remark}

%

\begin{remark}
\label{remarkriemannian}
In finite dimensions, the metric $\dist_\H$ may be defined as the length metric associated to a certain Riemannian metric on $\H$; cf. \cite[\62.2]{Quint}. The same procedure works in infinite dimensions; cf. \cite{Lang_differential_geometry} for an exposition of the theory of infinite dimensional manifolds. Although a detailed account of the theory of infinite-dimensional Riemannian manifolds would be too much of a digression, let us make the following points:
\begin{itemize}
\item An infinite-dimensional analogue of the Hopf--Rinow theorem is false \cite{Atkin}, i.e. there exists an infinite-dimensional Riemannian manifold such that some two points on that manifold cannot be connected by a geodesic. However, if an infinite-dimensional Riemannian manifold $X$ is nonpositively curved, then any two points of $X$ can be connected by a unique geodesic as a result of the infinite-dimensional Cartan--Hadamard theorem \cite[IX, Theorem 3.8]{Lang_differential_geometry}; moreover, this geodesic is length-minimizing. In particular, if one takes a Riemannian manifolds approach to defining infinite-dimensional algebraic hyperbolic spaces, then the second assertion of Proposition \ref{propositionROSSONCT} follows from the Cartan--Hadamard theorem.
\item A bijection between two infinite-dimensional Riemannian manifolds is an isometry with respect to the length metric if and only if it is a diffeomorphism which induces an isometry on the Riemannian metric \cite[Theorem 7]{GJR}. This theorem is commonly known as the Myers--Steenrod theorem, as S. B. Myers and N. E. Steenrod proved its finite-dimensional version \cite{MyersSteenrod}. The difficult part of this theorem is proving that any bijection which is an isometry with respect to the length metric is differentiable. In the case of algebraic hyperbolic spaces, however, this follows directly from Theorem \ref{theoremisometries} below.
\end{itemize}
\end{remark}

\bigskip
\section{Isometries of algebraic hyperbolic spaces}
\label{subsectionisometries}

We define the \emph{group of isometries} of a metric space $(X,\dist)$ to be the group
\[
\Isom(X) := \{g:X\to X: \text{$g$ is a bijection and } \dist(g(x),g(y)) = \dist(x,y)\all x,y\in X\}.
\]
In this section we will compute the group of isometries of an arbitrary algebraic hyperbolic space. Fix $\F\in\{\R,\C,\Q\}$ and an index set $J$, and let $\HH = \HH_\F^J$, $\LL = \LL_\F^{J\cup\{0\}}$, and $\H = \H_\F^J$. We begin with the following observation:

\begin{observation}
Let $\O_\F(\LL;\QQ)$ denote the group of $\QQ$-preserving $\F$-linear automorphisms of $\LL$. Then for all $T\in\O_\F(\LL;\QQ)$, the map $[T]:\H\to\H$ defined by the equation
\begin{equation}
\label{[T]def}
[T]([\xx]) = [T\xx]
\end{equation}
is an isometry of $(\H,\dist_\H)$.
\end{observation}
\begin{proof}
The map $[T]$ is well-defined by the associativity property $T(\xx a) = (T\xx) a$. Since $T$ is $\QQ$-preserving and $\F$-linear, the polarization identity (the three versions cover the three cases when the base field $\F = \R$, $\C$, and $\Q$ respectively)
\[
B_\QQ(\xx,\yy) = \begin{cases}
\frac 14 [\QQ(\xx + \yy) - \QQ(\xx - \yy)] \\
\frac 14 [\QQ(\xx + \yy) - \QQ(\xx - \yy) - i \QQ(\xx + \yy i) + i\QQ(\xx - \yy i)] \\
\frac 14 \left[\QQ(\xx + \yy) - \QQ(\xx - \yy) + \sum\limits_{\ell = i,j,k}\big( -\ell \QQ(\xx + \yy \ell) + \ell\QQ(\xx - \yy \ell) \big)\right] 
\end{cases}
\] 
guarantees that
\begin{equation}
\label{BQpreserving}
B_\QQ(T\xx,T\yy) = B_\QQ(\xx,\yy) \all \xx,\yy\in\HH.
\end{equation}
Comparing with \eqref{distanceinL} shows that $[T]$ is an isometry.
\end{proof}

The group $\O_\F(\LL;\QQ)$ is quite large. In addition to containing all maps of the form $T\oplus I$, where $T\in\O_\F(\HH;\EE)$ and $I:\F\to\F$ is the identity map, it also contains the so-called \emph{Lorentz boosts}
\begin{equation}
\label{lorentzboost}
T_{j,t}(\xx) = \left(\begin{cases}
x_i & i\neq 0,j\\
\cosh(t) x_j + \sinh(t) x_0 & i = j\\
\sinh(t) x_j + \cosh(t) x_0 & i = 0
\end{cases}\right)_{i\in J\cup\{0\}}, \;\; j\in J, t\in\R.
\end{equation}
We leave it as an exercise that $\O_\F(\HH;\EE)\oplus\{I\}$ and the Lorentz boosts in fact generate the group $\O_\F(\LL;\QQ)$.

\begin{observation}
\label{observationtransitivity}
The group
\[
\PO_\F(\LL;\QQ) = \{[T] : T\in\O_\F(\LL;\QQ)\} \leq \Isom(\H)
\]
acts transitively on $\H$.
\end{observation}
\begin{proof}
Let $\zero = [(1,\0)]$. The orbit of $\zero$ under $\PO_\F(\LL;\QQ)$ contains its image under the Lorentz boosts. Specifically, for every $t\in\R$ the orbit of $\zero$ contains the point $[(\cosh(t),\sinh(t),\0)]$. Applying maps of the form $[T\oplus I]$, $T\in \O_\F(\HH,\EE)$, shows that the orbit of $\zero$ is $\H$.
\end{proof}

We may ask the question of whether the group $\PO_\F(\LL;\QQ)$ is equal to $\Isom(\H)$ or is merely a subgroup. The answer turns out to depend on the division algebra $\F$:

\begin{theorem}
\label{theoremisometries}
If $\F\in\{\R,\Q\}$ then $\Isom(\H) = \PO_\F(\LL;\QQ)$. If $\F = \C$, then $\PO_\F(\LL;\QQ)$ is of index 2 in $\Isom(\H)$.
\end{theorem}
\begin{remark}
In finite dimensions, Theorem \ref{theoremisometries} is given as an exercise in \cite[Exercise II.10.21]{BridsonHaefliger}. Because of the importance of Theorem \ref{theoremisometries} to \thispaper, we provide a full proof.
\end{remark}
Before proving Theorem \ref{theoremisometries}, it will be convenient for us to introduce a group somewhat larger than $\O_\F(\LL;\QQ)$. Let $\Aut(\F)$ denote the group of automorphisms of $\F$ as an $\R$-algebra, i.e.
\[
\Aut(\F) = \left\{\sigma:\F\to\F \left|
\begin{split}
\text{$\sigma$ is an $\R$-linear bijection and}\\
\text{$\sigma(ab) = \sigma(a)\sigma(b)$ for all $a,b\in\F$}
\end{split}
\right.
\right\}.
\]
We will say that an $\R$-linear map $T:\LL\to\LL$ is \emph{$\F$-skew linear} if there exists $\sigma\in\Aut(\F)$ such that
\begin{equation}
\label{sigmaTdef}
T(\xx a) = T(\xx) \sigma(a) \text{ for all $\xx\in\HH$ and $a\in\F$}.
\end{equation}
The group of skew-linear bijections $T:\LL\to\LL$ which preserve $\QQ$ will be denoted $\O_\F^*(\LL;\QQ)$. For each $T$, the unique $\sigma\in\Aut(\F)$ satisfying \eqref{sigmaTdef} will be denoted $\sigma_T$. Note that the map $T\mapsto\sigma_T$ is a homomorphism.

\begin{warning*}
The associative law $(T\xx)a = T(\xx a)$ is \emph{not} valid for $T\in\O_\F^*(\LL;\QQ)$; rather, $T(\xx a) = (T\xx)\sigma_T(a)$ by \eqref{sigmaTdef}. Thus when discussing elements of $\O_\F^*(\LL;\QQ)$, we must be careful of parentheses.
\end{warning*}

\begin{example}
For each $\sigma\in\Aut(\F)$, the map
\[
\sigma^J(\xx) = \left(\sigma(x_i)\right)_{i\in J}
\]
is $\F$-skew-linear and $\QQ$-preserving, and $\sigma_{\sigma^J} = \sigma$.
\end{example}

\begin{observation}
\label{observationBQpreserving}
For $T\in\O_\F^*(\LL;\QQ)$,
\[
B_\QQ(T\xx,T\yy) = \sigma_T(B_\QQ(\xx,\yy)) \all \xx,\yy\in\LL.
\]
\end{observation}
\begin{proof}
By \eqref{BQpreserving}, the formula holds when $T\in\O_\F(\LL;\QQ)$, and direct calculation shows that it holds when $T = \sigma^J$ for some $\sigma\in\Aut(\F)$. Since $\O_\F^*(\LL;\QQ)$ is a semidirect product of the groups $\O_\F(\LL;\QQ)$ and $\{\sigma^J : \sigma\in\Aut(\F)\}$, this completes the proof.
\end{proof}

We observe that if $T\in\O_\F^*(\LL;\QQ)$, then \eqref{sigmaTdef} shows that $T$ preserves $\F$-lines, i.e. $T(\xx\F) = T(\xx)\F$ for all $\xx\in\LL\butnot\{\0\}$. Thus the equation \eqref{[T]def} defines a map $[T]:\H\to\H$, which is an isometry by Observation \ref{observationBQpreserving}. Thus if
\[
\PO_\F^*(\LL;\QQ) = \{[T] : T\in\O_\F^*(\LL;\QQ)\},
\]
then
\[
\PO_\F(\LL;\QQ) \leq \PO_\F^*(\LL;\QQ) \leq \Isom(\H).
\]
We are now ready to begin the
\begin{proof}[Proof of Theorem \ref{theoremisometries}]
The proof will consist of two parts. In the first, we show that $\PO_\F^*(\LL;\QQ) = \Isom(\H)$, and in the second we show that $\PO_\F(\LL;\QQ)$ is equal to $\PO_\F^*(\LL;\QQ)$ if $\F = \R,\Q$ and is of index $2$ in $\PO_\F^*(\LL;\QQ)$ if $\F = \C$.

Fix $g\in \Isom(\H)$; we claim that $g\in\PO_\F^*(\LL;\QQ)$. Let $\zz = (1,\0)$, and let $\zero = [\zz]$. By Observation \ref{observationtransitivity}, there exists $[T]\in \PO_\F(\LL;\QQ)$ such that $[T](\zero) = g(\zero)$. Thus, we may without loss of generality assume that $g(\zero) = \zero$.

We observe that $\zz^\perp = \HH$. Let $S(\HH)$ denote the unit sphere of $\HH$, i.e. $S(\HH) = \{\ww\in\HH : \QQ(\ww) = 1\}$. For each $\ww\in S(\HH)$, the embedding $\gamma_{\zz,\ww}:\R\to\H$ defined in the proof of Proposition \ref{propositionROSSONCT} is an isometry. By Proposition \ref{propositionROSSONCT}, its image under $g$ must also be an isometry. Specifically, there exists $f(\ww)\in S(\HH)$ such that
\begin{equation}
\label{tanht}
g([\zz + \tanh(t)\ww]) = [\zz + \tanh(t)f(\ww)] \all t\in\R.
\end{equation}
The fact that $g$ is a bijection implies that $f:S(\HH)\to S(\HH)$ is a bijection. Moreover, the fact that $g$ is an isometry means that for all $\ww_1,\ww_2\in S(\HH)$ and $t_1,t_2\in\R$, we have
\[
\dist([\zz + \tanh(t_1)\ww_1], [\zz + \tanh(t_2)\ww_2]) = \dist([\zz + \tanh(t_1)f(\ww_1)], [\zz + \tanh(t_2)f(\ww_2)]).
\]
Recalling that
\begin{align*}
\cosh\dist([\zz + \tanh(t_1)\ww_1] &, [\zz + \tanh(t_2)\ww_2])\\ &= \frac{|B_\QQ(\zz + \tanh(t_1)\ww_1, \zz + \tanh(t_2)\ww_2)|}{\sqrt{|\QQ(\zz + \tanh(t_1)\ww_1)|\cdot|\QQ(\zz + \tanh(t_2)\ww_2)|}}\\
&= \frac{|-1 + \tanh(t_1)\tanh(t_2)B_\QQ(\ww_1,\ww_2)|}{\sqrt{\left(1 - \tanh^2(t_1)\right)\left(1 - \tanh^2(t_2)\right)}},
\end{align*}
we see that
\[
|-1 + \tanh(t_1)\tanh(t_2)B_\QQ(\ww_1,\ww_2)| = |-1 + \tanh(t_1)\tanh(t_2)B_\QQ(f(\ww_1),f(\ww_2))|.
\]
Write $\theta = \tanh(t_1)\tanh(t_2)$. Squaring both sides gives
\begin{eqnarray}
\label{thetapolynomial}
\mathrm{LHS}^2 & = & \theta^2 |B_\QQ(\ww_1,\ww_2)|^2 - 2\theta \Re[B_\QQ(\ww_1,\ww_2)] + 1 \nonumber \\
|| &  & \\ 
\mathrm{RHS}^2 & = & \theta^2 |B_\QQ(f(\ww_1),f(\ww_2))|^2 - 2\theta \Re[B_\QQ(f(\ww_1),f(\ww_2))] + 1. \nonumber
\end{eqnarray}
We observe that for $\ww_1,\ww_2\in S(\HH)$ fixed, \eqref{thetapolynomial} holds for all $-1 < \theta < 1$. In particular, taking the first and second derivatives and plugging in $\theta = 0$ gives
\begin{align} \label{represerved}
\Re[B_\QQ(\ww_1,\ww_2)] &= \Re[B_\QQ(f(\ww_1),f(\ww_2))]\\ \label{abspreserved}
|B_\QQ(\ww_1,\ww_2)| &= |B_\QQ(f(\ww_1),f(\ww_2))|.
\end{align}
Extend $f$ to a bijection $f:\HH\to\HH$ by letting $f(\0) = \0$ and $f(t\ww) = tf(\ww)$ for $t > 0$, $\ww\in S(\HH)$. We observe that \eqref{represerved} and \eqref{abspreserved} hold also for the extended version of $f$.

\begin{claim}
\label{claimflinear}
$f$ is $\R$-linear.
\end{claim}
\begin{subproof}
Fix $\ww_1,\ww_2\in\HH$ and $c_1,c_2\in\R$. By \eqref{represerved}, the maps
\[
\ww\mapsto \Re[B_\QQ(f(c_1\ww_1 + c_2\ww_2),f(\ww))] \text{ and } \ww\mapsto \Re[B_\QQ(c_1 f(\ww_1) + c_2 f(\ww_2),f(\ww))] 
\]
are identical. By the surjectivity of $f$ together with the Riesz representation theorem, this implies that $f(c_1\ww_1 + c_2\ww_2) = c_1 f(\ww_1) + c_2 f(\ww_2)$.
\end{subproof}
\begin{claim}
\label{claimFlines}
$f$ preserves $\F$-lines.
\end{claim}
\begin{subproof}
For each $\xx\in\HH\butnot\{\0\}$, the $\F$-line $\xx\F$ may be defined using the quantity $|B_\QQ|$ via the formula
\[
\xx\F = \left\{\yy\in\HH : \forall \ww\in \HH, \;\; |B_\QQ(\xx,\ww)| = 0 \;\;\Leftrightarrow\;\; |B_\QQ(\yy,\ww)| = 0\right\}.
\]
The claim therefore follows from \eqref{abspreserved}.
\end{subproof}
From Claim \ref{claimFlines}, we see that for all $\xx\in\HH\butnot\{\0\}$ and $a\in\F$, there exists $\sigma_\xx(a)\in\F$ such that
\[
f(\xx a) = f(\xx) \sigma_\xx(a).
\]
\begin{claim}
For $\xx,\yy\in\HH\butnot\{\0\}$,
\[
\sigma_\xx(a) = \sigma_\yy(a).
\]
\end{claim}
\begin{subproof}
By Claim \ref{claimflinear},
\[
[f(\xx) + f(\yy)] \sigma_{\xx + \yy}(a) = f(\xx a + \yy a) = f(\xx)\sigma_\xx(a) + f(\yy)\sigma_\yy(a).
\]
Rearranging, we see that
\[
f(\xx)[\sigma_{\xx + \yy}(a) - \sigma_\xx(a)] + f(\yy)[\sigma_{\xx + \yy}(a) - \sigma_\yy(a)] = 0.
\]
If $\xx$ and $\yy$ are linearly independent, then $\sigma_{\xx + \yy}(a) - \sigma_\xx(a) = 0$ and $\sigma_{\xx + \yy}(a) - \sigma_\yy(a) = 0$, so $\sigma_\xx(a) = \sigma_\yy(a)$. But the general case clearly follows from the linearly independent case.
\end{subproof}
For $a\in\F$, denote the common value of $\sigma_\xx(a)$ ($\xx\in\HH\butnot\{\0\}$) by $\sigma(a)$. Then
\begin{equation}
\label{fsigma}
f(\xx a) = f(\xx)\sigma(a) \all \xx\in\HH \all a\in\F.
\end{equation}
\begin{claim}
$\sigma\in\Aut(\F)$.
\end{claim}
\begin{subproof}
The $\R$-linearity of $\sigma$ follows from Claim \ref{claimflinear}, and the bijectivity of $\sigma$ follows from the bijectivity of $f$. Fix $\xx\in\HH\butnot\{\0\}$ arbitrary. For $a,b\in\F$,
\[
f(\xx)\sigma(ab) = f(\xx ab) = f(\xx a)\sigma(b) = f(\xx) \sigma(a)\sigma(b),
\]
which proves that $\sigma$ is a multiplicative homomorphism.
\end{subproof}
Thus $f\in\O_\F^*(\HH;\EE)$, and so $T = f\oplus I\in\O_\F^*(\LL;\QQ)$. But $[T] = g$ by \eqref{tanht}, so $g\in\PO_\F^*(\LL;\QQ)$. This completes the first part of the proof, namely that $\PO_\F^*(\LL;\QQ) = \Isom(\H)$.

To complete the proof, we need to show that $\PO_\F(\LL;\QQ)$ is equal to $\PO_\F^*(\LL;\QQ)$ if $\F = \R,\Q$ and is of index $2$ in $\PO_\F^*(\LL;\QQ)$ if $\F = \C$. If $\F = \R$, this is obvious. If $\F = \C$, it follows from the semidirect product structure $\O_\F^*(\LL;\QQ) = \O_\F(\LL;\QQ) \ltimes \{\sigma^J : \sigma\in\Aut(\F)\}$ together with the fact that $\Aut(\F) = \{I,z\mapsto \bar z\} \equiv\Z_2$.

If $\F = \Q$, then $\Aut(\F) = \{\Phi_a : a\in S(\Q)\}$, where $\Phi_a(b) = aba^{-1}$. Here $S(\F) = \{a\in\F : |a| = 1\}$. So $\O_\Q^*(\LL;\QQ) \neq \O_\Q(\LL;\QQ)$; nevertheless, we will show that $\PO_\Q^*(\LL;\QQ) = \PO_\Q(\LL;\QQ)$. Fix $[T]\in \PO_\Q^*(\LL;\QQ)$, and fix $a\in S(\Q)$ for which $\sigma_T = \Phi_a$. Consider the map
\begin{equation}
\label{Tadef}
T_a(\xx) = \xx a.
\end{equation}
We have $T_a\in\O_\Q^*(\LL;\QQ)$ and $\sigma_{T_a} = \Phi_a^{-1}$. Thus $\sigma_{T_a T} = \Phi_a \Phi_a^{-1} = I$, so $T_a T$ is $\F$-linear. But
\[
[T_a T] = [T],
\]
so $[T]\in\PO_\Q(\LL;\QQ)$. The completes the proof of Theorem \ref{theoremisometries}.
\end{proof}

\begin{remark}
Using algebraic language, the automorphisms $\Phi_a$ of $\Q$ are \emph{inner} automorphisms, while the automorphism $z\mapsto \bar z$ of $\C$ is an \emph{outer} automorphism. Although both inner and outer automorphisms contribute to the quotient $\O_\F^*(\LL;\QQ)/\O_\F(\LL;\QQ)$, only the outer automorphisms contribute to the quotient $\PO_\F^*(\LL;\QQ)/\PO_\F(\LL;\QQ)$. This explains why the index $\#(\PO_\F^*(\LL;\QQ)/\PO_\F(\LL;\QQ))$ is smaller when $\F = \Q$ than when $\F = \C$: although the group $\Aut(\Q)$ is much larger than $\Aut(\C)$, it consists entirely of inner automorphisms, while $\Aut(\C)$ has an outer automorphism.
\end{remark}

\begin{definition}
The \emph{bordification} of $\H$ is its closure relative to the topological space $\proj(\LL)$, i.e.
\[
\bord\H = \{[\xx] : \QQ(\xx) \leq 0\}.
\]
The \emph{boundary} of $\H$ is its topological boundary relative to $\proj(\LL)$, i.e.
\[
\del\H = \bord\H\butnot\H = \{[\xx] : \QQ(\xx) = 0\}.
\]
\end{definition}

The following is a corollary of Theorem \ref{theoremisometries}:

\begin{corollary}
\label{corollaryextension}
Every isometry of $\H$ extends uniquely to a homeomorphism of $\bord\H$.
\end{corollary}
\begin{proof}
If $T\in\O_\F^*(\LL;\QQ)$, then the formula \eqref{[T]def} defines a homeomorphism of $\bord\H$ which extends the action of $[T]$ on $\H$. The uniqueness is automatic.
\end{proof}
\begin{remark}
Corollary \ref{corollaryextension} can also be proven independently of Theorem \ref{theoremisometries} via the theory of hyperbolic metric spaces; cf. Lemma \ref{lemmaisometryextension} and Proposition \ref{propositionboundariesequivalent}.
\end{remark}

The following observation will be useful in the sequel:
\begin{observation}
\label{observationBQnonzero}
Fix $[\xx],[\yy]\in\bord\H$. Then
\[
B_\QQ(\xx,\yy) = 0 \;\;\Leftrightarrow\;\; [\xx] = [\yy] \in\del\H.
\]
\end{observation}
\begin{proof}
If either $[\xx]$ or $[\yy]$ is in $\H$, this follows from Lemma \ref{lemmasylvester}. Suppose that $[\xx],[\yy]\in\del\H$, and that $B_\QQ(\xx,\yy) = 0$. Then $\QQ$ is identically zero on $\xx\F + \yy\F$. Thus $(\xx\F + \yy\F)\cap\HH = \{\0\}$, and so $\xx\F + \yy\F$ is one-dimensional. This implies $[\xx] = [\yy]$.
\end{proof}

\section{Totally geodesic subsets of algebraic hyperbolic spaces}
\label{subsectiontotallygeodesic}

Given two pairs $(X,\bord X)$ and $(Y,\bord Y)$, where $X$ and $Y$ are metric spaces contained in the topological spaces $\bord X$ and $\bord Y$ (and dense in these spaces), an \emph{isomorphism} between $(X,\bord X)$ and $(Y,\bord Y)$ is a homeomorphism between $\bord X$ and $\bord Y$ which restricts to an isometry between $X$ and $Y$.

\begin{proposition}
\label{propositiontotallygeodesic}
Let $\K\leq\F$ be an $\R$-subalgebra, and let $V\leq\LL$ be a closed (right) $\K$-module such that
\begin{equation}
\label{totallygeodesic}
B_\QQ(\xx,\yy)\in\K \all \xx,\yy\in V.
\end{equation}
Then either $[V]\cap\H = \emptyset$ and $\#([V]\cap\bord\H)\leq 1$, or $([V]\cap\H,[V]\cap\bord\H)$ is isomorphic to an algebraic hyperbolic space together with its closure.
\end{proposition}
\begin{proof}~
\begin{itemize}
\item[Case 1:] $[V]\cap\H\neq\emptyset$. In this case, fix $[\zz]\in [V]\cap\H$, and let $\zz$ be a representative of $[\zz]$ with $\QQ(\zz) = -1$. By Lemma \ref{lemmasylvester}, $\QQ$ is positive-definite on $\zz^\perp$. We leave it as an exercise that the quadratic forms $\QQ\given\zz^\perp$ and $\EE\given\zz^\perp$ agree up to a bounded multiplicative error factor, which implies that $\zz^\perp$ is complete with respect to the norm $\sqrt\QQ$.

From \eqref{totallygeodesic}, we see that $(V\cap\zz^\perp,B_\QQ)$ is a $\K$-Hilbert space. By the usual Gram--Schmidt process, we may construct an orthonormal basis $(\ee_i)_{i\in J'}$ for $V\cap\zz^\perp$, thus proving that $V\cap\zz^\perp$ is isomorphic to $\HH_\K^{J'}$ for some set $J'$. Thus $V$ is isomorphic to $\LL_\K^{J'\cup\{0\}}$, and so $([V]\cap\H,[V]\cap\bord\H)$ is isomorphic to $(\H_\K^{J'},\bord\H_\K^{J'})$.

\item[Case 2:] $[V]\cap\H = \emptyset$. We need to show that $\#([V]\cap\bord\H)\leq 1$. By contradiction fix $[\xx],[\yy]\in[V]$ distinct, and let $\xx,\yy\in V$ be representatives. By Observation \ref{observationBQnonzero}, $B_\QQ(\xx,\yy)\neq 0$. On the other hand, $\QQ(\xx) = \QQ(\yy) = 0$ since $[\xx],[\yy]\in\del\H$. Thus $\QQ(\xx - \yy B(\xx,\yy)^{-1}) = -2 < 0$. On the other hand, $\xx - \yy B(\xx,\yy)^{-1}\in V$ by \eqref{totallygeodesic}. Thus $[\xx - \yy B(\xx,\yy)^{-1}]\in [V]\cap\H$, a contradiction.
\end{itemize}
\end{proof}

\begin{definition}
\label{definitiontotallygeodesic}
A \emph{totally geodesic subset} of an algebraic hyperbolic space $\H$ is a set of the form $[V]\cap\bord\H$, where $V$ is as in Proposition \ref{propositiontotallygeodesic}. A totally geodesic subset is \emph{nontrivial} if it contains an element of $\H$.
\end{definition}

\begin{remark}
As with Definition \ref{definitionROSSONCT}, the terminology ``totally geodesic'' is motivated here by the finite-dimensional situation, where totally geodesic subsets correspond precisely with the closures of those submanifolds which are totally geodesic in the sense of Riemannian geometry; see \cite[Proposition A.4 and A.7]{Quint}. However, note that we consider both the empty set and singletons in $\del\H$ to be totally geodesic.
\end{remark}

\begin{remark}
\label{remarkVrotation}
If $V\leq\LL$ is a closed $\K$-module satisfying \eqref{totallygeodesic}, then for each $a\in\F\butnot\{0\}$, $Va$ is a closed $a^{-1}\K a$-module satisfying \eqref{totallygeodesic} (with $\K = a^{-1}\K a$).
\end{remark}

\begin{lemma}
\label{lemmatotallygeodesic}
The intersection of any collection of totally geodesic sets is totally geodesic.
\end{lemma}
\begin{proof}
Suppose that $(S_\alpha)_{\alpha\in A}$ is a collection of totally geodesic sets, and suppose that $S = \bigcap_\alpha S_\alpha\neq\emptyset$. Fix $[\zz]\in S$, and let $\zz$ be a representative of $[\zz]$. Then for each $\alpha\in A$, there exist (cf. Remark \ref{remarkVrotation}) an $\R$-algebra $\K_\alpha$ and a closed $\K_\alpha$-subspace $V_\alpha \leq \LL$ satisfying \eqref{totallygeodesic} (with $\K = \K_\alpha$) such that $\zz\in V_\alpha$ and $S_\alpha = [V_\alpha]\cap\bord\H$. Let $\K = \bigcap_\alpha\K_\alpha$ and $V = \bigcap_\alpha V_\alpha$. Clearly, $V$ is a $\K$-module and satisfies \eqref{totallygeodesic}.

We have $[V]\cap\bord\H\subset S$. To complete the proof, we must show the converse direction. Fix $[\xx]\in S\butnot\{[\zz]\}$. By Observation \ref{observationBQnonzero}, there exists a representative $\xx$ of $[\xx]$ such that $B_\QQ(\zz,\xx) = 1$. Then for each $\alpha$, we may find $a_\alpha\in\F\butnot\{0\}$ such that $\xx a_\alpha\in V_\alpha$. We have
\[
a_\alpha = B_\QQ(\zz,\xx)a_\alpha = B_\QQ(\zz,\xx a_\alpha)\in\K_\alpha.
\]
Since $V_\alpha$ is a $\K_\alpha$-module, this implies $\xx\in V_\alpha$. Since $\alpha$ was arbitrary, $\xx\in V$, and so $[\xx]\in[V]\cap\bord\H$.
\end{proof}
\begin{remark}
Given $K\subset\bord\H$, Lemma \ref{lemmatotallygeodesic} implies that there exists a smallest totally geodesic set containing $K$. If we are only interested in the geometry of $K$, then by Proposition \ref{propositiontotallygeodesic} we can assume that this totally geodesic set is really our ambient space. In such a situation, we may without loss of generality suppose that there is no proper totally geodesic subset of $\bord\H$ which contains $K$. In this case we say that $K$ is \emph{irreducible}.
\end{remark}

\begin{warning*}
Although the intersection of any collection of totally geodesic sets is totally geodesic, it is not necessarily the case that the decreasing intersection of nontrivial totally geodesic sets is nontrivial; cf. Remark \ref{remarkparabolictorsion}.
\end{warning*}


The main reason that totally geodesic sets are relevant to our development is their relationship with the group of isometries. Specifically, we have the following:
\begin{theorem}
\label{theoremtotallygeodesic}
Let $(g_n)_1^\infty$ be a sequence in $\Isom(\H)$, and let
\begin{equation}
\label{fixedROSSONCT}
S = \left\{[\xx]\in\bord\H : g_n([\xx]) \tendsto n [\xx]\right\}.
\end{equation}
Then either $S\subset\del\H$ and $\#(S) = 2$, or $S$ is a totally geodesic set.
\end{theorem}
\begin{remark}
\label{remarktotallygeodesic}
An important example is the case where the sequence $(g_n)_1^\infty$ is constant, say $g_n = g$ for all $n$. Then $S$ is precisely the \emph{fixed point set} of $g$:
\[
S = \Fix(g) := \left\{[\xx]\in\bord\H : g([\xx]) = [\xx]\right\}.
\]
If $\H$ is finite-dimensional, then it is possible to reduce Theorem \ref{theoremtotallygeodesic} to this special case by a compactness argument.
\end{remark}
\begin{proof}[Proof of Theorem \ref{theoremtotallygeodesic}]
If $S = \emptyset$, then the statement is trivial. Suppose that $S\neq\emptyset$, and fix $[\zz]\in S$.

{\it Step 1: Choosing representatives $T_n$.} From the proof of Theorem \ref{theoremisometries}, we see that each $g_n$ may be written in the form $[T_n]$ for some $T_n\in\O_\F^*(\LL;\QQ)$. We have some freedom in choosing the representatives $T_n$; specifically, given $a_n\in S(\F)$ we may replace $T_n$ by $T_n T_{a_n}$, where $T_{a_n}$ is defined by \eqref{Tadef}.

Since $g_n([\zz])\to [\zz]$, there exist representatives $\zz_n$ of $g_n([\zz])$ such that $\zz_n\to \zz$. For each $n$, there is a unique representative $T_n$ of $g_n$ such that
\[
(T_n \zz) c_n = \zz_n \text{ for some $c_n\in\R\butnot\{0\}$.}
\]
Then
\[
(T_n \zz) c_n \to \zz.
\]
\begin{remark}
If $\F = \Q$, it may be necessary to choose $T_n\in\O_\F^*(\LL;\QQ)\butnot\O_\F(\LL;\QQ)$, despite the fact that each $g_n$ can be represented by an element of $\O_\F(\LL;\QQ)$.
\end{remark}

{\it Step 2: A totally geodesic set.}
Write $\sigma_n = \sigma_{T_n}$, and let
\begin{align*}
\K &= \{a\in \F : \sigma_n(a)\to a\}\\
V &= \left\{\xx\in\LL : T_n \xx \tendsto n \xx\right\}.
\end{align*}
Then $\K$ is an $\R$-subalgebra of $\F$, and $V$ is a $\K$-module. Given $\xx,\yy\in V$, by Observation \ref{observationBQpreserving} we have
\[
\sigma_n(B_\QQ(\xx,\yy)) = B_\QQ(T_n\xx,T_n\yy) \tendsto n B_\QQ(\xx,\yy),
\]
so $B(\xx,\yy)\in \K$. Thus $V$ satisfies \eqref{totallygeodesic}. If $V$ is closed, then the above observations show that $[V]\cap\bord\H$ is totally geodesic. However, this issue is a bit delicate:

\begin{claim}
If $\#([V]\cap\bord\H)\geq 2$, then $V$ is closed.
\end{claim}
\begin{subproof}
Suppose that $\#([V]\cap\bord\H)\geq 2$. The proof of Proposition \ref{propositiontotallygeodesic} shows that $[V]\cap\H\neq\emptyset$. Thus, there exists $\xx\in V$ for which $[\xx]\in\H$. In particular, $g_n([\xx])\to [\xx]$. Letting $\zero = [(1,\0)]$, we have
\[
\dist_\H(\zero,g_n(\zero)) \leq 2\dist_\H(\zero,[\xx]) + \dist_\H([\xx],g_n([\xx])) \tendsto n 2\dist_\H(\zero,[\xx]).
\]
In particular $\dist_\H(\zero,g_n(\zero))$ is bounded, say $\dist_\H(\zero,g_n(\zero)) \leq C$.

\begin{lemma}
\label{lemmaoperatornorm}
Fix $T\in\O_\F^*(\LL;\QQ)$, and let $\|T\|$ denote the operator norm of $T$. Then
\[
\|T\| = e^{\dist_\H(\zero,[T](\zero))}.
\]
\end{lemma}
\begin{subproof}
Write $T = T_{j,t}(A\oplus I)$, where $T_{j,t}$ is a Lorentz boost (cf. \eqref{lorentzboost}) and $A\in\O_\F^*(\HH;\EE)$. Then
\[
[T](\zero) = [T_{j,t}](\zero) = [(\cosh(t),\sinh(t),\0)].
\]
Here the second entry represents the $j$th coordinate. In particular,
\[
\cosh\dist_\H(\zero,[T](\zero)) = \frac{|B_\QQ((1,\0),(\cosh(t),\sinh(t),\0))|}{\sqrt{|\QQ(1,\0)|\cdot|\QQ(\cosh(t),\sinh(t),\0)|}} = \frac{\cosh(t)}{1} = \cosh(t).
\]
On the other hand,
\[
\|T\| = \|T_{j,t}\| = \left\|
\left[\begin{array}{ccc}
\cosh(t) & \sinh(t) & \\
\sinh(t) & \cosh(t) &\\
&& I
\end{array}\right] \right\| = e^t.
\]
This completes the proof.
\end{subproof}
\noindent Thus $\|T_n\| \leq e^C$ for all $n$, and so the sequence $(T_n)_1^\infty$ is equicontinuous. It follows that $V$ is closed.
\end{subproof}
\noindent Since $\#([V]\cap\bord\H)\leq 1$ implies that $[V]\cap\bord\H$ is totally geodesic, we conclude that $[V]\cap\bord\H$ is totally geodesic, regardless of whether or not $V$ is closed.

\begin{remark}
When $\#([V]\cap\bord\H)\leq 1$, there seems to be no reason to think that $V$ should be closed.
\end{remark}

{\it Step 3: Relating $S$ to $[V]\cap\bord\H$.}
The object of this step is to show that $S = [V]\cap\bord\H$ unless $S\subset\del\H$ and $\#(S)\leq 2$. For each $[\xx]\in S\butnot\{[\zz]\}$, let $\xx$ be a representative of $[\xx]$ such that $B_\QQ(\zz,\xx) = 1$; this is possible by Observation \ref{observationBQnonzero}. It is possible to choose a sequence of scalars $(a_n^{([\xx])})_{n = 1}^\infty$ in $\F\butnot\{\0\}$ such that $(T_n\xx) a_n^{([\xx])}\to \xx$. Let $a_n^{([\zz])} = c_n$. For $[\xx],[\yy]\in S$, we have
\begin{equation}
\label{anxany}
\begin{split}
\wbar{a_n^{([\xx])}}\sigma_{T_n}(B_\QQ(\xx,\yy)) a_n^{([\yy])}
&= \wbar{a_n^{([\xx])}}B_\QQ(T_n\xx,T_n\yy) a_n^{([\yy])} \hspace{.5 in} \text{(by Observation \ref{observationBQpreserving})}\\
&= B_\QQ((T_n\xx) a_n^{([\xx])},(T_n\yy) a_n^{([\yy])})\\
&\tendsto n B_\QQ(\xx,\yy).
\end{split}
\end{equation}
In particular,
\begin{equation}
\label{anxanyabs}
|a_n^{([\xx])}|\cdot|a_n^{([\yy])}| \tendsto n 1 \text{ whenever }B_\QQ(\xx,\yy)\neq 0.
\end{equation}
\begin{claim}
Unless $S\subset\del\H$ and $\#(S)\leq 2$, then for all $[\xx]\in S$ we have
\begin{equation}
\label{anx}
|a_n^{([\xx])}|\tendsto n 1.
\end{equation}
\end{claim}
\begin{subproof}
We first observe that it suffices to demonstrate \eqref{anx} for one value of $\xx$; if \eqref{anx} holds for $\xx$ and $[\yy]\neq[\xx]$, then $B_\QQ(\xx,\yy)\neq 0$ by Observation \ref{observationBQnonzero} and so \eqref{anxanyabs} implies $|a_n^{([\yy])}|\to 1$.

Now suppose that $S\nsubset\del\H$, and choose $[\xx]\in S\cap\H$. Then $B_\QQ(\xx,\xx)\neq 0$, and so \eqref{anxanyabs} implies \eqref{anx}.

Finally, suppose that $\#(S)\geq 3$, and choose $[\xx],[\yy],[\zz]\in S$ distinct. By \eqref{anxanyabs} together with Observation \ref{observationBQnonzero}, we have $|a_n^{([\xx])}|\cdot|a_n^{([\yy])}| \to 1$, $|a_n^{([\xx])}|\cdot|a_n^{([\zz])}| \to 1$, and $|a_n^{([\yy])}|\cdot|a_n^{([\zz])}| \to 1$. Multiplying the first two formulas and dividing by the third, we see that $|a_n^{([\xx])}|\to 1$.
\end{subproof}

For the remainder of the proof we assume that either $S\nsubset\del\H$ or $\#(S)\geq 3$.

Plugging $\zz = \xx$ into \eqref{anx}, we see that $c_n\to 1$. In particular, $[\zz]\in[V]\cap\bord\H$. Now fix $[\xx]\in S\butnot\{[\zz]\}$. Since $c_n\to 1$ and $B_\QQ(\zz,\xx) = 1$, \eqref{anxany} becomes
\[
a_n^{([\xx])} \to 1.
\]
Thus $\xx\in V$, and so $[\xx]\in [V]\cap\bord\H$.

\end{proof}

\bigskip
\section{Other models of hyperbolic geometry}
\label{subsectionmodels}

Fix $\F\in\{\R,\C,\Q\}$ and a set $J$, and let $\H = \H_\F^J$. The pair $(\H,\bord\H)$ is known as the \emph{hyperboloid model} of hyperbolic geometry (over the division algebra $\F$ and in dimension $\#(J)$). In this section we discuss two other important models of hyperbolic geometry. Note that the Poincar\'e ball model, which many of the figures of later chapters are drawn in, is not discussed here. References for this section include \cite{CFKP, Goldman}.

\subsection{The (Klein) ball model}
\label{subsubsectionballmodel}
Let
\[
\B = \B_\F^J = \{\xx\in \HH := \HH_\F^J : \|\xx\| < 1\},
\]
and let $\bord\B$ denote the closure of $\B$ relative to $\HH$.
\begin{observation}
\label{observationHequivB}
The map $e_{\B,\H}:\bord\B\to\bord\H$ defined by the equation
\[
e_{\B,\H}(\xx) = [(1,\xx)]
\]
is a homeomorphism, and $e_{\B,\H}(\B) = \H$. Thus if we let
\begin{equation}
\label{distB}
\cosh\dist_\B(\xx,\yy) = \cosh\dist_\H(e_{\B,\H}(\xx),e_{\B,\H}(\yy)) = \frac{|1 - B_\EE(\xx,\yy)|}{\sqrt{1 - \|\xx\|^2}\sqrt{1 - \|\yy\|^2}},
\end{equation}
then $e_{\B,\H}$ is an isomorphism between $(\B,\bord\B)$ and $(\H,\bord\H)$.
\end{observation}
The pair $(\B,\bord\B)$ is called the \emph{ball model} of hyperbolic geometry. It is often convenient for computations, especially those for which a single point plays an important role: by Observation \ref{observationtransitivity}, such a point can be moved to the origin $\0\in\B$ via an isomorphism of $(\B,\bord\B)$.

\begin{remark}
We should warn that the ball model $\B_\R^J$ of real hyperbolic geometry is not the same as the well-known Poincar\'e model, rather, it is the same as the Klein model.
\end{remark}

\begin{observation}
For all $T\in \O_\F^*(\HH;\EE)$, $T\given\B$ is an isometry which stabilizes $\0$.
\end{observation}

\begin{proposition}
\label{propositionIsomB}
In fact,
\[
\Stab(\Isom(\B);\0) = \{T\given\B : T\in\O_\F^*(\HH;\EE)\}.
\]
\end{proposition}
\begin{proof}
This is an immediate consequence of Theorem \ref{theoremisometries}.
\end{proof}

\subsection{The half-space model}
\label{subsubsectionE}

Now suppose $\F = \R$.\Footnote{The appropriate analogue of the half-space model when $\amsbb F\in\{\amsbb C,\amsbb Q\}$ is the \emph{paraboloid model}; see e.g. \cite[Chapter 4]{Goldman}.} Assume that $1\in J$, and let
\[
\E = \E^J = \left\{\left.\xx\in\HH := \HH_\F^J \right\vert x_1 > 0 \right\}.
\]
We will view $\E$ as resting inside the larger space
\[
\what\HH := \HH \cup\{\infty\}.
\]
The topology on $\what\HH$ is defined as follows: a subset $U \subset \what\HH$ is open if and only if
\[
\text{$U \cap \HH$ is open and ($\infty\in U$ \implies $\HH \setminus U$ is bounded)}.
\]
The boundary and closure of $\E$ will be subsets of $\what\HH$ according to the topology defined above, i.e.
\begin{align*}
\del\E &= \{\xx\in\HH : x_1 = 0\}\cup\{\infty\}\\
\bord\E &= \{\xx\in\HH : x_1 \geq 0\}\cup\{\infty\}.
\end{align*}
\begin{proposition}
\label{propositionHequivE}
The map $e_{\E,\H}:\bord\E\to\bord\H$ defined by the formula
\begin{equation}
\label{paraboloidembedding}
e_{\E,\H}(\xx) = \begin{cases}
\left[\left(\begin{cases}
2x_i & i\neq 0,1\\
1 + \|\xx\|^2 & i = 0\\
1 - \|\xx\|^2 & i = 1
\end{cases}\right)_{i\in J\cup\{0\}}\right] & \xx\neq\infty\\
[(1,-1,\0)] & \xx = \infty
\end{cases}
\end{equation}
is a homeomorphism, and $e_{\E,\H}(\E) = \H$. Thus if we let
\begin{equation}
\label{distE}
\cosh\dist_\E(\xx,\yy) = \cosh\dist_\H(e_{\E,\H}(\xx),e_{\E,\H}(\yy)) = 1 + \frac{\|\yy - \xx\|^2}{2x_1 y_1},
\end{equation}
then $e_{\E,\H}$ is an isomorphism between $(\E,\bord\E)$ and $(\H,\bord\H)$.
\end{proposition}
\begin{proof}
For $\xx\in\bord\E\butnot\{\infty\}$,
\[
\QQ(e_{\E,\H}(\xx)) = -(1 + \|\xx\|^2)^2 + (1 - \|\xx\|^2)^2 + \sum_{i\in J\butnot\{1\}}(2x_i)^2 = -4x_1^2.
\]
It follows that $e_{\E,\H}(\E) \subset\H$ and $e_{\E,\H}(\del\E)\subset\del\H$. Calculation verifies that the map
\begin{equation}
\label{iotainverse}
e_{\H,\E}([\xx]) = \begin{cases}\left(\begin{cases}
x_i/2 & i\neq 1\\
\sqrt{-\QQ(\xx)}/2 & i = 1
\end{cases}\right)_{i\in J} &\text{ if }x_0 + x_1 = 2 \\
\infty &\text{ if } \xx = (1,-1,\0)
\end{cases}
\end{equation}
is both a left and a right inverse of $e_{\E,\H}$. Notice that it is defined in a way such that for each $[\xx]\in\bord\H$, there is a unique representative $\xx$ of $[\xx]$ for which the formula \eqref{iotainverse} makes sense. We leave it to the reader to verify that $e_{\E,\H}$ and $e_{\H,\E}$ are both continuous, and that \eqref{distE} holds.
\end{proof}

The point $\infty\in\del\E$, corresponding to the point $[(1,-1,\0)]\in\del\H$, plays a special role in the half-space model. In fact, the half-space model can be thought of as an attempt to understand the geometry of hyperbolic space when a single point on the boundary is fixed. Consequently, we are less interested in the set of all isometries of $\E$ than simply the set of all isometries which fix $\infty$.

\begin{observation}[Poincar\'e extension]
\label{observationpoincareextension}
Let $\BB = \del\E\butnot\{\infty\} = \HH_\R^{J\butnot\{1\}}$, and let $g:\BB\to\BB$ be a \emph{similarity}, i.e. a map of the form
\[
g(\xx) = \lambda T\xx + \bb,
\]
where $\lambda > 0$, $T\in\O_\R(\BB;\EE)$, and $\bb\in\BB$. Then the map $\what g:\bord\E\to\bord\E$ defined by the formula
\begin{equation}
\label{poincareextension}
\what g(\xx) = \begin{cases}
(\lambda x_1,g(\pi(\xx))) & \xx\neq \infty \\
\infty & \xx = \infty
\end{cases}
\end{equation}
is an isomorphism of $(\E,\bord\E)$; in particular, $\what g\given\E$ is an isometry of $\E$. Here $\pi:\HH\to\BB$ is the natural projection.
\end{observation}
\begin{proof}
This is immediate from \eqref{distE}.
\end{proof}
The isometry $\what g$ defined by \eqref{poincareextension} is called the \emph{Poincar\'e extension} of $g$ to $\E$.

\begin{remark}
Intuitively we shall think of the number $x_1$ as representing the \emph{height} of a point $\xx\in\bord\E$. Then \eqref{poincareextension} says that if $g:\BB\to\BB$ is an isometry, then the Poincar\'e extension of $g$ is an isometry of $\E$ which preserves the heights of points.
\end{remark}

\begin{proposition}
\label{propositionIsomE}
For all $g\in\Isom(\E)$ such that $g(\infty) = \infty$, there exists a similarity $h:\BB\to\BB$ such that $g = \what h$.
\end{proposition}
\begin{proof}
By Theorem \ref{theoremisometries}, there exists $T\in\O(\LL;\QQ)$ such that $[T] = e_{\E,\H}\circ g\circ e_{\E,\H}^{-1}$. This gives an explicit formula for $g$, and one must check that if $[T]$ preserves $[(1,-1,\0)]$, then $g$ is a Poincar\'e extension.
\end{proof}

\subsection{Transitivity of the action of $\Isom(\H)$ on $\del\H$}
Using the ball and half-plane models of hyperbolic geometry, it becomes easy to prove the following assertion:

\begin{proposition}
\label{propositiontransitivity}
If $\F = \R$, the group $\Isom(\H)$ acts triply transitively on $\del\H$.
\end{proposition}
This complements the fact that $\Isom(\H)$ acts transitively on $\H$ (Observation \ref{observationtransitivity}).
\begin{proof}
By Observation \ref{observationHequivB} and Proposition \ref{propositionHequivE}, we may switch between models as convenient. It is clear that $\Isom(\B)$ acts transitively on $\del\B$, and that $\Stab(\Isom(\E);\infty)$ acts doubly transitively on $\del\E\butnot\{\infty\}$. Therefore given any triple $(\xi_1,\xi_2,\xi_3)$, we may conjugate to $\B$, conjugate $\xi_1$ to a standard point, conjugate to $\E$ while conjugating $\xi_1$ to $\infty$, and then conjugate $\xi_2,\xi_3$ to standard points.
\end{proof}

We end this chapter with a convention:

\begin{convention}
\label{conventionH}
When $\alpha$ is a cardinal number, $\H_\F^\alpha$ will denote $\H_\F^J$ for any set $J$ of cardinality $\alpha$, but particularly $J = \{1,\ldots,n\}$ if $\alpha = n\in\Namer$ and $J = \Namer$ if $\alpha = \#(\Namer)$. Moreover, $\H_\F^\infty$ will always be used to denote $\H_\F^{\#(\Namer)} = \H_\F^{\Namer}$, the unique (up to isomorphism) infinite-dimensional separable algebraic hyperbolic space defined over $\F$. Finally, real hyperbolic spaces will be denoted without using $\R$ as a subscript, e.g. $\H^\infty = \H_\R^\infty$, $\B^J = \B_\R^J$, $\HH^\alpha = \HH_\R^\alpha$.
\end{convention}

\chapter{$\R$-trees, CAT(-1) spaces, and Gromov hyperbolic metric spaces} \label{sectiongeometry1}

In this chapter we review the theory of ``negative curvature'' in general metric spaces. A good reference for this subject is \cite{BridsonHaefliger}. We begin by defining the class of \emph{$\R$-trees}, the main class of examples we will talk about in this monograph other than the class of algebraic hyperbolic spaces, which we will discuss in more detail in Chapter \ref{sectionRtrees}. Next we will define CAT(-1) spaces, which are geodesic metric spaces whose triangles are ``thinner'' than the corresponding triangles in two-dimensional real hyperbolic space $\H^2$. Both algebraic hyperbolic spaces and $\R$-trees are examples of CAT(-1) spaces. The next level of generality considers Gromov hyperbolic metric spaces. After defining these spaces, we proceed to define the boundary $\del X$ of a hyperbolic metric space $X$, introducing the families of so-called \emph{visual metametrics} and \emph{extended visual metrics} on the bordification $\bord X := X\cup\del X$. We show that the bordification of an algebraic hyperbolic space $X$ is isomorphic to its closure $\bord X$ defined in Chapter \ref{sectionROSSONCTs}; under this isomorphism, the visual metric on $\del \B^\alpha$ is proportional to the Euclidean metric.

\bigskip
\section{Graphs and $\R$-trees}
\label{subsectiongraphs}

To motivate the definition of $\R$-trees we begin by defining simplicial trees, which requires first defining graphs.

\begin{definition}
\label{definitiongraphmetrization}
A \emph{weighted undirected graph} is a triple $(V,E,\ell)$, where $V$ is a nonempty set, $E\subset V\times V\butnot\{(x,x) : x\in V\}$ is invariant under the map $(x,y)\mapsto (y,x)$, and $\ell:E\to(0,\infty)$ is also invariant under $(x,y)\to (y,x)$. (If $\ell \equiv 1$, the graph is called \emph{unweighted}, and can be denoted simply $(V,E)$.) The graph is called \emph{connected} if for all $x,y\in V$, there exist $x = z_0,z_1,\ldots,z_n = y$ such that $(z_i,z_{i + 1}) \in E$ for all $i = 0,\ldots,n - 1$. If $(V,E,\ell)$ is connected, then the \emph{path metric} on $V$ is the metric
\begin{equation}
\label{pathmetric}
\dist_{E,\ell}(x,y) := \inf\left\{ \sum_{i = 0}^{n - 1} \ell(z_i,z_{i + 1}) \left|
\begin{split}
z_0 = x,\; z_n = y,\\
(z_i,z_{i + 1})\in E \all i = 0,\ldots,n - 1
\end{split}
\right.
\right\}
\end{equation}
The \emph{geometric realization} of the graph $(V,E,\ell)$ is the metric space
\[
X = X(V,E,\ell) = \left(V\cup\bigcup_{(v,w)\in E} [0,\ell(v,w)]\right)/\sim,
\]
where $\sim$ represents the following identifications:
\begin{align*}
v &\sim ((v,w),0) \all (v,w)\in E\\
((v,w),t) &\sim ((w,v),\ell(v,w) - t) \all (v,w)\in E \all t\in[0,\ell(v,w)]
\end{align*}
and the metric $\dist$ on $X$ is given by
\[
\dist\big(((v_0,v_1),t),((w_0,w_1),s)\big) = \min_{\substack{i\in\{0,1\}\\ j\in \{0,1\}}}\{ |t - i\ell(v_0,v_1)| + \dist(v_i,w_j) + |s - j\ell(w_0,w_1)|\}.
\]
(The geometric realization of a graph is sometimes also called a graph. In the sequel, we shall call it a \emph{geometric graph}.)
\end{definition}

\begin{example}[The Cayley graph of a group]
\label{examplecayleygraph}
Let $\Gamma$ be a group, and let $E_0\subset \Gamma$ be a generating set. (In most circumstances $E_0$ will be finite; there is an exception in Example \ref{examplenotstronglydiscrete} below.) Assume that $E_0 = E_0^{-1}$. The \emph{Cayley graph} of $\Gamma$ with respect to the generating set $E_0$ is the unweighted graph $(\Gamma,E)$, where
\begin{equation}
\label{Edef}
(\gamma,\beta)\in E \;\;\Leftrightarrow\;\; \gamma^{-1}\beta\in E_0.
\end{equation}
More generally, if $\ell_0:E_0\to (0,\infty)$ satisfies $\ell_0(g^{-1}) = \ell_0(g)$, the \emph{weighted Cayley graph} of $\Gamma$ with respect to the pair $(E_0,\ell_0)$ is the graph $(\Gamma,E,\ell)$, where $E$ is defined by \eqref{Edef}, and
\begin{equation}
\label{elldef}
\ell(\gamma,\beta) = \ell_0(\gamma^{-1}\beta).
\end{equation}
The path metric of a Cayley graph is called a \emph{Cayley metric}.
\end{example}

\begin{remark}
\label{remarknaturalaction}
The equations \eqref{Edef}, \eqref{elldef} guarantee that for each $\gamma\in\Gamma$, the map $\Gamma\ni\beta\to \gamma\beta\in\Gamma$ is an isometry of $\Gamma$ with respect to any Cayley metric. This isometry extends in a unique way to an isometry of the geometric Cayley graph $X = X(\Gamma,E,\ell)$. The map sending $\gamma$ to this isometry is a homomorphism from $\Gamma$ to $\Isom(X)$, and is called the \emph{natural action} of $\Gamma$ on $X$.
\end{remark}

\begin{remark}
\label{remarkpathmetricuniversal}
The path metric \eqref{pathmetric} satisfies the following \emph{universal property}: If $Y$ is a metric space and if $\phi:V\to Y$ satisfies $\dist(\phi(v),\phi(w)) \leq \ell(v,w)$ for every $(v,w)\in E$, then $\dist(\phi(v),\phi(w)) \leq \dist(v,w)$ for every $v,w\in V$.
\end{remark}

\begin{remark}
\label{remarkgeodesicmetric}
The main difference between the metric space $(V,\dist_{E,\ell})$ and the geometric graph $X = X(V,E,\ell)$ is that the latter is a \emph{geodesic} metric space. A metric space $X$ is said to be \emph{geodesic} if for every $p,q\in X$, there exists an isometric embedding $\pi:[t,s]\to X$ such that $\pi(t) = p$ and $\pi(s) = q$, for some $t,s\in\R$. The set $\pi([t,s])$ is denoted $\geo pq$ and is called a \emph{geodesic segment connecting $p$ and $q$}. The map $\pi$ is called a \emph{parameterization} of the geodesic segment $\geo pq$. (Note that although $\geo qp = \geo pq$, $\pi$ is not a parameterization of $\geo qp$.)

Warning: Although we may denote any geodesic segment connecting $p$ and $q$ by $\geo pq$, such a geodesic segment is not necessarily unique. A geodesic metric space $X$ is called \emph{uniquely geodesic} if for every $p,q\in X$, the geodesic segment connecting $p$ and $q$ is unique.
\end{remark}

\begin{notation}
\label{notationgeopqt}
If $\pi:[0,t_0]\to X$ is a parameterization of the geodesic segment $\geo pq$, then for each $t\in[0,t_0]$, $\geo pq_t$ denotes the point $\pi(t)$, i.e. the unique point on the geodesic segment $\geo pq$ such that $\dist(p,\geo pq_t) = t$.
\end{notation}

We are now ready to define the class of simplicial trees. Let $(V,E,\ell)$ be a weighted undirected graph. A \emph{cycle} in $(V,E,\ell)$ is a finite sequence of distinct vertices $v_1,\ldots,v_n\in V$, with $n\geq 3$, such that
\begin{equation}
\label{cycle}
(v_1,v_2), (v_2,v_3), \ldots, (v_{n - 1},v_n), (v_n,v_1)\in E.
\end{equation}
\begin{definition}
\label{definitionweightedtree}
A \emph{simplicial tree} is the geometric realization of a weighted undirected graph with no cycles. A \emph{$\Z$-tree} (or \emph{unweighted simplicial tree}, or just \emph{tree}\Footnote{However, in \cite{Tits}, the word ``trees'' is used to refer to what are now known as $\R$-trees.}) is the geometric realization of an unweighted undirected graph with no cycles.
\end{definition}

\begin{example}
Let $\F_2(\Z)$ denote the free group on two elements $\gamma_1,\gamma_2$. Let $E_0 = \{\gamma_1,\gamma_1^{-1},\gamma_2,\gamma_2^{-1}\}$. The geometric Cayley graph of $\F_2(\Z)$ with respect to the generating set $E_0$ is an unweighted simplicial tree.
\end{example}

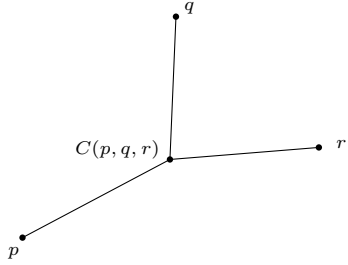
\begin{figure}
\begin{center}
\begin{tikzpicture}[line cap=round,line join=round,>=triangle 45,x=1.0cm,y=1.0cm]
\clip(0.22,0.49) rectangle (5.46,4.0);
\draw (0.8,0.8)-- (2.76,1.84);
\draw (2.76,1.84)-- (2.84,3.74);
\draw (2.76,1.84)-- (4.74,2);
\begin{scriptsize}
\fill [color=black] (0.8,0.8) circle (1.2pt);
\draw[color=black] (0.68,0.6) node {$p$};
\fill [color=black] (2.76,1.84) circle (1.2pt);
\draw[color=black] (2.06,1.98) node {$C(p,q,r)$};
\fill [color=black] (2.84,3.74) circle (1.2pt);
\draw[color=black] (3.02,3.86) node {$q$};
\fill [color=black] (4.74,2) circle (1.2pt);
\draw[color=black] (5.04,2.04) node {$r$};
\end{scriptsize}
\end{tikzpicture}
\caption{A geodesic triangle in an $\R$-tree}
\label{figureRtree}
\end{center}
\end{figure}

\begin{example}
\label{exampletreetriangle}
Let $V = \{\wbar C,\wbar p,\wbar q,\wbar r\}$, and fix $\ell_{\wbar p},\ell_{\wbar q},\ell_{\wbar r} > 0$. Let 
\begin{align*}
E &= \{(\wbar C,x), (x,\wbar C) : x = \wbar p,\wbar q,\wbar r\},\\
\ell(\wbar C,x) &= \ell(x,\wbar C) = \ell_x.
\end{align*} 
The geometric realization of the graph $(V,E,\ell)$ is a simplicial tree; see Figure \ref{figureRtree}. It will be denoted $\wbar \Delta = \Delta(\wbar p,\wbar q,\wbar r)$, and will be called a \emph{tree triangle}. For $x,y\in\{\wbar p,\wbar q,\wbar r\}$ distinct, the distance between $x$ and $y$ is given by
\[
\dist(x,y) = \ell_x + \ell_y.
\]
Solving for $\ell_{\wbar p}$ in terms of $\dist(\wbar p,\wbar q),\dist(\wbar p,\wbar r),\dist(\wbar q,\wbar r)$ gives
\begin{equation}
\label{gromovprecursor}
\ell_{\wbar p} = \dist(\wbar p,\wbar C) = \frac12[\dist(\wbar p,\wbar q) + \dist(\wbar p,\wbar r) - \dist(\wbar q,\wbar r)].
\end{equation}
\end{example}

\begin{definition}
\label{definitionRtree}
A metric space $X$ is an \emph{$\R$-tree} if for all $p,q,r\in X$, there exists a tree triangle $\wbar\Delta = \Delta(\wbar p,\wbar q,\wbar r)$ and an isometric embedding $\iota:\wbar\Delta\to X$ sending $\wbar p,\wbar q,\wbar r$ to $p,q,r$, respectively.
\end{definition}
\begin{definition}
\label{definitioncenterRtree}
Let $X$ be an $\R$-tree, fix $p,q,r\in X$, and let $\iota:\wbar\Delta\to X$ be as above. The point $C = C(p,q,r) := \iota(\wbar C)$ is called the \emph{center} of the geodesic triangle $\Delta = \Delta(p,q,r)$.
\end{definition}

As the name suggests, every simplicial tree is an $\R$-tree; the converse does not hold; see e.g. \cite[Example on p.50]{Chiswell}. Before we can prove that every simplicial tree is an $\R$-tree, we will need a lemma:

\begin{lemma}[Cf. {\cite[p.29]{Chiswell}}]
\label{lemmaRtreeequivalent}
Let $X$ be a metric space. The following are equivalent:
\begin{itemize}
\item[(A)] $X$ is an $\R$-tree.
\item[(B)] There exists a collection of geodesics $\GG$, with the following properties:
\begin{itemize}
\item[(BI)] For each $x,y\in X$, there is a geodesic $\geo xy\in\GG$ connecting $x$ and $y$.
\item[(BII)] Given $\geo xy\in\GG$ and $z,w\in \geo xy$, we have $\geo zw\in\GG$, where $\geo zw$ is interpreted as the set of points in $\geo xy$ which lie between $z$ and $w$.
\item[(BIII)] Given $x_1,x_2,x_3\in X$ distinct and geodesics $\geo{x_1}{x_2},\geo{x_1}{x_3},\geo{x_2}{x_3}\in\GG$, at least one pair of the geodesics $\geo{x_i}{x_j}$, $i\neq j$, has a nontrivial intersection. More precisely, there exist distinct $i,j,k\in\{1,2,3\}$ such that
\[
\geo{x_i}{x_j} \cap \geo{x_i}{x_k} \propersupset \{x_i\}.
\]
\end{itemize}
\end{itemize}
\end{lemma}
\begin{proof}[Proof of \text{(A)} \implies \text{(B)}]
Note that (BI) and (BII) are true for any uniquely geodesic metric space. Given $x_1,x_2,x_3$ distinct, let $C$ be the center. Then $x_i\neq C$ for some $i$; without loss of generality $x_1\neq C$. Then
\[
\geo{x_1}{x_2} \cap \geo{x_1}{x_3} = \geo{x_1}{C} \propersupset \{x_1\}.
\]
\end{proof}
\begin{proof}[Proof of \text{(B)} \implies \text{(A)}]
We first show that given points $x_1,x_2,x_3\in X$ and geodesics $\geo{x_1}{x_2},\geo{x_1}{x_3},\geo{x_2}{x_3}\in\GG$, the intersection $\bigcap_{i\neq j} \geo{x_i}{x_j}$ is nonempty. Indeed, suppose not. For $i = 2,3$ let $\gamma_i:[0,\dist(x_1,x_i)]\to X$ be a parameterization of $\geo{x_1}{x_i}$, and let
\[
t_1 = \max\{t\geq 0 : \gamma_2(t) = \gamma_3(t)\}.
\]
By replacing $x$ with $\gamma_2(t_1) = \gamma_3(t_1)$ and using (BII), we may without loss of generality assume that $t_1 = 0$, or equivalently that $\geo{x_1}{x_2}\cap \geo{x_1}{x_3} = \{x_1\}$. Similarly, we may without loss of generality assume that $\geo{x_2}{x_1}\cap \geo{x_2}{x_3} = \{x_2\}$ and $\geo{x_3}{x_1}\cap\geo{x_3}{x_2} = \{x_3\}$. But then (BIII) implies that $x_1,x_2,x_3$ cannot be all distinct. This immediately implies that $\bigcap_{i\neq j} \geo{x_i}{x_j} \neq \emptyset$.

To complete the proof, we must show that $X$ is uniquely geodesic. Indeed suppose that for some $x_1,x_2\in X$, there is more than one geodesic connecting $x_1$ and $x_2$. Let $\geo{x_1}{x_2}\in\GG$ be a geodesic connecting $x_1$ and $x_2$, and let $\geo{x_1}{x_2}'$ be another geodesic connecting $x_1$ and $x_2$. Then there exists $x_3\in\geo{x_1}{x_2}'\butnot\geo{x_1}{x_2}$. By the above paragraph, there exists $w\in\bigcap_{i\neq j} \geo{x_i}{x_j}$. Since $w\in \geo{x_i}{x_3}$, we have
\begin{equation}
\label{dxiw}
\dist(x_i,w) \leq \dist(x_i,x_3).
\end{equation}
On the other hand, since $w\in \geo{x_1}{x_2}$ and $x_3\in\geo{x_1}{x_2}'$, we have
\[
\dist(x_1,x_2) = \dist(x_1,w) + \dist(x_2,w) \leq \dist(x_1,x_3) + \dist(x_2,x_3) = \dist(x_1,x_2).
\]
It follows that equality holds in \eqref{dxiw}, i.e. $\dist(x_i,w) = \dist(x_i,x_3)$. Since $w\in \geo{x_i}{x_3}$, this implies $w = x_3$. But then $x_3 = w\in \geo{x_1}{x_2}$, a contradiction.
\end{proof}

\begin{corollary}
Every simplicial tree is an $\R$-tree.
\end{corollary}
\begin{proof}
Let $X = X(V,E,\ell)$ be a simplicial tree, and let $\GG$ be the collection of all geodesics; then (BI) and (BII) both hold. By contradiction, suppose that there exist points $x_1,x_2,x_3\in X$ such that $\geo{x_i}{x_j}\cap\geo{x_i}{x_k} = \{x_i\}$ for all distinct $i,j,k\in\{1,2,3\}$. Then the path $\geo{x_1}{x_2}\cup\geo{x_2}{x_3}\cup\geo{x_3}{x_1}$ is equal to the union of the edges of a cycle of the graph $(V,E,\ell)$. This is a contradiction.
\end{proof}

We shall investigate $\R$-trees in more detail in Chapter \ref{sectionRtrees}, where we will give various examples of $\R$-trees together with groups acting isometrically on them.





\bigskip
\section{CAT(-1) spaces}
\label{subsectionCAT}

The following definitions have been modified from \cite[p.158]{BridsonHaefliger}, to which the reader is referred for more details.

A \emph{geodesic triangle} in $X$ consists of three points $p,q,r\in X$ (the \emph{vertices} of the triangle) together with a choice of three geodesic segments $\geo pq$, $\geo qr$, and $\geo rp$ joining them (the \emph{sides}). Such a geodesic triangle will be denoted $\Delta(p,q,r)$, although we note that this could cause ambiguity if $X$ is not uniquely geodesic. Although formally $\Delta(p,q,r)$ is an element of $X^3\times \P(X)^3$, we will sometimes identify $\Delta(p,q,r)$ with the set $\geo pq \cup \geo qr\cap \geo rp \subset X$, writing $x\in\Delta(p,q,r)$ if $x\in \geo pq \cup \geo qr\cup \geo rp$.

A triangle $\wbar{\Delta} = \Delta(\wbar p,\wbar q,\wbar r)$ in $\H^2$ is called a \emph{comparison triangle} for $\Delta = \Delta(p,q,r)$ if $\dist(\wbar p, \wbar q) = \dist(p,q)$, $\dist(\wbar q,\wbar r) = \dist(q, r)$, and $\dist(\wbar p,\wbar r) = \dist(p,r)$. Any triangle admits a comparison triangle, unique up to isometry. For any point $x\in\geo pq$, we define its \emph{comparison point} $\wbar x\in\geo{\wbar p}{\wbar q}$ to be the unique point such that $\dist(\wbar x,\wbar p) = \dist(x,p)$ and $\dist(\wbar x,\wbar q) = \dist(x,q)$. In the notation above, the comparison point of $\geo pq_t$ is equal to $\geo{\wbar p}{\wbar q}_t$ for all $t\in[0,\dist(p,q)] = [0,\dist(\wbar p,\wbar q)]$. For $x\in \geo qr$ and $x\in \geo rp$, the comparison point is defined similarly.

Let $X$ be a metric space and let $\Delta$ be a geodesic triangle in $X$. We say that $\Delta$ \emph{satisfies the CAT(-1) inequality} if for all $x,y\in\Delta$,
\begin{equation}
\label{CAT}
\dist(x, y) \leq \dist(\wbar x, \wbar y),
\end{equation}
where $\wbar x$ and $\wbar y$ are any\Footnote{The comparison points $\wbar x$ and $\wbar y$ may not be uniquely determined if either $x$ or $y$ lies on two sides of the triangle simultaneously.} comparison points for $x$ and $y$, respectively. Intuitively, $\Delta$ satisfies the CAT(-1) inequality if it is ``thinner'' than its comparison triangle $\wbar{\Delta}$.

\begin{definition}
\label{definitionCAT}
$X$ is a \emph{CAT(-1) space} if it is a geodesic metric space and if all of its geodesic triangles satisfy the CAT(-1) inequality.
\end{definition}
\begin{observation}[{\cite[Proposition II.1.4(1)]{BridsonHaefliger}}]
\label{observationCATuniquelygeo}
CAT(-1) spaces are uniquely geodesic.
\end{observation}
\begin{proof}
Let $X$ be a CAT(-1) space, and suppose that two points $p,q\in X$ are connected by two geodesic segments $\geo pq$ and $\geo pq'$. Fix $t\in [0,\dist(p,q)]$ and let $x = \geo pq_t$, $x' = \geo pq_t'$. Consider the triangle $\Delta(p,q,x)$ determined by the geodesic segments $\geo pq'$, $\geo px$, and $\geo xq$, and a comparison triangle $\Delta(\wbar p,\wbar q,\wbar x)$. Then $x$ and $x'$ have the same comparison point $\wbar x$, so by the CAT(-1) inequality
\[
\dist(x,x')\leq \dist(\wbar x,\wbar x) = 0,
\]
and thus $x = x'$. Since $t$ was arbitrary, it follows that $\geo pq = \geo pq'$. Since $\geo pq'$ was arbitrary, $\geo pq$ is the unique geodesic segment connecting $p$ and $q$.
\end{proof}

\subsection{Examples of CAT(-1) spaces}
\label{subsubsectionexamplesCAT}
In this text we concentrate on two main examples of CAT(-1) spaces: algebraic hyperbolic spaces and $\R$-trees. We therefore begin by proving the following result which implies that algebraic hyperbolic spaces are CAT(-1):

\begin{proposition}
\label{propositioncurvatureCAT}
Any Riemannian manifold (finite- or infinite-dimensional) with sectional curvature bounded above by $-1$ is a CAT(-1) space.
\end{proposition}
\begin{proof}
The finite-dimensional case is proven in \cite[Theorem II.1A.6]{BridsonHaefliger}. The infinite-dimensional follows upon augmenting the finite-dimensional proof with the infinite-dimensional Cartan--Hadamard theorem \cite[IX, Theorem 3.8]{Lang_differential_geometry} to guarantee surjectivity of the exponential map.
\end{proof}

Since algebraic hyperbolic spaces have sectional curvature bounded between $-4$ and $-1$ (e.g. \cite[Corollary of Proposition 4]{Heintze}; see also \cite[Lemmas 2.3, 2.7, and 2.11]{Quint}
), the following corollary is immediate:

\begin{corollary}
\label{corollaryROSSONCTCAT}
Every algebraic hyperbolic space is a CAT(-1) space.
\end{corollary}

\begin{remark}
One can prove Corollary \ref{corollaryROSSONCTCAT} without using the full strength of Proposition \ref{propositioncurvatureCAT}. Indeed, any geodesic triangle in an algebraic hyperbolic space is isometric to a geodesic triangle in $\H_\F^2$ for some $\F\in\{\R,\C,\Q\}$. Since $\H_\F^2$ is finite-dimensional, thinness of its geodesic triangles follows from the finite-dimensional version of Proposition \ref{propositioncurvatureCAT}.
\end{remark}

\begin{observation}
\label{observationRtreesCAT}
$\R$-trees are CAT(-1).
\end{observation}
\begin{proof}
First of all, an argument similar to the proof of Observation \ref{observationCATuniquelygeo} shows that $\R$-trees are uniquely geodesic, justifying Figure \ref{figureRtree}. In particular, if $\Delta(p,q,r)$ is a geodesic triangle in an $\R$-tree and if $C = C(p,q,r)$ then $\geo pq = \geo pC\cup\geo qC$, $\geo qr = \geo qC\cup\geo rC$, and $\geo rp = \geo rC\cup\geo pC$. It follows that any two points $x,y\in \Delta$ share a side in common, without loss of generality say $x,y\in\geo pq$. Then
\[
\dist(x,y) = \dist(p,q) - \dist(x,p) - \dist(y,q) = \dist(\wbar p,\wbar q) - \dist(\wbar x,\wbar p) - \dist(\wbar y,\wbar q) \leq \dist(\wbar x,\wbar y).
\]
\end{proof}

In a sense $\R$-trees are the ``most negatively curved spaces''; although we did not define the notion of a CAT($\kappa$) space, $\R$-trees are CAT($\kappa$) for every $\kappa\in\R$.
\bigskip
\section{Gromov hyperbolic metric spaces}

We now come to the theory of Gromov hyperbolic metric spaces. In a sense, Gromov hyperbolic metric spaces are those which are ``approximately $\R$-trees''. A good reference for this section is \cite{Vaisala}.

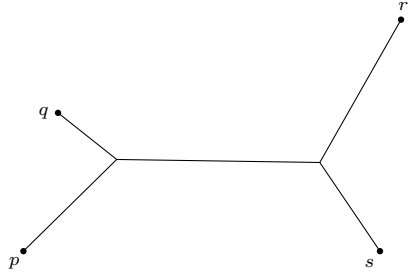
\begin{figure}
\begin{center}
\begin{tikzpicture}[line cap=round,line join=round,>=triangle 45,x=1.0cm,y=1.0cm]
\clip(0.56,1.2) rectangle (7.44,4.8);
\draw (2.02,3.3)-- (2.8,2.68);
\draw (2.8,2.68)-- (1.56,1.46);
\draw (2.8,2.68)-- (5.5,2.64);
\draw (5.5,2.64)-- (6.58,4.54);
\draw (5.5,2.64)-- (6.3,1.46);
\begin{scriptsize}
\fill [color=black] (2.02,3.3) circle (1.2pt);
\draw[color=black] (1.83,3.3) node {$q$};
\fill [color=black] (1.56,1.46) circle (1.2pt);
\draw[color=black] (1.44,1.3) node {$p$};
\fill [color=black] (6.58,4.54) circle (1.2pt);
\draw[color=black] (6.61,4.72) node {$r$};
\fill [color=black] (6.3,1.46) circle (1.2pt);
\draw[color=black] (6.16,1.3) node {$s$};
\end{scriptsize}
\end{tikzpicture}
\caption{A quadruple of points in an $\R$-tree}
\label{figurequadruple}
\end{center}
\end{figure}

For any three numbers $d_{pq},d_{qr},d_{rp}\geq 0$ satisfying the triangle inequality, there exists an $\R$-tree $X$ and three points $p,q,r\in X$ such that $\dist(p,q) = d_{pq}$, etc. Thus in some sense looking at triples ``does not tell you'' that you are looking at an $\R$-tree. Now let us look at quadruples. A quadruple $(p,q,r,s)$ in an $\R$-tree $X$ looks something like Figure \ref{figurequadruple}. Of course, the points $p,q,r,s\in X$ could be arranged in any order. However, let us consider them the way that they are arranged in Figure \ref{figurequadruple} and note that

\begin{equation}
\label{centersequal}
C(p,q,r) = C(p,q,s).
\end{equation}
In order to write this equality in terms of distances, we need some way of measuring the distance from the vertex of a geodesic triangle to its center.
\begin{observation}
If $\Delta(p,q,r)$ is a geodesic triangle in an $\R$-tree then the distance from the vertex $p$ to the center $C(p,q,r)$, i.e. $\dist(p,C(p,q,r))$, is equal to
\begin{equation}
\label{gromovproduct}
\lb q|r\rb_p := \frac12[\dist(p,q) + \dist(p,r) - \dist(q,r)].
\end{equation}
\end{observation}
The expression $\lb q|r\rb_p$ is called the \emph{Gromov product} of $q$ and $r$ with respect to $p$, and it makes sense in any metric space. It can be thought of as measuring the ``defect in the triangle inequality''; indeed, the triangle inequality is exactly what assures that $\lb q|r\rb_p\geq 0$ for all $p,q,r\in X$.

Now \eqref{centersequal} implies that
\[
\lb q|r\rb_p = \lb q|s\rb_p \leq \lb r|s\rb_p.
\]
(The last inequality does not follow from \eqref{centersequal} but it may be seen from Figure \ref{figurequadruple}.) However, since the arrangement of points was arbitrary we do not know which two Gromov products will be equal and which one will be larger. An inequality which captures all possibilities is
\begin{equation}
\label{gromovinRtrees}
\lb q|r\rb_p \geq \min(\lb q|s\rb_p,\lb r|s\rb_p).
\end{equation}
Now, as mentioned before, we will define hyperbolic metric spaces as those which are ``approximately $\R$-trees''. Thus they will satisfy \eqref{gromovinRtrees} with an asymptotic.
\begin{definition}
\label{definitiongromovhyperbolic}
A metric space $X$ is called \emph{hyperbolic} (or \emph{Gromov hyperbolic}) if for every four points $x,y,z,w\in X$ we have 
\begin{equation}
\label{gromov}
\lb x|z\rb_w \gtrsim_\plus \min(\lb x|y\rb_w,\lb y|z\rb_w),
\end{equation}
We refer to \eqref{gromov} as \emph{Gromov's inequality}.
\end{definition}
From the above discussion, every $\R$-tree is Gromov hyperbolic with an implied constant of $0$ in \eqref{gromov}. (This can also be deduced from Proposition \ref{propositionCATimpliesGromov} below.)

Note that many authors require $X$ to be a geodesic metric space in order to be hyperbolic; we do not. If $X$ is a geodesic metric space, then the condition of hyperbolicity can be reformulated in several different ways, including the \emph{thin triangles condition}; for details, see \cite[\6 III.H.1]{BridsonHaefliger} or Section \ref{subsectionrips} below.

It will be convenient for us to make a list of several identities satisfied by the Gromov product. For each $z\in X$, let $\busemann_z$ denote the \emph{Busemann function}
\begin{equation}
\label{busemanndef}
\busemann_z(x,y) := \dist(z,x) - \dist(z,y).
\end{equation}

\begin{proposition}
\label{propositionbasicidentities}
The Gromov product and Busemann function satisfy the following identities and inequalities:
\begin{align*}
\tag{a} \lb x|y\rb_z &= \lb y|x\rb_z \\
\tag{b} \dist(y,z) &= \lb y|x\rb_z + \lb z|x\rb_y\\
\tag{c} 0\leq \lb x|y\rb_z &\leq \min(\dist(x,z),\dist(y,z))\\
\tag{d} \lb x|y\rb_z &\leq \lb x|y\rb_w + \dist(z,w)\\
\tag{e} \lb x|y\rb_w &\leq \lb x|z\rb_w + \dist(y,z)\\
\tag{f} |\busemann_x(z,w)| &\leq \dist(z,w)\\
\tag{g} \lb x|y\rb_z &= \lb x|y\rb_w + \frac12[\busemann_x(z,w) + \busemann_y(z,w)]\\
\tag{h} \lb x|y\rb_z &= \frac12[\dist(x,z) + \busemann_y(z,x)]\\
\tag{j} \busemann_x(y,z) &= \lb z|x\rb_y - \lb y|x\rb_z \\
\tag{k} \lb x|y\rb_z &= \lb x|y\rb_w + \dist(z,w) - \lb x|z\rb_w - \lb y|z\rb_w\\
\tag{l} \lb x|y\rb_w &= \lb x|z\rb_w + \frac12[\busemann_w(y,z) - \busemann_x(y,z)]
\end{align*}
\end{proposition}
The proof is a straightforward computation. We remark that (a)-(e) may be found in \cite[Lemma 2.8]{Vaisala}.

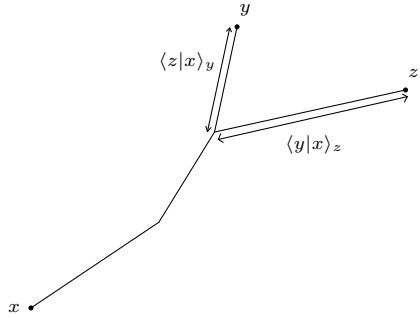
\begin{figure}
\begin{center}
\begin{tikzpicture}[line cap=round,line join=round,>=triangle 45,x=1.0cm,y=1.0cm]
\clip(-1.54,-0.63) rectangle (5.02,3.6);
\draw (-0.68,-0.48)-- (1.02,0.66);
\draw (1.02,0.66)-- (1.76,1.86);
\draw (1.76,1.86)-- (2.06,3.26);
\draw (1.76,1.86)-- (4.3,2.42);
\draw (1.9621294832739777,3.244449991852601)-- (1.6662459503620535,1.8663069374640353);
\draw (1.9621294832739777,3.244449991852601)-- (1.920651867546605,3.2029723761252282);
\draw (1.9621294832739777,3.244449991852601)-- (1.9850941143735463,3.194322722245245);
\draw (1.6662459503620535,1.8663069374640353)-- (1.645781680327839,1.9358422158191166);
\draw (1.6662459503620535,1.8663069374640353)-- (1.712620501296158,1.919132510577037);
\draw (1.811554396247592,1.7618958525911042)-- (4.322472096464471,2.330897992557095);
\draw (1.811554396247592,1.7618958525911042)-- (1.8414665959029357,1.8362946416296355);
\draw (1.811554396247592,1.7618958525911042)-- (1.8757672524481015,1.7219591198124153);
\draw (4.322472096464471,2.330897992557095)-- (4.258520422692855,2.3649896261428003);
\draw (4.322472096464471,2.330897992557095)-- (4.284871880487056,2.2595837949659954);
\begin{scriptsize}
\draw [fill=black] (-0.68,-0.48) circle (.75pt);
\draw[color=black] (-0.9047548310928106,-0.48743844483680343) node {$x$};
\draw [fill=black] (2.06,3.26) circle (.75pt);
\draw[color=black] (2.1685467646948235,3.4849673412078808) node {$y$};
\draw[color=black] (1.3983054592304987,2.7983786868297873) node {$\lb z|x\rb_y$};
\draw [fill=black] (4.3,2.42) circle (.75pt);
\draw[color=black] (4.408133565880492,2.6512525466059103) node {$z$};
\draw[color=black] (3.083998303865608,1.701190403090175) node {$\lb y|x\rb_z$};
\end{scriptsize}
\end{tikzpicture}
\caption[Expressing distance via Gromov products in an $\R$-tree]{An illustration of (b) of Proposition \ref{propositionbasicidentities} in an $\R$-tree.}
\end{center}
\end{figure}

\subsection{Examples of Gromov hyperbolic metric spaces}

\begin{proposition}[Proven in Section \ref{subsectionGromovinH}]
\label{propositionCATimpliesGromov}
Every CAT(-1) space (in particular every algebraic hyperbolic space) is Gromov hyperbolic. In fact, if $X$ is a CAT(-1) space then for every four points $x,y,z,w\in X$ we have 
\begin{equation}
\label{gromovforCAT}
e^{-\lb x|z\rb_w} \leq e^{-\lb x|y\rb_w} + e^{-\lb y|z\rb_w},
\end{equation}
and so $X$ satisfies \eqref{gromov} with an implied constant of $\log(2)$.
\end{proposition}
\begin{remark}
The first assertion of Proposition \ref{propositionCATimpliesGromov}, namely, that CAT(-1) spaces are Gromov hyperbolic, is \cite[Proposition III.H.1.2]{BridsonHaefliger}. The inequality \eqref{gromovforCAT} in the case where $x,y,z\in\del X$ and $w\in X$ can be found in \cite[Th\'eor\`eme 2.5.1]{Bourdon}.
\end{remark}

\begin{definition}
\label{definitionstronglyhyperbolic}
A space $X$ satisfying the conclusion of Proposition \ref{propositionCATimpliesGromov} is said to be \emph{strongly hyperbolic}.
\end{definition}
Note that
\[
\text{$\R$-tree \implies CAT(-1) \implies Strongly hyperbolic \implies Hyperbolic.}
\]

A large class of examples of hyperbolic metric spaces which are not CAT(-1) is furnished by the Cayley graphs of finitely presented groups. Indeed, we have the following:

\begin{theorem}[{\cite[p.78]{Gromov4}}, \cite{Olshanskii}; see also \cite{Champetier}]
\label{theoremolshanskii}
Fix $k\geq 2$ and an alphabet $A = \{a_1^{\pm 1},a_2^{\pm 1},\cdots,a_k^{\pm 1}\}$. Fix $i\in\Namer$ and a sequence of positive integers $(n_1,\cdots,n_i)$. Let $N = N(k,i,n_1,\cdots,n_i)$ be the number of group presentations $G = \lb a_1,\cdots, a_k|r_1,\cdots,r_i\rb$ such that $r_1,\cdots,r_i$ are reduced words in the alphabet $A$ such that the length of $r_j$ is $n_j$ for $j = 1,2,\cdots,i$. If $N_h$ is the number of groups in this collection whose Cayley graphs are hyperbolic and if $n = \min(n_1,\cdots,n_i)$ then $\lim_{n\to\infty}N_h/N = 1$.
\end{theorem}

This theorem says that in some sense, ``almost every'' finitely presented group is hyperbolic.


If one has a hyperbolic metric space $X$, there are two ways to get another hyperbolic metric space from $X$, one trivial and one nontrivial.

\begin{observation}
\label{observationsubspace}
Any subspace of a hyperbolic metric space is hyperbolic. Any subspace of a strongly hyperbolic metric space is strongly hyperbolic.
\end{observation}

To describe the other method we need to define the notion of a \emph{quasi-isometric embedding}.
\begin{definition}
\label{definitionquasiisometry}
Let $(X_1,\dist_1)$ and $(X_2,\dist_2)$ be metric spaces. A map $\Phi:X_1\to X_2$ is a \emph{quasi-isometric embedding} if for every $x,y\in X_1$
\[
\dist_2(\Phi(x),\Phi(y)) \asymp_{\plus,\times} \dist_1(x,y).
\]
A quasi-isometric embedding $\Phi$ is called a \emph{quasi-isometry} if its image $\Phi(X_1)$ is cobounded in $X_2$, that is, if there exists $R > 0$ such that $\Phi(X_1)$ is $R$-dense in $X_2$, meaning that $\min_{x\in X_2} \dist(x,\Phi(X_1)) \leq R$. In this case, the spaces $X_1$ and $X_2$ are said to be \emph{quasi-isometric}.
\end{definition}

\begin{theorem}[{\cite[Theorem III.H.1.9]{BridsonHaefliger}}]
\label{theoremquasiisometric}
Any geodesic metric space which can be quasi-isometrically embedded into a geodesic hyperbolic metric space is also a hyperbolic metric space.
\end{theorem}

\begin{remark}
Theorem \ref{theoremquasiisometric} is not true if the hypothesis of geodesicity is dropped. For example, $\R$ is quasi-isometric to $\Rplus\times\{0\}\cup\{0\}\times\Rplus\subset\R^2$, but the former is hyperbolic and the latter is not.
\end{remark}

There are many more examples of hyperbolic metric spaces which we will not discuss; cf. the list in \6\ref{subsubsectionlist}.

\bigskip
\section{The boundary of a hyperbolic metric space}
\label{subsectionboundary}

In this section we define the Gromov boundary of a hyperbolic metric space $X$. The construction will depend on a distinguished point $\zero\in X$, but the resulting space will be independent of which point is chosen. If $X$ is an $\R$-tree, then the boundary of $X$ will turn out to be the set of infinite branches through $X$, i.e. the set of all isometric embeddings $\pi:\Rplus\to X$ sending $0$ to $\zero$, where $\zero\in X$ is a distinguished fixed point. If $X$ is an algebraic hyperbolic space, then the boundary of $X$ will turn out to be isomorphic to the space $\del X$ defined in Chapter \ref{sectionROSSONCTs}.

\begin{figure}
\begin{center}
\begin{tikzpicture}[line cap=round,line join=round,>=triangle 45,x=1.0cm,y=1.0cm]
\clip(-2.74,-1.0) rectangle (7.98,3.16);
\draw (-1.8,-0.78)-- (0.52,0.72);
\draw (0.52,0.72)-- (0.7,2.3);
\draw (0.52,0.72)-- (2.28,0.52);
\draw (2.28,0.52)-- (2.96,-0.64);
\draw (2.28,0.52)-- (3.46,1.62);
\draw (3.46,1.62)-- (3.7,2.48);
\draw (3.46,1.62)-- (4.78,2.16);
\draw (4.78,2.16)-- (5.38,2.92);
\draw (4.78,2.16)-- (6.8,2.36);
\begin{scriptsize}
\fill [color=black] (-1.8,-0.78) circle (1pt);
\draw[color=black] (-2,-0.86) node {$\zero$};
\fill [color=black] (0.7,2.3) circle (1pt);
\draw[color=black] (0.9,2.48) node {$x_1$};
\fill [color=black] (2.96,-0.64) circle (1pt);
\draw[color=black] (3.32,-0.72) node {$x_2$};
\fill [color=black] (3.7,2.48) circle (1pt);
\draw[color=black] (3.96,2.74) node {$x_3$};
\fill [color=black] (5.38,2.92) circle (1pt);
\draw[color=black] (5.74,3.08) node {$x_5$};
\fill [color=black] (6.8,2.36) circle (1pt);
\draw[color=black] (7.36,2.4) node {$\ldots$};
\end{scriptsize}
\end{tikzpicture}
\caption{A Gromov sequence in an $\R$-tree}
\label{figuregromovseq}
\end{center}
\end{figure}
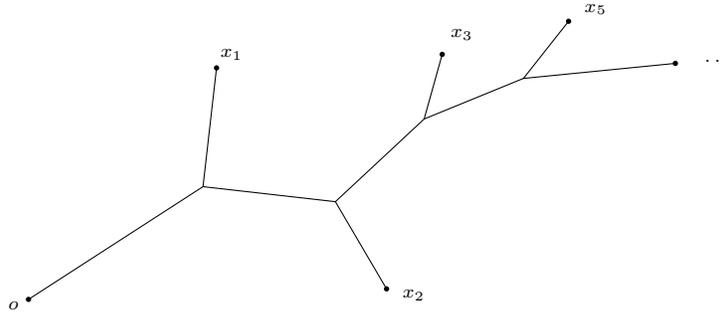

To motivate the definition of the boundary, suppose that $X$ is an $\R$-tree. An infinite branch through $X$ can be approximated by finite branches which agree on longer and longer segments. Suppose that $(\geo\zero{x_n})_1^\infty$ is a sequence of geodesic segments. For each $n,m\in\Namer$, the length of the intersection of $\geo\zero{x_n}$ and $\geo\zero{x_m}$ is equal to $\dist(\zero,C(\zero,x_n,x_m))$, which in turn is equal to $\lb x_n|x_m\rb_\zero$. Thus, the sequence $(\geo\zero{x_n})_1^\infty$ converges to an infinite geodesic if and only if
\begin{equation}
\label{gromovsequence}
\lb x_n|x_m\rb_\zero \tendsto{n,m} \infty.
\end{equation}
(Cf. Figure \ref{figuregromovseq}.) The formula \eqref{gromovsequence} is reminiscent of the definition of a Cauchy sequence. This intuition will be made explicit in Section \ref{subsectionmetametrics}, where we will introduce a \emph{metametric} on $X$ with the property that a sequence in $X$ satisfies \eqref{gromovsequence} if and only if it is Cauchy with respect to this metametric.

\begin{definition}
\label{definitiongromovsequence}
A sequence $(x_n)_1^\infty$ in $X$ for which \eqref{gromovsequence} holds is called a \emph{Gromov sequence}. Two Gromov sequences $(x_n)_1^\infty$ and $(y_n)_1^\infty$ are called \emph{equivalent} if
\[
\lb x_n|y_n\rb_\zero \tendsto n \infty,
\]
or equivalently if
\[
\lb x_n|y_m\rb_\zero \tendsto{n,m} \infty.
\]
In this case, we write $(x_n)_1^\infty \sim (y_n)_1^\infty$. It is readily verified using Gromov's inequality that $\sim$ is an equivalence relation on the set of Gromov sequences in $X$. We will denote the class of sequences equivalent to a given sequence $(x_n)_1^\infty$ by $[(x_n)_1^\infty]$.
\end{definition}

\begin{definition}
\label{definitiongromovboundary}
The \emph{Gromov boundary} of $X$ is the set of Gromov sequences modulo equivalence. It is denoted $\del X$. The \emph{Gromov closure} or \emph{bordification} of $X$ is the disjoint union $\bord X := X\cup\del X$. 
\end{definition}

\begin{remark}
\label{remarkambiguityboundary}
If $X$ is an algebraic hyperbolic space, then this notation causes some ambiguity, since it is not clear whether $\del X$ represents the Gromov boundary of $X$, or rather the topological boundary of $X$ as in Chapter \ref{sectionROSSONCTs}. This ambiguity will be resolved in \6\ref{subsubsectionboundariesequivalent} below when it is shown that the two bordifications are isomorphic.
\end{remark}

\begin{remark}
\label{remarkambiguityboundary2}
In the literature, the ideal boundary of a hyperbolic metric space is often taken to be the set of equivalence classes of geodesic rays under asymptotic equivalence, rather than the set of equivalence classes of Gromov sequences (e.g. \cite[p.427]{BridsonHaefliger}). If $X$ is proper and geodesic, then these two notions are equivalent \cite[Lemma III.H.3.13]{BridsonHaefliger}, but in general they may be different.
\end{remark}

\begin{remark}
\label{remarkbasepoint1}
By (d) of Proposition \ref{propositionbasicidentities}, the concepts of Gromov sequence and equivalence do not depend on the basepoint $\zero$. In particular, the Gromov boundary $\del X$ is independent of $\zero$.
\end{remark}

\subsection{Extending the Gromov product to the boundary}
We now wish to extend the Gromov product and Busemann function to the boundary ``by continuity''. Fix $\xi,\eta\in\del X$ and $z\in X$. Ideally, we would like to define $\lb \xi|\eta\rb_z$ to be
\begin{equation}
\label{gromovproductdef}
\lim_{n,m\to\infty}\lb x_n|y_m\rb_z ,
\end{equation}
where $(x_n)_1^\infty\in\xi$ and $(y_m)_1^\infty\in\eta$. (The definition would then have to be shown independent of which sequences were chosen.) The naive definition \eqref{gromovproductdef} does not work, because the limit \eqref{gromovproductdef} does not necessarily exist:
\begin{example}
\label{examplegromovnotcont}
Let
\[
X = \{\xx\in\R^2: x_2\in[0,1]\}
\]
be interpreted as a subspace of $\R^2$ with the $L^1$ metric. Then $X$ is a hyperbolic metric space, since it contains the cobounded hyperbolic metric space $\R\times\{0\}$. Its Gromov boundary consists of two points $-\infty$ and $+\infty$, which are the limits of $\xx$ as $x_1$ approaches $-\infty$ or $+\infty$, respectively. Let $\yy = (0,1)$ and $\zz = (1,0)$. Then for all $\xx\in X$, $\lb \xx|\yy\rb_\zz = x_2$. In particular, we can find a sequence $\xx_n\to +\infty$ such that $\lim_{n\to\infty} \lb \xx_n|\yy\rb_\zz$ does not exist.
\end{example}

Fortunately, the limit \eqref{gromovproductdef} ``exists up to a constant'':
\begin{lemma}
\label{lemmawelldefined}
Let $(x_n)_1^\infty$ and $(y_m)_1^\infty$ be Gromov sequences, and fix $y,z\in X$. Then
\begin{align} \label{welldefined1}
\liminf_{n,m\to\infty}\lb x_n|y_m\rb_z &\asymp_\plus \limsup_{n,m\to\infty}\lb x_n|y_m\rb_z\\ \label{welldefined2}
\liminf_{n\to\infty}\lb x_n|y\rb_z &\asymp_\plus \limsup_{n\to\infty}\lb x_n|y\rb_z ,
\end{align}
with equality if $X$ is strongly hyperbolic.
\end{lemma}
Note that except for the statement about strongly hyperbolic spaces, this lemma is simply \cite[Lemma 5.6]{Vaisala}.
\begin{proof}[Proof of Lemma \ref{lemmawelldefined}]
Fix $n_1,n_2,m_1,m_2\in\Namer$. By Gromov's inequality
\[
\lb x_{n_1}|y_{m_1}\rb_z \gtrsim_\plus \min(\lb x_{n_2}|y_{m_2}\rb_z,\lb x_{n_1}|x_{n_2}\rb_z,\lb y_{m_1}|y_{m_2}\rb_z).
\]
Taking the liminf over $n_1,m_1$ and the limsup over $n_2,m_2$ gives
\begin{align*}
&\liminf_{n,m\to\infty}\lb x_n|y_m \rb_z\\ 
&\gtrsim_\plus \min\left(\limsup_{n,m\to\infty}\lb x_n|y_m \rb_z,\liminf_{n_1,n_2\to\infty}\lb x_{n_1}|x_{n_2}\rb_z,\liminf_{m_1,m_2\to\infty}\lb y_{m_1}|y_{m_2}\rb_z\right) \noreason\\
&=_\pt \limsup_{n,m\to\infty}\lb x_n|y_m\rb_z, \phantom{XxX}\since{$(x_n)_1^\infty$ and $(y_m)_1^\infty$ are Gromov}
\end{align*}
demonstrating \eqref{welldefined1}. On the other hand, suppose that $X$ is strongly hyperbolic. Then by \eqref{gromovforCAT} we have
\[
\exp\big(-\lb x_{n_1}|y_{m_1}\rb_z\big) \leq \exp\big(-\lb x_{n_2}|y_{m_2}\rb_z\big) + \exp\big(-\lb x_{n_1}|x_{n_2}\rb_z\big) + \exp\big(-\lb y_{m_1}|y_{m_2}\rb_z\big);
\]
taking the limsup over $n_1,m_1$ and the liminf over $n_2,m_2$ gives
\begin{align*}
\exp \big(-\liminf_{n,m\to\infty}\lb x_n|y_m\rb_z\big) &\leq \exp\big(-\limsup_{n,m\to\infty}\lb x_n|y_m \rb_z\big) +\\
&\phantom{XxX}+ \exp\big(-\liminf_{n_1,n_2\to\infty}\lb x_{n_1}|x_{n_2}\rb_z\big) +\\ 
&\phantom{XxX}+ \exp\big(-\liminf_{m_1,m_2\to\infty}\lb y_{m_1}|y_{m_2}\rb_z\big) \noreason\\
&= \exp\big(-\limsup_{n,m\to\infty}\lb x_n|y_m\rb_z\big),\\ 
&\phantom{XxXxXx}\text{(since $(x_n)_1^\infty$ and $(y_m)_1^\infty$ are Gromov)}
\end{align*}
demonstrating equality in \eqref{welldefined1}. The proof of \eqref{welldefined2} is similar and will be omitted.
\end{proof}

\begin{remark}
Many of the statements in this monograph concerning strongly hyperbolic metric spaces are in fact valid for all hyperbolic metric spaces satisfying the conclusion of Lemma \ref{lemmawelldefined}.
\end{remark}

Now that we know that it does not matter too much whether we replace the limit in \eqref{gromovproductdef} by a liminf or a limsup, we make the following definition without fear:

\begin{definition}
For $\xi,\eta\in\del X$ and $y,z\in X$, let
\begin{align} \label{gromovboundarydef1}
\lb \xi|\eta\rb_z &:= \inf\left\{\liminf_{n,m\to\infty}\lb x_n|y_m\rb_z:(x_n)_1^\infty\in\xi,(y_m)_1^\infty\in\eta\right\}\\ \label{gromovboundarydef2}
\lb \xi|y\rb_z &:= \lb y|\xi\rb_z := \inf\left\{\liminf_{n\to\infty}\lb x_n|y\rb_z:(x_n)_1^\infty\in\xi\right\}\\ \label{gromovboundarydef3}
\busemann_\xi(y,z) &= \lb z|\xi\rb_y - \lb y|\xi\rb_z.
\end{align}
\end{definition}

As a corollary of Lemma \ref{lemmawelldefined}, we have the following:

\begin{lemma}
\label{lemmaboundaryasymptotic}
Fix $\xi,\eta\in\del X$ and $y,z\in X$. For all $(x_n)_1^\infty\in\xi$ and $(y_m)_1^\infty\in\eta$ we have 
\begin{align} \label{boundaryasymptotic1}
\lb x_n|y_m\rb_z &\tendsto{n,m,\plus} \lb \xi|\eta\rb_z\\ \label{boundaryasymptotic2}
\lb x_n|y\rb_z &\tendsto{\;\; n,\plus \;\;} \lb \xi|y\rb_z\\ \label{boundaryasymptotic3}
\busemann_{x_n}(y,z) &\tendsto{\;\; n,\plus \;\;} \busemann_\xi(y,z),
\end{align}
(cf. Convention \ref{conventiontendston}), with exact limits if $X$ is strongly hyperbolic.
\end{lemma}
Note that except for the statement about strongly hyperbolic spaces, this lemma is simply \cite[Lemma 5.11]{Vaisala}.
\begin{proof}[Proof of Lemma \ref{lemmaboundaryasymptotic}]
Say we are given $(x_n^{(i)})_1^\infty\in \xi$ and $(y_m^{(i)})_1^\infty\in\eta$ for each $i = 1, 2$.
Let
\[
x_n = \begin{cases}
x_{n/2}^{(1)} & n \text{ even}\\
x_{(n + 1)/2}^{(2)} & n \text{ odd}
\end{cases},
\]
and define $y_m$ similarly. It may be verified using Gromov's inequality that $(x_n)_1^\infty\in \xi$ and $(y_m)_1^\infty\in\eta$. Applying Lemma \ref{lemmawelldefined}, we have
\[
\min_{i = 1}^2 \min_{j = 1}^2 \liminf_{n,m\to\infty} \lb x_n^{(i)}|y_m^{(j)}\rb_z \asymp_\plus \max_{i = 1}^2 \max_{j = 1}^2 \limsup_{n,m\to\infty} \lb x_n^{(i)}|y_m^{(j)}\rb_z .
\]
In particular,
\begin{align*}
\liminf_{n,m\to\infty} \lb x_n^{(1)} | y_n^{(1)} \rb_z
&\lesssim_\plus \limsup_{n,m\to\infty} \lb x_n^{(1)} | y_n^{(1)} \rb_z\\
&\lesssim_\plus \liminf_{n,m\to\infty} \lb x_n^{(2)} | y_n^{(2)} \rb_z\\
&\lesssim_\plus \limsup_{n,m\to\infty} \lb x_n^{(2)} | y_n^{(2)} \rb_z
\lesssim_\plus \liminf_{n,m\to\infty} \lb x_n^{(1)} | y_n^{(1)} \rb_z.
\end{align*}
Taking the infimum over all $(x_n^{(2)})_1^\infty\in\xi$ and $(y_m^{(2)})_1^\infty\in\eta$ gives \eqref{boundaryasymptotic1}. A similar argument gives \eqref{boundaryasymptotic2}. Finally, \eqref{boundaryasymptotic3} follows from \eqref{boundaryasymptotic2}, \eqref{gromovboundarydef3}, and (j) of Proposition \ref{propositionbasicidentities}.

If $X$ is strongly hyperbolic, then all error terms are equal to zero, demonstrating that the limits converge exactly.
\end{proof}
\begin{remark}
In the sequel, the statement that ``if $X$ is strongly hyperbolic, then all error terms are zero'' will typically be omitted from our proofs.
\end{remark}

A simple but useful consequence of Lemma \ref{lemmaboundaryasymptotic} is the following:

\begin{corollary}
\label{corollaryboundaryasymptotic}
The formulas of Proposition \ref{propositionbasicidentities} together with Gromov's inequality hold for points on the boundary as well, if the equations and inequalities there are replaced by additive asymptotics.
 If $X$ is strongly hyperbolic, then we may keep the original formulas without adding an error term.
\end{corollary}
\begin{proof}
For each identity, choose a Gromov sequence representing each element of the boundary which appears in the formula. Replace each occurrence of this element in the formula by the general term of the chosen sequence. This yields a sequence of formulas, each of which is known to be true. Take a subsequence on which each term in these formulas converges. Taking the limit along this subsequence again yields a true formula, and by Lemma \ref{lemmaboundaryasymptotic} we may replace each limit term by the term which it stood for, with only bounded error in doing so, and no error if $X$ is strongly hyperbolic. Thus the formula holds as an additive asymptotic, and holds exactly if $X$ is strongly hyperbolic.
\end{proof}

\begin{remark}
\label{remarkacde}
In fact, (a), (c), (d), and (e) of Proposition \ref{propositionbasicidentities} hold in $\bord X$ in the usual sense, i.e. as exact formulas without additive constants. 
\end{remark}
\begin{proof}
These are the identities where there is at most one Gromov product on each side of the formula. For each element of the boundary, we may simply replace each occurence of that element with the general term of an arbitrary Gromov sequence, take the liminf, and then take the infimum over all Gromov sequences.
\end{proof}

\begin{observation}
\label{observationyzinfinity}
$\lb x|y\rb_z = \infty$ if and only if $x = y\in\del X$.
\end{observation}
\begin{proof}
This follows directly from \eqref{gromovboundarydef1} and \eqref{gromovboundarydef2}.
\end{proof}

\subsection{A topology on $\bord X$}
\label{subsubsectiontopologyclX}
One can endow the bordification $\bord X = X\cup\del X$ with a topological structure $\scrT$ as follows: Given $S\subset\bord X$, write $S\in\scrT$ (i.e. call $S$ open) if
\begin{itemize}
\item[(I)] $S\cap X$ is open, and
\item[(II)] For each $\xi\in S\cap\del X$ there exists $t\geq 0$ such that $N_t(\xi)\subset S$, where
\[
N_t(\xi) := N_{t,\zero}(\xi) := \{y\in \bord X:\lb y|\xi\rb_\zero > t\}
\]
\end{itemize}

\begin{remark}
The topology $\scrT$ may equivalently be defined to be the unique topology on $\bord X$ satisfying:
\begin{itemize}
\item[(I)] $\scrT\given X$ is compatible with the metric $\dist$, and
\item[(II)] For each $\xi\in\del X$, the collection
\begin{equation}
\label{Ntxibase}
\{N_t(\xi) : t \geq 0\}
\end{equation}
is a neighborhood base for $\scrT$ at $\xi$.
\end{itemize}
\end{remark}

\begin{remark}
\label{remarkNtxiopen}
It follows from Lemma \ref{lemmalowersemicontinuous} below that the sets $N_t(\xi)$ are open in the topology $\scrT$.
\end{remark}

\begin{remark}
\label{remarkbasepoint2}
By (d) of Proposition \ref{propositionbasicidentities} (cf. Remark \ref{remarkacde}), we have $N_{t,x}(\xi) \supset N_{t + \dist(x,y),y}(\xi)$ for all $x,y\in X$, $\xi\in\del X$, and $t\geq 0$. Thus the topology $\scrT$ is independent of the basepoint $\zero$.
\end{remark}

The topology $\scrT$ is quite nice. In fact, we have the following:

\begin{proposition}\label{propositionclX}
The topological space $(\bord X,\scrT)$ is completely metrizable. If $X$ is proper and geodesic, then $\bord X$ (and thus also $\del X$) is compact. If $X$ is separable, then $\bord X$ (and thus also $\del X$) is separable.
\end{proposition}
\begin{remark}
If $X$ is proper and geodesic, then Proposition \ref{propositionclX} is \cite[Exercise III.H.3.18(4)]{BridsonHaefliger}.
\end{remark}
\begin{proof}[Proof of Proposition \ref{propositionclX}]
We delay the proof of the complete metrizability of $\bord X$ until Section \ref{subsectionmetametrics}, where we will introduce a class of compatible complete metrics on $\bord X$ which are important from a geometric point of view, the so-called \emph{visual} metrics.

Since $X$ is dense in $\bord X$, the separability of $X$ implies the separability of $\bord X$. Moreover, since $\bord X$ is metrizable (as we will show in Section \ref{subsectionmetametrics}), the separability of $\bord X$ implies the separability of $\del X$.

Finally, assume that $X$ is proper and geodesic; we claim that $\bord X$ is compact. Let $(x_n)_1^\infty$ be a sequence in $X$. If $(x_n)_1^\infty$ contains a bounded subsequence, then since $X$ is proper it contains a convergent subsequence. Thus we assume that $(x_n)_1^\infty$ contains no bounded subsequence, i.e. $\dox{x_n}\to \infty$.

For each $n\in\Namer$ and $t\geq 0$ let
\[
x_{n,t} = \geo\zero{x_n}_{t\wedge\dox{x_n}},\Footnote{Here and from now on $A\wedge B = \min(A,B)$ and $A\vee B = \max(A,B)$.}
\]
where $\geo\zero{x_n}$ is any geodesic connecting $\zero$ and $x_n$. Since $X$ is proper, there exists a sequence $(n_k)_1^\infty$ such that for each $t\geq 0$, the sequence $(x_{n_k,t})_1^\infty$ is convergent, say
\[
x_{n_k,t}\tendsto k x_t.
\]
It is readily verified that the map $t\mapsto x_t$ is an isometric embedding from $\Rplus$ to $X$. Thus there exists a point $\xi\in\del X$ such that $x_t\to \xi$. We claim that $x_{n_k}\to \xi$. Indeed, for each $t\geq 0$,
\[
\limsup_{k\to\infty}\Dist(x_{n_k},x_{n_k,t}) \asymp_\times \limsup_{k\to\infty}\Dist(x_{n_k,t},x_t) \asymp_\times \limsup_{k\to\infty} \Dist(x_t,\xi) \asymp_\times b^{-t},
\]
and so the triangle inequality gives
\[
\limsup_{k\to\infty}\Dist(x_{n_k},\xi) \lesssim_\times b^{-t}.
\]
Letting $t\to\infty$ shows that $x_{n_k}\to\xi$.
\end{proof}

\begin{observation}
\label{observationboundaryconvergence}
A sequence $(x_n)_1^\infty$ in $\bord X$ converges to a point $\xi\in \del X$ if and only if
\begin{equation}\label{boundaryconvergence}
\lb x_n|\xi\rb_\zero \tendsto n \infty.
\end{equation}
\end{observation}

\begin{observation}
\label{observationboundaryconvergence2}
A sequence $(x_n)_1^\infty$ in $X$ converges to a point $\xi\in \del X$ if and only if $(x_n)_1^\infty$ is a Gromov sequence and $(x_n)_1^\infty\in\xi$.
\end{observation}


We now investigate the continuity properties of the Gromov product and Busemann function.

\begin{lemma}[Near-continuity of the Gromov product and Busemann function]
\label{lemmanearcontinuity}
The maps $(x,y,z)\mapsto \lb x|y\rb_z$ and $(x,z,w)\mapsto \busemann_x(z,w)$ are \underline{nearly continuous} in the following sense: Suppose that $(x_n)_1^\infty$ and $(y_n)_1^\infty$ are sequences in $\bord X$ which converge to points $x_n\to x\in\bord X$ and $y_n\to y\in\bord X$. Suppose that $(z_n)_1^\infty$ and $(w_n)_1^\infty$ are sequences in $X$ which converge to points $z_n\to z\in X$ and $w_n\to w\in X$. Then
\begin{align} \label{continuousextension1}
\lb x_n|y_n\rb_{z_n} &\tendsto{n,\plus} \lb x|y\rb_z\\ \label{continuousextension2}
\busemann_{x_n}(z_n,w_n) &\tendsto{n,\plus} \busemann_x(z,w),
\end{align}
with $\tendsto n$ if $X$ is strongly hyperbolic.
\end{lemma}
\begin{proof}
In the proof of \eqref{continuousextension1}, there are three cases:
\begin{itemize}
\item[Case 1:] $x,y\in X$. In this case, \eqref{continuousextension1} follows directly from (d) and (e) of Proposition \ref{propositionbasicidentities}.
\item[Case 2:] $x,y\in\del X$. In this case, for each $n\in\Namer$, choose $\what x_n\in X$ such that either
\begin{itemize}
\item[(1)] $\what x_n = x_n$ (if $x_n\in X$), or
\item[(2)] $\lb \what x_n|x_n\rb_z\geq n$ (if $x_n\in\del X$).
\end{itemize}
Choose $\what y_n$ similarly. Clearly, $\what x_n\to x$ and $\what y_n \to y$. By Observation \ref{observationboundaryconvergence2}, $(\what x_n)_1^\infty\in x$ and $(\what y_n)_1^\infty\in y$. Thus by Lemma \ref{lemmaboundaryasymptotic},
\begin{equation}
\label{whatxwhaty}
\lb \what x_n|\what y_n\rb_z \tendsto{n,\plus} \lb x|y\rb_z.
\end{equation}
Now by Gromov's inequality and (e) of Proposition \ref{propositionbasicidentities}, either
\begin{align*} \tag{1}
\lb \what x_n|\what y_n\rb_z &\asymp_\plus \lb x_n|y_n\rb_{z_n} \text{ or }\\ \tag{2}
\lb \what x_n|\what y_n\rb_z &\gtrsim_\plus n,
\end{align*}
with which asymptotic is true depending on $n$. But for $n$ sufficiently large, \eqref{whatxwhaty} ensures that the (2) fails, so (1) holds.
\item[Case 3:] $x\in X$, $y\in\del X$, or vice-versa. In this case, a straightforward combination of the above arguments demonstrates \eqref{continuousextension1}.
\end{itemize}
Finally, note that \eqref{continuousextension2} is an immediate consequence of \eqref{continuousextension1}, \eqref{gromovboundarydef3}, and (j) of Proposition \ref{propositionbasicidentities}.
\end{proof}

Although Lemma \ref{lemmanearcontinuity} is generally sufficient for applications, we include the following lemma which reassures us that the Gromov product does behave somewhat regularly even on an ``exact'' level.

\begin{lemma}
\label{lemmalowersemicontinuous}
The function $(x,y,z)\mapsto\lb x|y\rb_z$ is lower semicontinuous on $\bord X\times\bord X\times X$.
\end{lemma}
\begin{proof}
Since $\bord X$ is metrizable, it is enough to show that if $x_n\to x$, $y_n\to y$, and $z_n\to z$, then
\[
\liminf_{n\to\infty}\lb x_n|y_n\rb_{z_n} \geq \lb x|y\rb_z.
\]
Now fix $\epsilon > 0$.
\begin{claim}
For each $n\in\N$, there exist points $\what x_n,\what y_n\in X$ satisfying:
\begin{align} \label{LSC1}
\lb \what x_n|\what y_n\rb_{z_n} &\leq \lb x_n|y_n\rb_{z_n} + \epsilon,\\ \label{LSC2}
\lb \what x_n|x_n\rb_\zero &\geq n, \text{ or }\what x_n = x_n\in X,\\ \label{LSC3}
\lb \what y_n|y_n\rb_\zero &\geq n, \text{ or }\what y_n = y_n\in X.
\end{align}
\end{claim}
\begin{subproof}
Suppose first that $x_n,y_n\in\del X$. By the definition of $\lb x_n|y_n\rb_{z_n}$, there exist $(x_{n,k})_1^\infty\in x_n$ and $(y_{n,\ell})_1^\infty\in y_n$ such that
\[
\liminf_{k,\ell\to\infty} \lb x_{n,k}|y_{n,\ell}\rb_{z_n} \leq \lb x_n|y_n\rb_{z_n} + \epsilon/2.
\]
It follows that there exist arbitrarily large\Footnote{Here, of course, ``arbitrarily large'' means that $\min(k,\ell)$ can be made arbitrarily large.} pairs $(k,\ell)\in\N^2$ such that the points $\what x_n := x_{n,k}$ and $\what y_n := y_{n,\ell}$ satisfy \eqref{LSC1}. Since \eqref{LSC2} and \eqref{LSC3} are satisfied for all sufficiently large $(k,\ell)\in\N^2$, this completes the proof. Finally, if either $x_n\in X$, $y_n\in X$, or both, a straightforward adaptation of the above argument yields the claim.
\end{subproof}
The equations \eqref{LSC2} and \eqref{LSC3}, together with Gromov's inequality, imply that $\what x_n\to x$ and $\what y_n\to y$. Now suppose that $x,y\in \del X$. Then by Observation \ref{observationboundaryconvergence2}, $(\what x_n)_1^\infty\in x$ and $(\what y_n)_1^\infty\in y$. So by the definition of $\lb x|y\rb_z$, we have
\begin{align*}
\lb x|y\rb_z &\leq \liminf_{n\to\infty}\lb \what x_n|\what y_n\rb_z \by{the definition of $\lb x|y\rb_z$}\\
&= \liminf_{n\to\infty}\lb \what x_n|\what y_n\rb_{z_n} \by{(d) of Proposition \ref{propositionbasicidentities}}\\
&\leq \liminf_{n\to\infty}\lb x_n|y_n\rb_{z_n} + \epsilon. \by{\eqref{LSC1}}
\end{align*}
Letting $\epsilon$ tend to zero completes the proof. A similar argument applies to the case where $x\in X$, $y\in X$, or both.
\end{proof}

\begin{lemma}
\label{lemmaisometryextension}
If $g$ is an isometry of $X$, then it extends in a unique way to a continuous map $\w g:\bord X\to\bord X$.
\end{lemma}
\begin{proof}
This follows more or less directly from Remarks \ref{remarkbasepoint1} and \ref{remarkbasepoint2}; details are left to the reader.
\end{proof}
In the sequel we will omit the tilde from the extended map $\w g$.

\bigskip
\section{The Gromov product in algebraic hyperbolic spaces} \label{subsectionGromovinH}

In this section we analyze the Gromov product in an algebraic hyperbolic space $X$. We prove Proposition \ref{propositionCATimpliesGromov} which states that CAT(-1) spaces are strongly hyperbolic, and then we show that the Gromov boundary of $X$ is isomorphic to its topological boundary, justifying Remark \ref{remarkambiguityboundary}.

In what follows, we will switch between the hyperboloid model $\H = \H_\F^\alpha$ and the ball model $\B = \B_\F^\alpha$ according to convenience. In the following lemma, $\del\B$ and $\bord\B$ denote the topological boundary and closure of $\B$ as defined in Chapter \ref{sectionROSSONCTs}, not the Gromov boundary and closure as defined above.

\begin{figure}
\begin{center}
\definecolor{qqwuqq}{rgb}{0.0,0.39215686274509803,0.0}
\begin{tikzpicture}[line cap=round,line join=round,>=triangle 45,x=1.0cm,y=1.0cm]
\clip(-3.40,-3.1) rectangle (3.40,3.13);
\draw [shift={(0.0,-0.0)},color=qqwuqq,fill=qqwuqq,fill opacity=0.1] (0,0) -- (27.07967008819517:0.5963691863565421) arc (27.07967008819517:60.832386620422206:0.5963691863565421) -- cycle;
\draw(0.0,0.0) circle (3.0cm);
\draw (0.0,-0.0)-- (2.6711231603651164,1.3656870293596084);
\draw (0.0,-0.0)-- (1.4620984777924115,2.619593106044737);
\draw (1.4620984777924115,2.619593106044737)-- (2.6711231603651164,1.3656870293596084);
\begin{scriptsize}
\draw [fill=black] (1.4620984777924115,2.619593106044737) circle (.75pt);
\draw[color=black] (1.5997838888502904,2.7999980589688112) node {$\xx$};
\draw [fill=black] (2.6711231603651164,1.3656870293596084) circle (.75pt);
\draw[color=black] (2.9117960988346834,1.4482279032273154) node {$\yy$};
\draw [fill=black] (0.0,-0.0) circle (.75pt);
\draw[color=black] (-0.2290816159764389,-0.12160581857100267) node {$\0$};
\draw[color=qqwuqq] (0.6040198381009913,0.57476336778553957) node {$\theta$};
\end{scriptsize}
\end{tikzpicture}
\caption[Relating angle and the Gromov product]{If $\F = \R$ and $\xx,\yy\in\del\B$, then $e^{-\lb\xx|\yy\rb_\0} = \frac12\|\yy - \xx\| = \sin(\theta/2)$, where $\theta$ denotes the angle $\ang_\0(\xx,\yy)$ drawn in the figure.}
\end{center}
\end{figure}

\begin{lemma}
\label{lemmagromovextension}
The Gromov product $(\xx,\yy,\zz)\mapsto \lb \xx|\yy\rb_\zz:\B\times\B\times\B\to\Rplus$ extends uniquely to a continuous function $(\xx,\yy,\zz)\mapsto \lb \xx|\yy\rb_\zz:\bord\B\times\bord\B\times\B\to[0,\infty]$. Moreover, the extension satisfies the following:
\begin{itemize}
\item[(i)] $\lb\xx|\yy\rb_{\zz} = \infty$ if and only if $\xx = \yy\in\del\B$.
\item[(ii)] For all $\xx,\yy\in\bord\B$,
\begin{equation}
\label{euclideanball1}
e^{-\lb\xx|\yy\rb_\0} \geq \frac{1}{\sqrt 8}\|\yy - \xx\|.
\end{equation}
If $\F = \R$ and $\xx,\yy\in\del\B$, then
\begin{equation}
\label{euclideanball2}
e^{-\lb\xx|\yy\rb_\0} = \frac12\|\yy - \xx\|.
\end{equation}
\end{itemize}
\end{lemma}
\begin{proof}
We begin by making some computations in the hyperboloid model $\H$. For $[\xx],[\yy]\in\bord\H$ and $[\zz]\in\H$, let
\[
\alpha_{[\zz]}([\xx],[\yy]) = \frac{|\QQ(\zz)|\cdot |B_\QQ(\xx,\yy)|}{|B_\QQ(\xx,\zz)|\cdot |B_\QQ(\yy,\zz)|}\in\Rplus.
\]
By \eqref{distanceinL}, for $[\xx],[\yy],[\zz]\in\H$ we have
\begin{equation}
\label{lorentzlawofcosh}
\alpha_{[\zz]}([\xx],[\yy]) = \frac{\cosh\dist_\H([\xx],[\yy])}{\cosh\dist_\H([\xx],[\zz])\cosh\dist_\H([\yy],[\zz])}\cdot
\end{equation}
Let $\DD = \{(A,B,C)\in\Rplus^3: \cosh(A)\cosh(B)C \geq 1\}$, and define $F: \DD\to\Rplus$ by
\[
F(A,B,C) = \frac{\exp\left[\cosh^{-1}\left(\cosh(A)\cosh(B)C\right)\right]}{e^A e^B}\cdot
\]
Then by \eqref{lorentzlawofcosh}, we have
\[
e^{-2\lb [\xx]|[\yy]\rb_{[\zz]}} = F\left(\dist_\H([\zz],[\xx]),\dist_\H([\zz],[\yy]),\alpha_{[\zz]}([\xx],[\yy])\right)
\]
for all $[\xx],[\yy],[\zz]\in\H$. Now since $\lim_{t\to\infty}\frac{e^t}{\cosh(t)} = 2$, we have for all $A\geq 0$ and $C > 0$
\begin{align*}
\lim_{B\to\infty}F(A,B,C)
&= \lim_{B\to\infty}\frac{2\left(\cosh(A)\cosh(B)C\right)}{e^A e^B}\\
&= \frac{\cosh(A)}{e^A} C
\end{align*}
and 
\[
\lim_{A\to\infty}\frac{\cosh(A)}{e^A} C
= C/2.
\]
Let $\what \DD$ be the closure of $\DD$ relative to $[0,\infty]^2\times \Rplus$, i.e. 
\[
\what \DD = \DD\cup \left([0,\infty]^2\times \Rplus\butnot \Rplus^3\right).
\] 
If we let
\[
\what F(A,B,C) :=
\begin{cases}
F(A,B,C) & A,B < \infty\\
\frac{\cosh(A)}{e^A} C & A < B = \infty\\
\frac{\cosh(B)}{e^B} C & B < A = \infty\\
C/2 & A = B = \infty
\end{cases}
\]
then $\what F: \what \DD\to\Rplus$ is a continuous function.\Footnote{Technically, the calculations above do not prove the continuity of $\what F$; however, this continuity is easily verified using standard methods.} Thus, letting
\[
\lb [\xx]|[\yy]\rb_{[\zz]} := -\frac12\log\what F\left(\dist_\H([\zz],[\xx]),\dist_\H([\zz],[\yy]),\alpha_{[\zz]}([\xx],[\yy])\right)
\]
defines a continuous extension of the Gromov product to $\bord\H\times\bord\H\times\H$.

We now prove (i)-(ii):
\begin{itemize}
\item[(i)] Using the inequality $e^t/2 \leq \cosh(t) \leq e^t$, it is easily verified that
\begin{equation}
\label{C4}
\what F(A,B,C) \geq C/4
\end{equation}
for all $(A,B,C)\in\what\DD$. In particular, if $\what F(A,B,C) = 0$ then $C = 0$. Thus if $\lb [\xx]|[\yy]\rb_{[\zz]} = \infty$ then $\alpha_{[\zz]}([\xx],[\yy]) = 0$; since $[\zz]\in\H$ we have $B_\QQ(\xx,\yy) = 0$, and by Observation \ref{observationBQnonzero} we have $[\xx] = [\yy]\in\del\H$. Conversely, if $[\xx] = [\yy]\in\del\H$, then $\dist_\H([\zz],[\xx]) = \dist_\H([\zz],[\yy]) = \infty$ and $\alpha_{[\zz]}([\xx],[\yy]) = 0$, so $\lb [\xx]|[\yy]\rb_{[\zz]} = -\frac12 \log\what F(\infty,\infty,0) = \infty$.

\item[(ii)] Recall that in $\B$, $\zero = [(1,\0)]$. For $\xx,\yy\in\B$,
\begin{align*}
\alpha_\zero(e_{\B,\H}(\xx),e_{\B,\H}(\yy)) &= \frac{|\QQ(1,\0)|\cdot|B_\QQ((1,\xx),(1,\yy))|}{|B_\QQ((1,\0),(1,\xx))|\cdot |B_\QQ((1,\0),(1,\yy))|} \noreason\\
&= |1 - B_\EE(\xx,\yy)|\\
&\geq 1 - \Re B_\EE(\xx,\yy) &\text{(with equality if $\F = \R$)}\\
&\geq \frac12 [\|\xx\|^2 + \|\yy\|^2] - \Re B_\EE(\xx,\yy) &\text{(with equality if $\xx,\yy\in\del\B$)}\\
&= \frac12\|\yy - \xx\|^2.
\end{align*}
Combining with \eqref{C4} gives
\[
e^{-2\lb \xx|\yy\rb_\0} \geq \frac14\alpha_\zero(e_{\B,\H}(\xx),e_{\B,\H}(\yy)) \geq \frac18\|\yy - \xx\|^2.
\]
If $\F = \R$ and $\xx,\yy\in\del\B$, then
\[
e^{-2\lb \xx|\yy\rb_\0} = \frac12\alpha_\zero(e_{\B,\H}(\xx),e_{\B,\H}(\yy)) = \frac14\|\yy - \xx\|^2.
\]
\end{itemize}
\end{proof}

We now prove Proposition \ref{propositionCATimpliesGromov}, beginning with the following lemma:
\begin{lemma}
\label{lemmaHstrongly}
If $\F = \R$ then $\B$ is strongly Gromov hyperbolic.
\end{lemma}
\begin{proof}
By the transitivity of the isometry group (Observation \ref{observationtransitivity}), it suffices to check \eqref{gromovforCAT} for the special case $w = \zero$. So let us fix $x,y,z\in\B$, and by contradiction suppose that
\[
e^{-\lb x|z\rb_\zero} > e^{-\lb x|y\rb_\zero} + e^{-\lb y|z\rb_\zero},
\]
or equivalently that
\[
1 > e^{\lb x|z\rb_\zero - \lb x|y\rb_\zero} + e^{\lb x|z\rb_\zero - \lb y|z\rb_\zero}.
\]
Clearly, the above inequality implies that $x\neq z$ and $y\neq\zero$. Now let $\gamma_1$ and $\gamma_2$ be the unique bi-infinite geodesics extending the geodesic segments $\geo xz$ and $\geo\zero y$, respectively. Let $x_\infty,z_\infty\in\del\B$ be the appropriate endpoints of $\gamma_1$, and let $y_\infty$ be the endpoint of $\gamma_2$ which is closer to $y$ than to $\zero$. (See Figure \ref{figureHstrongly}.) For each $t\in\Rplus$, let
\[
x_t = \geo x{x_\infty}_t\in\gamma_1,
\]
and let $y_t\in\gamma_2$, $z_t\in\gamma_1$ be defined similarly.

\begin{figure}
\begin{center}
\begin{tikzpicture}[line cap=round,line join=round,>=triangle 45,x=0.4402812927177227cm,y=0.44028129271787425cm]
\clip(-9.71,-7.2) rectangle (11,7.5);
\draw(0,0) ellipse (3.11cm and 3.11cm);
\draw (5,5)-- (-5,-5);
\draw [shift={(12.49,9.91)}] plot[domain=3.35:4.27,variable=\t]({1*14.67*cos(\t r)+0*14.67*sin(\t r)},{0*14.67*cos(\t r)+1*14.67*sin(\t r)});
\begin{scriptsize}
\fill [color=black] (0,0) circle (1pt);
\draw[color=black] (0.45,-0.58) node {$\zero$};
\fill [color=black] (5,5) circle (1pt);
\draw[color=black] (5.9,4.96) node {$y_\infty$};
\fill [color=black] (-1.85,6.83) circle (1pt);
\draw[color=black] (-2.5,7.14) node {$x_\infty$};
\fill [color=black] (6.29,-3.33) circle (1pt);
\draw[color=black] (7.2,-3.6) node {$z_\infty$};
\fill [color=black] (-0.97,4.09) circle (1pt);
\draw[color=black] (-1.53,3.88) node {$x$};
\fill [color=black] (-1.36,5.08) circle (1pt);
\draw[color=black] (-1.84,4.81) node {$x_t$};
\fill [color=black] (4.34,-2.28) circle (1pt);
\draw[color=black] (3.99,-2.7) node {$z$};
\fill [color=black] (4.91,-2.64) circle (1pt);
\draw[color=black] (4.93,-3.3) node {$z_t$};
\fill [color=black] (3.12,3.12) circle (1pt);
\draw[color=black] (3.67,2.63) node {$y_t$};
\fill [color=black] (2,2) circle (1pt);
\draw[color=black] (2.42,1.7) node {$y$};
\end{scriptsize}
\end{tikzpicture}
\caption[$\B$ is strongly Gromov hyperbolic]{If Gromov's inequality fails for the quadruple $x,y,z,\zero$, then it also fails for the quadruple $x_\infty,y_\infty,z_\infty,\zero$.}
\label{figureHstrongly}
\end{center}
\end{figure}
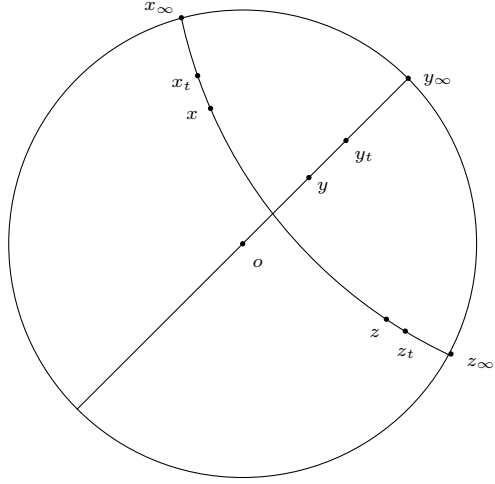

We observe that
\begin{align*}
\frac{\del}{\del t}\left[\lb x_t|z_t\rb_\zero - \lb x_t|y_t\rb_\zero\right]
&= \frac12 \frac{\del}{\del t}\bigl[\dist_\B(\zero,z_t) + \dist_\B(x_t,y_t) - \dist_\B(x_t,z_t) - \dist_\B(\zero,y_t)\bigr]\\
&= \frac12 \frac{\del}{\del t}\bigl[\dist_\B(\zero,z_t) + \dist_\B(x_t,y_t) - 2t - t\bigr]\\
&\leq \frac12 \frac{\del}{\del t}\bigl[t + 2t - 2t - t\bigr] = 0,
\end{align*}
i.e. the expression $\lb x_t|z_t\rb_\zero - \lb x_t|y_t\rb_\zero$ is nonincreasing with respect to $t$. Taking the limit as $t$ approaches infinity, we have
\[
\lb x_\infty|z_\infty\rb_\zero - \lb x_\infty|y_\infty\rb_\zero \leq \lb x|z\rb_\zero - \lb x|y\rb_\zero
\]
and a similar argument shows that
\[
\lb x_\infty|z_\infty\rb_\zero - \lb y_\infty|z_\infty\rb_\zero \leq \lb x|z\rb_\zero - \lb y|z\rb_\zero.
\]
Thus
\[
1 > e^{\lb x_\infty|z_\infty\rb_\zero - \lb x_\infty|y_\infty\rb_\zero} + e^{\lb x_\infty|z_\infty\rb_\zero - \lb y_\infty|z_\infty\rb_\zero}
\]
or equivalently,
\[
e^{-\lb x_\infty|z_\infty\rb_\zero} > e^{-\lb x_\infty|y_\infty\rb_\zero} + e^{-\lb y_\infty|z_\infty\rb_\zero}.
\]
But by \eqref{euclideanball2}, if we write $x_\infty = \xx$, $y_\infty = \yy$, and $z_\infty = \zz$, then
\[
\frac12\|\zz - \xx\| > \frac12\|\yy - \xx\| + \frac12\|\zz - \yy\|.
\]
This is a contradiction.
\end{proof}

We are now ready to prove
\begin{repproposition}{propositionCATimpliesGromov}
Every CAT(-1) space is strongly hyperbolic.
\end{repproposition}
\begin{proof}
Let $X$ be a CAT(-1) space, and fix $x,y,z,w\in X$. By \cite[Proposition II.1.11]{BridsonHaefliger}, there exist $\wbar x,\wbar y,\wbar z,\wbar w\in \H^2$ such that
\begin{align*}
\dist(x,y) &= \dist(\wbar x,\wbar y) &
\dist(y,z) &= \dist(\wbar y,\wbar z)\\
\dist(z,w) &= \dist(\wbar z,\wbar w) &
\dist(w,x) &= \dist(\wbar w,\wbar x)\\
\dist(x,z) &\leq \dist(\wbar x,\wbar z) &
\dist(y,w) &\leq \dist(\wbar y,\wbar w).
\end{align*}
It follows that
\[
e^{-\lb x|z\rb_w} \leq e^{-\lb \wbar x|\wbar z\rb_{\wbar w}} \leq e^{-\lb \wbar x|\wbar y\rb_{\wbar w}} + e^{-\lb \wbar y|\wbar z\rb_{\wbar w}} \leq e^{-\lb x|y\rb_w} + e^{-\lb y|z\rb_w}.
\]
\end{proof}

\bigskip
\subsection{The Gromov boundary of an algebraic hyperbolic space}
\label{subsubsectionboundariesequivalent}

Again let $X = \H = \H_\F^\alpha$ be an algebraic hyperbolic space. By Proposition \ref{propositionCATimpliesGromov}, $X$ is a (strongly) hyperbolic metric space. (If $\F = \R$, we can use Lemma \ref{lemmaHstrongly}.) In particular, $X$ has a Gromov boundary, defined in Section \ref{subsectionboundary}. On the other hand, $X$ also has a topological boundary, defined in Chapter \ref{sectionROSSONCTs}. For this subsection only, we will write
\begin{align*}
\del_G X &= \text{ Gromov boundary of $X$}, \\
\del_T X &= \text{ topological boundary of $X$}.
\end{align*}
We will now show that this distinction is in fact unnecessary.
\begin{proposition}
\label{propositionboundariesequivalent}
The identity map $\id:X\to X$ extends uniquely to a homeomorphism $\what\id:X\cup\del_G X\to X\cup\del_T X$. Thus the pairs $(X,X\cup\del_G X)$ and $(X,X\cup\del_T X)$ are isomorphic in the sense of Section \ref{subsectiontotallygeodesic}.
\end{proposition}

\begin{proof}[Proof of Proposition \ref{propositionboundariesequivalent}]
By Observation \ref{observationHequivB} and Proposition \ref{propositionHequivE}, it suffices to consider the case where $X = \B = \B_\F^\alpha$ is the ball model. Fix $\xi\in\del_G\B$. By definition, $\xi = [(\xx_n)_1^\infty]$ for some Gromov sequence $(\xx_n)_1^\infty$. By \eqref{euclideanball1}, the sequence $(\xx_n)_1^\infty$ is Cauchy in the metric $\|\cdot - \cdot\|$. Thus $\xx_n\to\xx$ for some $\xx\in\bord\B$; since $(\xx_n)_1^\infty$ is a Gromov sequence, we have
\[
\lb \xx|\xx\rb_\0 = \lim_{n,m\to\infty}\lb \xx_n|\xx_m\rb_\0 = \infty,
\]
and thus $\xx\in\del_T\B$ by (i) of Lemma \ref{lemmagromovextension}. Let
\[
\what\id(\xi) = \xx.
\]
To see that the map $\what\id$ is well-defined, note that if $(\yy_n)_1^\infty$ is another Gromov sequence equivalent to $(\xx_n)_1^\infty$, and if $\yy_n\to\yy\in\del_T \B$, then
\[
\lb \xx|\yy\rb_\0 = \lim_{n\to\infty}\lb \xx_n|\yy_n\rb_\0 = \infty,
\]
and so by (i) of Lemma \ref{lemmagromovextension} we have $\xx = \yy$.

We next claim that $\what\id:\del_G\B\to\del_T\B$ is a bijection. To demonstrate injectivity, we note that if $\what\id(\xi) = \what\id(\eta) = \xx$, then by (i) of Lemma \ref{lemmagromovextension}
\[
\lim_{n\to\infty}\lb \xx_n|\yy_n\rb_\0 = \lb \xx|\xx\rb_\0 = \infty,
\]
where $(\xx_n)_1^\infty$ and $(\yy_n)_1^\infty$ are Gromov sequences representing $\xi$ and $\eta$, respectively. Thus $(\xx_n)_1^\infty$ and $(\yy_n)_1^\infty$ are equivalent, and so $\xi = \eta$.

To demonstrate surjectivity, we observe that for $\xx\in\del_T\B$, we have
\[
\what\id\left(\left[\left(\frac{n - 1}{n}\xx\right)_1^\infty\right]\right) = \xx.
\]
Finally, we must demonstrate that $\what\id$ is a homeomorphism, or in other words that the topology defined in \6\ref{subsubsectiontopologyclX} is the usual topology on $\bord\B$ (i.e. the topology inherited from $\HH$). It suffices to demonstrate the following:
\begin{claim}
For any $\xx\in\del_T\B$, the collection \eqref{Ntxibase} (with $\xi = \xx$) is a neighborhood base of $\xx$ with respect to the usual topology.
\end{claim}
\begin{subproof}
By \eqref{euclideanball1}, we have
\[
N_t(\xx) \subset B(\xx,\sqrt 8 e^{-t}).
\]
On the other hand, the continuity of the Gromov product on $\bord\B$ guarantees that $N_t(\xx)$ contains a neighborhood of $\xx$ with respect to the usual topology.
\end{subproof}

\end{proof}

In the sequel, the following will be useful:

\begin{figure}
\begin{center}
\begin{tikzpicture}[line cap=round,line join=round,>=triangle 45,x=1.0cm,y=1.0cm]
\clip(-3.8,-0.0) rectangle (5.57,4.41);
\draw (-3.0,0.0)-- (3.0,0.0);
\draw [dash pattern=on 5pt off 5pt] (-1.52,2.22)-- (-1.52,0.0);
\draw [dash pattern=on 5pt off 5pt] (1.76,1.48)-- (1.7800000000000002,0.0);
\begin{scriptsize}
\draw [fill=black] (0.0,4.0) circle (.75pt);
\draw[color=black] (0.13316978586974992,4.270045192923707) node {$\infty$};
\draw [fill=black] (-1.52,2.22) circle (.75pt);
\draw[color=black] (-1.7346933526424422,2.5810200144818323) node {$\xx$};
\draw [fill=black] (1.76,1.48) circle (.75pt);
\draw[color=black] (1.961291155477427,1.7663137519392818) node {$\yy$};
\draw[color=black] (-1.834047774903729,1.1304454494670468) node {$x_1$};
\draw[color=black] (2.0209038088341993,0.8721239515877014) node {$y_1$};
\end{scriptsize}
\end{tikzpicture}
\caption[A formula for the Busemann function in the half-space model]{The value of the Busemann function $\busemann_\infty(\xx,\yy)$ depends on the heights of the points $\xx$ and $\yy$.}
\end{center}
\end{figure}

\begin{proposition}
\label{propositionbusemannE}
Let $\E = \E^\alpha$ be the half-space model of a real hyperbolic space. For $\xx,\yy\in\E$, we have
\[
\busemann_\infty(\xx,\yy) = -\log(x_1/y_1).
\]
\end{proposition}
\begin{proof}
By \eqref{distE} we have
\begin{align*}
e^{\busemann_\infty(\xx,\yy)}
= \lim_{\zz\to\infty} \frac{\exp\dist_\E(\zz,\xx)}{\exp\dist_\E(\zz,\yy)}
&= \lim_{\zz\to\infty} \frac{\cosh\dist_\E(\zz,\xx)}{\cosh\dist_\E(\zz,\yy)}\\
&= \lim_{\zz\to\infty} \frac{1 + \frac{\|\zz - \xx\|^2}{2 x_1 z_1}}{1 + \frac{\|\zz - \yy\|^2}{2 y_1 z_1}}
= \lim_{\zz\to\infty} \frac{\left(\frac{\|\zz - \xx\|^2}{2 x_1 z_1}\right)}{\left(\frac{\|\zz - \yy\|^2}{2 y_1 z_1}\right)}
= \frac{y_1}{x_1}\cdot
\end{align*}
\end{proof}

\bigskip
\section{Metrics and metametrics on $\bord X$}
\label{subsectionmetametrics}

\subsection{General theory of metametrics}

\begin{definition}
\label{definitionmetametric}
Recall that a \emph{metric} on a set $Z$ is a map $\Dist:Z\times Z\to\Rplus$ which satisfies:
\begin{itemize}
\item[(I)] Reflexivity: $\Dist(x,x) = 0$.
\item[(II)] Reverse reflexivity: $\Dist(x,y) = 0$ \implies $x = y$.
\item[(III)] Symmetry: $\Dist(x,y) = \Dist(y,x)$.
\item[(IV)] Triangle inequality: $\Dist(x,z) \leq \Dist(x,y) + \Dist(y,z)$.
\end{itemize}
Now we can define a \emph{metametric} on $Z$ to be a map $\Dist: Z\times Z\to\Rplus$ which satisfies (II), (III), and (IV), but not necessarily (I). This concept is not to be confused with the more common notion of a \emph{pseudometric}, which satisfies (I), (III), and (IV), but not necessarily (II). The term ``metametric'' was introduced by J. V\"ais\"al\"a in \cite{Vaisala}.

If $\Dist$ is a metametric, we define its \emph{domain of reflexivity} to be the set $Z_\refl := \{x\in Z:\Dist(x,x) = 0\}$.\Footnote{In the terminology of \cite[p.19]{Vaisala}, the domain of reflexivity is the set of ``small points''.} Obviously, $\Dist$ restricted to its domain of reflexivity is a metric.
\end{definition}

As in metric spaces, a sequence $(x_n)_1^\infty$ in a metametric space $(Z,\Dist)$ is called \emph{Cauchy} if $\Dist(x_n,x_m)\tendsto{n,m} 0$, and \emph{convergent} if there exists $x\in Z$ such that $\Dist(x_n,x)\to 0$. (However, see Remark \ref{remarktopologymetametrics} below.) The metametric space $(Z,\Dist)$ is called \emph{complete} if every Cauchy sequence is convergent. Using these definitions, the standard proof of the Banach contraction principle immediately yields the following:

\begin{theorem}[Banach contraction principle for metametric spaces]\label{banachcontractionprinciple}
Let $(Z,\Dist)$ be a complete metametric space. Fix $0 < \lambda < 1$. If $g:Z\to Z$ satisfies
\[
\Dist(g(z),g(w)) \leq \lambda\Dist(z,w) \all z,w\in Z,
\]
then there exists a unique point $z\in Z$ so that $g(z) = z$. Moreover, for all $w\in Z$, we have $g^n(w)\to z$ with respect to the metametric $\Dist$.
\end{theorem}

\begin{observation}
The fixed point coming $z$ coming from Theorem \ref{banachcontractionprinciple} must lie in the domain of reflexivity $Z_\refl$.
\end{observation}
\begin{proof}
\[
\Dist(z,z) = \Dist(g(z),g(z)) \leq \lambda \Dist(z,z),
\]
and thus $\Dist(z,z) = 0$.
\end{proof}

Recall that a metric is said to be \emph{compatible} with a topology if that topology is equal to the topology induced by the metric. We now generalize this concept by introducing the notion of compatibility between a topology and a metametric.

\begin{definition}
\label{definitioncompatiblemetametric}
Let $(Z,\Dist)$ be a metametric space. A topology $\scrT$ on $Z$ is \emph{compatible} with the metametric $\Dist$ if for every $\xi\in Z_\refl$, the collection
\begin{equation}
\label{balls}
\left\{B_\Dist(\xi,r) := \{y\in Z : \Dist(\xi,y) < r\} : r > 0\right\}
\end{equation}
forms a neighborhood base for $\scrT$ at $\xi$.
\end{definition}

Note that unlike a metric, a metametric may have multiple topologies with which it is compatible.\Footnote{The topology considered in \cite[p.19]{Vaisala} is the finest topology compatible with a given metametric.} The metametric is viewed as determining a neighborhood base for points in the domain of reflexivity; neighborhood bases for other points must arise from some other structure. In the case we are interested in, namely the case where the underlying space for the metametric is the Gromov closure of a hyperbolic metric space $X$, the topology on the complement of the domain of reflexivity will come from the original metric $\dist$ on $X$.

\begin{remark}
\label{remarktopologymetametrics}
If $(Z,\Dist)$ is a metametric space with a compatible topology $\scrT$, then there are two notions of what it means for a sequence $(x_n)_1^\infty$ in $Z$ to converge to a point $x\in Z$: the sequence may converge with respect to the topology $\scrT$, or it may converge with respect to the metametric (i.e. $\Dist(x_n,x)\to 0$). The relation between these notions is as follows: $x_n\to x$ with respect to the metametric $\Dist$ if and only if both of the following hold: $x_n\to x$ with respect to the topology $\scrT$, and $x\in Z_\refl$.
\end{remark}

\begin{remark}
If a metametric $\Dist$ on a set $Z$ is compatible with a topology $\scrT$, then the metric $\Dist\given Z_\refl$ is compatible with the topology $\scrT\given Z_\refl$. However, the converse does not necessarily hold.
\end{remark}

For the remainder of this chapter, we fix a hyperbolic metric space $X$, and we let $\scrT$ be the topology on $\bord X$ introduced in \6\ref{subsubsectiontopologyclX}. We will consider various metrics and metametrics on $\bord X$ which are compatible with the topology $\scrT$.

\subsection{The visual metametric based at a point $\notzero\in X$}
The first metametric that we will consider is designed to emulate the Euclidean or ``spherical'' metric on the boundary of the ball model $\B$. Recall from Lemma \ref{lemmagromovextension} that
\begin{repequation}{euclideanball2}
\frac12\|\yy - \xx\| = e^{-\lb \xx|\yy\rb_\0} \text{ for all $\xx,\yy\in\del\B$.}
\end{repequation}
The metric $(\xx,\yy)\mapsto \frac12\|\yy - \xx\|$ can be thought of as ``seen from $\0$''. The expression on the right hand side makes sense if $\xx,\yy\in\bord\B$, and defines a metametric on $\bord\B$. Moreover, the formula can be generalized to an arbitrary strongly hyperbolic metric space:
\begin{observation}
\label{observationDist}
If $X$ is a strongly hyperbolic metric space, then for each $\notzero\in X$ the map $\Dist_\notzero:\bord X\times\bord X\to \Rplus$ defined by
\begin{equation}
\label{distancedef}
\Dist_\notzero(x,y) := e^{-\lb x|y\rb_\notzero}
\end{equation}
is a complete metametric on $\bord X$. This metametric is compatible with the topology $\scrT$; moreover, its domain of reflexivity is $\del X$.
\end{observation}
\begin{proof}
Reverse reflexivity and the fact that $(\bord X)_\refl = \del X$ follow directly from Observation \ref{observationyzinfinity}; symmetry follows from (a) of Proposition \ref{propositionbasicidentities} together with Corollary \ref{corollaryboundaryasymptotic}; the triangle inequality follows from the definition of strong hyperbolicity together with Corollary \ref{corollaryboundaryasymptotic}.

To show that $\Dist_\notzero$ is complete, suppose that $(x_n)_1^\infty$ is a Cauchy sequence in $X$. Applying \eqref{distancedef}, we see that $\lb x_n | x_m\rb_\notzero \tendsto{n,m} \infty$, i.e. $(x_n)_1^\infty$ is a Gromov sequence. Letting $\xi = [(x_n)_1^\infty]$, we have $x_n\to \xi$ in the $\Dist_{b,\notzero}$ metametric. Thus every Cauchy sequence in $X$ converges in $\bord X$. Since $X$ is dense in $\bord X$, a standard approximation argument shows that $\bord X$ is complete.

Given $\xi\in (\bord X)_\refl = \del X$, the collection \eqref{balls} is equal to the collection \eqref{Ntxibase}, and is therefore a neighborhood base for $\scrT$ at $\xi$. Thus $\Dist_\notzero$ is compatible with $\scrT$.
\end{proof}

Next, we drop the assumption that $X$ is strongly hyperbolic. Fix $b > 1$ and $\notzero\in X$, and consider the function
\begin{equation}
\label{Distdef}
\Dist_{b,\notzero}(x,y) = \inf_{(x_i)_0^n} \sum_{i = 0}^{n - 1} b^{-\lb x_i | x_{i + 1}\rb_\notzero},
\end{equation}
where the infimum is taken over finite sequences $(x_i)_0^n$ satisfying $x_0 = x$ and $x_n = y$.

\begin{proposition}
\label{propositionDist}
If $b > 1$ is sufficiently close to $1$, then for each $\notzero\in X$, the function $\Dist_{b,\notzero}$ defined by \eqref{Distdef} is a complete metametric on $\bord X$ satisfying the following inequality:
\begin{equation}
\label{distanceasymptotic}
b^{-\lb x|y\rb_\notzero} / 4 \leq \Dist_{b,\notzero}(x,y) \leq b^{-\lb x|y\rb_\notzero}.
\end{equation}
This metametric is compatible with the topology $\scrT$; moreover, its domain of reflexivity is $\del X$.
\end{proposition}
We will refer to $\Dist_{b,\notzero}$ as the ``visual (meta)metric from the point $\notzero$ with respect to the parameter $b$''.

\begin{remark}
The metric $\Dist_{b,\notzero}\given\del X$ has been referred to in the literature as the \emph{Bourdon metric}.
\end{remark}

\begin{remark}
The first part of Proposition \ref{propositionDist} is \cite[Propositions 5.16 and 5.31]{Vaisala}.
\end{remark}

\begin{proof}[Proof of Proposition \ref{propositionDist}]
Let $\delta \geq 0$ be the implied constant in Gromov's inequality, and fix $1 < b \leq 2^{1/\delta}$. Then raising $b^{-1}$ to the power of both sides of Gromov's inequality gives
\[
b^{-\lb x|z\rb_\notzero} \leq 2\max\left(b^{-\lb x|y\rb_w}, b^{-\lb y|z\rb_w}\right),
\]
i.e. the function
\[
(x,y)\mapsto b^{-\lb x|y\rb_\notzero}
\]
satisfies the ``weak triangle inequality'' of \cite{Schroeder}. A straightforward adaptation of the proof of \cite[Theorem 1.2]{Schroeder} demonstrates \eqref{distanceasymptotic}. Condition (II) of being a metametric and the equality $(\bord X)_\refl = \del X$ now follow from Observation \ref{observationyzinfinity}. Conditions (III) and (IV) of being a metametric are immediate from \eqref{Distdef}.

The argument for completeness is the same as in the proof of Observation \ref{observationDist}.

Finally, given $\xi\in (\bord X)_\refl = \del X$, we observe that although the collections \eqref{balls} and \eqref{Ntxibase} are no longer equal, \eqref{distanceasymptotic} guarantees that the filters they generate are equal, which is enough to show that $\Dist_{b,\notzero}$ is compatible with $\scrT$.
\end{proof}

\begin{remark}
\label{remarkDist}
If $X$ is strongly hyperbolic, then Proposition \ref{propositionDist} holds for all $1 < b\leq e$; moreover, the metametric $\Dist_{e,\notzero}$ is equal to the metametric $\Dist_\notzero$ defined in Observation \ref{observationDist}.
\end{remark}

\begin{remark}
\label{remarkDistRtree}
If $(X,\dist)$ is an $\R$-tree, then for all $t > 0$, $(X,t\dist)$ is also an $\R$-tree and is therefore strongly hyperbolic (by Observation \ref{observationRtreesCAT} and Proposition \ref{propositionCATimpliesGromov}). It follows that Proposition \ref{propositionDist} holds for all $b > 1$.
\end{remark}

For the remainder of this chapter, we fix $b > 1$ close enough to $1$ so that Proposition \ref{propositionDist} holds.


\subsection{The extended visual metric on $\bord X$}
\label{subsectionvisual}
Although the metametric $\Dist_{b,\notzero}$ has the advantage of being directly linked to the Gromov product via \eqref{distanceasymptotic}, it is sometimes desirable to put a metric on $\bord X$, not just a metametric. We show now that such a metric can be constructed which agrees with $\Dist_{b,\notzero}$ on $\del X$.

In the following proposition, we use the convention that $\dist(x,y) = \infty$ if $x,y\in\bord X$ and either $x\in\del X$ or $y\in\del X$.

\begin{proposition}
\label{propositionwbarDist}
Fix $\notzero\in X$, and for all $x,y\in \bord X$ let
\[
\wbar\Dist_{b,\notzero}(x,y) = \min\big(\log(b)\dist(x,y), \Dist_{b,\notzero}(x,y)\big).
\]
Then $\wbar\Dist = \wbar\Dist_{b,\notzero}$ is a complete metric on $\bord X$ which agrees with $\Dist = \Dist_{b,\notzero}$ on $\del X$ and induces the topology $\scrT$.
\end{proposition}
We call the metric $\wbar\Dist$ an \emph{extended visual metric}.

As an immediate consequence we have the following result which was promised in \6\ref{subsubsectiontopologyclX}:
\begin{corollary}
The topological space $(\bord X,\scrT)$ is completely metrizable.
\end{corollary}
\begin{proof}[Proof of Proposition \ref{propositionwbarDist}]
Let us show that $\wbar\Dist$ is a metric. Conditions (I)-(III) are obvious. To demonstrate the triangle inequality, fix $x,y,z\in\bord X$.
\begin{itemize}
\item[(1)] If $\wbar\Dist(x,y) = \log(b)\dist(x,y)$ and $\wbar\Dist(y,z) = \log(b)\dist(y,z)$, then $\wbar\Dist(x,z) \leq \log(b)\dist(x,z) \leq \log(b)\dist(x,y) + \log(b)\dist(y,z) = \wbar\Dist(x,y) + \wbar\Dist(y,z)$. Similarly, if $\wbar\Dist(x,y) = \Dist(x,y)$ and $\wbar\Dist(y,z) = \Dist(y,z)$, then $\wbar\Dist(x,z) \leq \Dist(x,z) \leq \Dist(x,y) + \Dist(y,z) = \wbar\Dist(x,y) + \wbar\Dist(y,z)$.
\item[(2)] If $\wbar\Dist(x,y) = \log(b)\dist(x,y)$ and $\wbar\Dist(y,z) = \Dist(y,z)$, fix $\epsilon > 0$, and let $y = y_0,y_1,\ldots,y_n = z$ be a sequence such that
\[
\sum_{i = 0}^{n - 1}b^{-\lb y_i,y_{i + 1}\rb_\notzero} \leq \Dist(y,z) + \epsilon.
\]
Let $x_i = y_i$ for $i\geq 1$ but let $x_0 = x$. Then by (e) of Proposition \ref{propositionbasicidentities} and the inequality
\[
b^{-t} \leq s\log(b) + b^{-(t + s)} \hspace{.5 in} (s,t \geq 0),
\]
we have
\[
b^{-\lb x|y_1\rb_\notzero} \leq \log(b)\dist(x,y) + b^{-\lb y|y_1\rb_\notzero}.
\]
It follows that
\begin{align*}
\wbar\Dist(x,z) \leq \Dist(x,z) &\leq \sum_{i = 0}^{n - 1} b^{-\lb x_i|x_{i + 1}\rb_\notzero}\\
&= b^{-\lb x|y_1\rb_\notzero} + \sum_{i = 1}^{n - 1} b^{-\lb y_i|y_{i + 1}\rb_\notzero}\\
&\leq \log(b)\dist(x,y) + b^{-\lb y|y_1\rb_\notzero} + \sum_{i = 1}^{n - 1} b^{-\lb y_i|y_{i + 1}\rb_\notzero}\\
&\leq \log(b)\dist(x,y) + \Dist(y,z) + \epsilon = \wbar\Dist(x,y) + \wbar\Dist(y,z) + \epsilon.
\end{align*}
Taking the limit as $\epsilon$ goes to zero finishes the proof.

\item[(3)] The third case is identical to the second.
\end{itemize}

If a sequence $(x_n)_1^\infty$ is Cauchy with respect to $\wbar\Dist$, then Ramsey's theorem (for example) guarantees that some subsequence is Cauchy with respect to either $\dist$ or $\Dist$. This subsequence converges with respect to that metametric, and therefore also with respect to $\wbar\Dist$. It follows that the entire sequence converges, and therefore that $\wbar\Dist$ is complete.

Finally, to show that $\wbar\Dist$ induces the topology $\scrT$, suppose that $U\subset\del X$ is open in $\scrT$, and fix $x\in U$. If $x\in X$, then $B_\dist(x,r)\subset U$ for some $r > 0$. On the other hand, by the triangle inequality $\Dist(x,y) \geq \frac12\Dist(x,x) > 0$ for all $y\in\bord X$. Letting $\w r = \min(r, \frac12\Dist(x,x))$, we have $B_{\wbar\Dist}(x,\w r) \subset B_\dist(x,r) \subset U$. If $x\in\del X$, then $N_t(x)\subset U$ for some $t\geq 0$; letting $C$ be the implied constant of \eqref{distanceasymptotic}, we have $B_{\wbar\Dist}(x,e^{-t}/C) = B_\Dist(x,e^{-t}/C) \subset N_t(x) \subset U$. Thus $U$ is open in the topology generated by the $\wbar\Dist$ metric. The converse direction is similar but simpler, and will be omitted.
\end{proof}

\begin{remark}
\label{remarkDistunifcont}
The proof of Proposition \ref{propositionwbarDist} actually shows more, namely that
\[
\Dist(x,z) \leq \Dist(x,y) + \wbar\Dist(y,z) \all x,y,z\in\bord X.
\]
Since $\Dist(x,x) = b^{-\dox x} = \inf_{y\in\bord X} \Dist(x,y)$, plugging in $x = z$ gives
\[
b^{-\dox x} \leq b^{-\dox y} + \wbar\Dist(x,y) \all x,y\in\bord X.
\]
\end{remark}

\begin{remark}
Although the metric $\wbar\Dist$ is convenient since it induces the correct topology on $\bord X$, it is not a generalization of the Euclidean metric on the closure of an algebraic hyperbolic space. Indeed, when $X = \B^2$, then $\wbar\Dist$ is not bi-Lipschitz equivalent to the Euclidean metric on $\bord X$.
\end{remark}

\subsection{The visual metametric based at a point $\xi\in\del X$}
Our final metametric is supposed to generalize the Euclidean metric on the boundary of the half-space model $\E$. This metric should be thought of as ``seen from the point $\infty$''.

\begin{notation}
If $X$ is a hyperbolic metric space and $\xi\in\del X$, then let $\EE_\xi := \bord X\butnot\{\xi\}$.
\end{notation}

Since we have not yet introduced a formula analogous to \eqref{euclideanball2} for the Euclidean metric on $\del\E\butnot\{\infty\}$, we will instead motivate the visual metametric based at a point $\xi\in\del X$ by considering a sequence $(\notzero_n)_1^\infty$ in $X$ converging to $\xi$, and taking the limits of their visual metametrics.

In fact, $\Dist_{b,\notzero_n}(y_1,y_2)\to 0$ for every $y_1, y_2\in\EE_\xi$. Some normalization is needed.

\begin{lemma}
\label{lemmaeuclideanmetametric}
Fix $\zero\in X$, and suppose $\notzero_n\to\xi\in\del X$. Then for all $y_1, y_2\in\EE_\xi$,
\[
b^{\dox{\notzero_n}} \Dist_{b,\notzero_n} (y_1,y_2) \tendsto{n,\times} b^{-[\lb y_1 | y_2 \rb_\zero - \sum_{i = 1}^2 \lb y_i | \xi \rb_\zero]} .
\]
with $\tendsto n$ if $X$ is strongly hyperbolic.
\end{lemma}
\begin{proof}
\begin{align*}
b^{\dox{\notzero_n}} \Dist_{b,\notzero_n}(y_1,y_2)
&\asymp_\times b^{-[\lb y_1 | y_2 \rb_{\notzero_n} - \dox{\notzero_n}]} \by{\eqref{distanceasymptotic}}\\
&\asymp_\times b^{-[\lb y_1 | y_2 \rb_\zero - \sum_{i = 1}^2 \lb y_i | \notzero_n \rb_\zero]} \by{(k) of Proposition \ref{propositionbasicidentities}}\\
&\tendsto{n,\times} b^{-[\lb y_1 | y_2 \rb_\zero - \sum_{i = 1}^2 \lb y_i | \xi \rb_\zero]}. \by{Lemma \ref{lemmanearcontinuity}}
\end{align*}
In each step, equality holds if $X$ is strongly hyperbolic.
\end{proof}

We can now construct the visual metametric based at a point $\xi\in\del X$.

\begin{proposition}
\label{propositioneuclideanmetametric}
For each $\zero\in X$ and $\xi\in\del X$, there exists a complete metametric $\Dist_{b,\xi,\zero}$ on $\EE_\xi$ satisfying
\begin{equation}
\label{euclideanmetametric}
\Dist_{b,\xi,\zero}(y_1,y_2) \asymp_\times b^{-[\lb y_1 | y_2 \rb_\zero - \sum_{i = 1}^2 \lb y_i | \xi \rb_\zero]} \ ,
\end{equation}
with equality if $X$ is strongly hyperbolic. The metametric $\Dist_{b,\xi,\zero}$ is compatible with the topology $\scrT\given\EE_\xi$; moreover, a set $S\subset\EE_\xi$ is bounded in the metametric $\Dist_{b,\xi,\zero}$ if and only if $\xi\notin\cl S$.
\end{proposition}
\begin{remark}
The metric $\Dist_{b,\xi,\zero}\given\EE_\xi\cap\del X$ has been referred to in the literature as the \emph{Hamenst\"adt metric}.
\end{remark}

\begin{proof}[Proof of Proposition \ref{propositioneuclideanmetametric}]
Let
\begin{align*}
\Dist_{b,\xi,\zero}(y_1,y_2) 
&= \limsup_{\notzero\to\xi}b^{\dox{\notzero}}\Dist_\notzero(y_1,y_2)\\ 
&:= \sup\left\{\limsup_{n\to\infty} b^{\dox{\notzero_n}}\Dist_{\notzero_n}(y_1,y_2) : \notzero_n \tendsto n \xi\right\}
\end{align*}
Since the class of metametrics is closed under suprema and limits, it follows that $\Dist_{b,\xi,\zero}$ is a metametric. The asymptotic \eqref{euclideanmetametric} follows from Lemma \ref{lemmaeuclideanmetametric}.

For the remainder of this proof, we write $\Dist = \Dist_{b,\zero}$ and $\Dist_\xi = \Dist_{b,\xi,\zero}$.

For all $x\in\EE_\xi$,
\begin{equation}
\label{euclideancomparison}
\Dist_\xi(\zero,x) \asymp_\times b^{-[\lb \zero|x\rb_\zero - \lb \zero|\xi\rb_\zero - \lb x|\xi\rb_\zero]} = b^{\lb x|\xi\rb_\zero} \asymp_\times \frac{1}{\Dist(x,\xi)},
\end{equation}
with equality if $X$ is strongly hyperbolic. It follows that for any set $S\subset\EE_\xi$, the function $\Dist_\xi(\zero,\cdot)$ is bounded on $S$ if and only if the function $\Dist(\cdot,\xi)$ is bounded from below on $S$. This demonstrates that $S$ is bounded in the $\Dist_\xi$ metametric if and only if $\xi\notin\cl S$.

Let $(x_n)_1^\infty$ be a Cauchy sequence with respect to $\Dist_\xi$. Since $\Dist \lesssim_\times \Dist_\xi$, it follows that $(x_n)_1^\infty$ is also Cauchy with respect to the metametric $\Dist$, so it converges to a point $x\in\bord X$ with respect to $\Dist$. If $x\in\EE_\xi$, then we have
\[
\Dist_\xi(x_n,x) \asymp_\times b^{\lb x_n | \xi\rb_\zero + \lb x | \xi\rb_\zero} \Dist(x_n,x) \tendsto{n,\times} b^{2\lb x | \xi\rb_\zero} 0 = 0.
\]
On the other hand, if $x = \xi$, then the sequence $(x_n)_1^\infty$ is unbounded in the $\Dist_\xi$ metametric, which contradicts the fact that it is Cauchy. Thus $\Dist_\xi$ is complete.

Finally, given $\eta\in (\EE_\xi)_\refl = \EE_\xi\cap\del X$, consider the filters $\FF_1$ and $\FF_2$ generated by the collections $\{B_{\Dist}(\eta,r) : r > 0\}$ and $\{B_{\Dist_\xi}(\eta,r) : r > 0\}$, respectively. Since $\Dist \lesssim_\times \Dist_\xi$, we have $\FF_2\subset\FF_1$. Conversely, since $B_{\Dist_\xi}(\eta,1)$ is bounded in the $\Dist_\xi$ metametric, its closure does not contain $\xi$, and so the function $\lb \cdot | \xi\rb_\zero$ is bounded on this set. Thus $\Dist \asymp_{\times,\eta} \Dist_\xi$ on $B_{\Dist_\xi}(\eta,1)$. Letting $C$ be the implied constant of the asymptotic, we have $B_{\Dist_\xi}(\eta,\min(r,1)) \subset B_{\Dist}(\eta,Cr)$, which demonstrates that $\FF_1\subset\FF_2$. Thus $\Dist_\xi$ is compatible with the topology $\scrT\given\EE_\xi$.
\end{proof}

From Lemma \ref{lemmaeuclideanmetametric} and Proposition \ref{propositioneuclideanmetametric} it immediately follows that
\begin{equation}
\label{euclideanmetrictendsasymp}
b^{d(\zero,\notzero_n)} \Dist_{b,\notzero_n} (y_1,y_2) \tendsto{n,\times} \Dist_{b,\xi,\zero}(y_1,y_2) 
\end{equation}
whenever $(\notzero_n)_1^\infty\in\xi$.

\begin{remark}
It is not clear whether a result analogous to Proposition \ref{propositionwbarDist} holds for the metametric $\Dist_{b,\xi,\zero}$. A straightforward adaptation of the proof of Proposition \ref{propositioneuclideanmetametric} does not work, since
\begin{align*}
b^{\dox{\notzero_n}}\wbar\Dist_{b,\notzero_n}(x,y)
&= \min(b^{\dox{\notzero_n}}\dist(x,y)) , b^{\dox{\notzero_n}}\Dist_{b,\notzero_n}(x,y)\\
&\tendsto{n,\times} \min(\infty, \Dist_{b,\xi,\zero}(x,y))\\
&= \Dist_{b,\xi,\zero}(x,y).
\end{align*}
\end{remark}

We finish this chapter by describing the relation between the visual metametric based at $\infty$ and the Euclidean metric on the boundary of the half-space model $\E$.

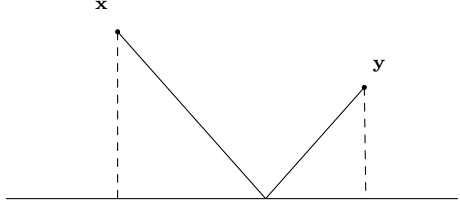
\begin{figure}[h!]
\begin{center}
\begin{tikzpicture}[line cap=round,line join=round,>=triangle 45,x=1.0cm,y=1.0cm]
\clip(-3.801,-0.18) rectangle (3.96,4.50);
\draw (-3.0,0.0)-- (3.0,0.0);
\draw (0.44800000000000006,-0.0)-- (1.76,1.48);
\draw [dashed] (-1.52,2.22)-- (-1.52,0);
\draw [dashed] (1.76,1.48)-- (1.78,0);
\draw (0.44800000000000006,-0.0)-- (-1.52,2.22);
\begin{scriptsize}
\draw[color=black] (0.13316978586974992,4.270045192923707) node {$\infty$};
\draw [fill=black] (-1.52,2.22) circle (.75pt);
\draw[color=black] (-1.7346933526424422,2.5810200144818323) node {$\xx$};
\draw [fill=black] (1.76,1.48) circle (.75pt);
\draw[color=black] (1.961291155477427,1.7663137519392818) node {$\yy$};
\draw [fill=black] (1.78,-6.938893903907228) circle (1.5pt);
\end{scriptsize}
\end{tikzpicture}
\caption[The Hamenst\"adt distance]{The Hamenst\"adt distance $\Dist_{e,\infty,\zero}(\xx,\yy)$ between two points $\xx,\yy\in\E$ is coarsely asymptotic to the maximum of the following three quantities: $x_1$, $y_1$, and $\|\yy - \xx\|$. Equivalently, $\Dist_{e,\infty,\zero}(\xx,\yy)$ is coarsely asymptotic to the length of the shortest path which both connects $\xx$ and $\yy$ and touches $\BB$.}
\label{figurehamenstadt}
\end{center}
\end{figure}

\begin{proposition}[Cf. Figure \ref{figurehamenstadt}]
\label{propositionhamenstadt}
Let $X = \E = \E^\alpha$, let $\zero = (1,\0)\in X$, and fix $\xx,\yy\in \EE_\infty = \E\cup\BB$. We have
\begin{equation}
\label{hamenstadt}
\Dist_{e,\infty,\zero}(\xx,\yy) \asymp_\times \max(x_1,y_1,\|\yy - \xx\|),
\end{equation}
with equality if $\xx,\yy\in\BB = \del\E\butnot\{\infty\}$.
\end{proposition}
\begin{proof}
First suppose that $\xx,\yy\in\E$. By (h) of Proposition \ref{propositionbasicidentities},
\begin{align*}
\Dist_{e,\infty,\zero}(\xx,\yy)
&=_\pt \exp\left(\frac12\big(\dist(\xx,\yy) + \busemann_\infty(\zero,x) + \busemann_\infty(\zero,\yy)\big)\right)\\
&=_\pt \sqrt{x_1 y_1} \exp\left(\frac12\left(\cosh^{-1}\left(1 + \frac{\|\yy - \xx\|^2}{2 x_1 y_1}\right)\right)\right)\\ 
\by{\eqref{distE} and Proposition \ref{propositionbusemannE}}\\
&\asymp_\times \sqrt{x_1 y_1} \sqrt{1 + \frac{\|\yy - \xx\|^2}{2 x_1 y_1}}\\ 
\since{$e^{t/2} \asymp_\times \sqrt{\cosh(t)}$}\\
&=_\pt \sqrt{x_1 y_1 + \|\yy - \xx\|^2}\\
&\asymp_\times \max(\sqrt{x_1 y_1}, \|\yy - \xx\|).
\end{align*}
Since $\sqrt{x_1 y_1}\leq \max(x_1,y_1)$, this demonstrates the $\lesssim$ direction of \eqref{hamenstadt}. Since $y_1 \leq x_1 + \|\yy - \xx\|$ and $x_1 \leq y_1 + \|\yy - \xx\|$, we have
\[
\max(x_1,y_1) \lesssim_\times \max(\min(x_1,y_1),\|\yy - \xx\|) \leq \max(\sqrt{x_1 y_1}, \|\yy - \xx\|)
\]
which demonstrates the reverse inequality. Thus \eqref{hamenstadt} holds for $\xx,\yy\in\E$; a continuity argument demonstrates \eqref{hamenstadt} for $\xx,\yy\in\EE_\infty$.

If $\xx,\yy\in\BB$, then
\begin{align*}
\Dist_{e,\infty,\zero}(\xx,\yy)
&= \lim_{a,b\to 0} \sqrt{a b} \exp\left(\frac12\left(\cosh^{-1}\left(1 + \frac{\|\yy - \xx\|^2}{2 a b}\right)\right)\right) \noreason\\
&= \lim_{a,b\to 0} \sqrt{a b} \sqrt{2\left(1 + \frac{\|\yy - \xx\|^2}{2 a b}\right)} \since{$\lim_{t\to\infty} e^{t/2}/\sqrt{2\cosh(t)} = 1$}\\
&= \lim_{a,b\to 0} \sqrt{2\left(ab + \frac{\|\yy - \xx\|^2}{2}\right)} = \sqrt{\|\yy - \xx\|^2} = \|\yy - \xx\|.\noreason
\end{align*}
\end{proof}

\begin{corollary}[Cf. {\cite[Fig. 5]{Sullivan_entropy}}]
\label{corollarysullivansformula}
For $\xx,\yy\in\E$, we have 
\[
e^{\dist(\xx,\yy)} \asymp_\times \frac{\max(x_1^2,y_1^2,\|\yy - \xx\|^2)}{x_1 y_1}\cdot
\]
\end{corollary}
\begin{proof}
The result follows from 
\begin{align*}
\max(x_1^2,y_1^2,\|\yy - \xx\|^2) \asymp_\times \Dist_{e,\infty,\zero}(\xx,\yy)^2
&= \exp\big(\dist(\xx,\yy) + \busemann_\infty(\zero,\xx) + \busemann_\infty(\zero,\yy)\big)\\
&= x_1 y_1 e^{\dist(\xx,\yy)}
\end{align*}
which may be easily rearranged to complete the proof.
\end{proof}




\chapter{More about the geometry of hyperbolic metric spaces}\label{sectiongeometry2}

In this chapter we discuss various topics regarding the geometry of hyperbolic metric spaces, including metric derivatives, the Rips condition for hyperbolicity, construction of geodesic rays and lines in CAT(-1) spaces, ``shadows at infinity'', and some functions which we call ``generalized polar coordinates''. We start by introducing some conventions to apply in the remainder of the paper.

\section{Gromov triples}
\label{standingassumptions2}

The following definition is made for convenience of notation:

\begin{definition}
\label{definitiongromovtriple}
A \emph{Gromov triple} is a triple $(X,\zero,b)$, where $X$ is a hyperbolic metric space, $\zero\in X$, and $b > 1$ is close enough to $1$ to guarantee for every $\notzero\in X$ the existence of a visual metametric $\Dist_{b,\notzero}$ via Proposition \ref{propositionDist} above.
\end{definition}

\begin{notation}
Let $(X,\zero,b)$ be a Gromov triple. Given $\notzero\in X$ and $\xi\in\del X$, we will let $\Dist_\notzero = \Dist_{b,\notzero}$ be the metametric defined in Proposition \ref{propositionDist}, we will let $\wbar\Dist_\notzero = \wbar\Dist_{b,\notzero}$ be the metric defined in Proposition \ref{propositionwbarDist}, and we will let $\Dist_{\xi,\notzero} = \Dist_{b,\xi,\notzero}$ be the metametric defined in Proposition \ref{propositioneuclideanmetametric}. If $\notzero = \zero$, then we use the further shorthand $\Dist = \Dist_\zero$, $\wbar\Dist = \wbar\Dist_\zero$, and $\Dist_\xi = \Dist_{\xi,\zero}$.

We will denote the diameter of a set $S$ with respect to the metametric $\Dist$ by $\Diam(S)$.
\end{notation}

\begin{convention}
\label{conventiontriples}
For the remainder of the paper, with the exception of Chapter \ref{sectiondiscreteness}, all statements should be assumed to be universally quantified over Gromov triples $(X,\zero,b)$ unless context indicates otherwise.
\end{convention}

\begin{convention}
For the remainder of the paper, whenever we make statements of the form ``Let $X = Y$'', where $Y$ is a hyperbolic metric space, we implicitly want to ``beef up'' $X$ into a Gromov triple $(X,\zero,b)$ whose underlying hyperbolic metric space is $Y$. For general $Y$, this may be done arbitrarily, but if $Y$ is strongly hyperbolic, we want to set $b = e$, and if $Y$ is an algebraic hyperbolic space, then we want to set $\zero = [(1,\0)]$, $\zero = \0$, or $\zero = (1,\0)$ depending on whether $Y$ is the hyperboloid model $\H$, the ball model $\B$, or the half-space model $\E$, respectively.

For example, when saying ``Let $X = \H = \H^\infty$'', we really mean ``Let $X = \H = \H^\infty$, let $\zero = [(1,\0)]$, and let $b = e$.''
\end{convention}

\begin{convention}
\label{conventionstandard}
The term ``Standard Case'' will always refer to the finite-dimensional situation where $X = \H^d$ for some $2\leq d < \infty$.
\end{convention}

\bigskip
\section{Derivatives}\label{subsectionderivatives}

\subsection{Derivatives of metametrics}

Let $(Z,\scrT)$ be a perfect topological space, and let $\Dist_1$ and $\Dist_2$ be two metametrics on $Z$. The \emph{metric derivative} of $\Dist_1$ with respect to $\Dist_2$ is the function $\Dist_1/\Dist_2:Z\to[0,\infty]$ defined by
\[
\frac{\Dist_1}{\Dist_2}(z) := \lim_{w\to z}\frac{\Dist_1(z,w)}{\Dist_2(z,w)},
\]
assuming the limit exists. If the limit does not exist, then we can speak of the \emph{upper} and \emph{lower} derivatives; these will be denoted $(\Dist_1/\Dist_2)^*(z)$ and $(\Dist_1/\Dist_2)_*(z)$, respectively. Note that the chain rule for metric derivatives takes the following form:
\[
\frac{\Dist_1}{\Dist_3} = \frac{\Dist_1}{\Dist_2}\frac{\Dist_2}{\Dist_3},
\]
assuming all limits exist.

We proceed to calculate the derivatives of the metametrics that were introduced in Section \ref{subsectionmetametrics}:

\begin{observation}
\label{observationGMVT}
Fix $y_1,y_2\in\bord X$.
\begin{itemize}
\item[(i)] For all $\notzero_1,\notzero_2\in X$, we have
\begin{equation}
\label{GMVT1}
\frac{\Dist_{\notzero_1}(y_1,y_2)}{\Dist_{\notzero_2}(y_1,y_2)} \asymp_\times b^{-\frac12[\busemann_{y_1}(\notzero_1,\notzero_2) + \busemann_{y_2}(\notzero_1,\notzero_2)]}.
\end{equation}
\item[(ii)] For all $\notzero\in X$ and $\xi\in\del X$, we have
\begin{equation}
\label{GMVT2}
\frac{\Dist_{\xi,\notzero}(y_1,y_2)}{\Dist_\notzero(y_1,y_2)} \asymp_\times b^{-[\lb y_1|\xi\rb_\notzero + \lb y_2|\xi\rb_\notzero]}.
\end{equation}
\item[(iii)] For all $\notzero_1,\notzero_2\in X$ and $\xi\in\del X$, we have
\begin{equation}
\label{GMVT3}
\frac{\Dist_{\xi,\notzero_1}(y_1,y_2)}{\Dist_{\xi,\notzero_2}(y_1,y_2)} \asymp_\times b^{\busemann_\xi(\notzero_1,\notzero_2)}.
\end{equation}
\end{itemize}
In each case, equality holds if $X$ is strongly hyperbolic.
\end{observation}
\begin{proof}
(i) follows from (g) of Proposition \ref{propositionbasicidentities}, while (ii) is immediate from \eqref{distanceasymptotic} and \eqref{euclideanmetametric}. (iii) follows from \eqref{euclideanmetrictendsasymp}.
\end{proof}

Combining with Lemma \ref{lemmanearcontinuity} yields the following:

\begin{corollary}
\label{corollarymetricderivative}
Suppose that $\bord X$ is perfect. Fix $y\in\bord X$.
\begin{itemize}
\item[(i)] For all $\notzero_1,\notzero_2\in X$, we have
\begin{equation}
\label{derivative1}
\left(\frac{\Dist_{\notzero_1}}{\Dist_{\notzero_2}}\right)^*(y) \asymp_\times \left(\frac{\Dist_{\notzero_1}}{\Dist_{\notzero_2}}\right)_*(y) \asymp_\times b^{-\busemann_y(\notzero_1,\notzero_2)}.
\end{equation}
\item[(ii)] For all $\notzero\in X$ and $\xi\in\del X$, we have
\begin{equation}
\label{derivative2}
\left(\frac{\Dist_{\xi,\notzero}}{\Dist_\notzero}\right)^*(y) \asymp_\times \left(\frac{\Dist_{\xi,\notzero}}{\Dist_\notzero}\right)_*(y) \asymp_\times b^{-2\lb y|\xi\rb_\notzero}.
\end{equation}
\item[(iii)] For all $\notzero_1,\notzero_2\in X$ and $\xi\in\del X$, we have
\begin{equation}
\label{derivative3}
\left(\frac{\Dist_{\xi,\notzero_1}}{\Dist_{\xi,\notzero_2}}\right)^*(y) \asymp_\times \left(\frac{\Dist_{\xi,\notzero_1}}{\Dist_{\xi,\notzero_2}}\right)_*(y) \asymp_\times b^{\busemann_\xi(\notzero_1,\notzero_2)}.
\end{equation}
\end{itemize}
In each case, equality holds if $X$ is strongly hyperbolic.
\end{corollary}
\begin{remark}
In case $\bord X$ is not perfect, \eqref{derivative1} - \eqref{derivative3} may be taken as definitions. We will ignore the issue henceforth.
\end{remark}
Combining Observation \ref{observationGMVT} with Corollary \ref{corollarymetricderivative} yields the following:
\begin{proposition}[Geometric mean value theorem]
\label{propositionGMVT}
Fix $y_1,y_2\in\bord X$.
\begin{itemize}
\item[(i)] For all $\notzero_1,\notzero_2\in X$, we have
\[
\frac{\Dist_{\notzero_1}(y_1,y_2)}{\Dist_{\notzero_2}(y_1,y_2)} \asymp_\times \left(\frac{\Dist_{\notzero_1}}{\Dist_{\notzero_2}}(y_1)\frac{\Dist_{\notzero_1}}{\Dist_{\notzero_2}}(y_2)\right)^{1/2}.
\]
\item[(ii)] For all $\notzero\in X$ and $\xi\in\del X$, we have
\[
\frac{\Dist_{\xi,\notzero}(y_1,y_2)}{\Dist_\notzero(y_1,y_2)} \asymp_\times \left(\frac{\Dist_{\xi,\notzero}}{\Dist_\notzero}(y_1)\frac{\Dist_{\xi,\notzero}}{\Dist_\notzero}(y_2)\right)^{1/2}.
\]
\item[(iii)] For all $\notzero_1,\notzero_2\in X$ and $\xi\in\del X$, we have
\[
\frac{\Dist_{\xi,\notzero_1}(y_1,y_2)}{\Dist_{\xi,\notzero_2}(y_1,y_2)} \asymp_\times \left(\frac{\Dist_{\xi,\notzero_1}}{\Dist_{\xi,\notzero_2}}(y_1)\frac{\Dist_{\xi,\notzero_1}}{\Dist_{\xi,\notzero_2}}(y_2)\right)^{1/2}.
\]
\end{itemize}
In each case, equality holds if $X$ is strongly hyperbolic.
\end{proposition}

\subsection{Derivatives of maps}
\label{subsubsectionderivativemaps}
As before, let $(Z,\scrT)$ be a perfect topological space, and now fix just one metametric $\Dist$ on $Z$. For any map $g:Z\to Z$, the \emph{metric derivative} of $G$ is the function $g':Z\to(0,\infty)$ defined by
\[
g'(z) := \frac{\Dist\circ g}{\Dist}(z) = \lim_{w\to z}\frac{\Dist(g(z),g(w))}{\Dist(z,w)}\cdot
\]
If the limit does not exist, the upper and lower metric derivatives will be denoted $\overline g'$ and $\underline g'$, respectively.

\begin{remark}
\label{remarkconfusion}
To avoid confusion, in what follows $g'$ will always denote the derivative of an isometry $g\in\Isom(X)$ with respect to the metametric $\Dist = \Dist_{b,\zero}$, rather than with respect to any other metametric.
\end{remark}

\begin{proposition}
\label{propositionderivativebusemann}
For all $g\in\Isom(X)$,
\begin{align*}
\overline g'(y) \asymp_\times \underline g'(y) &\asymp_\times b^{-\busemann_y(g^{-1}(\zero),\zero)} \all y\in\bord X\\
\frac{\Dist(g(y_1),g(y_2))}{\Dist(y_1,y_2)} &\asymp_\times \left(\overline g'(y_1) \overline g'(y_2)\right)^{1/2} \all y_1,y_2\in\bord X,
\end{align*}
with equality if $X$ is strongly hyperbolic.
\end{proposition}
\begin{proof}
This follows from (i) of Corollary \ref{corollarymetricderivative}, (i) of Proposition \ref{propositionGMVT}, and the fact that $\Dist\circ g = \Dist_{g^{-1}(\zero)}$.
\end{proof}

\begin{corollary}
\label{corollaryproductoneorig}
For any distinct $y_1,y_2\in\Fix(g)\cap\del X$ we have
\[
\overline g'(y_1) \overline g'(y_2) \asymp_\times 1,
\]
with equality if $X$ is strongly hyperbolic.
\end{corollary}

The next proposition shows the relation between the derivative of an isometry $g\in\Isom(X)$ at a point $\xi\in\Fix(g)$ and the action on the metametric space $(\EE_\xi,\Dist_\xi)$:

\begin{proposition}
\label{propositioneuclideansimilarity}
Fix $g\in\Isom(X)$ and $\xi\in\Fix(g)$. Then for all $y_1,y_2\in\EE_\xi$,
\[
\frac{\Dist_\xi(g(y_1), g(y_2))}{\Dist_\xi(y_1,y_2)} \asymp_\times \frac{1}{g'(\xi)},
\]
with equality if $X$ is strongly hyperbolic.
\end{proposition}
\begin{proof}
\begin{align*}
\frac{\Dist_\xi(g(y_1), g(y_2))}{\Dist_\xi(y_1,y_2)} &=_\pt \frac{\Dist_{\xi,g^{-1}(\zero)}(y_1, y_2)}{\Dist_{\xi,\zero}(y_1,y_2)} \\
&\asymp_\times b^{-\busemann_\xi(\zero,g^{-1}(\zero))} \by{\eqref{GMVT3}} \\
&\asymp_\times 1/g'(\xi). \by{Proposition \ref{propositionderivativebusemann}}
\end{align*}
\end{proof}

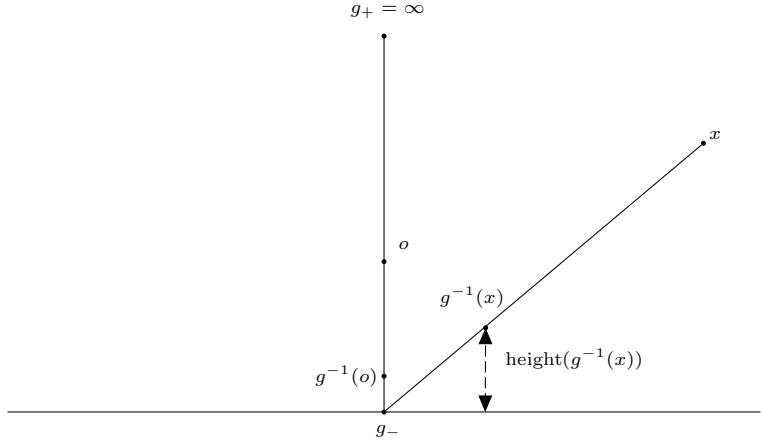
\begin{figure}
\begin{center}
\begin{tikzpicture}[line cap=round,line join=round,>=triangle 45,x=1.0cm,y=1.0cm]
\clip(-5.8,-0.36) rectangle(6.01,5.53);
\draw (-5.0,0.0)-- (5.0,0.0);
\draw [<->, dash pattern=on 4pt off 2pt] (1.35,1.12)-- (1.3452937374574478,0.0);
\draw (0.0,5.0)-- (0.0,0.0);
\draw (0.0,0.0)-- (4.246442867487024,3.5745642370946995);
\begin{scriptsize}
\draw[color=black] (0.057057757531416,-0.26212992254098266) node {$g_-$};
\draw [fill=black] (0.0,5.0) circle (.75pt);
\draw[color=black] (0.04808489249598222,5.322254959074134) node {$g_+ = \infty$};
\draw [fill=black] (1.35,1.12) circle (.75pt);
\draw[color=black] (1.17,1.53) node {$g^{-1}(x)$};
\draw [fill=black] (0,0.47635850682121) circle (.75pt);
\draw[color=black] (-0.5,0.47635850682121) node {$g^{-1}(\zero)$};
\draw[color=black] (2.5272758356701829,0.7026278596223968) node {$\mathrm{height}(g^{-1}(x))$};
\draw [fill=black] (0.0,0.0) circle (.75pt);
\draw [fill=black] (4.246442867487024,3.5745642370946995) circle (.75pt);
\draw[color=black] (4.4036034520894525,3.6885775833822156) node {$x$};
\draw [fill=black] (0.0,2.0) circle (.75pt);
\draw[color=black] (0.2618688112016899,2.2273323702788814) node {$\zero$};
\end{scriptsize}
\end{tikzpicture}
\caption[The derivative of $g$ at $\infty$]{The derivative of $g$ at $\infty$ is equal to the reciprocal of the dilatation ratio of $g$. In particular, $\infty$ is an attracting fixed point if and only if $g$ is expanding, and $\infty$ is a repelling fixed point if and only if $g$ is contracting.}
\end{center}
\end{figure}

\begin{remark}
Proposition \ref{propositioneuclideansimilarity} can be interpreted as a geometric mean value theorem for the action of $g$ on the metametric space $(\EE_\xi,\Dist_\xi)$. Specifically, it tells us that the derivative of $g$ on this metametric space is identically $1/g'(\xi)$.
\end{remark}

\begin{remark}
If $g'(\xi) = 1$, then Proposition \ref{propositioneuclideansimilarity} tells us that the bi-Lipschitz constant of $g$ is independent of $g$, and that $g$ is an isometry if $X$ is strongly hyperbolic. This special case will be important in Chapter \ref{sectionparabolic}.
\end{remark}

\begin{example}
\label{examplegprimeinfty}
Suppose that $X = \E = \E^\alpha$ is the half-space model of a real hyperbolic space, let $\BB = \del\E\butnot\{\infty\}$, let $g(\xx) = \lambda T(\xx) + \bb$ be a similarity of $\BB$, and consider the Poincar\'e extension $\what g\in\Isom(\E)$ defined in Observation \ref{observationpoincareextension}. Clearly $\what g$ acts as a similarity on the metametric space $(\EE_\infty,\Dist_\infty)$ in the following sense: For all $y_1,y_2\in \EE_\infty$,
\[
\Dist_\infty(g(y_1),g(y_2)) = \lambda \Dist_\infty(y_1,y_2).
\]
Comparing with Proposition \ref{propositioneuclideansimilarity} shows that $g'(\infty) = 1/\lambda$.
\end{example}

\subsection{The dynamical derivative} \label{subsubsectiondynamicalderivative}
We can interpret Corollary \ref{corollarymetricderivative} as saying that the metric derivative is well-defined only up to an asymptotic in a general hyperbolic metric space (although it is perfectly well defined in a strongly hyperbolic metric space). Nevertheless, if $\xi$ is a fixed point of the isometry $g$, then we can iterate in order to get arbitrary accuracy.

\begin{proposition}
\label{propositiondynamicalderivative}
Fix $g\in\Isom(X)$ and $\xi\in\Fix(g)$. Then
\[
g'(\xi) := \lim_{n\to\infty} \big( (\overline{g^n})'(\xi) \bigr)^{1/n} = \lim_{n\to\infty} \bigl( (\underline{g^n})'(\xi) \bigr)^{1/n}.
\]
Furthermore
\[
\underline g'(\xi)\leq g'(\xi) \leq \overline g'(\xi).
\]
\end{proposition}
The number $g'(\xi)$ will be called the \emph{dynamical derivative} of $g$ at $\xi$.
\begin{proof}[Proof of Proposition \ref{propositiondynamicalderivative}]
The limits converge due to the submultiplicativity and supermultiplicativity of the expressions inside the radicals, respectively. To see that they converge to the same number, note that by Corollary \ref{corollarymetricderivative}
\[
\lim_{n\to\infty} \left(\frac{(\overline{g^n})'(\xi) }{(\underline{g^n})'(\xi) }\right)^{1/n} \leq \lim_{n\to\infty} C^{1/n} = 1
\]
for some constant $C$ independent of $n$.
\end{proof}
\begin{remark}
Let $\beta_\xi$ denote the \emph{Busemann quasicharacter} of \cite[p.14]{CCMT}. Then $\beta_\xi$ is related to the dynamical derivative via the following formula: $g'(\xi) = b^{-\beta_\xi(g)}$.
\end{remark}
Note that although the dynamical derivative is ``well-defined'', it is not necessarily the case that the chain rule holds for any two $g,h\in\Stab(\Isom(X);\xi)$ (although it must hold up to a multiplicative coarse asymptotic). For a counterexample see \cite[Example 3.12]{CCMT}. Note that this counterexample includes the possibility of two elementa $g,h\in\Stab(\Isom(X);\xi)$ such that $g'(\xi) = h'(\xi) = 1$ but $(gh)'(\xi)\neq 1$. A sufficient condition for the chain rule to hold exactly is given in \cite[Corollary 3.9]{CCMT}.

Despite the failure of the chain rule, the following ``iteration'' version of the chain rule holds:

\begin{proposition}
\label{propositionderivativeiteration}
Fix $g\in\Isom(X)$ and $\xi\in\Fix(g)$. Then
\[
(g^n)'(\xi) = [g'(\xi)]^n \all n\in\Z.
\]
In particular
\begin{equation}
\label{derivativeiteration2}
(g^{-1})'(\xi) = \frac{1}{g'(\xi)}\cdot
\end{equation}
\end{proposition}
\begin{proof}
The only difficulty lies in establishing \eqref{derivativeiteration2}:
\begin{align*}
(g^{-1})'(\xi) = \lim_{n\to\infty} \big((\overline{g^{-n}})'(\xi)\big)^{1/n}
&= \exp_{1/b} \lim_{n\to\infty} \frac{1}{n} \busemann_\xi(\zero,g^n(\zero))\\
&= \exp_{1/b} \lim_{n\to\infty} \frac{1}{n} \busemann_\xi(g^{-n}(\zero),\zero)\\
&= \exp_{1/b} \left( -\lim_{n\to\infty} \frac{1}{n} \busemann_\xi(\zero,g^{-n}(\zero)) \right) \\
&= \frac 1{g'(\xi)}\cdot
\end{align*}
\end{proof}

Combining with Corollary \ref{corollaryproductoneorig} yields the following:

\begin{corollary}
\label{corollaryproductone}
For any distinct $y_1,y_2\in\Fix(g)\cap\del X$ we have
\[
g'(y_1) g'(y_2) = 1.
\]
\end{corollary}

We end this section with the following result relating the dynamical derivative with the Busemann function:

\begin{proposition}
\label{propositionbusemannpreserved}
Fix $g\in\Isom(X)$ and $\xi\in\Fix(g)$. Then for all $x\in X$ and $n\in\Z$,
\[
\busemann_\xi(x , g^{-n}(x)) \asymp_\plus n\log_b g'(\xi).
\]
with equality if $X$ is strongly hyperbolic.
\end{proposition}
\begin{proof}
If $x = \zero$, then
\begin{align*}
b^{-\busemann_\xi(\zero , g^{-n}(\zero))}
&\asymp_\times b^{\busemann_\xi(g^{-n}(\zero),\zero)}\\
&\asymp_\times (\overline{g^n})'(\xi) \by{Proposition \ref{propositionderivativebusemann}}\\
&\asymp_\times (g^n)'(\xi) = (g'(\xi))^n.
\end{align*}
For the general case, we note that
\begin{align*}
\busemann_\xi(x , g^{-n}(x))
&\asymp_\plus \busemann_\xi(x,\zero) + \busemann_\xi(\zero , g^{-n}(\zero)) + \busemann_\xi(g^{-n}(\zero), g^{-n}(x))\\
&\asymp_\plus \busemann_\xi(x,\zero) + n\log_b g'(\xi) + \busemann_\xi(\zero,x)\\
&\asymp_\plus n\log_b g'(\xi).
\end{align*}
\end{proof}

\bigskip
\section{The Rips condition}\label{subsectionrips}

In this section, in addition to assuming that $X$ is a hyperbolic metric space (cf. \6\ref{standingassumptions2}), we assume that $X$ is geodesic. Recall (Section \ref{subsectionCAT}) that $\geo xy$ denotes the geodesic segment connecting two points $x,y\in X$.

\begin{proposition}~
\label{propositionrips}
\begin{itemize}
\item[(i)] For all $x,y,z\in X$,
\[
\dist(z,\geo xy) \asymp_\plus \lb x|y\rb_z.
\]
\item[(ii)] (Rips' thin triangles condition) For all $x,y_1,y_2\in X$ and for any $z\in\geo{y_1}{y_2}$, we have
\[
\min_{i = 1}^2 \dist(z,\geo x{y_i}) \asymp_\plus 0.
\]
\end{itemize}
\end{proposition}

\begin{figure}
\begin{center}
\begin{tikzpicture}[line cap=round,line join=round,>=triangle 45,scale=0.9]
\clip(-6.66,-3.12) rectangle (6.24,3.31);
\draw(0.0,0.0) circle (3.0cm);
\draw [shift={(3.001597531764653,1.959673853196985)}] plot[domain=2.7288669538659853:3.9124414341882043,variable=\t]({1.0*2.037582193265268*cos(\t r)+-0.0*2.037582193265268*sin(\t r)},{0.0*2.037582193265268*cos(\t r)+1.0*2.037582193265268*sin(\t r)});
\draw (0.0,-0.0)-- (1.309924961105806,0.823898739896739);
\begin{scriptsize}
\draw[color=black] (0.75,0.23) node {$\nwarrow$};
\draw (.85,0.3) node[anchor=north west] {$\asymp_\plus \langle x,y \rangle_z$};
\draw [fill=black] (0.0,-0.0) circle (.75pt);
\draw[color=black] (-0.19452753572169953,-0.12825951916697592) node {$z$};
\draw [fill=black] (1.135109224427748,2.7769636383321687) circle (.75pt);
\draw[color=black] (1.252879052282513,3.0161065168421803) node {$x$};
\draw [fill=black] (1.54,0.54) circle (.75pt);
\draw[color=black] (1.7519847722839659,0.5372147741682952) node {$y$};
\end{scriptsize}
\end{tikzpicture}
\caption[The Rips condition]{An illustration of Proposition \ref{propositionrips}(i).}
\end{center}
\end{figure}
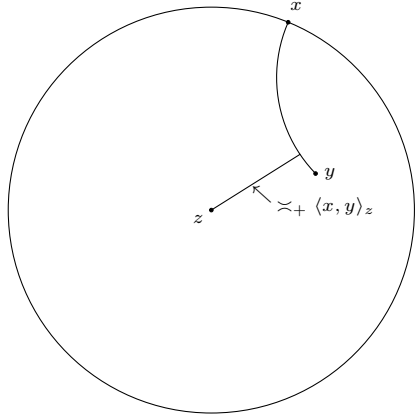

\noindent In fact, the thin triangles condition is equivalent to hyperbolicity; see e.g. \cite[Proposition III.H.1.22]{BridsonHaefliger}.
\begin{proof}~
\begin{itemize}
\item[(i)] By the intermediate value theorem, there exists $w\in\geo xy$ such that $\lb x|z\rb_w = \lb y|z\rb_w$. Applying Gromov's inequality gives $\lb x|z\rb_w = \lb y|z\rb_w \lesssim_\plus \lb x|y\rb_w = 0$. Now (k) of Proposition \ref{propositionbasicidentities} shows that 
\[
\dist(z,\geo xy) \leq \dist(z,w) \lesssim_\plus \lb x|y\rb_z.
\] 
The other direction is immediate, since for each $w\in\geo xy$, we have $\lb x|y\rb_w = 0$, and so (d) of Proposition \ref{propositionbasicidentities} gives $\lb x|y\rb_w \leq \dist(z,w)$.
\item[(ii)] This is immediate from (i), Gromov's inequality, and the equation 
\[
\lb y_1|y_2\rb_z = 0.
\]
\end{itemize}
\end{proof}

The next lemma demonstrates the correctness of the intuitive notion that if two points are close to each other, then the geodesic connecting them should not be very large.

\begin{lemma}
\label{lemmageodesicdiameter}
Fix $x_1,x_2\in \bord X$. We have
\[
\Diam(\geo{x_1}{x_2}) \asymp_\times \Dist(x_1,x_2).
\]
\end{lemma}
\begin{proof}
It suffices to show that if $y\in\geo{x_1}{x_2}$, then
\[
\Dist(y,\{x_1,x_2\}) \lesssim_\times \Dist(x_1,x_2).
\]
Indeed, by the thin triangles condition, we may without loss of generality suppose that $\dist(y, \geo{\zero}{x_1}) \asymp_\plus 0$. Write $\dist(y,z) \asymp_\plus 0$ for some $z\in \geo\zero{x_1}$. Then
\begin{align*}
\Dist(x_1,y) \asymp_\times \Dist(x_1,z) = e^{-\dox z} \asymp_\times e^{-\dox y} 
&\leq e^{-\dist(\zero,\geo{x_1}{x_2})}\\ 
&\asymp_\times e^{-\lb x_1|x_2\rb_\zero}\\
&\asymp_\times \Dist(x_1,x_2).
\end{align*}
\end{proof}

%
%

\bigskip
\section{Geodesics in CAT(-1) spaces} \label{subsectiongeodesicsCAT}

\begin{observation}
Any isometric embedding $\pi:\CO t\infty\to X$ extends uniquely to a continuous map $\pi:[t,\infty]\to \bord X$. Similarly, any isometric embedding $\pi:(-\infty,+\infty)\to X$ extends uniquely to a continuous map $\pi:[-\infty,+\infty]\to \bord X$. Abusing terminology, we will also call the extended maps ``isometric embeddings''.
\end{observation}

\begin{definition}
\label{definitiongeodesic}
Fix $x\in X$ and $\xi,\eta\in\del X$.
\begin{itemize}
\item A \emph{geodesic ray} connecting $x$ and $\xi$ is the image of an isometric embedding $\pi:[0,\infty]\to X$ satisfying
\[
\pi(0) = x,\hspace{.5 in}\pi(\infty) = \xi.
\]
\item A \emph{geodesic line} or \emph{bi-infinite geodesic} connecting $\xi$ and $\eta$ is the image of an isometric embedding $\pi:[-\infty,+\infty]\to X$ satisfying
\[
\pi(-\infty) = \xi,\hspace{.5 in}\pi(+\infty) = \eta.
\]
\end{itemize}
When we do not wish to distinguish between geodesic segments (cf. Section \ref{subsectionCAT}), geodesic rays, and geodesic lines, we shall simply call them geodesics. For $x,y\in\bord X$, any geodesic connecting $x$ and $y$ will be denoted $\geo xy$.
\end{definition}

\begin{notation}
\label{notationgeopqt2}
Extending Notation \ref{notationgeopqt}, if $\geo x\xi$ is the image of the isometric embedding $\pi:[0,\infty]\to X$, then for $t\in[0,\infty]$ we let $\geo x\xi_t = \pi(t)$, i.e. $\geo x\xi_t$ is the unique point on the geodesic ray $\geo x\xi$ such that $\dist(x,\geo x\xi_t) = t$.
\end{notation}


The main goal of this section is to prove the following:

\begin{proposition}
\label{propositiongeodesicconvergence}
Suppose that $X$ is a complete CAT(-1) space. Then:
\begin{itemize}
\item[(i)] For any two distinct points $x,y\in\bord X$, there is a unique geodesic $\geo xy$ connecting them.
\item[(ii)] Suppose that $(x_n)_1^\infty$ and $(y_n)_1^\infty$ are sequences in $\bord X$ which converge to points $x_n\to x\in\bord X$ and $y_n\to y\in\bord X$, with $x\neq y$. Then $\geo{x_n}{y_n}\to\geo xy$ in the Hausdorff metric on $(\bord X,\wbar\Dist)$. If $x = y$, then $\geo{x_n}{y_n}\to \{x\}$ in the Hausdorff metric.
\end{itemize}
\end{proposition}

\begin{definition}
\label{definitionregularlygeodesic}
A hyperbolic metric space $X$ satisfying the conclusion of Proposition \ref{propositiongeodesicconvergence} will be called \emph{regularly geodesic}.
\end{definition}

\begin{remark}
The existence of a geodesic connecting any two points in $\bord X$ was proven in \cite[Proposition 0.2]{Buyalo} under the weaker hypothesis that $X$ is a Gromov hyperbolic complete CAT(0) space. However, this weaker hypothesis does not imply the uniqueness of such a geodesic, nor does it imply (ii) of Proposition \ref{propositiongeodesicconvergence}, as shown by the following example:
\end{remark}

\begin{example}[A proper and uniquely geodesic hyperbolic CAT(0) space which is not regularly geodesic]
\label{examplestrip}
Let
\[
X = \{\xx\in\R^2: x_2\in[0,1]\}
\]
be interpreted as a subspace of $\R^2$ with the usual metric. Then $X$ is hyperbolic, proper, and uniquely geodesic, but is not regularly geodesic.
\end{example}
\begin{proof}
It is hyperbolic since it is roughly isometric to $\R$. It is uniquely geodesic since it is a convex subset of $\R^2$. It is proper because it is a closed subset of $\R^2$. It is not regularly geodesic because if we write $\del X = \{\xi_+,\xi_-\}$, then the two points $\xi_+$ and $\xi_-$ have infinitely many distinct geodesics connecting them: for each $t\in[0,1]$, $\R\times\{t\}$ is a geodesic connecting $\xi_+$ and $\xi_-$.
\end{proof}

The proof of Proposition \ref{propositiongeodesicconvergence} will proceed through several lemmas, the first of which is as follows:

\begin{lemma}
\label{lemmaxyz}
Fix $\epsilon > 0$. There exists $\delta = \delta_X(\epsilon) > 0$ such that if $\Delta = \Delta(x,y_1,y_2)$ is a geodesic triangle in $X$ satisfying
\begin{equation}
\label{y1y2delta}
\wbar\Dist(y_1,y_2) \leq \delta,
\end{equation}
then for all $t\in[0,\min_{i = 1}^2 \dist(x,y_i)]$, if $z_i = \geo x{y_i}_t$, then
\begin{equation}
\label{z1z2epsilon}
\wbar\Dist(z_1,z_2) \leq \epsilon.
\end{equation}
\end{lemma}
\begin{proof}
We prove the assertion first for $X = \H^2$ and then in general:
\begin{itemize}
\item[If $X = \H^2$:] Let $\epsilon > 0$, and by contradiction, suppose that for each $\delta = \frac 1n > 0$ there exists a 5-tuple $(x^{(n)},y_1^{(n)},y_2^{(n)},z_1^{(n)},z_2^{(n)})$ satisfying the hypotheses but not the conclusion of the theorem. Since $\bord\H^2$ is compact, there exists a convergent subsequence
\[
(x^{(n_k)},y_1^{(n_k)},y_2^{(n_k)},z_1^{(n_k)},z_2^{(n_k)}) \to (x,y_1,y_2,z_1,z_2)\in\left(\bord\H^2\right)^5.
\]
Taking the limit of \eqref{y1y2delta} as $k\to\infty$ shows that $\wbar\Dist(y_1,y_2) = 0$, so $y_1 = y_2$. Conversely, taking the limit of \eqref{z1z2epsilon} shows that $\wbar\Dist(z_1,z_2) \geq \epsilon > 0$, so $z_1\neq z_2$. Write $y = y_1 = y_2$.

We will take for granted that Proposition \ref{propositiongeodesicconvergence} holds when $X = \H^2$. (This can be proven using the explicit form of geodesics in this space.) It follows that $z_i\in \geo xy$ if $x\neq y$, and $z_i = x$ if $x = y$. The second case is clearly a contradiction, so we assume that $x\neq y$.

Writing $z_i^{(n_k)} = \geo{x^{(n_k)}}{y_i^{(n_k)}}_{t_k}$, we observe that
\begin{align*}
t_k - \dox{x^{(n_k)}}
&= \lb z_i^{(n_k)} | y_i^{(n_k)} \rb_\zero - \lb z_i^{(n_k)} | x^{(n_k)} \rb_\zero - \lb x^{(n_k)} | y_i^{(n_k)} \rb_\zero\\
&\tendsto k \lb z_i | y \rb_\zero - \lb z_i | x \rb_\zero - \lb x |y\rb_\zero.
\end{align*}
Since the left hand side is independent of $i$, so is the right hand side. But the function
\[
z\mapsto \lb z | y \rb_\zero - \lb z | x \rb_\zero - \lb x | y \rb_\zero
\]
is an isometric embedding from $\geo xy$ to $[-\infty,+\infty]$; it is therefore injective. Thus $z_1 = z_2$, a contradiction.

\item[In general:] Let $\epsilon > 0$, and fix $\w\epsilon > 0$ to be determined, depending on $\epsilon$. Let $\w\delta = \delta_{\H^2}(\w\epsilon)$, and fix $\delta > 0$ to be determined, depending on $\w\delta$. Now suppose that $\Delta = \Delta(x,y_1,y_2)$ is a geodesic triangle in $X$ satisfying \eqref{y1y2delta}, fix $t\geq 0$, and let $z_i = \geo x{y_i}_t$. To complete the proof, we must show that $\wbar\Dist(z_1,z_2) \leq \epsilon$.

By contradiction suppose not, i.e. suppose that $\wbar\Dist(z_1,z_2) > \epsilon$. Then $\Dist(x,z_i) > \epsilon/2$ for some $i = 1,2$; without loss of generality suppose $\Dist(x,z_1) > \epsilon/2$. By Proposition \ref{propositionrips} this implies $\dist(\zero,\geo x{z_1}) \asymp_{\plus,\epsilon} 0$; fix $w_1\in\geo x{z_1}$ with $\dox{w_1}\asymp_{\plus,\epsilon} 0$. Let $s = \dist(x,w_1) \leq t$, and let $w_2 = \geo x{z_2}_s$. (See Figure \ref{figurexyz}.)

\begin{figure}
\begin{center}
\begin{tikzpicture}[line cap=round,line join=round,>=triangle 45,x=1.0cm,y=1.0cm]
\clip(-1.539,-1.576) rectangle (3.98,1.8);
\draw [shift={(4.7207100633754004,-4.251144523671373)}] plot[domain=1.948435252123557:2.6141429865240005,variable=\t]({1.0*5.902960522581992*cos(\t r)+-0.0*5.902960522581992*sin(\t r)},{0.0*5.902960522581992*cos(\t r)+1.0*5.902960522581992*sin(\t r)});
\draw [shift={(2.6440776402056905,-5.2367125169303135)}] plot[domain=1.5361702930450398:2.2233807737101374,variable=\t]({1.0*4.980022039672619*cos(\t r)+-0.0*4.980022039672619*sin(\t r)},{0.0*4.980022039672619*cos(\t r)+1.0*4.980022039672619*sin(\t r)});
\draw (0.0,-0.0)-- (0.48334729355344397,-0.1414365427938682);
\draw [shift={(2.802578610044664,0.5103691534918943)}] plot[domain=1.9130096151383342:4.730441794850227,variable=\t]({1.0*0.7701702590181908*cos(\t r)+-0.0*0.7701702590181908*sin(\t r)},{0.0*0.7701702590181908*cos(\t r)+1.0*0.7701702590181908*sin(\t r)});
\draw [shift={(7.3783710087828664,2.3193354663807444)}] plot[domain=3.4843932225971725:3.560165863983011,variable=\t]({1.0*7.320980188110974*cos(\t r)+-0.0*7.320980188110974*sin(\t r)},{0.0*7.320980188110974*cos(\t r)+1.0*7.320980188110974*sin(\t r)});
\begin{scriptsize}
\draw [fill=black] (0.0,-0.0) circle (1pt);
\draw[color=black] (-0.20573038068032204,-0.16484898217738642) node {$\zero$};
\draw [fill=black] (-0.38,-1.28) circle (1pt);
\draw[color=black] (-0.6011417551376018,-1.450757318050246) node {$x$};
\draw [fill=black] (2.544130372359072,1.2358805154383208) circle (1pt);
\draw[color=black] (2.8569090399799677,1.294830750435049) node {$y_1$};
\draw [fill=black] (2.8164815956355223,-0.25967560810823964) circle (1pt);
\draw[color=black] (3.13494327476329,-0.3038660995690469) node {$y_2$};
\draw [fill=black] (0.48334729355344397,-0.1414365427938682) circle (1pt);
\draw[color=black] (0.458863764973814,0.0892241456218923) node {$w_1$};
\draw [fill=black] (0.6894125111628517,-0.6563300275601325) circle (1pt);
\draw[color=black] (0.8759151171487972,-0.8641974071794341) node {$w_2$};
\draw [fill=black] (1.4043001305961753,0.6321297967948407) circle (1pt);
\draw[color=black] (1.3277207486716958,0.8430251189121523) node {$z_1$};
\draw [fill=black] (1.7083588447287483,-0.34538858046728116) circle (1pt);
\draw[color=black] (1.796903519868552,-0.6861631723961132) node {$z_2$};
\end{scriptsize}
\end{tikzpicture}
\caption[The triangle $\Delta(x,y_1,y_2)$]{The triangle $\Delta(x,y_1,y_2)$.}
\label{figurexyz}
\end{center}
\end{figure}
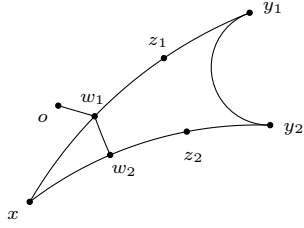

Now let $\wbar{\Delta} = \Delta(\wbar x,\wbar{y_1},\wbar{y_2})$ be a comparison triangle for $\Delta(x,y_1,y_2)$, and let $\wbar{z_1},\wbar{z_2},\wbar{w_1},\wbar{w_2}$ be the corresponding comparison points. Note that $\wbar{z_i} = \geo{\wbar x}{\wbar{y_i}}_t$ and $\wbar{w_i} = \geo{\wbar x}{\wbar{y_i}}_s$. Without loss of generality, suppose that $\wbar{w_1} = \zero_\H$. Then $\dox{y_2} \leq \dox{w_1} + \dist(w_1,y_2) \asymp_{\plus,\epsilon} \dist(\zero_\H,\wbar{y_2})$, and so $\lb y_1|y_2\rb_\zero \lesssim_{\plus,\epsilon} \lb \wbar{y_1}|\wbar{y_2}\rb_{\zero_\H}$, and thus
\[
\wbar\Dist(\wbar{y_1},\wbar{y_2}) \lesssim_{\times,\epsilon} \wbar\Dist(y_1,y_2) \leq \delta.
\]
Setting $\delta$ equal to $\w\delta$ divided by the implied constant, we have
\[
\wbar\Dist(\wbar{y_1},\wbar{y_2}) \leq \w\delta = \delta_\H(\w\epsilon).
\]
Thus $\wbar\Dist(\wbar{z_1},\wbar{z_2})\leq\w\epsilon$ and $\wbar\Dist(\wbar{w_1},\wbar{w_2})\leq\w\epsilon$.
\begin{itemize}
\item If $\dist(\wbar{z_1},\wbar{z_2})\leq\w\epsilon$, then the CAT(-1) inequality finishes the proof (as long as $\w\epsilon\leq \epsilon$). Thus, suppose that
\begin{equation}
\label{z1z2}
\Dist(\wbar{z_1},\wbar{z_2}) \leq \w\epsilon.
\end{equation}
\item If $\Dist(\wbar{w_1},\wbar{w_2})\leq\w\epsilon$, then $0 = \lb \wbar{w_1}|\wbar{w_2}\rb_{\zero_\H} \geq -\log(\w\epsilon)$, a contradiction for $\w\epsilon$ sufficiently small. Thus, suppose that
\begin{equation}
\label{w1w2}
\dist(\wbar{w_1},\wbar{w_2}) \leq\w\epsilon.
\end{equation}
\end{itemize}
By \eqref{z1z2}, we have $\dist(\zero_\H,\wbar{z_i}) \geq -\log(\w\epsilon)$. Applying \eqref{w1w2} gives
\[
\lb z_i|y_i\rb_{w_i} = \dist(w_i,z_i) = \dist(\wbar{w_i},\wbar{z_i}) \gtrsim_\plus -\log(\w\epsilon).
\]
Applying \eqref{w1w2}, the coarse asymptotic $\dox{w_1} \asymp_{\plus,\epsilon} 0$, and the CAT(-1) inequality, we have
\[
\lb z_i|y_i\rb_\zero \gtrsim_{\plus,\epsilon} -\log(\w\epsilon),
\]
and thus $\wbar\Dist(z_i,y_i)\lesssim_{\times,\epsilon} \w\epsilon$. Using the triangle inequality together with the assumption $\wbar\Dist(y_1,y_2)\leq\delta$, we have
\[
\wbar\Dist(z_1,z_2) \lesssim_{\times,\epsilon} \max(\delta,\w\epsilon).
\]
Setting $\w\epsilon$ equal to $\epsilon$ divided by the implied constant, and decreasing $\delta$ if necessary, completes the proof.
\end{itemize}

%
%
%
%
\end{proof}


\begin{notation}
\label{notationwpi}
If the map $\pi:[t,s]\to X$ is an isometric embedding, then the map $\w\pi:[-\infty,+\infty]\to X$ is defined by the equation
\[
\w\pi(r) = \pi(t\vee r\wedge s).
\]
\end{notation}

\begin{corollary}
\label{corollaryxyz}
If $\epsilon$, $\delta$, and $\Delta(x,y_1,y_2)$ are as in Lemma \ref{lemmaxyz}, and if $\pi_1:[t,s_1]\to \geo x{y_1}$ and $\pi_2:[t,s_2]\to\geo x{y_2}$ are isometric embeddings, then
\[
\wbar\Dist(\w\pi_1(r),\w\pi_2(r)) \lesssim_\times \epsilon \all r\in [-\infty,+\infty].
\]
\end{corollary}
\begin{proof}
If $r\leq t$, then $\w\pi_1(r) = x = \w\pi_2(r)$. If $t\leq r \leq \min_{i = 1}^2 s_i$, then $\w\pi_i(r) = \geo x{y_i}_{r - t}$, allowing us to apply Lemma \ref{lemmaxyz} directly. Finally, suppose $r\geq r_0 := \min_{i = 1}^2 s_i$. Without loss of generality suppose that $s_1 \leq s_2$, so that $r_0 = s_1$. Applying the previous case to $r_0$, we have
\[
\wbar\Dist(y_1,w_2) \leq \epsilon,
\]
where $w_2 = \pi_2(s_1)$. Now $\w\pi_1(r) = y_1$, and $\w\pi_2(r)\in\geo{w_2}{y_2}$, so Lemma \ref{lemmageodesicdiameter} completes the proof.
\end{proof}

\begin{lemma}
\label{lemmageodesicconvergence}
Suppose that $(x_n)_1^\infty$ and $(y_n)_1^\infty$ are sequences in $X$ which converge to points $x_n\to x\in\bord X$ and $y_n\to y\in\bord X$, with $x\neq y$. Then there exists a geodesic $\geo xy$ connecting $x$ and $y$ such that $\geo{x_n}{y_n} \to \geo xy$ in the Hausdorff metric. If $x = y$, then $\geo{x_n}{y_n}\to \{x\}$ in the Hausdorff metric.
\end{lemma}
\begin{proof}
We observe first that if $x = y$, then the conclusion follows immediately from Lemma \ref{lemmageodesicdiameter}. Thus we assume in what follows that $x\neq y$.

For any pair $p,q\in X$, we define the \emph{standard parameterization} of the geodesic $\geo pq$ to be the unique isometry $\pi:[-\lb \zero|q\rb_p , \lb \zero|p\rb_q] \to \geo pq$ sending $-\lb \zero|q\rb_p$ to $p$ and $\lb \zero|p\rb_q$ to $q$. For each $n$ let $\pi_n:[t_n,s_n]\to \geo{x_n}{y_n}$ be the standard parameterization, and for each $m,n\in\N$ let $\pi_{m,n}:[t_{m,n},s_{m,n}]\to \geo{x_n}{y_m}$ be the standard parameterization. Let $\w\pi_n:[-\infty,+\infty]\to \geo{x_n}{y_n}$ and $\w\pi_{m,n}:[-\infty,+\infty]\to \geo{x_n}{y_m}$ be as in Notation \ref{notationwpi}. Note that
\[
t_n - t_{m,n} = \lb \zero|y_m\rb_{x_n} - \lb \zero|y_n\rb_{x_n} = \lb x_n|y_n\rb_\zero - \lb x_n|y_m\rb_\zero \tendsto{m,n} \lb x|y\rb_\zero - \lb x|y\rb_\zero = 0.
\]
(We have $\lb x|y\rb_\zero < \infty$ since $x\neq y$.) Thus
\[
\wbar\Dist(\w\pi_n(r), \w\pi_n(r - t_n + t_{m,n})) \leq \dist(\w\pi_n(r), \w\pi_n(r - t_n + t_{m,n})) \leq |t_n - t_{m,n}| \tendsto n 0.
\]
Here and below, the limit converges uniformly for $r\in[-\infty,+\infty]$. On the other hand, Corollary \ref{corollaryxyz} implies that
\[
\wbar\Dist(\w\pi_n(r - t_n + t_{m,n}) , \w\pi_{m,n}(r)) \tendsto{m,n} 0,
\]
so the triangle inequality gives
\[
\wbar\Dist(\w\pi_n(r) , \w\pi_{m,n}(r)) \tendsto{m,n} 0.
\]
A similar argument shows that
\[
\wbar\Dist(\w\pi_{m,n}(r) , \w\pi_m(r)) \tendsto{m,n} 0,
\]
so the triangle inequality gives
\[
\wbar\Dist(\w\pi_n(r) , \w\pi_m(r)) \tendsto{m,n} 0,
\]
i.e. the sequence of functions $(\w\pi_n)_1^\infty$ is uniformly Cauchy. Since $(\bord X,\wbar\Dist)$ is complete, they converge uniformly to a function $\w\pi:[-\infty,+\infty]\to X$.

Clearly, $\geo{x_n}{y_n} = \w\pi_n([-\infty,+\infty]) \to \w\pi([-\infty,+\infty])$ in the Hausdorff metric. We claim that $\w\pi([-\infty,+\infty])$ is a geodesic connecting $x$ and $y$. Indeed,
\[
t_n \tendsto n t := \lb x|y\rb_\zero - \dox x \text{ and } s_n \tendsto n s := \dox y - \lb x|y\rb_\zero.
\]
For all $t < r_1 < r_2 < s$, we have $t_n < r_1 < r_2 < s_n$ for all sufficiently large $n$, which implies that
\[
\dist(\w\pi(r_1) , \w\pi(r_2)) = \lim_{n\to\infty} \dist(\w\pi_n(r_1) , \w\pi_n(r_2)) = \lim_{n\to\infty} (r_2 - r_1) = r_2 - r_1,
\]
i.e. $\w\pi\given(t,s)$ is an isometric embedding. Since $\w\pi$ is continuous (being the uniform limit of continuous functions), $\pi := \w\pi\given[t,s]$ is also an isometric embedding. A similar argument shows that $\w\pi(r) = \pi(t)$ for all $r\leq t$, and $\w\pi(r) = \pi(s)$ for all $r\geq s$; thus $\w\pi([-\infty,+\infty]) = \pi([t,s])$ is a geodesic. To complete the proof, we must show that $\pi(t) = x$ and $\pi(s) = y$. Indeed,
\[
\pi(t) = \w\pi(-\infty) = \lim_{n\to\infty} \w\pi_n(-\infty) = \lim_{n\to\infty} x_n = x,
\]
and a similar argument shows that $\pi(s) = y$. Thus the geodesic $\pi([t,s])$ connects $x$ and $y$.
\end{proof}

Using Lemma \ref{lemmageodesicconvergence}, we prove Proposition \ref{propositiongeodesicconvergence}.

\begin{proof}[Proof of Proposition \ref{propositiongeodesicconvergence}]~
\begin{itemize}
\item[(i)] Given distinct points $x,y\in\bord X$, we may find sequences $X\ni x_n\to x$ and $X\ni y_n\to y$. Applying Lemma \ref{lemmageodesicconvergence} proves the existence of a geodesic connecting $x$ and $y$. To show uniqueness, suppose that $\geo xy_1$ and $\geo xy_2$ are two geodesics connecting $x$ and $y$. Fix sequences $\geo xy_1\ni x_n^{(1)}\to x$, $\geo xy_2\ni x_n^{(2)}\to x$, $\geo xy_1\ni y_n^{(1)}\to y$,and $\geo xy_2\ni y_n^{(2)}\to y$. By considering the intertwined sequences 
\[
x_1^{(1)},x_1^{(2)},x_2^{(1)},x_2^{(2)},\ldots 
\]
and 
\[
y_1^{(1)},y_1^{(2)},y_2^{(1)},y_2^{(2)},\ldots
\] 
Lemma \ref{lemmageodesicconvergence} shows that both sequences $(\geo{x_n^{(1)}}{y_n^{(1)}})_1^\infty$ and $(\geo{x_n^{(2)}}{y_n^{(2)}})_1^\infty$ converge in the Hausdorff metric to a common geodesic $\geo xy$. But clearly the former tend to $\geo xy_1$, and the latter tend to $\geo xy_2$; we must have $\geo xy_1 = \geo xy_2$.

\item[(ii)] Suppose that $\bord X\ni x_n\to x$ and $\bord X\ni y_n\to y$. For each $n$, choose $\what x_n,\what y_n\in \geo{x_n}{y_n}\cap X$ such that $\wbar\Dist(\what x_n,x_n) , \wbar\Dist(\what y_n,y_n) \leq 1/n$. Then $\what x_n\to x$ and $\what y_n\to y$, so by Lemma \ref{lemmageodesicconvergence} we have $\geo{\what x_n}{\what y_n}\to \geo xy$ in the Hausdorff metric, or $\geo{\what x_n}{\what y_n} \to \{x\}$ if $x = y$. To complete the proof it suffices to show that the Hausdorff distance between $\geo{x_n}{y_n}$ and $\geo{\what x_n}{\what y_n}$ tends to zero as $n$ tends to infinity. Indeed, $\geo{\what x_n}{\what y_n} \subset \geo{x_n}{y_n}$, and for each $z\in \geo{x_n}{y_n}$, either $z\in \geo{x_n}{\what x_n}$, $z\in \geo{\what x_n}{\what y_n}$, or $z\in \geo{\what y_n}{y_n}$. In the first case, Lemma \ref{lemmageodesicdiameter} shows that $\wbar\Dist(z,\geo{\what x_n}{\what y_n}) \leq \wbar\Dist(z,\what x_n) \lesssim_\times \wbar\Dist(x_n,\what x_n) \leq 1/n \to 0$; the third case is treated similarly.

\end{itemize}
\end{proof}

Having completed the proof of Proposition \ref{propositiongeodesicconvergence}, in the remainder of this section we prove that a version of the CAT(-1) equality holds for ideal triangles.

\begin{definition}
\label{definitiontriangles}
A \emph{geodesic triangle} $\Delta = \Delta(x,y,z)$ consists of three distinct points $x,y,z\in\bord X$ together with the geodesics $\geo xy$, $\geo yz$, and $\geo zx$.

A geodesic triangle $\wbar{\Delta} = \Delta(\wbar x,\wbar y,\wbar z)$ is called a \emph{comparison triangle} for $\Delta$ if
\[
\lb \wbar x|\wbar y\rb_{\wbar z} = \lb x|y\rb_z, \text{ etc.}
\]
For any point $p\in\geo xy$, its \emph{comparison point} is defined to be the unique point $\wbar p\in\geo {\wbar x}{\wbar y}$ such that
\[
\lb x|z\rb_p - \lb y|z\rb_p = \lb \wbar x|\wbar z\rb_{\wbar p} - \lb \wbar y|\wbar z\rb_{\wbar p}.
\]
We say that the geodesic triangle $\Delta$ \emph{satisfies the CAT(-1) inequality} if for all points $p,q\in\Delta$ and for any comparison points $\wbar p,\wbar q\in\wbar\Delta$, we have $\dist(p,q) \leq \dist(\wbar p,\wbar q)$.
\end{definition}

It should be checked that these definitions are consistent with those given in Section \ref{subsectionCAT}.

\begin{proposition}
\label{propositionidealCAT}
Any geodesic triangle (including ideal triangles) satisfies the CAT(-1) inequality.
\end{proposition}
\begin{proof}
Let $\Delta = \Delta(x,y,z)$ be a geodesic triangle, and fix $p,q\in \Delta$. Choose sequences $x_n\to x$, $y_n\to y$, and $z_n\to z$. By Proposition \ref{propositiongeodesicconvergence}, we have $\Delta_n = \Delta(x_n,y_n,z_n) \to \Delta$ in the Hausdorff metric, so we may choose $p_n,q_n\in\Delta_n$ so that $p_n\to p$, $q_n\to q$. For each $n$, let $\wbar\Delta_n = \Delta(\wbar x_n,\wbar y_n,\wbar z_n)$ be a comparison triangle for $\Delta_n$. Without loss of generality, we may assume that
\begin{equation}
\label{zerocenter}
\zero\in\geo{\wbar x_n}{\wbar y_n} \text{ and } \lb \wbar x_n|\wbar z_n \rb_\zero = \lb \wbar y_n | \wbar z_n \rb_\zero \asymp_\plus 0.
\end{equation}
By extracting a convergent subsequence, we may without loss of generality assume that $\wbar x_n\to \wbar x$, $\wbar y_n\to \wbar y$, and $\wbar z_n\to \wbar z$ for some points $\wbar x,\wbar y,\wbar z\in\bord \H^2$. By \eqref{zerocenter}, the points $\wbar x,\wbar y,\wbar z$ are distinct. Thus $\wbar\Delta = \Delta(\wbar x,\wbar y,\wbar z)$ is a geodesic triangle, and is in fact a comparison triangle for $\Delta$. If $\wbar p,\wbar q$ are comparison points for $p,q$, then $\wbar p_n\to \wbar p$ and $\wbar q_n\to q$. It follows that
\[
\dist(p,q) = \lim_{n\to\infty} \dist(p_n,q_n) \leq \lim_{n\to\infty} \dist(\wbar p_n,\wbar q_n) = \dist(\wbar p,\wbar q).
\]
\end{proof}

\bigskip
\section{The geometry of shadows}\label{subsectionshadows}

\subsection{Shadows in regularly geodesic hyperbolic metric spaces}
Suppose that $X$ is regularly geodesic. For each $z\in X$ we consider the relation $\pi_z\subset X\times\del X$ defined by
\[
(x,\xi)\in \pi_z \Leftrightarrow x\in \geo z\xi
\]
(see Definition \ref{definitiongeodesic} for the definition of $\geo z\xi$). Note that if $X$ is an algebraic hyperbolic space, then the relation $\pi_z$ is a function when restricted to $X\butnot\{z\}$; in particular, for $\xx\in \B = \B_\F^\alpha$ with $\xx\neq\0$ we have
\[
\pi_\0(\xx) = \frac{\xx}{\|\xx\|}\cdot
\]
However, in general the relation $\pi_z$ is not necessarily a function; $\R$-trees provide a good counterexample. The reason is that in an $\R$-tree, there may be multiple ways to extend a geodesic segment to a geodesic ray.

For any set $S$, we define its \emph{shadow} with respect to the light source $z$ to be the set
\[
\pi_z(S) := \{\xi\in\del X:\exists x\in S \;\; (x,\xi)\in \pi_z\}.
\]

\begin{figure}
\begin{center}
\begin{tikzpicture}[line cap=round,line join=round,>=triangle 45,x=1.0cm,y=1.0cm]
\clip(-3.85,-3.03) rectangle (3.99,3.1);
\draw(0,0) circle (3cm);
\draw[line width=0.81pt](1.7,1.4) circle (0.44cm);
\draw (0,0)-- (1.89,2.33);
\draw (0,0)-- (2.65,1.4);
\draw [shift={(0,0)},line width=1.2pt] plot[domain=0.49:0.89,variable=\t]({1*3*cos(\t r)+0*3*sin(\t r)},{0*3*cos(\t r)+1*3*sin(\t r)});
\draw [shift={(1.83,2.51)}] plot[domain=4.19:4.78,variable=\t]({1*0.96*cos(\t r)+0*0.96*sin(\t r)},{0*0.96*cos(\t r)+1*0.96*sin(\t r)});
\begin{scriptsize}
\fill [color=black] (0,0) circle (1pt);
\draw[color=black] (-0.1,-0.16) node {$z$};
\draw[color=black] (0.72,1.66) node {$B(x,\sigma)$};
\fill [color=black] (1.89,1.55) circle (1pt);
\draw[color=black] (2,1.44) node {$x$};
\draw[color=black] (1.59,1.46) node {$\sigma$};
\draw[color=black] (3.2,1.98) node {$\pi_z\big(B(x,\sigma)\big)$};
\end{scriptsize}
\end{tikzpicture}
\caption[Shadows in regularly geodesic hyperbolic metric spaces]{The set $\pi_z(B(x,\sigma))$. Although this set is not equal to $\Shad_z(x,\sigma)$, they are approximately the same in regularly geodesic spaces by Corollary \ref{corollaryshadowsrelation}. In our drawings, we will draw the set $\pi_z(B(x,\sigma))$ to indicate the set $\Shad_z(x,\sigma)$ (since the latter is hard to draw).}
\label{figureshadow}
\end{center}
\end{figure}
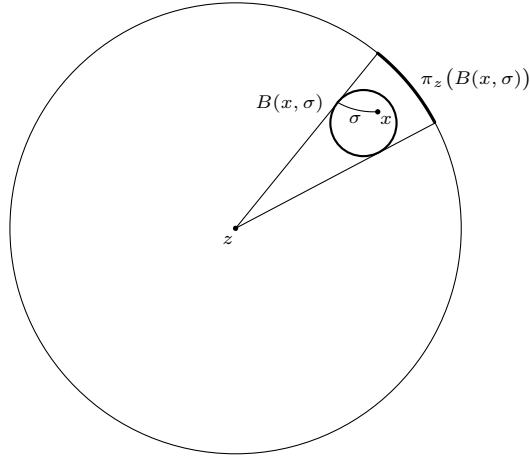

\subsection{Shadows in hyperbolic metric spaces}
In regularly geodesic hyperbolic metric spaces, it is particularly useful to consider $\pi_z(B(x,\sigma))$ where $x\in X$ and $\sigma > 0$. We would like to have an analogue for this set in the Gromov hyperbolic setting.

\begin{definition}
\label{definitionshadow}
For each $\sigma > 0$ and $x,z\in X$, let
\[
\Shad_z(x,\sigma) = \{\eta\in \del X:\lb z|\eta\rb_x\leq\sigma\}.
\]
We say that $\Shad_z(x,\sigma)$ is the \emph{shadow cast by $x$ from the light source $z$, with parameter $\sigma$}. For shorthand we will write $\Shad(x,\sigma) = \Shad_\zero(x,\sigma)$.
\end{definition}
The relation between $\pi_z(B(x,\sigma))$ and $\Shad_z(x,\sigma)$ in the case where $X$ is a regularly geodesic hyperbolic metric space will be made explicit in Corollary \ref{corollaryshadowsrelation} below.

Let us establish up front some geometric properties of shadows.

\begin{observation}
\label{observationshadowsclosed}
For each $x,z\in X$ and $\sigma > 0$ the set $\Shad_z(x,\sigma)$ is closed.
\end{observation}
\begin{proof}
This follows directly from Lemma \ref{lemmalowersemicontinuous}.
\end{proof}

\begin{observation}\label{observationmetricshadow}
If $\eta\in \Shad_z(x,\sigma)$, then
\[
\lb x|\eta\rb_z \asymp_{\plus,\sigma}
 \dist(z,x).
\]
\end{observation}
\begin{proof}
Follows directly from (b) of Proposition \ref{propositionbasicidentities} together with the definition of $\Shad_z(x,\sigma)$.
\end{proof}

\begin{lemma}[Intersecting Shadows Lemma]
\label{lemmatau}
For each $\sigma > 0$, there exists $\tau = \tau_\sigma > 0$ such that for all $x,y,z\in X$ satisfying $\dist(z,y)\geq \dist(z,x)$ and $\Shad_z(x,\sigma)\cap \Shad_z(y,\sigma)\neq\emptyset$, we have
\begin{equation}
\label{shadcontainment}
\Shad_z(y,\sigma) \subset\Shad_z(x,\tau)
\end{equation}
and
\begin{equation}
\label{distbusemann}
\dist(x,y) \asymp_{\plus,\sigma} \dist(z,y) - \dist(z,x).
\end{equation}
\end{lemma}

\begin{figure}
\begin{center}
\definecolor{ttqqqq}{rgb}{0.2,0,0}
\definecolor{ffttww}{rgb}{1,0.1,0.1}
\begin{tikzpicture}[line cap=round,line join=round,>=triangle 45,x=1.0cm,y=1.0cm]
\clip(-4.88,-4.1) rectangle (5.4,4.1);
\draw(0,0) circle (4.03cm);
\draw(1.84,2.21) circle (0.76cm);
\draw(2.92,1.93) circle (0.25cm);
\draw [shift={(0,0)},line width=2pt,color=ffttww] plot[domain=0.61:1.14,variable=\t]({1*4.03*cos(\t r)+0*4.03*sin(\t r)},{0*4.03*cos(\t r)+1*4.03*sin(\t r)});
\draw [shift={(0,0)},line width=1.2pt,color=yellow] plot[domain=0.66:0.61,variable=\t]({1*4.03*cos(\t r)+0*4.03*sin(\t r)},{0*4.03*cos(\t r)+1*4.03*sin(\t r)});
\draw [shift={(0,0)},dash pattern=on 1pt off 2pt,line width=1pt,color=blue] plot[domain=0.51:0.66,variable=\t]({1*4.03*cos(\t r)+0*4.03*sin(\t r)},{0*4.03*cos(\t r)+1*4.03*sin(\t r)});
\draw (1.15,2.53)-- (1.67,3.67);
\draw (2.28,1.59)-- (3.31,2.31);
\draw (2.76,2.13)-- (3.19,2.47);
\draw (3.04,1.71)-- (3.51,1.98);
\begin{scriptsize}
\fill [color=black] (0,0) circle (1pt);
\draw[color=black] (-0.08,-0.34) node {$z$};
\fill [color=black] (2.07,2.49) circle (1pt);
\draw[color=black] (0.5,2.13) node {$B(x,\sigma)$};
\fill [color=black] (3.01,1.99) circle (1pt);
\draw[color=black] (3.1,1.4) node {$B(y,\sigma)$};
\draw[color=ffttww] (3.22,3.36) node {$\Shad_z(x,\sigma)$};
\draw[color=blue] (4.38,2.1) node {$\Shad_z(y,\sigma)$};
\end{scriptsize}
\end{tikzpicture}
\caption[The Intersecting Shadows Lemma]{In this figure, $\dist(z,y)\geq \dist(z,x)$ and $\Shad_z(x,\sigma)\cap \Shad_z(y,\sigma)\neq\emptyset$. The Intersecting Shadows Lemma (Lemma \ref{lemmatau}) provides a $\tau_\sigma > 0$ such that the shadow cast from $z$ about $B(x,\tau_\sigma)$ will capture $\Shad_z(y,\sigma)$.}
\label{figureintersectingshadows}
\end{center}
\end{figure}
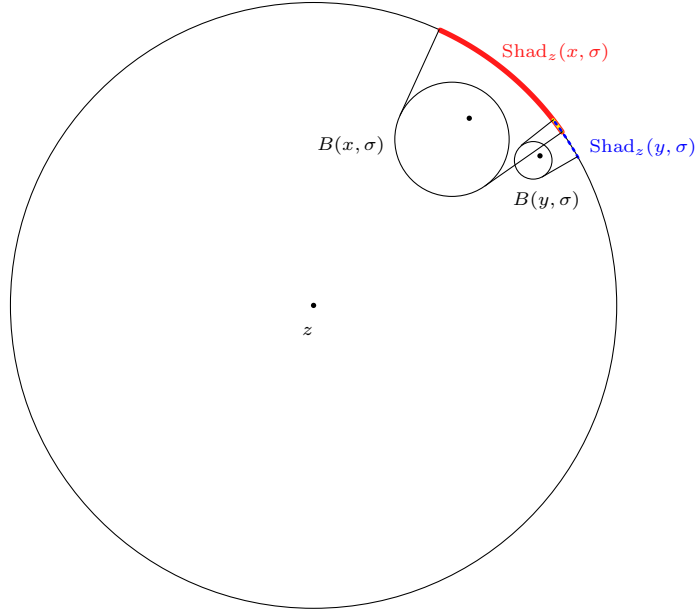

\begin{proof}
Fix $\eta\in\Shad_z(x,\sigma)\cap \Shad_z(y,\sigma)$, so that by Observation \ref{observationmetricshadow}
\[
\lb x|\eta\rb_z \asymp_{\plus,\sigma} \dist(z,x)\ \text{ and } \
\lb y|\eta\rb_z \asymp_{\plus,\sigma} \dist(z,y) \geq \dist(z,x).
\]
Gromov's inequality along with (c) of Proposition \ref{propositionbasicidentities} then gives
\begin{equation}\label{distbusemann2}
\lb x|y\rb_z \asymp_{\plus,\sigma} \dist(z,x).
\end{equation}
Rearranging yields \eqref{distbusemann}. In order to show \eqref{shadcontainment}, fix $\xi\in\Shad_z(y,\sigma)$, so that $\lb y|\xi\rb_z \asymp_{\plus,\sigma} \dist(z,y)\geq \dist(z,x)$. Gromov's inequality and \eqref{distbusemann2} then give
\[
\lb x|\xi\rb_z \asymp_{\plus,\sigma} \dist(z,x),
\]
i.e. $\xi\in\Shad_z(x,\tau)$ for some $\tau > 0$ sufficiently large (depending on $\sigma$).
\end{proof}

\begin{corollary}
\label{corollaryshadowsrelation}
Suppose that $X$ is regularly geodesic. For every $\sigma > 0$, there exists $\tau = \tau_\sigma > 0$ such that for any $x,z\in X$ we have
\begin{equation}
\label{shadowsrelation}
\pi_z(B(x,\sigma)) \subset \Shad_z(x,\sigma) \subset \pi_z(B(x,\tau)).
\end{equation}
\end{corollary}
\begin{proof}
Suppose $\xi\in\pi_z(B(x,\sigma))$. Then there exists a point $y\in B(x,\sigma)\cap \geo z\xi$. By (d) of Proposition \ref{propositionbasicidentities}
\begin{align*}
\lb z|\xi\rb_x \leq \lb z|\xi\rb_y + \dist(x,y)
\leq \lb z|\xi\rb_y + \sigma = \sigma,
\end{align*}
i.e. $\xi\in\Shad_z(x,\sigma)$. This demonstrates the first inclusion of \eqref{shadowsrelation}. On the other hand, suppose that $\xi\in\Shad_z(x,\sigma)$. Let $y\in\geo z\xi$ be the unique point so that $\dist(z,y) = \dist(z,x)$. Clearly $\xi\in\Shad_z(y,\sigma)$, so $\Shad_z(x,\sigma)\cap\Shad_z(y,\sigma)\neq\emptyset$; by the Intersecting Shadows Lemma \ref{lemmatau} we have
\[
\dist(x,y) \asymp_{\plus,\sigma} \busemann_\xi(y,x) = 0
\]
i.e. $\dist(x,y)\leq\tau$ for some $\tau = \tau_\sigma > 0$ depending only on $\sigma$. Then $y\in B(x,\tau)\cap\geo z\xi$, which implies that $\xi = \pi_z(y)\in \pi_z(B(x,\tau))$. This finishes the proof.
\end{proof}

\begin{lemma}[Bounded Distortion Lemma]
\label{lemmaboundeddistortion}
Fix $\sigma > 0$. Then for every $g\in\Isom(X)$ and for every $y\in\Shad_{g^{-1}(\zero)}(\zero,\sigma)$ we have
\begin{equation}
\label{boundeddistortion1}
\overline g'(y) \asymp_{\times,\sigma} b^{-\dogo g}.
\end{equation}
Moreover, for every $y_1,y_2\in \Shad_{g^{-1}(\zero)}(\zero,\sigma)$, we have
\begin{equation}
\label{boundeddistortion2}
\frac{\Dist (g(y_1),g(y_2) )}{\Dist(y_1,y_2)} \asymp_{\times,\sigma} b^{-\dogo g}.
\end{equation}
\end{lemma}
\begin{proof}
We have $\overline g'(y) \asymp_\times b^{\busemann_y(\zero,g^{-1}(\zero))}
\asymp_\times b^{2\lb g^{-1}(\zero)|y\rb_\zero - \dogo g}
\asymp_{\times,\sigma} b^{-\dogo g}$, giving \eqref{boundeddistortion1}.
Now \eqref{boundeddistortion2} follows from \eqref{boundeddistortion1} and the geometric mean value theorem (Proposition \ref{propositionGMVT}).
\end{proof}

\begin{lemma}[Big Shadows Lemma]
\label{lemmabigshadow}
For every $\epsilon > 0$, for every $\sigma > 0$ sufficiently large (depending on $\epsilon$), and for every $z\in X$, we have
\begin{equation}
\label{bigshadow}
\Diam(\del X\butnot\Shad_z(\zero,\sigma)) \leq \epsilon.
\end{equation}
\end{lemma}

\begin{figure}
\begin{center}
\begin{tikzpicture}[line cap=round,line join=round,>=triangle 45,scale=0.85]
\clip(-5.4,-4.4) rectangle (5.4,4.4);
\draw(0,0) circle (4.12cm);
\draw(0,0) circle (3.07cm);
\draw [shift={(-3.49,-2.35)}] plot[domain=-0.28:1.95,variable=\t]({1*1.14*cos(\t r)+0*1.14*sin(\t r)},{0*1.14*cos(\t r)+1*1.14*sin(\t r)});
\draw [shift={(-1.98,-3.76)}] plot[domain=-0.24:1.94,variable=\t]({1*1.18*cos(\t r)+0*1.18*sin(\t r)},{0*1.18*cos(\t r)+1*1.18*sin(\t r)});
\draw [shift={(0,0)},line width=1.2pt] plot[domain=3.46:4.51,variable=\t]({1*4.12*cos(\t r)+0*4.12*sin(\t r)},{0*4.12*cos(\t r)+1*4.12*sin(\t r)});
\draw (0,0)-- (2.56,1.7);
\begin{scriptsize}
\fill [color=black] (0,0) circle (1.0pt);
\draw[color=black] (-0.23,0.12) node {$\zero$};
\fill [color=black] (-2.4,-2.66) circle (1.0pt);
\draw[color=black] (-2.19,-2.84) node {$z$};
\draw[color=black] (-4.09,-3.21) node {$\del X\butnot\Shad_z(\zero,\sigma)$};
\draw[color=black] (1.97,1.08) node {$\sigma$};
\end{scriptsize}
\end{tikzpicture}
\caption[The Big Shadows Lemma]{The Big Shadows Lemma \ref{lemmabigshadow} tells us that for any $\epsilon > 0$, we may choose $\sigma > 0$ sufficiently large so that $\Diam(\del X\butnot\Shad_z(\zero,\sigma)) \leq \epsilon$ for every $z\in X$.}
\label{figurebigshadows}
\end{center}
\end{figure}

\begin{proof}
If $\xi, \eta\in\del X \setminus \Shad_z(\zero, \sigma)$, then $\lb z | \xi \rb_\zero > \sigma$ and $\lb z | \eta \rb_\zero > \sigma$. Thus by Gromov's inequality we have
\[
\lb \xi | \eta \rb_\zero \gtrsim_\plus \sigma.
\]
Exponentiating gives $\Dist(\xi,\eta) \lesssim_\times b^{-\sigma}$. Thus
\[
\Diam (\del X \setminus\Shad_z (\zero,\sigma)) \lesssim_\times b^{-\sigma} \tendsto\sigma 0,
\]
and the convergence is uniform in $z$.
\end{proof}

\begin{figure}
\begin{center}
\begin{tikzpicture}[line cap=round,line join=round,>=triangle 45,x=1.0cm,y=1.0cm]
\clip(-3.54,-3.1) rectangle (3.943,3.38);
\draw(0.0,0.0) circle (3.0cm);
\draw(1.5,1.5) circle (0.75cm);
\draw (0.0,-0.0)-- (1.2343134832984426,2.7343134832984424);
\draw (0.0,-0.0)-- (2.7343134832984424,1.2343134832984426);
\draw (0.0,-0.0)-- (1.7415599609449397,1.7700977248889884);
\draw [shift={(1.317194030436857,2.707806971266285)}] plot[domain=4.204270499810774:5.13736643041908,variable=\t]({1.0*1.0292643361729137*cos(\t r)+-0.0*1.0292643361729137*sin(\t r)},{0.0*1.0292643361729137*cos(\t r)+1.0*1.0292643361729137*sin(\t r)});
\draw [shift={(0.0,0.0)},line width=1.2pt] plot[domain=0.4240310394907405:1.146765287304156,variable=\t]({1.0*3.0*cos(\t r)+-0.0*3.0*sin(\t r)},{0.0*3.0*cos(\t r)+1.0*3.0*sin(\t r)});
\begin{scriptsize}
\draw [fill=black] (0.0,-0.0) circle (.75pt);
\draw[color=black] (-0.19797292010819517,-0.09284349186632893) node {$z$};
\draw [fill=black] (1.7415599609449397,1.7700977248889884) circle (.75pt);
\draw[color=black] (1.8727829942209163,1.50) node {$g(\zero)$};

\draw[color=black] (1.5727829942209163,0.71) node {$\Big\uparrow$};
\draw[color=black] (1.9727829942209163,0.31) node {$B(g(\zero),\sigma)$};

\draw[color=black] (1.2421411016356007,1.814123570799388) node {$\sigma$};
\draw[color=black] (3.0115304609745355,2.3508938419017653) node {$\Shad_z(g(\zero),\sigma)$};
\end{scriptsize}
\end{tikzpicture}
\caption[The Diameter of Shadows Lemma]{The Diameter of Shadows Lemma \ref{lemmadiameterasymptotic} says that the diameter of $\Shad(g(\zero),\sigma)$ is coarsely asymptotic to $b^{-\dogo g}$.}
\end{center}
\end{figure}

\begin{lemma}[Diameter of Shadows Lemma]
\label{lemmadiameterasymptotic}
For all $\sigma > 0$ sufficiently large, we have for all $g\in\Isom(X)$ and for all $z\in X$
\begin{equation}
\label{diameterasymptotic}
\Diam_z(\Shad_z(g(\zero),\sigma)) \lesssim_{\times,\sigma} b^{-\dist(z,g(\zero))},
\end{equation}
with $\asymp$ if $\#(\del X)\geq 3$. Moreover, for every $C > 0$ there exists $\sigma > 0$ such that
\begin{equation}
\label{ShadcontainsB}
B_z(x,C e^{-\dist(z,x)}) \subset \Shad_z(x,\sigma) \all x,z\in X.
\end{equation}
\end{lemma}
\begin{proof}
Let $x = g(\zero)$. For any $\xi,\eta\in \Shad_z(x,\sigma)$, we have
\begin{align*}
\Dist_z(\xi,\eta) \asymp_\times b^{-\lb\xi|\eta\rb_z}
\lesssim_{\times\phantom{,\sigma}} b^{-\min\left(\lb x|\xi\rb_z,\lb x|\eta\rb_z\right)}
\lesssim_{\times,\sigma} b^{-\dist(z,x)}
\end{align*}
which demonstrates \eqref{diameterasymptotic}.

Now let us prove the converse of \eqref{diameterasymptotic}, assuming $\#(\del X)\geq 3$. Fix $\xi_1,\xi_2,\xi_3\in\del X$, let $\epsilon = \min_{i\neq j}\Dist(\xi_i,\xi_j)/2$, and fix $\sigma > 0$ large enough so that \eqref{bigshadow} holds for every $z\in X$. By \eqref{bigshadow} we have
\[
\Diam(\del X\butnot \Shad_{g^{-1}(z)}(\zero,\sigma)) \leq \epsilon,
\]
and thus
\[
\#\big\{i = 1,2,3: \xi_i\in \Shad_{g^{-1}(z)}(\zero,\sigma)\big\}\geq 2.
\]
Without loss of generality suppose that $\xi_1,\xi_2\in \Shad_{g^{-1}(z)}(\zero,\sigma)$. By applying $g$, we have $g(\xi_1),g(\xi_2)\in \Shad_z(x,\sigma)$. Then
\begin{align*}
\Diam_z(\Shad_z(x,\sigma))
&\geq_\pt \Dist_z(g(\xi_1),g(\xi_2))\\
&\asymp_\times b^{-\lb g(\xi_1)|g(\xi_2)\rb_z}
= b^{-\lb \xi_1|\xi_2\rb_{g^{-1}(z)}}\\
&\gtrsim_\times b^{-\lb\xi_1|\xi_2\rb_\zero}b^{-\dox{g^{-1}(z)}}
\asymp_{\times,\xi_1,\xi_2}b^{-\dist(z,x)}.
\end{align*}
Finally, given $y\in B_z(x,C b^{-\dist(z,x)})$, we have
\[
\lb x|y\rb_z \gtrsim_\plus -\log_b(C b^{-\dox x}) \asymp_\plus \dist(z,x)
\]
and thus $\lb z|y\rb_x \asymp_\plus 0$, demonstrating \eqref{ShadcontainsB}.
\end{proof}

\bigskip
\section{Generalized polar coordinates} \label{subsectionpolar}

Suppose that $X = \E = \E^\alpha$ is the half-space model of a real hyperbolic space. Fix a point $\xx\in\E$, and consider the numbers $\|\xx\|$ and $\ang(\xx) := \cos^{-1}(x_1 / \|\xx\|)$, i.e. the radial and unsigned angular coordinates of $\xx$. (The angular coordinate is computed with respect to the ray $\{(t,\0) : t\in\Rplus\}$; cf. Figure \ref{figurepolarcoordinates}.) These ``polar coordinates'' of $\xx$ do not completely determine $\xx$, but they are enough to compute certain important quantities depending on $\xx$, e.g. $\dist_\E(\zero,\xx)$, $\busemann_\infty(\zero,\xx)$, and $\busemann_\0(\zero,\xx)$. (We omit the details.) In this section we consider a generalization, in a loose sense, of these coordinates to an arbitrary hyperbolic metric space.

\begin{figure}
\begin{center}
\begin{tikzpicture}[line cap=round,line join=round,>=triangle 45,x=1.0cm,y=1.0cm]
\clip(-5.2,-0.33) rectangle (5.38,5.43);
\draw [shift={(0,0)},color=black,fill=black,fill opacity=0.1] (0,0) -- (42.62:0.79) arc (42.62:90:0.79) -- cycle;
\draw (-4,0)-- (4,0);
\draw (0,0)-- (0,5);
\draw (0,0)-- (2.76,2.54);
\draw [dash pattern=on 2pt off 2pt] (2.76,2.54)-- (2.78,0);
\begin{scriptsize}
\fill [color=black] (0,0) circle (1pt);
\draw[color=black] (-0.1,-0.2) node {$\0$};
\fill [color=black] (2.76,2.54) circle (1pt);
\draw[color=black] (3,2.6) node {$\xx$};
\draw[color=black] (1.44,1.82) node {$\|\xx\|$};
\draw[color=black] (3.06,1.09) node {$x_1$};
\draw[color=black] (0.43,0.9) node {$\ang(\xx)$};
\end{scriptsize}
\end{tikzpicture}
\caption[Polar coordinates in the half-space model]{The quantities $\|\xx\|$ and $\ang(\xx)$ can be interpreted as ``polar coordinates'' of $\xx$.}
\label{figurepolarcoordinates}
\end{center}
\end{figure}

Let us note that the isometries of $\E$ which preserve the polar coordinate functions defined above are exactly those of the form $\what T$ where $T\in\O(\EE)$. Equivalently, these are the members of $\Isom(\E)$ which preserve $\0$, $\zero = (1,\0)$, and $\infty$. This suggests that our ``coordinate system'' is fixed by choosing a point in $\E$ and two distinct points in $\del\E$.

We now return to the general case of \6\ref{standingassumptions2}. Fix two distinct points $\xi_1,\xi_2\in\del X$.
\begin{definition}
\label{definitionpolarcoordinates}
The \emph{generalized polar coordinate functions} are the functions $r = r_{\xi_1,\xi_2,\zero}$ and $\theta = \theta_{\xi_1,\xi_2,\zero}:X\to\R$ defined by
\begin{align*}
r(x) &= \frac12\left[\busemann_{\xi_1}(x,\zero) - \busemann_{\xi_2}(x,\zero)\right]\\
\theta(x) &= \frac12\left[\busemann_{\xi_1}(x,\zero) + \busemann_{\xi_2}(x,\zero)\right] \asymp_\plus \lb \xi_1|\xi_2\rb_x - \lb \xi_1|\xi_2\rb_\zero.
\end{align*}
\end{definition}

The connection between generalized polar coordinates and classical polar coordinates is given in Proposition \ref{propositionrtheta} below. For now, we list some geometrical facts about generalized polar coordinates. Our first lemma says that the hyperbolic distance from a point to the origin is essentially the sum of the ``radial'' distance and the ``angular'' distance.

\begin{lemma}
\label{lemmafabs}
For all $x\in X$ we have
\[
\dox x \asymp_{\plus,\zero,\xi_1,\xi_2} \max_{i = 1}^2\busemann_{\xi_i}(x,\zero) = |r(x)| + \theta(x).
\]
\end{lemma}
\begin{proof}
The equality is trivial, so we concentrate on the asymptotic. The $\gtrsim$ direction follows directly from (f) of Proposition \ref{propositionbasicidentities}. On the other hand, by Gromov's inequality
\begin{align*}
\dox x - \max_{i = 1}^2\busemann_{\xi_i}(x,\zero)
\asymp_\plus 2\min_{i = 1}^2 \lb x|\xi_i\rb_\zero
\lesssim_\plus 2\lb \xi_1|\xi_2\rb_\zero \asymp_{\plus,\zero,\xi_1,\xi_2} 0.
\end{align*}
\end{proof}

Our next lemma describes the effect of isometries on generalized polar coordinates.

\begin{lemma}
\label{lemmatranslationdistance}
Fix $g\in \Isom(X)$ such that $\xi_1,\xi_2\in\Fix(g)$. For all $x\in X$ we have
\begin{align} \label{rtranslation}
r(g(x)) &\asymp_\plus r(x) + \log_b g'(\xi_1) = r(x) - \log_b g'(\xi_2)\\ \label{thetatranslation}
\theta(g(x)) &\asymp_\plus \theta(x),
\end{align}
with equality if $X$ is strongly hyperbolic. The implied constants are independent of $g$, $\xi_1$, and $\xi_2$.
\end{lemma}
\begin{proof}
\begin{align*}
2[r(g(x)) - r(x)]
&=_\pt \big[\busemann_{\xi_1}(g(x),\zero) - \busemann_{\xi_2}(g(x),\zero)\big] - \big[\busemann_{\xi_1}(x,\zero) + \busemann_{\xi_2}(x,\zero)\big] \noreason\\
&=_\pt \big[\busemann_{\xi_1}(x,g^{-1}(\zero)) - \busemann_{\xi_2}(x,g^{-1}(\zero))\big] - \big[\busemann_{\xi_1}(x,\zero) - \busemann_{\xi_2}(x,\zero)\big] \noreason\\
&\asymp_\plus \busemann_{\xi_1}(\zero,g^{-1}(\zero)) - \busemann_{\xi_2}(\zero,g^{-1}(\zero))\\
&\asymp_\plus \log_b g'(\xi_1) - \log_b g'(\xi_2). & \by{Proposition \ref{propositionbusemannpreserved}}
\end{align*}
Now \eqref{rtranslation} follows from Corollary \ref{corollaryproductone}.

On the other hand, by (g) of Proposition \ref{propositionbasicidentities},
\begin{align*}
\theta(g(x)) - \theta(x) &\asymp_\plus \big[\lb \xi_1|\xi_2\rb_{g(x)} - \lb \xi_1|\xi_2\rb_\zero\big] - \big[\lb \xi_1|\xi_2\rb_x - \lb \xi_1|\xi_2\rb_\zero\big]\\
&=_\pt \lb g^{-1}(\xi_1)|g^{-1}(\xi_2)\rb_x - \lb \xi_1|\xi_2\rb_x = 0,
\end{align*}
proving \eqref{thetatranslation}.
\end{proof}

%

We end this chapter by describing the relation between generalized polar coordinates and classical polar coordinates.

\begin{proposition}
\label{propositionrtheta}
If $X = \E$, $\zero = (1,\0)$, $\xi_1 = \0$, and $\xi_2 = \infty$, then
\begin{align*}
r(\xx) &= \log \|\xx\| \\
\theta(\xx) &= -\log(x_1/\|\xx\|) = -\log \cos(\ang(\xx)).
\end{align*}
\end{proposition}
Thus the notations $r$ and $\theta$ are slightly inaccurate as they really represent the logarithm of the radius and the negative logarithm of the cosine of the angle, respectively.
\begin{proof}[Proof of Proposition \ref{propositionrtheta}]
We consider first the case $\|\xx\| = 1$. Let us set $g(\yy) = \yy/\|\yy\|^2$, and note that $g\in\Isom(\E)$, $g(\zero) = \zero$, and $g(\xi_i) = \xi_{3 - i}$. On the other hand, since $\|\xx\| = 1$ we have $g(\xx) = \xx$, and so
\[
\busemann_{\xi_1}(\xx,\zero) = \busemann_{g(\xi_1)}(g(\xx),g(\zero)) = \busemann_{\xi_2}(\xx,\zero).
\]
It follows that $r(\xx) = 0$ and $\theta(\xx) = \busemann_{\xi_2}(\xx,\zero) = \busemann_\infty(\xx,\zero)$. By Proposition \ref{propositionbusemannE}, we have $\busemann_\infty(\xx,\zero) = -\log(x_1/\zero_1) = -\log(x_1/\|\xx\|) = -\log\cos(\ang(\xx))$.

The general case follows upon applying Lemma \ref{lemmatranslationdistance} to maps of the form $g_\lambda(\xx) = \lambda\xx$, $\lambda > 0$.
\end{proof}


\chapter{Discreteness} \label{sectiondiscreteness}

Let $X$ be a metric space. In this chapter we discuss several different notions of what it means for a group or semigroup $G\prec\Isom(X)$ to be \emph{discrete}. We show that these notions are equivalent in the Standard Case. Finally, we give examples to show that these notions are no longer equivalent when $X = \H^\infty$.

Throughout this chapter, the standing assumptions that $X$ is a (not necessarly hyperbolic) metric space and that $\zero\in X$ replace the paper's overarching standing assumption that $(X,\zero,b)$ is a Gromov triple (cf. \6\ref{standingassumptions2}). Of course, if $(X,\zero,b)$ is a Gromov triple then $X$ is a metric space and $\zero\in X$, and therefore all theorems in this chapter can be used in other chapters without comment.

\section{Topologies on $\Isom(X)$}
\label{subsectiontopologies}

In this section we discuss different topologies that may be put on the isometry group of the metric space $X$.

In the Standard Case, the most natural topology is the \emph{compact-open topology (COT)}, i.e. the topology whose subbasic open sets are of the form
\[
\GG(K,U) = \{f\in \Isom(X):f(K)\subset U\}
\]
where $K\subset X$ is compact and $U\subset X$ is open. When we replace $X$ by a metric space which is not proper, it is tempting to replace the compact-open topology with a ``bounded-open'' topology. However, it is hard to define such a topology in a way that does not result in pathologies. It turns out that the compact-open topology is still the ``right'' topology for many applications in an arbitrary metric space. But we are getting ahead of ourselves.

Let's start by considering the case where $X$ is an algebraic hyperbolic space, i.e., $X = \H = \H_\F^\alpha$, and figure out what topology or topologies we can put on $\Isom(\H)$. Recall from Theorem \ref{theoremisometries} that
\begin{equation}
\label{IsomH}
\Isom(\H) \equiv \PO^*(\LL;\QQ) \equiv \O^*(\LL;\QQ)/\sim
\end{equation}
where $\LL = \HH_\F^{\alpha + 1}$, $\QQ$ is the quadratic form \eqref{Qdef}, and $T_1\sim T_2$ means that $[T_1] = [T_2]$ (in the notation of Section \ref{subsectionisometries}). Thus $\Isom(\H)$ is isomorphic to a quotient of a subspace of $L(\LL)$, the set of bounded linear maps from $\LL$ to itself. This indicates that to define a topology or topologies on $\Isom(\H)$, it may be best to start from the functional analysis point of view and look for topologies on $L(\LL)$. In particular, we will be interested in the following widely used topologies on $L(\LL)$:
\begin{itemize}
\item The \emph{uniform operator topology (UOT)} is the topology on $L(\LL)$ which comes from looking it as a metric space with the metric
\[
\dist(T_1,T_2) = \|T_1 - T_2\| = \sup\{\|(T_1 - T_2)\xx\|:\xx\in\LL,\|\xx\| = 1\}.
\]
\item The \emph{strong operator topology (SOT)} is the topology on $L(\LL)$ which comes from looking at it as a subspace of the product space $\LL^\LL$. Note that in this topology,
\[
T_n\tendsto n T \;\;\;\; \Leftrightarrow \;\;\;\; T_n \xx \tendsto n T\xx \all \xx\in\LL.
\]
The strong operator topology is weaker than the uniform operator topology.
\end{itemize}

\begin{remark}
There are many other topologies used in functional analysis, for example the weak operator topology, which we do not consider here.
\end{remark}

Starting with either the uniform operator topology or the strong operator topology, we may restrict to the subspace $\O^*(\LL;\QQ)$ and then quotient by $\sim$ to induce a topology on $\Isom(\H)$ using the identification \eqref{IsomH}. For convenience, we will also call these induced topologies the uniform operator topology and the strong operator topology, respectively.

We now return to the general case of a metric space $X$. Define the \emph{Tychonoff topology} to be the topology on $\Isom(X)$ inherited from the product topology on $X^X$.

\begin{proposition}
\label{propositionCOTSOT}~
\begin{itemize}
\item[(i)] The Tychonoff topology and the compact-open topology on $\Isom(X)$ are identical.
\item[(ii)] If $X$ is an algebraic hyperbolic space, then the strong operator topology is identical to the Tychonoff topology (and thus also to the compact-open topology).
\end{itemize}
\end{proposition}
\begin{proof}~
\begin{itemize}
\item[(i)] Since subbasic sets in the Tychonoff topology take the form $\GG(\{x\},U)$, it is clear that the compact-open topology is at least as fine as the Tychonoff topology. Conversely, suppose that $\GG(K,U)$ is a subbasic open set in the Tychonoff topology, and fix $f\in \GG(K,U)$. Let $\epsilon = \dist(f(K),X\butnot U) > 0$, and let $(x_i)_1^n$ be a set of points in $K$ such that $K\subset \bigcup_1^n B(x_i,\epsilon/3)$. Then
\[
f\in\UU := \bigcap_{i = 1}^n \GG(\{x_i\},\thicken_{\epsilon/3}(f(K))).\Footnote{Here and elsewhere $\thicken_\epsilon(S) = \{x\in X : \dist(x,S) \leq \epsilon\}$.}
\]
The set $\UU$ is open in the Tychonoff topology; we claim that $\UU\subset \GG(K,U)$. Indeed, suppose that $\w f\in\UU$. Then for $x\in K$, fix $i$ with $x\in B(x_i,\epsilon/3)$; since $\w f$ is an isometry, $\dist(\w f(x),f(K)) \leq \dist(\w f(x),\w f(x_i)) + \dist(\w f(x_i),f(K)) \leq 2\epsilon/3 < \epsilon$. It follows that $\w f(x)\in U$; since $x\in K$ was arbitrary, $\w f\in \GG(K,U)$.

\item[(ii)] It is clear that the strong operator topology is at least as fine as the Tychonoff topology. Conversely, suppose that a set $\UU\subset\Isom(\H)$ is open in the strong operator topology, and fix $[T]\in\UU$. Let $T\in\O^*(\LL;\QQ)$ be a representative of $[T]$. There exist $(\vv_i)_1^n$ in $\LL$ and $\epsilon > 0$ such that for all $\w T\in \O^*(\LL;\QQ)$ satisfying $\|(\w T - T)\vv_i\| \leq \epsilon\all i$, we have $[\w T]\in \UU$. Let $\ff_0 = \ee_0$, and let $V = \lb \ff_0,\vv_1,\ldots,\vv_n\rb$. Extend $\{\ff_0\}$ to an $\F$-basis $\{\ff_0,\ff_1,\ldots,\ff_k\}$ of $V$ with the property that $B_\QQ(\ff_{j_1},\ff_{j_2}) = 0$ for all $j_1\neq j_2$. Without loss of generality, suppose that $k\geq 1$. For each $i = 1,\ldots,n$ we have $\vv_i = \sum_j \ff_j c_{i,j}$ for some $c_{i,j}\in\F$, so there exists $\epsilon_2 > 0$ such that for all $\w T\in\O^*(\LL;\QQ)$ satisfying $\|(\w T - T)\ff_j\| \leq \epsilon_2 \all j$ and $\|\sigma_{\w T} - \sigma_T\| \leq \epsilon_2$, we have $[\w T]\in \UU$.

Let
\begin{equation}
\label{IFdef}
I_\F = \begin{cases}
\{1\} & \F = \R\\
\{1,i\} & \F = \C\\
\{1,i,j,k\}& \F = \Q
\end{cases},
\end{equation}
and let
\[
F = \{\ee_0\} \cup \{\ee_0 \pm (1/2)\ff_1 \ell : j = 1,\ldots,k , \ell\in I_\F\}.
\]
Fix $\epsilon_3 > 0$ small to be determined, and for the remainder of this proof write $A\sim B$ if $\|A - B\|$ is bounded by a constant which tends to zero as $\epsilon\to 0$. Let
\[
\VV = \left\{[\w T]\in\Isom(\H): \forall \xx\in F, \;\;\exists \yy_\xx\in [\w T]([\xx]) \text{ such that } \|\yy_\xx - T\xx\| < \epsilon_3\right\}.
\]
For each $\xx\in F$, we have $[\xx]\in\H$, so the set 
\[
\{[\yy]\in\H : \exists \yy\in[\yy] \text{ such that } \|\yy - T\xx\| < \epsilon_3\}
\] 
is open in the natural topology on $\H$. It follows that $\VV$ is open in the Tychonoff topology. Moreover, $[T]\in\VV$. To complete the proof we show that $\VV\subset\UU$. Indeed, fix $[\w T]\in\VV$, and let $\yy = \yy_{\ee_0}$. There exists a representative $\w T\in\O^*(\LL;\QQ)$ such that $\w T\ee_0 = \lambda \yy$ for some $\lambda > 0$. Since
\[
-1 = \QQ(\ee_0) \sim \QQ(\yy) = \lambda^{-2}\QQ(\lambda\yy) = -\lambda^{-2},
\]
we have $\lambda\sim 1$ and thus $\w T\ee_0 \sim T\ee_0$.

Now for each $\xx\in F\butnot\{\ee_0\}$, there exists $a_\xx\in\F$ such that $\yy_\xx = \w T(\xx a_\xx)$. Fix $j = 1,\ldots,k$ and $\ell\in I_\F$. Writing $a_\pm = a_{\ee_0 \pm (1/2)\ff_j \ell}$, we have
\[
\left\|T\left(\ee_0 \pm \frac12\ff_j \ell\right) - \w T\left(\left(\ee_0 \pm \frac12\ff_j \ell\right)a_\pm\right)\right\| < \epsilon_3,
\]
i.e. $T(\ee_0 \pm (1/2)\ff_j \ell) \sim \w T((\ee_0 \pm (1/2)\ff_j \ell) a_\pm)$. Substituting $\pm = +$ and $\pm = -$ and adding the resulting equations gives
\[
2T\ee_0 \sim \w T(\ee_0 (a_+ + a_-)) + \frac12 \w T(\ff_j \ell (a_+ - a_-));
\]
using $T\ee_0 \sim \w T\ee_0$ and rearranging gives
\[
\w T(\ee_0 (2 - a_+ - a_-)) \sim \frac12 \w T(\ff_j \ell (a_+ - a_-)).
\]
Now by Lemma \ref{lemmaoperatornorm}, we have $\|\w T\|\sim 1$, and thus $\ee_0(2 - a_+ - a_-) \sim (1/2)\ff_j \ell(a_+ - a_-)$. Since $\|\ee_0 a + \ff_j \ell b\| \asymp_\times \max(|a|,|b|)$ for all $a,b\in\F$, it follows that
\[
2 - a_+ - a_- \sim \ell(a_+ - a_-) \sim 0,
\]
from which we deduce $a_+\sim a_- \sim 1$. Thus
\[
T\left(\ee_0 \pm \frac12\ff_j \ell\right) \sim \w T\left(\ee_0 \pm \frac12\ff_j \ell\right).
\]
Substituting $\pm = +$ and $\pm = -$, subtracting the resulting equations, and using the fact that $T\ee_0\sim \w T\ee_0$ gives
\[
T(\ff_j \ell) \sim \w T(\ff_j \ell).
\]
In particular, letting $\ell = 1$ we have $T\ff_j \sim \w T\ff_j$. Thus
\[
(T\ff_j) (\sigma_T\ell) \sim (\w T\ff_j) (\sigma_{\w T}\ell) \sim (T\ff_j) (\sigma_{\w T}\ell).
\]
Since this holds for all $\ell\in I_\F$, we have $\sigma_T\sim \sigma_{\w T}$. By the definition of $\sim$, this means that we can choose $\epsilon_3$ small enough so that $\|T\ff_j \ell - \w T\ff_j \ell\| \leq \epsilon_2 \all j$ and $\|\sigma_{\w T} - \sigma_T\| \leq \epsilon_2$. Then $[\w T]\in\UU$, completing the proof.
\end{itemize}
\end{proof}

\begin{proposition}
The compact-open topology makes $\Isom(X)$ into a topological group, i.e. the maps
\[
(g,h) \mapsto gh, \hspace{.5 in}
g \mapsto g^{-1}
\]
are continuous.
\end{proposition}
\begin{proof}
Fix $g_0,h_0\in\Isom(X)$, and let $\GG(\{x\},U)$ be a neighborhood of $g_0 h_0$. For some $\epsilon > 0$, we have $B(g_0 h_0(x),\epsilon)\subset U$. We claim that
\[
\GG\Big(\{h_0(x)\},B(g_0 h_0(x),\epsilon/2)\Big) \GG\Big(\{x\},B(h_0(x),\epsilon/2)\Big) \subset \GG(\{x\},U).
\]
Indeed, fix $g\in \GG(\{h_0(x)\},B(g_0 h_0(x),\epsilon/2))$ and $h\in \GG(\{x\},B(h_0(x),\epsilon/2))$. Then
\[
\dist(gh(x),g_0 h_0(x)) \leq \dist(h(x),h_0(x)) + \dist(gh_0(x),g_0 h_0(x)) \leq \epsilon/2 + \epsilon/2 = \epsilon,
\]
demonstrating that $gh(x)\in U$, and thus that the map $(g,h) \mapsto gh$ is continuous.

Now fix $g_0\in\Isom(X)$, and let $\GG(\{x\},U)$ be a neighborhood of $g_0^{-1}$. For some $\epsilon > 0$, we have $B(g_0^{-1}(x),\epsilon)\subset U$. We claim that
\[
\GG\big(\{g_0^{-1}(x)\},B(x,\epsilon)\big)^{-1} \subset \GG(\{x\},U).
\]
Indeed, fix $g\in \GG(\{g_0^{-1}(x)\},B(x,\epsilon))$. Then
\[
\dist(g^{-1}(x),g_0^{-1}(x)) = \dist(x,g g_0^{-1}(x)) \leq \epsilon,
\]
demonstrating that $g^{-1}(x)\in U$, and thus that the map $g \mapsto g^{-1}$ is continuous.
\end{proof}

\begin{remark}[{\cite[9.B(9), p.60]{Kechris}}]
\label{remarkpolish}
If $X$ is a separable complete metric space, then the group $\Isom(X)$ with the compact-open topology is a Polish space.
\end{remark}

\bigskip
\section{Discrete groups of isometries}
\label{subsectiondiscrete}

In this section we discuss several different notions of what it means for a group $G\leq\Isom(X)$ to be \emph{discrete}, and then we show that they are equivalent in the Standard Case. However, each of our notions will be distinct when $X = \H = \H_\F^\alpha$ for some infinite cardinal $\alpha$.

\begin{definition}
\label{definitiondiscreteness}
Fix $G\leq\Isom(X)$.
\begin{itemize}
\item $G$ is called \emph{strongly discrete (SD)} if for every bounded set $B \subset X$, we have
\[
\#\{g\in G: g(B) \cap B \neq \emptyset\} < \infty.
\]
\item $G$ is called \emph{moderately discrete (MD)} if for every $x\in X$, there exists an open set $U\ni x$ such that
\[
\#\{g\in G: g(U) \cap U \neq \emptyset\} < \infty.
\]
\item $G$ is called \emph{weakly discrete (WD)} if for every $x\in X$, there exists an open set $U \ni x$ such that
\[
g(U) \cap U \neq \emptyset \Rightarrow g(x) = x.
\]
\end{itemize}
\end{definition}
\begin{remark}
Strongly discrete groups are known in the literature as \emph{metrically proper}, and moderately discrete groups are known as \emph{wandering}.
\end{remark}
\begin{remark}
\label{remarkdiscretenessequivalent}
We may equivalently give the definitions as follows:
\begin{itemize}
\item $G$ is \emph{strongly discrete (SD)} if for every $R > 0$ and $x\in X$,
\begin{equation}
\label{SD2}
\#\{g\in G: \dist(x,g(x)) \leq R\} < \infty.
\end{equation}
\item $G$ is \emph{moderately discrete (MD)} if for every $x\in X$, there exists $\epsilon > 0$ such that
\begin{equation}
\label{MD2}
\#\{g\in G: \dist(x,g(x)) \leq \epsilon\} < \infty.
\end{equation}
\item $G$ is \emph{weakly discrete (WD)} if for every $x\in X$, there exists $\epsilon > 0$ such that
\begin{equation}
\label{WD2}
G(x)\cap B(x,\epsilon) = \{x\}.
\end{equation}
\end{itemize}
\end{remark}

As our naming suggests, the condition of strong discreteness is stronger than the condition of moderate discreteness, which is in turn stronger than the condition of weak discreteness.
\begin{proposition}
\label{propositionSDMDWD}
Any strongly discrete group is moderately discrete, and any moderately discrete group is weakly discrete.
\end{proposition}
\begin{proof}
It is clear from the second formulation that strongly discrete groups are moderately discrete. Let $G\leq\Isom(X)$ be a moderately discrete group. Fix $x\in X$, and let $\epsilon > 0$ be such that \eqref{MD2} holds. Letting $\epsilon' = \epsilon \wedge \min\{\dist(x,g(x)) : g(x)\neq x, g(x)\in B(x,\epsilon)\}$, we see that \eqref{WD2} holds.
\end{proof}

The reverse directions, WD \implies MD and MD \implies SD, both fail in infinite dimensions. Examples \ref{examplevalette} and \ref{exampletranslations}-\ref{examplenotstronglydiscreteBIM} are moderately discrete groups which are not strongly discrete, and Examples \ref{exampleMDWD} and \ref{exampleAutT} are weakly discrete groups which are not moderately discrete.

If $X$ is a proper metric space, then the classes MD and SD coincide, but are still distinct from WD. Example \ref{exampleAutT} is a weakly discrete group acting on a proper metric space which is not moderately discrete. We show now that MD $\Leftrightarrow$ SD when $X$ is proper:

\begin{proposition}
\label{propositionproperMDSD}
Suppose that $X$ is proper. Then a subgroup of $\Isom(X)$ is moderately discrete if and only if it is strongly discrete.
\end{proposition}
\begin{proof}
Let $G\leq\Isom(X)$ be a moderately discrete subgroup. Fix $x\in X$, and let $\epsilon > 0$ satisfy \eqref{MD2}. Fix $R > 0$ and let $K = \cl{G(\zero)}\cap B(x,R)$; $K$ is compact since $X$ is proper. The collection $\{B(g(x),\epsilon): g\in G\}$ covers $K$, so there is a finite subcover $\{B(g_i(x),\epsilon): i = 1,\ldots,n\}$. Now
\[
\#\{g\in G: \dist(x,g(x)\leq R)\} \leq \sum_{i = 1}^n \#\{g\in G: g(x)\in B(g_i(x),\epsilon)\} < \infty,
\]
i.e. \eqref{SD2} holds.
\end{proof}

\subsection{Topological discreteness}
\begin{definition}
\label{definitionparametricdiscreteness}
Let $\scrT$ be a topology on $\Isom(X)$. A group $G\leq\Isom(X)$ is \emph{$\scrT$-discrete} if it is discrete as a subspace of $\Isom(X)$ in the topology $\scrT$.
\end{definition}

Most of the time, we will let $\scrT$ be the compact-open topology (COT). The relation between COT-discreteness and our previous notions of discreteness is as follows:
\begin{proposition}
\label{propositionparametricdiscreteness}~
\begin{itemize}
\item[(i)] Any moderately discrete group is COT-discrete.
\item[(ii)] Any weakly discrete group that is acting on an algebraic hyperbolic space is COT-discrete.
\item[(iii)] Any COT-discrete group that is acting on a proper metric space is strongly discrete.
\end{itemize}
\end{proposition}
\begin{proof}~
\begin{itemize}
\item[(i)]
Let $G\leq\Isom(X)$ be moderately discrete, and let $\epsilon > 0$ satisfy \eqref{MD2}. Then the set $\UU := \GG(\{\zero\},B(\zero,\epsilon))\subset\Isom(X)$ satisfies $\#(\UU\cap G) < \infty$. But $\UU$ is a neighborhood of $\id$ in the compact-open topology. It follows that $G$ is COT-discrete.

\item[(ii)]
Suppose that $X = \H = \H_\F^\alpha$. Let $G\leq\Isom(\H)$ be weakly discrete, and by contradiction suppose it is not COT-discrete. For any finite set $F\subset\H$, let $\epsilon > 0$ be small enough so that \eqref{WD2} holds for all $x\in F$; since $G$ is not COTD, there exists $g = g_F\in G\butnot\{\id\}$ such that $\dist(x,g(x))\leq \epsilon$ for all $x\in F$, and it follows that $g(x) = x$ for all $x\in F$. Now suppose that $J$ is a finite set of indices, and let $F = \{[\ee_0]\} \cup \{[\ee_0 \pm (1/2)\ee_i] \ell : i\in J, \ell\in I_\F\}$, where $I_\F$ is as in \eqref{IFdef}. Then if $T_I$ is a representative of $g_F$ satisfying $T_J\ee_0 = \ee_0$, an argument similar to the proof of Proposition \ref{propositionCOTSOT}(ii) shows that $\sigma_{T_J} = I$ and $T_J \ee_i = \ee_i$ for all $i\in J$.

Now we define an infinite sequence of indices $(i_n)_1^\infty$ as follows: If $i_1,\ldots,i_{n - 1}$ have been defined, let $T_n = T_{\{i_1,\ldots,i_{n - 1}\}}$, and let $i_n$ be such that $\ee_{i_n}\notin \Fix(T_n)$.

Choose a nonnegative summable sequence $(t_n)_1^\infty$, and let $\xx = \ee_0 + \sum_{n = 1}^\infty t_n \ee_{i_n}$. Then $T_n\xx\to \xx$; since $G$ is weakly discrete, it follows that $T_n\xx = \xx$ for all $n$ sufficiently large. Fix such an $n$, and observe that
\[
0 = T_n\xx - \xx = t_n (T_n(\ee_n) - \ee_n) + \sum_{m > n} t_m (T_n(\ee_m) - \ee_m);
\]
the triangle inequality gives
\[
t_n \leq \frac{\sum_{m > n}2t_m}{\|T_n\ee_n - \ee_n\|}\cdot
\]
By choosing the sequence $(t_n)_1^\infty$ to satisfy
\[
t_{n + 1} < \frac14 \|T_n\ee_n - \ee_n\| t_n \leq \frac12 t_n ,
\]
we arrive at a contradiction.

%
%
%

\item[(iii)]
Let $G$ be a COT-discrete group acting by isometries on a proper metric space $X$. By contradiction, suppose that $G$ is not strongly discrete. Then there exists an infinite set $A\subset G$ such that the set $A(\zero)$ is bounded. Without loss of generality we may suppose that $A^{-1} = A$. Note that for each $x\in X$, the set $A(x)$ is bounded and therefore precompact. Now since $X$ is a proper metric space, it is $\sigma$-compact and therefore separable. Let $\SS$ be a countable dense subset of $X$. Then
\[
\KK := \left(\prod_{q\in\SS}\cl{A(q)}\right)^2
\]
is a compact metrizable space. For each $g\in A$ let
\[
\phi_g := \left((g(q))_{q\in\SS},(g^{-1}(q))_{q\in\SS}\right)\in \KK.
\]
Since $A$ is infinite, there exists an infinite sequence $(g_n)_1^\infty$ in $A$ such that
\[
\phi_{g_n}\to \left((y_q^{(+)})_{q\in\SS},(y_q^{(-)})_{q\in\SS}\right)\in\KK.
\]
Thus
\[
g_n^\pm(q) \tendsto n y_q^{(\pm)} \all q\in\SS.
\]
The density of $\SS$ and the equicontinuity of the sequences $(g_n)_1^\infty$ and $(g_n^{-1})_1^\infty$ imply that for all $x\in X$, there exist $y_x^{(\pm)}$ such that $g_n^\pm(y)\to y_x^{(\pm)}$. Thus, the sequence $(g_n)_1^\infty$ converges in the Tychonoff topology to some $g^{(+)}\in X^X$. Similarly, the sequence $(g_n^{-1})_1^\infty$ converges to some $g^{(-)}\in X^X$. Again, the equicontinuity of the sequences $(g_n)_1^\infty$ and $(g_n^{-1})_1^\infty$ allows us to take limits and deduce that
\[
g^{(+)} g^{(-)} = \lim_{n\to\infty} g_n g_n^{-1} = \id.
\]
Similarly, $g^{(-)}g^{(+)} = \id$. Thus $g^{(+)}$ and $g^{(-)}$ are inverses, and in particular $g^{(+)}\in\Isom(X)$. Since $g_n\to g^{(+)}$ in the compact-open topology, the proof is completed by the following lemma from topological group theory:
\begin{lemma}
\label{lemmadiscretesubgroup}
Let $H$ be a topological group, and let $G$ be a subgroup of $H$. Suppose there is a sequence $(g_n)_1^\infty$ of distinct elements in $G$ which converges to an element of $H$. Then $G$ is not discrete in the topology inherited from $H$.
\end{lemma}
\begin{subproof}
Suppose $g_n\to h\in H$. Then
\[
g_n g_{n + 1}^{-1} \to h h^{-1} = \id,
\]
while on the other hand $g_n g_{n + 1}^{-1}\neq \id$ (since the sequence $(g_n)_1^\infty$ consists of distinct elements). This demonstrates that $G$ is not discrete in the inherited topology.
\end{subproof}
\end{itemize}
\end{proof}
If $X$ is not an algebraic hyperbolic space, then it is possible for a weakly discrete group to not be COT-discrete; see Example \ref{exampleAutT}. Conversely, it is possible for a COT-discrete group to not be weakly discrete; see Examples \ref{examplepoincareextension} amd \ref{examplepoincareextensiontwo}.

On the other hand, suppose that $X$ is an algebraic hyperbolic space. The uniform operator topology (abbreviated as UOT) is finer than the COT, i.e. it has more open sets, and therefore it is easier for every subset of $G$ to be relatively open in that topology, which is exactly what it means to be discrete. Notice that there is an ``order switch'' here; the UOT is finer than the COT, but the condition of being COT-discrete is stronger than the condition of being UOT-discrete. We record this for later use as the following 
\begin{observation}
\label{observationCOTUOT}
Let $X$ be an algebraic hyperbolic space. If a subgroup $G\leq\Isom(X)$ is COT-discrete, then it is also UOT-discrete.
\end{observation}

The inclusion in the previous observation is strict. A significant example of a group acting on $\H^\infty$ which is UOT-discrete but not COT-discrete is described in Example \ref{exampleAutTBIM}.

The various relations between the distinct shades of discreteness are somewhat subtle when first discerned. We speculate that it may be fruitful to study such distinctions with a finer lens. For the reader's ease, we summarize the relations between our different notions of discreteness in Table \ref{figurediscreteness} below.

\subsection{Equivalence in finite dimensions}

\begin{proposition}\label{propositionfinitedimequivalent}
Suppose that $X$ is a finite-dimensional Riemannian manifold. Then the notions of strong discreteness, moderate discreteness, weak discreteness, and COT-discreteness agree. If $X$ is an algebraic hyperbolic space, these notions also agree with the notion of UOT-discreteness.
\end{proposition}
\begin{proof}
By Propositions \ref{propositionSDMDWD} and \ref{propositionparametricdiscreteness}, the conditions of strong discreteness, moderate discreteness, and COT-discreteness agree and imply weak discreteness. Conversely, suppose that $G\leq\Isom(X)$ is weakly discrete, and by contradiction suppose that $G$ is not COT-discrete. Since $X$ is separable, so is $\Isom(X)$, and thus there exists a sequence $\Isom(X)\butnot\{\id\}\ni g_n\to\id$ in the compact-open topology. For each $n$ let $F_n = \{x\in X : g_n(x) = x\}$. Since $G$ is weakly discrete, $X = \bigcup_1^\infty F_n$, so by the Baire category theorem, $F_n$ has nonempty interior for some $n$. But then $g_n = \id$ on an open set; in particular there exists a point $x_0\in X$ such that $g_n(x_0) = x_0$ and $g_n'(x_0)$ is the identity map on the tangent space of $x_0$. By the naturality of the exponential map, this implies that $g_n$ is the identity map, a contradiction.

Finally, suppose $X = \H = \H_\F^\alpha$ is an algebraic hyperbolic space, and let $\LL = \LL_\F^{\alpha + 1}$. Since $\LL$ is finite-dimensional, the SOT and UOT topologies on $L(\LL)$ are equivalent. This in turn  demonstrates that the notions of COT-discreteness and UOT-discreteness agree.
\end{proof}
In such a setting, we shall call a group satisfying any of these equivalent definitions simply \emph{discrete}.

\subsection{Proper discontinuity}
\begin{definition}
\label{definitionproperdiscontinuity}
A group $G\leq\Isom(X)$ \emph{acts properly discontinuously (PrD)} on $X$ if for every $x\in X$, there exists an open set $U \ni x$ with 
\[
g(U) \cap U \neq \emptyset \Rightarrow g = \id,
\]
or equivalently, if
\[
\dist(x,\{g(x):g\neq\id\}) > 0.
\]
\end{definition}

Let us discuss the relations between proper discontinuity and some of our notions of discreteness. We begin by noting that even in finite dimensions, the notion of proper discontinuity is not the same as the notion of discreteness; instead, a group acts properly discontinuously if and only if it both discrete and torsion-free. We also remark that in finite dimensions Selberg's lemma (see e.g. \cite{Alperin}) can be used to pass from a discrete group to a finite-index subgroup that acts properly discontinuously. However, it is impossible to do this in infinite dimensions; cf. Example \ref{exampleparabolictorsion}.

Although no notion of discreteness implies proper discontinuity, the reverse is true for certain types of discreteness. Namely, since $\#\{\id\} = 1 < \infty$, we have:
\begin{observation}
\label{observationproperlydiscontinuous}
Any group which acts properly discontinuously is moderately discrete.
\end{observation}

In particular, by combining with Proposition \ref{propositionproperMDSD} we see that if $X$ is proper then any group which acts properly discontinuously is strongly discrete. This provides a connection between our results, in which strong discreteness is often a hypothesis, and many results from the literature in which proper discontinuity and properness are both hypotheses.

Observation \ref{observationproperlydiscontinuous} admits the following partial converse, which generalizes the fact that in finite dimensions every discrete torsion-free group acts properly discontinuously: 
\begin{remark}
If $X$ is a proper CAT(0) space, then a group acts properly discontinously if and only if it is moderately discrete and torsion free.
\end{remark}
\begin{proof}
Suppose that $G\leq\Isom(X)$ acts properly discontinuously. If $g\in G\butnot\{\id\}$ is a torsion element, then by Cartan's lemma \cite[II.2.8(1)]{BridsonHaefliger}, $g$ has a fixed point. This contradicts $G$ acting properly discontinuously. Thus $G$ is torsion-free.

Conversely, suppose that $G\leq\Isom(X)$ is moderately discrete and torsion-free. Given $x\in X$, let $\epsilon > 0$ be as in \eqref{WD2}, and by contradiction suppose that there exists $g\neq\id$ such that $\dist(x,g(x)) < \epsilon$. By \eqref{WD2}, $g(x) = x$. But then by \eqref{MD2}, the set $\{g^n : n\in\Z\}$ is finite, i.e. $g$ is a torsion element. This is a contradiction, so $G$ acts properly discontinuously.
\end{proof}

We summarize the relations between our various notions of discreteness, together with proper discontinuity, in the following table:

\begin{table}[h!]
\begin{tabular}{ | c | c c c c c c c | }
 \hline
 Finite dimensional			& SD & $\leftrightarrow$ & MD & $\leftrightarrow$ & WD && \\
 Riemannian manifold 	& $\uparrow$ & & & & $\updownarrow$ && \\
	 			& PrD & & & & COTD & $\leftrightarrow$ & UOTD\\
 \hline
	 		& SD & $\rightarrow$ & MD & $\rightarrow$ & WD && \\
 General metric space		& & $\nearrow$ & & $\searrow$ & && \\
				& PrD & & & & COTD && \\
 \hline
 Infinite dimensional				& SD & $\rightarrow$ & MD & $\rightarrow$ & WD && \\
 algebraic hyperbolic space 	& & $\nearrow$ & & & $\downarrow$ && \\
			& PrD & & & & COTD & $\rightarrow$ & UOTD \\
 \hline
	 		& SD & $\leftrightarrow$ & MD & $\leftrightarrow$ & COTD && \\
 Proper metric  space 		& $\uparrow$ & & & & $\downarrow$ &&\\
			& PrD & & & & WD && \\
 \hline
\end{tabular}
\vskip0.12cm
\caption{The relations between different notions of discreteness. COTD and UOTD stand for discrete with respect to the compact-open and uniform operator topologies respectively. All implications not listed have counterexamples; see Chapter \ref{sectionexamples}.}
\label{figurediscreteness}
\end{table}

\subsection{Behavior with respect to restrictions}
Fix $G\leq\Isom(X)$, and suppose $Y\subset X$ is a subspace of $X$ preserved by $G$, i.e. $g(Y) = Y$ for all $g\in G$. Then $G$ can be viewed as a group acting on the metric space $(Y,\dist\given_Y)$.

\begin{observation}~
\label{observationrestrictions}
\begin{itemize}
\item[(i)] $G$ is strongly discrete $\Leftrightarrow$ $G\given Y$ is strongly discrete
\item[(ii)] $G$ is moderately discrete \implies $G\given Y$ is moderately discrete
\item[(iii)] $G$ is weakly discrete \implies $G\given Y$ is weakly discrete
\item[(iv)] $G$ is $\scrT$-discrete $\Leftarrow$ $G\given Y$ is $\scrT\given Y$-discrete
\item[(v)] $G$ acts properly discontinuously on $X$ \implies $G$ acts properly discontinuously on $Y$.
\end{itemize}
\end{observation}
In particular, strong discreteness is the only concept which is ``independent of the space being acted on''. It is thus the most robust of all our definitions.

Note that for the notions of topological discreteness like COTD and UOTD, the order of implication reverses; restricting to a subspace may cause a group to no longer be discrete. Example \ref{examplepoincareextension} is an example of this phenomenon.

\subsection{Countability of discrete groups}
In finite dimensions, all discrete groups are countable. In general, it depends on what type of discreteness you are considering.

\begin{proposition}
Fix $G\leq\Isom(X)$, and suppose that either
\begin{itemize}
\item[(1)] $G$ is strongly discrete, or
\item[(2)] $X$ is separable and $G$ is COT-discrete.
\end{itemize}
Then $G$ is countable.
\end{proposition}
\begin{proof}
If $G$ is strongly discrete, then
\[
\#(G) \leq \sum_{n\in\N}\#\{g\in G:\dogo g\leq n\} \leq \sum_{n\in\N}\#(\N) = \#(\N).
\]
On the other hand, if $X$ is a separable metric space, then by Remark \ref{remarkpolish} $\Isom(X)$ is separable metrizable, so it contains no uncountable discrete subspaces.
\end{proof}
\begin{remark}
An example of an uncountable UOT-discrete subgroup of $\Isom(\H^\infty)$ is given in Example \ref{exampleAutTBIM}, and an example of an uncountable weakly discrete group acting on a separable $\R$-tree is given in Example \ref{exampleAutT}. An example of an uncountable moderately discrete group acting on a (non-separable) $\R$-tree is given in Remark \ref{remarkMDuncountable}.
\end{remark}


\draftnewpage
\chapter{Classification of isometries and semigroups} \label{sectionclassification}

In this chapter we classify subsemigroups $G\prec\Isom(X)$ into six categories, depending on the behavior of the orbit of the basepoint $\zero\in X$. We start by classifying individual isometries, although it will turn out that the category into which an isometry is classified is the same as the category of the cyclic group that it generates.

We remark that if $X$ is geodesic and $G\leq\Isom(X)$ is a group, then the main results of this chapter were proven in \cite{Hamann}. Moreover, our terminology is based on \cite[\63.A]{CCMT}, where a similar classification was given based on \cite[\6 3.1]{Gromov3}.

\bigskip
\section{Classification of isometries}
Fix $g\in\Isom(X)$, and let
\[
\Fix(g) := \{x\in\bord X:g(x) = x\}.
\]
Consider $\xi\in\Fix(g)\cap\del X$. Recall that $g'(\xi)$ denotes the dynamical derivative of $g$ at $\xi$ (see \6\ref{subsubsectiondynamicalderivative}). 
\begin{definition}
\label{definitionderivativefixedpoint}
$\xi$ is said to be
\begin{itemize}
\item a \emph{neutral} or \emph{indifferent} fixed point if $g'(\xi) = 1$,
\item an \emph{attracting} fixed point if $g'(\xi) < 1$, and
\item a \emph{repelling} fixed point if $g'(\xi) > 1$.
\end{itemize}
\end{definition}

\begin{definition}
\label{definitionclassification1}
An isometry $g\in\Isom(X)$ is called
\begin{itemize}
\item \emph{elliptic} if the orbit $\{g^n(\zero):n\in\N\}$ is bounded,
\item \emph{parabolic} if it is not elliptic and has a unique fixed point in $\del X$, which is neutral, and
\item \emph{loxodromic} if it has exactly two fixed points in $\del X$, one of which is attracting and the other of which is repelling.
\end{itemize}
\end{definition}
\begin{remark}
We use the terminology ``loxodromic'' rather than the more common ``hyperbolic'' to avoid confusion with the many other meanings of the word ``hyperbolic''. In particular, when we get to classification of groups it would be a bad idea to call any group ``hyperbolic'' if it is not hyperbolic in the sense of Gromov.
\end{remark}

The categories of elliptic, parabolic, and loxodromic are clearly mutually exclusive.\Footnote{Proposition \ref{propositionbusemannpreserved} can be used to show that loxodromic isometries are not elliptic.} In the converse direction we have the following:
\begin{theorem}
\label{classificationofisometries}
Any isometry is either elliptic, parabolic, or loxodromic.
\end{theorem}

The proof of Theorem \ref{classificationofisometries} will proceed through several lemmas.

\begin{lemma}[A corollary of {\cite[Proposition 5.1]{Karlsson}}]
\label{lemmakarlsson}
If $g\in\Isom(X)$ is not elliptic, then $\Fix(g)\cap\del X\neq\emptyset$.
\end{lemma}
We include the proof for completeness.

\begin{proof}
For each $t\in\N$, let $n_t$ be the smallest integer such that
\[
\dogo{g^{n_t}} \geq n_t.
\]
The sequence $(n_t)_1^\infty$ is nondecreasing. Given $s,t\in\N$ with $s < t$, we have
\[
\dist(g^{n_s}(\zero) , g^{n_t}(\zero)) = \dogo{g^{n_t - n_s}} < n_t,
\]
and thus
\[
\lb g^{n_s}(\zero) | g^{n_t}(\zero)\rb_\zero > \frac12 [ n_s + n_t - n_t] = \frac12 n_s \tendsto{s,t} \infty,
\]
i.e. $(g^{n_t}(\zero))_t$ is a Gromov sequence. Let $\xi = [(g^{n_t}(\zero))_t]$, and note that
\[
\lb \xi|g(\xi)\rb_\zero = \lim_{t\to\infty} \lb g^{n_t}(\zero) | g^{n_t + 1}(\zero) \rb \geq \lim_{t\to\infty} \left[\dogo{g^{n_t}} - \dist(g^{n_t}(\zero), g^{n_t + 1}(\zero))\right] = \infty.
\]
Thus $g(\xi) = \xi$, i.e. $\xi\in\Fix(g)\cap\del X$.
\end{proof}

\begin{remark}[{\cite[Proposition 5.2]{Karlsson}}]
If $g\in\Isom(X)$ is elliptic and if $X$ is CAT(0), then $\Fix(g)\cap X\neq\emptyset$ due to Cartan's lemma (Theorem \ref{Cartan's lemma} below). Thus if $X$ is a CAT(0) space, then any isometry of $X$ has a fixed point in $\bord X$.
\end{remark}

\begin{lemma}
\label{lemmarepellingloxodromic}
If $g\in\Isom(X)$ has an attracting or repelling periodic point, then $g$ is loxodromic.
\end{lemma}
\begin{proof}
Suppose that $\xi\in\del X$ is a repelling fixed point for $g\in\Isom(X)$, i.e. $g'(\xi) > 1$. Recall from Proposition \ref{propositioneuclideansimilarity} that
\[
\Dist_\xi(g^n(y_1), g^n(y_2)) \leq C g'(\xi)^{-n} \Dist_\xi(y_1, y_2) \all y_1,y_2\in\EE_\xi\all n\in\Z
\]
for some constant $C > 0$. Now let $n$ be large enough so that $g'(\xi)^n > C$; then the above inequality shows that the map $g^n$ is a strict contraction of the complete metametric space $(\EE_\xi,\Dist_\xi)$ (cf. Proposition \ref{propositioneuclideanmetametric}). Then by Theorem \ref{banachcontractionprinciple}, $g$ has a unique fixed point $\eta\in(\EE_\xi)_\refl = \del X\butnot\{\xi\}$. By Corollary \ref{corollaryproductone}, $\eta$ is an attracting fixed point. Corollary \ref{corollaryproductone} also implies that $g$ cannot have a third fixed point. Thus $g$ is loxodromic.

On the other hand, if $g$ has an attracting fixed point, then by Proposition \ref{propositionderivativeiteration}, $g^{-1}$ has a repelling fixed point. Thus $g^{-1}$ is loxodromic, so applying Proposition \ref{propositionderivativeiteration} again, we see that $g$ is loxodromic.
\end{proof}

\begin{proof}[Proof of Theorem \ref{classificationofisometries}]
By contradiction suppose that $g$ is not elliptic or loxodromic, and we will show that it is parabolic. By Lemma \ref{lemmakarlsson}, we have $\Fix(g)\cap\del X\neq\emptyset$; on the other hand, by Lemma \ref{lemmarepellingloxodromic}, every fixed point of $g$ in $\del X$ is neutral. It remains to show that $\#(\Fix(g)) = 1$. By contradiction, suppose otherwise. Since $g$ is not elliptic, we clearly have $\Fix(g)\cap X = \emptyset$. Thus we may suppose that there are two distinct neutral fixed points $\xi_1,\xi_2\in\del X$.

Now for each $n\in\N$, we have
\[
\busemann_{\xi_i}(\zero,g^n(\zero)) \asymp_\plus n\log_b(g'(\xi_i)) = 0, \hspace{.5 in} i = 1,2
\]
by Proposition \ref{propositionbusemannpreserved}. Let $r = r_{\xi_1,\xi_2,\zero}$ and $\theta = \theta_{\xi_1,\xi_2,\zero}$ be as in Section \ref{subsectionpolar}. Then by Lemma \ref{lemmatranslationdistance} we have $r(g^n(\zero))\asymp_\plus \theta(g^n(\zero)) \asymp_\plus 0$. Thus by Lemma \ref{lemmafabs} we have
\[
\dogo{g^n} \asymp_\plus |r(g^n(\zero))| + \theta(g^n(\zero)) \asymp_\plus 0,
\]
i.e. the sequence $\{g^n(\zero):n\in\N\}$ is bounded. Thus $g$ is elliptic, contradicting our hypothesis.
\end{proof}

\begin{remark}[Cf. {\cite[Chapter 3, Theorem 1.4]{Chiswell}}]
\label{remarkparabolicRtree}
For $\R$-trees, parabolic isometries are impossible, so Theorem \ref{classificationofisometries} shows that every isometry is elliptic or loxodromic.
\end{remark}
\begin{proof}
By contradiction suppose that $X$ is an $\R$-tree and that $g\in\Isom(X)$ is a parabolic isometry with fixed point $\xi\in\del X$. Let $x = C(\zero,g(\zero),\xi)\in X$; then $x = \geo\zero\xi_t$ for some $t\geq 0$. Now,
\[
\dist(g(\zero),x) = \dox x + \busemann_\xi(g(\zero),\zero) = t + 0 = t.
\]
It follows that $g(x) = \geo{g(\zero)}{\xi}_t = x$. Thus $g$ is elliptic, a contradiction.
\end{proof}

\subsection{More on loxodromic isometries}

\begin{notation}
Suppose $g\in\Isom(X)$ is loxodromic. Then $g_+$ and $g_-$ denote the attracting and repelling fixed points of $g$, respectively.
\end{notation}

\begin{theorem}
\label{theoremloxodromicextra}
Let $g\in\Isom(X)$ be loxodromic. Then
\begin{equation}
\label{loxodromicinversederivative}
g'(g_+) = \frac{1}{g'(g_-)}\cdot
\end{equation}
Furthermore, for every $x\in \bord X\butnot\{g_-\}$ and for every $n\in\N$ we have
\begin{equation}
\label{loxodromiccontracting}
\Dist(g^n(x), g_+) \lesssim_\times \frac{[g'(g_+)]^n}{\Dist(g_-,g_+)\Dist(x, g_-)},
\end{equation}
with $\leq$ if $X$ is strongly hyperbolic. In particular
\[
x\neq g_- \;\;\Rightarrow\;\; g^n(x)\tendsto n g_+,
\]
and the convergence is uniform over any set whose closure does not contain $g_-$. Finally,
\begin{equation}
\label{dognoloxodromic}
\dogo{g^n} \asymp_\plus |n| \log_b g'(g_-) = |n| \log_b \frac{1}{g'(g_+)}\cdot
\end{equation}
\end{theorem}
\begin{proof}
\eqref{loxodromicinversederivative} follows directly from Corollary \ref{corollaryproductone}.

To demonstrate \eqref{loxodromiccontracting}, note that
\begin{align*}
&\lb x|g_-\rb_\zero + \lb g^n(x)|g_+\rb_\zero\\
&\gtrsim_\plus \busemann_{g_-}(\zero,x) + \busemann_{g_+}(\zero,g^n(x)) \by{(j) of Proposition \ref{propositionbasicidentities}}\\
&\asymp_\plus \busemann_{g_-}(\zero,x) + \busemann_{g_+}(\zero,x) - n \log_b g'(g_+) \by{Proposition \ref{propositionbusemannpreserved}}\\
&\asymp_\plus \lb g_-|g_+\rb_x - \lb g_-|g_+\rb_\zero - n \log_b g'(g_+) \by{(g) of Proposition \ref{propositionbasicidentities}}\\
&\geq_\pt -\lb g_-|g_+\rb_\zero - n\log_b g'(g_+).
\end{align*}
Exponentiating and rearranging yields \eqref{loxodromiccontracting}.

Finally, \eqref{dognoloxodromic} follows directly from Lemmas \ref{lemmafabs} and \ref{lemmatranslationdistance}.
\end{proof}

%
%

\subsection{The story for real hyperbolic spaces}
If $X$ is a real hyperbolic space, then we may conjugate each $g\in\Isom(X)$ to a ``normal form'' whose geometrical significance is clearer. The normal form will depend on the classification of $g$ as elliptic, parabolic, or hyperbolic.

\begin{proposition}\label{class}
Let $X$ be a real hyperbolic space, and fix $g\in\Isom(X)$.
\begin{itemize}
\item[(i)] If $g$ is elliptic, then $g$ is conjugate to a map of the form $T \given \B$ for some linear isometry $T\in\O^*(\HH)$.
\item[(ii)] If $g$ is parabolic, then $g$ is conjugate to a map of the form $\xx\mapsto \what T \xx + \pp:\E\to\E$, where $T\in\O^*(\BB)$ and $\pp\in\BB$. Here $\BB = \del\E\butnot\{\infty\} = \HH^{\alpha - 1}$.
\item[(iii)] If $g$ is hyperbolic, then $g$ is conjugate to a map of the form $\xx\mapsto \lambda \what T \xx: \E \to \E$, where $0< \lambda < 1$ and $T\in\O^*(\BB)$.
\end{itemize}
\end{proposition}
\begin{proof}~
\begin{itemize}
\item[(i)] If $g$ is elliptic, then by Cartan's lemma (Theorem \ref{Cartan's lemma} below), $g$ has a fixed point $x\in X$. Since $\Isom(X)$ acts transitively on $X$ (Observation \ref{observationtransitivity}), we may conjugate to $\B$ in a way such that $g(\0) = \0$. But then by Proposition \ref{propositionIsomB}, $g$ is of the form (i).
\item[(ii)] Let $\xi$ be the neutral fixed point of $g$. Since $\Isom(X)$ acts transitively on $\del X$ (Proposition \ref{propositiontransitivity}), we may conjugate to $\E$ in a way such $g(\infty) = \infty$. Then by Proposition \ref{propositionIsomE} and Example \ref{examplegprimeinfty}, $g$ is of the form (ii).
\item[(iii)] Since $\Isom(X)$ acts doubly transitively on $\del X$ (Proposition \ref{propositiontransitivity}), we may conjugate to $\E$ in a way such that $g_+ = \0$ and $g_- = \infty$. Then by Proposition \ref{propositionIsomE} and Example \ref{examplegprimeinfty}, $g$ is of the form (iii). (We have $\pp = \0$ since $\0\in\Fix(g)$.)
\end{itemize}
\end{proof}


\begin{remark}
If $g\in\Isom(X)$ is elliptic or loxodromic, then the orbit $(g^n(\zero))_1^\infty$ exhibits some ``regularity'' - either it remains bounded forever, or it diverges to the boundary. On the other hand, if $g$ is parabolic then the orbit can oscillate, both accumulating at infinity and returning infinitely often to a bounded region. This is in sharp contrast to finite dimensions, where such behavior is impossible. We discuss such examples in detail in \6\ref{subsubsectionedelstein}.
\end{remark}
\bigskip
\section{Classification of semigroups}

\begin{notation}
\label{notationglobalfixedpoints}
We denote the set of \emph{global fixed points} of a semigroup $G\prec\Isom(X)$ by
\[
\Fix(G) := \bigcap_{g\in G}\Fix(g).
\]
\end{notation}

\begin{definition}
\label{definitionclassification2}
$G$ is
\begin{itemize}
\item \emph{elliptic} if $G(\zero)$ is a bounded set.
\item \emph{parabolic} if $G$ is not elliptic and has a global fixed point $\xi\in \Fix(G)$ such that 
\[
g'(\xi) = 1 \all g\in G,
\]
i.e. $\xi$ is neutral with respect to every element of $G$.
\item \emph{loxodromic} if it contains a loxodromic isometry.
\end{itemize}
\end{definition}

Below we shall prove the following theorem:

\begin{theorem}
\label{classificationofsemigroups}
Every semigroup of isometries of a hyperbolic metric space is either elliptic, parabolic, or loxodromic.
\end{theorem}

\begin{observation}
An isometry $g$ is elliptic, parabolic, or loxodromic according to whether the cyclic group generated by it is elliptic, parabolic, or loxodromic. A similar statement holds if ``group'' is replaced by ``semigroup''. Thus, Theorem \ref{classificationofisometries} is a special case of Theorem \ref{classificationofsemigroups}.
\end{observation}

Before proving Theorem \ref{classificationofsemigroups}, let us say a bit about each of the different categories in this classification.

\subsection{Elliptic semigroups}
Elliptic semigroups are the least interesting of the semigroups we consider. Indeed, we observe that any strongly discrete elliptic semigroup is finite. We now consider the question of whether every elliptic semigroup has a global fixed point.

\begin{theorem}[Cartan's lemma] \label{Cartan's lemma}
If $X$ is a CAT(0) space (and in particular if $X$ is a CAT(-1) space), then every elliptic subsemigroup $G\prec\Isom(X)$ has a global fixed point.
\end{theorem}
We remark that if $G$ is a group, then this result may be found as \cite[Corollary II.2.8(1)]{BridsonHaefliger}.
\begin{proof}
Since $G(\zero)$ is a bounded set, it has a unique circumcenter \cite[Proposition II.2.7]{BridsonHaefliger}, i.e. the minimum
\[
\min_{x\in X}\sup_{g\in G}\dist(x,g(\zero))
\]
is achieved at a single point $x\in X$. We claim that $x$ is a global fixed point of $G$. Indeed, for each $h\in G$ we have
\[
\sup_{g\in G}\dist(h^{-1}(x),g(\zero)) = \sup_{g\in G}\dist(x,h g(\zero)) \leq \sup_{g\in G}\dist(x,g(\zero));
\]
since $x$ is the circumcenter we deduce that $h^{-1}(x) = x$, or equivalently that $h(x) = x$.
\end{proof}

On the other hand, if we do not restrict to CAT(0) spaces, then it is possible to have an elliptic group with no global fixed point. We have the following simple example:
\begin{example}
Let $X = \B\butnot B_\B(\0,1)$ and let $g(\xx) = -\xx$. Then $X$ is a hyperbolic metric space, $g$ is an isometry of $X$, and $G = \{\id,g\}$ is an elliptic group with no global fixed point.
\end{example}

\subsection{Parabolic semigroups}
Parabolic semigroups will be important in Chapter \ref{sectionGF} when we consider geometrically finite semigroups. In particular, we make the following definition:
\begin{definition}
\label{definitionparabolic}
Let $G\prec\Isom(X)$. A point $\xi\in\del X$ is a \emph{parabolic fixed point} of $G$ if the semigroup
\[
G_\xi := \Stab(G;\xi) = \{g\in G:g(\xi) = \xi\}
\]
is a parabolic semigroup.
\end{definition}
In particular, if $G$ is a parabolic semigroup then the unique global fixed point of $G$ is a parabolic fixed point.
\begin{warning}
A parabolic group does not necessarily contain a parabolic isometry; see Example \ref{exampleparabolictorsion}.
\end{warning}

Note that Proposition \ref{propositioneuclideansimilarity} yields the following observation:
\begin{observation}
\label{observationuniformlyLipschitz}
Let $G\prec\Isom(X)$, and let $\xi$ be a parabolic fixed point of $G$. Then the action of $G_\xi$ on $(\EE_\xi,\Dist_\xi)$ is uniformly Lipschitz, i.e.
\[
\Dist_\xi(g(y_1),g(y_2)) \asymp_\times \Dist_\xi (y_1,y_2) \all y_1,y_2\in\EE_\xi\all g\in G,
\]
and the implied constant is independent of $g\in G$. Furthermore, if $X$ is strongly hyperbolic, then $G$ acts isometrically on $\EE_\xi$.
\end{observation}

\begin{observation}
\label{observationeuclideanparabolicasymp}
Let $G\prec\Isom(X)$, and let $\xi$ be a parabolic fixed point of $G$. Then for all $g\in G_\xi$,
\begin{equation}
\label{euclideanparabolicasymp}
\Dist_\xi(\zero,g(\zero)) \asymp_\times b^{(1/2)\dogo g},
\end{equation}
with equality if $X$ is strongly hyperbolic.
\end{observation}
\begin{proof}
This is a direct consequence of \eqref{euclideancomparison}, (h) of Proposition \ref{propositionbasicidentities}, and Proposition \ref{propositionbusemannpreserved}.
\end{proof}

As a corollary we have the following:
\begin{observation}
\label{observationparabolic}
Let $G\prec\Isom(X)$, and let $\xi$ be a parabolic fixed point of $G$. Then for any sequence $(g_n)_1^\infty$ in $G_\xi$,
\[
\dogo{g_n}\tendsto n \infty \Leftrightarrow g_n(\zero)\tendsto n\xi.
\]
\end{observation}
\begin{proof}
Indeed,
\[
g_n(\zero)\tendsto n\xi \Leftrightarrow \Dist_\xi(\zero,g_n(\zero))\tendsto n \infty \Leftrightarrow \dogo{g_n}\tendsto n \infty.
\]
\end{proof}

\begin{remark}
\label{remarkparabolicRtree2}
If $X$ is an $\R$-tree, then any parabolic group must be infinitely generated. This follows from a straightforward modification of the proof of Remark \ref{remarkparabolicRtree}.
\end{remark}

\subsection{Loxodromic semigroups}
\label{subsubsectionloxodromic}
We now come to loxodromic semigroups, which are the most diverse out of these classes. In fact, they are so diverse that we separate them into three subclasses.
\begin{definition}[\cite{CCMT}]
\label{definitionCCMT}
Let $G\prec\Isom(X)$ be a loxodromic semigroup. $G$ is
\begin{itemize}
\item \emph{lineal} if $\Fix(g) = \Fix(h)$ for all loxodromic $g,h\in G$.
\item \emph{of general type} if it has two loxodromic elements $g,h\in G$ with 
\[
\Fix(g)\cap\Fix(h) = \emptyset.
\]
\item \emph{focal} if $\#(\Fix(G)) = 1$.
\end{itemize}
(We remark that focal groups were called \emph{quasiparabolic} by Gromov \cite[\63, Case 4']{Gromov3}.
\end{definition}
We observe that any cyclic loxodromic group or semigroup is lineal, so this refined classification does not give any additional information for individual isometries.

\begin{proposition}
Any loxodromic semigroup is either lineal, focal, or of general type.
\end{proposition}
\begin{proof}
Clearly, $\#(\Fix(G))\leq 2$ for any loxodromic semigroup $G$; moreover, $\#(\Fix(G)) = 2$ if and only if $G$ is lineal. So to complete the proof, it suffices to show that $\#(\Fix(G)) = 0$ if and only if $G$ is of general type. The backward direction is obvious. Suppose that $\#(\Fix(G)) = 0$, but that $G$ is not of general type. Combinatorial considerations show that there exist three points $\xi_1,\xi_2,\xi_3\in\del X$ such that $\Fix(g)\subset\{\xi_1,\xi_2,\xi_3\}$ for all $g\in G$. But then the set $\{\xi_1,\xi_2,\xi_3\}$ would have to be preserved by every element of $g$, which contradicts the definition of a loxodromic isometry.
\end{proof}

Let $G$ be a focal semigroup, and let $\xi$ be the global fixed point of $G$. The dynamics of $G$ will be different depending on whether or not $g'(\xi) > 1$ for any $g\in G$.
\begin{definition}
\label{definitionfocal}
$G$ will be called \emph{outward focal} if $g'(\xi) > 1$ for some $g\in G$, and \emph{inward focal} otherwise.
\end{definition}
Note that an inward focal semigroup cannot be a group.

\begin{proposition}
\label{propositionfocalequivalent}
For $G\leq\Isom(X)$, the following are equivalent:
\begin{itemize}
\item[(A)] $G$ is focal.
\item[(B)] $G$ has a unique global fixed point $\xi\in\del X$, and $g'(\xi)\neq 1$ for some $g\in G$.
\item[(C)] $G$ has a unique global fixed pont $\xi\in\del X$, and there are two loxodromic isometries $g,h\in G$ so that $g_+ = h_+ = \xi$, but $g_- \neq h_-$.
\end{itemize}
\end{proposition}
\begin{proof}
The implications (C) \implies (A) \iff (B) are straightforward. Suppose that $G$ is focal, and let $g\in G$ be a loxodromic isometry. Since $G$ is a group, we may without loss of generality suppose that $g'(\xi) < 1$, so that $g_+ = \xi$. Let $j\in G$ be such that $g_-\notin\Fix(j)$. By choosing $n$ sufficiently large, we may guarantee that $(jg^n)'(\xi) < 1$. Then if $h = jg^n$, then $h$ is loxodromic and $h_+ = \xi$. But $g_-\notin \Fix(h)$, so $g_-\neq h_-$.
\end{proof}

\bigskip
\section{Proof of the Classification Theorem}

We begin by recalling the following definition from Section \ref{subsectionshadows}:
\begin{repdefinition}{definitionshadow}
For each $\sigma > 0$ and $x,y\in X$, let
\[
\Shad_y(x,\sigma) = \{\eta\in\del X:\lb y|\eta\rb_x\leq\sigma\}.
\]
We say that $\Shad_y(x,\sigma)$ is the \emph{shadow cast by $x$ from the light source $y$, with parameter $\sigma$}. For shorthand we will write $\Shad(x,\sigma) = \Shad_\zero(x,\sigma)$.
\end{repdefinition}

\begin{lemma}
\label{lemmaloxodromic}
For every $\sigma > 0$, there exists $r > 0$ such that for every $g\in\Isom(X)$ with $\dogo g \geq r$, if there exists a nonempty closed set
\[
Z \subset \Shad_{g^{-1}(\zero)}(\zero,\sigma)
\]
satisfying $g(Z) \subset Z$, then $g$ is loxodromic and $g_+\in Z$.
\end{lemma}

\begin{proof}
Recall from the Bounded Distortion Lemma \ref{lemmaboundeddistortion} that
\begin{equation}
\label{equationloxodromic}
\frac{\Dist(g(y_1),g(y_2))}{\Dist(y_1,y_2)} \leq C b^{-\dogo g} \all y_1,y_2\in Z
\end{equation}
for some $C > 0$ independent of $g$. Now choose $r > 0$ large enough so that $C b^{-r} < 1$. If $g\in\Isom(X)$ satisfies $\dogo g\geq r$, we can conclude that $g:Z\to Z$ is a strict contraction of the complete metametric space $(Z,\Dist)$. Then by Theorem \ref{banachcontractionprinciple}, $g$ has a unique fixed point $\xi\in Z_\refl = Z\cap \del X$.

To complete the proof we must show that $g'(\xi) < 1$ to prove that $g$ is not parabolic and that $\xi = g_+$. Indeed, by the Bounded Distortion Lemma, we have $g'(\xi) \lesssim_\times b^{-\dogo g} \leq b^{-r}$, so choosing $r$ sufficiently large completes the proof.
\end{proof}

\begin{corollary}
\label{corollaryloxodromic}
For every $\sigma > 0$, there exists $r = r_\sigma > 0$ such that for every $g\in\Isom(X)$ with $\dogo g\geq r$, if $g$ is not loxodromic, then
\begin{equation}
\label{nonloxodromic}
\lb g(\zero)|g^{-1}(\zero)\rb_\zero \geq \sigma.
\end{equation}
\end{corollary}
\begin{proof}
Fix $\sigma > 0$, and let $\sigma' = \sigma + \delta$, where $\delta$ is the implied constant in Gromov's inequality. Apply Lemma \ref{lemmaloxodromic} to get $r' > 0$. Let $r = \max(r',2\sigma')$. Now suppose that $g\in\Isom(X)$ satisfies $\dogo g\geq r\geq r'$ but is not loxodromic. Then by Lemma \ref{lemmaloxodromic}, we have
\[
\Shad(g(\zero),\sigma')\butnot \Shad_{g^{-1}(\zero)}(\zero,\sigma')\neq \emptyset.
\]
Let $x$ be a member of this set. By definition this means that
\[
\lb \zero|x\rb_{g(\zero)} \leq \sigma' < \lb g^{-1}(\zero)|x\rb_\zero.
\]
Since $\dogo g\geq r\geq 2\sigma'$, we have
\[
\lb g(\zero)|x\rb_\zero = \dogo g - \lb \zero|x\rb_{g(\zero)} \geq 2\sigma' - \sigma' = \sigma'.
\]
Now by Gromov's inequality we have
\[
\lb g(\zero)|g^{-1}(\zero)\rb_\zero \geq \min(\lb g(\zero)|x\rb_\zero,\lb g^{-1}(\zero)|x\rb_\zero) - \delta
\geq \sigma' - \delta = \sigma.
\]
\end{proof}

\begin{lemma}\label{keylemma}
Let $G\prec\Isom(X)$ be a semigroup which is not loxodromic, and let $(g_n)_1^\infty$ be a sequence in $G$ such that $\dogo{g_n} \to \infty$. Then $(g_n(\zero))_1^\infty$ is a Gromov sequence.
\end{lemma}
\begin{proof}
Fix $\sigma > 0$ large, and let $r = r_\sigma$ be as in Corollary \ref{corollaryloxodromic}. Since $G$ is not loxodromic, \eqref{nonloxodromic} holds for every $g\in G$ for which $\dogo g\geq r$.

Fix $n,m\in\N$ with $\dogo{g_n},\dogo{g_m}\geq r$; Corollary \ref{corollaryloxodromic} gives
\begin{align} \label{gngn1}
\lb g_n(\zero)|g_n^{-1}(\zero)\rb_\zero &\geq \sigma\\ \label{gmgm1}
\lb g_m(\zero)|g_m^{-1}(\zero)\rb_\zero &\geq \sigma.
\end{align}
By contradiction, suppose that $\lb g_n(\zero)|g_m(\zero)\rb_\zero \leq \sigma/2$; then Gromov's inequality together with \eqref{gngn1} gives
\begin{equation}
\label{gn1gm}
\lb g_n^{-1}(\zero)|g_m(\zero)\rb_\zero \asymp_\plus 0.
\end{equation}
It follows that
\[
\dogo{g_n g_m} = \dist(g_n^{-1}(\zero),g_m(\zero)) \geq 2r - \lb g_n^{-1}(\zero)|g_m(\zero)\rb_\zero \asymp_\plus 2r.
\]
Choosing $r$ sufficiently large, we have $\dogo{g_n g_m} \geq r$. So by Corollary \ref{corollaryloxodromic},
\begin{equation}
\label{gngmgm1gn1}
\lb g_n g_m(\zero)|g_m^{-1} g_n^{-1}(\zero)\rb_\zero \geq \sigma.
\end{equation}
Now
\begin{align*}
\lb g_n(\zero)|g_n g_m(\zero)\rb_\zero &=_\pt \lb \zero|g_m(\zero)\rb_{g_n^{-1}(\zero)}\\
&=_\pt \dogo{g_n} - \lb g_n^{-1}(\zero)|g_m(\zero)\rb_\zero\\
&\asymp_\plus \dogo{g_n} \by{\eqref{gn1gm}}\\
&\geq_\pt r,
\end{align*}
i.e.
\begin{equation}
\label{gngngm}
\lb g_n(\zero)|g_n g_m(\zero)\rb_\zero \gtrsim_\plus r.
\end{equation}

A similar argument yields
\begin{equation}
\label{gm1gm1gn1}
\lb g_m^{-1}(\zero)|g_m^{-1} g_n^{-1}(\zero)\rb_\zero \gtrsim_\plus r.
\end{equation}
Combining \eqref{gmgm1}, \eqref{gm1gm1gn1}, \eqref{gngmgm1gn1}, and \eqref{gngngm}, together with Gromov's inequality, yields
\[
\lb g_n|g_m\rb_\zero \gtrsim_\plus \min(\sigma,r).
\]
This completes the proof.
\end{proof}

\begin{proof}[Proof of Theorem \ref{classificationofsemigroups}]
Suppose that $G$ is neither elliptic nor loxodromic, and we will show that it is parabolic. Since $G$ is not elliptic, there is a sequence $(g_n)_1^\infty$ in $G$ such that $\dogo{g_n} \to \infty$. By Lemma \ref{keylemma}, $(g_n(\zero))_1^\infty$ is a Gromov sequence; let $\xi\in\del X$ be the limit point.

Note that $\xi$ is uniquely determined by $G$; if $(h_n(\zero))_1^\infty$ were another Gromov sequence, then we could let
\[
j_n := \begin{cases}
g_{n/2} & n\text{ even}\\
h_{(n - 1)/2} & n\text{ odd}
\end{cases}.
\]
The sequence $(j_n(\zero))_1^\infty$ would tend to infinity, so by Lemma \ref{keylemma} it would be a Gromov sequence. But that exactly means that the Gromov sequences $(g_n(\zero))_1^\infty$ and $(h_n(\zero))_1^\infty$ are equivalent. Moreover, it is easy to see that $\xi$ does not depend on the choice of the basepoint $\zero\in X$.

In particular, the fact that $\xi$ is canonically determined by $G$ implies that $\xi$ is a global fixed point of $G$. To complete the proof, we need to show that $g'(\xi) = 1$ for all $g\in G$. Suppose we have $g\in G$ such that $g'(\xi)\neq 1$. Then $g$ is loxodromic by Lemma \ref{lemmarepellingloxodromic}, a contradiction.
\end{proof}

\bigskip
\section{Discreteness and focal groups}

\begin{proposition}
\label{propositionnonfocal}
Fix $G\leq\Isom(X)$, and suppose that either
\begin{itemize}
\item[(1)] $G$ is strongly discrete,
\item[(2)] $X$ is CAT(-1) and $G$ is moderately discrete, or
\item[(3)] $X$ admits unique geodesic extensions (e.g. $X$ is an algebraic hyperbolic space) and $G$ is weakly discrete.
\end{itemize}
Then $G$ is not focal.
\end{proposition}
\begin{proof}[Strongly discrete case]
Suppose that $G$ is a focal group. Let $\xi\in\del X$ be its global fixed point, and let $g,h\in G$ be as in (C) of Proposition \ref{propositionfocalequivalent}. Since $h^{-n}(\zero)\to h_- \neq \xi$, we have
\[
\lb h^{-n}(\zero)|\xi\rb_\zero \asymp_{\plus,h} 0
\]
and thus
\[
\lb h^n(\zero)|\xi\rb_\zero \asymp_\plus \dogo{h^n} - \lb\zero|\xi\rb_{h^n(\zero)} \asymp_{\plus,h} \dogo{h^n}.
\]
Applying $g$ we have
\[
\lb gh^n(\zero)|\xi\rb_\zero
= \lb h^n(\zero)|\xi\rb_{g^{-1}(\zero)}
\asymp_{\plus,g} \lb h^n(\zero)|\xi\rb_\zero
\asymp_{\plus,h} \dogo{h^n}
\]
and applying Gromov's inequality we have
\[
\lb h^n(\zero) | gh^n(\zero)\rb_\zero \asymp_{\plus,g,h} \dogo{h^n} \asymp_{\plus,g} \dogo{g h^n}.
\]
Now
\begin{align*}
\dogo{h^{-n} gh^n}
&= \dist(h^n(\zero),gh^n(\zero))\\
&= \dogo{h^n} + \dogo{g h^n} - 2\lb h^n(\zero),gh^n(\zero)\rb_\zero
\asymp_{\plus,g,h} 0.
\end{align*}
Since $G$ is strongly discrete, this implies that the collection $\{h^{-n} gh^n:n\in\N\}$ is finite, and so for some $n_1 < n_2$ we have
\[
h^{-n_1} gh^{n_1} = h^{-n_2} gh^{n_2}
\]
or
\[
h^{n_2 - n_1} g = g h^{n_2 - n_1},
\]
i.e. $h^{n_2 - n_1}$ commutes with $g$. But then $h^{n_2 - n_1}(g_-) = g_-$, contradicting that $g_-\neq h_-$. This completes the proof of Proposition \ref{propositionnonfocal}$(1)$.
\end{proof}

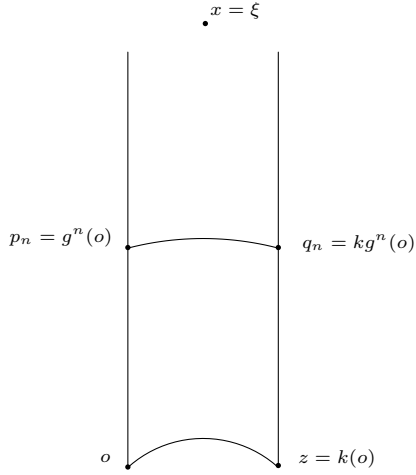
\begin{figure}
\begin{center}
\begin{tikzpicture}[line cap=round,line join=round,>=triangle 45,x=1.0cm,y=1.0cm]
\clip(-5.13,-0.16) rectangle (6.57,8.26);
\draw (0,1.08)-- (0,6.6);
\draw (2,1.1)-- (2,6.6);
\draw [shift={(1,0.14)}] plot[domain=1.32:1.82,variable=\t]({1*3.98*cos(\t r)+0*3.98*sin(\t r)},{0*3.98*cos(\t r)+1*3.98*sin(\t r)});
\draw [shift={(1,-0.03)}] plot[domain=0.84:2.31,variable=\t]({1*1.49*cos(\t r)+0*1.49*sin(\t r)},{0*1.49*cos(\t r)+1*1.49*sin(\t r)});
\draw (-4,0)-- (6,0);
\begin{scriptsize}
\fill [color=black] (2,1.1) circle (1pt);
\draw[color=black] (2.78,1.2) node {$z=k(\zero)$};
\fill [color=black] (1.03,6.98) circle (1pt);
\draw[color=black] (1.43,7.15) node {$x=\xi$};
\fill [color=black] (0,4) circle (1pt);
\draw[color=black] (-0.89,4.14) node {$p_n = g^n(\zero)$};
\fill [color=black] (2,4) circle (1pt);
\draw[color=black] (3.07,4.06) node {$q_n = kg^n(\zero)$};
\fill [color=black] (0,1.08) circle (1pt);
\draw[color=black] (-0.3,1.2) node {$\zero$};
\end{scriptsize}
\end{tikzpicture}
\caption[High altitude implies small displacement in the half-space model]{The higher the point $g^n(\zero)$ is, the smaller its displacement under $k$ is.}
\label{figurexyzpq}
\end{center}
\end{figure}

\begin{proof}[Moderately discrete case]
Suppose that $G$ is a focal group. Let $\xi\in\del X$ be its global fixed point, and let $g,h\in G$ be as in (C) of Proposition \ref{propositionfocalequivalent}. Let
\[
k = [g,h] = g^{-1} h^{-1} gh\in G.
\]
We observe first that
\begin{equation}
\label{kprimeequalsone}
k'(\xi) = \frac{1}{g'(\xi)}\frac{1}{h'(\xi)}g'(\xi) h'(\xi) = 1.
\end{equation}
Note that strong hyperbolicity is necessary to deduce equality in \eqref{kprimeequalsone} rather than merely a coarse asymptotic.

Next, we claim that $k(g_-)\neq g_-$. Indeed, $g_-\notin\Fix(h)$, so $h(g_-)\neq g_-$. This in turn implies that $h(g_-)\notin\Fix(g)$, so $gh(g_-)\neq h(g_-)$. Now applying $g^{-1} h^{-1}$ to both sides shows that $k(g_-)\neq g_-$.

\begin{claim}
\label{claimfocalMD}
$g^{-n} k g^n(\zero)\to \zero$.
\end{claim}
\begin{subproof}
Indeed,
\[
\dogo{g^{-n} k g^n} = \dist(g^n(\zero),k g^n(\zero)).
\]
Let
\begin{align*}
x &= \xi\\
y &= \zero\\
z &= k(\zero)\\
p_n &= g^n(\zero)\\
q_n &= k g^n(\zero).
\end{align*}
(See Figure \ref{figurexyzpq}.) Then $p_n,q_n\in\Delta := \Delta(x,y,z)$. By Proposition \ref{propositionidealCAT} $\dist(p_n,q_n) \leq \dist(\wbar p_n,\wbar q_n)$, where $\wbar p_n,\wbar q_n$ are comparison points for $p_n,q_n$ on the comparison triangle $\wbar \Delta = \Delta(\wbar x,\wbar y,\wbar z)$. Now notice that
\[
\busemann_{\wbar x}(\wbar p_n,\wbar q_n) = \busemann_\xi(g^n(\zero),k g^n(\zero)) = 0
\]
by Proposition \ref{propositionbusemannpreserved} and \eqref{kprimeequalsone}. On the other hand, $\wbar p_n,\wbar q_n\to \wbar x$. An easy calculation based on \eqref{distE} and Proposition \ref{propositionbusemannE} (letting $\wbar x = \infty$) shows that $\dist(\wbar p_n,\wbar q_n)\to 0$, and thus that $\dogo{g^{-n} k g^n}\to 0$ i.e. $g^{-n} k g^n(\zero)\to \zero$.
\end{subproof}
Since $G$ is moderately discrete, this implies that the collection $\{g^{-n} k g^n:n\in\N\}$ is finite. As before (in the proof of the strongly discrete case), this implies that $g^n$ and $k$ commute for some $n\in\N$. But $(g^n)_- = g_-$, and $k(g_-) \neq g_-$, which contradicts that $g^n$ and $k$ commute. This completes the proof of Proposition \ref{propositionnonfocal}$(2)$.
\end{proof}

\begin{proof}[Weakly discrete case]
Suppose that $G$ is a focal group. Let $\xi$, $g$, $h$, and $k$ be as above. Without loss of generality, supposet that $\zero\in \geo{g_-}{\xi}$.

\begin{claim}
$g^{-n} k g^n(\zero)\neq \zero$ for all $n\in\N$.
\end{claim}
\begin{subproof}
Fix $n\in\N$. As observed above, $k(g_-)\neq g_-$. On the other hand, $k(\xi) = \xi$, and $g^n(\zero)\in \geo{g_-}{\xi}$. Since $X$ admits unique geodesic extensions, it follows that $k(g^n(\zero))\neq g^n(\zero)$, or equivalently that $g^{-n} k g^n(\zero)\neq \zero$.
\end{subproof}
Together with Claim \ref{claimfocalMD}, this contradicts that $G$ is weakly discrete. This completes the proof of Proposition \ref{propositionnonfocal}$(3)$.
\end{proof}


\draftnewpage
\chapter{Limit sets} \label{sectionlimitsets}

Throughout this chapter, we fix a subsemigroup $G\prec\Isom(X)$. We define the \emph{limit set} of $G$, along with various subsets. We then define several concepts in terms of the limit set including elementariness and compact type, while relating other concepts to the limit set, such as the quasiconvex core and irreducibility of a group action. We also prove that the limit set is minimal in an approprate sense (Proposition \ref{propositionminimal} - Proposition \ref{propositionminimal2}).

\bigskip
\section{Modes of convergence to the boundary}
\label{subsectionmodes}
We recall (Observation \ref{observationboundaryconvergence}) that a sequence $(x_n)_1^\infty$ in $X$ converges to a point $\xi\in\del X$ if and only if
\[
\lb x_n|\xi\rb_\zero \tendsto n \infty.
\]
In this section we define more restricted modes of convergence. To get an intuition let us consider the case where $X = \E = \E^\alpha$ is the half-space model of a real hyperbolic space. Consider a sequence $(\xx_n)_1^\infty$ in $\E$ which converges to a point $\xi\in\BB := \del\E\butnot\{\infty\} = \HH^{\alpha - 1}$. We say that $\xx_n \to \xi$ \emph{conically} if there exists $\theta > 0$ such that if we let
\[
C(\xi,\theta) = \{\xx\in\E: x_1 \geq \sin(\theta)\|\xx - \xi\|\}
\]
then $\xx_n\in C(\xi,\theta)$ for all $n\in\N$. We call $C(\xi,\theta)$ the \emph{cone centered at $\xi$ with inclination $\theta$}; see Figure \ref{figurecone}.

\definecolor{qqwuqq}{rgb}{0,0.39,0}
\begin{figure}
\begin{center}
\begin{tikzpicture}[line cap=round,line join=round,>=triangle 45,x=1.0cm,y=1.0cm]
\clip(-5.69,-0.26) rectangle (5.35,3.56);
\draw [shift={(0,0)},color=qqwuqq,fill=qqwuqq,fill opacity=0.1] (0,0) -- (147.02:0.44) arc (147.02:180:0.44) -- cycle;
\draw (0,0)-- (4.26,2.76);
\draw (-4.26,2.76)-- (0,0);
\draw (-5,0)-- (5,0);
\begin{scriptsize}
\fill [color=black] (0,0) circle (1.2pt);
\draw[color=black] (0.09,-0.17) node {$\xi$};
\draw[color=black] (-0.60,0.18) node {$\theta$};
\fill [color=black] (-1.62,2.69) circle (0.81pt);
\fill [color=black] (0.35,1.2) circle (0.81pt);
\fill [color=black] (-1.67,1.77) circle (0.81pt);
\fill [color=black] (-0.26,1.24) circle (0.81pt);
\fill [color=black] (0.04,0.64) circle (0.81pt);
\draw[color=black] (0.35,2.2) node {$C(\xi,\theta)$};
\end{scriptsize}
\end{tikzpicture}
\caption[Conical convergence to the boundary]{A sequence converging conically to $\xi$. For each point $\xx$, the height of $\xx$ is greater than $\sin(\theta)$ times the distance from $\xx$ to $\xi$.}
\label{figurecone}
\end{center}
\end{figure}
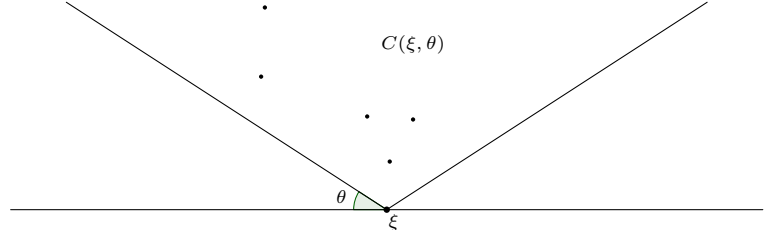


\begin{proposition}
\label{propositionradialconvergence}
Let $(\xx_n)_1^\infty$ be a sequence in $\E$ converging to a point $\xi\in\BB$. Then the following are equivalent:
\begin{itemize}
\item[(A)] $(\xx_n)_1^\infty$ converges conically to $\xi$.
\item[(B)] The sequence $(\xx_n)_1^\infty$ lies within a bounded distance of the geodesic ray $\geo\zero\xi$.
\item[(C)] There exists $\sigma > 0$ such that for all $n\in\N$,
\[
\lb\zero|\xi\rb_{\xx_n}\leq\sigma,
\]
or equivalently
\begin{equation}
\label{sigmaradially}
\xi\in\Shad(\xx_n,\sigma).
\end{equation}
\end{itemize}
Moreover, the equivalence of (B) and (C) holds in all geodesic hyperbolic metric spaces.
\end{proposition}
\begin{proof}
The equivalence of (B) and (C) follows directly from (i) of Proposition \ref{propositionrips}. Moreover, conditions (B) and (C) are clearly independent of the basepoint $\zero$. Thus, in proving (A) \iff (B) we may without loss of generality suppose that $\xi = \0$ and $\zero = (1,\0)$. Note that if $\theta > 0$ is fixed, then
\[
C(\0,\theta) = \{\xx\in\E : \ang(\xx) \leq \pi/2 - \theta\}
= \{\xx\in\E : \theta(\xx) \leq -\log\cos(\pi/2 - \theta)\},
\]
where $\theta = \theta_{\0,\infty,\zero}$ is as in Proposition \ref{propositionrtheta}. Since $-\log\cos(\pi/2 - \theta) \to \infty$ as $\theta \to 0$, we have (A) if and only if the sequence $(\theta(\xx_n))_1^\infty$ is bounded. But
\begin{align*}
\theta(\xx_n) = \lb \0 | \infty \rb_{\xx_n}
&\asymp_\plus \dist(\xx_n,\geo\0\infty) \by{(i) of Proposition \ref{propositionrips}}\\
&=_\pt \dist(\xx_n,\geo\zero\0), &\text{(for $n$ sufficiently large)}
\end{align*}
which completes the proof.
\end{proof}

Condition (B) of Proposition \ref{propositionradialconvergence} motivates calling this kind of convergence \emph{radial}; we shall use this terminology henceforth. However, condition (C) is best suited to a general hyperbolic metric space.

\begin{definition}
\label{definitionradialconvergence}
Let $(x_n)_1^\infty$ be a sequence in $X$ converging to a point $\xi\in\del X$. We will say that $(x_n)_1^\infty$ converges to $\xi$
\begin{itemize}
\item \emph{$\sigma$-radially} if \eqref{sigmaradially} holds for all $n\in\N$,
\item \emph{radially} if it converges $\sigma$-radially for some $\sigma > 0$,
\item \emph{$\sigma$-uniformly radially} if it converges $\sigma$-radially, $x_1 = \zero$, and
\[
\dist(x_n,x_{n + 1})\leq \sigma \all n\in\N,
\]
\item \emph{uniformly radially} if it converges $\sigma$-uniformly radially for some $\sigma > 0$.
\end{itemize}
\end{definition}
Note that a sequence can converge $\sigma$-radially and uniformly radially without converging $\sigma$-uniformly radially.

We next define horospherical convergence. Again, we motivate the discussion by considering the case of a real hyperbolic space $X = \E = \E^\alpha$. This time, however, we will let $\xi = \infty$, and we will say that a sequence $(\xx_n)_1^\infty$ converges \emph{horospherically} to $\xi$ if
\[
\height(\xx_n) \tendsto n \infty,
\]
where the \emph{height} of a point $\xx\in\E$ is its first coordinate $x_1$. This terminology comes from defining a \emph{horoball centered at $\infty$} to be a set of the form $H_{\infty,t} = \{\xx : \height(\xx) > e^t\}$; then $\xx_n\to\infty$ horospherically if and only if for every horoball $H_{\infty,t}$ centered at infinity, we have $\xx_n\in H_{\infty,t}$ for all sufficiently large $n$. (See also Definition \ref{definitionhoroball} below.)

Recalling (cf. Proposition \ref{propositionbusemannE}) that
\[
\height(x) = b^{\busemann_\infty(\zero,x)},
\]
the above discussion motivates the following definition:
\begin{definition}
\label{definitionhorosphericalconvergence}
A sequence $(x_n)_1^\infty$ in $X$ converges \emph{horospherically} to a point to $\xi\in\del X$ if
\[
\busemann_\xi(\zero,x_n)\tendsto n +\infty.
\]
\end{definition}

\begin{observation}
If $x_n\to\xi$ radially, then $x_n\to\xi$ horospherically.
\end{observation}
\begin{proof}
Indeed,
\[
\busemann_\xi(\zero,x_n) \asymp_\plus \dox{x_n} - 2\lb\zero|\xi\rb_{x_n} \asymp_\plus \dox{x_n} \tendsto n \infty.
\]
\end{proof}

\begin{figure}
\begin{center}
\begin{tikzpicture}[line cap=round,line join=round,>=triangle 45,x=1.0cm,y=1.0cm]
\clip(-5.69,-0.39) rectangle (5.35,3.76);
\draw (0,0)-- (4.26,2.76);
\draw (-4.26,2.76)-- (0,0);
\draw (-5,0)-- (5,0);
\draw(0.03,1.86) circle (1.86cm);
\begin{scriptsize}
\fill [color=black] (0,0) circle (1.2pt);
\draw[color=black] (0.12,-0.26) node {$\xi$};
\fill [color=black] (-1.61,2.26) circle (0.81pt);
\fill [color=black] (-0.91,0.42) circle (0.81pt);
\fill [color=black] (-1.61,1.39) circle (0.81pt);
\fill [color=black] (-1.32,0.75) circle (0.81pt);
\fill [color=black] (-0.5,0.19) circle (0.81pt);
\end{scriptsize}
\end{tikzpicture}
\caption[Converging horospherically but not radially to the boundary]{A sequence converging horospherically but not radially to $\xi$.}
\label{figurehorospherical}
\end{center}
\end{figure}
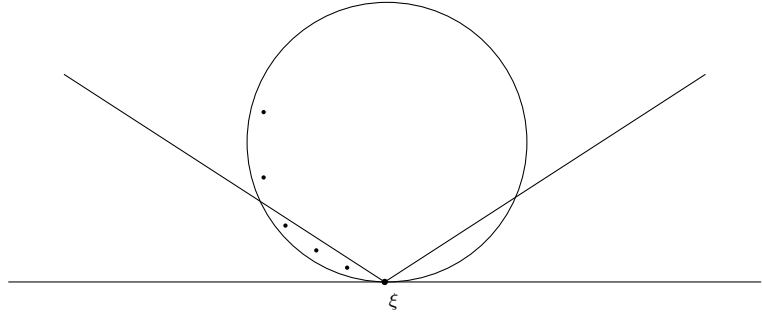

The converse is false, as illustrated in Figure \ref{figurehorospherical}.

\begin{observation}
\label{observationdependenceonbasepoint}
The concepts of convergence, radial convergence, uniformly radial convergence, and horospherical convergence are independent of the basepoint $\zero$, whereas the concepts of $\sigma$-radial convergence and $\sigma$-uniformly radial convergence depend on the basepoint. (Regarding $\sigma$-radial convergence, this dependence on basepoint is not too severe; see Proposition \ref{propositionnearinvariance} below.)
\end{observation}

\bigskip
\section{Limit sets}

We define the limit set of $G$, a subset of $\del X$ which encodes geometric information about $G$. We also define a few important subsets of the limit set.

\begin{definition}
\label{definitionlimitset}
Let
\begin{align*}
\Lambda(G) &:= \{\eta\in\del X: g_n(\zero)\to \eta \text{ for some $(g_n)_1^\infty \in G^\N$}\}\\
\Lr(G) &:= \{\eta\in\del X: g_n(\zero)\to \eta \text{ radially for some $(g_n)_1^\infty \in G^\N$}\}\\
\Lrsigma(G) &:= \{\eta\in\del X: g_n(\zero)\to \eta \text{ $\sigma$-radially for some $(g_n)_1^\infty \in G^\N$}\}\\
\Lur(G) &:= \{\eta\in\del X: g_n(\zero)\to \eta \text{ uniformly radially for some $(g_n)_1^\infty \in G^\N$}\}\\
\Lursigma(G) &:= \{\eta\in\del X: g_n(\zero)\to \eta \text{ $\sigma$-uniformly radially for some $(g_n)_1^\infty \in G^\N$}\}\\
\Lh(G) &:= \{\eta\in\del X: g_n(\zero)\to \eta \text{ horospherically for some $(g_n)_1^\infty \in G^\N$}\}.
\end{align*}
These sets are respectively called the \emph{limit set}, \emph{radial limit set}, \emph{$\sigma$-radial limit set}, \emph{uniformly radial limit set}, \emph{$\sigma$-uniformly radial limit set}, and \emph{horospherical limit set} of the semigroup $G$.
\end{definition}
Note that
\begin{align*}
\Lr &= \bigcup_{\sigma > 0}\Lrsigma\\
\Lur &= \bigcup_{\sigma > 0}\Lursigma\\
\Lur &\subset \Lr \subset \Lh \subset \Lambda.
\end{align*}

\begin{observation}
\label{observationlimitsets}
The sets $\Lambda$, $\Lr$, $\Lur$, and $\Lh$ are invariant\Footnote{By \emph{invariant} we always mean \emph{forward invariant}.} under the action of $G$, and are independent of the basepoint $\zero$. The set $\Lambda$ is closed.
\end{observation}
\begin{proof}
The first assertion follows from Observation \ref{observationdependenceonbasepoint} and the second follows directly from the definition of $\Lambda$ as the intersection of $\del X$ with the set of accumulation points of the set $G(\zero)$.
\end{proof}

\begin{proposition}[Near-invariance of the sets $\Lrsigma$]
\label{propositionnearinvariance}
For every $\sigma > 0$, there exists $\tau > 0$ such that for every $g\in G$, we have
\begin{equation}
\label{nearinvariance}
g(\Lrsigma) \subset \Lrtau.
\end{equation}
If $X$ is strongly hyperbolic, then \eqref{nearinvariance} holds for all $\tau > \sigma$.
\end{proposition}
\begin{proof}
Fix $\xi\in\Lrsigma$. There exists a sequence $(h_n)_1^\infty$ so that $h_n(\zero)\to \xi$ $\sigma$-radially, i.e.
\[
\lb \zero|\xi\rb_{h_n(\zero)} \leq \sigma\all n\in\N
\]
and $h_n(\zero)\to \xi$. Now
\[
\lb \zero|g^{-1}(\zero)\rb_{h_n(\zero)} \geq \dogo{h_n} - \dogo{g^{-1}} \tendsto n +\infty.
\]
Thus, for $n$ sufficiently large, Gromov's inequality gives
\begin{equation}
\label{ghnsigma}
\lb g^{-1}(\zero)|\xi\rb_{h_n(\zero)} \lesssim_\plus \sigma
\end{equation}
i.e.
\[
\lb \zero|g(\xi)\rb_{gh_n(\zero)} \lesssim_\plus \sigma.
\]
So $gh_n(\zero)\to g(\xi)$ $\tau$-radially, where $\tau$ is equal to $\sigma$ plus the implied constant of this asymptotic. Thus, $g(\xi)\in\Lrtau$.

If $X$ is strongly hyperbolic, then by using \eqref{gromovforCAT} instead of Gromov's inequality, the implied constant of \eqref{ghnsigma} can be made arbitrarily small. Thus $\tau$ may be taken arbitrarily close to $\sigma$.
\end{proof}

\bigskip
\section{Cardinality of the limit set}

In this section we characterize the cardinality of the limit set according to the classification of the semigroup $G$.

\begin{proposition}[Cardinality of the limit set by classification]
\label{propositioncardinalitylimitset}
Fix $G\prec\Isom(X)$.
\begin{itemize}
\item[(i)] If $G$ is elliptic, then $\LambdaG = \emptyset$.
\item[(ii)] If $G$ is parabolic or inward focal with global fixed point $\xi$, then $\LambdaG = \{\xi\}$.
\item[(iii)] If $G$ is lineal with fixed pair $\{\xi_1,\xi_2\}$, then $\LambdaG \subset \{\xi_1,\xi_2\}$, with equality if $G$ is a group.
\item[(iv)] If $G$ is outward focal or of general type, then $\#(\LambdaG) \geq \#(\R)$. Equality holds if $X$ is separable.
\end{itemize}
\end{proposition}
Case (i) is immediate, while case (iv) requires the theory of Schottky groups and will be proven in Chapter \ref{sectionschottky} (see Proposition \ref{propositionnonelementaryequivalent}).
\begin{proof}[Proof of \text{(ii)}]
For $g\in G$, $g'(\xi) \leq 1$, so by Proposition \ref{propositionbusemannpreserved}, we have $\busemann_\xi(g(\zero),\zero) \lesssim_\plus 0$. In particular, by (h) of Proposition \ref{propositionbasicidentities} we have
\[
\lb x | \xi \rb_\zero \gtrsim_\plus \frac12 \dox x \all x\in G(\zero).
\]
This implies that $x_n\to \xi$ for any sequence $(x_n)_1^\infty$ in $G(\zero)$ satisfying $\dox{x_n}\to \infty$. It follows that $\LambdaG = \{\xi\}$.
\end{proof}
\begin{proof}[Proof of \text{(iii)}]
By Lemma \ref{lemmatranslationdistance} we have
\[
\theta(g(\zero)) \asymp_\plus \theta(\zero) = \zero \all g\in G,
\]
where $\theta = \theta_{\xi_1,\xi_2,\zero} = \theta_{\xi_2,\xi_1,\zero}$ is as in Section \ref{subsectionpolar}. Thus
\[
\lb \xi_1 | \xi_2 \rb_x \asymp_\plus 0 \all x\in G(\zero).
\]
Fix a sequence $G(\zero)\ni x_n \to \xi\in\LambdaG$. By Gromov's inequality, there exists $i = 1,2$ such that
\[
\lb \zero | \xi_i \rb_{x_n} \asymp_\plus 0 \text{ for infinitely many $n$.}
\]
It follows that $x_n \to \xi_i$ radially along some subsequence, and in particular $\xi = \xi_i$. Thus $\LambdaG \subset \{\xi_1,\xi_2\}$.
\end{proof}

\begin{definition}
\label{definitionelementary}
Fix $G\prec\Isom(X)$. $G$ is called \emph{elementary} if $\#(\Lambda) < \infty$ and \emph{nonelementary} if $\#(\Lambda) = \infty$.
\end{definition}

Thus, according to Proposition \ref{propositioncardinalitylimitset}, elliptic, parabolic, lineal, and inward focal semigroups are elementary while outward focal semigroups and semigroups of general type are nonelementary.

\begin{remark}
In the Standard Case, some authors (e.g. \cite[\65.5]{Ratcliffe}) define a subgroup of $\Isom(X)$ to be elementary if there is a global fixed point or a global fixed geodesic line. According to this definition, focal groups are considered elementary. By contrast, we follow \cite{CCMT} and others in considering them to be nonelementary.

Another common definition in the Standard Case is that a group is elementary if it is virtually abelian. This agrees with our definition, but beyond the Standard Case this equivalence no longer holds (cf. Observation \ref{observationmargulislemma} and Remark \ref{remarkmargulislemma}).
\end{remark}

\bigskip
\section{Minimality of the limit set}

Observation \ref{observationlimitsets} identified the limit set $\Lambda$ as a closed $G$-invariant subset of the Gromov boundary $\del X$. In this section, we give a characterization of $\Lambda$ depending on the classification of $G$.

\begin{proposition}[Cf. {\cite[Th\'eor\`eme 5.1]{Coornaert}}]
\label{propositionminimal}
Fix $G\prec\Isom(X)$. Then any closed $G$-invariant subset of $\del X$ containing at least two points contains $\Lambda$.
\end{proposition}
\begin{proof}
We begin with the following lemma, which will also be useful later:
\begin{lemma}
\label{lemmaminimal}
Let $(x_n)_1^\infty$, $(y_n^{(1)})_1^\infty$, $(y_n^{(2)})_1^\infty$ be sequences in $\bord X$ satisfying
\[
\lb y_n^{(1)} | y_n^{(2)} \rb_{x_n} \asymp_\plus 0
\]
and
\[
x_n \to \xi\in\del X.
\]
Then $\xi\in\cl{\{y_n^{(i)} : n\in\N , i = 1,2\}}$.
\end{lemma}
\begin{subproof}
For $n\in\N$ fixed, by Gromov's inequality there exists $i_n = 1,2$ such that
\[
\lb \zero | y_n^{(i_n)} \rb_{x_n} \asymp_\plus 0.
\]
It follows that
\[
\lb x_n | y_n^{(i_n)}\rb_\zero \asymp_\plus \dox{x_n} \tendsto n \infty.
\]
On the other hand
\[
\lb x_n | \xi\rb_\zero \tendsto n \infty,
\]
so by Gromov's inequality
\[
\lb y_n^{(i_n)} | \xi \rb_\zero \tendsto n \infty,
\]
i.e. $y_n^{(i_n)} \to \xi$.
\end{subproof}
Now let $F$ be a closed $G$-invariant subset of $\del X$ containing two points $\xi_1\neq\xi_2$, and let $\eta\in\Lambda$. Then there exists a sequence $(g_n)_1^\infty$ so that $g_n(\zero)\to\eta$. Applying Lemma \ref{lemmaminimal} with $x_n = g_n(\zero)$ and $y_n^{(i)} = g_n(\xi_i)\in F$ completes the proof.

\end{proof}
The proof of Proposition \ref{propositionminimal} may be compared to the proof of \cite[Theorem 3.1]{FSU4}, where a quantitative convergence result is proven assuming that $\eta$ is in the radial limit set (and assuming that $G$ is a group).

\begin{corollary}
\label{corollaryminimal}
Let $G\prec\Isom(X)$ be nonelementary.
\begin{itemize}
\item[(i)] If $G$ is outward focal with global fixed point $\xi$, then $\LambdaG$ is the smallest closed $G$-invariant subset of $\del X$ which contains a point other than $\xi$.
\item[(ii)] (Cf. {\cite[Theorem 5.3.7]{Beardon_book}}) If $G$ is of general type, then $\LambdaG$ is the smallest nonempty closed $G$-invariant subset of $\del X$.
\end{itemize}
\end{corollary}
\begin{proof}
Any $G$-invariant set containing a point which is not fixed by $G$ contains two points.
\end{proof}

\begin{corollary}
Let $G\prec\Isom(X)$ be nonelementary. Then
\[
\Lambda = \cl{\Lr} = \cl{\Lur}.
\]
\end{corollary}
\begin{proof}
The implications $\supset$ are clear. On the other hand, for each loxodromic $g\in G$ we have $g_+\in\Lur$. Thus $\Lur\neq\emptyset$, and $\Lur\nsubset\{\xi\}$ if $G$ is outward focal with global fixed point $\xi$. By Proposition \ref{propositioncardinalitylimitset}, $G$ is either outward focal or of general type. Applying Corollary \ref{corollaryminimal}, we have $\cl{\Lur}\supset\LambdaG$.
\end{proof}

\begin{remark}
If $G$ is elementary, it is easily verified that $\Lambda = \Lr = \Lur$ unless $G$ is parabolic, in which case $\Lr = \Lur = \emptyset \propersubset \Lambda$.
\end{remark}

If $G$ is a nonelementary group, then Corollary \ref{corollaryminimal} immediately implies that the set of loxodromic fixed points of $G$ is dense in $\Lambda$. However, if $G$ is not a group then this conclusion does not follow, since the set of attracting loxodromic fixed points is not necessarily $G$-invariant. (The set of attracting fixed points is the right set to consider, since the set of repelling fixed points is not necessarily a subset of $\Lambda$.) Nevertheless, we have the following:

\begin{proposition}
\label{propositionminimal2}
Let $G\prec\Isom(X)$ be nonelementary. Then the set
\[
\Lambda_+ := \{g_+ : g\in G \text{ is loxodromic}\}.
\]
is dense in $\Lambda$.
\end{proposition}
\begin{proof}
First note that it suffices to show that $\cl{\Lambda_+}$ contains all elements of $\Lambda$ which are not global fixed points. Indeed, if this is true, then $\cl{\Lambda_+}$ is $G$-invariant, and applying Corollary \ref{corollaryminimal} completes the proof.

Fix $\xi\in\Lambda$ which is not a global fixed point of $G$, and choose $h\in G$ such that $h(\xi)\neq \xi$. Fix $\epsilon > 0$ small enough so that $\Dist(B,h(B)) > \epsilon$, where $B = B(\xi,\epsilon)$. Let $\sigma > 0$ be large enough so that the Big Shadows Lemma \ref{lemmabigshadow} holds. Since $\xi\in\Lambda$, there exists $g\in G$ such that
\[
\Shad(g(\zero),\sigma) \subset B.
\]
Let $Z = g^{-1}(\Shad(g(\zero),\sigma)) = \Shad_{g^{-1}(\zero)}(\zero,\sigma)$. Then by Lemma \ref{lemmabigshadow}, $\Diam(\del X\butnot Z)\leq\epsilon$. Thus $\del X\butnot Z$ can intersect at most one of the sets $B$, $h(B)$. So $B\subset Z$ or $h(B)\subset Z$. If $B\subset Z$ then
\[
g(B) \subset B \text{ and } B \subset \Shad_{g^{-1}(\zero)}(\zero,\sigma),
\]
whereas if $h(B)\subset Z$ then
\[
gh(B) \subset B \text{ and } B \subset \Shad_{(gh)^{-1}(\zero)}(\zero,\sigma + \dogo h).
\]
So by Lemma \ref{lemmaloxodromic}, we have $j_+\in B$, where $j = g$ or $j = gh$ is a loxodromic isometry.
\end{proof}

The following improvement over Proposition \ref{propositionminimal2} has a quite intricate proof:

\begin{proposition}[Cf. {\cite[Theorem 5.3.8]{Beardon_book}}, {\cite[p.349]{Koebe}}]
\label{propositionminimal3}
Let $G\prec\Isom(X)$ be of general type. Then
\[
\{(g_+,g_-) : g\in G \text{ is loxodromic}\}
\]
is dense in $\Lambda(G)\times\Lambda(G^{-1})$. Here $G^{-1} = \{g^{-1} : g\in G\}$.
\end{proposition}
\begin{proof}

\begin{claim}
\label{claimi04}
Let $g$ be a loxodromic isometry and fix $\epsilon > 0$. There exists $\delta = \delta(\epsilon,g)$ such that for all $\xi_1,\xi_2\in\del X$ with $\Dist(\xi_2,\Fix(g)) \geq \epsilon$,
\[
\#\{i = 0,\ldots,4 : \Dist(g^i(\xi_1),\xi_2) \leq \delta\} \leq 1.
\]
\end{claim}
\begin{subproof}
Suppose that $\Dist(g^i(\xi_1),\xi_2) \leq \delta$ for two distinct values of $i$. Then $\Dist(g^{i_1}(\xi_1),g^{i_2}(\xi_1)) \leq 2\delta$. For every $n$, we have
\[
\Dist(g^{n + i_1}(\xi_1),g^{n + i_2}(\xi_1)) \lesssim_\times b^{\dogo{g^n}}\delta
\]
and thus by the triangle inequality
\[
\Dist(g^{i_1}(\xi_1),g^{n(i_2 - i_1) + i_1}(\xi_1)) \lesssim_{\times,n} \delta.
\]
By Theorem \ref{theoremloxodromicextra}, if $n$ is sufficiently large then $\Dist(g^{n(i_2 - i_1) + i_1}(\xi_1), g_+) \leq \epsilon/2$, which implies that
\[
\epsilon/2 \leq \Dist(\xi_2,\Fix(g)) - \Dist(g^{n(i_2 - i_1) + i_1}(\xi_1), g_+) \leq \Dist(\xi_2, g^{n(i_2 - i_1) + i_1}(\xi_1)) \lesssim_{\times,n} \delta,
\]
which is a lower bound on $\delta$ independent of $\xi_1,\xi_2$. Choosing $\delta$ less than this lower bound yields a contradiction.
\end{subproof}

\begin{claim}
\label{claimxi1234}
There exist $\epsilon,\rho > 0$ such that for all $\xi_1,\xi_2,\xi_3,\xi_4\in\Lambda$, there exists $j\in G$ such that
\begin{equation}
\label{xi1234}
\Dist(j(\xi_k),\xi_\ell) \geq \epsilon \all k = 1,2 \all \ell = 3,4 \text{ and } \dogo j \leq \rho.
\end{equation}
\end{claim}
\begin{subproof}
Fix $g,h\in G$ loxodromic with $\Fix(g)\cap\Fix(h) = \emptyset$, and let
\[
\rho = \max_{i = 0}^4 \max_{j = 0}^4 \dogo{g^i h^j}.
\]
Now fix $\xi_1,\xi_2,\xi_3,\xi_4\in\Lambda$. By Claim \ref{claimi04}, for each $k = 1,2$ and $\eta\in\Fix(g)$, we have
\[
\#\{j = 0,\ldots,4 : \Dist(h^j(\xi_k),\eta) \leq \delta_1 := \delta(\Dist(\Fix(g),\Fix(h)),h)\} \leq 1.
\]
It follows that there exists $j\in\{0,\ldots,4\}$ such that $\Dist(h^j(\xi_k),\eta) \geq \delta_1$ for all $k = 1,2$ and $\eta\in\Fix(g)$. Applying Claim \ref{claimi04} again, we see that for each $k = 1,2$ and $\ell = 3,4$, we have
\[
\#\{i = 0,\ldots,4 : \Dist(g^{-i}(\xi_\ell),h^j(\xi_k)) \leq \delta_2 := \delta(\delta_1,g^{-1})\} \leq 1.
\]
It follows that there exists $i\in\{0,\ldots,4\}$ such that $\Dist(g^{-i}(\xi_\ell),h^j(\xi_k)) \geq \delta_2$ for all $k = 1,2$ and $\ell = 3,4$. But then
\[
\Dist(g^i h^j(\xi_k) , \xi_\ell) \gtrsim_\times \delta_2,
\]
completing the proof.
\end{subproof}

Now fix $\xi_+\in\Lambda$, $\xi_-\in\Lambda(G^{-1})$ distinct, and fix $\delta > 0$ arbitrarily small. By the definition of $\Lambda$, there exist $g,h\in G$ such that
\[
\Dist(g(\zero),\xi_+) , \Dist(h^{-1}(\zero),\xi_-) \leq \delta.
\]
Let $\sigma > 0$ be large enough so that the Big Shadows Lemma \ref{lemmabigshadow} holds for $\epsilon = \epsilon/2$, where $\epsilon$ is as in Claim \ref{claimxi1234}. Then
\[
\Diam(\del X\butnot\Shad_{g^{-1}(\zero)}(\zero,\sigma)) , \Diam(\del X\butnot \Shad_{h(\zero)}(\zero,\sigma)) \leq \epsilon/2.
\]
On the other hand,
\[
\Diam(\Shad(g(\zero),\sigma)) , \Diam(\Shad(h^{-1}(\zero),\sigma)) \lesssim_\times \delta.
\]
Since $\Shad(g(\zero),\sigma)$ is far from $h^{-1}(\zero)$ and $\Shad(h^{-1}(\zero),\sigma)$ is far from $g(\zero)$, the Bounded Distortion Lemma \ref{lemmaboundeddistortion} gives
\[
\Diam(h(\Shad(g(\zero),\sigma))) , \Diam(g^{-1}(\Shad(h^{-1}(\zero),\sigma))) \lesssim_\times \delta.
\]
Choose $\xi_1\in h(\Shad(g(\zero),\sigma))$, $\xi_2\in g^{-1}(\Shad(h^{-1}(\zero),\sigma))$, $\xi_3\in \del X\butnot\Shad_{g^{-1}(\zero)}(\zero,\sigma)$ and $\xi_4\in \del X\butnot\Shad_{h(\zero)}(\zero,\sigma)$. By Claim \ref{claimxi1234}, there exists $j\in G$ such that \eqref{xi1234} holds. Then
\[
\Diam(jh(\Shad(g(\zero),\sigma))) , \Diam(j^{-1} g^{-1}(\Shad(h^{-1}(\zero),\sigma))) \lesssim_\times \delta,
\]
and by choosing $\delta$ sufficiently small, we can make these diameters less than $\epsilon/2$. It follows that
\begin{align*}
jh(\Shad(g(\zero),\sigma)) &\subset \Shad_{g^{-1}(\zero)}(\zero,\sigma), \text{ and }\\ 
j^{-1} g^{-1}(\Shad(h^{-1}(\zero),\sigma)) &\subset \Shad_{h(\zero)}(\zero,\sigma)
\end{align*}
or equivalently that
\begin{align*}
gjh(\Shad(g(\zero),\sigma)) &\subset \Shad(g(\zero),\sigma), \text{ and }\\ 
(gjh)^{-1}(\Shad(h^{-1}(\zero),\sigma)) &\subset \Shad(h^{-1}(\zero),\sigma).
\end{align*}
By Lemma \ref{lemmaloxodromic}, it follows that $gjh$ is a loxodromic isometry satisfying
\[
(gjh)_+\in\Shad(g(\zero),\sigma) ,\; (gjh)_-\in\Shad(h^{-1}(\zero),\sigma).
\]
In particular $\Dist((gjh)_+,\xi_+) , \Dist((gjh)_-,\xi_-) \lesssim_\times \delta$. Since $\delta$ was arbitrary, this completes the proof.
\end{proof}

\bigskip
\section{Convex hulls}
\label{subsectionconvexhulls}

In this section, we assume that $X$ is regularly geodesic (see Section \ref{subsectiongeodesicsCAT}). Recall that for points $x,y\in\bord X$, the notation $\geo xy$ denotes the geodesic segment, line, or ray joining $x$ and $y$.

\begin{definition}
\label{definitionconvexhull}
Given $S\subset\bord X$, let
\begin{align*}
\Hull_1(S) &:= \bigcup_{x,y\in S} \geo xy,\\
\Hull_n(S) &:= \underbrace{\Hull_1\cdots \Hull_1}_{\text{$n$ times}}(S)\\
\Hull_\infty(S) &:= \cl{\bigcup_{n\in\N}\Hull_n(S)}.
\end{align*}
The set $\Hull_n(S)$ will be called the \emph{$n$th convex hull} of $S$. Moreover, $\Hull_\infty(S)$ will be called the \emph{convex hull} of $S$, and $\Hull_1(S)$ will be called the \emph{quasiconvex hull} of $S$.
\end{definition}

The terminology ``convex hull'' comes from the following fact:

\begin{proposition}
\label{propositionconvexhull}
$\Hull_\infty(S)$ is the smallest closed set $F\subset \bord X$ such that $S\subset F$ and
\begin{equation}
\label{convex}
\geo xy\subset F \all x,y\in F.
\end{equation}
\end{proposition}
A set $F$ satisfying \eqref{convex} will be called \emph{convex}.
\begin{proof}
It is clear that $S\subset\Hull_\infty(S)\subset\bord X$. To show that $\Hull_\infty(S)$ is convex, fix $x,y\in\Hull_\infty(S)$. Then there exist sequences $A\ni x_n\to x$ and $A\ni y_n\to y$, where $A = \bigcup_{n\in\N}\Hull_n(S)$. For each $n$, $\geo{x_n}{y_n}\subset A\subset \Hull_\infty(S)$. But since $X$ is regularly geodesic, $\geo{x_n}{y_n}\to \geo xy$ in the Hausdorff metric on $\bord X$. Since $\Hull_\infty(S)$ is closed, it follows that $\geo xy\subset\Hull_\infty(S)$.

Conversely, if $S\subset F\subset \bord X$ is a closed convex set, then an induction argument shows that $F\supset \Hull_n(S)$ for all $n$. Since $F$ is closed, we have $F\supset \Hull_\infty(S)$.
\end{proof}

Another connection between the operations $\Hull_1$ and $\Hull_\infty$ is given by the following proposition:

\begin{proposition}
\label{propositionhull}
Suppose that $X$ is a algebraic hyperbolic space. Then there exists $\tau > 0$ such that for every set $S\subset\bord X$ we have
\[
X\cap \Hull_1(S) \subset X\cap \Hull_\infty(S) \subset \thicken_\tau(X\cap\Hull_1(S)).
\]
(Recall that $\thicken_\tau(S)$ denotes the $\tau$-thickening of a set with respect to the hyperbolic metric $\dist$.)
\end{proposition}
\begin{proof}
The proof will proceed using the ball model $X = \B = \B_\F^\alpha$. We will need the following lemma:
\begin{lemma}
\label{lemmaF}
There exists a closed convex set $F\propersubset\bord\B$ whose interior intersects $\del\B$.
\end{lemma}
\begin{subproof}
If $\alpha < \infty$, this is a consequence of \cite[Theorem 3.3]{Anderson}.

We will use the finite-dimensional case to prove the infinite-dimensional case. Suppose that $\alpha$ is infinite. Let $Y = \B_\F^3 \subset X$. Then by the $\alpha < \infty$ case of Lemma \ref{lemmaF}, there exists a closed convex set $F_2 \propersubset\bord Y$ whose interior intersects $\del Y$, say $\xi\in \interior(F_2)\cap\del Y$. Choose $\epsilon > 0$ such that $B_Y(\xi,\epsilon) \subset F_2$. Then
\[
F_1 := \Hull_\infty(B_Y(\xi,\epsilon)) \subset F_2 \propersubset \bord Y
\]
by Proposition \ref{propositionconvexhull}. On the other hand, $F_1$ is invariant under the action of the group
\[
G_1 := \{g\in\Isom(Y) : g(\0) = \0, g(\xi) = \xi\}.
\]
Let
\[
G = \{g\in\Isom(X) : g(\0) = \0, g(\xi) = \xi\},
\]
and note that $G(\bord Y) = \bord X$. Let $F = G(F_1)$, and note that $F\cap \bord Y = F_1$. We claim that $F$ is convex. Indeed, suppose that $x,y\in F$; then there exists $g\in G$ such that $g(x),g(y)\in\bord Y$. (Note that in this step, we need all three dimensions of $Y$.) Then $g(x),g(y)\in F\cap \bord Y = F_1$, so by the convexity of $F_1$ we have $g(\geo xy) = \geo{g(x)}{g(y)} \subset F_1 \subset F$. Since $F$ is $G$-invariant, we have $\geo xy\subset F$.

In addition to being convex, $F$ is also closed and contains the set $G(B_Y(\xi,\epsilon)) = B_X(\xi,\epsilon)$. Thus, $\xi\in\interior(F)$. Finally, since $F\cap \bord Y = F_1 \propersubset \bord Y$, it follows that $F\propersubset\bord X$.
\end{subproof}
Let $F$ be as in Lemma \ref{lemmaF}. Since $F\propersubset\bord\B$ is a closed set, it follows that $\B\butnot F\neq\emptyset$. By the transitivity of $\Isom(\B)$ (Observation \ref{observationtransitivity}), we may without loss of generality assume that $\0\in\B\butnot F$. By the transitivity of $\Stab(\Isom(\B);\0)$ on $\del\B$, we may without loss of generality assume that $\ee_1\in\interior(F)$. Fix $\epsilon > 0$ such that $B(\ee_1,\epsilon) \subset F$.

We now proceed with the proof of Proposition \ref{propositionhull}. It is clear from the definitions that $\B\cap\Hull_1(S) \subset \B\cap\Hull_\infty(S)$. To prove the second inclusion, fix $\zz\in\B \butnot \thicken_\tau(\Hull_1(S))$ and we will show that $\zz\notin\Hull_\infty(S)$. By the transitivity of $\Isom(\B)$, we may without loss of generality assume that $\zz = \0$. Now for every $\xx, \yy\in S$, we have $\zz = \0 \notin \thicken_\tau(\geo\xx\yy)$. By (i) of Proposition \ref{propositionrips}, we have
\[
\lb \xx | \yy \rb_\0 \gtrsim_\plus \tau
\]
and thus by \eqref{euclideanball1},
\[
\|\yy - \xx\| \lesssim_\times e^{-\tau}.
\]
By choosing $\tau$ sufficiently large, this implies that
\[
\|\yy - \xx\| \leq \epsilon/2 \all \xx,\yy\in S.
\]
Moreover, since $\dist(\0,\xx) = 2\lb \xx | \xx \rb_\0 \gtrsim_\plus \tau$, by choosing $\tau$ sufficiently large we may guarantee that
\[
\|\xx\| \geq 1 - \epsilon/2 \all \xx\in S.
\]
Since the claim is trivial if $S = \emptyset$, assume that $S\neq \emptyset$ and choose $\xx\in S$. Without loss of generality, assume that $\xx = \lambda \ee_1$ for some $\lambda \geq 2/3$. Then $S \subset B_\EE(\xx,\epsilon/2) \subset B(\ee_1,\epsilon) \subset F$. But then $F$ is a closed convex set containing $S$, so by Proposition \ref{propositionconvexhull} $\Hull_\infty(S) \subset F$. Since $\zz = \0 \notin F$, it follows that $\zz\notin\Hull_\infty(S)$.
\end{proof}
\begin{corollary}
\label{corollaryanderson}
Suppose that $X$ is an algebraic hyperbolic space. Then for every closed set $S\subset\bord X$, we have
\[
\Hull_\infty(S) \cap \del X = S\cap \del X.
\]
\end{corollary}
\begin{proof}
The inclusion $\supset$ is immediate. Suppose that $\xi\in\Hull_\infty(S)\cap\del X$, and find a sequence $X\cap\Hull_\infty(S)\ni x_n\to \xi$. By Proposition \ref{propositionhull}, for each $n$ there exist $y_n^{(1)},y_n^{(2)}\in S$ such that $x_n\in\thicken_\tau(\geo{y_n^{(1)}}{y_n^{(2)}})$; by Proposition \ref{propositionrips} we have $\lb y_n^{(1)} | y_n^{(2)}\rb_{x_n} \asymp_\plus 0$. Applying Lemma \ref{lemmaminimal} gives $\xi\in S$.
\end{proof}

\begin{remark}
Corollary \ref{corollaryanderson} was proven for the case where $X$ is a pinched (finite-dimensional) Hadamard manifold and $S\subset\del X$ by M. T. Anderson \cite[Theorem 3.3]{Anderson}. It was conjectured to hold whenever $X$ is ``strictly convex'' by Gromov \cite[p.11]{Gromov1}, who observed that it holds in the Standard Case. However, this conjecture was proven to be false independently by A. Ancona \cite[Corollary C]{Ancona} and A. Borb\'ely \cite[Theorem 1]{Borbely}, who each constructed a three-dimensional CAT(-1) manifold $X$ and a point $\xi\in\del X$ such that for every neighborhood $U$ of $\xi$, $\Hull_\infty(U) = \bord X$.


Thus, although the $\infty$-convex hull has more geometric and intuitive appeal based on Proposition \ref{propositionconvexhull}, without more hypotheses there is no way to restrain its geometry. The $1$-convex hull is thus more useful for our applications. Proposition \ref{propositionhull} indicates that in the case of an algebraic hyperbolic space, we are not losing too much by the change.
\end{remark}

\begin{definition}
\label{definitionconvexcore}
The \emph{convex core} of a semigroup $G\prec\Isom(X)$ is the set
\[
\CC_\Lambda := X\cap \Hull_\infty(\LambdaG),
\]
and the \emph{quasiconvex core} is the set
\[
\CC_\zero := X\cap \cl{\Hull_1(G(\zero))}.
\]
\end{definition}

\begin{observation}
The convex core and quasiconvex core are both closed $G$-invariant sets. The quasiconvex core depends on the distinguished point $\zero$. However:
\end{observation}

\begin{proposition}
\label{propositionC0comparison}
Fix $x,y\in X$. Then
\[
\CC_x \subset \thicken_R(\CC_y)
\]
for some $R > 0$.
\end{proposition}
\begin{proof}
Fix $z\in\CC_x$. Then $z\in \geo{g(x)}{h(x)}$ for some $g,h\in G$. It follows that
\[
\lb g(y)|h(y)\rb_z \asymp_\plus \lb g(x)|h(x)\rb_z = 0.
\]
So by Proposition \ref{propositionrips}, $\dist(z,\geo{g(y)}{h(y)}) \asymp_\plus 0$. But $\geo{g(y)}{h(y)} \subset \CC_y$, so $\dist(z,\CC_y) \asymp_\plus 0$. Letting $R$ be the implied constant completes the proof.
\end{proof}

\begin{remark}
In many cases, we can get information about the action of $G$ on $X$ by looking just at its restriction to $\CC_\Lambda$ or to $\CC_\zero$. We therefore also remark that if $X$ is a CAT(-1) space, then $\CC_\Lambda$ is also a CAT(-1) space.
\end{remark}

In the sequel the following notation will be useful:
\begin{notation}
For a set $S\subset\bord X$ let
\begin{equation}
\label{primenotation}
S' = \cl S \cap \del X.
\end{equation}
\end{notation}

\begin{observation}
\label{observationboundaryofconvexcore}
$(\CC_\zero)' = \LambdaG$.
\end{observation}
\begin{proof}
Since $\LambdaG = (G(\zero))'$ and $G(\zero) \subset\CC_\zero$, we have $(\CC_\zero)' \supset \LambdaG$. Suppose that $\xi\in(\CC_\zero)'$, and let $\CC_\zero\ni x_n\to \xi$. By definition, for each $n$ there exist $y_n^{(1)},y_n^{(2)}\in G(\zero)$ such that $x_n\in \geo{y_n^{(1)}}{y_n^{(2)}}$. Lemma \ref{lemmaminimal} completes the proof.
\end{proof}

\bigskip
\section{Semigroups which act irreducibly on algebraic hyperbolic spaces}
\label{subsectionirreducibly}

\begin{definition}
\label{definitionreducibly}
Suppose that $X$ is an algebraic hyperbolic space, and fix $G\prec\Isom(X)$. We shall say that $G$ \emph{acts reducibly} on $X$ if there exists a nontrivial totally geodesic $G$-invariant subset $S\propersubset \bord X$. Otherwise, we shall say that $G$ \emph{acts irreducibly} on $X$.
\end{definition}

\begin{remark}
A parabolic or focal subsemigroup of $\Isom(X)$ may act either reducibly or irreducibly on $X$.
\end{remark}

%
%

\begin{proposition}
\label{propositionactsreducibly}
Let $G\prec\Isom(X)$ be nonelementary. Then the following are equivalent:
\begin{itemize}
\item[(A)] $G$ acts reducibly on $X$.
\item[(B)] There exists a nontrivial totally geodesic subset $S\propersubset \bord X$ such that $\LambdaG\subset S$.
\item[(C)] There exists a nontrivial totally geodesic subset $S\propersubset \bord X$ such that $\CC_\Lambda\subset S$.
\item[(D)] There exists a nontrivial totally geodesic subset $S\propersubset \bord X$ such that $\CC_\zero\subset S$ for some $\zero\in X$.
\end{itemize}
\end{proposition}
\begin{proof}[Proof of \text{(A) \implies (B)}]
Let $S\propersubset \bord X$ be a nontrivial totally geodesic $G$-invariant subset. Fix $\zero\in S\cap X$. Then $\LambdaG\subset\cl{G(\zero)}\subset S$.
\end{proof}
\begin{proof}[Proof of \text{(B) \implies (C)}]
If $S$ is any totally geodesic set which contains $\LambdaG$, then $S$ is a closed convex set containing $\LambdaG$, so by Proposition \ref{propositionconvexhull}, $\CC_\Lambda\subset S$.
\end{proof}
\begin{proof}[Proof of \text{(C) \implies (D)}]
Since $G$ is nonelementary, $\CC_\Lambda \neq \emptyset$. Fix $\zero\in\CC_\Lambda$; then $\CC_\zero\subset \CC_\Lambda$.
\end{proof}
\begin{proof}[Proof of \text{(D) \implies (A)}]
Let $S$ be the smallest totally geodesic subset of $X$ which contains $\CC_\zero$, i.e.
\[
S := \bigcap\{W:W\supseteq \CC_\zero\text{ totally geodesic}\}.
\]
Then our hypothesis implies that $S\propersubset \bord X$. Since $\zero\in S$, $S$ is nontrivial. It is obvious from the definition that $S$ is $G$-invariant. This completes the proof.
\end{proof}

\begin{remark}
If $G\prec\Isom(X)$ is nonelementary, then Proposition \ref{propositionactsreducibly} gives us a way to find a nontrivial totally geodesic set on which $G$ acts reducibly; namely, the smallest totally geodesic set containing $\LambdaG$, or equivalently $\CC_G$, will have this property (cf. Lemma \ref{lemmatotallygeodesic}). On the other hand, there exists a parabolic group $G\leq\Isom(\H^\infty)$ such that $G$ does not act irreducibly on any nontrivial totally geodesic subset $S\subset\bord \H^\infty$ (Remark \ref{remarkparabolictorsion}).
\end{remark}

\bigskip
\section{Semigroups of compact type}

\begin{definition}
\label{definitioncompacttype}
We say that a semigroup $G\prec\Isom(X)$ is of \emph{compact type} if its limit set $\Lambda$ is compact.
\end{definition}

\begin{proposition}
\label{propositioncompacttype}
For $G\prec\Isom(X)$, the following are equivalent:
\begin{itemize}
\item[(A)] $G$ is of compact type.
\item[(B)] Every sequence $(x_n)_1^\infty$ in $G(\zero)$ with $\dox{x_n}\to \infty$ has a convergent subsequence.
\end{itemize}
Furthermore, if $X$ is regularly geodesic, then \text{(A)-(B)} are equivalent to:
\begin{itemize}
\item[(C)] The set $\CC_\zero$ is a proper metric space.
\end{itemize}
and if $X$ is an algebraic hyperbolic space, then they are equivalent to:
\begin{itemize}
\item[(D)] The set $\CC_\Lambda$ is a proper metric space.
\end{itemize}
\end{proposition}
\begin{proof}[Proof of \text{(A)} \implies \text{(B)}]
Fix a sequence $(g_n)_1^\infty$ in $G$ with $\dogo{g_n}\to \infty$. The existence of such a sequence implies that $G$ is not elliptic. If $G$ is parabolic or inward focal, then the proof of Proposition \ref{propositioncardinalitylimitset}(ii) shows that $g_n(\zero)\to \xi$, where $\Lambda = \{\xi\}$. So we may assume that $G$ is lineal, outward focal, or of general type, in which case Proposition \ref{propositioncardinalitylimitset} gives $\#(\Lambda)\geq 2$.

Fix distinct $\xi_1,\xi_2\in\LambdaG$, and let $(n_k)_1^\infty$ be a sequence such that $(g_{n_k}(\xi_i))_1^\infty$ converges for $i = 1,2$, and such that
\[
\lb g_{n_k}^{-1}(\zero)|\xi_1\rb_\zero \leq \lb g_{n_k}^{-1}(\zero)|\xi_2\rb_\zero
\]
for all $k$. (If this is not possible, switch $\xi_1$ and $\xi_2$.) We have
\[
0 \asymp_{\plus,\xi_1,\xi_2} \lb \xi_1|\xi_2\rb_\zero
\gtrsim_\plus \min\left(\lb g_{n_k}^{-1}(\zero)|\xi_1\rb_\zero,\lb g_{n_k}^{-1}(\zero)|\xi_2\rb_\zero\right)
= \lb g_{n_k}^{-1}(\zero)|\xi_1\rb_\zero
\]
and thus
\[
\lb g_{n_k}(\zero)|g_{n_k}(\xi_1)\rb_\zero \asymp_{\plus,\xi_1,\xi_2} \dogo{g_{n_k}} \tendsto n \infty.
\]
On the other hand, there exists $\eta\in\LambdaG$ such that $g_{n_k}(\xi_1)\tendsto k\eta$, and thus
\[
\lb g_{n_k}(\xi_1)|\eta\rb_\zero\tendsto n \infty.
\]
Applying Gromov's inequality yields
\[
\lb g_{n_k}(\zero)|\eta\rb_\zero \tendsto n \infty
\]
and thus $g_{n_k}(\zero)\tendsto k\eta$. This completes the proof.
\end{proof}
\begin{proof}[Proof of \text{(B)} \implies \text{(A)}]
Fix a sequence $(\xi_n)_1^\infty$ in $\LambdaG$. For each $n\in\N$, choose $g_n\in G$ with
\[
\lb g_n(\zero)|\xi_n\rb_\zero \geq n.
\]
In particular $\dogo{g_n}\geq n\tendsto n \infty$. Thus by our hypothesis, there exists a convergent subsequence $g_{n_k}(\zero)\tendsto k\eta\in\LambdaG$. Now
\[
\Dist(\xi_{n_k},\eta) \leq \Dist(g_{n_k}(\zero),\xi_{n_k}) + \Dist(g_{n_k}(\zero),\eta)
\lesssim_\times b^{-n_k} + \Dist(g_{n_k}(\zero),\eta) \tendsto k 0,
\] 
i.e. $\xi_{n_k}\tendsto k\eta$.
\end{proof}
\begin{proof}[Proof of \text{(A)} \implies \text{(C)}]
Let $(x_n)_1^\infty$ be a bounded sequence in $\CC_\zero$. For each $n\in\N$, there exist $y_n^{(1)},y_n^{(2)}\in G(\zero)$ such that $x_n\in\thicken_{1/n}(\geo{y_n^{(1)}}{y_n^{(2)}})$. Choose a sequence $(n_k)_1^\infty$ on which $y_{n_k}^{(1)}\tendsto k\alpha$ and $y_{n_k}^{(2)}\tendsto k\beta$. Since $X$ is regularly geodesic we have
\[
\geo{y_{n_k}^{(1)}}{y_{n_k}^{(2)}}\tendsto k\geo{y^{(1)}}{y^{(2)}}.
\]
For each $k$, choose $z_k\in \geo{y_{n_k}^{(1)}}{y_{n_k}^{(2)}}$ with $\dist(x_{n_k},z_k)\leq 1/n_k$. Since the sequence $(z_k)_1^\infty$ is bounded, it must have a subsequence which converges to a point in $\geo{y^{(1)}}{y^{(2)}}$; it follows that the corresponding subsequence of $(x_{n_k})_1^\infty$ is also convergent. Thus every bounded sequence in $\CC_\zero$ has a convergent subsequence, so $\CC_\zero$ is proper.
\end{proof}
\begin{proof}[Proof of \text{(C)} \implies \text{(B)}]
Obvious since $G(\zero)\subset\CC_\zero$.
\end{proof}
\begin{proof}[Proof of \text{(A)} \implies \text{(D)}]
Note first of all that we cannot get (A) \implies (D) immediately from Proposition \ref{propositionhull}, since the $\tau$-thickening of a compact set is no longer compact.

By \cite[Proposition 1.5]{Bowditch_convexity}, there exists a metric $\rho$ on $\bord X$ compatible with the topology such that the map $F\mapsto \Hull_1(F)$ is a semicontraction with respect to the Hausdorff metric of $(\bord X,\rho)$. (Finite-dimensionality is not used in any crucial way in the proof of \cite[Proposition 1.5]{Bowditch_convexity},\Footnote{One should keep in mind that the Cartan--Hadamard theorem \cite[IX, Theorem 3.8]{Lang_differential_geometry} can be used as a substitute for the Hopf--Rinow theorem in most circumstances.} and in any case for algebraic hyperbolic spaces it can be proven by looking at finite-dimensional subsets, as we did in the proof of Proposition \ref{propositionhull}.) We remark that if $\F = \R$, then such a metric $\rho$ can be prescribed explicitly: if $X = \B$ is the ball model, then the Euclidean metric on $\bord\B\subset\HH$ has this property, due to the fact that geodesics in the ball model are line segments in $\HH$ (cf. \eqref{geodesic}).\Footnote{Recall that our ``ball model'' $\amsbb B$ is the Klein model rather than the Poincar\'e model.}

Now let us demonstrate (D). It suffices to show that $\bord\CC_\Lambda = \Hull_\infty(\Lambda)$ is compact. Since $\Hull_\infty(\Lambda)$ is by definition closed, it suffices to show that $\Hull_\infty(\Lambda)$ is totally bounded with respect to the $\rho$ metric. Indeed, fix $\epsilon > 0$. Since $\Lambda$ is compact, there is a finite set $F_\epsilon\subset\Lambda$ such that
\[
\Lambda\subset \thicken_{\epsilon/2}(F_\epsilon).
\]
(In this proof, all neighborhoods are taken with respect to the $\rho$ metric.) Let $X_\epsilon \subset X$ be a finite-dimensional totally geodesic set containing $F_\epsilon$. Then $\Lambda\subset \thicken_{\epsilon/2}(X_\epsilon)$. On the other hand, since $X_\epsilon$ is compact, there exists a finite set $F_\epsilon'\subset X_\epsilon$ such that $X_\epsilon\subset \thicken_{\epsilon/2}(F_\epsilon')$.

Now, our hypothesis on $\rho$ implies that
\[
\Hull_1(\thicken_{\epsilon/2}(X_\epsilon)) \subset \thicken_{\epsilon/2}(\Hull_1(X_\epsilon)) = \thicken_{\epsilon/2}(X_\epsilon),
\]
and thus that $\thicken_{\epsilon/2}(X_\epsilon)$ is convex. But $\Lambda\subset \thicken_{\epsilon/2}(X_\epsilon)$, so $\Hull_\infty(\Lambda) \subset \thicken_{\epsilon/2}(X_\epsilon)$. Thus
\[
\Hull_\infty(\Lambda)\subset \thicken_{\epsilon/2}(X_\epsilon) \subset \thicken_\epsilon(F_\epsilon').
\]
Since $\epsilon$ was arbitrary, this shows that $\Hull_\infty(\Lambda)$ is totally bounded, completing the proof.
\end{proof}
\begin{proof}[Proof of \text{(D)} \implies \text{(B)}]
Since property (B) is clearly basepoint-independent, we may without loss of generality suppose $\zero\in\CC_\Lambda$. Then (D) \implies (C) \implies (B).
\end{proof}

As an example of an application we prove the following corollary.
\begin{corollary}
Suppose that $X$ is regularly geodesic. Then any moderately discrete subgroup of $\Isom(X)$ of compact type is strongly discrete.
\end{corollary}
\begin{proof}
If $G$ is a moderately discrete group, then $G\given\CC_\zero$ is moderately discrete by Observation \ref{observationrestrictions}, and therefore strongly discrete by Propositions \ref{propositionproperMDSD} and \ref{propositioncompacttype}. Thus by Observation \ref{observationrestrictions}, $G$ is strongly discrete.
\end{proof}

A well-known characterization of the complement of the limit set in the Standard Case is that it is the set of points where the action of $G$ is discrete. We extend this characterization to hyperbolic metric spaces for groups of compact type:

\begin{proposition}
\label{propositionSDonclX}
Let $G\leq\Isom(X)$ be a strongly discrete group of compact type. Then the action of $G$ on $\bord X\butnot\LambdaG$ is strongly discrete in the following sense: For any set $S\subset\bord X\butnot\LambdaG$ satisfying
\begin{equation}
\label{Lambdabounded}
\Dist(S,\LambdaG) > 0,
\end{equation}
we have
\[
\#\{g\in G : g(S)\cap S \neq \emptyset\} < \infty.
\]
\end{proposition}
\begin{proof}
By contradiction, suppose that there exists a sequence of distinct $(g_n)_1^\infty$ such that $g_n(S)\cap S\neq \emptyset$ for all $n\in\N$. Since $G$ is strongly discrete, we have $\dogo{g_n}\to\infty$, and since $G$ is of compact type there exist an increasing sequence $(n_k)_1^\infty$ and $\xi_+,\xi_-\in\LambdaG$ such that $g_{n_k}(\zero)\to \xi_+$ and $g_{n_k}^{-1}(\zero)\to\xi_-$. In the remainder of the proof we restrict to this subsequence, so that $g_n(\zero)\to \xi_+$ and $g_n^{-1}(\zero)\to\xi_-$.

For each $n$, fix $x_n\in g_n^{-1}(g_n(S)\cap S)$, so that $x_n,g_n(x_n)\in S$. Then
\[
\Dist(x_n,\xi_-), \Dist(g_n(x_n),\xi_+) \geq \Dist(S,\LambdaG) \asymp_\times 1,
\]
and so
\[
\lb x_n|\xi_-\rb_\zero , \lb g_n(x_n)|\xi_+ \rb_\zero \asymp_\plus 0.
\]
On the other hand, $\lb g_n^{-1}(\zero)|\xi_-\rb_\zero, \lb g_n(\zero)|\xi_+\rb_\zero \to \infty$. Applying Gromov's inequality gives
\[
\lb x_n|g_n^{-1}(\zero)\rb_\zero , \lb g_n(x_n)|g_n(\zero) \rb_\zero \asymp_\plus 0
\]
for all $n$ sufficiently large. But then
\[
\dogo{g_n} = \lb g_n(x_n)|\zero\rb_{g_n(\zero)} + \lb g_n(x_n)|g_n(\zero) \rb_\zero \asymp_\plus 0,
\]
a contradiction.
\end{proof}

\part{The Bishop--Jones theorem}
\label{partbishopjones}
This part will be divided as follows: In Chapter \ref{sectionmodified}, we motivate and define the \emph{modified} Poincar\'e exponent of a semigroup, which is used in the statement of Theorem \ref{theorembishopjonesmodified}. In Chapter \ref{sectionbishopjones} we prove Theorem \ref{theorembishopjonesmodified} and deduce Theorem \ref{theorembishopjonesregular} from Theorem \ref{theorembishopjonesmodified}.

\chapter{The modified Poincar\'e exponent} \label{sectionmodified}
In this chapter we define the modified Poincar\'e exponent of a semigroup. We first recall the classical notion of the Poincar\'e exponent, introduced in the Standard Case by A. F. Beardon in \cite{Beardon1}. Although it is usually defined only for groups, the generalization to semigroups is trivial.

\section{The Poincar\'e exponent of a semigroup}
\label{subsectionpoincare}

\begin{definition}
\label{definitionpoincareexponent}
Fix $G\prec\Isom(X)$. For each $s\geq 0$, the series
\[
\Sigma_s(G) := \sum_{g\in G}b^{-s\dogo g}
\]
is called the \emph{Poincar\'e series} of the semigroup $G$ in dimension $s$ (or ``evaluated at $s$'') relative to $b$. The number
\[
\delta_G = \delta(G) := \inf\{s\geq 0:\Sigma_s(G) < \infty\}
\]
is called the \emph{Poincar\'e exponent} of the semigroup $G$ relative to $b$. Here, we let $\inf\emptyset = \infty$.
\end{definition}

\begin{remark}
The Poincar\'e series is usually defined with a summand of $e^{-s\dogo g}$ rather than $b^{-s\dogo g}$. The change of exponents here is important because it relates the Poincar\'e exponent to the metric $\Dist = \Dist_{b,\zero}$ defined in Proposition \ref{propositionDist}. In the Standard Case, and more generally for CAT(-1) spaces, we have made the convention that $b = e$ (see \6\ref{standingassumptions2}), so in this case our series reduces to the classical one.
\end{remark}

\begin{remark}
\label{remarkorbitalcounting}
Given $G\prec\Isom(X)$, we may define the \emph{orbital counting function} of $G$ to be the function
\[
\NN_{X,G}(\rho) = \#\{g\in G : \dogo g \leq \rho\}.
\]
The Poincar\'e series may be written as an integral over the orbital counting function as follows:
\begin{equation}
\label{poincareseriesalternate}
\begin{split}
\Sigma_s(G) &= \log(b^s) \sum_{g\in G} \int_{\dogo g}^\infty b^{-s\rho} \;\dee \rho\\
&= \log(b^s) \int_0^\infty b^{-s\rho} \sum_{g\in G}[\dogo g\leq \rho] \;\dee \rho\\
&= \log(b^s) \int_0^\infty b^{-s\rho} \NN_{X,G}(\rho) \;\dee \rho.
\end{split}
\end{equation}
The Poincar\'e exponent is written in terms of the orbital counting function as
\begin{equation}
\label{poincarealternate}
\delta_G = \limsup_{\rho\to \infty} \frac1\rho\log_b \NN_{X,G}(\rho)
\end{equation}
\end{remark}

\begin{definition}
\label{definitiondivergencetype}
A semigroup $G\prec\Isom(X)$ with $\delta_G < \infty$ is said to be of \emph{convergence type} if $\Sigma_{\delta_G}(G) < \infty$. Otherwise, it is said to be of \emph{divergence type}. In the case where $\delta_G = \infty$, we say that the semigroup is neither of convergence type nor of divergence type.
\end{definition}

The most basic question about the Poincar\'e exponent is whether it is finite. For groups, the finiteness of the Poincar\'e exponent is related to strong discreteness:

\begin{observation}\label{observationnotstronglydiscrete}
Fix $G\leq\Isom(X)$. If $G$ is not strongly discrete, then $\delta_G = \infty$.
\end{observation}
\begin{proof}
Fix $\rho > 0$ such that $\#\{g\in G : \dogo g \leq \rho\} = \infty$. Then for all $s\geq 0$ we have
\[
\Sigma_s(G) \geq \sum_{\substack{g\in G\\ \dogo g\leq \rho}} b^{-s \dogo g} \geq \sum_{\substack{g\in G\\ \dogo g\leq \rho}} b^{-s\rho} = \infty .
\]
Since $s$ was arbitrary, we have $\delta_G = \infty$.
\end{proof}

\begin{remark}
Although the converse to Observation \ref{observationnotstronglydiscrete} holds in the Standard Case, it fails for infinite-dimensional algebraic hyperbolic spaces; see Example \ref{exampleinfinitepoincareBIM}.
\end{remark}

\begin{notation}
\label{notationpoincareset}
The Poincar\'e exponent and type can be conveniently combined into a single mathematical object, the \emph{Poincar\'e set}
\[
\Delta_G := \{s\geq 0 : \Sigma_s(G) = \infty\} = \begin{cases}
[0,\delta_G] & \text{$G$ is of divergence type}\\
\CO 0{\delta_G} & \text{$G$ is of convergence type}\\
\Rplus & \delta_G = \infty
\end{cases}.
\]
\end{notation}

\section{The modified Poincar\'e exponent of a semigroup}
\label{subsectionmodified}

From a certain perspective, Observation \ref{observationnotstronglydiscrete} indicates a flaw in the Poincar\'e exponent: If $G\leq\Isom(X)$ is not strongly discrete, then the Poincar\'e exponent is always infinity even though there may be more geometric information to capture. In this section we introduce a modification of the Poincar\'e exponent which agrees with the Poincar\'e exponent in the case where $G$ is strongly discrete, but can be finite even if $G$ is not strongly discrete.

We begin by defining the modified Poincar\'e exponent of a locally compact group $G\leq\Isom(X)$. Let $\mu$ be a Haar measure on $G$, and for each $s$ consider the \emph{Poincar\'e integral}
\begin{equation}
\label{IsG}
I_s(G) := \int b^{-s\dogo g}\;\dee\mu(g).
\end{equation}
\begin{definition}
\label{definitionmodified1}
The \emph{modified Poincar\'e exponent} of a locally compact group $G\leq\Isom(X)$ is the number
\[
\w\delta_G = \w\delta(G) := \inf\{s\geq 0: I_s(G) < \infty\},
\]
where $I_s(G)$ is defined by \eqref{IsG}.
\end{definition}
\begin{example}
\label{exampleLie}
Let $X = \H^d$ for some $2\leq d < \infty$, and let $G\leq\Isom(X)$ be a positive-dimensional Lie subgroup. Then $G$ is locally compact, but not strongly discrete. Although the Poincar\'e series diverges for every $s$, the exponent of convergence of the Poincar\'e integral (or ``modified Poincar\'e exponent'') is equal to the Hausdorff dimension of the limit set of $G$ (Theorem \ref{theorembishopjonesmodified} below), and so in particular the Poincar\'e integral converges whenever $s > d - 1$.
\end{example}

We now proceed to generalize Definition \ref{definitionmodified1} to the case where $G\leq\Isom(X)$ is not necessarily locally compact. Fix $\rho > 0$, and consider a maximal $\rho$-separated\Footnote{Here, as usual, a \emph{$\rho$-separated} subset of a metric space $X$ is a set $S\subset X$ such that $\dist(x,y)\geq \rho$ for any distinct $x,y\in S$. The existence of a maximal $\rho$-separated subset of any metric space is guaranteed by Zorn's lemma.\label{footnoterhoseparated}} subset $S_\rho\subset G(\zero)$. Then we have
\[
\bigcup_{x\in S_\rho} B(x,\rho/2) \subset G(\zero) \subset \bigcup_{x\in S_\rho} B(x,\rho),
\]
and the former union is disjoint. Now suppose that $G$ is in fact locally compact, and let $\nu$ denote the image of Haar measure on $G$ under the map $g\mapsto g(\zero)$. Then if $f$ is a positive function on $X$ whose logarithm is uniformly continuous, we have
\[
\sum_{x\in S_\rho} f(x)
\asymp_{\times,\rho,f} \sum_{x\in S_\rho} \int_{B(x,\rho/2)} f\;\dee\nu
\leq \int f\;\dee\nu
\leq \sum_{x\in S_\rho} \int_{B(x,\rho)} f\;\dee\nu
\asymp_{\times,\rho,f} \sum_{x\in S_\rho} f(x).
\]
Thus in some sense, the counting measure on $S_\rho$ is a good approximation to the measure $\nu$. In particular, taking $f(x) = b^{-\dox x}$ gives
\[
I_s(G) \asymp_{\times,\rho} \sum_{x\in S_\rho} b^{-\dox x}.
\]
Thus the integral $I_s(G)$ converges if and only if the series $\sum_{x\in S_\rho} b^{-\dox x}$ converges. But the latter series is well-defined even if $G$ is not locally compact. This discussion shows that the definition of the ``modified Poincar\'e exponent'' given in Definition \ref{definitionmodified1} agrees with the following definition:

\begin{definition}
\label{definitionmodifiedexponent}
Fix $G\prec\Isom(X)$.
\begin{itemize}
\item For each set $S\subset X$ and $s\geq 0$, let
\begin{align*}
\Sigma_s(S) &= \sum_{x\in S} b^{-s\dox x}\\
\Delta(S) &= \{s\geq 0: \Sigma_s(S) = \infty\}\\
\delta(S) &= \sup\Delta(S).
\end{align*}
\item Let
\begin{equation}
\label{modifiedpoincaredef}
\w\Delta_G = \bigcap_{\rho > 0} \bigcap_{S_\rho} \Delta(S_\rho),
\end{equation}
where the second intersection is taken over all maximal $\rho$-separated sets $S_\rho$.
\item The number $\w\delta_G = \sup\w\Delta_G$ is called the \emph{modified Poincar\'e exponent} of $G$. If $\w\delta_G\in\w\Delta_G$, we say that $G$ is of \emph{generalized divergence type},\Footnote{We use the adjective ``generalized'' rather than ``modified'' because all groups of convergence/divergence type are also of generalized convergence/divergence type; see Corollary \ref{corollarygeneralizedtypes} below.} while if $\w\delta_G\in\Rplus\butnot\Delta_G$, we say that $G$ is of \emph{generalized convergence type}. Note that if $\w\delta_G = \infty$, then $G$ is neither of generalized convergence type nor of generalized divergence type.
\end{itemize}
\end{definition}

The basic properties of the modified Poincar\'e exponent are summarized as follows:

\begin{proposition}
\label{propositionbasicmodified}
Fix $G\prec\Isom(X)$.
\begin{itemize}
\item[(i)] $\w\Delta_G \subset \Delta_G$. (In particular $\w\delta_G\leq\delta_G$.)
\item[(ii)] If $G$ satisfies
\begin{equation}
\label{SD3}
\sup_{x\in X} \#\{g\in G:\dist(g(\zero),x)\leq \rho\} < \infty \all \rho > 0,
\end{equation}
then $\w\Delta_G = \Delta_G$. (In particular $\w\delta_G = \delta_G$.)
\item[(iii)] If $\w\delta_G < \infty$, then there exist $\rho > 0$ and a maximal $\rho$-separated set $S_\rho\subset G(\zero)$ such that $\#(S_\rho\cap B) < \infty$ for every bounded set $B$.
\item[(iv)] For all $\rho > 0$ sufficiently large and for every maximal $\rho$-separated set $S_\rho\subset G(\zero)$, we have $\Delta(S_\rho) = \w\Delta_G$. (In particular $\delta(S_\rho) = \w\delta(G)$.)
\end{itemize}
\end{proposition}
\begin{remark}
\label{remarkstronglydiscrete}
If $G$ is a group, then it is clear that \eqref{SD3} is equivalent to the assertion that $G$ is strongly discrete. If $G$ is not a group, then by analogy we will say that $G$ is \emph{strongly discrete} if \eqref{SD3} holds. (Recall that in Chapter \ref{sectiondiscreteness}, the various notions of discreteness are defined only for groups.)
\end{remark}

\begin{proof}[Proof of Proposition \ref{propositionbasicmodified}]
~
\begin{itemize}
\item[(i)] Indeed, for every $s\geq 0$, $\rho > 0$, and maximal $\rho$-separated set $S_\rho$ we have
\[
\Sigma_s(S_\rho) \leq \Sigma_s(G) \text{ and thus }\w\Delta(G) \subset \Delta(S_\rho) \subset \Delta(G).
\]
\item[(ii)] Fix $\rho > 0$, and let $S_\rho\subset G(\zero)$ be a maximal $\rho$-separated set. For every $x\in G(\zero)$ there exists $y_x\in S_\rho$ with $\dist(x,y_x) \leq \rho$. Then for each $y\in S_\rho$, we have
\[
\#\{ x\in G(\zero) : y_x = y \} \leq M_\rho,
\]
where $M_\rho$ is the value of the supremum \eqref{SD3}. Therefore for each $s\geq 0$ we have
\[
\Sigma_s(G)
= \sum_{x\in G(\zero)} b^{-s \dox x}
\asymp_\times \sum_{x\in G(\zero)} b^{-s\dox{y_x}}
\leq M_\rho \sum_{y\in S_\rho} b^{-s \dox y}
= M_\rho\Sigma_s(S_\rho).
\]
In particular, $\Sigma_s(G) < \infty$ if and only if $\Sigma_s(S_\rho) < \infty$, i.e. $\Delta(G) = \Delta(S_\rho)$. Intersecting over $\rho > 0$ and $S_\rho\subset G(\zero)$ yields $\Delta(G) = \w\Delta(G)$.
\item[(iii)] Take $\rho$ and $S_\rho$ such that $\delta(S_\rho) < \infty$.
\end{itemize}

Before proving (iv), we need a lemma:
\begin{lemma}
\label{lemma2tau}
Fix $\rho_1,\rho_2 > 0$ with $\rho_2\geq 2\rho_1$. Let $S_1\subset G(\zero)$ be a $\rho_1$-net,\Footnote{Here, as usual, a \emph{$\rho$-net} in a metric space $X$ is a subset $S\subset X$ such that $X = \thicken_\rho(S)$. Note that every maximal $\rho$-separated set is a $\rho$-net (but not conversely). \label{footnoterhonet}} and let $S_2\subset G(\zero)$ be a $\rho_2$-separated set. Then
\begin{equation}
\label{2tau}
\Delta(S_2) \subset \Delta(S_1).
\end{equation}
\end{lemma}
\begin{subproof}
Since $S_1$ is a $\rho_1$-net, for every $y\in S_2$, there exists $x_y\in S_1$ with $\dist(y,x_y) < \rho_1$. If $x_y = x_z$ for some $y,z\in S_2$, then $\dist(y,z) < 2\rho_1 \leq \rho_2$ and since $S_2$ is $\rho_2$-separated we have $y = z$. Thus the map $y\mapsto x_y$ is injective. It follows that for every $s\geq 0$, we have
\[
\Sigma_s(S_2)
= \sum_{y\in S_2} b^{-s \dox y}
\asymp_\times \sum_{y\in S_2} b^{-s \dox{x_y}}
\leq \sum_{x\in S_1} b^{-s \dox x}
= \Sigma_s(S_1),
\]
demonstrating \eqref{2tau}.
\end{subproof}

\begin{itemize}
\item[(iv)] The statement is trivial if $\w\delta_G = \infty$. So suppose that $\w\delta_G < \infty$, and let $\rho$, $S_\rho$ be as in (iii). Fix $\rho'\geq 2\rho$ and a maximal $\rho'$-separated set $S_{\rho'}\subset G(\zero)$, and we will show that $\Delta(S_{\rho'}) = \w\Delta_G$. The inclusion $\supset$ follows by definition. To prove the reverse direction, fix $\rho'' > 0$ and a maximal $\rho''$-separated set $S_{\rho''}$, and we will show that $\Delta(S_{\rho'})\subset \Delta(S_{\rho''})$.

Let $F = S_\rho\cap B(\zero,\rho'' + \rho)$; then $\#(F) < \infty$. We then set 
\[
S_\rho' := \bigcup_{x\in S_{\rho''}} g_x(F)
\] 
where for each $x\in S_{\rho''}$, $x = g_x(\zero)$. Then for all $s\geq 0$,
\begin{align*}
\Sigma_s(S_{\rho''}) = \sum_{x\in S_{\rho''}} b^{-s\dox x} 
&\asymp_\times \sum_{x\in S_{\rho''}} \sum_{y\in F} b^{-s\dox x}\\ 
&\asymp_\times \sum_{x\in S_{\rho''}} \sum_{y\in F} b^{-s\dox{g_x(y)}}\\ 
&= \Sigma_s(S_\rho')
\end{align*}
and therefore $\Delta(S_{\rho''}) = \Delta(S_\rho')$. But $S_\rho'$ is a $\rho$-net, so by Lemma \ref{lemma2tau}, we have $\Delta(S_{\rho'}) \subset \Delta(S_\rho')$. This completes the proof.
\end{itemize}
\end{proof}

Combining with Observation \ref{observationnotstronglydiscrete} yields the following:
\begin{corollary}
\label{corollarydeltawdelta}
Suppose that $G$ is a group. If $\Delta\neq\w\Delta$ then
\[
\w\delta < \delta = \infty.
\]
\end{corollary}

\begin{corollary}
\label{corollarygeneralizedtypes}
If a group $G$ is of convergence or divergence type, then it is also of generalized convergence or divergence type, respectively.
\end{corollary}

We will call a group $G\leq\Isom(X)$ \emph{Poincar\'e regular} if $\w\Delta_G = \Delta_G$, and \emph{Poincar\'e irregular} otherwise. A list of sufficient conditions for Poincar\'e regularity is given in Proposition \ref{propositionpoincareregular} below. Conversely, several examples of Poincar\'e irregular groups may be found in Section \ref{subsectionpoincareirregular}.


\chapter{Generalization of the Bishop--Jones theorem}
\label{sectionbishopjones}

In this chapter we prove Theorem \ref{theorembishopjonesmodified}, the first part of which states that if $G\prec\Isom(X)$ is a nonelementary semigroup, then
\begin{repequation}{bishopjonesmodified}
\HD(\Lr) = \HD(\Lur) = \HD(\Lur\cap\Lrsigma) = \w\delta
\end{repequation}
for some $\sigma > 0$. Our strategy is to prove that $\HD(\Lur\cap\Lrsigma)\leq \HD(\Lur) \leq \HD(\Lr) \leq \w\delta \leq \HD(\Lur\cap\Lrsigma)$ for some $\sigma > 0$. The first two inequalities are obvious. The third we prove now, and the proof of the fourth inequality will occupy \6\6\ref{subsectionpartitionstructures}-\ref{subsectionlemmastructure}.

\begin{lemma}
\label{lemmaonedirection}
For $G\prec\Isom(X)$, we have
\[
\HD(\Lr)\leq \w\delta.
\]
\end{lemma}
\begin{proof}
It suffices to show that for each $\sigma > 0$ and for each $s > \w\delta$, we have $\HD(\Lrsigma)\leq s$. Fix $\sigma > 0$ and $s > \w\delta$. Then there exists $\rho > 0$ and a maximal $\rho$-separated set $S_\rho \subset G(\zero)$ such that $s > \delta(S_\rho)$, which implies that $\Sigma_s(S_\rho) < \infty$. For each $x\in S_\rho$ let $\PP_x = \Shad(x,\sigma + \rho)$.
\begin{claim}
\label{claimonedirection}
\[
\xi\in \Lrsigma \Rightarrow \xi\in \PP_x\text{ for infinitely many }x\in S_\rho.
\]
\end{claim}
\begin{subproof}
Fix $\xi\in\Lrsigma$. Then there exists a sequence $g_n(\zero)\to\xi$ such that for all $n\in\N$ we have $\xi\in \Shad(g_n(\zero),\sigma)$. For each $n$, let $x_n\in S_\rho$ be such that $\dist(g_n(\zero),x_n)\leq \rho$; such an $x_n$ exists since $S_\rho$ is maximal $\rho$-separated. Then by (d) of Proposition \ref{propositionbasicidentities} we have $\xi\in \PP_{x_n} = \Shad(x_n,\sigma + \rho)$.

To complete the proof of Claim \ref{claimonedirection}, we need to show that the collection $(x_n)_1^\infty$ is infinite. Indeed, if $x_n\in F$ for some finite $F$ and for all $n\in\N$, then we would have $\dist(g_n(\zero),F)\leq\rho$ for all $n\in\N$. This would imply that the sequence $(g_n(\zero))_1^\infty$ is bounded, contradicting that $g_n(\zero)\to\xi$.
\end{subproof}
We next observe that by the Diameter of Shadows Lemma \ref{lemmadiameterasymptotic} we have
\[
\sum_{x\in S_\rho}\Diam^s(\PP_x)
\lesssim_{\times,\sigma,\rho} \sum_{x\in S_\rho}b^{-s\dox x}
= \Sigma_s(S_\rho) < \infty.
\]
Thus by the Hausdorff--Cantelli lemma \cite[Lemma 3.10]{BernikDodson}, we have $\HH^s(\Lrsigma) = 0$, and thus $\HD(\Lrsigma)\leq s$.
\end{proof}

\section{Partition structures} \label{subsectionpartitionstructures}
In this section we introduce the notion of a partition structure, an important technical tool for proving Theorem \ref{theorembishopjonesmodified}. We state some theorems about these structures, which will be proven in subsequent sections, and then use them to prove Theorem \ref{theorembishopjonesmodified}.
\Footnote{
\label{footnoteFSUcomparison}
Much of the material for this section has been taken (with modifications) from \cite[\65]{FSU4}. In \cite{FSU4} we also included as standing assumptions that $G$ was strongly discrete and of general type (see Definitions \ref{definitiondiscreteness} and \ref{definitionCCMT}). Thus some propositions which appear to have the exact same statement are in fact stronger in this monograph than in \cite{FSU4}. Specifically, this applies to Proposition \ref{propositionstructure} and Lemmas \ref{lemmaDSU} and \ref{sublemmaDSU}.
}

Throughout this section, $(Z,\Dist)$ denotes a metric space. We will constantly have in mind the special case $Z = \del X$, $\Dist = \Dist_{b,\zero}$.

\begin{notation}
Let
\[
\N^* = \bigcup_{n = 0}^\infty \N^n.
\]
If $\omega\in\N^*\cup\N^\N$, then we denote by $|\omega|$ the unique element of $\N\cup \{\infty\}$ such that $\omega\in \N^{|\omega|}$ and call $|\omega|$ the \emph{length} of $\omega$. For each $r\in\N$, we denote the initial segment of $\omega$ of length $r$ by
\[
\omega_1^r := (\omega_n)_1^r\in\N^r.
\]
For two words $\omega,\tau\in\N^\N$, let $\omega\wedge\tau$ denote their longest common initial segment, and let
\[
\dist_2(\omega,\tau) = 2^{-|\omega\wedge\tau|}.
\]
Then $(\N^\N,\dist_2)$ is a metric space.
\end{notation}
\begin{definition}
\label{definitiontree}
A \emph{tree} on $\N$ is a set $T\subset\N^*$ which is closed under initial segments. (Not to be confused with the various notions of ``trees'' introduced in Section \ref{subsectiongraphs}.)
\end{definition}
\begin{notation}
If $T$ is a tree on $\N$, then we denote its set of infinite branches by
\[
T(\infty) := \{\omega\in\N^\N:\omega_1^n\in T\all n\in\N\}.
\]
On the other hand, for $n\in\N$ we let
\[
T(n) := T\cap \N^n.
\]
For each $\omega\in T$, we denote the set of its children by
\[
T(\omega) := \{a\in\N:\omega a\in T\}.
\]
\end{notation}

\begin{definition}
\label{definitionpartitionstructure}
A \emph{partition structure} on $Z$ consists of a tree $T\subset\N^*$ together with a collection of closed subsets $(\PP_\omega)_{\omega\in T}$ of $Z$, each having positive diameter and enjoying the following properties:
\begin{itemize}
\item[(I)] If $\omega\in T$ is an initial segment of $\tau\in T$ then $\PP_\tau\subset \PP_\omega$. If neither $\omega$ nor $\tau$ is an initial segment of the other then $\PP_\omega\cap \PP_\tau = \emptyset$.
\item[(II)] For each $\omega\in T$ let
\[
D_\omega = \Diam(\PP_\omega).
\]
There exist $\kappa > 0$ and $0 < \lambda < 1$ such that for all $\omega\in T$ and for all $a\in T(\omega)$, we have 
\begin{equation}\label{kappa}
\Dist(\PP_{\omega a},Z\setminus \PP_\omega) \geq \kappa D_\omega
\end{equation}
and
\begin{equation}\label{lambda}
\kappa D_\omega \leq D_{\omega a} \leq \lambda D_\omega.
\end{equation}
\end{itemize}
Fix $s > 0$. The partition structure $(\PP_\omega)_{\omega\in T}$ is called \emph{$s$-thick} if for all $\omega\in T$, 
\begin{equation}
\label{redistribution}
\sum_{a\in T(\omega)}D_{\omega a}^s \geq D_\omega^s.
\end{equation}
\end{definition}

\begin{definition}
\label{definitionpartitionsubstructure}
Given a partition structure $(\PP_\omega)_{\omega\in T}$, a \emph{substructure} of $(\PP_\omega)_{\omega\in T}$ is a partition structure of the form $(\PP_\omega)_{\omega\in \w T}$, where $\w T\subset T$ is a subtree.
\end{definition}

\begin{observation}
Let $(\PP_\omega)_{\omega\in T}$ be a partition structure on a
complete metric space $(Z,\Dist)$. For each $\omega\in T(\infty)$, the set 
\[
\bigcap_{n\in\N}\PP_{\omega_1^n}
\]
is a singleton. If we define $\pi(\omega)$ to be the unique member of this set, then the map $\pi:T(\infty)\to Z$ is continuous. (In fact, it was shown in \cite[Lemma 5.11]{FSU4} that $\pi$ is quasisymmetric.)
\end{observation}

\begin{definition}
\label{definitionpartitionlimitset}
The set $\pi(T(\infty))$ is called the \emph{limit set} of the partition structure.
\end{definition}

We remark that a large class of examples of partition structures comes from the theory of conformal iterated function systems \cite{MauldinUrbanski1} (or in fact even graph directed Markov systems \cite{MauldinUrbanski2}) satisfying the strong separation condition (also known as the disconnected open set condition \cite{RiediMandelbrot}; see also \cite{FalconerMarsh}, where the limit sets of iterated function systems satisfying the strong separation condition are called \emph{dust-like}). Indeed, the notion of a partition structure was intended primarily to generalize these examples. The difference is that in a partition structure, the sets $(\PP_\omega)_\omega$ do not necessarily have to be defined by dynamical means. We also note that if $Z = \R^d$ for some $d\in\N$, and if $(\PP_\omega)_{\omega\in T}$ is a partition structure on $Z$, then the tree $T$ has bounded degree, meaning that there exists $N < \infty$ such that $\#(T(\omega))\leq N$ for every $\omega\in T$.

We will now state two propositions about partition structures and then use them to prove Theorem \ref{theorembishopjonesmodified}. Theorem \ref{theoremahlforsgeneral} will be proven below, and Proposition \ref{propositionstructure} will be proven in the following section.

\begin{theorem}[{\cite[Theorem 5.12]{FSU4}}]
\label{theoremahlforsgeneral}
Fix $s > 0$. Then any $s$-thick partition structure $(\PP_\omega)_{\omega\in T}$ on a complete metric space $(Z,\Dist)$ has a substructure $(\PP_\omega)_{\omega\in \w T}$ whose limit set is Ahlfors $s$-regular. Furthermore the tree $\w T$ can be chosen so that for each $\omega\in \w T$, we have that $\w T(\omega)$ is an initial segment of $T(\omega)$, i.e. $\w T(\omega) = T(\omega)\cap\{1,\ldots,N_\omega\}$ for some $N_\omega\in\N$. 
\end{theorem}

After these theorems about partition structures on an abstract metric space, we return to our more geometric setting of a Gromov triple $(X,\zero,b)$:

\begin{proposition}[Cf. {\cite[Lemma 5.13]{FSU4}}, Footnote \ref{footnoteFSUcomparison}]
\label{propositionstructure}
Let $G\prec\Isom(X)$ be nonelementary. Then for all $\sigma > 0$ sufficiently large and for every $0 < s < \w\delta_G$, there exist $\tau > 0$, a tree $T$ on $\N$, and an embedding $T\ni\omega\mapsto x_\omega\in G(\zero)$ such that if 
\[
\PP_\omega := \Shad(x_\omega,\sigma),
\]
then $(\PP_\omega)_{\omega\in T}$ is an $s$-thick partition structure on $(\del X,\Dist)$, whose limit set is a subset of $\Lurtau\cap\Lrsigma$.\end{proposition}

\begin{proof}[Proof of Theorem \ref{theorembishopjonesmodified} using Theorem \ref{theoremahlforsgeneral} \& Proposition \ref{propositionstructure}]
\hspace{0 in}

\noindent We first demonstrate the ``moreover'' clause. Fix $\sigma > 0$ large enough such that Proposition \ref{propositionstructure} holds. Fix $0 < s < \w\delta$, and let $(\PP_\omega)_{\omega\in T}$ be the partition structure guaranteed by Proposition \ref{propositionstructure}. Since this structure is $s$-thick, applying Theorem \ref{theoremahlforsgeneral} yields a substructure $(\PP_\omega)_{\omega\in \w T}$ whose limit set $\limitset_s\subset\Lurtau\cap\Lrsigma$ is Ahlfors $s$-regular, where $\tau > 0$ is as in Proposition \ref{propositionstructure}. Since $0 < s < \w\delta$ was arbitrary, this completes the proof of the ``moreover'' clause.

To demonstrate \eqref{bishopjonesmodified}, note that the inequality $\HD(\Lr)\leq\w\delta$ has already been established (Lemma \ref{lemmaonedirection}), and that the inequalities 
\[
\HD(\Lur\cap\Lrsigma) \leq \HD(\Lur) \leq \HD(\Lr)\] 
are obvious. Thus it suffices to show that 
\[
\HD(\Lur\cap\Lrsigma)\geq\w\delta.
\] 
But the mass distribution principle guarantees that 
\[
\HD(\Lur\cap\Lrsigma)\geq \HD(\limitset_s)\geq s
\] 
for each $0 < s < \w\delta$. This completes the proof.
\end{proof}

\begin{proof}[Proof of Theorem \ref{theoremahlforsgeneral}]
We will recursively define a sequence of maps 
\[
\mu_n:T(n)\to[0,1]
\]
with the following consistency property:
\begin{equation}
\label{consistency}
\mu_n(\omega) = \sum_{a\in T(\omega)}\mu_{n + 1}(\omega a).
\end{equation}
The Kolmogorov consistency theorem will then guarantee the existence of a measure $\w\mu$ on $T(\infty)$ satisfying 
\begin{equation}
\label{Kolmogorov}
\w\mu([\omega]) = \mu_n(\omega)
\end{equation}
for each $\omega\in T(n)$.

Let $c = 1 - \lambda^s > 0$, where $\lambda$ is as in \eqref{lambda}. For each $n\in\N$, we will demand of our function $\mu_n$ the following property: for all $\omega\in T(n)$, if $\mu_n(\omega) > 0$, then 
\begin{equation}
\label{regularity}
cD_\omega^s \leq \mu_n(\omega) < D_\omega^s.
\end{equation}

We now begin our recursion. For the case $n = 0$, let $\mu_0(\emptyset):= c D_\smallemptyset^s$; \eqref{regularity} is clearly satisfied.

For the inductive step, fix $n\in\N$ and suppose that $\mu_n$ has been constructed satisfying \eqref{regularity}. Fix $\omega\in T(n)$, and suppose that $\mu_n(\omega) > 0$. Formulas \eqref{redistribution} and \eqref{regularity} imply that
\[
\sum_{a\in T(\omega)}D_{\omega a}^s > \mu_n(\omega).
\]
Let $N_\omega\in T(\omega)$ be the smallest integer such that
\begin{equation}\label{regularitya}
\sum_{a\leq N_\omega}D_{\omega a}^s > \mu_n(\omega).\Footnote{Obviously, this and similar sums are restricted to $T(\omega)$.}
\end{equation}
Then the minimality of $N_\omega$ says precisely that
\[
\sum_{a\leq N_\omega - 1}D_{\omega a}^s \leq \mu_n(\omega).
\]
Using the above, \eqref{regularitya}, and \eqref{lambda}, we have
\begin{equation}
\label{sumaNbounds}
\mu_n(\omega) < \sum_{a\leq N_\omega}D_{\omega a}^s \leq \mu_n(\omega) +
D_{\omega N_\omega}^s \leq \mu_n(\omega) + \lambda^s D_\omega^s. 
\end{equation}
For each $a\in T(\omega)$ with $a > N_\omega$, let $\mu_{n + 1}(\omega a) =
0$, and for each $a \leq N_\omega$, let 
\[
\mu_{n + 1}(\omega a) = \frac{D_{\omega a}^s\mu_n(\omega)}{\sum_{b\leq N_\omega}D_{\omega b}^s}. 
\]
Obviously, $\mu_{n + 1}$ defined in this way satisfies \eqref{consistency}. Let us prove that \eqref{regularity} holds (of course, with $n = n + 1$). The second inequality follows directly from the definition of $\mu_{n + 1}$ and from \eqref{regularitya}. Using \eqref{sumaNbounds}, \eqref{regularity} (with $n = n$), and the equation $c = 1 - \lambda^s$, we deduce the first inequality as follows:
\begin{align*}
\mu_{n + 1}(\omega a)
&\geq \frac{D_{\omega a}^s\mu_n(\omega)}{\mu_n(\omega) + \lambda^sD_\omega^s}
= D_{\omega a}^s\left[1 - \frac{\lambda^sD_\omega^s}{\mu_n(\omega) + \lambda^sD_\omega^s}\right]\\ 
&\geq D_{\omega a}^s\left[1 - \frac{\lambda^s}{c + \lambda^s}\right] \\
&= cD_{\omega a}^s.
\end{align*}
The proof of \eqref{regularity} (with $n = n + 1$) is complete. This completes the recursive step.

Let
\[
\w T = \bigcup_{n = 1}^\infty\{\omega\in T(n):\mu_n(\omega) > 0\}.
\] 
Clearly, the limit set of the partition structure $(\PP_\omega)_{\omega\in \w T}$ is exactly the topological support of $\mu := \pi[\w\mu]$, where $\w\mu$ is defined by \eqref{Kolmogorov}. Furthermore, for each $\omega\in \w T$, we have $\w T(\omega) = T(\omega)\cap\{1,\ldots,N_\omega\}$. Thus, to complete the proof of Theorem \ref{theoremahlforsgeneral} it suffices to show that the measure $\mu$ is Ahlfors $s$-regular.

To this end, fix $z = \pi(\omega)\in\Supp(\mu)$ and $0 < r \leq \kappa D_\smallemptyset$, where $\kappa$ is as in \eqref{kappa} and \eqref{lambda}. For convenience of notation let 
\[
\PP_n := \PP_{\omega_1^n}, \;\; D_n := \Diam(\PP_n),
\]
and let $n\in\N$ be the largest integer such that $r < \kappa D_n$. We have
\begin{equation}\label{1fsdk9}
\kappa^2 D_n \leq \kappa D_{n + 1} \leq r < \kappa D_n.
\end{equation}
(The first inequality comes from \eqref{lambda}, whereas the latter two come from the definition of $r$.)

We now claim that
\[
B(z,r)\subseteq \PP_n.
\]
Indeed, by contradiction suppose that $w\in B(z,r)\butnot \PP_n$. By \eqref{kappa} we have
\[
\Dist(z,w)\geq \Dist(z,Z\butnot \PP_n)\geq \kappa D_n > r
\]
which contradicts the fact that $w\in B(z,r)$.

Let $k\in\N$ be large enough so that $\lambda^k\leq\kappa^2$. It follows from \eqref{1fsdk9} and repeated applications of the second inequality of \eqref{lambda} that
\[
D_{n + k} \leq \lambda^k D_n \leq \kappa^2 D_n \leq r,
\]
and thus
\[
\PP_{n + k} \subseteq B(z,r)\subseteq \PP_n.
\]
Thus, invoking \eqref{regularity}, we get
\begin{equation}
\label{muomegarbounds}
(1 - \lambda^s)D_{n + k}^s \leq \mu(\PP_{n + k}) \leq \mu(B(z,r))
\leq \mu(\PP_n) \leq D_n^s. 
\end{equation}

On the other hand, it follows from \eqref{1fsdk9} and repeated applications of the first inequality of \eqref{lambda} that
\begin{equation}
\label{above}
D_{n + k} \geq \kappa^k D_n \geq \kappa^{k - 1}r.
\end{equation}
Combining \eqref{1fsdk9}, \eqref{muomegarbounds}, and \eqref{above} yields
\[
(1 - \lambda^s)\kappa^{s(k - 1)}r^s \leq \mu(B(z,r)) \leq \kappa^{-2s}r^s,
\]
i.e. $\mu$ is Ahlfors $s$-regular. This completes the proof of Theorem \ref{theoremahlforsgeneral}.
\end{proof}

\section{A partition structure on $\del X$}
\label{subsectionlemmastructure}

We begin by stating our key lemma.

\begin{lemma}[Construction of children; cf. {\cite[Lemma 5.14]{FSU4}}, Footnote \ref{footnoteFSUcomparison}]
\label{lemmaDSU}
Let $G\prec\Isom(X)$ be nonelementary. Then for all $\sigma > 0$ sufficiently large, for every $0 < s < \w\delta_G$, for every $0 < \lambda < 1$, and for every $w\in G(\zero)$, there exists a finite subset $T(w)\subset G(\zero)$ (the \underline{children} of $w$) such that if we let
\begin{align*}
\PP_x &:= \Shad(x,\sigma)\\
D_x &:= \Diam(\PP_x)
\end{align*}
then the following hold:
\begin{itemize}
\item[(i)] The family $(\PP_x)_{x\in T(w)}$ consists of pairwise disjoint shadows contained in $\PP_w$.
\item[(ii)] There exists $\kappa > 0$ independent of $w$ such that for all $x\in T(w)$,
\begin{align*}
\Dist(\PP_x,\del X\butnot \PP_w) &\geq \kappa D_w\\
\kappa D_w \leq D_x &\leq \lambda D_w.
\end{align*}
\item[(iii)]
\[
\sum_{x\in T(w)}D_x^s \geq D_w^s.
\]
\end{itemize}
\end{lemma}

\begin{figure}
\begin{center}
\begin{tabular}{@{}ll|@{}}

\begin{tikzpicture}[line cap=round,line join=round,>=triangle 45,scale=0.85]
\clip(-3.3,-2.7) rectangle (3.5,2.7);
\draw(0.0,-0.0) circle (2.6129676614916217cm);
\draw(1.2444064831035602,0.7132854804669598) circle (0.3679391468439666cm);
\draw(1.7882027204795168,1.4262864278568905) circle (0.1351026446829843cm);
\draw(2.0629778937273007,1.1515112546091066) circle (0.11246824164327565cm);
\draw(2.1202227214872558,0.7965933224973857) circle (0.12035035832634534cm);
\draw (0,0)-- (1.85778232681272,1.8374560746269037);
\draw (0,0)-- (2.5244367854319023,0.6744026366787416);
\draw (2.156504665055924,0.6818421598881995)-- (2.491401614304384,0.7877296466697858);
\draw (2.0719679454251,0.9068461471164001)-- (2.3937366645507177,1.0476758949153115);
\draw (2.113056993311799,1.0508077831538636)-- (2.3396374875948704,1.1634846052448466);
\draw (2.003548920208303,1.2469958213453562)-- (2.218386998924586,1.3807096447129514);
\draw (1.86606101634726,1.315874425045114)-- (2.1354355680057964,1.5058269936804682);
\draw (1.6978674333274448,1.5267466712145967)-- (1.9429624621724708,1.747139625384784);
\begin{scriptsize}
\draw [fill=black] (0.0,-0.0) circle (1pt);
\draw[color=black] (-0.16181839193637584,-0.10332614941794697) node {$\zero$};
\draw [fill=black] (1.3742581787476496,0.7888486360608391) circle (0.75pt);
\draw[color=black] (1.3,-0.85) node {$g(\zero)$};
\draw[thin,->,>=stealth] (1.37,-0.58) -- (1.37,0.67);
\draw [fill=black] (1.8429159652907778,1.4693036915224922) circle (0.5pt);
\draw [fill=black] (2.122308107268412,1.1854052246742528) circle (0.5pt);
\draw [fill=black] (2.1711576420553347,0.8167939855719956) circle (0.5pt);
\end{scriptsize}
\end{tikzpicture}

\begin{tikzpicture}[line cap=round,line join=round,>=triangle 45,scale=0.36]
\clip(-1.59,-2.97) rectangle (1.57,2.72);
\draw [shift={(0.0,0.0)}] plot[domain=0.7354397676755055:2.3211211596591195,variable=\t]({1.0*1.6991762710207554*cos(\t r)+-0.0*1.6991762710207554*sin(\t r)},{0.0*1.6991762710207554*cos(\t r)+1.0*1.6991762710207554*sin(\t r)});
\draw (1.212984563827671,1.3292510070809178)-- (1.26,1.14);
\draw (1.0801070470812233,1.1485375843057488)-- (1.26,1.14);
\begin{scriptsize}
\draw[color=black] (0.36,2.34) node {$g^{-1}$};
\end{scriptsize}
\end{tikzpicture}

\begin{tikzpicture}[line cap=round,line join=round,>=triangle 45,scale=0.85]
\clip(-3.1,-2.7) rectangle (3.2,2.8);
\draw(0.0,-0.0) circle (2.611235899552007cm);
\draw(0.0,-0.0) circle (1.1901105590689844cm);
\draw [shift={(-3.3794741494806733,0.7796181967120758)}] plot[domain=-0.7176127868835138:0.6257318419197768,variable=\t]({1.0*2.279326680308009*cos(\t r)+-0.0*2.279326680308009*sin(\t r)},{0.0*2.279326680308009*cos(\t r)+1.0*2.279326680308009*sin(\t r)});
\draw [shift={(-1.7623065030786575,-3.0487890331884815)}] plot[domain=0.2114358851577187:1.5278846809760946,variable=\t]({1.0*2.3316956611045625*cos(\t r)+-0.0*2.3316956611045625*sin(\t r)},{0.0*2.3316956611045625*cos(\t r)+1.0*2.3316956611045625*sin(\t r)});
\draw(0.16511694177641764,2.046955552778598) circle (0.30789402561067564cm);
\draw(1.6180556349350896,-1.1660470857511598) circle (0.31806155761471966cm);
\draw(1.8661183386451068,1.0310797185375715) circle (0.27576581560448304cm);
\draw [shift={(-3.923946814679157,2.3225255909878184)}] plot[domain=-0.06724121069882294:0.07518735908343538,variable=\t]({1.0*3.7904448753248996*cos(\t r)+-0.0*3.7904448753248996*sin(\t r)},{0.0*3.7904448753248996*cos(\t r)+1.0*3.7904448753248996*sin(\t r)});
\draw [shift={(-4.64426430538874,3.825227515133297)}] plot[domain=5.9298153627837396:6.043103905803468,variable=\t]({1.0*5.4354771773979484*cos(\t r)+-0.0*5.4354771773979484*sin(\t r)},{0.0*5.4354771773979484*cos(\t r)+1.0*5.4354771773979484*sin(\t r)});
\draw [shift={(-1.060629324509099,-5.220932382140508)}] plot[domain=0.8585679463877474:0.9872882702824342,variable=\t]({1.0*4.5417198944649*cos(\t r)+-0.0*4.5417198944649*sin(\t r)},{0.0*4.5417198944649*cos(\t r)+1.0*4.5417198944649*sin(\t r)});
\draw [shift={(-0.21908948347075194,-6.549804754789044)}] plot[domain=1.1275322215585417:1.2404427016586532,variable=\t]({1.0*6.006520366955727*cos(\t r)+-0.0*6.006520366955727*sin(\t r)},{0.0*6.006520366955727*cos(\t r)+1.0*6.006520366955727*sin(\t r)});
\draw [shift={(-30.82463631611864,52.86866068223544)}] plot[domain=5.274920930454794:5.282870997078253,variable=\t]({1.0*61.009065644834784*cos(\t r)+-0.0*61.009065644834784*sin(\t r)},{0.0*61.009065644834784*cos(\t r)+1.0*61.009065644834784*sin(\t r)});
\draw [shift={(-51.465499934296915,125.44723710829089)}] plot[domain=5.117413840983604:5.120957170067834,variable=\t]({1.0*135.64048286553125*cos(\t r)+-0.0*135.64048286553125*sin(\t r)},{0.0*135.64048286553125*cos(\t r)+1.0*135.64048286553125*sin(\t r)});
\begin{scriptsize}
\draw [fill=black] (0.0,-0.0) circle (1pt);
\draw[color=black] (0.08094925045079697,0.16792479370770388) node {$\zero$};
\draw [fill=black] (-1.6622803095029393,-0.719239845767749) circle (0.75pt);
\draw[color=black] (-1.8311186992417487,-0.9150924600470811) node {$g^{-1}(\zero)$};
\draw [fill=black] (0.17980173196364405,2.1717367802041783) circle (0.5pt);
\draw [fill=black] (2.0000614646072825,1.1099186028287176) circle (0.5pt);
\draw [fill=black] (1.7294585171156684,-1.2410351367538701) circle (0.5pt);
\end{scriptsize}
\end{tikzpicture}
\end{tabular}
\caption[The construction of children]{The strategy for the proof of Lemma \ref{lemmaDSU}. To construct a collection of ``children'' of the point $w = g(\zero)$, we ``pull back'' the entire picture via $g^{-1}$. In the pulled-back picture, the Big Shadows Lemma \ref{lemmabigshadow} guarantees the existence of many points $x\in G(\zero)$ such that $\Shad_z(x,\sigma)\subset \Shad_z(\zero,\sigma)$, where $z = g^{-1}(\zero)$. (Cf. Lemma \ref{sublemmaDSU} below.) These children can then be pushed forward via $g$ to get children of $w$.}
\label{figureexperiment}
\end{center}
\end{figure}
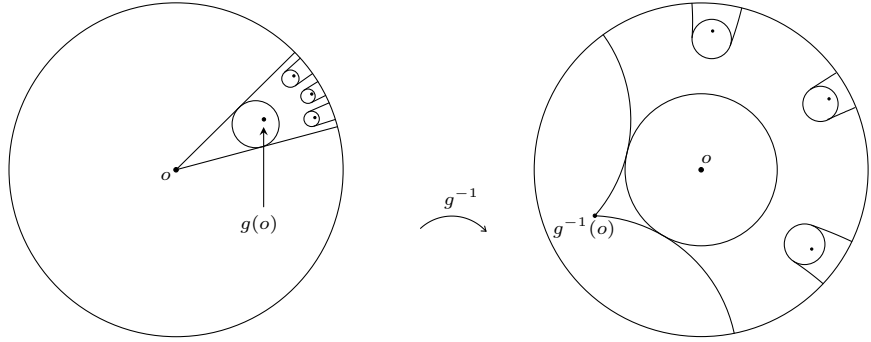

It is not too hard to deduce Proposition \ref{propositionstructure} from Lemma \ref{lemmaDSU}. We do it now:
\begin{proof}[Proof of Proposition \ref{propositionstructure} assuming Lemma \ref{lemmaDSU}]
Let $\sigma > 0$ be large enough so that Lemma \ref{lemmaDSU} holds. Fix $0 < s < \w\delta$, and let $\lambda = 1/2$. For each $w\in G(\zero)$, let $(y_n(w))_{n = 1}^{N(w)}$ be an enumeration of $T(w)$. Define a tree $T\subset\N^*$ and a collection $(x_\omega)_{\omega\in T}$ inductively as follows:
\begin{align*}
x_\smallemptyset &= \zero\\
T(\omega) &= \{1,\ldots,N(x_\omega)\}\\
x_{\omega a} &= y_a(x_\omega).
\end{align*}
Then the conclusion of Lemma \ref{lemmaDSU} precisely implies that $(\PP_\omega := \PP_{x_\omega})_{\omega\in T}$ is an $s$-thick partition structure on $(\del X,\Dist)$.

To complete the proof, we must show that the limit set of the partition structure $(\PP_\omega)_{\omega\in T}$ is contained in $\Lurtau\cap\Lrsigma$ for some $\tau > 0$. Indeed, fix $\omega\in T(\infty)$. Then for each $n\in\N$, $\pi(\omega)\in \PP_{\omega_1^n}=\Shad(x_{\omega_1^n},\sigma)$ and $\dox{x_{\omega_1^n}}\to \infty$. So, the sequence $(x_{\omega_1^n})_1^\infty$ converges $\sigma$-radially to $\pi(\omega)$. On the other hand,
\begin{align*}
\dist(x_{\omega_1^n},x_{\omega_1^{n + 1}}) &\asymp_{\plus,\sigma} \busemann_\zero(x_{\omega_1^{n + 1}},x_{\omega_1^n}) \by{\eqref{distbusemann}}\\
&\asymp_{\plus,\sigma} -\log_b\left(\frac{D_{\omega_1^{n + 1}}}{D_{\omega_1^n}}\right) \by{Lemma \ref{lemmadiameterasymptotic}}\\
&\leq_{\phantom{\plus,\sigma}} -\log_b(\kappa) \asymp_{\plus,\kappa} 0. \by {\eqref{lambda}}
\end{align*}
Thus the sequence $(x_{\omega_1^n})_1^\infty$ converges to $\pi(\omega)$ $\tau$-uniformly radially, where $\tau$ depends only on $\sigma$ and $\kappa$ (which in turn depends on $s$).
\end{proof}

The proof of Lemma \ref{lemmaDSU} will proceed through a series of lemmas.



\begin{lemma}[Cf. {\cite[Lemma 5.15]{FSU4}}]
\label{lemmatdivergent}
Fix $\tau > 0$, and let $S_\tau\subset G(\zero)$ be a maximal $\tau$-separated subset. Let $B\subset\bord X$ be an open set which intersects $\Lambda$. Then for every $0 < t < \w\delta$, the series
\[
\Sigma_t(S_\tau\cap B)
\]
diverges.
\end{lemma}
\begin{proof}

By Proposition \ref{propositionminimal2}, there exists a loxodromic isometry $g\in G$ such that $g_+\in B$. Let $\ell(g) = \log_b g'(g_-) = -\log_b g'(g_+) > 0$, and let the functions $r = r_{g_+,g_-,\zero}, \theta = \theta_{g_+,g_-,\zero}$ be as in Section \ref{subsectionpolar}. Fix $N\in\N$ large to be determined, let $\kappa = N\ell(g)$, and for each $n\in\Z$ let
\[
C_n = \{x\in X: n\kappa \leq r(x) < (n + 1)\kappa\}
\]
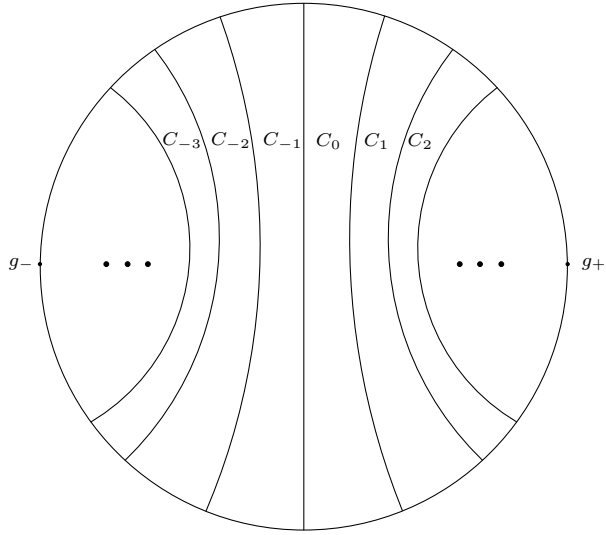
\begin{figure}
\begin{center}
\begin{tikzpicture}[line cap=round,line join=round,>=triangle 45,scale=1.2]
\clip(-3.28,-3.04) rectangle (3.35,2.93);
\draw(0.0,-0.0) circle (2.92cm);
\draw (0.0,-2.92)-- (0.0,2.9199999999999995);
\draw [shift={(-7.9807385594024085,0.22846692826476303)}] plot[domain=-0.40293664452450617:0.34569769451385884,variable=\t]({1.0*7.497384482406696*cos(\t r)+-0.0*7.497384482406696*sin(\t r)},{0.0*7.497384482406696*cos(\t r)+1.0*7.497384482406696*sin(\t r)});
\draw [shift={(-4.358450318907789,0.3145920313756659)}] plot[domain=-0.8023891554417188:0.6582794941876258,variable=\t]({1.0*3.4231925825980647*cos(\t r)+-0.0*3.4231925825980647*sin(\t r)},{0.0*3.4231925825980647*cos(\t r)+1.0*3.4231925825980647*sin(\t r)});
\draw [shift={(-3.50200819090357,0.20509253639486452)}] plot[domain=-1.0351538170621417:0.9181589198842997,variable=\t]({1.0*2.240662432657915*cos(\t r)+-0.0*2.240662432657915*sin(\t r)},{0.0*2.240662432657915*cos(\t r)+1.0*2.240662432657915*sin(\t r)});
\draw [shift={(8.568849511994815,0.3086758429165517)}] plot[domain=2.830100098137627:3.5251001108755213,variable=\t]({1.0*8.061553144764865*cos(\t r)+-0.0*8.061553144764865*sin(\t r)},{0.0*8.061553144764865*cos(\t r)+1.0*8.061553144764865*sin(\t r)});
\draw [shift={(4.358450318907789,0.3145920313756659)}] plot[domain=2.4833131594021673:3.943981809031512,variable=\t]({1.0*3.4231925825980647*cos(\t r)+-0.0*3.4231925825980647*sin(\t r)},{0.0*3.4231925825980647*cos(\t r)+1.0*3.4231925825980647*sin(\t r)});
\draw [shift={(3.50200819090357,0.20509253639486452)}] plot[domain=2.2234337337054937:4.176746470651935,variable=\t]({1.0*2.2406624326579156*cos(\t r)+-0.0*2.2406624326579156*sin(\t r)},{0.0*2.2406624326579156*cos(\t r)+1.0*2.2406624326579156*sin(\t r)});
\begin{scriptsize}
\draw[color=black] (0.273556687080068,1.39) node {$C_0$};
\draw[color=black] (-0.23313923665707456,1.39) node {$C_{-1}$};
\draw[color=black] (-0.81,1.39) node {$C_{-2}$};
\draw[color=black] (-1.35,1.39) node {$C_{-3}$};
\draw[color=black] (0.8049694851458516,1.39) node {$C_1$};
\draw[color=black] (1.29,1.39) node {$C_2$};
\draw [fill=black] (1.727358295319231,0.028247325459469327) circle (0.75pt);
\draw [fill=black] (1.9520571528523318,0.028247325459469327) circle (0.75pt);
\draw [fill=black] (2.1879909532620876,0.028247325459469327) circle (0.75pt);
\draw [fill=black] (-1.727358295319231,0.028247325459469327) circle (0.75pt);
\draw [fill=black] (-1.9520571528523318,0.028247325459469327) circle (0.75pt);
\draw [fill=black] (-2.1879909532620876,0.028247325459469327) circle (0.75pt);
\draw [fill=black] (-2.923,0.028247325459469327) circle (0.5pt);
\draw[color=black] (-3.127,0.028247325459469327) node {$g_-$};
\draw [fill=black] (2.923,0.028247325459469327) circle (0.5pt);
\draw[color=black] (3.217,0.028247325459469327) node {$g_+$};
\end{scriptsize}
\end{tikzpicture}
\caption[The sets $C_n$, for $n\in\Z$]{The sets $C_n$, for $n\in\Z$.}
\label{figurepartitionintofour}
\end{center}
\end{figure}
(cf. Figure \ref{figurepartitionintofour}). Let
\begin{align*}
C_{+,0} = \bigcup_{\substack{n\geq 0 \\ \text{even}}}C_n,&\hspace{.5 in}
C_{+,1} = \bigcup_{\substack{n\geq 0 \\ \text{odd}}}C_n\\
C_{-,0} = \bigcup_{\substack{n < 0 \\ \text{even}}}C_n,&\hspace{.5 in}
C_{-,1} = \bigcup_{\substack{n < 0 \\ \text{odd}}}C_n.
\end{align*}
Fix $\rho > 0$, and let $S_\rho\subset G(\zero)$ be a maximal $\rho$-separated set. Since $\Sigma_t(S_\rho) = \infty$, one of the series $\Sigma_t(S_\rho\cap C_{+,0})$, $\Sigma_t(S_\rho\cap C_{+,1})$, $\Sigma_t(S_\rho\cap C_{-,0})$, and $\Sigma_t(S_\rho\cap C_{-,1})$ must diverge. By way of illustration let us consider the case where
\[
\Sigma_t(S_\rho\cap C_{-,0}) = \infty.
\]
Let
\[
A_\rho = \bigcup_{\substack{n < 0 \\ \text{even}}}g^{2N |n|}(C_n\cap S_\rho).
\]
\begin{claim}
$\Sigma_t(A_\rho) = \infty$.
\end{claim}
\begin{subproof}
Fix $n = -m < 0$ even and $x\in C_n$. Then by \eqref{rtranslation},
\[
r(g^{2Nm}(x)) \asymp_+ 2Nm\ell(g) + r(x) = 2m\kappa + r(x) \asymp_{\plus,\kappa} 2m\kappa -m\kappa = m\kappa
\]
and thus
\[
|r(g^{2Nm}(x))| \asymp_{\plus,\kappa} |r(x)|.
\]
On the other hand, by \eqref{thetatranslation} we have $\theta(g^{2Nm}(x)) \asymp_\plus \theta(x)$. Combining with Lemma \ref{lemmafabs} gives
\[
\dist(0,x) \asymp_{\plus,\kappa} \dist(0,g^{2Nm}(x)).
\]
Thus
\begin{align*}
\Sigma_t(A_\rho)
&= \sum_{\substack{m > 0 \\ \text{even}}}\sum_{x\in C_{-m}\cap S_\rho} b^{-t\dox{g^{2Nm}(x)}}\\
&\asymp_{\times,\kappa} \sum_{\substack{m > 0 \\ \text{even}}}\sum_{x\in C_{-m}\cap S_\rho} b^{-t\dox x}\\
&= \Sigma_t(C_{-,0}\cap S_\rho) = \infty.
\end{align*}
\end{subproof}
\begin{claim}
$A_\rho$ is a $\rho$-separated set.
\end{claim}
\begin{subproof}
Fix $y_1,y_2\in A_\rho$. Then for some $m_1,m_2 > 0$ even, we have $x_i := g^{-2Nm_i}(y_i)\in C_{-m_i}$ ($i = 1,2$). If $n_1 = n_2$, then we have
\[
\dist(y_1,y_2) = \dist(x_1,x_2) \geq \rho,
\]
since $x_1,x_2\in S_\rho$ and $S_\rho$ is $\rho$-separated. So suppose $n_1\neq n_2$; without loss of generality we may assume $n_1 > n_2$. Then by \eqref{rtranslation} we have
\begin{align*}
r(y_1) - r(y_2) &\asymp_\plus 2Nm_1\ell(g) + r(x_1) - (2Nm_2\ell(g) + r(x_2))\\
&=_\pt 2\kappa(m_1 - m_2) + r(x_1) - r(x_2)\\
&\geq_\pt 2\kappa[m_1 - m_2] + \kappa(-m_1) - \kappa(-m_2 + 1)\\
&=_\pt \kappa(m_1 - m_2 - 1)\\
&\geq_\pt \kappa = N\ell(g).
\end{align*}
By choosing $N$ sufficiently large, we may guarantee that $r(y_1) - r(y_2)\geq \rho$, which implies $\dist(y_1,y_2)\geq \rho$.
\end{subproof}

For all $x\in \thicken_\rho(A_\rho)$, we have $r(x)\gtrsim_\plus 0$. Thus $g_-\notin \cl{\thicken_\rho(A_\rho)}$. So by Theorem \ref{theoremloxodromicextra}, we can find $n\in\N$ such that $\thicken_\rho(g^n(A_\rho))\subset B$.

Let $S_{\rho/2}\subset G(\zero)$ be a maximal $\rho/2$-separated set. By Lemma \ref{lemma2tau}, we have
\[
\Sigma_t(S_{\rho/2}\cap B) \gtrsim_\times \Sigma_t(g^n(A_\rho)) \asymp_\times \Sigma_t(A_\rho) = \infty.
\]
Since $\rho > 0$ was arbitrary, this completes the proof.
\end{proof}

\begin{lemma}[Cf. {\cite[Sublemma 5.17]{FSU4}}, Footnote \ref{footnoteFSUcomparison}]
\label{sublemmaDSU}
Let $B\subset\bord X$ be an open set which intersects $\Lambda$. For all $\sigma > 0$ sufficiently large and for all $0 < s < \w\delta$, there exists a set $S_B\subset G(\zero)\cap B$ such that for all $z\in X\butnot B$,
\begin{itemize}
\item[(i)] If
\[
\PP_{z,x} := \Shad_z(x,\sigma),
\]
then the family $(\PP_{z,x})_{x\in S_B}$ consists of pairwise disjoint shadows contained in $\PP_{z,\zero}\cap B$.
\item[(ii)] There exists $\kappa > 0$ independent of $z$ (but depending on $s$) such that for all $x\in S_B$,
\begin{align} \label{epsilonnew}
\Dist_{b,z}(\PP_{z,x},\del X\butnot \PP_{z,\zero}) &\geq \kappa\Diam_z(\PP_{z,\zero})\\ \label{lambdanew}
\kappa\Diam_z(\PP_{z,\zero}) \leq \Diam_z(\PP_{z,x}) &\leq \lambda\Diam_z(\PP_{z,\zero}).
\end{align}
\item[(iii)]
\[
\sum_{x\in S_B}\Diam_z^s(\PP_{z,x}) \geq \Diam_z^s(\PP_{z,\zero}).
\]
\end{itemize}
\end{lemma}
\begin{proof}
Let $B\subset\bord X$ be an open set which contains a point $\eta\in\Lambda$. Choose $\rho > 0$ large enough so that
\[
\{x\in\bord X:\lb x|\eta\rb_\zero\geq \rho\} \subset B.
\]
Then fix $\sigma > 0$ large to be determined, depending only on $\rho$. Fix $\w\rho\geq\rho$ large to be determined, depending only on $\rho$ and $\sigma$.

Fix $0 < s < \w\delta$ and $z\in X\butnot B$. For all $x\in X$ we have
\[
0 \asymp_{\plus,\rho} \lb z|\eta\rb_\zero \gtrsim_\plus \min(\lb x|\eta\rb_\zero,\lb x|z\rb_\zero).
\]
Let
\[
\w B = \{x\in X:\lb x|\eta\rb_\zero \geq \w\rho\ew\}.
\]
If $\w\rho$ is chosen large enough, then we have
\begin{equation}\label{xz0}
\lb x|z\rb_\zero \asymp_{\plus,\rho} 0,
\end{equation}
for all $x\in \w B$. We emphasize that the implied constants of these asymptotics are independent of both $z$ and $s$.

For each $n\in\N$ let 
\[
A_n := B(\zero,n)\butnot B(\zero,n - 1)
\] 
be the $n$th annulus centered at $\zero$. We shall need the following variant of the Intersecting Shadows Lemma:

\begin{claim}
\label{claimtau}
There exists $\tau > 0$ depending on $\rho$ and $\sigma$ such that for all $n\in\N$ and for all $x,y\in A_n\cap \w B$, if
\[
\PP_{z,x}\cap \PP_{z,y}\neq\emptyset,
\]
then 
\[
\dist(x,y) < \tau.
\] 
\end{claim}
\begin{subproof}
Without loss of generality suppose $\dist(z,y)\geq \dist(z,x)$. Then by the Intersecting Shadows Lemma \ref{lemmatau} we have
\[
\dist(x,y) \asymp_{\plus,\sigma} \busemann_z(y,x)
= \busemann_\zero(y,x) + 2\lb x|z\rb_\zero - 2\lb y|z\rb_\zero.
\]
Now $|\busemann_\zero(y,x)| \leq 1$ since $x,y\in A_n$. On the other hand, since $x,y\in\w B$, we have
\[
\lb x|z\rb_\zero \asymp_{\plus,\rho} \lb y|z\rb_\zero \asymp_{\plus,\rho} 0.
\]
Combining gives
\[
\dist(x,y) \asymp_{\plus,\rho,\sigma} 0,
\]
and letting $\tau$ be the implied constant finishes the proof.
\end{subproof}

Fix $M > 0$ large to be determined, depending on $\rho$ and $\tau$ (and thus implicitly on $\sigma$). Let $S_\tau\subset G(\zero)$ be a maximal $\tau$-separated set. Fix $t\in (s,\w\delta)$; then by Lemma \ref{lemmatdivergent} we have
\begin{align*}
\infty
= \Sigma_t(S_\tau\cap \w B)
&=_\pt \sum_{n = 1}^\infty\Sigma_t(S_\tau\cap\w B\cap A_n)\\
&=_\pt \sum_{n = 1}^\infty\sum_{x\in S_\tau\cap\w B\cap A_n}b^{-t\dox x}\\
&\asymp_\times \sum_{n = 1}^\infty b^{-(t - s)n}\sum_{x\in S_\tau\cap\w B\cap A_n} b^{-s\dox x}.
\end{align*}
It follows that there exist arbitrarily large numbers $n\in\N$ such that
\begin{equation}
\label{ndef}
\sum_{x\in S_\tau\cap\w B\cap A_n} b^{-s\dox x} \geq M.
\end{equation}
Fix such an $n$, also to be determined, depending on $\lambda$, $\rho$, $\w\rho$, and $M$ (and thus implicitly on $\tau$ and $\sigma$), and let $S_B = S_\tau\cap \w B\cap A_n$. To complete the proof, we demonstrate (i)-(iii).
\begin{subproof}[Proof of \textup{(i)}]
In order to see that the shadows $(\PP_{z,x})_{x\in S_B}$ are pairwise disjoint, suppose that $x,y\in S_B$ are such that $\PP_{z,x}\cap \PP_{z,y}\neq\emptyset$. By Claim \ref{claimtau} we have $\dist(x,y) < \tau$. Since $S_B$ is $\tau$-separated, this implies $x = y$.

Fix $x\in S_B$. Using \eqref{xz0} and the fact that $x\in A_{n}$, we have
\[
\lb\zero|z\rb_x \asymp_\plus \dox x - \lb x|z\rb_\zero\asymp_{\plus,\rho} \dox x\asymp_\plus n.
\]
Thus for all $\xi\in\PP_{z,x}$,
\[
0
\asymp_{\plus,\sigma} \lb z|\xi\rb_x
\gtrsim_\plus \min(\lb\zero|z\rb_x,\lb\zero|\xi\rb_x)
\asymp_\plus \min(n,\lb\zero|\xi\rb_x);
\]
taking $n$ sufficiently large (depending on $\sigma$), this gives
\[
\lb\zero|\xi\rb_x \asymp_{\plus,\sigma} 0,
\]
from which it follows that
\[
\lb x|\xi\rb_\zero
\asymp_\plus d(\zero,x)-\lb\zero|\xi\rb_x
\asymp_{\plus,\sigma} n.
\]
Therefore, since $x\in\w B$, we get
\[
\lb\xi|\eta\rb_\zero
\gtrsim_\plus \min(\lb x|\xi\rb_\zero,\lb x|\eta\rb_\zero)
\gtrsim_{\plus,\sigma} \min(n,\w\rho).
\]
Thus $\xi\in B$ as long as $\w\rho$ and $n$ are large enough (depending on $\sigma$). Thus $\PP_{z,x}\subset B$.

Finally, note that we do not need to prove that $\PP_{z,x}\subset \PP_{z,\zero}$, since it is implied by \eqref{epsilonnew} which we prove below.

\ignore{
Next we shall prove that
\begin{equation}\label{BaB0}
\PP_{z,x}\subset \PP_{z,\zero}.
\end{equation}
Indeed, fix $\xi\in \PP_{z,x}$. Then by \eqref{1fsdk13} and \eqref{2fsdk13} we get that
\[
\sigma
\gtrsim_\plus \lb\zero|\xi\rb_x
\gtrsim_\plus \min(\lb\zero|z\rb_x,\lb z|\xi\rb_x)
\asymp_\plus \min(n,\lb z|\xi\rb_x).
\]
Thus if $n$ is large enough (depending on $\sigma$) then
\[
\lb\zero|\xi\rb_x\gtrsim_\plus \lb z|\xi\rb_x;
\]
it follows from this, \eqref{xz0}, and (z) of Proposition \ref{propositionbasicidentities} that
\begin{equation}\label{2fsdk15}
\lb z|\xi\rb_\zero
\asymp_\plus \lb z|\xi\rb_x + \lb x|z\rb_\zero - \lb\zero|\xi\rb_x
\asymp_\plus \lb z|\xi\rb_x - \lb\zero|\xi\rb_x
\lesssim_\plus 0.
\end{equation}
So, if $\sigma > 0$ is taken large enough (not depending on anything), then $\lb z|\xi\rb_\zero < \sigma$, meaning that $\xi\in \PP_{z,\zero}$. The inclusion \eqref{BaB0} is proved, and in consequence the whole item (i).
}
\end{subproof}
\begin{subproof}[Proof of \textup{(ii)}]
Take any $x\in S_B$. Then by \eqref{xz0}, we have
\begin{equation}
\label{asympan1}
\dist(x,z) - \dox z = \dox x - 2\lb x|z\rb_\zero \asymp_{\plus,\rho} \dox x \asymp_\plus n.
\end{equation}
Combining with the Diameter of Shadows Lemma \ref{lemmadiameterasymptotic} gives
\begin{equation}
\label{asympan1new}
\frac{\Diam_z(\PP_{z,x})}{\Diam_z(\PP_{z,\zero})}
\asymp_{\times,\sigma} \frac{b^{-\dist(z,x)}}{b^{-\dist(z,\zero)}}
\asymp_{\times,\rho} b^{-n}.
\end{equation}
Thus by choosing $n$ sufficiently large depending on $\sigma$, $\lambda$, and $\rho$ (and satisfying \eqref{ndef}), we guarantee that the second inequality of \eqref{lambdanew} holds. On the other hand, once $n$ is chosen, \eqref{asympan1new} guarantees that if we choose $\kappa$ sufficiently small, then the first inequality of \eqref{lambdanew} holds.

In order to prove \eqref{epsilonnew}, let $\xi\in \PP_{z,x}$ and let $\gamma\in \del X\setminus \PP_{z,\zero}$. We have
\begin{align*}
\lb x|\xi\rb_z &\asymp_\plus \dist(x,z) - \lb z|\xi\rb_x \geq \dist(x,z) - \sigma\\
\lb \zero|\gamma\rb_z &\asymp_\plus \dox z - \lb\zero|\xi\rb_x \leq \dox z - \sigma.
\end{align*}
Also, by \eqref{xz0} we have
\[
\lb \zero|x\rb_z \asymp_\plus \dox z - \lb x|z\rb_\zero \asymp_{\plus,\rho} \dox z.
\]
Applying Gromov's inequality twice and then applying \eqref{asympan1} gives
\begin{align*}
\dox z - \sigma \gtrsim_\plus \lb \zero|\gamma\rb_z
&\gtrsim_{\plus\phantom{,\rho}} \min\left(\lb\zero|x\rb_z,\lb x|\xi\rb_z,\lb \xi|\gamma\rb_z\right)\\
&\gtrsim_{\plus,\rho} \min\left(\dox z,\dist(x,z) - \sigma,\lb \xi|\gamma\rb_z\right)\\
&\asymp_{\plus\phantom{,\rho}} \min\left(\dox z,\dox z + n - \sigma,\lb \xi|\gamma\rb_z\right).
\end{align*}
By choosing $n$ and $\sigma$ sufficiently large (depending on $\rho$), we can guarantee that neither of the first two expressions can represent the minimum without contradicting the inequality. Thus
\[
\dox z - \sigma \gtrsim_{\plus,\rho} \lb\xi|\gamma\rb_z;
\]
exponentiating and the Diameter of Shadows Lemma \ref{lemmadiameterasymptotic} give
\[
\Dist_{b,z}(\xi,\gamma)
\gtrsim_{\times,\rho} b^{-(\dox z - \sigma)}
\asymp_{\times,\sigma} b^{-\dox z}
\asymp_{\times,\sigma} \Diam_z(\PP_{z,\zero}).
\]
Thus we may choose $\kappa$ small enough, depending on $\rho$ and $\sigma$, so that \eqref{epsilonnew} holds.
\ignore{
Note that $\lb\gamma|z\rb_\zero\lesssim_\plus 0$ as $\gamma\in \del X\setminus \PP_{z,\zero}$. Therefore, using \eqref{2fsdk15}, we get
\[
0
\asymp_\plus \lb z|\xi\rb_\zero
\gtrsim_\plus \min(\lb\xi|\gamma\rb_\zero,\lb z|\gamma\rb_\zero)
\gtrsim_\plus \min(\lb\xi|\gamma\rb_\zero,\sigma).
\]
So, $\lb\xi|\gamma\rb_\zero\lesssim_\plus 0$. Combining with (d) of Proposition \ref{propositionbasicidentities}, we have
\[
\lb\xi|\gamma\rb_z \lesssim_\plus d(\zero,z).
\]
Inserting this in turn to \eqref{1fsdk19} and applying the Diameter of Shadows Lemma \ref{lemmadiameterasymptotic} gives
\[
\Dist_{b,z}(\xi,\gamma)\gtrsim_\times b^{-d(\zero,z)}\asymp_\times\Diam_z(\PP_{z,\zero}).
\]
Thus,
\[
\Dist_{b,z}(\PP_{z,x},\del X\setminus \PP_{z,\zero})\gtrsim_\times\Diam_z(\PP_{z,\zero}),
\]
finishing the proof of \eqref{epsilonnew} and thus of item (ii). 
}
\end{subproof}
\begin{subproof}[Proof of \textup{(iii)}]
\begin{align*}
\sum_{x\in S_B}\Diam_z^s(\PP_{z,x})
&\asymp_{\times\phantom{,\rho}} \sum_{x\in S_B}b^{-s\dist(z,x)} \by{the Diameter of Shadows Lemma}\\
&\asymp_{\times,\rho} \;\; b^{-s\dist(z,\zero)}\sum_{x\in S_B}b^{-\dox x} \hspace{-0.85in}\by{\eqref{asympan1}}\\
&\geq_{\phantom{\times,\rho}}Mb^{-s\dist(z,\zero)} \by{\eqref{ndef}}\\
&\asymp_{\times\phantom{,\rho}} M\Diam_z^s(\PP_{z,\zero}). \by{the Diameter of Shadows Lemma}
\end{align*}
Letting $M$ be larger than the implied constant yields the result.
\end{subproof}
\end{proof}

We may now complete the proof of Lemma \ref{lemmaDSU}:

\begin{proof}[Proof of Lemma \ref{lemmaDSU}]
Let $\eta_1,\eta_2\in\Lambda$ be distinct points, and let $B_1$ and $B_2$ be disjoint neighborhoods of $\eta_1$ and $\eta_2$, respectively. Let $\sigma > 0$ be large enough so that Lemma \ref{sublemmaDSU} holds for both $B_1$ and $B_2$. Fix $0 < s < \w\delta$, and let $S_1\subset G(\zero)\cap B_1$ and $S_2\subset G(\zero)\cap B_2$ be the sets guaranteed by Lemma \ref{sublemmaDSU}. Now suppose that $w = g_w(\zero)\in G(\zero)$. Let $z = g_w^{-1}(\zero)$. Then either $z\notin B_1$ or $z\notin B_2$; say $z\notin B_i$. Let $T(w) = g_w(S_i)$; then (i)-(iii) of Lemma \ref{sublemmaDSU} exactly guarantee (i)-(iii) of Lemma \ref{lemmaDSU}.
\end{proof}

\section{Sufficient conditions for Poincar\'e regularity}
\label{subsectionpropositionpoincareregular}
We end this chapter by relating the modified Poincar\'e exponent to the classical Poincar\'e exponent under certain additional assumptions, thus completing the proof of Theorem \ref{theorembishopjonesregular}.

\begin{proposition}
\label{propositionpoincareregular}
Let $G\leq\Isom(X)$ be nonelementary, and assume either that
\begin{itemize}
\item[(1)] $X$ is regularly geodesic and $G$ is moderately discrete,
\item[(2)] $X$ is an algebraic hyperbolic space and $G$ is weakly discrete, or that
\item[(3)] $X$ is an algebraic hyperbolic space and $G$ is COT-discrete and acts irreducibly.
\end{itemize}
Then $G$ is Poincar\'e regular.
\end{proposition}

\begin{remark}
Example \ref{exampleAutTBIM} shows that Proposition \ref{propositionpoincareregular} cannot be improved by replacing ``COT'' with ``UOT'', Example \ref{examplepoincareextension} shows that Proposition \ref{propositionpoincareregular} cannot be improved by removing the assumption that $G$ acts irreducibly, Example \ref{exampleAutT} shows that Proposition \ref{propositionpoincareregular} cannot be improved by removing the hypothesis that $X$ is an algebraic hyperbolic space from (ii), and Example \ref{exampleAutTparttwo} shows that Proposition \ref{propositionpoincareregular} cannot be improved by removing the assumption that $X$ is regularly geodesic.
\end{remark}

We begin with the following observation:
\begin{observation}
If (3) implies that $G$ is Poincar\'e regular, then (2) does as well.
\end{observation}
\begin{proof}
Suppose (2) holds, and let $S$ be the smallest totally geodesic subset of $\bord X$ which contains $\LambdaG$ (cf. Lemma \ref{lemmatotallygeodesic}). Since $G$ is nonelementary, $V := S\cap X$ is nonempty; it is clear that $V$ is $G$-invariant. By Observation \ref{observationrestrictions}, the action $G\given V$ is weakly discrete. By Proposition \ref{propositionparametricdiscreteness}, $G\given V$ is COT-discrete. Furthermore, $G$ acts irreducibly on $V$ because of the way $V$ was defined (cf. Proposition \ref{propositionactsreducibly}). Thus (3) holds for the action $G\given V$, which by our hypothesis implies $\Delta_G = \w\Delta_G$ (since the Poincar\'e set and modified Poincar\'e set are clearly stable under restrictions).
\end{proof}

We now proceed to prove that (1) and (3) each imply that $G$ is Poincar\'e regular.

By contradiction, let us suppose that $G$ is Poincar\'e irregular. By Proposition \ref{propositionbasicmodified}(ii), we have that $G$ is not strongly discrete and thus
\[
\w\delta_G < \delta_G = \infty.
\]
This gives us two contrasting behaviors: On one hand, by Proposition \ref{propositionbasicmodified}(iii), there exist $\rho > 0$ and a maximal $\rho$-separated set $S_\rho\subset G(\zero)$ so that $S_\rho$ does not contain an bounded infinite set. On the other hand, since $G$ is not strongly discrete, there exists $\sigma > 0$ such that $\#(G_\sigma) = \infty$, where
\[
G_\sigma := \{g\in G:g(\zero)\in B(\zero,\sigma)\}.
\]
\begin{claim}
\label{claimprecompact}
For every $\xi\in\LambdaG$, the orbit $G_\sigma(\xi)$ is precompact.
\end{claim}
\begin{proof}
Suppose not. Then the set $\cl{G_\sigma(\xi)}$ is complete (with respect to the metric $\Dist$) but not compact. It follows that $\cl{G_\sigma(\xi)}$, and thus also $G_\sigma(\xi)$, is not totally bounded. So there exists $\epsilon > 0$ and an infinite $\epsilon$-separated subset $(g_n(\xi))_1^\infty$.

Fix $L$ large to be determined. Since $\xi\in\Lambda$, we can find $x\in G(\zero)$ such that $\lb x | \xi \rb_\zero \geq L$.

\begin{subclaim}
By choosing $L$ large enough we can ensure
\[
\dist(g_m(x),g_n(x)) \geq 2\rho \all m,n\in\N.
\]
\end{subclaim}
\begin{subproof}
By (d) of Proposition \ref{propositionbasicidentities},
\[
\lb g_n(x) | g_n(\xi) \rb_\zero \asymp_{\plus,\sigma} \lb g_n(x) | g_n(\xi) \rb_{g_n(\zero)} = \lb x|\xi\rb_\zero \geq L ,
\]
and thus 
\[
\Dist(g_n(x),g_n(\xi)) \lesssim_{\times,\sigma} b^{-L} .
\]
If $L$ is large enough, then this implies
\[
\Dist(g_n(x),g_n(\xi)) \leq \epsilon/3.
\]
Since by construction the sequence $(g_n(\xi))_1^\infty$ is $\epsilon$-separated, we also have
\[
\Dist(g_m(\xi),g_n(\xi)) \geq \epsilon
\]
and then the triangle inequality gives
\[
\Dist(g_m(x),g_n(x)) \geq \epsilon/3,
\]
or, taking logarithms,
\[
\lb g_m(x) | g_n(x) \rb_\zero \lesssim_\plus \ -\log_b(\epsilon/3).
\]
Now we also have
\[
\dox{g_n(x)} \asymp_{\plus,\sigma} \dox x \geq \lb x|\xi\rb_\zero \geq L
\]
and therefore
\begin{align*}
\dist(g_m(x),g_n(x))
&=_{\phantom{\plus,\sigma}} \dox{g_m(x)} + \dox{g_n(x)} - 2\lb g_m(x)|g_n(x)\rb_\zero\\
&\gtrsim_{\plus,\sigma} 2L - 2(- \log_b (\epsilon/3)).
\end{align*}
Thus by choosing $L$ sufficiently large, we ensure that $\dist(g_m(x),g_n(x)) \geq 2\rho$.

\end{subproof}
Recall that $S_\rho$ is a maximal $\rho$-separated set. Thus for each $n\in\N$, we can find $y_n\in S_\rho$ with $\dist(g_n(x),y_n) < \rho$. Then the subclaim implies $y_n\neq y_m$ for $n\neq m$. But on the other hand
\[
\dox{y_n} \leq \dox x + \sigma + \rho \all n\in\N,
\]
which implies that $S_\rho$ contains a bounded infinite set, a contradiction.
\end{proof}

We now proceed to disprove the hypotheses (1) and (3) of Proposition \ref{propositionpoincareregular}. Thus if either of these hypotheses are assumed, we have a contradiction which finishes the proof.

\begin{proof}[Proof that \text{(1)} cannot hold]
Since $G$ is assumed to be nonelementary, we can find distinct points $\xi_1,\xi_2\in \LambdaG$. By Claim \ref{claimprecompact}, there exist a sequence $(g_n)_1^\infty$ in $G_\sigma$ and points $\eta_1,\eta_2\in\LambdaG$ such that
\[
g_n(\xi_i)\tendsto n \eta_i.
\]
Next, choose a point $x\in\geo{\xi_1}{\xi_2}$. For each $n\in\N$, we have
\[
g_n(x)\in\geo{g_n(\xi_1)}{g_n(\xi_2)}.
\]
Thus since $X$ is regularly geodesic there exist a sequence $(n_k)_1^\infty$ and a point $z\in \geo{\eta_1}{\eta_2}$ such that
\[
g_{n_k}(x)\tendsto k z.
\]
Since $g_n\in G_\sigma \all n$, the sequence $(g_n(x))_1^\infty$ is bounded and thus $z\in X$. By contradiction, suppose that $G$ is moderately discrete, and fix $\epsilon > 0$ satisfying \eqref{MD2}. For all $m,n\in\N$ large enough so that $g_m(x),g_n(x)\in B(z,\epsilon/2)$, we have $\dist(x,g_m^{-1} g_n(x)) = \dist(g_m(x),g_n(x)) \leq \epsilon$. Thus for some $N\in\N$, we have 
\[
\#\{g_m^{-1} g_n : m,n\geq N\} < \infty.
\] 
This is clearly a contradiction.
\end{proof}

\begin{proof}[Proof that \text{(3)} cannot hold]
Now we assume that $X$ is an algebraic hyperbolic space, say $X = \H = \H_\F^\alpha$, and that $G$ acts irreducibly on $X$. Using the identification
\[
\Isom(\H) \equiv \O^*(\LL;\QQ)/\sim,
\]
(Theorem \ref{theoremisometries}), for each $g\in G_\sigma$ let $T_g\in\O^*(\LL;\QQ)$ be a representative of $g$. Recall (Lemma \ref{lemmaoperatornorm}) that
\[
\|T_g\| = \|T_g^{-1}\| = e^{\dogo g},
\]
so since $g\in G_\sigma$ we have $\|T_g\| = \|T_g^{-1}\| \leq b^\sigma$. In particular, the family $(T_g)_{g\in G_\sigma}$ acts equicontinuously on $\LL$.

For simplicity of exposition, in the following proof we will assume that $X$ is separable. (In the non-separable case, the reader should use nets instead of sequences.) It follows that $\LambdaG\subset\del X$ is also separable; let $(\xi_k = [\xx_k])_1^\infty$ be a dense sequence, with $\xx_k\in\LL$, $\|\xx_k\| = 1$.

\begin{claim}
There exists a sequence of distinct elements $(g_n)_1^\infty$ in $G_\sigma$ such that the following hold:
\begin{align*}
T_{g_n}[\xx_k] &\;\tendsto n\; \yy_k^{(+)}\in\LL\butnot\{\0\} \\
T_{g_n}^{-1}[\xx_k] &\;\tendsto n\; \yy_k^{(-)}\in\LL\butnot\{\0\}\\
\sigma(T_{g_n}) &\tendsto n \sigma\in\Aut(\F).
\end{align*} 
\end{claim}
\begin{subproof}
For each $k\in\N$ let
\[
\KK_k = \{\yy\in\LL\butnot\{\0\}: [\yy]\in \cl{G_\sigma(\xi_k)} \text{ and } b^{-\sigma} \leq \|\yy\| \leq b^\sigma\},
\]
and let
\[
\KK := \left(\prod_{k\in\N}\KK_k\right)^2 \times \Aut(\F).
\]
Then by Claim \ref{claimprecompact} (and general topology), $\KK$ is a compact metrizable space, and is in particular sequentially compact. Now for each $g\in G_\sigma$,
\[
b^{-\sigma} \leq \|T_g[\xx_k]\| \leq b^\sigma \text{ and } b^{-\sigma} \leq \|T_g^{-1}[\xx_k]\| \leq b^\sigma 
\]
and thus
\[
\phi_g := \left((T_g(\xx_k))_1^\infty, (T_g^{-1}(\xx_k))_1^\infty, \sigma(T_g)\right)\in\KK,
\]
and so since $\#(G_\sigma) = \infty$, there exists a sequence of distinct elements $(g_n)_1^\infty$ in $G_\sigma$ so that the sequence $(\phi_{g_n})_1^\infty$ converges to a point
\[
\left((\yy_k^{(+)})_1^\infty, (\yy_k^{(-)})_1^\infty, \sigma\right)\in\KK.
\]
Writing out what this means yields the claim.
\end{subproof}

Let $T_n = T_{g_n}$ and $\sigma_n = \sigma(T_n) \to \sigma$. We claim that the sequence $(T_n)_1^\infty$ is convergent in the strong operator topology. Let
\begin{align*}
\K &= \{a\in\F : \sigma(a) = a\}\\
V &= \{\xx\in\LL:\text{the sequence $(T_n[\xx])_1^\infty$ converges}\}.
\end{align*}
Then $\K$ is an $\R$-subalgebra of $\F$, and $V$ is a $\K$-module. Given $\xx,\yy\in V$, by Observation \ref{observationBQpreserving} we have
\[
\sigma_n(B_\QQ(\xx,\yy)) = B_\QQ(T_n\xx,T_n\yy) \tendsto n B_\QQ(\xx,\yy),
\]
so $B(\xx,\yy)\in \K$. Thus $V$ satisfies \eqref{totallygeodesic}. On the other hand, since the family $(T_n)_1^\infty$ acts equicontinuously on $\LL$, the set $V$ is closed. Thus $[V]\cap\bord X$ is totally geodesic. But by construction, $\xi_k\in[V]$ for all $k$, and so
\[
\LambdaG \subset [V]
\]
Therefore, since by hypothesis $G$ acts irreducibly, it follows that $[V] = X$, i.e. $V = \LL$ (Proposition \ref{propositionactsreducibly}). So for every $\xx\in\LL$, the sequence $(T_n(\xx))_1^\infty$ converges. Thus
\[
T_n \tendsto n T^{(+)}\in L(\LL)
\]
in the strong operator topology. (The boundedness of the operator $T^{(+)}$ follows from the uniform boundedness of the operators $(T_n)_1^\infty$.) We do not yet know that $T^{(+)}$ is invertible. But a similar argument yields that
\[
T_n^{-1} \tendsto n T^{(-)}\in L(\LL),
\]
and since the sequences $(T_n)_1^\infty$ and $(T_n^{-1})_1^\infty$ are equicontinuous, we have
\[
T^{(+)}T^{(-)} = \lim_{n\to\infty}T_n T_n^{-1} = I,
\]
and similarly $T^{(-)}T^{(+)} = I$. Thus $T^{(+)}$ and $T^{(-)}$ are inverses of each other and in particular
\[
T^{(+)}\in \O^*(\LL;\QQ).
\]
Let $h = [T^{(+)}]\in\Isom(X)$. By Proposition \ref{propositionCOTSOT}, we have $g_n\to h$ in the compact-open topology. Thus, Lemma \ref{lemmadiscretesubgroup} completes the proof.
\end{proof}

\part{Examples}
\label{partexamples}

This part will be divided as follows: In Chapter \ref{sectionschottky} we consider semigroups of isometries which can be written as the ``Schottky product'' of two subsemigroups. Next, we analyze in detail the class of parabolic groups of isometries in Chapter \ref{sectionparabolic}. 
In Chapter \ref{sectionGF}, we define a subclass of the class of groups of isometries which we call \emph{geometrically finite}, generalizing known results from the Standard Case.
In Chapter \ref{sectionexamples}, we provide a list of examples whose main importance is that they are counterexamples to certain implications; however, these examples are often geometrically interesting in their own right. 
Finally, in Chapter \ref{sectionRtrees}, we consider methods of constructing $\R$-trees which admit natural group actions, including what we call the ``stapling method''.

%
%

\chapter{Schottky products} \label{sectionschottky}

An important tool for constructing examples of discrete subgroups of $\Isom(X)$ is the technique of \emph{Schottky products}. Schottky groups are a special case of Schottky products; cf. Definition \ref{definitionschottkygroup}. In this chapter we explain the basics of Schottky products on hyperbolic metric spaces, and give several important examples. We intend to give a more comprehensive account of Schottky products in \cite{DSU2}, where we will study their relation to pseudo-Markov systems (defined in \cite{StratmannUrbanski3}).

\begin{remark}
Throughout this chapter, $E$ denotes an index set with at least two elements. There are no other restrictions on the cardinality of $E$; in particular, $E$ may be infinite.
\end{remark}

\section{Free products}
\label{subsectionfreeproducts}
We provide a brief review of the theory of free products, mainly to fix notation. Let $(\Gamma_a)_{a\in E}$ be a collection of nontrivial abstract semigroups. Let
\[
\Gamma_E = \coprod_{a\in E} (\Gamma_a\butnot\{e\}) = \bigcup_{a\in E} \{a\}\times (\Gamma_a\butnot\{e\}).
\]
Let $(\Gamma_E)^*$ denote the set of finite words with letters in $\Gamma_E$, including the empty word, which we denote by $\emptyset$. The \emph{free product} of $(\Gamma_a)_{a\in E}$, denoted $\ast_{a\in E}\Gamma_a$, is the set
\[
\left\{\gg = (a_1,\gamma_1)\cdots(a_n,\gamma_n)\in (\Gamma_E)^* : a_i\neq a_{i + 1} \all i = 1,\ldots,n - 1, \; n \geq 0\right\},
\]
together with the operation of multiplication defined as follows: To multiply two words $\gg,\hh\in \ast_{a\in E} \Gamma_a$, begin by concatenating them. The concatenation may no longer satisfy $a_i\neq a_{i + 1}$ for all $i$; namely, this condition may fail at the point where the two words are joined. \emph{Reduce} the concatenated word $\gg\ast\hh$ using the rule
\[
(a,\gamma_1)(a,\gamma_2) = \begin{cases}
(a,\gamma_1 \gamma_2) & \gamma_1 \gamma_2 \neq e\\
\emptyset & \gamma_1 \gamma_2 = e
\end{cases}.
\]
The word may require multiple reductions in order to satisfy $a_i\neq a_{i + 1}$. The reduced form of $\gg\ast\hh$ will be denoted $\gg\hh$.

One verifies that the operation of multiplication defined above is associative, so that the free product $\ast_{a\in E}\Gamma_a$ is a semigroup. If $(\Gamma_a)_{a\in E}$ are groups, then $\ast_{a\in E}\Gamma_a$ is a group. The inclusion maps $\iota_a: \Gamma_a\to \ast_{a\in E}\Gamma_a$ defined by $\iota_a(\gamma) = (a,\gamma)$ are homomorphisms, and $\ast_{a\in E}\Gamma_a = \lb \iota_a(\Gamma_a)\rb_{a\in E}$.

An important fact about free products is their \emph{universal property}: Given any semigroup $\Gamma$ and any collection of homomorphisms $(\pi_a: \Gamma_a\to \Gamma)$, there exists a unique homomorphism $\pi:\ast_{a\in E}\Gamma_a\to \Gamma$ such that $\pi_a = \pi\circ\iota_a$ for all $a$. For example, if $(\Gamma_a)_{a\in E}$ are subsemigroups of $\Gamma$ and $(\pi_a)_{a\in E}$ are the identity inclusions, then $\pi((a_1,\gamma_1)\cdots(a_n,\gamma_n)) = \gamma_1\cdots \gamma_n$. We will call the map $\pi$ the \emph{natural map} from $\ast_{a\in E} \Gamma_a$ to $\Gamma$.

\begin{remark}
\label{remarkfreegroup}
We will use the notation $\Gamma_1\ast\cdots\ast \Gamma_n$ to denote $\ast_{a\in\{1,\ldots,n\}}\Gamma_a$. The semigroups
\[
\F_n(\Z) = \underbrace{\Z\ast \cdots\ast \Z}_{n\text{ times}} \text{ and } \F_n(\N) = \underbrace{\N\ast\cdots\ast\N}_{n\text{ times}}
\]
are called the \emph{free group on $n$ elements} and the \emph{free semigroup on $n$ elements}, respectively.
\end{remark}

\section{Schottky products}
\label{subsectionschottkyproducts}
Given a collection of semigroups $G_a \prec \Isom(X)$, we can ask whether the semigroup $\lb G_a\rb_{a\in E}\prec\Isom(X)$ is isomorphic to the free product $\ast_{a\in E} G_a$. A sufficient condition for this is that the groups $(G_a)_{a\in E}$ are in \emph{Schottky position}.

\begin{definition}
\label{definitionschottkyproduct}
A collection of nontrivial semigroups $(G_a \prec\Isom(X) )_{a\in E}$ is in \emph{Schottky position} if there exist disjoint open sets $U_a\subset\bord X$ satisfying:
\begin{itemize}
\item[(I)] For all $a,b\in E$ distinct and $g\in G_a\butnot\{\id\}$,
\[
g(U_b) \subset U_a.
\]
\item[(II)] There exists $\zero\in X\butnot \cl{\bigcup_{a\in E} U_a}$ satisfying
\begin{equation}
\label{zerodef}
g(\zero)\in U_a \all a\in E \all g\in G_a\butnot\{\id\}.
\end{equation}
\end{itemize}
Such a collection $(U_a)_{a\in E}$ is called a \emph{Schottky system} for $(G_a)_{a\in E}$. If the collection $(G_a)_{a\in E}$ is in Schottky position, then we will call the semigroup $G = \lb G_a\rb_{a\in E}$ the \emph{Schottky product} of $(G_a)_{a\in E}$.

A Schottky system will be called \emph{global} if for all $a\in E$ and $g\in G_a\butnot\{\id\}$,
\begin{equation}
\label{schottkyglobal}
g(\bord X\butnot U_a) \subset \schottkycl{U_a}.
\end{equation}
\end{definition}

\begin{remark*}
In most references (e.g. \cite[\65]{DPPS}), \eqref{schottkyglobal} or a similar hypothesis is taken as the definition of Schottky position. So what these references call a ``Schottky group'', we would call a ``global Schottky group''. There are important examples of Schottky semigroups which are not global; see e.g. (B) of Proposition \ref{propositionnonelementaryequivalent}. It should be noted that such examples tend to be semigroups rather than groups, which explains why references which consider only groups can afford to include globalness in their definition of Schottky position.
\end{remark*}

\begin{remark*}
The above definition may be slightly confusing to someone familiar with classical Schottky groups, since in that context the sets $U_a$ in the above definition are not half-spaces but rather unions of pairs of half-spaces; cf. Definition \ref{definitionschottkygroup}.
\end{remark*}

The basic properties of Schottky products are summarized in the following lemma:

\begin{lemma}
\label{lemmapingpong}
Let $G = \lb G_a\rb_{a\in E}$ be a Schottky product. Then:
\begin{itemize}
\item[(i)] (Ping-Pong Lemma) The natural map $\pi:\ast_{a\in E} G_a\to G$ is an injection (and therefore an isomorphism).
\item[(ii)] Fix $\gg = (a_1,g_1)(a_2,g_2)\cdots(a_n,g_n)\in\ast_{a\in E}G_a$, and let $g = \pi(\gg)$. Then
\begin{equation}
\label{pingpong}
g(\zero)\in U_{a_1}.
\end{equation}
Moreover, for all $b\neq a_n$
\begin{equation}
\label{pingpong2}
g(U_b) \subset U_{a_1},
\end{equation}
and if the system $(U_a)_{a\in E}$ is global
\begin{equation}
\label{pingpong3}
g(\bord X\butnot U_{a_n}) \subset \schottkycl{U_{a_1}}.
\end{equation}
\item[(iii)] If $G$ is a group, then $G$ is COT-discrete.
\end{itemize}
\end{lemma}
\begin{proof}
\eqref{pingpong}-\eqref{pingpong3} may be proven by an easy induction argument. Now \eqref{pingpong} immediately demonstrates (i), since it implies that $\pi(\gg)\neq\id$. (iii) also follows from \eqref{pingpong}, since it shows that $\dogo g$ is bounded from below for all $g\in G\butnot\{\id\}$.
\end{proof}
\begin{remark}
Lemma \ref{lemmapingpong}(i) says that Schottky products are (isomorphic to) free products. However, we warn the reader that the converse is not necessarily true; cf. Lemma \ref{lemmafreenotschottky}.
\end{remark}

Two important classes of Schottky products are \emph{Schottky groups} and \emph{Schottky semigroups}.

\begin{definition}
\label{definitionschottkygroup}
A \emph{Schottky group} is the Schottky product of cyclic groups $G_a = (g_a)^\Z$ with the following property: For each $a\in E$, $U_a$ may be written as the disjoint union of two sets $U_a^+$ and $U_a^-$ satisfying
\[
g_a(\bord X\butnot U_a^-) \subset \schottkycl{U_a^+}.
\]
A \emph{Schottky semigroup} is simply the Schottky product of cyclic semigroups; no additional hypotheses are needed.
\end{definition}

\begin{remark}
\label{remarkhalfspace}
In the classical theory of Schottky groups, the sets $U_a^\pm$ are required to be half-spaces. A \emph{half-space} in $\bord\H^\alpha$ is a connected component of the complement of a totally geodesic subset of $\bord\H^\alpha$ of codimension one. Requiring the sets $U_a^\pm$ to be half-spaces has interesting effects on the geometry of Schottky groups.

Although the notion of half-spaces cannot be generalized to hyperbolic metric spaces in general or even to nonreal algebraic hyperbolic spaces (since a totally geodesic subspace of an algebraic hyperbolic space over $\F = \C$ or $\Q$ always has real codimension at least 2, so deleting it yields a connected set), it at least makes sense over real hyperbolic spaces and in the context of $\R$-trees. A \emph{half-space} in an $\R$-tree $X$ is a connected component of the complement of a point in $X$.

We hope to study the effect of requiring the sets $U_a^\pm$ to be half-spaces, both in the case of real (but infinite-dimensional) algebraic hyperbolic spaces and in the case of $\R$-trees, in more detail in \cite{DSU2}.
\end{remark}

\section{Strongly separated Schottky products}
Many questions about Schottky products cannot be answered without some additional information. For example, one can ask whether or not the Schottky product of strongly (resp. moderately, weakly) discrete groups is strongly (resp. moderately, weakly) discrete. One can also ask about the relation between the Poincar\'e exponent of a Schottky group and the Poincar\'e exponent of its factors. 

For the purposes of this monograph, we will be interested in Schottky products which satisfy the following condition:

\begin{definition}
\label{definitionstrongseparation}
A Schottky product $G = \lb G_a\rb_{a\in E}$ is said to be \emph{strongly separated} (with respect to a Schottky system $(U_a)_{a\in E}$) if there exists $\epsilon > 0$ such that for all $a,b\in E$ distinct and $g\in G_a\butnot\{\id\}$,
\begin{equation}
\label{strongseparation}
\wbar\Dist\big(U_a\cup g^{-1}(\bord X\butnot U_a),U_b\big) \geq \epsilon.
\end{equation}
Here $\wbar\Dist$ is as in Proposition \ref{propositionwbarDist}. Abusing terminology, we will also call the semigroup $G$ and the Schottky system $(U_a)_{a\in E}$ strongly separated.

The product $G = \lb G_a\rb_{a\in E}$ is \emph{weakly separated} if \eqref{strongseparation} holds for a constant $\epsilon > 0$ which depends on $a$ and $b$ (but not on $g$).
\end{definition}

\begin{remark}
There are many important examples of Schottky products which are not strongly separated, and we hope to analyze these in more detail in \cite{DSU2}. Some examples of Schottky products that do satisfy the condition are given in Section \ref{subsectionschottkyexamples}.
\end{remark}

\begin{standingassumptions}
For the remainder of this chapter, 
\[
G = \lb G_a\rb_{a\in E}
\] 
denotes a strongly separated Schottky product and $(U_a)_{a\in E}$ denotes the corresponding Schottky system. Moreover, from now on we assume that the hyperbolic metric space $X$ is \emph{geodesic}.
\end{standingassumptions}

\begin{notation}
\label{notationVa}
Let $\Gamma$ denote the free product $\Gamma = \ast_{a\in E} G_a$, and let $\pi:\Gamma\to G$ denote the natural isomorphism. Whenever we have specified an element $\gg\in\Gamma$, we denote its length by $|\gg|$ and we write $\gg = (a_1,g_1)\cdots(a_{|\gg|},g_{|\gg|})$. For $\hh\in\Gamma$, we write $\hh = (b_1,h_1)\cdots(b_{|\hh|},h_{|\hh|})$.

Let $\zero\in X$ satisfy \eqref{zerodef}. Let $\epsilon\leq \dist(\zero,\bigcup_a U_a)$ satisfy \eqref{strongseparation}, and for each $a\in E$ let $V_a$ denote the closed $\epsilon/4$-thickening of $U_a$ with respect to the $\wbar\Dist$ metric. Then the sets $(\Int(V_a))_{a\in E}$ are also a Schottky system for $(G_a)_{a\in E}$; they are strongly separated with $\epsilon = \epsilon/2$; moreover,
\begin{equation}
\label{UaVa}
\wbar\Dist(U_a,\bord X\butnot V_a) \geq \epsilon/2 \all a\in E.
\end{equation}
Finally, let
\[
X_a = \begin{cases}
\bord X \butnot \Int(V_a) & \text{$(U_a)_{a\in E}$ is global}\\
\{\zero\}\cup\bigcup_{b\neq a} V_b & \text{otherwise}
\end{cases},
\]
so that
\begin{equation}
\label{XaUa}
g(X_a) \subset \schottkycl{U_a} \all a\in E.
\end{equation}
Note that since the sets $(V_a)_{a\in E}$ are $\epsilon/2$-separated, they have no accumulation points and thus $X_a$ is closed for all $a\in E$.
\end{notation}

The strong separation condition will allow us to relate the discreteness of the groups $G_a$ to the discreteness of their Schottky product $G$. It will also allow us to relate the Poincar\'e exponents of $G_a$ with the Poincar\'e exponent of $G$. The underlying fact which will allow us to prove both of these relations is the following lemma:

\begin{lemma}
\label{lemmastrongseparation}
There exist constants $C,\epsilon > 0$ such that for all $\gg\in\Gamma$,
\begin{equation}
\label{quasiisometry1}
\sum_{i = 1}^{|\gg|} (\dogo{g_i} - C) \vee \epsilon \leq \dist(X\butnot V_{a_1},\pi(\gg)(X_{a_{|\gg|}})) \leq \sum_{i = 1}^{|\gg|} \dogo{g_i}.
\end{equation}
In particular
\begin{equation}
\label{quasiisometry2}
\sum_{i = 1}^{|\gg|} (\dogo{g_i} - C) \vee \epsilon \leq \dogo{\pi(\gg)} \leq \sum_{i = 1}^{|\gg|} \dogo{g_i}
\end{equation}
and thus
\begin{equation}
\label{quasiisometry3}
\dogo{\pi(\gg)} \asymp_\times \sum_{i = 1}^{|\gg|} 1\vee\dogo{g_i}
\end{equation}
\end{lemma}
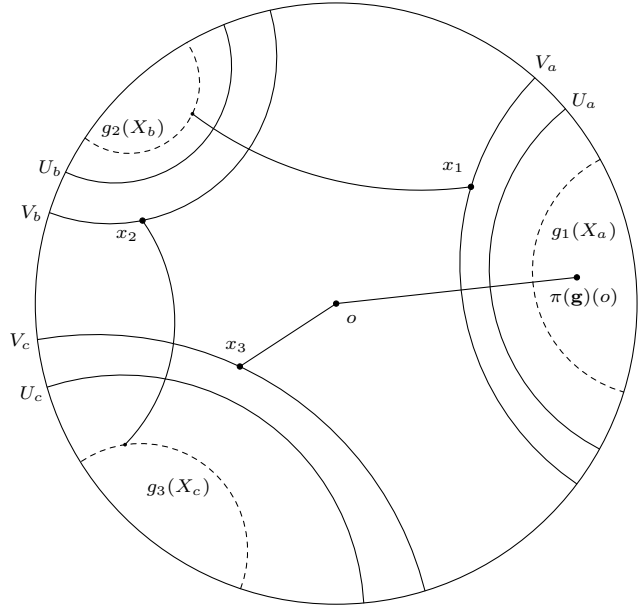
\begin{figure}
\begin{center}
\begin{tikzpicture}[line cap=round,line join=round,>=triangle 45,x=1.0cm,y=1.0cm]
\clip(-4.895,-4.02) rectangle (5.00,4.05);
\draw(0.0,0.0) circle (4.0cm);
\draw [shift={(4.283099034513885,0.4329696912868889)},dash pattern=on 2pt off 2pt] plot[domain=2.0507039028046368:4.4339728284114255,variable=\t]({1.0*1.6748784256043165*cos(\t r)+-0.0*1.6748784256043165*sin(\t r)},{0.0*1.6748784256043165*cos(\t r)+1.0*1.6748784256043165*sin(\t r)});
\draw [shift={(4.748353951077781,0.47506938161209916)}] plot[domain=2.247889155975885:4.234731034706145,variable=\t]({1.0*2.7142499830912525*cos(\t r)+-0.0*2.7142499830912525*sin(\t r)},{0.0*2.7142499830912525*cos(\t r)+1.0*2.7142499830912525*sin(\t r)});
\draw [shift={(5.192066955590817,0.5369298934377099)}] plot[domain=2.3728917598059724:4.116387978350911,variable=\t]({1.0*3.5487868285549222*cos(\t r)+-0.0*3.5487868285549222*sin(\t r)},{0.0*3.5487868285549222*cos(\t r)+1.0*3.5487868285549222*sin(\t r)});
\draw [shift={(-3.012083040210043,3.2835854841554135)}] plot[domain=-1.9392257232190016:0.28223225373733796,variable=\t]({1.0*2.222612337905064*cos(\t r)+-0.0*2.222612337905064*sin(\t r)},{0.0*2.222612337905064*cos(\t r)+1.0*2.222612337905064*sin(\t r)});
\draw [shift={(-2.9425720821605914,3.156105058241987)}] plot[domain=-2.00634379344784:0.36555006299211396,variable=\t]({1.0*1.5526360140016449*cos(\t r)+-0.0*1.5526360140016449*sin(\t r)},{0.0*1.5526360140016449*cos(\t r)+1.0*1.5526360140016449*sin(\t r)});
\draw [shift={(-2.754185174565803,2.935836630029507)},dash pattern=on 2pt off 2pt] plot[domain=-2.2439490561832196:0.6093252827785631,variable=\t]({1.0*0.9385331439412178*cos(\t r)+-0.0*0.9385331439412178*sin(\t r)},{0.0*0.9385331439412178*cos(\t r)+1.0*0.9385331439412178*sin(\t r)});
\draw [shift={(-3.194110354106763,-4.9225766162334885)}] plot[domain=0.24650964162935474:1.7439226966111823,variable=\t]({1.0*4.511504567019594*cos(\t r)+-0.0*4.511504567019594*sin(\t r)},{0.0*4.511504567019594*cos(\t r)+1.0*4.511504567019594*sin(\t r)});
\draw [shift={(-2.8366865058196393,-4.158283214595901)}] plot[domain=0.05462014450100066:1.8896444693721808,variable=\t]({1.0*3.2086081775267683*cos(\t r)+-0.0*3.2086081775267683*sin(\t r)},{0.0*3.2086081775267683*cos(\t r)+1.0*3.2086081775267683*sin(\t r)});
\draw [shift={(-2.6011377721289524,-3.285711019315482)},dash pattern=on 2pt off 2pt] plot[domain=-0.3641423721372501:2.1664758878545376,variable=\t]({1.0*1.42391741251383*cos(\t r)+-0.0*1.42391741251383*sin(\t r)},{0.0*1.42391741251383*cos(\t r)+1.0*1.42391741251383*sin(\t r)});
\draw (-1.2804269894594213,-0.8370546865272379)-- (0.0,-0.0);
\draw (0.0,-0.0)-- (3.1999910454812306,0.348915424894208);
\draw [shift={(1.11888072943081,6.533679518855568)}] plot[domain=4.06530989835595:4.846727654763512,variable=\t]({1.0*5.025532715624298*cos(\t r)+-0.0*5.025532715624298*sin(\t r)},{0.0*5.025532715624298*cos(\t r)+1.0*5.025532715624298*sin(\t r)});
\draw [shift={(-4.486159245286336,-0.24644631226103308)}] plot[domain=-0.7705046356929692:0.6152247592940139,variable=\t]({1.0*2.3406407569929826*cos(\t r)+-0.0*2.3406407569929826*sin(\t r)},{0.0*2.3406407569929826*cos(\t r)+1.0*2.3406407569929826*sin(\t r)});
\begin{scriptsize}
\draw [fill=black] (0.0,-0.0) circle (1pt);
\draw[color=black] (0.2,-0.2) node {$\zero$};
\draw[color=black] (3.2956551912579393,0.9655415669897534) node {$g_1(X_a)$};

\draw[color=black] (3.3,2.7) node {$U_a$};
\draw[color=black] (2.8,3.2) node {$V_a$};

\draw[color=black] (-2.6964002891957583,2.33515619473106552) node {$g_2(X_b)$};

\draw[color=black] (-3.8,1.8) node {$U_b$};
\draw[color=black] (-4.05,1.2) node {$V_b$};

\draw [fill=black] (3.1999910454812306,0.348915424894208) circle (1pt);
\draw[color=black] (3.2935235250000745,0.0801927097768155) node {$\pi(\gg)(\zero)$};
\draw [fill=black] (1.7919753189387375,1.5534262372055987) circle (1pt);
\draw[color=black] (1.5580967616703456,1.8177511393644147) node {$x_1$};
\draw [fill=black] (-1.9108045456133846,2.5240677907713494) circle (.5pt);
\draw [fill=black] (-2.574689601306593,1.104435974779788) circle (1pt);
\draw[color=black] (-2.780470146653977,0.9077896623446498) node {$x_2$};
\draw [fill=black] (-2.8066107407049277,-1.8766965918362397) circle (.5pt);
\draw[color=black] (-2.094694956189467,-2.450970721037196) node {$g_3(X_c)$};

\draw[color=black] (-4.2,-0.5) node {$V_c$};
\draw[color=black] (-4.05,-1.2) node {$U_c$};

\draw [fill=black] (-1.2804269894594213,-0.8370546865272379) circle (1pt);
\draw[color=black] (-1.3295395740340399,-0.57201686259784564) node {$x_3$};
\end{scriptsize}
\end{tikzpicture}
\caption[The strong separation lemma for Schottky products]{The geodesic segment $\geo\zero{\pi(\gg)(\zero})$ splits up naturally into four subsegments, which can then be rearranged by the isometry group to form geodesic segments which connect $\zero$ with $g_1(X_a)$, $V_a$ with $g_2(X_b)$, $V_b$ with $g_3(X_c)$, and $V_c$ with $\zero$, respectively. Here $\gg = (a,g_1)(b,g_2)(c,g_3)$.}
\label{figureschottkylemma}
\end{center}
\end{figure}
\begin{proof}
The second inequality of \eqref{quasiisometry1} is immediate from the triangle inequality. For the first inequality, fix $\gg\in\Gamma$, $x\in X\butnot V_{a_1}$, and $y\in X_{a_{|\gg|}}$. Write $n = |\gg|$. We have
\begin{align*}
\pi(\gg)(y)
\in g_1\cdots g_n(X_{a_n})
&\subset g_1\cdots g_{n - 1}(V_{a_n}) 
\subset g_1\cdots g_{n - 1}(X_{a_{n - 1}})\\ 
&\subset \cdots \subset g_1(V_{a_2})
\subset g_1(X_{a_1}) \subset V_{a_1} \not\ni x.
\end{align*}
Consequently, the geodesic $\geo{x}{\pi(\gg)(y)}$ intersects the sets
\[
\del V_{a_1},\; g_1(\del X_{a_1}),\; g_1(\del V_{a_2}), \ldots,\; g_1\cdots g_{n - 1}(\del V_{a_n}),\; g_1\cdots g_n(\del X_{a_n})
\]
in their respective orders. Thus
\begin{equation}
\label{d0pig}
\begin{split}
\dist(x,\pi(\gg)(y))
&\geq \sum_{i = 1}^n \dist(g_1\cdots g_{i - 1}(\del V_{a_i}),g_1\cdots g_i(\del X_{a_i}))\\
&= \sum_{i = 1}^n \dist(\del V_{a_i},g_i(\del X_{a_i}))\\
&\geq \sum_{i = 1}^n \dist(X\butnot V_{a_i}, g_i(X_{a_i})).
\end{split}
\end{equation}
(Cf. Figure \ref{figureschottkylemma}.) Now fix $i = 1,\ldots,n$, and we will estimate the distance $\dist(X\butnot V_{a_i}, g_i(X_{a_i}))$. For convenience of notation write $a = a_i$ and $g = g_i$.

Fix $z\in X\butnot V_a$ and $w\in g(X_a)$. Combining \eqref{UaVa} and \eqref{XaUa} gives
\[
\wbar\Dist(z,w) \geq \epsilon/2
\]
and in particular
\begin{equation}
\label{XVagVb}
\dist(z,w) \geq \epsilon/2.
\end{equation}
On the other hand, converting the inequality $\Dist(z,w) \geq \epsilon/2$ into a statement about Gromov products shows that
\[
\dist(z,w) \asymp_{\plus,\epsilon} \dist(\zero,z) + \dist(\zero,w) \geq \dist(\zero,w) \geq \dist(g^{-1}(\zero),X_a).
\]
Since $\wbar\Dist(g^{-1}(\zero),X_a) \geq \wbar\Dist(g^{-1}(\bord X\butnot V_a),X_a) \geq \epsilon/2$, we have
\[
\dist(g^{-1}(\zero),X_a) \gtrsim_{\plus,\epsilon} \dist(g^{-1}(\zero),\zero) = \dogo g.
\]
Combining with \eqref{XVagVb} gives
\[
\dist(z,w) \geq (\dogo g - C) \vee (\epsilon/2)
\]
for some $C > 0$ depending only on $\epsilon$. Taking the infimum over all $z,w$ gives
\[
\dist(X\butnot V_{a_i}, g_i(X_{a_i})) \geq (\dogo{g_i} - C) \vee (\epsilon/2).
\]
Summing over all $i = 1,\ldots,n$ and combining with \eqref{d0pig} yields \eqref{quasiisometry1}. Since $\zero\in X_{a_{|\gg|}}$ and $\zero\in X\butnot V_{a_1}$, \eqref{quasiisometry2} follows immediately. Finally, the coarse asymptotic
\[
(\dogo{g_i} - C)\vee \epsilon \asymp_{\times,C,\epsilon} 1\vee \dogo{g_i}.
\]
implies \eqref{quasiisometry3}.
\end{proof}

\begin{corollary}
\label{corollaryschottkySD}
Suppose that $\#\{a\in E : \dist(\zero,U_a) \leq \rho\} < \infty$ for all $\rho > 0$. If the groups $(G_a)_{a\in E}$ are strongly discrete, then $G$ is strongly discrete.

In fact, this corollary holds even if $G$ is only weakly separated and not strongly separated.
\end{corollary}
\begin{proof}
Since $\dogo g \geq \dist(\zero,U_a)$ for all $a\in E$ and $g\in G_a$, our hypotheses implies that
\[
\#\{(a,g)\in \Gamma_E : \dogo g \leq \rho\} < \infty \all \rho.
\]
It follows that for all $N\in\N$,
\begin{equation}
\label{Gamman}
\begin{split}
\#&\left\{\gg\in\Gamma : \sum_{i = 1}^{|\gg|} 1\vee\dogo{g_i} \leq N\right\}\\
&\leq \sum_{n = 0}^N \#\left\{\gg\in (\Gamma_E)^n : \dogo{g_i} \leq N \all i = 1,\ldots,n\right\}\\
&\leq \sum_{n = 0}^N \#\left\{(a,g)\in\Gamma_E : \dogo g \leq N\right\}^n
< \infty.
\end{split}
\end{equation}
Applying \eqref{quasiisometry3} completes the proof. If $G$ is only weakly separated, then for all $\rho > 0$ the Schottky product $\lb G_a\rb_{\dist(\zero,U_a)\leq \rho}$ is still stronglly separated, which is enough to apply \eqref{quasiisometry3} in this context.
\end{proof}

\begin{proposition}
\label{propositionSproductexponent}
~
\begin{itemize}
\item[(i)] If $\#(E) < \infty$ and the groups $G_a$ satisfy $\delta_{G_a} < \infty$, then $\delta_G < \infty$.
\item[(ii)] Suppose that for some $a\in E$, $G_a$ is of divergence type. Then $\delta_G > \delta_{G_a}$.
\item[(iii)] Suppose that $G$ is a group. If $\delta_{G_a} = \infty$ for some $a$, and if $G_b$ is infinite for some $b\neq a$, then $\w\delta_G = \infty$.
\item[(iv)] If $E = \{a,b\}$ and $G_b = g^\Z$, then
\[
\lim_{n\to\infty} \delta(G_a\ast g^{n\Z}) = \delta(G_a).
\]
Moreover, if $G_a$ is of convergence type, then for all $n$ sufficiently large, $G_a\ast g^{n\Z}$ is of convergence type.
\end{itemize}
Moreover, (ii) holds for any free product $G = \lb G_a\rb_{a\in E}$, even if the product is not Schottky.
\end{proposition}


\begin{remark}
Property (iii) tells us that an analogue of property (i) cannot hold for the modified Poincar\'e exponent: if we take the Schottky product of two groups $G_1,G_2$ with $\w\delta(G_i) < \infty$ but $\delta(G_i) = \infty$, then the product $G$ will have $\w\delta(G) = \infty$.
\end{remark}

\begin{proof}[Proof of \text{(i)}]
\eqref{Gamman} shows that for some $C > 0$,
\[
\#\{g\in G : \dogo g \leq \rho\} \leq \#\left\{(a,g)\in\Gamma_E : \dogo g \leq C\rho \right\}^{C\rho} \all \rho > 0.
\]
Applying \eqref{poincarealternate} completes the proof.
\end{proof}
\begin{proof}[Proof of \text{(ii)}]
For all $s\geq 0$,
\begin{align*}
\Sigma_s(G) = \sum_{\gg\in\Gamma} b^{-s\dogo{\pi(\gg)}}
&\geq \sum_{\gg\in\Gamma} b^{-s\sum_1^{|\gg|} \dogo{g_i}}\\
&= \sum_{\gg\in\Gamma} \prod_{i = 1}^{|\gg|} b^{-s\dogo{g_i}}\\
&= \sum_{n = 0}^\infty \sum_{a_1\neq \cdots \neq a_n} \sum_{g_1\in G_{a_1}\butnot\{\id\}}\cdots \sum_{g_n\in G_{a_n}\butnot\{\id\}}\prod_{i = 1}^n b^{-s\dogo{g_i}}\\
&= \sum_{n = 0}^\infty \sum_{a_1\neq \cdots \neq a_n} \prod_{i = 1}^n \sum_{g\in G_{a_i}\butnot\{\id\}} b^{-s\dogo g}\\
&= \sum_{n = 0}^\infty \sum_{a_1\neq \cdots \neq a_n} \prod_{i = 1}^n (\Sigma_s(G_{a_i}) - 1).
\end{align*}
To simplify further calculations, we will assume that $\#(E) = 2$; specifically we will let $E = \{1,2\}$. Then
\begin{align*}
\Sigma_s(G)
&\geq_\pt \sum_{n = 0}^\infty \begin{cases}
2\left(\displaystyle\prod_{a\in E}\big(\Sigma_s(G_a) - 1\big)\right)^{n/2} & n \text{ even}\\
\left(\displaystyle\prod_{a\in E}\big(\Sigma_s(G_a) - 1\big)\right)^{(n - 1)/2}\left(\displaystyle\sum_{a\in E}\big(\Sigma_s(G_a) - 1\big)\right) & n \text{ odd}
\end{cases}\\
&\asymp_\times \sum_{n = 0}^\infty \left(\prod_{a\in E}\big(\Sigma_s(G_a) - 1\big)\right)^{n/2}.
\end{align*}
This series diverges if and only if
\begin{equation}
\label{prodSsGa}
\prod_{a\in E}(\Sigma_s(G_a) - 1) \geq 1.
\end{equation}
Now suppose that $G_1$ is of divergence type, and let $\delta_1 = \delta(G_1)$. By the monotone convergence theorem,
\[
\lim_{s\searrow \delta_1} \prod_{a\in E}(\Sigma_s(G_a) - 1) = \prod_{a\in E}(\Sigma_{\delta_1}(G_a) - 1) = \infty(\Sigma_{\delta_1}(G_2) - 1) = \infty.
\]
(The last equality holds since $G_2$ is nontrivial, see Definition \ref{definitionschottkyproduct}.) So for $s$ sufficiently close to $\delta_1$, \eqref{prodSsGa} holds, and thus $\Sigma_s(G) = \infty$.
\end{proof}
\begin{proof}[Proof of \text{(iii)}]
Fix $\rho > 0$, and let $h\in G_b$ satisfy $\dist(h(\zero),U_a) \geq \rho$. (This is possible since $G_b$ is non-elliptic and $\dist(h(\zero),U_a) \asymp_\plus \dogo h \all h\in G_b$.) Then the set
\[
S_\rho = \{gh(\zero) : g\in G_a\}
\]
is $\rho$-separated, but $\delta(S_\rho) = \delta(G_a) = \infty$. Since $\rho$ was arbitrary, it follows from Proposition \ref{propositionbasicmodified}(iv) that $\w\delta(G) = \infty$.
\end{proof}
\begin{proof}[Proof of \text{(iv)}]
We will in fact show the following more general result:
\begin{proposition}
Suppose $E = \{a,b\}$, and fix $s\notin \Delta(G_a)\cup\Delta(G_b)$. Then there exists a finite set $F\subset G_b$ such that for all $H\leq G_b$, if $H\cap F = \{\id\}$, then $s\notin\Delta(G_a\ast H)$.
\end{proposition}
\noindent Indeed, for such an $s$, the Poincar\'e series $\Sigma_s(G_a\ast H)$ can be estimated using \eqref{quasiisometry2} as follows:
\[
\Sigma_s(G_a\ast H) = \sum_{\gg\in\Gamma} b^{-s\dogo{\pi(\gg)}}
\leq \sum_{\gg\in\Gamma} b^{-s\sum_1^{|\gg|} (\dogo{g_i} - C)} = \sum_{\gg\in\Gamma} b^{sC|\gg|} b^{-s\sum_1^{|\gg|} \dogo{g_i}}.
\]
Continuing as in the proof of part (ii), we get
\[
\Sigma_s(G_a\ast H) \lesssim_\times \sum_{n = 0}^\infty b^{sCn} \Big(\big(\Sigma_s(G_a) - 1\big)\big(\Sigma_s(H) - 1\big)\Big)^{n/2}.
\]
Since $\Sigma_s(H) - 1 \leq \Sigma_s(G_b\butnot F)$, to show that $\Sigma_s(G_a\ast H) < \infty$ it suffices to show that
\begin{equation}
\label{ETSSproductexponent4}
b^{sC} \Big(\big(\Sigma_s(G_a) - 1\big)\big(\Sigma_s(G_b\butnot F)\big)\Big)^{1/2} < 1.
\end{equation}
But since the series $\Sigma_s(G_a)$ and $\Sigma_s(G_b)$ both converge by assumption, \eqref{ETSSproductexponent4} holds for all $F\subset G_b$ sufficiently large.
\end{proof}

We will sometimes find the following variant of Proposition \ref{propositionSproductexponent}(ii) more useful than the original:

\begin{proposition}[Cf. {\cite[Proposition 2]{DOP}}]
\label{propositiondeltagtrdeltap}
Fix $H\leq G\leq\Isom(X)$, and suppose that
\begin{itemize}
\item[(I)] $\Lambda_H \propersubset \Lambda_G$,
\item[(II)] $G$ is of general type, and
\item[(III)] $H$ is of compact type and of divergence type.
\end{itemize}
Then $\delta_G > \delta_H$.
\end{proposition}
\begin{proof}
Fix $\xi\in\Lambda_G\butnot\Lambda_H$, and fix $\epsilon > 0$ small enough so that $B(\xi,\epsilon)\cap\Lambda_H = \emptyset$. Since $G$ is of general type, by Proposition \ref{propositionminimal3}, there exists a loxodromic isometry $g\in G$ such that $g_+,g_-\in B(\xi,\epsilon/4)$. After replacing $g$ by an appropriate power, we may assume that $g^n(\bord X\butnot B(\xi,\epsilon/2)) \subset B(\xi,\epsilon/2)$ for all $n\in\Z\butnot\{0\}$. Now let
\begin{align*}
U_1 &= B(\xi,\epsilon/2)\\
U_2 &= \thick{\Lambda_H}{\epsilon/4}.
\end{align*}
Since $H$ is of compact type (and strongly discrete by Observation \ref{observationnotstronglydiscrete}), Proposition \ref{propositionSDonclX} shows that there exists a finite set $F\subset H$ such that for all $h\in H\butnot F$, $h(\bord X\butnot U_2)\subset U_2$. Let
\[
S = \coprod_{n\geq 0} \big((H\butnot F)\times(g^\Z\butnot\{\id\})\big)^n,
\]
and define $\pi:S\to G$ via the formula
\[
\pi\big((h_i,j_i)_{i = 1}^n\big) = h_1 j_1 \cdots h_n j_n.
\]
A variant of the Ping-Pong Lemma shows that $\pi$ is injective. On the other hand, for all $(h_i,j_i)_{i = 1}^n\in S$, the triangle inequality gies
\[
\dist\left(\zero,\pi\big((h_i,j_i)_{i = 1}^n\big)(\zero)\right) \leq \sum_{i = 1}^n [\dogo{h_i} + \dogo{j_i}].
\]
Thus for all $s > \delta_H$,
\begin{align*}
\Sigma_s(G) &\geq \sum_{g\in \pi(S)} e^{-s\dogo g}\\
&\geq \sum_{(h_i,j_i)_{i = 1}^n \in S} \prod_{i = 1}^n e^{-s[\dogo{h_i} + \dogo{j_i}]}\\
&= \sum_{n\geq 0} \left(\sum_{h\in H\butnot F}\sum_{j\in g^\Z\butnot\{\id\}} e^{-s[\dogo h + \dogo j]}\right)^n\\
&= \sum_{n\geq 0} \left(\Sigma_s(H\butnot F) \Sigma_s(g^\Z\butnot\{\id\})\right)^n\\
&\begin{cases}
= \infty & \Sigma_s(H\butnot F) \Sigma_s(g^\Z\butnot\{\id\}) \geq 1\\
< \infty & \Sigma_s(H\butnot F) \Sigma_s(g^\Z\butnot\{\id\}) < 1
\end{cases}.
\end{align*}
Now since $H$ is of divergence type, by the monotone convergence theorem,
\begin{align*}
\lim_{s\searrow\delta_H} \Sigma_s(H\butnot F) \Sigma_s(g^\Z\butnot\{\id\}) &= \Sigma_{\delta_H}(H\butnot F) \Sigma_{\delta_H}(g^\Z\butnot\{\id\})\\ 
&= \infty \cdot (\text{positive constant})\\
&= \infty.
\end{align*}
Thus, for $s$ sufficiently close to $\delta_H$, $\Sigma_s(G) = \infty$. This shows that $\delta_G > \delta_H$.
\end{proof}
\begin{remark}
The reason that we couldn't deduce Proposition \ref{propositiondeltagtrdeltap} directly from Proposition \ref{propositionSproductexponent}(ii) is that the group $\lb H,g^\Z\rb$ considered in the proof of Proposition \ref{propositiondeltagtrdeltap} is not necessarily a free product due to the existence of the finite set $F$. In the Standard Case, this could be solved by taking a finite-index subgroup of $H$ whose intersection with $F$ is trivial, but in general, it is not clear that such a subgroup exists.
\end{remark}

\section{A partition-structure--like structure}
\label{subsectionWg}
For each $\gg \in\Gamma$, let
\[
W_\gg = \pi(\gg)(X_{a_{|\gg|}}),
\]
unless $\gg = \emptyset$, in which case let $W_\smallemptyset = \bord X$.

\begin{standingassumption}
\label{assumptionGa}
In what follows, we assume that for each $a\in E$, either
\begin{itemize}
\item[(1)] $G_a$ is a group, or
\item[(2)] $G_a\equiv\N$.
\end{itemize}
\end{standingassumption}

For $\gg,\hh\in\Gamma$, write $\gg \leq \hh$ if $\hh = \gg\ast\kk$ for some $\kk\in\Gamma$.

\begin{lemma}
\label{lemmaincomparable}
Fix $\gg,\hh\in\Gamma$. If $\gg\leq \hh$, then $W_\hh\subset W_\gg$. On the other hand, if $\gg$ and $\hh$ are incomparable ($\gg\not\leq\hh$ and $\hh\not\leq\gg$), then $W_\gg\cap W_\hh = \emptyset$.
\end{lemma}
\begin{proof}
The first assertion follows from Lemma \ref{lemmapingpong}. For the second assertion, it suffices to show that if $(a,g),(b,h)\in \Gamma_E$ are distinct, then $W_{(a,g)}\cap W_{(b,h)} = \emptyset$. Since $W_{(a,g)} \subset \schottkycl{U_a}$ and $(\schottkycl{U_a})_{a\in E}$ are disjoint, if $a\neq b$ then $W_{(a,g)}\cap W_{(b,h)} = \emptyset$. So suppose $a = b$. Assumption \ref{assumptionGa} guarantees that either $g^{-1}h\in G_a$ or $h^{-1}g\in G_a$; without loss of generality assume that $g^{-1}h\in G_a$. Then
\[
W_{(a,g)}\cap W_{(b,h)} = g(X_a)\cap h(X_a) = g(X_a\cap g^{-1}h(X_a)) \subset g(X_a\cap \schottkycl{U_a}) = \emptyset.
\]
\end{proof}

\begin{lemma}
\label{lemmastrongseparation2}
There exists $\sigma > 0$ such that for all $\gg\in \Gamma$,
\begin{equation}
\label{WgShad}
W_\gg \subset \Shad(\pi(\gg)(\zero),\sigma).
\end{equation}
In particular
\begin{equation}
\label{WgDiam}
\Diam(W_\gg) \lesssim_\times b^{-\dogo{\pi(\gg)}}.
\end{equation}
\end{lemma}
\begin{proof}
Let $n = |\gg|$, $g = g_n$, $a = a_n$, and $z = \pi(\gg)^{-1}(\zero)$. Observe that if $g(z)\in V_a$, then Lemma \ref{lemmapingpong} implies that $\zero\in V_{a_1}$, a contradiction. Thus $z\in g^{-1}(X\butnot V_a)$. If $X$ is not global, then $\eqref{strongseparation}_{U_a = V_a, \epsilon = \epsilon/2}$ implies that $\Dist(z,X_a) \geq \epsilon/2$. On the other hand, if $X$ is global then we have $z\in U_a$, so \eqref{UaVa} implies that $\Dist(z,X_a) \geq \epsilon/2$. Either way, we have $\Dist(z,X_a) \geq \epsilon/2$.

Let $\sigma > 0$ be large enough so that the Big Shadows Lemma \ref{lemmabigshadow} holds; then we have $X_a \subset \Shad_z(\zero,\sigma)$. Applying $\pi(\gg)$ yields \eqref{WgShad}, and combining with the Diameter of Shadows Lemma \ref{lemmadiameterasymptotic} yields \eqref{WgDiam}.
\end{proof}

Let $\del\Gamma$ denote the set of all infinite words with letters in $\Gamma_E$ such that $a_i \neq a_{i + 1}$ for all $i$. Given $\gg\in\del\Gamma$, for each $n$, $\gg\given n\in\Gamma$. Then Lemmas \ref{lemmaincomparable} and \ref{lemmastrongseparation2} show that the sequence $(W_{\gg\given n})_0^\infty$ is an infinite descending sequence of closed sets with diameters tending to zero; thus there exists a unique point $\xi\in \bigcap_0^\infty W_{\gg\given n}$, which will be denoted $\pi(\gg)$.

\begin{lemma}
\label{lemmaschottkyradial}
For all $\gg\in\del\Gamma$, $\pi(\gg\given n)(\zero)\to \pi(\gg)$ radially. In particular $\pi(\del\Gamma)\subset\Lr(G)$.
\end{lemma}
\begin{proof}
This is immediate from \eqref{WgShad}, since by definition $\pi(\gg)\in W_{\gg\given n}$ for all $n$.
\end{proof}

\begin{lemma}[Cf. Klein's combination theorem {\cite[Theorem 1.1]{LOW}}, \cite{Klein_combination_theorem}]
\label{lemmaschottkydomain}
The set
\begin{equation}
\label{schottkydomain}
\DD = \bord X \butnot \bigcup_{a\in E}\bigcup_{g\in G_a} g(X_a) = \bord X \butnot \bigcup_{(a,g)\in \Gamma_E} W_{(a,g)}
\end{equation}
satisfies $G(\DD) = \bord X\butnot\pi(\del\Gamma)$.
\end{lemma}
\noindent Before we begin the proof of this lemma, note that since $\DD\cap X$ is open (Lemma \ref{lemmaDopen} below), the connectedness of $X$ implies that $g(\DD)\cap \DD\neq\emptyset$ for some $g\in G$. Thus $\DD$ is not a fundamental domain.
\begin{proof}
Fix $x\in\bord X\butnot \pi(\del\Gamma)$, and consider the set
\[
\Gamma_x := \left\{\gg\in \Gamma : x\in W_\gg\right\}.
\]
By Lemma \ref{lemmaincomparable}, $\Gamma_x$ is totally ordered as a subset of $\Gamma$. If $\Gamma_x$ is infinite, let $\gg\in\del\Gamma$ be the unique word such that $\Gamma_x = \{\gg\given n : n\in\Neur\}$; Lemma \ref{lemmastrongseparation2} implies that $x = \pi(\gg)\in \pi(\del\Gamma)$, contradicting our hypothesis. Thus $\Gamma_x$ is finite. If $\Gamma_x = \emptyset$, we are done. Otherwise, let $\gg$ be the largest element of $\Gamma_x$. Then $x\in W_\gg$, so $\pi(\gg)^{-1}(x)\in X_a$, where $a = a_{|\gg|}$. The maximality of $\gg$ implies that
\[
\pi(\gg)^{-1}(x) \notin W_{(b,h)} = h(X_b) \all b\in E\butnot \{a\} \all h\in G_b\butnot\{\id\},
\]
but on the other hand $\pi(\gg)^{-1}(x)\in X_a\subset \bord X\butnot U_a$ implies that $\pi(\gg)^{-1}(x)\notin W_{(a,h)}$ for all $h\in G_a\butnot\{\id\}$. Thus $\pi(\gg)^{-1}(x)\in\DD$.
\end{proof}

\begin{lemma}
\label{lemmaDopen}
Suppose that for each $a\in E$, $G_a$ is strongly discrete. Then
\[
\DD\butnot\Int(\DD)\subset \bigcup_{a\in E} \Lambda_a,
\]
where $\Lambda_a = \Lambda(G_a)$. In particular, $\DD\cap X$ is open.
\end{lemma}
\begin{proof}
Fix $x\in\DD\butnot\Int(\DD)$, and find a sequence $(a_n,g_n)\in \Gamma_E$ such that $\wbar\Dist(x,g_n(X_{a_n})) \to 0$. Since $g_n(X_{a_n}) \subset \schottkycl{U_{a_n}}$, \eqref{strongseparation} implies that $a_n$ is constant for all sufficiently large $n$, say $a_n = a$. On the other hand, if there is some $g\in G_a$ such that $g_n = g$ for infinitely many $n$, then since $g(X_a)$ is closed we would have $x\in g(X_a)$, contradicting that $x\in\DD$. Since $G_a$ is strongly discrete, it follows that $\dogo{g_n}\to\infty$, and thus $\Diam(g_n(X_a))\to 0$ by Lemma \ref{lemmastrongseparation2}. Since $g_n(\zero)\in g_n(X_a)$, it follows that $g_n(\zero)\to x$, and thus $x\in \Lambda_a$.
\end{proof}

\begin{theorem}
\label{theoremschottkylimitset}
\[
\Lambda = \pi(\del\Gamma)\cup\bigcup_{g\in G}\bigcup_{a\in E} g(\Lambda_a).
\]
\end{theorem}
\begin{proof}
The $\supset$ direction follows from Lemma \ref{lemmaschottkyradial}, so let us show $\subset$. It suffices to show that $\Lambda\cap \DD\subset \bigcup_{a\in E} \Lambda_a$. Indeed, for all $\gg\in\Gamma\butnot\{\emptyset\}$, Lemma \ref{lemmapingpong} gives $\pi(\gg)(\zero)\in g_1(X_{a_1}) \subset \bord X\butnot\DD$. Thus $\Lambda\cap\DD = \DD\cap\del\DD \subset \bigcup_{a\in E}\Lambda_a$ by Lemma \ref{lemmaDopen}.
\end{proof}

\begin{corollary}
\label{corollaryschottkycompacttype}
If $E$ is finite and each $G_a$ is strongly discrete and of compact type, then $G$ is strongly discrete and of compact type.
\end{corollary}
\begin{proof}
Strong discreteness follows from Corollary \ref{corollaryschottkySD}, so let us show that $G$ is of compact type. Let $(\xi_n)_1^\infty$ be a sequence in $\Lambda$. For each $n\in\N$, if $\xi_n\in \pi(\del\Gamma)$, write $\xi_n = \pi(\gg_n)$ for some $\gg_n\in\del\Gamma$; otherwise, write $\xi_n = \pi(\gg_n)(\eta_n)$ for some $\gg_n\in\Gamma$ and $\eta_n\in\Lambda_{a_n}$. Either way, note that $\xi_n\in W_\hh$ for all $\hh\leq\gg$.

For each $\hh\in\Gamma$, let
\[
S_\hh = \{n\in\N : \hh\leq\gg_n\}.
\]
Since $\Gamma$ is countable, by extracting a subsequence we may without loss of generality assume that for all $\hh\in\Gamma$, either $n\in S_\hh$ for all but finitely many $n$, or $n\in S_\hh$ for only finitely many $n$. Let
\[
\Gamma' = \{\hh\in\Gamma : \text{$n\in S_\hh$ for all but finitely many $n$}\}.
\]
Then by Lemma \ref{lemmaincomparable}, the set $\Gamma'$ is totally ordered. Moreover, $\emptyset\in \Gamma'$. If $\Gamma'$ is infinite, then choose $\gg\in\del\Gamma$ such that $\Gamma' = \{\gg\given m : m\geq 0\}$; by Lemma \ref{lemmastrongseparation2}, we have $\xi_n\to\pi(\gg)$. Otherwise, let $\gg$ be the largest element of $\Gamma'$. For each $n$, either $\xi_n\in W_{\gg(b_n,h_n)}$ for some $(b_n,h_n)\in\Gamma_E$, or $\xi_n = \pi(\gg)(\eta_n)$ for some $a_n\in E$ and $\eta_n\in\Lambda_{a_n}$. By extracting another subsequence, we may assume that either the first case holds for all $n$, or the second case holds for all $n$. Suppose the first case holds for all $n$. The maximality of $\gg$ implies that for each $(b,h)\in\Gamma_E$, there are only finitely many $n$ such that $(b_n,h_n) = (b,h)$. Since $E$ is finite, by extracting a further subsequence we may assume that $b_n = b$ for all $n$. Since $G_b$ is strongly discrete and of compact type, by extracting a further subsequence we may assume that $h_n(\zero)\to \eta$ for some $\eta\in \Lambda_b$. But then $\xi_n \to \pi(\gg)(\eta)\in\Lambda$.

Suppose the second case holds for all $n$. Since $\Lambda_E = \bigcup_{a\in E} \Lambda_a$ is compact, by extracting a further subsequence we may assume that $\eta_n\to \eta$ for some $\eta\in\Lambda_E$. But then $\xi_n \to \pi(\gg)(\eta)\in\Lambda$.
\end{proof}

\begin{corollary}
\label{corollaryuncountable}
If $\#(\Gamma_E) \geq 3$, then $\#(\Lambda)\geq \#(\R)$.
\end{corollary}
\begin{proof}
The hypothesis $\#(\Gamma_E)\geq 3$ implies that for each $\gg\in\Gamma$,
\[
\#\{(a,g)(b,h)\in\Gamma_E^2 : \gg(a,g)(b,h)\in \Gamma\} = \sum_{b\neq a\neq a_{|\gg|}} (\#(G_a) - 1)(\#(G_b) - 1) \geq 2.
\]
Thus, the tree $\Gamma$ has no terminal nodes or infinite isolated paths. It follows that $\del\Gamma$ is perfect and therefore has cardinality at least $\#(\R)$; since $\pi(\del\Gamma)\subset\Lambda$, we have $\#(\Lambda) \geq \#(\del\Gamma)\geq \#(\R)$.
\end{proof}


\begin{proposition}
\label{propositionSproductMDWDPD}
Suppose that the Schottky system $(U_a)_{a\in E}$ is global. Then if $(G_a)_{a\in E}$ are moderately (resp. weakly) discrete, then $G$ is moderately (resp. weakly) discrete. If $(G_a)_{a\in E}$ act properly discontinuously, then $G$ acts properly discontinuously.
\end{proposition}
\begin{proof}
Let $\DD$ be as in \eqref{schottkydomain}. Fix $x\in \DD$ and $\gg \in\Gamma$, let $n = |\gg|$, and suppose that $\pi(\gg)(x)\in \DD$. Then:
\begin{itemize}
\item[(A)] For all $i = 1,\ldots,n$, if $g_{i + 1}\cdots g_n(x)\in X_{a_i}$ then Lemma \ref{lemmapingpong} would give $\pi(\gg)(x)\in g_1(X_{a_1})$, so $g_{i + 1}\cdots g_n(x)\in V_{a_i}$.
\item[(B)] For all $i = 0,\ldots,n - 1$, if $g_{i + 1}\cdots g_n(x)\in X_{a_{i + 1}}$ then Lemma \ref{lemmapingpong} would give $x\in g_n^{-1}(X_{a_n})$, so $g_{i + 1}\cdots g_n(x)\in V_{a_{i + 1}}$.
\end{itemize}
If $n\geq 2$, letting $i = 1$ in (A)-(B) yields a contradiction, so $n = 0$ or $1$. Moreover, if $n = 1$, plugging in $i = 1$ in (A) gives $x\in V_{a_1}$.

To summarize, if we let
\[
G_x = \begin{cases}
G_a & x\in V_a\\
\{\id\} & x\notin \bigcup_{a\in E} V_a
\end{cases}
\]
then
\[
g(x)\in \DD \;\Rightarrow\; g\in G_x \all g\in G.
\]
More concretely,
\[
\dist(x,g(x)) < \dist(x,X\butnot \DD) \;\Rightarrow\; g\in G_x \all g\in G.
\]
Comparing with the definitions of moderate and weak discreteness (Definition \ref{definitiondiscreteness}) and the definition of proper discontinuity (Definition \ref{definitionproperdiscontinuity}) completes the proof.
\end{proof}

\section{Existence of Schottky products}
\label{subsectionschottkyexamples}

\begin{proposition}
\label{propositionschottkyproduct}
Suppose that $\Lambda_{\Isom(X)} = \del X$, and let $G_1,G_2\leq\Isom(X)$ be groups with the following property: For $i = 1,2$, there exist $\xi_i\in\del X$ and $\epsilon > 0$ such that
\begin{equation}
\label{xiidef}
\Dist(\xi_i,g(\xi_i)) \geq \epsilon \all g\in G_i\butnot\{\id\}.
\end{equation}
Then there exists $\phi\in\Isom(X)$ such that the product $\lb G_1,\phi(G_2)\rb$ is a global strongly separated Schottky product.
\end{proposition}
\begin{proof} We begin the proof with the following
\begin{claim}
For each $i = 1,2$, there exists an open set $A_i\ni \xi_i$ such that $g(A_i)\cap A_i = \emptyset$ for all $g\in G_i\butnot\{\id\}$.
\end{claim}
\begin{subproof}
Fix $i = 1,2$. Clearly, \eqref{xiidef} implies that $\xi_i\notin \Lambda(G_i)$. Thus, there exists $\delta > 0$ such that $\Dist(g(\zero),\xi_i)\geq \delta$ for all $g\in G_i$. Applying the Big Shadows Lemma \ref{lemmabigshadow}, there exists $\sigma > 0$ such that $B(\xi_i,\delta/2) \subset \Shad_{g^{-1}(\zero)}(\zero,\sigma)$ for all $g\in G$. But then by the Bounded Distortion Lemma \ref{lemmaboundeddistortion}, we have
\[
\Diam(g(B(\xi_i,\gamma))) \asymp_{\times,\sigma} b^{-\dogo g} \Diam(B(\xi_i,\gamma)) \leq 2\gamma \all g\in G \all 0 < \gamma \leq \delta/2.
\]
Choosing $\gamma$ appropriately gives $\Diam(g(B(\xi_i,\gamma))) < \epsilon/2$ for all $g\in G$. Letting $A_i = B(\xi_i,\gamma)$ completes the proof of the claim.
\end{subproof}
For each $i = 1,2$, let $A_i$ be as above, and fix an open set $B_i\ni\xi_i$ such that $\wbar\Dist(B_i,\bord X\butnot A_i) > 0$. Since $\Lambda_{\Isom(X)} = \del X$, there exists a loxodromic isometry $\phi\in\Isom(X)$ such that $\phi_-\in B_1$ and $\phi_+\in B_2$ (Proposition \ref{propositionminimal3}). Then by Theorem \ref{theoremloxodromicextra}, $\phi^n\to \phi_+$ uniformly on $\bord X\butnot B_1$, so $\phi^n(B_1)\cup B_2 = \bord X$ for all $n\in\N$ sufficiently large. Fix such an $n$, and let $U_1 = \phi^n(\bord X\butnot A_1)$, $U_2 = \bord X\butnot A_2$, $V_1 = \phi^n(\bord X\butnot B_1)$, $V_2 = \bord X\butnot B_2$. Then $(V_1,V_2)$ is a global Schottky system for $(G_1,\phi(G_2))$, which implies that $(U_1,U_2)$ is a global strongly separated Schottky system for $(G_1,\phi(G_2))$. This concludes the proof of the proposition.
\end{proof}

\begin{remark}
\label{remarkschottkyproduct}
The hypotheses of the above theorem are satisfied if $X$ is an algebraic hyperbolic space and for each $i = 1,2$, $G_i$ is strongly discrete and of compact type and $\Lambda_i = \Lambda(G_i) \propersubset\del X$.
\end{remark}
\begin{proof}
We have $\Lambda_{\Isom(X)} = \del X$ by Observation \ref{observationtransitivity}. Fix $i = 1,2$. Since $\Lambda_i \propersubset\del X$, $\del X\butnot\Lambda(G_i)$ is a nonempty open set. For each $g\in G_i\butnot\{\id\}$, the set $\Fix(g)$ is totally geodesic (Theorem \ref{theoremtotallygeodesic}) and therefore nowhere dense; since $G_i$ is countable, it follows that $\bigcup_{g\in G_i\butnot\{\id\}} \Fix(g)$ is a meager set, so the set
\[
S_i = \big(\del X\butnot\Lambda(G_i)\big)\butnot\bigcup_{g\in G_i\butnot\{\id\}} \Fix(g)
\]
is nonempty. Fix $\xi_i\in S_i$. By Proposition \ref{propositionSDonclX}, 
\[
\liminf_{g\in G_i} \Dist(\xi_i,g(\xi_i)) \geq \Dist(\xi_i,\Lambda(G_i)) > 0. 
\]
On the other hand, for all $g\in G_i\butnot\{\id\}$ we have $\xi_i \notin \Fix(g)$ and therefore $\Dist(\xi_i,g(\xi_i)) > 0$. Combining yields $\inf_{g\in G_i\butnot\{\id\}} \Dist(\xi_i,g(\xi_i)) > 0$ and thus \eqref{xiidef}. This concludes the proof of the remark.
\end{proof}

\begin{proposition}
\label{propositionnonelementaryequivalent}
For a semigroup $G\prec\Isom(X)$, the following are equivalent:
\begin{itemize}
\item[(A)] $G$ is either outward focal or of general type.
\item[(B)] $G$ contains a strongly separated Schottky subsemigroup.
\item[(C)] $\w\delta(G) > 0$.
\item[(D)] $\#(\Lambda_G) \geq \#(\R)$.
\item[(E)] $\#(\Lambda_G) \geq 3$, i.e. $G$ is nonelementary.
\end{itemize}
If $G$ is a group, then these are also equivalent to:
\begin{itemize}
\item[(F)] $G$ contains a global strongly separated Schottky subgroup.
\end{itemize}
\end{proposition}
The implications (C) \implies (E) \implies (A) have been proven elsewhere in the paper; see Proposition \ref{propositioncardinalitylimitset} and Theorem \ref{theorembishopjonesmodified}. The implication (B) \implies (D) is an immediate consequence of Corollary \ref{corollaryuncountable}, and (D) \implies (E) and (F) \implies (B) are both trivial. So it remains to prove (A) \implies (B) \implies (C), and that (A) \implies (F) if $G$ is a group.
\begin{proof}[Proof of \text{(A) \implies (B)}]
Suppose first that $G$ is outward focal with global fixed point $\xi$. Then there exists $g\in G$ with $g'(\xi) > 1$, and there exists $h\in G$ such that $h_+ \neq g_+$. If we let $j = g^n h$, then $j'(\xi) > 1$ (after choosing $n$ sufficiently large), and $j_+ \neq g_+$.

So regardless of whether $G$ is outward focal or of general type, there exist loxodromic isometries $g,h\in G$ such that $g_+ \notin \Fix(h)$ and $h_+\notin\Fix(g)$. It follows that there exists $\epsilon > 0$ such that
\[
B(g_+,\epsilon) \cap h^n\big(B(g_+,\epsilon)\big) = B(h_+,\epsilon) \cap g^n\big(B(h_+,\epsilon)\big) = \emptyset \all n\geq 1.
\]
Let $U_1 = B(g_+,\epsilon/2)$, $U_2 = B(h_+,\epsilon/2)$, $V_1 = B(g_+,\epsilon)$, and $V_2 = B(h_+,\epsilon)$. By Theorem \ref{theoremloxodromicextra}, for all sufficiently large $n$ we have $g^n(V_1\cup V_2) \subset U_1$ and $h^n(V_1\cup V_2) \subset U_2$. It follows that $(V_1,V_2)$ is a Schottky system for $((g^n)^\N,(h^n)^\N)$, and that $(U_1,U_2)$ is a strongly separated Schottky system for $((g^n)^\N,(h^n)^\N)$.
\end{proof}
\begin{proof}[Proof of \text{(B) \implies (C)}]
Since a cyclic loxodromic semigroup is of divergence type (an immediate consequence of \eqref{dognoloxodromic}), Proposition \ref{propositionSproductexponent}(i),(ii) shows that $0 < \delta(H) < \infty$, where $H\prec G$ is a Schottky subsemigroup. Thus $\w\delta(H) > 0$, and so $\w\delta(G) > 0$.
\end{proof}
\begin{proof}[Proof of \text{(A) \implies (F)} for groups]
Fix loxodromic isometries $g,h\in G$ with $\Fix(g)\cap \Fix(h) = \emptyset$. Choose $\epsilon > 0$ such that
\[
B(\Fix(g),\epsilon) \cap h^n\big(B(\Fix(g),\epsilon)\big) = B(\Fix(h),\epsilon) \cap g^n\big(B(\Fix(h),\epsilon)\big) = \emptyset \all n\geq 1.
\]
Let $U_1 = B(\Fix(g),\epsilon/2)$, $U_2 = B(\Fix(h),\epsilon/2)$, $V_1 = B(\Fix(g),\epsilon)$, and $V_2 = B(\Fix(h),\epsilon)$. By Theorem \ref{theoremloxodromicextra}, for all sufficiently large $n$ we have $g^n(\bord X\butnot B(g_-,\epsilon/2)) \subset B(g_+,\epsilon/2)$ and $h^n(\bord X\butnot B(h_-,\epsilon/2)) \subset B(h_+,\epsilon/2)$. It follows that $(V_1,V_2)$ is a global Schottky system for $((g^n)^\Z,(h^n)^\Z)$, and that $(U_1,U_2)$ is a global strongly separated Schottky system for $((g^n)^\Z,(h^n)^\Z)$.

\end{proof}

\chapter{Parabolic groups}\label{sectionparabolic}

In this chapter we study parabolic groups. We begin with a list of several examples of parabolic groups acting on $\E^\infty$, the half-space model of infinite-dimensional real hyperbolic geometry. These examples include a counterexample to the infinite-dimensional analogue of Margulis's lemma, as well as a parabolic isometry that generates a cyclic group which is not discrete in the ambient isometry group. The latter example is the Poincar\'e extension of an example due to M. Edelstein. After giving these examples of parabolic groups, we prove a lower bound on the Poincar\'e exponent of a parabolic group in terms of its algebraic structure (Theorem \ref{theoremparaboliclowerbound}). We show that it is optimal by constructing explicit examples of parabolic groups acting on $\E^\infty$ which come arbitrarily close to this bound.

\section{Examples of parabolic groups acting on $\E^\infty$}
\label{subsectionparabolic}
Let $X = \E = \E^\infty$ be the half-space model of infinite-dimensional real hyperbolic geometry (\6\ref{subsubsectionE}). Recall that $\BB = \del \E\butnot\{\infty\}$ is an infinite-dimensional Hilbert space, and that \emph{Poincar\'e extension} is the homomorphism $\what\cdot:\Isom(\BB) \to \Isom(\E)$ given by the formula
\[
\what\cdot \;(g)(t,\xx) = \what g(t,\xx) = (t,g(\xx))
\]
(Observation \ref{observationpoincareextension}). The image of $\what\cdot$ is the set $\{g\in\Stab(\Isom(\E);\infty) : g'(\infty) = 1\}$. Thus, Poincar\'e extension provides a bijection between the class of subgroups of $\Isom(\BB)$ and the class of subgroups of $\Isom(\E)$ for which $\infty$ is a neutral global fixed point. Given a group $G\leq\Isom(\BB)$, we will denote its image under $\what\cdot$ by $\what G$. We may summarize the relation between $\what G$ and $G$ as follows:

\begin{observation}
\label{observationisomE}
~
\begin{itemize}
\item[(i)] $\what G$ is parabolic if and only if $G(\0)$ is unbounded; otherwise $\what G$ is elliptic.
\item[(ii)] $\what G$ is strongly (resp. moderately, weakly, COT) discrete if and only if $G$ is. $\what G$ acts properly discontinuously if and only if $G$ does.
\item[(iii)] Write $\Isom(\BB) = \O(\BB)\ltimes \BB$. Then the preimage of the uniform operator topology under $\what\cdot$ is equal to the product of the uniform operator topology on $\O(\BB)$ with the usual topology on $\BB$. Thus if we denote this topology by UOT*, then $\what G$ is UOT-discrete if and only if $G$ is UOT*-discrete.
\item[(iv)]
For all $g\in G$, we have
\[
e^{\dogo{\what g}} \asymp_\times \cosh\dogo{\what g} = 1 + \frac{\|(1,g(\0)) - (1,\0)\|^2}{2} \asymp_\times 1 \vee \|g(\0)\|^2
\]
and thus for all $s\geq 0$,
\begin{equation}
\label{poincareextensionpoincareexponent}
\Sigma_s(\what G) \asymp_\times \Sigma_s(G) := \sum_{g\in G} (1\vee \|g(\0)\|)^{-2s}.
\end{equation}
\end{itemize}
\end{observation}

In what follows, we let $\delta(G) = \inf\{s : \Sigma_s(G) < \infty\} = \delta(\what G)$, and we say that $G$ is of convergence or divergence type if $\what G$ is.

\subsection{The Haagerup property and the absence of a Margulis lemma}
\label{subsubsectionhaagerup}
One question which has been well studied in the literature is the following: For which abstract groups $\Gamma$ can $\Gamma$ be embedded as a strongly discrete subgroup of $\Isom(\BB)$? Such a group is said to have the \emph{Haagerup property}.\Footnote{The Haagerup property can also be defined for locally compact groups, by replacing ``finite'' with ``precompact'' in the definition of strong discreteness. However, in \thispaper\ we consider only discrete groups.} For a detailed account, see \cite{CCJJV}.

\begin{remark}
\label{remarkhaagerup}
The following groups have the Haagerup property:
\begin{itemize}
\item \cite[pp.73-74]{HarpeValette} Groups which admit a cocompact action on a proper $\R$-tree. In particular this includes $\F_n(\Z)$ for every $n$.
\item \cite{Jolissaint} Amenable groups. This includes solvable and nilpotent groups.
\end{itemize}
A class of examples of groups without the Haagerup property is the class of infinite groups with Kazdan's property (T). For example, if $n\geq 3$ then $\SL_n(\Z)$ does not have the Haagerup property \cite[\64.2]{BHV}.
\end{remark}

The example of (virtually) nilpotent groups will be considered in more detail in \6\ref{subsubsectiontheoremnilpotent}, since it turns out that every parabolic subgroup of $\Isom(\E)$ which has finite Poincar\'e exponent is virtually nilpotent.

Recall that \emph{Margulis's lemma} is the following lemma:

\begin{proposition}[Margulis's lemma, {\cite[p.126]{Harpe}} or {\cite[p.101]{BGS}}]
\label{propositionmargulislemma}
Let $X$ be a Hadamard manifold with curvature bounded away from $-\infty$. Then there exists $\epsilon = \epsilon_X > 0$ with the following property: For every discrete group $G\leq\Isom(X)$ and for every $x\in X$, the group
\[
G_\epsilon(x) := \lb g\in G : \dist(x,g(x)) \leq \epsilon\rb
\]
is virtually nilpotent.
\end{proposition}
For convenience, we will say that \emph{Margulis's lemma holds} on a metric space $X$ if the conclusion of Proposition \ref{propositionmargulislemma} holds, i.e. if there exists $\epsilon > 0$ such that for every strongly discrete group $G\leq\Isom(X)$ and for every $x\in X$, $G_\epsilon(x)$ is virtually nilpotent. It was proven recently by E. Breuillard, B. Green, and T. C. Tao \cite[Corollary 1.15]{BGT} that Margulis's lemma holds on all metric spaces with \emph{bounded packing} in the sense of \cite{BGT}. This result includes Proposition \ref{propositionmargulislemma} as a special case.

By contrast, in infinite dimensions we have the following:
\begin{observation}
\label{observationmargulislemma}
Margulis's lemma does not hold on the space $X = \E = \E^\infty$.
\end{observation}
\begin{proof}
Since $\F_2(\Z)$ has the Haagerup property, there exists a strongly discrete group $G\leq\Isom(\BB)$ isomorphic to $\F_2(\Z)$, say $G = (g_1)^\Z\ast (g_2)^\Z$. Let $\what G$ be the Poincar\'e extension of $G$. Fix $\epsilon > 0$, and let
\[
x = (t,\0)\in\E
\]
for $t > 0$ large to be determined. Then by \eqref{distE},
\[
\dist(x,\what g_i(x)) = \dist\big((t,\0), (t,g_i(\0))\big) \asymp_\times \|g_i(\0)\|/t.
\]
So if $t$ is large enough, then $\dist(x,\what g_i(x)) \leq \epsilon$. It follows that $\what g_1,\what g_2\in \what G_\epsilon(x)$, and so $\what G_\epsilon(x) = \what G \equiv \F_2(\Z)$ is not virtually nilpotent.
\end{proof}
\begin{remark}
In view of the fact that in the finite-dimensional Margulis's lemma, $\epsilon_{\E^d}$ depends on the dimension $d$ and tends to zero as $d\to\infty$ (see e.g. \cite[Proposition 5.2]{Belolipetsky_survey}), we should not be surprised that the lemma fails in infinite dimensions.
\end{remark}
\begin{remark}
\label{remarkmargulislemma}
In some references (e.g. \cite[Theorem 12.6.1]{Ratcliffe}), the conclusion of Margulis's lemma states that $G_\epsilon(x)$ is elementary rather than virtually nilpotent. The above example shows that the two statements should not be confused with each other. We will show (Example \ref{examplemargulis} below) that the alternative formulation of Margulis's lemma which states that $G_\epsilon(x)$ is elementary also fails in infinite dimensions.
\end{remark}

\begin{remark}
Parabolic groups acting on proper CAT(-1) spaces must be amenable \cite[Proposition 1.6]{BurgerMozes}. Therefore the existence of a parabolic subgroup of $\Isom(\H^\infty)$ isomorphic to $\F_2(\Z)$ also distinguishes $\H^\infty$ from the class of proper CAT(-1) spaces.
\end{remark}

\subsection{Edelstein examples}
\label{subsubsectionedelstein}
One of the oldest results in the field of groups acting by isometries on Hilbert space is the following example due to M. Edelstein:

\begin{proposition}[{\cite[Theorem 2.1]{Edelstein}}]
\label{propositionedelstein}
There exist an isometry $g$ belonging to $\Isom(\ell^2(\N;\C))$ and sequences $(n_k^{(1)})_1^\infty$, $(n_k^{(2)})_1^\infty$ such that
\begin{equation}
\label{edelstein}
g^{n_k^{(1)}}(\0) \tendsto k \0 \text{ but } \|g^{n_k^{(2)}}(\0)\| \tendsto k \infty.
\end{equation}
\end{proposition}
Since the specific form of Edelstein's example will be important to us, we recall the proof of Proposition \ref{propositionedelstein}, in a modified form suitable for generalization:
\begin{proof}[Proof of Proposition \ref{propositionedelstein}]
For each $k\in\N$ let $a_k = 1/k!$, let $b_k = 1$, and let
\begin{equation}
\label{edelsteindef1}
c_k = e^{2\pi i a_k}, \;\; d_k = b_k (1 - c_k).
\end{equation}
Then
\[
\sum_{k\in\N}|d_k|^2 \lesssim_\times \sum_{k\in\N}|a_k b_k|^2 = \sum_{k\in\N} \left(\frac{1}{k!}\right)^2 < \infty,
\]
so $\dd = (d_k)_1^\infty\in\ell^2(\N;\C)$. Let $g\in\Isom(\ell^2(\N;\C))$ be given by the formula
\begin{equation}
\label{edelsteindef2}
g(\xx)_k = c_k x_k + d_k.
\end{equation}
Then
\begin{equation}
\label{homomorphism}
g^n(\xx)_k = c_k^n x_k + \sum_{i = 0}^{n - 1} c_k^i d_k = c_k^n x_k + \frac{1 - c_k^n}{1 - c_k} d_k = c_k^n x_k + b_k(1 - c_k^n).
\end{equation}
In particular, $g^n(\0)_k = b_k(1 - c_k^n)$. So
\begin{equation}
\label{gn0norm}
\|g^n(\0)\|^2 = \sum_{k = 1}^\infty |b_k(1 - c_k^n)|^2
= \sum_{k = 1}^\infty |b_k(1 - e^{2\pi i na_k})|^2
\asymp_\times \sum_{k = 1}^\infty |b_k|^2\dist(n a_k,\Z)^2.
\end{equation}
Now for each $k\in\N$, let
\begin{align*}
n_k^{(1)} &= k! \\
n_k^{(2)} &= \frac12\sum_{j = 1}^k j!
\end{align*}
Then
\begin{align*}
\|g^{n_k^{(1)}}(\0)\|^2 \asymp_\times \sum_{j = 1}^\infty \dist\left(\frac{k!}{j!},\Z\right)^2
&=_\pt (k!)^2\sum_{j = k + 1}^\infty \left(\frac{1}{j!}\right)^2\\
&\asymp_\times \left(\frac{k!}{(k+1)!}\right)^2\\ 
&= \frac{1}{(k + 1)^2} \tendsto k 0,
\end{align*}
but on the other hand
\begin{align*}
\|g^{n_k^{(2)}}(\0)\|^2 \gtrsim_\times \sum_{j = 1}^k \dist\left(\frac{n_k^{(2)}}{(j + 1)!},\Z\right)^2
&= \sum_{j = 1}^k \frac 14\left[1 - \sum_{i = 1}^j \frac{i!}{(j + 1)!}\right]^2\\
&\geq \frac 14 \sum_{j = 1}^k \left[1 - \frac{2}{j + 1}\right]^2 \tendsto k \infty.
\end{align*}
This demonstrates \eqref{edelstein}.
\end{proof}
\begin{remark}
\label{remarkedelsteinsignificance}
Let us explain the significance of Edelstein's example from the point of view of hyperbolic geometry. Let $g\in\Isom(\BB)$ be as in Proposition \ref{propositionedelstein}, and let $\what g\in \Isom(\E^\infty)$ be its Poincar\'e extension. By Observation \ref{observationisomE}, $\what g$ is a parabolic isometry. But the orbit of $\zero = (1,\0)$ (cf. \6\ref{standingassumptions2}) is quite irregular; indeed,
\[
\what g^{\;n_{k}^{(1)}}(\zero) \tendsto k \zero \text{ but } \what g^{\;n_{k}^{(2)}}(\zero) \tendsto k \infty\in\del\E.
\]
So the orbit $(g^n(\zero))_1^\infty$ simultaneously tends to infinity on one subsequence while remaining bounded on another subsequence. Such a phenomenon cannot occur in proper metric spaces, as we demonstrate now:
\end{remark}
\begin{theorem}
If $X$ is a proper metric space and if $G\leq\Isom(X)$ is cyclic, then either $G$ has bounded orbits or $G$ is strongly discrete.
\end{theorem}
\begin{proof}
Write $G = g^\Z$ for some $g\in\Isom(X)$, and fix a point $\zero\in X$. For each $n\in\Z$ write $\|n\| = \dogo{g^n}$. Then $\|-n\| = \|n\|$, and $\|m + n\| \leq \|m\| + \|n\|$.

Suppose that $G$ is not strongly discrete. Then there exists $R > 0$ such that
\begin{equation}
\label{Rdef}
\#\{n\in\N : \|n\| \leq R\} = \infty.
\end{equation}
Now let $g^S(\zero) \subset g^\Z(\zero)\cap B(\zero,2R)$ be a maximal $R$-separated set. Since $X$ is proper, $S$ is finite. For each $k\in S$, choose $\ell_k > k$ such that $\|\ell_k\| \leq R$; such an $\ell_k$ exists by \eqref{Rdef}.

Now let $n\in\N$ be arbitrary. Let $0\leq m \leq n$ be the largest number for which $\|m\| \leq 2R$. Since $g^S(\zero)$ is a maximal $R$-separated subset of $g^\Z(\zero)\cap B(\zero,2R)\ni g^m(\zero)$, there exists $k\in S$ for which $\|m - k\| \leq R$. Then
\[
\|m - k + \ell_k\| \leq R + R = 2R.
\]
On the other hand, $m - k + \ell_k > m$ since $\ell_k > k$ by construction. Thus by the maximality of $m$, we have $m - k + \ell_k > n$. So
\[
n - m < \ell_k - k \leq C := \max_{k\in S} (\ell_k - k).
\]
It follows that
\[
\|n\| \leq \|m\| + \|n - m\| \leq 2R + C\dogo g,
\]
i.e. $\|n\|$ is bounded independent of $n$. Thus $G$ has bounded orbits.
\end{proof}

At this point, we shall use the different notions of discreteness introduced in Chapter \ref{sectiondiscreteness} to distinguish between different variations of Edelstein's example. To this end, we make the following definition:

\begin{definition}
\label{definitionedelstein}
An isometry $g\in\Isom(\ell^2(\N;\C))$ is said to be \emph{Edelstein-type} if it is of the form \eqref{edelsteindef2}, where the sequences $(c_k)_1^\infty$ and $(d_k)_1^\infty$ are of the form \eqref{edelsteindef1}, where $(a_k)_1^\infty$ and $(b_k)_1^\infty$ are sequences of positive real numbers satisfying
\[
\sum_{k = 1}^\infty |a_k b_k|^2 < \infty.
\]
\end{definition}

Our proof of Proposition \ref{propositionedelstein} shows that the isometry $g$ is always well-defined and satisfies \eqref{gn0norm}. On the other hand, the conclusion of Proposition \ref{propositionedelstein} does not hold for all Edelstein-type isometries; it is possible that the cyclic group $G = g^\Z$ is strongly discrete, and it is also possible that this group has bounded orbits. (But the two cannot happen simultaneously unless $g$ is a torsion element.) In the sequel, we will be interested in Edelstein-type isometries for which $G$ has unbounded orbits but is not necessarily strongly discrete. We will be able to distinguish between the examples using our more refined notions of discreteness.

\begin{edelsteinexample}
\label{exampleedelstein}
In our notation, Edelstein's original example can be described as the Edelstein-type isometry $g$ defined by the sequences $a_k = 1/k!$, $b_k = 1$. Edelstein's proof shows that $G = g^\Z$ has unbounded orbits and is not weakly discrete. However, we can show more:
\end{edelsteinexample}
\begin{proposition}
\label{propositionedelsteinUOT}
The cyclic group generated by the isometry in Edelstein's example is not UOT-discrete.
\end{proposition}
\begin{proof}
As in the proof of Proposition \ref{propositionedelstein}, we let $n_k = k!$, so that $g^{n_k}(\0)\to \0$. But if $T^n$ denotes the linear part of $g^n$, then
\[
T^{n_k}(\xx) = (e^{2\pi i k!/j!} x_j)_{j = 1}^\infty
\]
and so
\[
\|T^{n_k} - I\| \leq \sum_{j = k + 1}^\infty |1 - e^{2\pi i k!/j!}| \asymp_\times \frac{1}{k + 1} \to 0.
\]
Thus $T^{n_k}\to I$ in the uniform operator topology, so by Observation \ref{observationisomE}(iii), $\what g^{n_k}\to\id$ in the uniform operator topology. Thus $\what g^\Z$ is not UOT-discrete.
\end{proof}

\begin{edelsteinexample}
\label{examplevalette}
The Edelstein-type isometry $g$ defined by the sequences $a_k = 1/2^k$, $b_k = 1$. This example was considered by A. Valette \cite[Proposition 1.7]{Valette}. It has unbounded orbits, and is moderately discrete (in fact properly discontinuous) but not strongly discrete.
\end{edelsteinexample}
\begin{proof}
Letting $n_k^{(1)} = 2^k$, we have by \eqref{gn0norm}
\[
\|g^{n_k^{(1)}}(\0)\|^2 \asymp_\times \sum_{j = 1}^\infty \dist\left(\frac{2^k}{2^j},\Z\right)^2
= \sum_{j = k + 1}^\infty (2^{k - j})^2 = \frac 13
\]
and so $g^\Z$ is not strongly discrete. Letting $n_k^{(2)} = \lfloor 2^k/3\rfloor$, we have
\begin{align*}
\|g^{n_k^{(2)}}(\0)\|^2 \asymp_\times \sum_{j = 1}^\infty \dist\left(\frac{\lfloor 2^k/3\rfloor}{2^j},\Z\right)^2
&\geq \sum_{j = 1}^k \left[\dist\left(\frac{2^{k - j}}{3},\Z\right) - \frac{1}{2^j}\right]^2\\ 
&\geq \sum_{j = 7}^k \frac{1}{100}\\ 
&\asymp_{\plus,\times} k \tendsto k \infty
\end{align*}
and therefore $g^\Z$ has unbounded orbits.

Finally, we show that $g^\Z$ acts properly discontinuously. To begin with, we observe that for all $n\in\Namer$, we may write $n = 2^k (2j + 1)$ for some $j,k\geq 0$; then
\[
\|g^n(\0)\|^2 \asymp_\times \sum_{i = 1}^\infty \dist\left(\frac{2^k (2j + 1)}{2^i},\Z\right)^2
\geq \dist\left(\frac{2^k (2j + 1)}{2^{k + 1}},\Z\right)^2 = 1/4,
\]
i.e. $\0$ is an isolated point of $g^\Z(\0)$. So for some $\epsilon > 0$,
\[
\|g^n(\0)\| \geq \epsilon \all n\in\Namer.
\]
Now let $\xx\in\ell^2(\N;\C)$ be arbitrary, and let $N$ be large enough so that 
\[
\|(x_{N + 1},\ldots)\| \leq \epsilon/3.
\] 
Now for all $n\in\Namer$, we have 
\begin{align*}
\|g^{2^N n}(\xx) - \xx\| &= \|g^{2^N n}(0,\ldots,0,x_{N + 1},\ldots) - (0,\ldots,0,x_{N + 1},\ldots)\|\\ 
&\geq \|g^{2^N n}(\0)\| - 2\epsilon/3 \geq \epsilon/3
\end{align*}
which implies that the set $g^{2N \Z}(\xx)$ is discrete. But $g^\Z(\xx)$ is the union of finitely many isometric images of $g^{2^N\Z}(\xx)$, so it must also be discrete.
\end{proof}
\begin{remark}
It is not possible to differentiate further between unbounded Edelstein-type isometries by considering separately the conditions of weak discreteness, moderate discreteness, and proper discontinuity. Indeed, if $X$ is any metric space and if $G\leq\Isom(X)$ is any cyclic group with unbounded orbits, then the following are equivalent: $G$ is weakly discrete, $G$ is moderately discrete, $G$ acts properly discontinuously. This can be seen as follows: every nontrivial subgroup of $G$ is of finite index, and therefore also has unbounded orbits; it follows that no element of $G\butnot\{\id\}$ has a fixed point in $X$; it follows from this that the three notions of discreteness are equivalent.
\end{remark}

\begin{example}
Let $g\in\Isom(\ell^2(\N;\C))$ be as in Proposition \ref{propositionedelstein}, let $\sigma:\ell^2(\Z;\C)\to\ell^2(\Z;\C)$ be the shift map $\sigma(\xx)_k = x_{k + 1}$, and let $T:\ell^2(\N;\C)\to\ell^2(\N;\C)$ be given by the formula
\[
T(\xx)_k = e^{2\pi i/k} x_k.
\]
Then $g_1 = g\oplus\sigma$ has unbounded orbits and is COT-discrete but not weakly discrete; $g_2 = g\oplus T$ has unbounded orbits and is UOT-discrete but not COT-discrete.
\end{example}
\begin{proof}
Since $g$ has unbounded orbits and is not weakly discrete, the same is true for both $g_1$ and $g_2$. Since the sequence $(\sigma^n(\xx))_1^\infty$ diverges for every $\xx\in\ell^2(\Z;\C)$, the group generated by $\sigma$ is COT-discrete, which implies that $g_1$ is as well. Since the sequence $(\|T^n - I\|)_1^\infty$ is bounded from below, the group generated by $T$ is UOT-discrete, which implies that $g_2$ is as well. On the other hand, if we let $n_k = k!$, then $T^{n_k}(\xx)\to \xx$ for all $\xx\in\ell^2(\Z;\C)$. But we showed in Proposition \ref{propositionedelsteinUOT} that $g^{n_k}(\xx)\to \xx$ for all $\xx\in\ell^2(\N;\C)$; it follows that $g\oplus T$ is not COT-discrete.
\end{proof}

\begin{remark}
One might object to the above examples on the grounds that the isometries $g_1$ and $g_2$ do not act irreducibly. However, Edelstein-type isometries never act irreducibly: if $g$ is defined by \eqref{edelsteindef1} and \eqref{edelsteindef2} for some sequences $(a_k)_1^\infty$ and $(b_k)_1^\infty$, then for every $k$ the affine hyperplane $H_k = \{\xx\in \ell^2(\N;\C) : x_k = b_k\}$ is invariant under $g$. In general, it is not even possible to find a minimal subspace on which the restricted action of $g$ is irreducible, since such a minimal subspace would be given by the formula
\[
\bigcap_k H_k = \begin{cases}
\emptyset & \sum_1^\infty |b_k|^2 = \infty\\
\{(b_k)_1^\infty\} & \sum_1^\infty |b_k|^2 < \infty
\end{cases},
\]
and if $g$ has unbounded orbits (as in Examples \ref{exampleedelstein} and \ref{examplevalette}), the first case must hold.
\end{remark}

We conclude this section with one more Edelstein-type example:

\begin{edelsteinexample}
\label{exampleparabolicinfinite}
The Edelstein-type isometry $g$ defined by the sequences $a_k = 1/2^k$, $b_k = \log(1 + k)$. In this example, $g^\Z$ is strongly discrete but has infinite Poincar\'e exponent.
\end{edelsteinexample}
\begin{proof}
To show that $g^\Z$ is strongly discrete, fix $n\geq 1$, and let $k$ be such that $2^k\leq n < 2^{k + 1}$. Then $1/4 \leq n/2^{k + 2} < 1/2$, so by \eqref{gn0norm},
\[
\|g^n(\0)\|^2 \gtrsim_\times b_{k + 2} \,\dist\left(\frac{n}{2^{k + 2}},\Z\right)^2 \geq \frac{b_{k + 2}}{16} \tendsto n \infty.
\]
To show that $\delta(g^\Z) = \infty$, fix $\ell\geq 0$, and note that by \eqref{gn0norm},
\[
\|g^{2^\ell}(\0)\|^2 \asymp_\times \sum_{k = \ell + 1}^\infty \frac{4^\ell}{4^k} |b_k|^2 \asymp_\times |b_{\ell + 1}|^2 = \log^2(2 + \ell).
\]
It follows that
\[
\Sigma_s(g^\Z) \geq \sum_{\ell = 0}^\infty (1\vee \|g^{2^\ell}(\0)\|)^{-2s} \asymp_\times \sum_{\ell = 0}^\infty \log^{-2s}(2 + \ell) = \infty \all s \geq 0.
\]
\end{proof}

%
%

%
%

%
%

\section{The Poincar\'e exponent of a finitely generated parabolic group}
\label{subsectionparabolicexponent}
In this section, we relate the Poincar\'e exponent $\delta_G$ of a parabolic group $G$ with its algebraic structure. We will show below that $\delta_G$ is infinite unless $G$ is virtually nilpotent (Theorem \ref{theoremparaboliclowerbound} below), so we begin with a digression on the coarse geometry of finitely generated virtually nilpotent groups.

\subsection{Nilpotent and virtually nilpotent groups}
\label{subsubsectionnilpotent}
Recall that the \emph{lower central series} of an abstract group $\Gamma$ is the sequence $(\Gamma_i)_1^\infty$ defined recursively by the equations
\[
\Gamma_1 = \Gamma \text{ and } \Gamma_{i + 1} = [\Gamma,\Gamma_i].
\]
Here $[A,B]$ denotes the commutator of two sets $A,B\subset \Gamma$, i.e. $[A,B] = \lb aba^{-1}b^{-1} : a\in A, b\in B\rb$. The group $\Gamma$ is \emph{nilpotent} if its lower central series terminates, i.e. if $\Gamma_{k + 1} = \{\id\}$ for some $k\in\N$. The smallest integer $k$ for which this equality holds is called the \emph{nilpotency class} of $\Gamma$, and we will denote it by $k$.

Note that a group is abelian if and only if it is nilpotent of class $1$. The fundamental theorem of finitely generated abelian groups says that if $\Gamma$ is a finitely generated abelian group, then $\Gamma \equiv \Z^d \times F$ for some $d\in\Neur$ and some finite abelian group $F$. The number $d$ will be called the \emph{rank} of $\Gamma$, denoted $\rank(\Gamma)$. Note that the large-scale structure of $\Gamma$ depends only on $d$ and not on the finite group $F$. Specifically, if $\dist_\Gamma$ is a Cayley metric on $\Gamma$ then
\begin{equation}
\label{bassguivarchabelian}
\NN_\Gamma(R) \asymp_\times R^d \all R \geq 1.
\end{equation}
Here $\NN_\Gamma(R) = \#\{\gamma\in\Gamma : \dist(e,\gamma)\leq R\}$ is the orbital counting function of $\Gamma$ interpreted as acting on the metric space $(\Gamma,\dist_\Gamma)$ (cf. Remark \ref{remarkorbitalcounting}).

The following analogue of \eqref{bassguivarchabelian} was proven by H. Bass and independently by Y. Guivarch:

\begin{theorem}[\cite{Bass, Guivarch}]
\label{theorembassguivarch}
Let $\Gamma$ be a finitely generated nilpotent group with lower central series $(\Gamma_i)_1^\infty$ and nilpotency class $k$, and let
\[
\alpha_\Gamma = \sum_{i = 1}^k i\rank(\Gamma_i/\Gamma_{i + 1}).
\]
Let $\dist_\Gamma$ be a Cayley metric on $\Gamma$. Then for all $R\geq 1$,
\begin{equation}
\label{bassguivarch}
\NN_\Gamma(R) \asymp_\times R^{\alpha_\Gamma}.
\end{equation}
\end{theorem}
The number $\alpha_\Gamma$ will be called the \emph{(polynomial) growth rate} of $\NN_\Gamma$.

A group is \emph{virtually nilpotent} if it has a nilpotent subgroup of finite index. The following is an immediate corollary of Theorem \ref{theorembassguivarch}:
\begin{corollary}
\label{corollarybassguivarch}
Let $\Gamma$ be a finitely generated virtually nilpotent group. Let $\Gamma'\leq\Gamma$ be a nilpotent subgroup of finite index, and let $\dist_\Gamma$ be a Cayley metric. Let $\alpha_\Gamma = \alpha_{\Gamma'}$. Then for all $R\geq 1$,
\begin{equation}
\label{bassguivarch2}
\NN_\Gamma(R) \asymp_\times R^{\alpha_\Gamma}.
\end{equation}
\end{corollary}

\begin{example}
If $\Gamma$ is abelian, then \eqref{bassguivarch} reduces to \eqref{bassguivarchabelian}.
\end{example}

\begin{example}
\label{exampleheisenberg}
Let $\Gamma$ be the discrete Heisenberg group, i.e.
\[
\Gamma = \left\{\left[\begin{array}{ccc}
1 & a & c\\
& 1 & b\\
&& 1
\end{array}\right] : a,b,c\in\Z\right\}.
\]
We compute the growth rate of $\NN_\Gamma$. Note that $\Gamma$ is nilpotent of class $2$, and its lower central series is given by $\Gamma_1 = \Gamma$,
\[
\Gamma_2 = \left\{\left[\begin{array}{ccc}
1 & & c\\
& 1 &\\
&& 1
\end{array}\right] : c\in\Z\right\}.
\]
Thus
\[
\alpha_\Gamma = \rank(\Gamma_1/\Gamma_2) + 2\rank(\Gamma_2) = 2 + 2\cdot 1 = 4.
\]
\end{example}


Corollary \ref{corollarybassguivarch} implies that finitely generated virtually nilpotent groups have \emph{polynomial growth}, meaning that the growth rate
\begin{equation}
\label{growthrate}
\alpha_\Gamma := \lim_{R\to\infty}\frac{\log \NN_\Gamma(R)}{\log(R)}
\end{equation}
exists and is finite. The converse assertion is a deep theorem of M. Gromov:

\begin{theorem}[\cite{ShalomTao}]
\label{theoremgromov}
A finitely generated group $\Gamma$ has polynomial growth if and only if $\Gamma$ is virtually nilpotent. Moreover, if $\Gamma$ does not have polynomial growth then the limit \eqref{growthrate} exists and equals $\infty$.
\end{theorem}

Thus the limit \eqref{growthrate} exists in all circumstances, so we may refer to it unambiguously.

\subsection{A universal lower bound on the Poincar\'e exponent}
Now let $G\leq\Isom(X)$ be a parabolic group. Recall that in the Standard Case, if a group $G$ is discrete then it is virtually abelian. Moreover, in this case $\delta_G = \frac12\rank(G)$.

If $G$ is virtually nilpotent, then it is natural to replace this formula by the formula $\delta_G = \frac12\alpha_G$. However, in general equality does not hold in this formula, as we will see below (Theorem \ref{theoremnilpotentembedding}). We show now that the $\geq$ direction always holds. Precisely:

\begin{theorem}
\label{theoremparaboliclowerbound}
Let $G\leq\Isom(X)$ be a finitely generated parabolic group. Let $\alpha_G$ be as in \eqref{growthrate}. Then
\begin{equation}
\label{paraboliclowerbound}
\delta_G \geq \frac{\alpha_G}{2}\cdot
\end{equation}
Moreover, if equality holds and $\delta_G < \infty$, then $G$ is of divergence type.
\end{theorem}

Before proving Theorem \ref{theoremparaboliclowerbound}, we make a few remarks:

\begin{remark}
In this theorem, it is crucial that $b > 1$ is chosen close enough to $1$ so that Proposition \ref{propositionDist} holds (cf. \6\ref{standingassumptions2}). Indeed, by varying the parameter $b$ one may vary the Poincar\'e exponent at will (cf. \eqref{poincarealternate}); in particular, by choosing $b$ large, one could make $\delta_G$ arbitrarily small. If $X$ is strongly hyperbolic, then of course we may let $b = e$.
\end{remark}
\begin{remark}
Expanding on the above remark, we recall that if $X$ is an $\R$-tree, then any value of $b$ is permitted in Proposition \ref{propositionDist} (Remark \ref{remarkDistRtree}). This demonstrates that if a finitely generated parabolic group acting on an $\R$-tree has finite Poincar\'e exponent, then its growth rate is zero. This may also be seen more directly from Remark \ref{remarkparabolicRtree2}.
\end{remark}
\begin{remark}
Let $G\leq\Isom(X)$ be a group of general type, and let $H\leq G$ be a finitely generated parabolic subgroup. Then combining Theorem \ref{theoremparaboliclowerbound} with Proposition \ref{propositiondeltagtrdeltap} shows that $\delta_G > \alpha_H/2$. This generalizes a well-known theorem of A. F. Beardon \cite[Theorem 7]{Beardon1}.
\end{remark}

Combining Theorems \ref{theoremgromov} and \ref{theoremparaboliclowerbound} gives the following corollary:

\begin{corollary}
\label{corollaryfiniteimpliesnilpotent}
Any finitely generated parabolic group with finite Poincar\'e exponent is virtually nilpotent.
\end{corollary}
\noindent This corollary can be viewed very loosely as a generalization of Margulis's lemma (Proposition \ref{propositionmargulislemma}). As we have seen above (Observation \ref{observationmargulislemma}), a strict analogue of Margulis's lemma fails in infinite dimensions.
\begin{proof}[Proof of Theorem \ref{theoremparaboliclowerbound}]
Let $g_1,\ldots,g_n$ be a set of generators for $G$, and let $\dist_G$ denote the corresponding Cayley metric. Let $\xi\in\del X$ denote the unique fixed point of $G$. Fix $g\in G$, and write $g = g_{i_1}\cdots g_{i_m}$. By the universal property of path metrics (Remark \ref{remarkpathmetricuniversal}), we have
\[
\Dist_\xi(\zero,g(\zero)) \lesssim_\times \dist_G(\id,g).
\]
Now we apply Observation \ref{observationeuclideanparabolicasymp} to get
\[
b^{(1/2)\dogo g} \lesssim_\times \dist_G(\id,g).
\]
Letting $C > 0$ be the implied constant, we have
\begin{equation}
\label{NXGbound}
\NN_{X,G}(\rho) \geq \NN_G(b^{\rho/2}/C) \all \rho > 0
\end{equation}
(cf. Remark \ref{remarkorbitalcounting}). In particular, by \eqref{poincarealternate}
\[
\delta_G = \lim_{\rho\to\infty}\frac{\log_b\NN_{X,G}(\rho)}{\rho} \geq \lim_{R\to\infty}\frac{\log_b\NN_G(R)}{2\log_b(R)} = \frac{\alpha_G}{2}\cdot
\]
To demonstrate the final assertion of Theorem \ref{theoremparaboliclowerbound}, suppose that equality holds in \eqref{paraboliclowerbound} and that $\delta_G < \infty$. Then by Theorem \ref{theoremgromov}, $G$ is virtually nilpotent. Combining \eqref{NXGbound} with \eqref{bassguivarch} and then plugging into \eqref{poincareseriesalternate} gives us that $\Sigma_\delta(G) = \infty$, completing the proof.
\end{proof}

\subsection{Examples with explicit Poincar\'e exponents}
\label{subsubsectiontheoremnilpotent}
Theorem \ref{theoremparaboliclowerbound} raises a natural question: do the exponents allowed by this theorem actually occur as the Poincar\'e exponent of some parabolic group? More precisely, given a finitely generated abstract group $\Gamma$ and a number $\delta \geq \alpha_\Gamma/2$, does there exist a hyperbolic metric space $X$ and an injective homomorphism $\Phi:\Gamma\to\Isom(X)$ such that $G = \Phi(\Gamma)$ is a parabolic group satisfying $\delta_G = \delta$? If $\delta = \alpha_\Gamma/2$, then the problem appears to be difficult; cf. Remark \ref{remarknilpotentembedding}. However, we can provide a complete answer when $\delta > \alpha_\Gamma/2$ by embedding $\Gamma$ into $\Isom(\BB)$ and then using Poincar\'e extension to get an embedding into $\Isom(\E^\infty)$. Specifically, we have the following:

\begin{theorem}
\label{theoremnilpotentembedding}
Let $\Gamma$ be a virtually nilpotent group, and let $\alpha = \alpha_\Gamma$ be the growth rate of $\NN_\Gamma$. Then for all $\delta > \alpha_\Gamma/2$, there exists an injective homomorphism $\Phi:\Gamma\to\Isom(\BB)$ such that
\[
\delta(\Phi(\Gamma)) = \delta.
\]
Moreover, $\Phi(\Gamma)$ may be taken to be either of convergence type or of divergence type.
\end{theorem}

\begin{remark}
\label{remarknilpotentembedding}
Theorem \ref{theoremnilpotentembedding} raises the question of whether there exists an injective homomorphism $\Phi:\Gamma\to\Isom(\BB)$ such that
\begin{equation}
\label{nilpotentequality}
\delta(\Phi(\Gamma)) = \alpha_\Gamma/2.
\end{equation}
It is readily computed that if the map $\gamma\mapsto \Phi(\gamma)(\0)$ is bi-Lipschitz, then \eqref{nilpotentequality} holds. In particular, if $\Gamma = \Z^d$ for some $d\in\N$, then such a $\Phi$ is given by $\Phi(\nn)(\xx) = \xx + (\nn,\0)$. By contrast, if $\Gamma$ is a virtually nilpotent group which is not virtually abelian, then it is known \cite[Theorem 1.3]{CTV1} that there is no quasi-isometric embedding $\iota:\Gamma\to\BB$ (see also \cite[Theorem A]{Pauls} for an earlier version of this result which applied to nilpotent Lie groups). In particular, there is no homomorphism $\Phi:\Gamma\to\Isom(\BB)$ such that $\gamma\mapsto \Phi(\gamma)(\0)$ is a bi-Lipschitz embedding. So this approach of constructing an injective homomorphism $\Phi$ satisfying \eqref{nilpotentequality} is doomed to failure. However, it is possible that another approach will work. We leave the question as an open problem.
\end{remark}
\begin{remark}
\label{remarknilpotentembedding2}
Letting $\Gamma = \Z$ in Theorem \ref{theoremnilpotentembedding}, we have the following corollary: For any $\delta > 1/2$, there exists an isometry $g_\delta\in\Isom(\BB)$ such that the cyclic group $G_\delta = (g_\delta)^\Z$ satisfies $\delta(G_\delta) = \delta$, and may be taken to be either of convergence type or of divergence type. The isometries $(g_\delta)_{\delta > 1/2}$ exhibit ``intermediate'' behavior between the isometry $g_{1/2}(\xx) = \xx + \ee_1$ (which has Poincar\'e exponent $1/2$ as noted above) and the isometries described in the Edelstein-type isometries \ref{exampleedelstein}, \ref{examplevalette}, and \ref{exampleparabolicinfinite}: since $\delta > 1/2$, the sequence $(g_\delta^n(\0))_1^\infty$ converges to infinity much more slowly than the sequence $(g_{1/2}^n(\0))_1^\infty$, but since $\delta < \infty$, the sequence converges faster than in Example \ref{exampleparabolicinfinite}, not to mention Examples \ref{exampleedelstein} and \ref{examplevalette} where the sequence $(g_\delta^n(\0))_1^\infty$ does not converge to infinity at all (although it converges along a subsequence).
\end{remark}

\begin{remark}
\label{remarkheisenbergdeltainfty}
Theorem \ref{theoremnilpotentembedding} leaves open the question of whether there is a homomorphism $\Phi:\Gamma\to\Isom(\BB)$ such that $\Phi(\Gamma)$ is strongly discrete but $\delta(\Phi(\Gamma)) = \infty$. If $\Gamma = \Z$, this is answered affirmatively by Example \ref{exampleparabolicinfinite}, and if $\Gamma$ contains $\Z$ as a direct summand, i.e. $\Gamma \equiv \Z\times\Gamma'$ for some $\Gamma'\leq\Gamma$, then the answer can be achieved by taking the direct sum of Example \ref{exampleparabolicinfinite} with an arbitrary strongly discrete homomorphism from $\Gamma'$ to $\Isom(\BB)$. However, the Heisenberg group does not contain $\Z$ as a direct summand. Thus, it is unclear whether or not there is a a homomorphism from the Heisenberg group to $\Isom(\BB)$ whose image is strongly discrete with infinite Poincar\'e exponent.
\end{remark}

\begin{proof}[Proof of Theorem \ref{theoremnilpotentembedding}]
We will need the following variant of the Assouad embedding theorem:


\begin{theorem}
\label{theoremassouad}
Let $X$ be a doubling metric space,\Footnote{Recall that a metric space $X$ is \emph{doubling} if there exists $M > 0$ such that for all $x\in X$ and $\rho > 0$, the ball $B(x,\rho)$ can be covered by $M$ balls of radius $\rho/2$.\label{footnotedoubling}} and let $F:(0,\infty)\to(0,\infty)$ be a nondecreasing function such that
\begin{equation}
\label{assouadhypotheses}
0 < \lexp(F) \leq \uexp(F) < 1.
\end{equation}
Here
\begin{align*}
\lexp(F) &:= \liminf_{\lambda\to\infty} \inf_{R > 0} \frac{\log F(\lambda R) - \log F(R)}{\log(\lambda)}\\
\uexp(F) &:= \limsup_{\lambda\to\infty} \sup_{R > 0} \frac{\log F(\lambda R) - \log F(R)}{\log(\lambda)}\cdot
\end{align*}
Then there exist $d\in\N$ and a map $\iota:X\to\R^d$ such that for all $x,y\in X$,
\begin{equation}
\label{assouad}
\|\iota(y) - \iota(x)\| \asymp_\times F(\dist(x,y)).
\end{equation}
\end{theorem}
\begin{subproof}
The classical Assouad embedding theorem (see e.g. \cite[Theorem 12.2]{Heinonen}) gives the special case of Theorem \ref{theoremassouad} where $F(t) = t^\epsilon$ for some $0 < \epsilon < 1$. It is possible to modify the standard proof of the classical version in order to accomodate more general functions $F$ satisfying \eqref{assouadhypotheses}; however, we prefer to prove Theorem \ref{theoremassouad} directly as a consequence of the classical version.

Fix $\epsilon \in (\uexp(F),1)$, and let
\[
\what F(t) = t^\epsilon \inf_{s \leq t} \frac{F(s)}{s^\epsilon}\cdot
\]
The inequality $\epsilon > \lexp(f)$ implies that $\what F\asymp_\times F$, so we may replace $F$ by $\what F$ without affecting either the hypotheses or the conclusion of the theorem. Thus, we may without loss of generality assume that the function $t\mapsto F(t)/t^\epsilon$ is nonincreasing.

Let $G(t) = F(t)^{1/\epsilon}$, so that $t\mapsto G(t)/t$ is nonincreasing. It follows that
\[
G(t + s) \leq G(t) + G(s).
\]
Combining with the fact that $G$ is nondecreasing shows that $G\circ\dist$ is a metric on $X$. On the other hand, since $\lexp(G) = \lexp(F)/\epsilon > 0$, there exists $\lambda > 0$ such that $G(\lambda t) \geq 2G(t)$ for all $t > 0$. It follows that the metric $G\circ\dist$ is doubling. Thus we may apply the classical Assouad embedding theorem to the metric space $(X,G\circ\dist)$ and the function $t\mapsto t^\epsilon$, giving a map $\iota:X\to \R^d$ satisfying
\[
\|\iota(y) - \iota(x)\| \asymp_\times G^\epsilon \circ\dist(x,y) = F(\dist(x,y)).
\]
This completes the proof.
\end{subproof}

Now let $\Gamma$ be a virtually nilpotent group, and let $\dist_\Gamma$ be a Cayley metric on $\Gamma$.
\begin{lemma}
$(\Gamma,\dist_\Gamma)$ is a doubling metric space.
\end{lemma}
\begin{subproof}
For all $\gamma\in\Gamma$ and $R > 0$, we have by Corollary \ref{corollarybassguivarch}
\[
\#(B(\gamma,R)) = \#(\gamma(B(e,R))) = \#(B(e,R)) \asymp_\times (1\vee R)^{\alpha_\Gamma}.
\]
Now let $S\subset B(\gamma,2R)$ be a maximal $R$-separated set. Then $\{B(\beta,R) : \beta\in S\}$ is a cover of $B(\gamma,2R)$. On the other hand, $\{B(\beta,R/2) : \beta\in S\}$ is a disjoint collection of subsets of $B(\gamma,3R)$, so
\begin{align*}
\sum_{\beta\in S}\#(B(\beta,R/2)) &\leq_\pt \#(B(\gamma,3R))\\
\#(S) \leq \frac{\#(B(\gamma,3R))}{\min_{\beta\in S}\#(B(\beta,R/2))}
&\asymp_\times \frac{(1\vee 3R)^{\alpha_\Gamma}}{(1\vee R/2)^{\alpha_\Gamma}} \asymp_\times 1,
\end{align*}
i.e. $\#(S) \leq M$ for some $M$ independent of $\gamma$ and $R$. But then $B(\gamma,2R)$ can be covered by $M$ balls of radius $R$, proving that $\Gamma$ is doubling.
\end{subproof}

Now let $f:\CO 1\infty\to \CO 1\infty$ be a continuous increasing function satisfying
\begin{equation}
\label{fhypothesis}
\alpha < \lexp(f) \leq \uexp(f) < \infty
\end{equation}
and $f(1) = 1$. Let
\[
F(R) = \begin{cases}
f^{-1}(R^\alpha) & R \geq 1\\
R^{1/2} & R \leq 1
\end{cases}.
\]
Then
\[
0 < \lexp(F) = \min\left(\frac12,\frac{\alpha}{\uexp(f)}\right) \leq \uexp(F) = \max\left(\frac12,\frac{\alpha}{\lexp(f)}\right) < 1.
\]
Thus $F$ satisfies the hypotheses of Theorem \ref{theoremassouad}, so there exists an embedding $\iota:\Gamma\to\HH$ satisfying \eqref{assouad}. By \cite[Proposition 4.4]{CTV1}, we may without loss of generality assume that $\iota(\gamma) = \Phi(\gamma)(\0)$ for some homomorphism $\Phi:\Gamma\to\Isom(\BB)$. Now for all $R\geq 1$,
\begin{align*}
\NN_{\BB,\Phi(\Gamma)}(R) &= \#\{\gamma\in\Gamma : \Dist_\xi(\zero,\Phi(\gamma)(\zero)) \leq R\}\\
&= \#\{\gamma\in\Gamma : F(\dist_\Gamma(e,\gamma)) \leq R\}\\
&= \NN_\Gamma(F^{-1}(R)) \asymp_\times \big(F^{-1}(R)\big)^\alpha = f(R).
\end{align*}
In particular, given $\delta > \alpha_\Gamma/2$ and $k\in \{0,2\}$, we can let 
\[
f(R) = R^{2\delta} (1 + \log(R))^{-k}.
\] 
It is readily verified that $\alpha < \dexp(f) = 2\delta < \infty$, so in particular \eqref{fhypothesis} holds. By \eqref{poincarealternate}, $\delta(\Phi(\Gamma)) = \delta$ and by \eqref{poincareseriesalternate}, $\Phi(\Gamma)$ is of divergence type if and only if $k = 0$.
\end{proof}
\begin{remark}
\label{remarknilpotentembedding3}
The above proof shows a little more that what was promised; namely, it has been shown that
\begin{itemize}
\item[(i)] for every function $F:(0,\infty)\to(0,\infty)$ satisfying \eqref{assouadhypotheses}, there exists an injective homomorphism $\Phi:\Gamma\to \BB$ such that $\|\Phi(\gamma)(\0)\| \asymp_\times F(\dist(e,\gamma))$ for all $\gamma\in \Gamma$, and that
\item[(ii)] for every function $f:\CO 1\infty\to \CO 1\infty$ satisfying \eqref{fhypothesis}, there exists a group $G\leq\Isom(\BB)$ isomorphic to $\Gamma$ such that $\NN_{\BB,G}(R) \asymp_\times f(R)$ for all $R\geq 1$.
\end{itemize}
The latter will be of particular interest in Chapter \ref{sectionGFmeasures}, in which the orbital counting function of a parabolic subgroup of a geometrically finite group is shown to have implications for the geometry of the Patterson--Sullivan measure via the Global Measure Formula (Theorem \ref{theoremglobalmeasure}).
\end{remark}

We conclude this chapter by giving two examples of how the Poincar\'e exponents of infinitely generated parabolic groups behave somewhat erratically.

\begin{example}[A class of infinitely generated parabolic torsion groups]
\label{exampleparabolictorsion}
Let $(b_n)_1^\infty$ be an increasing sequence of positive real numbers, and for each $n\in\N$, let $g_n\in\Isom(\BB)$ be the reflection across the hyperplane $H_n := \{\xx : x_n = b_n\}$. Then $G := \lb g_n:n\in\N\rb$ is a strongly discrete subgroup of $\Isom(\BB)$ consisting of only torsion elements. It follows that its Poincar\'e extension $\what G$ is a strongly discrete parabolic subgroup of $\Isom(\H^\infty)$ with no parabolic element. 
To compute the Poincar\'e exponent of $\what G$, we use \eqref{poincareextensionpoincareexponent}:
\begin{align*}
\Sigma_s(G) = \sum_{g\in G} (1\vee \|g(\0)\|)^{-2s}
&= \sum_{\substack{S\subset\N \\ \text{finite}}} \left(1\vee \left\|\left(\prod_{n\in S} g_n\right)(\0)\right\|\right)^{-2s}\\
&= \sum_{\substack{S\subset\N \\ \text{finite}}} \left(1\vee \sum_{n\in S} (2b_n)^2\right)^{-s}.
\end{align*}
The special case $b_n = n$ gives
\[
\Sigma_s(G) \geq \sum_{S\subset \{1,\ldots,N\}} \left(\sum_{n = 1}^N (2n)^2\right)^{-s} \asymp_\times 2^N N^{-3s} \tendsto N \infty \all s \geq 0
\]
and thus $\delta = \infty$, while the special case $b_n = n^n$ gives
\[
\Sigma_s(G) \leq \sum_{n = 1}^\infty \sum_{\substack{S\subset\N \\ \max(S) = n}} (n^n)^{-2s}
= \sum_{n = 1}^\infty 2^{n - 1} (n^n)^{-2s} < \infty \all s > 0
\]
and thus $\delta = 0$. Intermediate values of $\delta$ can be achieved either by setting $b_n = 2^{n/(2\delta)}$, which gives a group of divergence type:
\begin{align*}
\Sigma_s(G)
\asymp_\times \sum_{\substack{S\subset\N \\ \text{finite}}} \left(1\vee \max_{n\in S} (2b_n)^2\right)^{-s}
&\asymp_\times \sum_{n = 1}^\infty \sum_{\substack{S\subset\N \\ \max(S) = n}} b_n^{-2s}\\
&= \sum_{n = 1}^\infty 2^{n - 1} 2^{-ns/\delta}
\begin{cases}
= \infty & \text{ for $s \leq \delta$}\\
< \infty & \text{ for $s > \delta$}
\end{cases}
\end{align*}
or by setting $b_n = 2^{n/(2\delta)} n^{1/\delta}$, which gives a group of convergence type:
\[
\Sigma_s(G)
\asymp_\times \sum_{n = 1}^\infty \sum_{\substack{S\subset\N \\ \max(S) = n}} b_n^{-2s}
= \sum_{n = 1}^\infty 2^{n - 1} 2^{-ns/\delta} n^{-2s/\delta}
\begin{cases}
= \infty & \text{ for $s < \delta$}\\
< \infty & \text{ for $s \geq \delta$}
\end{cases}
\]
\end{example}

\begin{remark}
\label{remarkparabolictorsion}
In Example \ref{exampleparabolictorsion}, for each $n$ the hyperplane $H_n$ is a totally geodesic subset of $\E^\infty$ which is invariant under $G$. However, the intersection $\bigcap_n H_n$ is trivial, since no point $\xx\in\bord \E^\infty\butnot\{\infty\}$ can satisfy $x_n = b_n$ for all $n$. In particular, $G$ does not act irreducibly on any nontrivial totally geodesic set $S\subset \bord\H^\infty$.
\end{remark}

\begin{example}[A torsion-free infinitely generated parabolic group with finite Poincar\'e exponent]
\label{exampleQparabolic}
Let $\Gamma = \{n/2^k : n\in\Z, k\geq 0\}$. Then $\Gamma$ is an infinitely generated abelian group. For each $k\in\N$ let $B_k = k^k$, and define an action $\Phi:\Gamma\to\Isom(\ell^2(\N;\C))$ by the following formula:
\[
\Phi(q)(x_0,\xx) = \big(x_0 + q, \big(e^{2\pi i 2^k q}(x_k - B_k) + B_k\big)_k\big),
\]
i.e. $\Phi(q)$ is the direct sum of the Edelstein-type example (cf. Definition \ref{definitionedelstein}) defined by the sequences $a_k = 2^k q$, $b_k = B_k$ with the map $\R\ni x_0\mapsto x_0 + q$. It is readily verified that $\Phi$ is a homomorphism (cf. \eqref{homomorphism}). We have
\begin{align*}
\|\Phi(q)(\0)\|^2 &= |q|^2 + \sum_k B_k^2 |e^{2\pi i 2^k q} - 1|^2\\
&\asymp_\times |q|^2 + \sum_k B_k^2 \dist(2^k q,\Z)\\
&\asymp_\times \max(|q|^2, B_{k_q}^2),
\end{align*}
where $k_q$ is the largest integer such that $2^{k_q} q\notin\Z$. Equivalently, $k_q$ is the unique integer such that $q = n/2^{k_q + 1}$ for some $k$.

To compute the Poincar\'e exponent of $G = \Phi(\Gamma)$, fix $s > 1/2$ and observe that
\begin{align*}
\Sigma_s(G) &=_\pt \sum_{g\in G} (1\vee \|g(\0)\|)^{-2s}\\
&=_\pt \sum_{q\in \Gamma} (|q|\vee B_{k_q})^{-2s}\\
&\leq_\pt \sum_{k\in\N} \sum_{n\in\Z} (|n|/2^{k + 1} \vee B_k)^{-2s}\\
&\asymp_\times \sum_{k\in\N} \int_0^\infty \left(\frac{x}{2^{k + 1}} \vee B_k\right)^{-2s} \;\dee x\\
&=_\pt \sum_{k\in\N} \left[\int_0^{2^{k + 1} B_k} B_k^{-2s} \;\dee x + \int_{2^{k + 1} B_k}^\infty \left(\frac{x}{2^{k + 1}}\right)^{-2s} \;\dee x\right]\\
&=_\pt \sum_{k\in\N} \left[2^{k + 1} B_k^{1 - 2s} + \left.\left((2^{k + 1})^{2s}\frac{x^{1 - 2s}}{1 - 2s}\right)\right\vert_{x = 2^{k + 1} B_k}^\infty\right]\\
&=_\pt \sum_{k\in\N} \left[2^{k + 1} B_k^{1 - 2s} + \frac{1}{2s - 1} 2^{k + 1} B_k^{1 - 2s}\right]\\
&\asymp_\times \sum_{k\in\N} 2^k B_k^{1 - 2s}
= \sum_{k\in\N} 2^k (k^k)^{1 - 2s} < \infty.
\end{align*}
Thus $\delta(G) \leq 1/2$, but Theorem \ref{theoremparaboliclowerbound} guarantees that $\delta(G) \geq \delta(\Phi(\Z)) \geq 1/2$. So $\delta(G) = 1/2$.
\end{example}
%
%

\chapter{Geometrically finite and convex-cobounded groups}\label{sectionGF}

In this chapter we generalize the notion of geometrically finite groups to regularly geodesic strongly hyperbolic metric spaces, mainly CAT(-1) spaces. We generalize finite-dimensional theorems such as the Beardon--Maskit theorem \cite{BeardonMaskit} 
and Tukia's isomorphism theorem \cite[Theorem 3.3]{Tukia2}.

\begin{standingassumptions}
Throughout this chapter, we assume that
\begin{itemize}
\item[(I)] $X$ is regularly geodesic and strongly hyperbolic, and that
\item[(II)] $G\leq\Isom(X)$ is strongly discrete.
\end{itemize}
Recall that for $x,y\in\bord X$, $\geo xy$ denotes the geodesic segment, ray, or line connecting $x$ and $y$.
\end{standingassumptions}
Note that we do not assume that $G$ is nonelementary.

\bigskip
\section{Some geometric shapes}
\label{subsectionGFgeometry}

To define geometrically finite groups requires three geometric concepts. The first, the quasiconvex core $\CC_\zero$ of the group $G$, has already been introduced in Section \ref{subsectionconvexhulls}. The remaining two concepts are horoballs and Dirichlet domains.

\subsection{Horoballs}
\label{subsubsectionhoroballs}

\begin{figure}
\begin{center}
\begin{tabular}{@{}ll@{}}
\begin{tikzpicture}[line cap=round,line join=round,>=triangle 45,scale=0.9]
\clip(-3.16,-3.06) rectangle (3.45,3.06);
\draw(0.0,0.0) circle (3.0cm);
\draw [fill=black,pattern=north east lines,pattern color=black] (2.0,0.0) circle (1.0cm);
\draw (0.03714500501000815,0.06519516118026925)-- (0.9681603188456883,0.060396113170806974);
\draw (0.03714500501000815,0.06519516118026925)-- (0.09979385703413489,0.11506179560495189);
\draw (0.03714500501000815,0.06519516118026925)-- (0.10974098801916482,0.018077268500909828);
\draw (0.9681603188456883,0.060396113170806974)-- (0.9204321632991048,0.10511466461992194);
\draw (0.9681603188456883,0.060396113170806974)-- (0.9204321632991048,0.023050833993424807);
\draw (0,0)-- (1,0);
\begin{scriptsize}
\draw[color=black] (0.7,1.0) node {$H_{\xi,t}$};
\draw[thin,->,>=stealth] (0.9,0.8) -- (1.6,0.1) ;
\draw [fill=black] (3.0,0.0) circle (.75pt);
\draw[color=black] (3.2,0.13) node {$\xi$};
\draw [fill=black] (0.0,-0.0) circle (.75pt);
\draw[color=black] (-0.1,-0.1) node {$\zero$};
\draw[color=black] (0.4469212289714138,0.1875703566282652) node {$t$};
\end{scriptsize}
\end{tikzpicture}

\begin{tikzpicture}[line cap=round,line join=round,>=triangle 45,scale=0.85]
\clip(-3.68,-0.75) rectangle (3.85,4.36);
\draw (-3.0,0.0)-- (3.0,0.0);
\draw (-3.0,1.5)-- (3.0,1.5);
\draw (-2.958916691465665,2.5399933831077948)-- (-1.646103861448389,3.7681086111884694);
\draw (-2.7895214875924683,1.8624125676150085)-- (-1.1590926503129475,3.4081438029579267);
\draw (-1.8366734658057353,1.8624125676150085)-- (-0.16389582755791549,3.4928414048945253);
\draw (-0.8626510435348529,1.8200637666467094)-- (0.7677777937446677,3.450492603926226);
\draw (0.09019697825187993,1.7777149656784101)-- (1.7418002160155501,3.4716670044103757);
\draw (1.1065682014910616,1.8200637666467094)-- (2.6099506358656845,3.3869694024737775);
\draw (2.038241822793645,1.735366164710111)-- (3.5416242571682677,3.3234462010213286);
\draw (0.0,1.5)-- (0.0,1.0);
\begin{scriptsize}
\draw [fill=black] (0.0,1.0) circle (.75pt);
\draw[color=black] (0.2942550361666618,0.9810500926903516) node {$\zero$};
\draw [fill=black] (0.0,4.0) circle (.75pt);
\draw[color=black] (0.1139837141782084,4.225933888482493) node {$\xi$};
\draw[color=black] (-0.9512559157535614,2.177396138613717) node {$H_{\xi,t}$};
\draw[color=black] (-0.16127704379135044,1.2449331126798539) node {$t$};
\end{scriptsize}
\end{tikzpicture}
\end{tabular}
\end{center}
\caption[Visualizing horoballs in the ball and half-space models]{Two pictures of the same horoball, in the ball model and half-space model, respectively.}
\label{figurehoroball}
\end{figure}
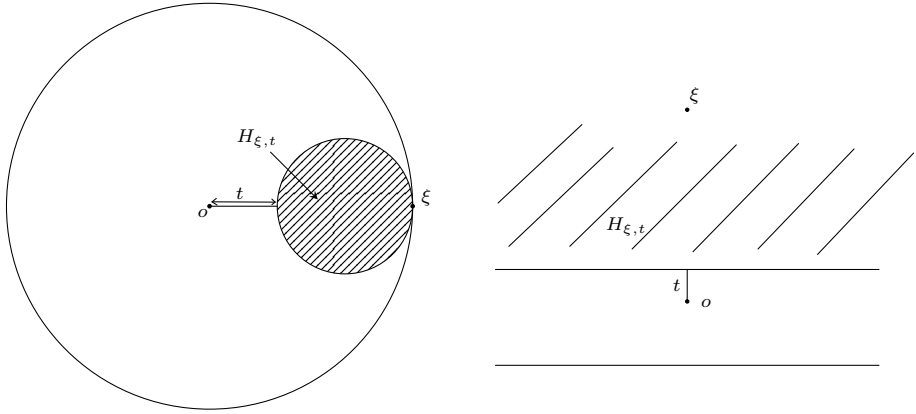

\begin{definition}
\label{definitionhoroball}
A \emph{horoball} is a set of the form
\[
H_{\xi,t} = \{x\in X : \busemann_\xi (\zero, x) > t\},
\]
where $\xi\in\del X$ and $t\in\R$. The point $\xi$ is called the \emph{center} of a horoball $H_{\xi,t}$, and will be denoted $\Center(H_{\xi,t})$. Note that for any horoball $H$, we have
\[
\cl H \cap \del X = \{\Center(H)\}.
\]
(Cf. Figure \ref{figurehoroball}.)
\end{definition}

\begin{lemma}\label{lemmaHdiameter}
For every horoball $H\subset X$, we have
\[
\Diam(H) \asymp_\times b^{-\dist(\zero,H)}.
\]
\end{lemma}
\begin{proof}
Write $H = H_{\xi,t}$ for some $\xi\in\del X$, $t\in\R$. If $t < 0$, then $\zero\in H$, so $\dist(\zero,H) = 0$ and $\Diam(H) = 1$. So suppose $t\geq 0$. Then the intersection $\geo\zero\xi\cap \del H$ consists of a single point $x_0$ satisfying $\dox{x_0} = t$. It follows that $\dist(\zero,H) \leq \dox{x_0} = t$ and $\Diam(H) \geq \Dist(x_0,x_0) = b^{-t}$. For the reverse directions, fix $x\in H$. Since $\busemann_\xi(\zero, x) > t$, we have
\[
\dox x > t
\]
and
\[
\Dist(x,\xi) = b^{-\lb x|\xi\rb_\zero} = b^{-[\busemann_\xi(\zero,x) + \lb\zero|\xi\rb_x]}
\leq b^{-\busemann_\xi(\zero,x)}
< b^{-t}.
\]
It follows that $\Diam(H) \asymp_\times b^{-t} = b^{-\dist(\zero,H)}$.
\end{proof}

\begin{lemma}[Cf. Figure \ref{figureDSC05176disc}]
\label{lemmaHcutBdiameter}
Suppose that $H$ is a horoball not containing $\zero$. Then
\[
\Diam(H\butnot B(\zero,\rho)) \leq 2e^{-(1/2)\rho}.
\]
\end{lemma}
\begin{proof}
Write $H = H_{\xi,t}$ for some $\xi\in\del X$ and $t\in\R$; we have $t\geq 0$ since $\zero\notin H$. Then for all $x\in H\butnot B(\zero,\rho)$,
\[
\lb x|\xi\rb_\zero = \frac 12[\dox x + \busemann_\xi(\zero,x)] \geq \frac 12[\rho + t] \geq \frac 12 \rho
\]
and so $\Dist(x,\xi) \leq e^{-(1/2)\rho}$.
\end{proof}

\begin{figure}
\begin{center}
\begin{tikzpicture}[line cap=round,line join=round,>=triangle 45,x=1.0cm,y=1.0cm]
\clip(-3.98,-3.06) rectangle (4.44,3.06);
\draw(0.0,0.0) circle (3.0cm);
\draw(2.0,0.0) circle (1.0cm);
\draw [shift={(0.0,-0.0)}] plot[domain=-0.49916630796940265:0.5863917531910169,variable=\t]({1.0*2.6608551854448765*cos(\t r)+-0.0*2.6608551854448765*sin(\t r)},{0.0*2.6608551854448765*cos(\t r)+1.0*2.6608551854448765*sin(\t r)});
\draw (0.0,-0.0)-- (2.2141953354369233,1.4709837919838016);
\draw (2.725648009077071,0.6128495628882283)-- (2.6440345586902594,0.5176338707702813);
\draw (2.8299318623491088,0.48136122615392063)-- (2.6984435256148007,0.3226684059573425);
\draw (2.8979430710047853,0.3226684059573425)-- (2.716579847922981,0.10956661883622314);
\draw (2.9432838767752365,0.11863477999031333)-- (2.702977606191846,-0.12167149059307653);
\draw (2.9523520379293267,-0.12167149059307653)-- (2.6848412838836655,-0.3665118417535115);
\draw (2.879806748696605,-0.4118526475239624)-- (2.6485686392673045,-0.5932158706057661);
\draw (2.6032278334968533,-0.7428405296482541)-- (2.712045767345936,-0.6657611598384875);
\begin{scriptsize}
\draw[color=black] (1.3244253846810095,-0.4609545880123348) node {$H$};
\draw [fill=black] (0.0,-0.0) circle (.75pt);
\draw[color=black] (-0.12690497269437848,-0.15753256418730905) node {$\zero$};
\draw[color=black] (1.1175813370309577,0.8973467018588107) node {$\rho$};
\draw [fill=black] (3.0,0.0) circle (.75pt);
\draw[color=black] (3.223917845298176,0.14576257479200725) node {$\xi$};
\draw[color=black] (3.7,0.75) node {$H\butnot B(\zero,\rho)$};
\draw[color=black] (2.87,0.63) node {$\swarrow$};
\end{scriptsize}
\end{tikzpicture}
\caption[Diameter decay of a ball complement inside a horoball]{The set $H\butnot B(\zero,\rho)$ decreases in diameter as $\rho\to\infty$.}
\label{figureDSC05176disc}
\end{center}
\end{figure}

\subsection{Dirichlet domains}

\begin{definition}
\label{definitiondirichletdomain}
Let $G$ be a group acting by isometries on a metric space $X$. Fix $z\in X$. We define the \emph{Dirichlet domain for $G$ centered at $z$} by
\begin{equation}
\label{dirichletdomain}
\DD_z := \{ x : \dist(z,x) \leq \dist(z,g(x)) \all g\in G\} = \{x : \busemann_x(z,g^{-1}(z))\leq 0 \all g\in G\}.
\end{equation}
\end{definition}

The idea is that the Dirichlet domain is a ``tile'' whose iterates under $G$ tile the space $X$. This is made explicit in the following proposition:

\begin{proposition}
\label{propositiondirichletdomain}
For all $z\in X$, $G(\DD_z) = X$.
\end{proposition}
\begin{proof} 
Fix $x\in X$. Since the group $G$ is strongly discrete, the minimum $\min_{g\in G}\{ \dist(x,g(z)) \}$ is attained at some $g\in G$. Now for every $h\in G$, we have $\dist(x,g(z)) \leq \dist(x,h(z))$. Replacing $h$ by $gh$, it follows that for every $h\in G$ we have $\dist(x,g(z)) \leq \dist(x,gh(z))$ which is the same as $\dist(g^{-1}(x),z) \leq \dist(g^{-1}(x),h(z))$. Thus $g^{-1}(x)\in\DD_z$, i.e. $x\in g(\DD_z)$.
\end{proof}

\begin{corollary}
\label{corollarydirichletdomain}
Let $S\subset X$ be a $G$-invariant set. The following are equivalent:
\begin{itemize}
\item[(A)] There exists a bounded set $S_0\subset X$ such that $S\subset G(S_0)$.
\item[(B)] The set $S\cap\DD_z$ is bounded.
\end{itemize}
\end{corollary}
\begin{proof}[Proof of \text{(A)} \implies \text{(B)}]
Given $x\in S\cap\DD_z$, fix $g\in G$ with $x\in g(S_0)$. Then $\dist(z,x) \leq \dist(z,g^{-1}(x)) \asymp_\plus 0$, i.e. $x$ is in a bounded set.
\end{proof}
\begin{proof}[Proof of \text{(B)} \implies \text{(A)}]
The set $S_0 = S\cap\DD_z$ is such a set. Specifically, given $x\in S$ by Proposition \ref{propositiondirichletdomain} there exists $g\in G$ such that $x\in g(\DD_z)$. Since $S$ is $G$-invariant, $g^{-1}(x)\in S\cap \DD_z = S_0$.
\end{proof}

\begin{figure}
\begin{center}
\begin{tabular}{@{}ll@{}}

\begin{tikzpicture}[line cap=round,line join=round,>=triangle 45,scale=0.45]
\clip(-6.705221014292816,-5.620235521243468) rectangle (6.556270148687806,5.50828853160739);
\draw (0.0,3.0)-- (0.0,-3.0);
\draw (-3.0,0.0)-- (3.0,0.0);
\draw (-1.5,3.0)-- (0.0,3.0);
\draw (0.0,3.0)-- (1.5,3.0);
\draw (0.0,3.0)-- (0.0,4.5);
\draw (-3.0,0.0)-- (-3.0,1.5);
\draw (-3.0,0.0)-- (-4.5,0.0);
\draw (-3.0,0.0)-- (-3.0,-1.5);
\draw (0.0,-3.0)-- (-1.5,-3.0);
\draw (0.0,-3.0)-- (1.5,-3.0);
\draw (0.0,-3.0)-- (0.0,-4.5);
\draw (3.0,0.0)-- (3.0,1.5);
\draw (3.0,0.0)-- (3.0,-1.5);
\draw (0.0,4.5)-- (0.0,5.5);
\draw (0.0,4.5)-- (-1.0,4.5);
\draw (0.0,4.5)-- (1.0,4.5);
\draw (-1.5,3.0)-- (-1.5,4.0);
\draw (-1.5,3.0)-- (-2.5,3.0);
\draw (-1.5,3.0)-- (-1.5,2.0);
\draw (1.5,3.0)-- (1.5,4.0);
\draw (1.5,3.0)-- (2.5,3.0);
\draw (1.5,3.0)-- (1.5,2.0);
\draw (3.0,1.5)-- (3.0,2.5);
\draw (3.0,1.5)-- (2.0,1.5);
\draw (3.0,1.5)-- (4.0,1.5);
\draw (3.0,-1.5)-- (2.0,-1.5);
\draw (3.0,-1.5)-- (4.0,-1.5);
\draw (3.0,-1.5)-- (3.0,-2.5);
\draw (1.5,-3.0)-- (1.5,-2.0);
\draw (1.5,-3.0)-- (2.5,-3.0);
\draw (1.5,-3.0)-- (1.5,-4.0);
\draw (0.0,-4.5)-- (1.0,-4.5);
\draw (0.0,-4.5)-- (0.0,-5.5);
\draw (0.0,-4.5)-- (-1.0,-4.5);
\draw (-1.5,-3.0)-- (-2.5,-3.0);
\draw (-1.5,-3.0)-- (-1.5,-2.0);
\draw (-1.5,-3.0)-- (-1.5,-4.0);
\draw (-3.0,-1.5)-- (-2.0,-1.5);
\draw (-3.0,-1.5)-- (-3.0,-2.5);
\draw (-3.0,-1.5)-- (-4.0,-1.5);
\draw (-4.5,0.0)-- (-4.5,1.0);
\draw (-4.5,0.0)-- (-4.5,-1.0);
\draw (-4.5,0.0)-- (-5.5,0.0);
\draw (-3.0,1.5)-- (-4.0,1.5);
\draw (-3.0,1.5)-- (-3.0,2.5);
\draw (-3.0,1.5)-- (-2.0,1.5);
\draw (4.5,0.0)-- (4.5,1.0);
\draw (4.5,0.0)-- (4.5,-1.0);
\draw (4.5,0.0)-- (5.5,0.0);
\draw (3.0,0.0)-- (4.5,0.0);
\draw [line width=2.0pt] (0.0,0.0)-- (3.0,0.0);
\begin{scriptsize}
\draw [fill=black] (1.5,0) circle (1.0pt);
\draw[color=black] (1.5,0.3) node {$z$};
\draw [fill=black] (0.0,3.0) circle (1.0pt);
\draw[color=black] (0.33,2.63) node {$\gamma_2$};
\draw [fill=black] (0.0,-3.0) circle (1.0pt);
\draw[color=black] (0.57,-3.37) node {$\gamma_2^{-1}$};
\draw [fill=black] (-3.0,0.0) circle (1.0pt);
\draw[color=black] (-2.4,0.4) node {$\gamma_1^{-1}$};
\draw [fill=black] (3.0,0.0) circle (1.0pt);
\draw[color=black] (3.4,0.3) node {$\gamma_1$};
\draw [fill=black] (-1.5,3.0) circle (1.0pt);
\draw [fill=black] (1.5,3.0) circle (1.0pt);
\draw [fill=black] (0.0,4.5) circle (1.0pt);
\draw [fill=black] (-3.0,1.5) circle (1.0pt);
\draw [fill=black] (-4.5,0.0) circle (1.0pt);
\draw [fill=black] (-3.0,-1.5) circle (1.0pt);
\draw [fill=black] (-1.5,-3.0) circle (1.0pt);
\draw [fill=black] (1.5,-3.0) circle (1.0pt);
\draw [fill=black] (0.0,-4.5) circle (1.0pt);
\draw [fill=black] (3.0,1.5) circle (1.0pt);
\draw [fill=black] (3.0,-1.5) circle (1.0pt);
\draw [fill=black] (4.5,0.0) circle (1.0pt);
\draw [fill=black] (0.0,5.4) circle (1.0pt);
\draw [fill=black] (-1.0,4.5) circle (1.0pt);
\draw [fill=black] (1.0,4.5) circle (1.0pt);
\draw [fill=black] (-1.5,4.0) circle (1.0pt);
\draw [fill=black] (-2.5,3.0) circle (1.0pt);
\draw [fill=black] (-1.5,2.0) circle (1.0pt);
\draw [fill=black] (1.5,4.0) circle (1.0pt);
\draw [fill=black] (2.5,3.0) circle (1.0pt);
\draw [fill=black] (1.5,2.0) circle (1.0pt);
\draw [fill=black] (3.0,2.5) circle (1.0pt);
\draw [fill=black] (2.0,1.5) circle (1.0pt);
\draw [fill=black] (4.0,1.5) circle (1.0pt);
\draw [fill=black] (2.0,-1.5) circle (1.0pt);
\draw [fill=black] (4.0,-1.5) circle (1.0pt);
\draw [fill=black] (3.0,-2.5) circle (1.0pt);
\draw [fill=black] (1.5,-2.0) circle (1.0pt);
\draw [fill=black] (2.5,-3.0) circle (1.0pt);
\draw [fill=black] (1.5,-4.0) circle (1.0pt);
\draw [fill=black] (1.0,-4.5) circle (1.0pt);
\draw [fill=black] (0.0,-5.5) circle (1.0pt);
\draw [fill=black] (-1.0,-4.5) circle (1.0pt);
\draw [fill=black] (-2.5,-3.0) circle (1.0pt);
\draw [fill=black] (-1.5,-2.0) circle (1.0pt);
\draw [fill=black] (-1.5,-4.0) circle (1.0pt);
\draw [fill=black] (-2.0,-1.5) circle (1.0pt);
\draw [fill=black] (-3.0,-2.5) circle (1.0pt);
\draw [fill=black] (-4.0,-1.5) circle (1.0pt);
\draw [fill=black] (-4.5,1.0) circle (1.0pt);
\draw [fill=black] (-4.5,-1.0) circle (1.0pt);
\draw [fill=black] (-5.5,0.0) circle (1.0pt);
\draw [fill=black] (-4.0,1.5) circle (1.0pt);
\draw [fill=black] (-3.0,2.5) circle (1.0pt);
\draw [fill=black] (-2.0,1.5) circle (1.0pt);
\draw [fill=black] (4.5,1.0) circle (1.0pt);
\draw [fill=black] (4.5,-1.0) circle (1.0pt);
\draw [fill=black] (5.5,0.0) circle (1.0pt);
\draw [fill=black] (0.0,0.0) circle (1.0pt);
\draw[color=black] (0.27329094384908476,0.3) node {$e$};
\draw[color=black] (1.9339097257857616,-0.48055544635952215) node {$\cl{\DD_z^*}$};
\end{scriptsize}
\end{tikzpicture}

\begin{tikzpicture}[line cap=round,line join=round,>=triangle 45,scale=0.45]
\clip(-6.705221014292816,-5.620235521243468) rectangle (6.556270148687806,5.50828853160739);
\draw (0.0,3.0)-- (0.0,-3.0);
\draw (-3.0,0.0)-- (3.0,0.0);
\draw [line width=2.0pt] (-1.5,3.0)-- (0.0,3.0);
\draw [line width=2.0pt] (0.0,3.0)-- (1.5,3.0);
\draw (0.0,3.0)-- (0.0,4.5);
\draw (-3.0,0.0)-- (-3.0,1.5);
\draw (-3.0,0.0)-- (-4.5,0.0);
\draw (-3.0,0.0)-- (-3.0,-1.5);
\draw [line width=2.0pt] (0.0,-3.0)-- (-1.5,-3.0);
\draw [line width=2.0pt] (0.0,-3.0)-- (1.5,-3.0);
\draw (0.0,-3.0)-- (0.0,-4.5);
\draw (3.0,0.0)-- (3.0,1.5);
\draw (3.0,0.0)-- (3.0,-1.5);
\draw (0.0,4.5)-- (0.0,5.4);
\draw [line width=2.0pt] (0.0,4.5)-- (-1.0,4.5);
\draw [line width=2.0pt] (0.0,4.5)-- (1.0,4.5);
\draw (-1.5,3.0)-- (-1.5,4.0);
\draw [line width=2.0pt] (-1.5,3.0)-- (-2.5,3.0);
\draw (-1.5,3.0)-- (-1.5,2.0);
\draw (1.5,3.0)-- (1.5,4.0);
\draw [line width=2.0pt] (1.5,3.0)-- (2.5,3.0);
\draw (1.5,3.0)-- (1.5,2.0);
\draw (3.0,1.5)-- (3.0,2.5);
\draw [line width=2.0pt] (3.0,1.5)-- (2.0,1.5);
\draw [line width=2.0pt] (3.0,1.5)-- (4.0,1.5);
\draw [line width=2.0pt] (3.0,-1.5)-- (2.0,-1.5);
\draw [line width=2.0pt] (3.0,-1.5)-- (4.0,-1.5);
\draw (3.0,-1.5)-- (3.0,-2.5);
\draw (1.5,-3.0)-- (1.5,-2.0);
\draw [line width=2.0pt] (1.5,-3.0)-- (2.5,-3.0);
\draw (1.5,-3.0)-- (1.5,-4.0);
\draw [line width=2.0pt] (0.0,-4.5)-- (1.0,-4.5);
\draw (0.0,-4.5)-- (0.0,-5.5);
\draw [line width=2.0pt] (0.0,-4.5)-- (-1.0,-4.5);
\draw [line width=2.0pt] (-1.5,-3.0)-- (-2.5,-3.0);
\draw (-1.5,-3.0)-- (-1.5,-2.0);
\draw (-1.5,-3.0)-- (-1.5,-4.0);
\draw [line width=2.0pt] (-3.0,-1.5)-- (-2.0,-1.5);
\draw (-3.0,-1.5)-- (-3.0,-2.5);
\draw [line width=2.0pt] (-3.0,-1.5)-- (-4.0,-1.5);
\draw (-4.5,0.0)-- (-4.5,1.0);
\draw (-4.5,0.0)-- (-4.5,-1.0);
\draw (-4.5,0.0)-- (-5.5,0.0);
\draw [line width=2.0pt] (-3.0,1.5)-- (-4.0,1.5);
\draw (-3.0,1.5)-- (-3.0,2.5);
\draw [line width=2.0pt] (-3.0,1.5)-- (-2.0,1.5);
\draw (4.5,0.0)-- (4.5,1.0);
\draw (4.5,0.0)-- (4.5,-1.0);
\draw (4.5,0.0)-- (5.5,0.0);
\draw (3.0,0.0)-- (4.5,0.0);
\draw [line width=2.0pt] (-5.5,0.0)-- (5.44,0.0);
\begin{scriptsize}
\draw [fill=black] (1.5,0) circle (1.0pt);
\draw[color=black] (1.5,0.3) node {$z$};
\draw [fill=black] (0.0,3.0) circle (1.0pt);
\draw[color=black] (0.33,2.68) node {$\gamma_2$};
\draw [fill=black] (0.0,-3.0) circle (1.0pt);
\draw[color=black] (0.574,-3.471) node {$\gamma_2^{-1}$};
\draw [fill=black] (-3.0,0.0) circle (1.0pt);
\draw[color=black] (-2.4,0.45) node {$\gamma_1^{-1}$};
\draw [fill=black] (3.0,0.0) circle (1.0pt);
\draw[color=black] (3.36,0.3) node {$\gamma_1$};
\draw [fill=black] (-1.5,3.0) circle (1.0pt);
\draw [fill=black] (1.5,3.0) circle (1.0pt);
\draw [fill=black] (0.0,4.5) circle (1.0pt);
\draw [fill=black] (-3.0,1.5) circle (1.0pt);
\draw [fill=black] (-4.5,0.0) circle (1.0pt);
\draw [fill=black] (-3.0,-1.5) circle (1.0pt);
\draw [fill=black] (-1.5,-3.0) circle (1.0pt);
\draw [fill=black] (1.5,-3.0) circle (1.0pt);
\draw [fill=black] (0.0,-4.5) circle (1.0pt);
\draw [fill=black] (3.0,1.5) circle (1.0pt);
\draw [fill=black] (3.0,-1.5) circle (1.0pt);
\draw [fill=black] (4.5,0.0) circle (1.0pt);
\draw [fill=black] (0.0,5.4) circle (1.0pt);
\draw [fill=black] (-1.0,4.5) circle (1.0pt);
\draw [fill=black] (1.0,4.5) circle (1.0pt);
\draw [fill=black] (-1.5,4.0) circle (1.0pt);
\draw [fill=black] (-2.5,3.0) circle (1.0pt);
\draw [fill=black] (-1.5,2.0) circle (1.0pt);
\draw [fill=black] (1.5,4.0) circle (1.0pt);
\draw [fill=black] (2.5,3.0) circle (1.0pt);
\draw [fill=black] (1.5,2.0) circle (1.0pt);
\draw [fill=black] (3.0,2.5) circle (1.0pt);
\draw [fill=black] (2.0,1.5) circle (1.0pt);
\draw [fill=black] (4.0,1.5) circle (1.0pt);
\draw [fill=black] (2.0,-1.5) circle (1.0pt);
\draw [fill=black] (4.0,-1.5) circle (1.0pt);
\draw [fill=black] (3.0,-2.5) circle (1.0pt);
\draw [fill=black] (1.5,-2.0) circle (1.0pt);
\draw [fill=black] (2.5,-3.0) circle (1.0pt);
\draw [fill=black] (1.5,-4.0) circle (1.0pt);
\draw [fill=black] (1.0,-4.5) circle (1.0pt);
\draw [fill=black] (0.0,-5.5) circle (1.0pt);
\draw [fill=black] (-1.0,-4.5) circle (1.0pt);
\draw [fill=black] (-2.5,-3.0) circle (1.0pt);
\draw [fill=black] (-1.5,-2.0) circle (1.0pt);
\draw [fill=black] (-1.5,-4.0) circle (1.0pt);
\draw [fill=black] (-2.0,-1.5) circle (1.0pt);
\draw [fill=black] (-3.0,-2.5) circle (1.0pt);
\draw [fill=black] (-4.0,-1.5) circle (1.0pt);
\draw [fill=black] (-4.5,1.0) circle (1.0pt);
\draw [fill=black] (-4.5,-1.0) circle (1.0pt);
\draw [fill=black] (-5.5,0.0) circle (1.0pt);
\draw [fill=black] (-4.0,1.5) circle (1.0pt);
\draw [fill=black] (-3.0,2.5) circle (1.0pt);
\draw [fill=black] (-2.0,1.5) circle (1.0pt);
\draw [fill=black] (4.5,1.0) circle (1.0pt);
\draw [fill=black] (4.5,-1.0) circle (1.0pt);
\draw [fill=black] (5.44,0.0) circle (1.0pt);
\draw [fill=black] (0.0,0.0) circle (1.0pt);
\draw[color=black] (0.27329094384908476,0.23676420562238503) node {$e$};
\draw[color=black] (1.8339097257857616,-0.5) node {$G(\cl{\DD_z^*})$};
\end{scriptsize}
\end{tikzpicture}
\end{tabular}
\caption[The Cayley graph of $\Gamma = \F_2(\Z) = \lb \gamma_1,\gamma_2\rb$]{The Cayley graph of $\Gamma = \F_2(\Z) = \lb \gamma_1,\gamma_2\rb$. The closure of the naive Dirichlet domain $\DD_z^*$ is the geodesic segment $\cl{\DD_z^*} = \geo e{\gamma_1}$. Its orbit $G(\cl{\DD_z^*})$ is the union of all geodesic segments which appear as horizontal lines in this picture.}
\label{figurefree2}
\end{center}
\end{figure}
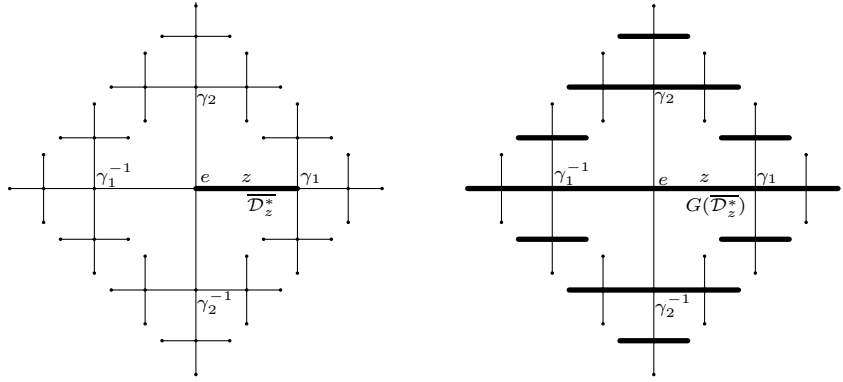

\begin{remark}
It is tempting to define the Dirichlet domain of $G$ centered at $z$ to be the set
\[
\DD_z^* := \{ x : \dist(z,x) < \dist(z,g(x)) \all g\in G \text{ such that }g(z)\neq z\},
\]
and then to try to prove that $G(\cl{\DD_z^*}) = X$. However, there is a simple example which disproves this hypothesis. Let $X$ be the Cayley graph of $\Gamma = \F_2(\Z) = \lb \gamma_1,\gamma_2\rb$, let $\Phi:\Gamma\to\Isom(X)$ be the natural action, and let $G = \Phi(\Gamma)$. If we let $z = ((e,\gamma_1),1/2)$, then $\DD_z^* = \{((e,\gamma_1),t) : t\in (0,1)\}$, and
\[
G(\cl{\DD_z^*}) = \{((\gg,\gg\gamma_1),t) : \gg\in\Gamma,t\in[0,1]\}.
\]
This set excludes all elements of the form $((\gg,\gg\gamma_2),t)$, $t\in (0,1)$. (Cf. Figure \ref{figurefree2}.)
\end{remark}

\begin{remark}
The assumption that $G$ is strongly discrete is crucial for Proposition \ref{propositiondirichletdomain}. In general, tiling Hilbert spaces turns out to be a very subtle problem and has been studied (among others) by Klee \cite{Klee1,Klee2}, Fonf and Lindenstrauss \cite{FonfLindenstrauss} and most recently by Preiss \cite{Preiss}.
\end{remark}

\bigskip
\section{Cobounded and convex-cobounded groups}
\label{subsectionCBCCB}

Before studying geometrically finite groups, we begin by considering the simpler case of cobounded and convex-cobounded groups. The theory of these groups will provide motivation for the theory of geometrically finite groups.

\begin{definition}
\label{definitionCB}
Let $G$ be a group acting by isometries on a metric space $X$. We say that $G$ is \emph{cobounded} if there exists $\sigma > 0$ such that $X = G(B(\zero,\sigma))$.
\end{definition}

It has been a long-standing conjecture to prove or disprove the existence of cobounded subgroups of $\Isom(\H^\infty)$ that are discrete in an appropriate sense. To the best of our knowledge, this conjecture was first stated explicitly by D. P. Sullivan in his IH\'ES seminar on conformal dynamics \cite[p.17]{Sullivan_seminar}. We give here two partial answers to this question, both negative. Our first partial answer is as follows:

\begin{proposition}
\label{propositionnocoboundedgroups}
A strongly discrete subgroup of $\Isom(\H^\infty)$ cannot be cobounded.
\end{proposition}
\begin{proof}
Let us work in the ball model $\B^\infty$. Suppose that $G\leq\Isom(\B^\infty)$ is a strongly discrete cobounded group, and choose $\sigma > 0$ so that $\B^\infty = G(B_\B(\0,\sigma))$. Since $G$ is strongly discrete, we have $\#(F) < \infty$ where
\[
F := \{\xx\in G(\0):\dist_\B(\0,\xx) \leq 2\sigma + 1\}.
\]
Choose $\vv\in \del\B^\infty$ such that $B_\EE(\vv,\zz) = 0$ for all $\zz\in F$, and let $\xx = t\vv$, where $0 < t < 1$ is chosen to make
\[
\dist_\B(\0,\xx) = \sigma + 1.
\]
Since $\xx\in\B$, we have $\xx\in B_\B(\yy,\sigma)$ for some $\yy\in G(\0)$. But then $\dist(\0,\yy)\leq 2\sigma + 1$, which implies $\yy\in F$, and thus $B_\EE(\xx,\yy) = 0$. On the other hand
\[
\dist_\B(\xx,\yy) \leq \sigma < \sigma + 1 = \dist_\B(\0,\xx),
\]
which contradicts \eqref{distB}.
\end{proof}

Proposition \ref{propositionnocoboundedgroups} leaves open the question of whether there exist cobounded subgroups of $\Isom(\H)$ which satisfy a weaker discreteness condition than strong discreteness. One way that we could try to construct such a group would be to take the direct limit of a sequence cobounded subgroups of $\Isom(\H^d)$ as $d\to\infty$. The most promising candidate for such a direct limit has been the direct limit of a sequence of \emph{arithmetic} cocompact subgroups of $\Isom(\H^d)$. (See e.g. \cite{Belolipetsky_survey} for the definition of an arithmetic subgroup of $\Isom(\H^d)$.) Nevertheless, such innocent hopes are dashed by the following result:

\begin{proposition}
\label{propositionarithmeticvolume}
If $G_d\leq\Isom(\H^d)$ is a sequence of arithmetic subgroups, then the codiameter of $G_d$ tends to infinity, that is, there is no $\sigma > 0$ such that $G_d(B(\zero,\sigma)) = \H^d$ for every $d$.
\end{proposition}
\begin{proof}
It is known \cite[Corollary 3.3]{Belolipetsky_survey} that the covolume of $G_d$ tends to infinity superexponentially fast as $d\to\infty$.
On the other hand, the volume of $B(\zero,\sigma)$ in $\H^d$ tends to zero superexponentially fast. Indeed, it is equal to 
\[
(2\pi^{d/2}/\Gamma(d/2))\int_0^\sigma \sinh^{d - 1}(r)\;\dee r \asymp_\times \pi^{d/2} \sigma^{d - 1}/\Gamma(d/2).
\] 
Thus, for sufficiently large $d$, the volume of $B(\zero,\sigma)$ is less than the covolume of $G_d$, which implies that $G_d(B(\zero,\sigma)) \propersubset \H^d$.
\end{proof}
\begin{remark}
Proposition \ref{propositionarithmeticvolume} strongly suggests, but does not prove, that it is impossible to get a cobounded subgroup of $\Isom(\H^\infty)$ as the direct limit of arithmetic subgroups of $\Isom(\H^d)$. One might ask whether one can get a cobounded subgroup of $\Isom(\H^\infty)$ as the direct limit of non-arithmetic subgroups of $\Isom(\H^d)$; the analogous known lower bounds on volume \cite{AdeboyeWei,Kellerhalls} are insufficient to disprove this. However, this approach seems unlikely to work, for two reasons: first of all, the much worse lower bounds for the covolumes of non-arithmetic groups may just be a failure of technique; there are no known examples of non-arithmetic groups with volume lower than the bound which holds for arithmetic groups, and it is conjectured that there are no such examples \cite[p.9]{Belolipetsky_survey}. 
Second of all, even if such groups exist, they are of no use to the problem unless an entire sequence of groups may be found, each one of which is a subgroup of all its higher dimensional analogues. Such structure exists in the arithmetic case but it is unclear whether or not it will also exist in the non-arithmetic case.
\end{remark}

From Propositions \ref{propositionnocoboundedgroups} and \ref{propositionarithmeticvolume}, we see that the theory of cobounded groups acting on $\H^\infty$ will be rather limited. Consequently we focus on the weaker condition of \emph{convex-coboundedness}.

For the remainder of this chapter, we return to our standing assumption that the group $G$ is strongly discrete.

\begin{definition}
\label{definitionCCB}
We say that $G\leq\Isom(X)$ is \emph{convex-cobounded} if its restriction to the quasiconvex core $\CC_\zero$ is cobounded, or equivalently if there exists $\sigma > 0$ such that
\[
\CC_\zero\subset G(B(\zero,\sigma)).
\]
\end{definition}
We remark that whether or not $G$ is convex-cobounded is independent of the base point $\zero$ (cf. Proposition \ref{propositionC0comparison}).

From Proposition \ref{propositionhull} we immediately deduce the following:
\begin{observation}
\label{observationCCBwCCB}
If $X$ is an algebraic hyperbolic space and if $G$ is nonelementary, then the following are equivalent:
\begin{itemize}
\item[(A)] $G$ is convex-cobounded.
\item[(B)] There exists $\sigma > 0$ such that $\CC_\Lambda\subset G(B(\zero,\sigma))$.
\end{itemize}
In particular, when $X$ is finite-dimensional, we see that the notion of convex-coboundedness coincides with the standard notion of convex-cocompactness.
\end{observation}

\subsection{Characterizations of convex-coboundedness}
The property of convex-coboundedness can be characterized in terms of the limit set. Precisely:

\begin{theorem}\label{theoremCCBcompact}
The following are equivalent:
\begin{itemize}
\item[(A)] $G$ is convex-cobounded.
\item[(B)] $G$ is of compact type and any of the following hold:
\begin{itemize}
\item[(B1)] $\Lambda(G) = \Lursigma(G)$ for some $\sigma > 0$.
\item[(B2)] $\Lambda(G) = \Lur(G)$.
\item[(B3)] $\Lambda(G) = \Lr(G)$.
\item[(B4)] $\Lambda(G) = \Lh(G)$.
\end{itemize}
\end{itemize}
\end{theorem}
\begin{remark}
\label{remarkanyBn}
(B1)-(B4) should be regarded as equivalent conditions which also assume that $G$ is of compact type, so that there are a total of $5$ equivalent conditions in this theorem.
\end{remark}
The implications (B1) \implies (B2) \implies (B3) \implies (B4) follow immediately from the definitions. We therefore proceed to prove (A) \implies (B1) and (B4) \implies (A).

\begin{proof}[Proof of \text{(A) \implies (B1)}]
The proof consists of two parts: showing that $\Lambda(G) = \Lursigma(G)$ for some $\sigma > 0$, and showing that $\Lambda(G)$ is compact.

\begin{subproof}[Proof that $\Lambda(G) = \Lursigma(G)$ for some $\sigma > 0$]
Fix $\xi\in\Lambda(G)$, so that 
\[
\geo\zero\xi \subset \CC_\zero \subset G(B(\zero,\sigma)).
\] 
For each $n\in\N$, let $x_n = \geo\zero\xi_n$, so that $x_n\to\xi$ and $\dist(x_n,x_{n + 1}) = 1$. Then for each $n$, there exists $g_n\in G$ satisfying $\dist(g_n(\zero),x_n) \leq \sigma$. Then
\[
\lb \zero|\xi\rb_{g_n(\zero)} \leq \lb \zero|\xi\rb_{x_n} + \sigma = \sigma;
\]
moreover,
\[
\dist(g_n(\zero),g_{n + 1}(\zero)) \leq \dist(x_n,x_{n + 1}) + 2\sigma = 2\sigma + 1.
\]
Thus the convergence $g_n(\zero)\to \xi$ is $(2\sigma + 1)$-uniformly radial, so $\xi\in\Lambda_{\text{ur},2\sigma + 1}(G)$.
\end{subproof}
\begin{subproof}[Proof that $G$ is of compact type]
By contradiction, suppose that $G$ is not of compact type. Then $\Lambda$ is a complete metric space which is not compact, which implies that there exist $\epsilon > 0$ and an infinite $\epsilon$-separated set $I\subset\Lambda$. Fix $\rho > 0$ large to be determined. For each $\xi\in I$, let $x_\xi = \geo\zero\xi_\rho$. Then $x_\xi\in\CC_\zero \subset G(B(\zero,\sigma))$, so there exists $g_\xi\in G$ such that $\dist(g_\xi(\zero),x_\xi)\leq \sigma$.
\begin{claim}
\label{claiminjectiveCCB}
For $\rho$ sufficiently large, the function $\xi \mapsto g_\xi (\zero)$ is injective.
\end{claim}
\begin{subproof}
Fix $\xi_1,\xi_2\in I$ distinct, and suppose $g_{\xi_1}(\zero) = g_{\xi_2}(\zero)$. Then (cf. Figure \ref{figureDSC05177A}) we have that 
\begin{figure}
\begin{center}
\begin{tikzpicture}[line cap=round,line join=round,>=triangle 45,scale=1.2]
\clip(-1.237,-0.3) rectangle(4.17,3.31);
\draw (0.0,-0.0)-- (2.0,3.0);
\draw (0.0,-0.0)-- (3.14,1.94);
\draw [shift={(2.5257374476052528,1.3202338754377372)}]
plot[domain=2.6968619804241416:3.3014726936094725,variable=\t]({1.0*1.3833805462278395*cos(\t r)+-0.0*1.3833805462278395*sin(\t r)},{0.0*1.3833805462278395*cos(\t r)+1.0*1.3833805462278395*sin(\t r)});
\draw [shift={(1.3158446112993265,2.7637714989592004)}]
plot[domain=4.618992018202531:5.148642321848299,variable=\t]({1.0*1.67105450049362*cos(\t r)+-0.0*1.67105450049362*sin(\t r)},{0.0*1.67105450049362*cos(\t r)+1.0*1.67105450049362*sin(\t r)});
\draw [shift={(2.0259885237973716,2.003438282067574)}]
plot[domain=3.258607010908317:4.707025332181733,variable=\t]({1.0*0.7542231048683804*cos(\t r)+-0.0*0.7542231048683804*sin(\t r)},{0.0*0.7542231048683804*cos(\t r)+1.0*0.7542231048683804*sin(\t r)});
\begin{scriptsize}
\draw [fill=black] (0.0,-0.0) circle (.75pt);
\draw[color=black] (-0.11563980272751029,-0.1993420851850086) node {$\zero$};
\draw [fill=black] (2.0,3.0) circle (.75pt);
\draw[color=black] (2.2403553977854944,3.1888033936479783) node {$\xi_1$};
\draw [fill=black] (3.14,1.94) circle (.75pt);
\draw[color=black] (3.373476898984607,2.1342150657993) node {$\xi_2$};
\draw [fill=black] (1.276923076923077,1.9153846153846157) circle (.75pt);
\draw[color=black] (1.0633290204081213,2.010805793391476) node {$x_1$};
\draw [fill=black] (2.021943155793059,1.249226026190616) circle (.75pt);
\draw[color=black] (2.2515744225498424,1.1581599113010557) node {$x_2$};
\draw [fill=black] (1.16,1.1) circle (.75pt);
\draw[color=black] (1.000262673707254,0.8889033169567123) node
{$g_\xi(\zero)$};
\draw[color=black] (1.511842209096808,1.2010935646001886) node {$\sigma$};
\draw[color=black] (1.2516143918733356,1.3601023570593131) node {$\sigma$};
\draw[color=black] (1.7028134319759368,1.6283877285245278) node {$2\sigma$};
\end{scriptsize}
\end{tikzpicture}
\caption[Proving that convex-cobounded groups are of compact type]{If $g_{\xi_1}(\zero) = g_{\xi_2}(\zero)$, then $\xi_1$ and $\xi_2$ must be close to each other.}
\label{figureDSC05177A}
\end{center}
\end{figure}
\[
\lb \xi_1 | \xi_2 \rb_\zero \geq \lb x_1 | x_2 \rb_\zero
= \frac12 [2\rho - \dist(x_1,x_2)]
\geq \rho - \sigma.
\]
On the other hand, since $I$ is $\epsilon$-separated we have $\lb \xi_1|\xi_2\rb_\zero \leq -\log(\epsilon)$. This is a contradiction if $\rho > \sigma - \log(\epsilon)$.
\end{subproof}
\noindent The strong discreteness of $G$ therefore implies
\[
\#(I) \leq \# \{ g\in G : \dogo g \leq \rho + \sigma \} < \infty,
\]
which is a contradiction since $\#(I) = \infty$ by assumption.
\end{subproof}
\noindent This completes the proof of (A) \implies (B1).
\end{proof}
\begin{proof}[Proof of \text{(B4) \implies (A)}]
We use the notation \eqref{primenotation}.
\begin{lemma}
\label{lemmahorosphericaldirichlet}
$\Lh\cap\DD_\zero' = \emptyset$.
\end{lemma}
(Lemma \ref{lemmahorosphericaldirichlet} is true even without assuming (B4); this fact will be used in the proof of Theorem \ref{theoremGFcompact} below.)
\begin{subproof}
By contradiction fix $\xi\in\Lh\cap\DD_\zero'$. Since $\xi\in(\DD_\zero)'$, \eqref{dirichletdomain} gives $\busemann_\xi(\zero,g(\zero)) \leq 0$ for all $g\in G$ (cf. Lemma \ref{lemmanearcontinuity}). But then $\xi \notin \Lh$, since by definition $\xi\in \Lh$ if and only if there exists a sequence $(g_n)_1^\infty$ satisfying $\busemann_\xi(\zero,g_n(\zero)) \to \ +\infty$.
\end{subproof}

Now by (B4) and Observation \ref{observationboundaryofconvexcore}, we have $(\CC_\zero\cap \DD_\zero)' \subset \Lambda\cap\DD_\zero' = \Lh\cap \DD_\zero'$, and so $(\CC_\zero\cap\DD_\zero)' = \emptyset$. By (C) of Proposition \ref{propositioncompacttype}, we get that $\CC_\zero\cap\DD_\zero$ is bounded, and Corollary \ref{corollarydirichletdomain} finishes the proof.
\end{proof}
\noindent The proof of Theorem \ref{theoremCCBcompact} is now complete.

\begin{remark}
(B4) \implies (A) may also be deduced as a consequence of Theorem \ref{theoremGFcompact}(B3)$\Rightarrow$(A) below; cf. Remark \ref{remarkGFCCBcompact}. However, the above prove is much shorter. Alternatively, the above proof may be viewed as the ``skeleton'' of the proof of Theorem \ref{theoremGFcompact}(B3)$\Rightarrow$(A), which is made more complicated by the presence of parabolic points.
\end{remark}

\subsection{Consequences of convex-coboundedness}
The notion of convex-coboundedness also has several important consequences. In the following theorem, $G$ is endowed with an arbitrary Cayley metric (cf. Example \ref{examplecayleygraph}).

\begin{theorem}[Cf. {\cite[Proposition I.8.19]{BridsonHaefliger}}]
\label{theoremCCB}
Suppose that $G$ is convex-cobounded. Then:
\begin{itemize}
\item[(i)] $G$ is finitely generated.
\item[(ii)] The orbit map $g\mapsto g(\zero)$ is a quasi-isometric embedding (cf. Definition \ref{definitionquasiisometry}).
\item[(iii)] 
$\delta_G < \infty$.
\end{itemize}
\end{theorem}
We shall prove Theorem \ref{theoremCCB} as a corollary of a similar statement about geometrically finite groups; cf. Theorem \ref{theoremGF} and Observation \ref{observationGFCCB} below. For now, we list some corollaries of Theorem \ref{theoremCCB}.

\begin{corollary}
Suppose that $G$ is convex-cobounded. Then $G$ is word-hyperbolic, i.e. $G$ is a hyperbolic metric space with respect to any Cayley metric.
\end{corollary}
\begin{proof}
This follows from Theorem \ref{theoremCCB}(ii) and Theorem \ref{theoremquasiisometric}.
\end{proof}

\begin{corollary}
Suppose that $G$ is convex-cobounded. Then 
\[
\HD(\LambdaG) = \delta < \infty.
\]
\end{corollary}
\begin{proof}
This follows from Theorem \ref{theoremCCB}(iii), Theorem \ref{theorembishopjonesregular}, and Theorem \ref{theoremCCBcompact}.
\end{proof}

%
%
%
%
%

\bigskip
\section{Bounded parabolic points}

The difference between groups that are geometrically finite and those that are convex-cobounded is the potential presence of \emph{bounded parabolic points} in the former. 
In the Standard Case, a parabolic fixed point $\xi$ in the limit set of a geometrically finite group $G$, is said to be \emph{bounded} if $(\Lambda\butnot\{\xi\})/\Stab(G;\xi)$ is compact \cite[p.272]{Bowditch_geometrical_finiteness}. We will have to modify this definition a bit to make it work for arbitrary hyperbolic metric spaces, but we show that in the usual case, our definition coincides with the standard one (Remark \ref{remarkLbpvsstandarddef}).

Fix $\xi\in\del X$. Recall that $\EE_\xi$ denotes the set $\bord X\butnot\{\xi\}$.

\begin{definition}
\label{definitionxibounded}
A set $S\subset\EE_\xi$ is \emph{$\xi$-bounded} if $\xi\notin\cl S$.
\end{definition}
The motivation for this definition is that if $X = \H^d$ and $\xi = \infty$, then $\xi$-bounded sets are exactly those which are bounded in the Euclidean metric. Actually, this can be generalized as follows:

\begin{observation}
\label{observationxibounded}
Fix $S\subset\EE_\xi$. The following are equivalent:
\begin{itemize}
\item[(A)] $S$ is $\xi$-bounded.
\item[(B)] $\lb x|\xi\rb_\zero\asymp_\plus 0$ for all $x\in X$.
\item[(C)] $\Dist_\xi(\zero,x)\lesssim_\times 1$ for all $x\in X$.
\item[(D)] $S$ has bounded diameter in the $\Dist_\xi$ metametric.
\end{itemize}
\end{observation}
Condition (D) motivates the terminology ``$\xi$-bounded''.

\begin{proof}[Proof of Observation \ref{observationxibounded}]
(A) \iff (B) follows from the definition of the topology on $\bord X$, (B) \iff (C) follows from \eqref{euclideancomparison}, and (C) \iff (D) is obvious.
\end{proof}

Now fix $G\leq\Isom(X)$, and let $G_\xi$ denote the stabilizer of $\xi$ relative to $G$. Recall (Definition \ref{definitionparabolic}) that $\xi$ is said to be a \emph{parabolic fixed point} of $G$ if $G_\xi$ is a parabolic group, i.e. if $G_\xi(\zero)$ is unbounded and
\[
g\in G_\xi \Rightarrow g'(\xi) = 1.
\]
(Here $g'(\xi)$ is the dynamical derivative of $g$ at $\xi$; cf. Proposition \ref{propositiondynamicalderivative}.)

\begin{observation}
If $\xi$ is a parabolic point then $\xi\in\LambdaG$.
\end{observation}
\begin{proof}
This follows directly from Observation \ref{observationparabolic}.
\end{proof}

\begin{definition}
\label{definitionboundedparabolic}
A parabolic point $\xi\in\LambdaG$ is a \emph{bounded} parabolic point if there exists a $\xi$-bounded set $S\subset\EE_\xi$ such that
\begin{equation}
\label{boundedparabolic}
G(\zero) \subset G_\xi(S).
\end{equation}
We denote the set of bounded parabolic points by $\Lbp$.
\end{definition}

\begin{lemma}
Let $G\leq\Isom(X)$, and fix $\xi\in\del X$. The following are equivalent:
\begin{itemize}
\item[(A)] $\xi$ is a bounded parabolic point.
\item[(B)] All three of the following hold:
\begin{itemize}
\item[(BI)] $\xi\in\LambdaG$, 
\item[(BII)] $g'(\xi) = 1 \all g\in G_\xi$, and
\item[(BIII)] there exists a $\xi$-bounded set $S\subset \EE_\xi$ satisfying \eqref{boundedparabolic}.
\end{itemize}
\end{itemize}
\end{lemma}
\begin{proof}
The only thing to show is that if (B) holds, then $G_\xi(\zero)$ is unbounded. By contradiction suppose otherwise. Let $S$ be a $\xi$-bounded set satisfying \eqref{boundedparabolic}. Then for all $x\in G(\zero)$, we have $x\in h(S)$ for some $h\in G_\xi$, and so
\begin{align*}
\lb x|\xi\rb_\zero = \lb h^{-1}(x)|\xi\rb_{h^{-1}(\zero)} &\asymp_\plus \lb h^{-1}(x)|\xi\rb_\zero \since{$G_\xi(\zero)$ is bounded}\\
&\asymp_\plus 0. \since{$h^{-1}(x)\in S$}
\end{align*}
By Observation \ref{observationxibounded}, the set $G(\zero)$ is $\xi$-bounded and so $\xi\notin\LambdaG$, contradicting (BI).
\end{proof}


We now prove a lemma that summarizes a few geometric properties about bounded parabolic points.

\begin{lemma}
\label{lemmaboundedparabolic}
Let $\xi$ be a parabolic limit point of $G$. The following are equivalent:
\begin{itemize}
\item[(A)] $\xi$ is a bounded parabolic point, i.e. there exists a $\xi$-bounded set $S\subset\EE_\xi$ such that
\begin{equation}
\label{LbpA}
G(\zero)\subset G_\xi(S).
\end{equation}
\item[(B)] There exists a $\xi$-bounded set $S\subset\EE_\xi\cap\del X$ such that
\begin{equation}
\label{LbpB}
\LambdaG\setminus\{\xi\}\subset G_\xi(S).
\end{equation}
\end{itemize}
Moreover, if $H$ is a horoball centered at $\xi$ satisfying $G(\zero)\cap H = \emptyset$, then \text{(A)-(B)} are moreover equivalent to the following:
\begin{itemize}
\item[(C)] There exists a $\xi$-bounded set $S\subset\EE_\xi$ such that
\begin{equation}
\label{LbpC}
\CC_\zero \setminus H \subset G_\xi(S).
\end{equation}
\item[(D)] There exists $\rho > 0$ such that
\begin{equation}
\label{LbpD}
\CC_\zero\cap\del H \subset G_\xi(B(\zero,\rho)).
\end{equation}
\end{itemize}
\end{lemma}

\begin{remark}
\label{remarkLbpvsstandarddef}
The equivalence of conditions (A) and (B) implies that in the Standard Case, our definition of a bounded parabolic point coincides with the usual one.
\end{remark}

\begin{proof}[Proof of \text{(A)} \implies \text{(B)}]
This is immediate since $\LambdaG\butnot\{\xi\} \subset G(\zero)^{(1)_\euc}$. Here $\thicken_{1,\euc}(S)$ denotes the $1$-thickening of $S$ with respect to the Euclidean metametric $\Dist_\xi$.
\end{proof}
\begin{proof}[Proof of \text{(B)} \implies \text{(A)}]
If $\#(\Lambda) = 1$, then $G = G_\xi$ and there is nothing to prove. Otherwise, let $\eta_1,\eta_2\in\Lambda$ be distinct points.

Let $S$ be as in \eqref{LbpB}. Fix $x = g_x(\zero)\in G$. Since $\lb g_x(\eta_1)|g_x(\eta_2)\rb_{g_x(\zero)} \asymp_\plus 0$, Gromov's inequality implies that there exists $i = 1,2$ such that $\lb g_x(\eta_i)|\xi\rb_x \asymp_\plus 0$. By \eqref{LbpB}, there exists $h_x\in G_\xi$ such that $h_x^{-1}g_x(\eta_i)\in S$. We have
\[
\lb h_x^{-1}g_x(\eta_i)|\xi\rb_\zero \asymp_\plus \lb h_x^{-1}g_x(\eta_i)|\xi\rb_{h_x^{-1}(x)} \asymp_\plus 0.
\]
By Proposition \ref{propositionrips}(i), this means that $\zero$ and $y_x := h_x^{-1}(x)$ are both within a bounded distance of the geodesic line $\geo{h_x^{-1}g_x(\eta_i)}\xi$. Since one of these two points must lie closer to $\xi$ then the other, we have either
\begin{equation}
\label{twocases}
\lb y_x|\xi\rb_\zero \asymp_\plus 0 \text{ or } \lb \zero|\xi\rb_{y_x} \asymp_\plus 0.
\end{equation}
By contradiction, let us suppose that there exists a sequence $x_n\in G(\zero)$ such that $\Dist_\xi(\zero,y_{x_n}) \to \infty$. (If no such sequence exists, then for some $N\in\N$ the set $S = \{y\in X : \Dist_\xi(\zero,y) \leq N\}$ is a $\xi$-bounded set satisfying \eqref{LbpA}.) For $n$ sufficiently large, the first case of \eqref{twocases} cannot hold, so the second case holds. It follows that $y_n := y_{x_n} \to \xi$ radially. So $\xi$ is a radial limit point of $G$. In the remainder of the proof, we show that this yields a contradiction.

By Proposition \ref{propositionrips}(i), for each $n\in\N$ there exists a point $z_n\in \geo\zero\xi$ satisfying
\begin{equation}
\label{dxyx}
\dist(y_n,z_n) \asymp_\plus 0.
\end{equation}
Now let $\rho$ be the implied constant of \eqref{dxyx}, and let $\delta$ be the implied constant of Proposition \ref{propositionrips}(ii). Since $G$ is strongly discrete, $M := \#\{g\in G : \dogo g \leq 2\rho + 2\delta\} < \infty$. Let $F\subset G_\xi$ be a finite set with cardinality strictly greater than $M$. By Proposition \ref{propositionrips}(ii), there exists $t > 0$ such that for all $y\in\geo\zero\xi$ with $y > t$, then $\dist(y,\geo{h(\zero)}\xi) \leq \delta$ for all $h\in F$.

Suppose $z_n > t$. Then for all $h\in F$, we have $\dist(z_n,\geo{h(\zero)}\xi) \leq \delta$. On the other hand, $h(z_n)\in \geo{h(\zero)}\xi$ and $\busemann_\xi(z_n,h(z_n)) = 0$; this implies that $\dist(z_n,h(z_n)) \leq 2\delta$ and thus $\dist(y_n,h(y_n)) \leq 2\rho + 2\delta$. But $y_n = g_n(\zero)$ for some $g_n\in G$, so we have $\dogo{g_n^{-1}hg_n} \leq 2\rho + 2\delta$. But since $\#(F) > M$, this contradicts the definition of $M$.

It follows that $z_n \leq t$. But then $\dox{y_n} \leq \dox{z_n} + \rho \leq t + \rho$, implying that the sequence $(y_n)_1^\infty$ is bounded, a contradiction.
\end{proof}

For the remainder of the proof, we fix a horoball $H = H_{\xi,t} \subset X$ disjoint from $G(\zero)$.

\begin{proof}[Proof of \text{(A)} \implies \text{(C)}]
Let $S$ be as in \eqref{LbpA}. Fix $x\in\CC_\zero \setminus H$. Then there exist $g_1,g_2\in G$ with $x\in\geo{g_1(\zero)}{g_2(\zero)}$. We have $\lb g_1(\zero) | g_2(\zero) \rb_x = 0$, so by Gromov's inequality there exists $i = 1,2$ such that $\lb g_i(\zero) | \xi \rb_x \asymp_\plus 0$. By \eqref{euclideancomparison}, we have $\Dist_{\xi,x}(x,g_i(\zero)) \asymp_\times 1$, and combining with \eqref{derivative3} gives
\[
\Dist_\xi(x,g_i(\zero)) \asymp_\times e^{\busemann_\xi(\zero,x)} \leq e^t \asymp_{\times,H} 1.
\]
Now by \eqref{LbpA}, there exists $h\in G_\xi$ such that $h^{-1}(g_i(\zero))\in S$. Then by Observation \ref{observationuniformlyLipschitz},
\[
\Dist_\xi(\zero,h^{-1}(x)) \leq \Dist_\xi(\zero,h^{-1}(g_i(\zero))) + \Dist_\xi(x,g_i(\zero)) \lesssim_\times 1.
\]
Thus $h^{-1}(x)$ lies in some $\xi$-bounded set which is independent of $x$.
\end{proof}
\begin{proof}[Proof of \text{(C)} \implies \text{(D)}]
Let $S$ be a $\xi$-bounded set satisfying \eqref{LbpC}. Then for all $x\in S\cap\del H$, by (h) of Proposition \ref{propositionbasicidentities} we have
\[
\dox x = 2\underbrace{\lb x|\xi\rb_\zero}_{\asymp_\plus 0 \text{ since }x\in S} - \underbrace{\busemann_\xi(\zero,x)}_{ = t \text{ since }x\in\del H} \asymp_{\plus,H} 0.
\]
Thus $S\cap\del H\subset B(\zero,\rho)$ for sufficiently large $\rho$. Applying $G_\xi$ demonstrates \eqref{LbpD}.
\end{proof}
\begin{proof}[Proof of \text{(D)} \implies \text{(A)}]
Let $\rho$ be as in \eqref{LbpD}, and fix $g\in G$. Since by assumption $G(\zero)\cap H = \emptyset$, we have $\busemann_\xi(\zero,g(\zero)) \leq t$. Let $x = \geo{g(\zero)}{\xi}_{t - \busemann_\xi(\zero,g(\zero))}$, so that $x\in \geo{g(\zero)}{\xi}\cap \del H$ (cf. Figure \ref{figureDSC05151}). By \eqref{LbpD}, there exists $h\in G_\xi$ such that $x\in B(h(\zero),\rho)$. Then
\begin{align*}
\lb h^{-1}g(\zero) | \xi \rb_\zero
= \lb g(\zero)|\xi\rb_{h(\zero)}
&\leq \lb g(\zero) | \xi \rb_x + \dist(h(\zero),x)\\
&= \dist(h(\zero),x) \since{$x\in\geo{g(\zero)}\xi$}\\
&\leq \rho.
\end{align*}
This demonstrates that $g(\zero)\in h(S)$ for some $\xi$-bounded set $S$.
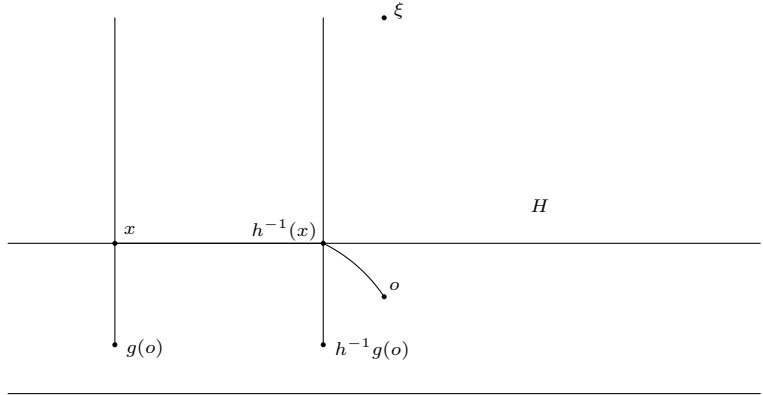
\begin{figure}
\begin{center} 
\begin{tikzpicture}[line cap=round,line join=round,>=triangle 45,x=1.0cm,y=1.0cm]
\clip(-5.74,-0.18) rectangle (6.03,5.39);
\draw (-5.0,0.0)-- (5.0,0.0);
\draw (-5.0,2.0)-- (5.0,2.0);
\draw (-3.58,2.0)-- (-0.81,2.0);
\draw (-0.81,2.0)-- (-0.81,5.0);
\draw (-3.58,2.0)-- (-3.58,5.0);
\draw [shift={(-1.7378227225855478,0.1269240644762515)}] plot[domain=0.5892046179331553:1.1108779656108625,variable=\t]({1.0*2.0902795183382112*cos(\t r)+-0.0*2.0902795183382112*sin(\t r)},{0.0*2.0902795183382112*cos(\t r)+1.0*2.0902795183382112*sin(\t r)});
\draw (-3.58,2.0)-- (-3.58,0.65);
\draw (-0.81,2.0)-- (-0.81,0.65);
\begin{scriptsize}
\draw[color=black] (2.069012341908649,2.5) node {$H$};
\draw [fill=black] (0.0,5.0) circle (.75pt);
\draw[color=black] (0.19765287332143666,5.1) node {$\xi$};
\draw [fill=black] (0.0,1.2884922500123306) circle (.75pt);
\draw[color=black] (0.13525283424696695,1.426597976040203) node {$\zero$};
\draw [fill=black] (-3.58,2.0) circle (.75pt);
\draw[color=black] (-3.3799461107423507,2.177237725963598) node {$x$};
\draw [fill=black] (-0.81,2.0) circle (.75pt);
\draw[color=black] (-1.327946892231745,2.177237725963598) node {$h^{-1}(x)$};
\draw [fill=black] (-3.58,0.65) circle (.75pt);
\draw[color=black] (-3.1799461107423507,0.6) node {$g(\zero)$};
\draw [fill=black] (-0.81,0.65) circle (.75pt);
\draw[color=black] (-0.157946892231745,0.6) node {$h^{-1}g(\zero)$};
\end{scriptsize}
\end{tikzpicture}
\caption[The geometry of bounded parabolic points]{By moving $x$ close to $\zero$ with respect to the $\dist$ metric, $h^{-1}$ also moves $g(\zero)$ close to $\zero$ with respect to the $\Dist_\xi$ metametric.}
\label{figureDSC05151}
\end{center}
\end{figure}
\end{proof}
\begin{remark}
\label{remarkLbpLr}
The proof of (B) \implies (A) given above shows a little more than asked for, namely that a parabolic point of a strongly discrete group cannot also be a radial limit point.
\end{remark}

It will also be useful to rephrase the above equivalent conditions in terms of a Dirichlet domain of $G_\xi$. Indeed, letting $\DD_\zero(G_\xi)$ denote such a Dirichlet domain, we have the following analogue of Corollary \ref{corollarydirichletdomain}:
\begin{lemma}
Let $\xi$ be a parabolic point of $G$, and let $S\subset\EE_\xi$ be a $G_\xi$-invariant set. The following are equivalent:
\begin{itemize}
\item[(A)] There exists a $\xi$-bounded set $S_0\subset\EE_\xi$ such that $S\subset G_\xi(S_0)$.
\item[(B)] The set $S\cap\cl{\DD_\zero(G_\xi)}$ is $\xi$-bounded.
\end{itemize}
\end{lemma}
\begin{proof}
We first observe that for all $x\in \EE_\xi$ and $h\in G_\xi$, (g) of Proposition \ref{propositionbasicidentities} gives
\[
\lb x|\xi\rb_\zero - \lb x|\xi\rb_{h(\zero)} = \frac12 \left[ \busemann_x(\zero,h(\zero)) + \busemann_\xi(\zero,h(\zero)) \right] = \frac12 \busemann_x(\zero,h(\zero)).
\]
In particular
\begin{align*}
x\in \cl{\DD_\zero(G_\xi)}
&\Leftrightarrow \lb x|\xi\rb_\zero \leq \lb x|\xi\rb_{h(\zero)} \all h\in G_\xi\\
&\Leftrightarrow \Dist_\xi(x,\xi) \leq \Dist_\xi(h(x),\xi) \all h\in G_\xi,
\end{align*}
i.e. $\cl{\DD_\zero(G_\xi)}$ is the Dirichlet domain of $\zero$ for the action of $G_\xi$ on the metametric space $(\EE_\xi,\DD_\xi)$. Note that this action is isometric (Observation \ref{observationuniformlyLipschitz}) and strongly discrete (Proposition \ref{propositionSDonclX}). Modifying the proof of Corollary \ref{corollarydirichletdomain} now yields the conclusion.

%

\end{proof}

\begin{corollary}
\label{corollaryboundedparabolic}
In Lemma \ref{lemmaboundedparabolic}, the equivalent conditions \text{(A)-(D)} are also equivalent to:
\begin{itemize}
\item[(A$'$)] $G(\zero)\cap \cl{\DD_\zero(G_\xi)}$ is $\xi$-bounded.
\item[(B$'$)] $\cl{\DD_\zero(G_\xi)}\cap\Lambda\butnot\{\xi\}$ is $\xi$-bounded.
\item[(C$'$)] $\CC_\zero\cap\cl{\DD_\zero(G_\xi)}\butnot H$ is $\xi$-bounded.
\end{itemize}
\end{corollary}

\bigskip
\section{Geometrically finite groups}

\begin{definition}\label{definitionGF}
We say that $G$ is \emph{geometrically finite} if there exists a disjoint $G$-invariant collection of horoballs $\scrH$ satisfying $\zero\notin\bigcup\scrH$ such that
\begin{itemize}
\item[(I)] for every $\rho > 0$, the set
\begin{equation}
\label{Hrhodef}
\scrH_\rho := \{ H\in\scrH : \dist(\zero,H) \leq \rho \}
\end{equation}
is finite, and
\item[(II)] there exists $\sigma > 0$ such that
\begin{equation}
\label{sigmadef}
\CC_\zero \subset G(B(\zero,\sigma)) \cup \bigcup\scrH.
\end{equation}
\end{itemize}
\end{definition}

\begin{observation}
Notice that the following implications hold:
\[
\text{$G$ cobounded} \Rightarrow \text{$G$ convex-cobounded} \Rightarrow \text{$G$ geometrically finite} .
\]
Indeed, $G$ is convex-cobounded if and only if it satisfies Definition \ref{definitionGF} with $\scrH = \emptyset$.
\end{observation}

\begin{remark}
It is not immediately obvious that the definition of geometrical finiteness is independent of the basepoint $\zero$, but this follows from Theorems \ref{theoremGFcompact} and \ref{theoremGF} below.
\end{remark}

\begin{remark}
Geometrical finiteness is closely related to the notion of \emph{relative hyperbolicity} of a group; see e.g. \cite{Bowditch_relatively_hyperbolic}. The main differences are:
\begin{itemize}
\item[1.] Relative hyperbolicity is a property of an abstract group, whereas geometrical finiteness is a property of an isometric group action (equivalently, of a subgroup of an isometry group)
\item[2.] The maximal parabolic subgroups of relatively hyperbolic groups are assumed to be finitely generated, whereas we do not make this assumption (cf. Corollary \ref{corollaryGF}(i)).
\item[3.] The relation between relative hyperbolicity and geometrical finiteness is only available in retrospect, once one proves that both are equivalent to a decomposition of the limit set into radial and bounded parabolic limit points plus auxiliary assumptions (compare Theorem \ref{theoremGFcompact} with \cite[Definition 1]{Bowditch_relatively_hyperbolic}).
\end{itemize}
\end{remark}

\subsection{Characterizations of geometrical finiteness}
We now state and prove an analogue of Theorem \ref{theoremCCBcompact} in the setting of geometrically finite groups. In the Standard Case, the equivalence (A) \iff (B2) of the following theorem was proven by A. F. Beardon and B. Maskit \cite{BeardonMaskit}. Note that while in Theorem \ref{theoremCCBcompact}, one of the equivalent conditions involved the uniformly radial limit set, no such characterization exists for geometrically finite groups. This is because for many geometrically finite groups, the typical point on the limit set is neither parabolic nor uniformly radial. (For example, the set of uniformly radial limit points of the geometrically finite Fuchsian group $\SL_2(\Z)$ is equal to the set of badly approximable numbers; cf. e.g. \cite[Observation 1.15 and Proposition 1.21]{FSU4}.)

\begin{theorem}[Generalization of the Beardon--Maskit Theorem; see also {\cite[Proposition 1.10]{Roblin1}}]
\label{theoremGFcompact}
The following are equivalent:
\begin{itemize}
\item[(A)] $G$ is geometrically finite.
\item[(B)] $G$ is of compact type and any of the following hold (cf. Remark \ref{remarkanyBn}):
\begin{itemize}
\item[(B1)] $\Lambda(G) = \Lrsigma(G) \cup \Lbp(G)$ for some $\sigma > 0$.
\item[(B2)] $\Lambda(G) = \Lr(G) \cup \Lbp(G)$.
\item[(B3)] $\Lambda(G) = \Lh(G) \cup \Lbp(G)$.
\end{itemize}
\end{itemize}
\end{theorem}
\begin{remark}
\label{remarkbowditchcomparision}
Of the equivalent definitions of geometrical finiteness discussed in \cite{Bowditch_geometrical_finiteness}, it seems the above definitions most closely correspond with (GF1) and (GF2).\Footnote{Cf. Remark \ref{remarkLbpvsstandarddef} above regarding (GF2).} It seems that definitions (GF3) and (GF5) cannot be generalized to our setting. Indeed, (GF5) depend on the notion of volume, which does not exist in infinite dimensional spaces, while (GF3) already fails in the case of variable curvature; cf. \cite{Bowditch_geometrical_finiteness2}. It seems plausible that a version of (GF4) could be made to work at least in the setting of algebraic hyperbolic spaces, but we do not study the issue at this stage.
\end{remark}

The implications (B1) \implies (B2) \implies (B3) follow immediately from the definitions. We therefore proceed to prove (A) \implies (B1) and then the more difficult (B3) \implies (A).

\begin{proof}[Proof of \text{(A) \implies (B1)}]
The proof consists of two parts: showing that $\Lambda(G) = \Lrsigma(G) \cup \Lbp(G)$ for some $\sigma > 0$, and showing that $G$ is of compact type.

\begin{subproof}[Proof that $\Lambda(G) = \Lrsigma(G) \cup \Lbp(G)$ for some $\sigma > 0$]
Let $\scrH$ be the disjoint $G$-invariant collection of horoballs as defined in Definition \ref{definitionGF}, and let $\sigma > 0$ be large enough so that \eqref{sigmadef} holds. Fix $\xi\in\LambdaG$, and we will show that $\xi\in\Lrsigma\cup\Lbp$. For each $t\geq 0$, recall that $\geo\zero\xi_t$ denotes the unique point on $\geo\zero\xi$ so that $\dist(\zero,\geo\zero\xi_t) = t$; since $\geo\zero\xi_t\in\CC_\zero$, by \eqref{sigmadef} either $\geo\zero\xi_t\in G(B(\zero,\sigma))$ or $\geo\zero\xi_t\in\bigcup\scrH$.

Now if there exists a sequence $t_n\to \infty$ satisfying $\geo\zero\xi_{t_n}\in G(B(\zero,\sigma))$, then $\xi\in\Lrsigma$ (Corollary \ref{corollaryshadowsrelation}). Assume not; then there exists $t_0$ such that $\geo\zero\xi_t\in\bigcup\scrH$ for all $t > t_0$. This in turn implies that the collection
\[
\big\{\{t > t_0: \geo\zero\xi_t\in H\}: H\in\scrH\big\}
\]
is a disjoint open cover of $(t_0,\infty)$. Since $(t_0,\infty)$ is connected, we have $(t_0,\infty) = \{t > t_0: \geo\zero\xi_t\in H\}$ for some $H\in\scrH$, or equivalently
\[
\geo\zero\xi_t\in H\all t > t_0.
\]
Therefore $\xi = \Center(H)$. Now it suffices to show 
\begin{lemma}
\label{lemmacentersLbp}
For every $H\in\scrH$, if $\Center(H)\in\LambdaG$, then $\Center(H)\in\Lbp$. 
\end{lemma}
\begin{subproof}
Let $\xi = \Center(H)$. For every $g\in G_\xi$, we have $g(H)\cap H\neq\emptyset$. Since $\scrH$ is disjoint, this implies $g(H) = H$ and thus $g'(\xi) = 1$. Thus $\xi$ is neutral with respect to every element of $G_\xi$.

We will demonstrate equivalent condition (D) of Lemma \ref{lemmaboundedparabolic}. First of all, we observe that $G(\zero)$ is disjoint from $H$ since $\zero\notin\bigcup\scrH$. Fix $x\in\CC_\zero\cap\del H \subset \CC_\zero\butnot\bigcup\scrH$. Then by \eqref{sigmadef}, we have $x\in g_x(B(\zero,\sigma))$ for some $g_x\in G$. It follows that $g_x^{-1}(x)\in B(\zero,\sigma)$ and so $g_x^{-1}(H) \cap B(\zero,\sigma + \epsilon) \neq \emptyset$ for every $\epsilon >0$. Equivalently, $g_x^{-1}(H)\in\scrH_{\sigma + \epsilon}$, where $\scrH_\rho$ is defined as in \eqref{Hrhodef}. Therefore, by (I) of Definition \ref{definitionGF}, the set
\[
\{g_x^{-1}(H) : x\in \CC_\zero\cap\del H\}
\]
is finite. Let $(g_{x_i}^{-1}(H))_1^n$ be an enumeration of this set. Then for any $x\in\CC_\zero\cap\del H$ there exists $i = 1,\ldots,n$ with $g_x^{-1}(H) = g_{x_i}^{-1}(H)$. Then $g_x g_{x_i}^{-1} (H) = H$ and so $g_x g_{x_i}^{-1}(\xi) = \xi$. Equivalently, $h_x := g_x g_{x_i}^{-1}\in G_\xi$. Thus
\[
\dist(x,G_\xi(\zero)) \leq \dist(h_x(\zero), x) = \dist(g_{x_i}^{-1}(\zero),g_x^{-1}(x)) \leq \dogo{g_{x_i}^{-1}} + \dox{g_x^{-1}(x)} \leq \sigma + \max_{i = 1}^n \dogo{g_{x_i}}.
\]
Letting $\rho = \sigma + \max_{i = 1}^n \dogo{g_{x_i}}$, we have \eqref{LbpD}, which completes the proof.
\end{subproof}
\noindent The identity $\Lambda(G) = \Lrsigma(G) \cup \Lbp(G)$ has been proven.
\end{subproof}

\begin{subproof}[Proof that $G$ is of compact type]
By contradiction, suppose otherwise. Then $\Lambda$ is a complete metric space which is not compact, which implies that there exist $\epsilon > 0$ and an infinite $\epsilon$-separated set $I\subset\Lambda$. Fix $\rho > 0$ large to be determined. For each $\xi\in I$, let $x_\xi = \geo\zero\xi_\rho$. Then $x_\xi\in\CC_\zero \subset G(B(\zero,\sigma)) \cup \bigcup\scrH$, so either
\begin{itemize}
\item[(1)] there exists $g_\xi\in G$ such that $\dist(g_\xi(\zero),x_\xi)\leq \sigma$, or
\item[(2)] there exists $H_\xi\in\scrH$ such that $x_\xi\in H_\xi$.
\end{itemize}

\begin{claim}
\label{claiminjectiveGF}
For $\rho$ sufficiently large, the partial functions $\xi \mapsto g_\xi (\zero)$ and $\xi \mapsto H_\xi$ are injective.
\end{claim}
\begin{figure}
\begin{center} 
\begin{tikzpicture}[line cap=round,line join=round,>=triangle 45,scale=0.7]
\clip(-6.028,-5.1) rectangle (6.161,4.39);
\draw (0.0,-0.0)-- (2.8,3.38);
\draw (0.0,-0.0)-- (4.0,2.0);
\draw(2.53983892263005,1.8532030631106415) circle (1.291991690158774cm);
\draw(-0.09397058823529784,-0.3482352941176511) circle (4.719619070844117cm);
\begin{scriptsize}
\draw [fill=black] (0.0,-0.0) circle (.75pt);
\draw[color=black] (-0.264,-0.12399999999999917) node {$\zero$};
\draw [fill=black] (2.8,3.38) circle (.75pt);
\draw[color=black] (2.978,3.5820000000000017) node {$\xi_1$};
\draw [fill=black] (4.0,2.0) circle (.75pt);
\draw[color=black] (4.258,2.21) node {$\xi_2$};
\draw [fill=black] (1.7594402109590748,2.123895683229169) circle (.75pt);
\draw[color=black] (1.6940000000000002,2.362) node {$x_1$};
\draw [fill=black] (2.464,1.232) circle (.75pt);
\draw[color=black] (2.706,1.0961) node {$x_2$};
\draw[color=black] (2.706,0.296) node {$H_{\xi_i}$};
\end{scriptsize}
\end{tikzpicture}
\caption[Proving that geometrically finite groups are of compact type]{If $H_{\xi_1} = H_{\xi_2}$, then $\xi_1$ and $\xi_2$ must be close to each other.}
\label{figureDSC05177B}
\end{center}
\end{figure}
\begin{subproof}
For the first partial function $\xi\mapsto g_\xi(\zero)$, see Claim \ref{claiminjectiveCCB}. Now fix $\xi_1,\xi_2\in I$ distinct, and suppose that $H_{\xi_1} = H_{\xi_2}$ (cf. Figure \ref{figureDSC05177B}). Then $x_i := x_{\xi_i}\in H_{\xi_i} \setminus B(\zero,\rho)$. By Lemma \ref{lemmaHcutBdiameter}, this implies that
\[
\epsilon \leq \Dist (\xi_1,\xi_2) \leq \Dist (x_1,x_2) \leq 2 e^{-(1/2)\rho}.
\]
For $\rho > 2(\log(2) - \log(\epsilon))$, this is a contradiction. Thus the second partial function $\xi \mapsto H_\xi$ is also injective.
\end{subproof}
\noindent The strong discreteness of $G$ and \eqref{Hrhodef} therefore imply
\[
\#(I) \leq \# \{ H\in\scrH : \dist(\zero,H) \leq \rho \} + \# \{ g\in G : \dogo g \leq \rho + \sigma \} < \infty \ ,
\]
which is a contradiction since $\#(I) = \infty$ by assumption.
\end{subproof}
\noindent This completes the proof of (A) \implies (B1).
\end{proof}
\begin{proof}[Proof of \text{(B3)} \implies \text{(A)}]
Let $F := (\CC_\zero \cap \DD_\zero)'$, where we use the notation \eqref{primenotation}. By Lemma \ref{lemmahorosphericaldirichlet} Observation \ref{observationboundaryofconvexcore}, and our hypothesis (B3), we have
\begin{equation}
\label{FLbp}
F \subset \Lambda \setminus \Lh \subset \Lbp.
\end{equation}
\begin{claim}
\label{claimFfinite}
$\#(F) < \infty$.
\end{claim}
\begin{subproof}
Note that $F$ is compact since $G$ is of compact type and so it is enough to show that $F$ has no accumulation points. By contradiction, suppose there exists $\xi\in F$ such that $\xi\in\cl{F \setminus \{ \xi \}}$. Then by \eqref{FLbp}, $\xi\in\Lbp$, so by (B$'$) of Corollary \ref{corollaryboundedparabolic}, $\cl{\DD_\zero(G_\xi)} \cap \LambdaG \setminus \{ \xi \}$ is $\xi$-bounded. But $F \setminus \{ \xi \} \subset \DD_\zero' \cap \LambdaG \setminus \{ \xi \} \subset \cl{\DD_\zero(G_\xi)} \cap \LambdaG \setminus \{ \xi \}$, contradicting that $\xi\in\cl{F\butnot\{\xi\}}$.
\end{subproof}
Let $P$ be a transversal of the partition of $F$ into $G$-orbits. Fix $t > 0$ large to be determined. For each $\bp\in P$ let
\[
H_\bp = H_{\bp,t} = \{ x : \busemann_\bp(\zero,x) > t \},
\]
and let
\begin{equation}
\label{scrHdef}
\scrH := \{ g(H_\bp) : \bp\in P, g\in G\} \ .
\end{equation}
Clearly, $\scrH$ is a $G$-invariant collection of horoballs. To finish the proof, we need to show that:
\begin{itemize}
\item[(i)] $\zero\notin\bigcup\scrH$.
\item[(ii)] For $t$ sufficiently large, $\scrH$ is a disjoint collection.
\item[(iii)] ((I) of Definition \ref{definitionGF}) For every $\rho > 0$ we have $\#(\scrH_\rho) < \infty$.
\item[(iv)] ((II) of Definition \ref{definitionGF}) There exists $\sigma > 0$ satisfying \eqref{sigmadef}.
\end{itemize}
It turns out that (ii) is the hardest, so we prove it last.
\begin{subproof}[Proof of \text{(i)}]
Fix $g\in G$ and $\bp\in P$. Since $\bp\in P \subset \DD_\zero'$, we have
\[
\busemann_\bp(\zero,g^{-1}(\zero)) \leq 0 < t.
\]
It follows that $g^{-1}(\zero)\notin H_\bp$, or equivalently $\zero\notin g(H_\bp)$.
\end{subproof}
\begin{subproof}[Proof of \text{(iii)}]
Fix $H = g(H_\bp)\in\scrH_\sigma$ for some $\bp\in P$. Consider the point $x_H = \geo\zero{g(\bp)}_{\dist(\zero,H)}\in\del H$, and note that $\dist (\zero, x_H) = \dist(\zero,H) \leq \sigma$. Now $g^{-1}(x_H)\in H_\bp$, so by (D) of Lemma \ref{lemmaboundedparabolic} there exists $h\in G_\bp$ such that
\[
\dist(h(\zero),g^{-1}(x_H)) \asymp_\plus 0.
\]
(Cf. Figure \ref{figureDSC05164}.) Letting $C$ be the implied constant, we have
\[
\dogo{gh} \leq \dist(\zero,x_H) + \dist(x_H,gh(\zero)) \leq \rho + C.
\]
On the other hand, $gh(H_\bp) = g(H_\bp) = H$ since $h\in G_\bp$. Summarizing, we have
\[
\scrH_\rho \subset \{g(H_\bp) : \bp\in P,\;\; \dogo g \leq \rho + C\}.
\]
But this set is finite because $G$ is strongly discrete and because of Claim \ref{claimFfinite}. Thus $\#(\scrH_\rho) < \infty$.
\begin{figure}
\begin{center} 
\begin{tikzpicture}[line cap=round,line join=round,>=triangle 45,x=1.0cm,y=1.0cm]
\clip(-5.746,-0.18) rectangle (6.031,5.39);
\draw (-5.0,0.0)-- (5.0,0.0);
\draw (-5.0,2.0)-- (5.0,2.0);
\draw [shift={(-1.7378227225855478,0.1269240644762515)}] plot[domain=0.5892046179331553:1.1108779656108625,variable=\t]({1.0*2.0902795183382112*cos(\t r)+-0.0*2.0902795183382112*sin(\t r)},{0.0*2.0902795183382112*cos(\t r)+1.0*2.0902795183382112*sin(\t r)});
\begin{scriptsize}
\draw[color=black] (2.069012341908649,2.5) node {$H_\bp$};
\draw [fill=black] (0.0,5.0) circle (.75pt);
\draw[color=black] (0.19765287332143666,5.17) node {$p$};
\draw [fill=black] (0.0,1.2884922500123306) circle (.75pt);
\draw[color=black] (0.13525283424696695,1.426597976040203) node {$\zero$};
\draw [fill=black] (-3.58,2.0) circle (.75pt);
\draw[color=black] (-3.3799461107423507,2.277237725963598) node {$g^{-1}(x_H)$};
\draw [fill=black] (-0.81,2.0) circle (.75pt);
\draw[color=black] (-1.027946892231745,2.277237725963598) node {$h^{-1}g^{-1}(x_H)$};
\end{scriptsize}
\end{tikzpicture}
\caption[Local finiteness of the horoball collection]{Since $g^{-1}(x_H)$ lies on the boundary of the horoball $H_\bp$, an element of $G_\bp$ can move it close to $\zero$.}
\label{figureDSC05164}
\end{center}
\end{figure}
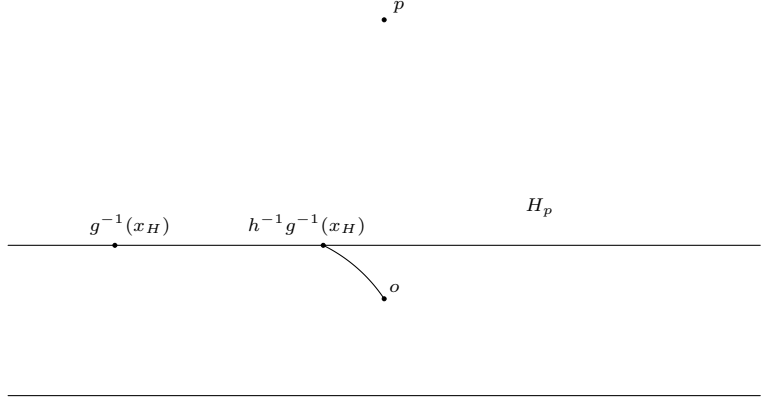
\end{subproof}
\begin{subproof}[Proof of \text{(iv)}]
\begin{claim}
\label{claimC0D0H}
\[
\Bigl( \CC_\zero \cap \DD_\zero \setminus \bigcup\scrH \Bigr)' = \emptyset.
\]
\end{claim}
\begin{subproof}
By contradiction, suppose that there exists
\begin{equation}
\label{CC0DzeroscrH}
\xi\in\Bigl( \CC_\zero \cap \DD_\zero \setminus \bigcup\scrH \Bigr)' \subset F \subset \Lbp .
\end{equation}
By the definition of $P$, there exist $\bp\in P$ and $g\in G$ so that $g(\bp) = \xi$. Then $H_\xi := g(H_\bp)\in\scrH$ is centered at $\xi$, and so by (C$'$) of Corollary \ref{corollaryboundedparabolic}, 
\[
\CC_\zero\cap\DD_\zero\butnot H_\xi \subset \cl{\DD_\zero(G_\xi)}\cap\CC_\zero\butnot H_\xi
\] 
is $\xi$-bounded, contradicting \eqref{CC0DzeroscrH}.
\end{subproof}
Since $G$ is of compact type, Claim \ref{claimC0D0H} implies that the set $\CC_\zero \cap \DD_\zero \setminus \bigcup\scrH$ is bounded (cf. (C) of Proposition \ref{propositioncompacttype}), and Corollary \ref{corollarydirichletdomain} finishes the proof.
\end{subproof}
\begin{subproof}[Proof of \text{(ii)}]
Fix $H_1,H_2\in\scrH$ distinct, and write $H_i = g_i(H_{\xi_i})$ for $i = 1, 2$. The distinctness of $H_1$ and $H_2$ implies that they have different centers, i.e. $g_1(\xi_1) \neq g_2(\xi_2)$. (This is due to the inequivalence of distinct points in $P$.) By contradiction, suppose that $H_1\cap H_2\neq\emptyset$. Without loss of generality, we may suppose that $g_1 = \id$ and that $g_2(\xi_2)\in\cl{\DD_\zero(G_{\xi_1})}$. Otherwise, let $h\in G_{\xi_1}$ be such that $hg_1^{-1}g_2(\xi_2)\in\cl{\DD_\zero(G_{\xi_1})}$ (such an $h$ exists by Proposition \ref{propositiondirichletdomain}), and we have $H_{\xi_1}\cap hg_1^{-1}g_2(H_{\xi_2}) \neq \emptyset$.

By (B$'$) of Corollary \ref{corollaryboundedparabolic}, we have
\[
\lb \xi_1 | g_2(\xi_2)\rb_\zero \asymp_\plus 0,
\]
where the implied constant depends on $\xi_1$. Since there are only finitely many choices for $\xi_1$, we may ignore this dependence.

Fix $x\in H_1\cap H_2$. We have
\begin{align*}
\busemann_{g_2(\xi_2)}(\zero,x) &= \busemann_{\xi_2}(g_2^{-1}(\zero),\zero) + \busemann_{\xi_2}(\zero,g_2^{-1}(x))\\
&> \busemann_{\xi_2}(g_2^{-1}(\zero),\zero) + t \since{$x\in H_2 = g_2(H_{\xi_2})$}\\
&\geq 0 + t. \since{$\xi_2\in \DD_\zero'$}
\end{align*}
On the other hand, $\busemann_{\xi_1}(\zero,x) > t$ since $x\in H_1$. Thus (g) of Proposition \ref{propositionbasicidentities} gives
\begin{align*}
0\leq \lb \xi_1 | g_2(\xi_2)\rb_x &= \lb \xi_1 | g_2(\xi_2)\rb_\zero - \frac12\left[\busemann_{\xi_1}(\zero,x) + \busemann_{g_2(\xi_2)}(\zero,x)\right]\\
&\leq \lb \xi_1 | g_2(\xi_2)\rb_\zero - \frac12\left[t + t\right] \asymp_\plus -t.
\end{align*}
This is a contradiction for sufficiently large $t$.
\end{subproof}
\noindent The implication (B3) \implies (A) has been proven.
\end{proof}
\noindent The proof of Theorem \ref{theoremGFcompact} is now complete.

\begin{remark}
\label{remarkGFCCBcompact}
The implication (B4) \implies (A) of Theorem \ref{theoremCCBcompact} follows directly from the proof of the implication (B3) \implies (A) of Theorem \ref{theoremGFcompact}, since if there are no parabolic points then we have $F = \emptyset$ and so no horoballs will be defined in \eqref{scrHdef}.
\end{remark}

\begin{observation}
\label{observationPfinite}
The proof of Theorem \ref{theoremGFcompact} shows that if $G\leq\Isom(X)$ is geometrically finite, then the set $G\bs\Lbp(G)$ is finite. When $X = \H^3$, this is a special case of Sullivan's Cusp Finiteness Theorem \cite{Sullivan_cusp_finiteness}, which applies to all finitely generated subgroups of $\Isom(\H^3)$ (not just the geometrically finite ones). However, the Cusp Finiteness Theorem does not generalize to higher dimensions \cite{Kapovich3}.
\end{observation}
\begin{proof}
Let $\scrH$ be the collection of horoballs defined in the proof of (B3) \implies (A), i.e. $\scrH = \{g(H_\bp) : \bp\in P\}$ for some finite set $P$. We claim that $\Lbp = G(P)$. Indeed, fix $\xi\in\Lbp$. By the proof of (A) \implies (B1), either $\xi\in\Lr$ or $\xi = \Center(H)$ for some $H\in\scrH$. Since $\Lbp\cap\Lr = \emptyset$ (Remark \ref{remarkLbpLr}), the latter possibility holds. Write $H = g(H_\bp)$; then $\xi = g(\bp)\in G(P)$.
\end{proof}
The set $G\bs\Lbp(G)$ is called the set of \emph{cusps} of $G$.

\begin{definition}
\label{definitioncompleteset}
A \emph{complete set of inequivalent parabolic points} for a geometrically finite group $G$ is a transversal of $G\bs\Lbp(G)$, i.e. a set $P$ such that $\Lbp = G(P)$ but $G(\bp_1)\cap G(\bp_2) = \emptyset$ for all $\bp_1,\bp_2\in P$ distinct.
\end{definition}

Then Observation \ref{observationPfinite} can be interpreted as saying that any complete set of inequivalent parabolic points for a geometrically finite group is finite.

\subsection{Consequences of geometrical finiteness}

Geometrical finiteness, like convex-coboundedness,  has some 
further geometric consequences. Recall (Theorem \ref{theoremCCB}) that if $G$ is convex-cobounded, then $G$ is finitely generated, and for any Cayley graph of $G$, the orbit map $g\mapsto g(\zero)$ is a quasi-isometric embedding. If $G$ is only geometrically finite rather than convex-cobounded, then in general neither of these things is true.\Footnote{For examples of infinitely generated strongly discrete parabolic groups, see Examples \ref{exampleparabolictorsion} and \ref{exampleQparabolic}; these examples can be extended to nonelementary examples by taking a Schottky product with a lineal group. Theorem \ref{theoremparaboliclowerbound} guarantees that the orbit map of a parabolic group is never a quasi-isometric embedding.} Nevertheless, by considering a certain \emph{weighted} Cayley metric with infinitely many generators, we can recover the rough metric structure of the orbit $G(\zero)$.

Recall that the \emph{weighted Cayley metric} of $G$ with respect to a generating set $E_0$ and a weight function $\ell_0:E_0\to(0,\infty)$ is the metric
\[
\dist_G(g_1,g_2) := \inf_{\substack{(h_i)_1^n\in (E\cup F)^n \\ g_1 = g_2 h_1\cdots h_n}} \sum_{i = 1}^n \ell_0(h_i).
\]
(Example \ref{examplecayleygraph}). To describe the generating set and weight function that we want to use, let $P$ be a complete set of inequivalent parabolic points of $G$, and consider the set
\[
E := \bigcup_{\bp\in P} G_\bp.
\]
We will show that there exists a finite set $F$ such that $G$ is generated by $E\cup F$. Without loss of generality, we will assume that this set is symmetric, i.e. $h^{-1}\in F$ for all $h\in F$. For each $h\in E\cup F$ let
\begin{equation}
\label{l0def}
\ell_0(h) := 1\vee\dogo h.
\end{equation}
We then claim that when $G$ is endowed with its weighted Cayley metric with respect to $(E\cup F,\ell_0)$, then the orbit map will be a quasi-isometric embedding. Specifically:

\begin{theorem}
\label{theoremGF}
If $G$ is geometrically finite, then
\begin{itemize}
\item[(i)] There exists a finite set $F$ such that $G$ is generated by $E\cup F$.
\item[(ii)] With the metric $\dist_G$ as above, the orbit map $g\mapsto g(\zero)$ is a quasi-isometric embedding.
\end{itemize}
\end{theorem}

\begin{observation}
\label{observationGFCCB}
Theorem \ref{theoremCCB} follows directly from Theorem \ref{theoremGF}, since by Theorem \ref{theoremCCBcompact} we have $\Lbp = \emptyset$ if $G$ is convex-cobounded.
\end{observation}


We now begin the proof of Theorem \ref{theoremGF}. Of course, part (i) has been proven already (Theorem \ref{theoremGFcompact}).

\begin{proof}[Proof of \text{(i)} and \text{(ii)}]
Let $\scrH$ and $\sigma$ be as in Definition \ref{definitionGF}. Without loss of generality, we may suppose that $\scrH = \{k(H_{\bp,t}) : k\in G,\; \bp\in P\}$ for some $t > 0$ (cf. the proof of Theorem \ref{theoremGFcompact}).

Fix $\rho > 2\sigma + 1$ large to be determined, and let $F = \{g\in G : \dogo g \leq \rho\}$. Then $F$ is finite since $G$ is strongly discrete.
\begin{claim}
\label{claimGcutF}
For all $g\in G\butnot F$, there exist $h_1,h_2\in E\cup F$ such that
\[
\dogo g - \dist(h_1 h_2(\zero),g(\zero)) \gtrsim_{\times,\rho} 1\vee \dogo{h_1}\vee \dogo{h_2} \asymp_\times \ell_0(h_1) + \ell_0(h_2).
\]
\end{claim}
\begin{subproof}
Let $\gamma:[0,\dogo g]\to \geo\zero{g(\zero)}$ be the unit speed parameterization. Let $I = [\sigma + 1,\rho - \sigma]$. Then $\gamma(I)\subset\CC_\zero$, so by \eqref{sigmadef}, either $\gamma(I)\cap h(B(\zero,\sigma))\neq\emptyset$ for some $h\in G$, or $\gamma(I)\subset\bigcup\scrH$.

\begin{itemize}
\item[Case 1:] $\gamma(I) \cap h(B(\zero,\sigma)) \neq \emptyset$ for some $h\in G$. In this case, fix $x\in \gamma(I) \cap h(B(\zero,\sigma))$. Then
\[
\dogo h \leq \dox x + \dist(x,h(\zero)) \leq (\rho - \sigma) + \sigma = \rho,
\]
so $h\in F$. On the other hand,
\begin{align*}
\dist(h(\zero),g(\zero)) &\leq \dist(h(\zero),x) + \dist(x, g(\zero)) \\
&= \dist(h(\zero),x) + \dogo g - \dox x\\
&\leq \sigma + \dogo g - (\sigma + 1)\\
&= \dogo g - 1,
\end{align*}
so
\[
\dogo g - \dist(h(\zero),g(\zero)) \geq 1 \asymp_{\times,\rho} \dogo h.
\]
The claim follows upon letting $h_1 = h$ and $h_2 = \id$.

\item[Case 2:] $\gamma(I) \subset \bigcup\scrH$. In this case, since $\gamma(I)$ is connected and $\scrH$ is a disjoint open cover of $\gamma(I)$, there exists $H\in\scrH$ such that $\gamma(I) \subset H$. Since $\gamma(0),\gamma(\dogo g)\in G(\zero)\subset X\butnot H$, there exist
\[
0 < t_1 < \sigma + 1 < \rho - \sigma < t_2 < \dogo g
\]
so that $\gamma(t_1),\gamma(t_2)\in\del H$. Let $x_i = \gamma(t_i)$ for $i = 1,2$ (cf. Figure \ref{figureDSC05180}).

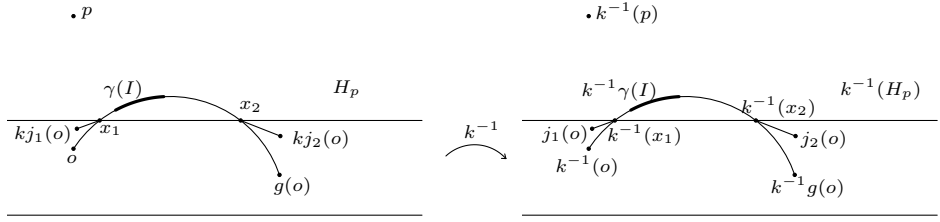
\begin{figure}
\begin{center}
\begin{tabular}{@{}ll|@{}}

\begin{tikzpicture}[line cap=round,line join=round,>=triangle 45,scale=0.84]
\clip(-3.2616777565392376,-0.7167337528936658) rectangle (3.3677361988958543,3.4064383413403285);
\draw (-3.3,0.0)-- (3.3,0.0);
\draw (-3.3,1.5)-- (3.3,1.5);
\draw [shift={(-0.6866012280010922,0.05351828745891047)}] plot[domain=0.32579805555008284:2.563050322940295,variable=\t]({1.0*1.8270110460256037*cos(\t r)+-0.0*1.8270110460256037*sin(\t r)},{0.0*1.8270110460256037*cos(\t r)+1.0*1.8270110460256037*sin(\t r)});
\draw [shift={(-0.6866012280010922,0.05351828745891047)},line width=1.2000000000000002pt] plot[domain=1.62579805555008284:2.063050322940295,variable=\t]({1.0*1.8270110460256037*cos(\t r)+-0.0*1.8270110460256037*sin(\t r)},{0.0*1.8270110460256037*cos(\t r)+1.0*1.8270110460256037*sin(\t r)});
\draw (-2.162165198076832,1.3691062477188254)-- (-1.8026926299856103,1.5);
\draw (1.0576640552230778,1.2529812582555515)-- (0.4294901739834259,1.5);
\begin{scriptsize}
\draw[color=black] (-1.4304721747735372,1.99) node {$\gamma(I)$};
\draw [fill=black] (-2.162165198076832,1.3691062477188254) circle (.75pt);
\draw[color=black] (-2.7,1.28) node {$k j_1 (\zero)$};
\draw [fill=black] (-2.2162855901742153,1.0525353615985757) circle (.75pt);
\draw[color=black] (-2.2430117097284796,0.8840271778089437) node {$\zero$};

\draw [fill=black] (-2.2162855901742153,3.1525353615985757) circle (.75pt);
\draw[color=black] (-2.0030117097284796,3.2) node {$\bp$};

\draw [fill=black] (1.0443009954459461,0.6382805085074923) circle (.75pt);
\draw[color=black] (1.2495575936226908,0.4546727591663085) node {$g(\zero)$};
\draw [fill=black] (-1.8026926299856103,1.5) circle (.75pt);
\draw[color=black] (-1.6304721747735372,1.34895081083073896) node {$x_1$};
\draw [fill=black] (0.4294901739834259,1.5) circle (.75pt);
\draw[color=black] (0.5866161980791817,1.7218467129680482) node {$x_2$};
\draw [fill=black] (1.0576640552230778,1.2529812582555515) circle (.75pt);
\draw[color=black] (1.7,1.2) node {$k j_2 (\zero)$};

\draw[color=black] (2.1057268959530203,2.0) node {$H_\bp$};

\end{scriptsize}
\end{tikzpicture}

\begin{tikzpicture}[line cap=round,line join=round,>=triangle 45,scale=0.33]
\clip(-1.594946279998675,-2.975509096397911) rectangle (1.5782665910292664,2.720001184934291);
\draw [shift={(0.0,0.0)}] plot[domain=0.7354397676755055:2.3211211596591195,variable=\t]({1.0*1.6991762710207554*cos(\t r)+-0.0*1.6991762710207554*sin(\t r)},{0.0*1.6991762710207554*cos(\t r)+1.0*1.6991762710207554*sin(\t r)});
\draw (1.212984563827671,1.3292510070809178)-- (1.26,1.14);
\draw (1.0801070470812233,1.1485375843057488)-- (1.26,1.14);
\begin{scriptsize}
\draw[color=black] (0.3,2.3) node {$k^{-1}$};
\end{scriptsize}
\end{tikzpicture}

\begin{tikzpicture}[line cap=round,line join=round,>=triangle 45,scale=0.84]
\clip(-3.2616777565392376,-0.7167337528936658) rectangle (3.3677361988958543,3.4064383413403285);
\draw (-3.3,0.0)-- (3.3,0.0);
\draw (-3.3,1.5)-- (3.3,1.5);
\draw [shift={(-0.6866012280010922,0.05351828745891047)}] plot[domain=0.32579805555008284:2.563050322940295,variable=\t]({1.0*1.8270110460256037*cos(\t r)+-0.0*1.8270110460256037*sin(\t r)},{0.0*1.8270110460256037*cos(\t r)+1.0*1.8270110460256037*sin(\t r)});
\draw [shift={(-0.6866012280010922,0.05351828745891047)},line width=1.2000000000000002pt] plot[domain=1.62579805555008284:2.063050322940295,variable=\t]({1.0*1.8270110460256037*cos(\t r)+-0.0*1.8270110460256037*sin(\t r)},{0.0*1.8270110460256037*cos(\t r)+1.0*1.8270110460256037*sin(\t r)});
\draw (-2.162165198076832,1.3691062477188254)-- (-1.8026926299856103,1.5);
\draw (1.0576640552230778,1.2529812582555515)-- (0.4294901739834259,1.5);
\begin{scriptsize}
\draw[color=black] (-1.7304721747735372,1.99) node {$k^{-1}\gamma(I)$};
\draw [fill=black] (-2.162165198076832,1.3691062477188254) circle (.75pt);
\draw[color=black] (-2.5885354307014453,1.2852755500491547) node {$j_1 (\zero)$};
\draw [fill=black] (-2.2162855901742153,1.0525353615985757) circle (.75pt);
\draw[color=black] (-2.2430117097284796,0.7840271778089437) node {$k^{-1}(\zero)$};

\draw [fill=black] (-2.2162855901742153,3.1525353615985757) circle (.75pt);
\draw[color=black] (-1.6030117097284796,3.2) node {$k^{-1}(\bp)$};

\draw [fill=black] (1.0443009954459461,0.6382805085074923) circle (.75pt);
\draw[color=black] (1.2495575936226908,0.4546727591663085) node {$k^{-1}g(\zero)$};
\draw [fill=black] (-1.8026926299856103,1.5) circle (.75pt);
\draw[color=black] (-1.3304721747735372,1.27895081083073896) node {$k^{-1}(x_1)$};
\draw [fill=black] (0.4294901739834259,1.5) circle (.75pt);
\draw[color=black] (0.7866161980791817,1.7218467129680482) node {$k^{-1}(x_2)$};
\draw [fill=black] (1.0576640552230778,1.2529812582555515) circle (.75pt);
\draw[color=black] (1.5057268959530203,1.2) node {$j_2 (\zero)$};

\draw[color=black] (2.4057268959530203,2.0) node {$k^{-1}(H_\bp)$};
\end{scriptsize}
\end{tikzpicture}

\end{tabular}
\caption[Orbit maps of geometrically finite groups are QI embeddings]{Since $j_1^{-1} j_2\in E$ and $k j_1\in F$, the points $\zero$, $k j_1(\zero)$, and $k j_2(\zero)$ are connected to each other by edges in the weighted Cayley graph. Since the distance from $k j_2(\zero)$ to $g(\zero)$ are both significantly less than the distance from $\zero$ to $g(\zero)$, our recursive algorithm will eventually halt.}
\label{figureDSC05180}
\end{center}
\end{figure}

Since $H\in\scrH$, we have $H = k(H_\bp)$ for some $\bp\in P$ and $k\in G$. By (D) of Lemma \ref{lemmaboundedparabolic}, there exist $j_1, j_2\in G_\bp$ with
\[
\dist(k^{-1}(x_i), j_i(\zero)) \leq \rho_\bp \;\; (i = 1,2)
\] 
for some $\rho_\bp > 0$ depending only on $\bp$. Letting $\rho_0 = \max_{\bp\in P} \rho_\bp$, we have
\[
\dogo{k j_1} \leq \dox{x_1} + \dist(x_1, k j_1(\zero))
\leq (\sigma + 1) + \rho_0.
\]
Letting $\rho = \max(\rho_0 + \sigma + 1,2\sigma + 2)$, we see that $\dogo{k j_1} \leq \rho$, so $h_1 := k j_1\in F$. On the other hand, $h_2 := j_1^{-1} j_2\in E$ by construction, since $j_1, j_2\in G_\bp$. Observe that $h_1 h_2 = k j_2$. Now
\begin{align*}
\dist(h_1 h_2(\zero),g(\zero)) &\leq \dist(g(\zero),x_2) + \dist(x_2, kj_2(\zero)) \\
&\leq (\dogo g - t_2) + \rho_0,
\end{align*}
and so
\begin{equation}
\label{t2rhoxi}
\dogo g - \dist(h_1 h_2(\zero),g(\zero)) \geq t_2 - \rho_0.
\end{equation}
Now
\begin{align*}
t_2 \geq t_2 - t_1 &= \dist(x_1,x_2)\\
&\geq \dist(j_1(\zero),j_2(\zero)) - \dist(k^{-1}(x_1),j_1(\zero)) - \dist(k^{-1}(x_2),j_2(\zero))\\
&\geq \dogo{j_1^{-1}j_2} - 2\rho_0 = \dogo{h_2} - 2\rho_0
\end{align*}
and on the other hand
\[
t_2 \geq \rho - \sigma \geq \rho_0 + 1.
\]
Combining with \eqref{t2rhoxi}, we see that
\begin{align*}
\dogo g - \dist(h_1 h_2(\zero),g(\zero))
&\geq_\pt (\dogo{h_2} - 2\rho_0)\vee(\rho_0 + 1) - \rho_0\\
&=_\pt (\dogo{h_2} - 3\rho_0)\vee 1\\
&\asymp_\times 1\vee\dogo{h_1}\vee\dogo{h_2}.
\end{align*}
\end{itemize}
\end{subproof}
Fix $j\in G$, and define the sequence $(h_i)_1^n$ in $E\cup F$ inductively as follows: If $h_1,\ldots,h_{2i}$ have been defined for some $i\geq 0$, then let
\[
g = g_{2i} = h_{2i}^{-1}\cdots h_1^{-1} j = (h_1\cdots h_{2i})^{-1} j.
\]
(Note that $g_0 = j$.) If $g\in F$, then let $h_{2i + 1} = g$ and let $n = 2i + 1$ (i.e. stop the sequence here). Otherwise, by Claim \ref{claimGcutF} there exist $h_{2i + 1},h_{2i + 2}\in E\cup F$ such that
\begin{equation}
\label{hi1hi2}
\dogo{g_{2i}} - \dist(h_{2i + 1} h_{2i + 2}(\zero),g_{2i}(\zero)) \gtrsim_{\times,\rho} \ell_0(h_{2i + 1}) + \ell_0(h_{2i + 2}).
\end{equation}
This completes the inductive step, as now $h_1,\ldots,h_{2(i + 1)}$ have been defined. We remark that \emph{a priori}, this process could be infinite and so we could have $n = \infty$; however, it will soon be clear that $n$ is always finite.

We observe that \eqref{hi1hi2} may be rewritten:
\[
\dogo{g_{2i}} - \dogo{g_{2(i + 1)}} \gtrsim_{\times,\rho} \ell_0(h_{2i + 1}) + \ell_0(h_{2i + 2}).
\]
Iterating yields
\begin{equation}
\label{h1h2m}
\dogo j - \dogo{g_{2m}} \gtrsim_\times \sum_{i = 1}^{2m} \ell_0(h_i) \all m \leq n/2.
\end{equation}
In particular, since $\ell_0(h_i) \geq 1$ for all $i$, we have
\[
\dogo j \gtrsim_\times \lfloor n/2\rfloor \asymp_\times n,
\]
and thus $n < \infty$. This demonstrates that the sequence $(h_i)_1^n$ is in fact a finite sequence. In particular, since the only way the sequence can terminate is if $g_{2i}\in F$ for some $i\geq 0$, we have $g_{n - 1}\in F$ and $h_n = g_{n - 1}$. From the definition of $g_{n - 1}$, it follows that $j = h_1\cdots h_n$. Since $j$ was arbitrary and $h_1,\ldots,h_n\in E\cup F$, this demonstrates that $E\cup F$ generates $G$, completing the proof of (i).

To demonstrate (ii), we observe that by \eqref{h1h2m} we have
\begin{align*}
\dogo j &\gtrsim_\times \sum_{i = 1}^{n - 1} \ell_0(h_i)\\
&\asymp_\plus \sum_{i = 1}^n \ell_0(h_i) \since{$h_n\in F$}\\
&\geq_\pt \dist_G(\id,j),
\end{align*}
where $\dist_G$ denotes the weighted Cayley metric. Conversely, if $(h_i)_1^n$ is any sequence satisfying $j = h_1\cdots h_n$, then
\[
\dogo j \leq \sum_{i = 1}^n \dist(h_1\cdots h_{i - 1}(\zero), h_1\cdots h_i(\zero)) = \sum_{i = 1}^n \dogo{h_i} \leq \sum_{i = 1}^n \ell_0(h_i),
\]
and taking the infimum gives $\dogo j \leq \dist_G(\id,j)$.
\end{proof}
This finishes the proof of Theorem \ref{theoremGF}.

\begin{corollary}
\label{corollaryGF}
If $G$ is geometrically finite, then
\begin{itemize}
\item[(i)] If for every $\xi\in\Lbp$, $G_\xi$ is finitely generated, then $G$ is finitely generated.
\item[(ii)] If for every $\xi\in\Lbp$, $\delta(G_\xi) < \infty$, then $\delta(G) < \infty$.
\end{itemize}
\end{corollary}
\begin{proof}[Proof of \text{(i)}]
This is immediate from Theorem \ref{theoremGF}(i) and Observation \ref{observationPfinite}.
\end{proof}
\begin{proof}[Proof of \text{(ii)}]
Call a sequence $(h_i)_1^n\in E^n$ \emph{minimal} if
\begin{equation}
\label{quasiisometry}
\sum_{i = 1}^n \ell_0(h_i) = \dist_G(\id,h_1\cdots h_n).
\end{equation}
Then for each $g\in G\butnot\{\id\}$, there exists a minimal sequence $(h_i)_1^n\in (E\cup F)^n$ so that $g = h_1\cdots h_n$.

Let $C$ be the implied multiplicative constant of \eqref{quasiisometry}, so that for every minimal sequence $(h_i)_1^n$, we have
\[
\sum_{i = 1}^n \ell_0(h_i) \gtrsim_\plus \frac1C \dogo{h_1\cdots h_m}.
\]
Fix $s > 0$. Then
\begin{align*}
\Sigma_s(G) - 1 &\leq_\pt \sum_{g\in G\butnot\{\id\}}\sum_{n\in\N}\sum_{\substack{(h_i)_1^n\in (E\cup F)^n \\ \text{minimal} \\ g = h_1\cdots h_n}} e^{-s\dogo g}\\
&=_\pt \sum_{n\in\N}\sum_{\substack{(h_i)_1^n\in (E\cup F)^n \\ \text{minimal}}} e^{-s\dogo g}\\
&\lesssim_\times \sum_{n\in\N}\sum_{\substack{(h_i)_1^n\in (E\cup F)^n \\ \text{minimal}}} \exp\left(-\frac sC\sum_{i = 1}^n \ell_0(h_i)\right)\\
&\leq_\pt \sum_{n\in\N}\sum_{(h_i)_1^n\in (E\cup F)^n} \exp\left(-\frac sC\sum_{i = 1}^n \ell_0(h_i)\right)\\
&=_\pt \sum_{n\in\N}\sum_{(h_i)_1^n\in (E\cup F)^n}\prod_{i = 1}^n e^{-(s/C)\ell_0(h_i)}\\
&=_\pt \sum_{n\in\N}\prod_{i = 1}^n \sum_{h\in E\cup F} e^{-(s/C)\ell_0(h)}\\
&=_\pt \sum_{n\in\N}\left(\sum_{h\in E\cup F} e^{-(s/C)\ell_0(h)}\right)^n.
\end{align*}
In particular, if
\[
\lambda_s := \sum_{h\in E\cup F} e^{-(s/C)\ell_0(h)} < 1,
\]
then $\Sigma_s(G) < \infty$. Now when $s/C > \max_{\bp\in P}\delta(G_\bp)$, we have $\lambda_s < \infty$. On the other hand, each term of the sum defining $\lambda_s$ tends to zero as $s\to\infty$. Thus $\lambda_s\to 0$ as $s\to \infty$, and in particular there exists some value of $s$ for which $\lambda_s < 1$. For this $s$, $\Sigma_s(G) < \infty$ and so $\delta_G \leq s < \infty$.
\end{proof}

\bigskip
\subsection{Examples of geometrically finite groups}

We conclude this section by giving some basic examples of geometrically finite groups. We begin with the following observation:

\begin{observation}~
\label{observationclassificationGF}
\begin{itemize}
\item[(i)] Any elliptic or lineal group is convex-cobounded.
\item[(ii)] Any parabolic group is geometrically finite and is not convex-cobounded.
\end{itemize}
\end{observation}
\begin{proof}
This follows directly from Theorems \ref{theoremCCBcompact} and \ref{theoremGFcompact}. It may also be proven directly; we leave this as an exercise to the reader.
\end{proof}

\begin{proposition}
\label{propositionSproductGF}
The strongly separated Schottky product $G = \lb G_a\rb_{a\in E}$ of a finite collection of geometrically finite groups is geometrically finite. Moreover, if $P_1$ and $P_2$ are complete sets of inequivalent parabolic points for $G_1$ and $G_2$ respectively, then $P = P_1\cup P_2$ is a complete set of inequivalent parabolic points for $G$. In particular, if the groups $(G_a)_{a\in E}$ are convex-cobounded, then $G$ is convex-cobounded.
\end{proposition}
\begin{proof}
This follows direction from Lemma \ref{lemmaschottkyradial}, Theorem \ref{theoremschottkylimitset}, Corollary \ref{corollaryschottkycompacttype}, and Theorem \ref{theoremGFcompact}.
\end{proof}

Combining Observation \ref{observationclassificationGF} and Proposition \ref{propositionSproductGF} yields the following:

\begin{corollary}
\label{corollaryschottkyGF}
The Schottky product of finitely many parabolic and/or lineal groups is geometrically finite. If only lineal groups occur in the product, then it is convex-cobounded.
\end{corollary}

%


\ignore{
\bigskip
\subsection{Tukia's isomorphism theorem}
\label{subsectiontukia}

As an application of Theorem \ref{theoremGF}, we prove Theorem \ref{theoremtukiaintro} from the introduction:

\begin{definition}
\label{definitiontypepreserving}
An isomorphism between two groups acting on hyperbolic metric spaces is \emph{type-preserving} if the image of a loxodromic (resp. parabolic, elliptic) isometry is loxodromic (resp. parabolic, elliptic).
\end{definition}

\begin{definition}
\label{definitionquasisymmetric}
Let $(Z,\Dist)$ and $(\w Z,\w\Dist)$ be metric spaces. A homeomorphism $\phi:Z\to\w Z$ is said to be \emph{quasisymmetric} if there exists an increasing homeomorphism $f:(0,\infty)\to(0,\infty)$ such that
\[
\frac{\w\Dist(\phi(z),\phi(y))}{\w\Dist(\phi(z),\phi(x))} \leq f\left(\frac{\Dist(z,y)}{\Dist(z,x)}\right) \all x,y,z\in Z.
\]
\end{definition}


\begin{theorem}[Generalization of Tukia's isomorphism theorem; cf. Theorem \ref{theoremtukiaintro}]
\label{theoremtukia}
Let $X$, $\w X$ be CAT(-1) spaces (or more generally, regularly geodesic strongly hyperbolic metric spaces), let $G\leq\Isom(X)$ and $\w G\leq\Isom(\w X)$ be two geometrically finite groups, and let $\Phi:G\to \w G$ be a type-preserving isomorphism. Let $P$ be a complete set of inequivalent parabolic points for $G$.
\begin{itemize}
\item[(i)] If for every $\bp\in P$ we have
\begin{equation}
\label{tukia}
\dogo{\Phi(h)} \asymp_{\plus,\times,\bp} \dogo h \all h\in G_\bp,
\end{equation}
then there is an equivariant homeomorphism between $\Lambda := \Lambda(G)$ and $\w\Lambda := \Lambda(\w G)$.
\item[(ii)] If for every $\bp\in P$ there exists $\alpha_\bp > 0$ such that
\begin{equation}
\label{tukia2}
\dogo{\Phi(h)}  \asymp_{\plus,\bp} \alpha_\bp \dogo h \all h\in G_\bp,
\end{equation}
then the homeomorphism of \text{(i)} is quasisymmetric.
\end{itemize}
\end{theorem}

When $X$ and $\w X$ are finite-dimensional real \ROSSs, Theorem \ref{theoremtukia} was proven by P. Tukia \cite[Theorem 3.3]{Tukia2}. Note that in this case, the hypothesis \eqref{tukia2} always holds with $\alpha_\bp = 1$ (Corollary \ref{corollaryrealROSSONCTquasi}; see also \cite[Theorem 5.4.3]{Ratcliffe}). This is why Tukia's original theorem does not need to mention the conditions \eqref{tukia} and \eqref{tukia2}.

A natural question is then whether the assumptions \eqref{tukia} and/or \eqref{tukia2} are really necessary. In the case of finite-dimensional nonreal \ROSSs, we show that \eqref{tukia} holds automatically (Corollary \ref{corollaryCROSSONCTquasi}), and that \eqref{tukia2} holds assuming both that (A) one of the groups $G$, $\w G$ is a lattice, and that (B) the underlying base fields of $X$ and $\w X$ are the same (Corollary \ref{corollarylatticequasi}). Without these assumptions, it is easy to construct examples of groups $G$, $\w G$ satisfying the hypotheses of the theorem but for which the equivariant homeomorphism is not quasisymmetric (Example \ref{exampletukia} and Remark \ref{remarktukia}). This shows that the assumption \eqref{tukia2} cannot be omitted from the second assertion of Theorem \ref{theoremtukia}.

For the remainder of this subsection, the notation will be as in Theorem \ref{theoremtukia}.

Observe that a subgroup of $G$ is parabolic if and only if it is infinite and consists only of parabolic and elliptic elements. Since $\Phi$ is type-preserving, it follows that $\Phi$ preserves the class of parabolic subgroups, and also the class of maximal parabolic subgroups. But all maximal parabolic subgroups of $G$ are of the form $G_\xi$, where $\xi$ is a parabolic fixed point of $G$. It follows that there is a bijection $\phi:\Lbp(G)\to\Lbp(\w G)$ such that $\Phi(G_\xi) = \w G_{\phi(\xi)}$ for all $\xi\in\Lbp(G)$. The equivariance of $\phi$ implies that $\w P := \phi(P)$ is a complete set of inequivalent parabolic points for $\w G$.

Let $\dist_G$ and $\dist_{\w G}$ denote the weighted Cayley metrics on $G$ and $\w G$, respectively.

\begin{lemma}
\label{lemmatukia1}
$\dist_G \asymp_\times \dist_{\w G}\circ\Phi$.
\end{lemma}
\begin{proof}
Let $E$ and $F$ be as in Theorem \ref{theoremGF}, and let $\w E$ and $\w F$ be the corresponding sets for $\w G$. Since $\w P = \phi(P)$, we have $\w E = \Phi(E)$. On the other hand, for all $h\in E$, we have $\ell_0(h) \asymp_\times \ell_0(\Phi(h))$ by \eqref{tukia}. Thus, edges in the weighted Cayley graph of $G$ have roughly the same weight as their corresponding edges in the weighted Cayley graph of $\w G$. (The sets $F$ and $\w F$ are both finite, and so their edges are essentially irrelevant.) The lemma follows.
\end{proof}

Thus, the map $\Phi(g(\zero)) := \Phi(g)(\zero)$ is a quasi-isometry between $G(\zero)$ and $\w G(\zero)$. At this point, we would like to extend $\Phi$ to an equivariant homeomorphism between $\Lambda$ and $\w\Lambda$. However, all known theorems which give such extensions, e.g. \cite[Theorem 6.5]{BonkSchramm}, require the spaces in question to be geodesic or at least roughly geodesic -- for the good reason that the extension theorems are false without this hypothesis\Footnote{A counterexample is given by letting $X_1 = X_2 = \amsbb R$, $\dist_1(x,y) = \log(1 + |y - x|)$, $\dist_2(x,y) = \begin{cases} \dist_1(x,y) & xy\geq 0 \\ \dist_1(0,x) + \dist_1(0,y) & xy\leq 0 \end{cases}$, and $\Phi:X_1\to X_2$ the identity map -- since $\#(\del X_1) = 1 < 2 = \#(\del X_2)$, $\Phi$ cannot be extended to a homeomorphism between $\del X_1$ and $\del X_2$. On the other hand, if one of the spaces in question is geodesic, then the extension theorem can be proven by isometrically embedding the other space into a geodesic hyperbolic metric space via \cite[Theorem 4.1]{BonkSchramm} -- a fact which however has no relevance to the present situation.} -- but the spaces $G(\zero)$ and $\w G(\zero)$ are not roughly geodesic. They are, however, embedded in the roughly geodesic metric spaces $\CC_\zero$ and $\w\CC_\zero$, which suggests the strategy of extending the map $\Phi$ to a quasi-isometry between $\CC_\zero$ and $\w\CC_\zero$. It turns out that this strategy works if we assume \eqref{tukia2}, and thus proves the existence of a \emph{quasisymmetric} equivariant homeomorphism between $\Lambda$ and $\w\Lambda$ in that case. Since we know that the equivariant homeomorphism is not necessarily quasisymmetric if \eqref{tukia2} fails (Example \ref{exampletukia} and Remark \ref{remarktukia}), this strategy can't be used to prove part (i) of Theorem \ref{theoremtukia}. Thus the proof splits into two parts at this point, depending on whether we have the stronger assumption \eqref{tukia2} which guarantees quasisymmetry, or only the weaker assumption \eqref{tukia}.

\subsection{Completion of the proof assuming \eqref{tukia2}}

\begin{lemma}
\label{lemmatukia2}
Fix $\bp\in P$ and let $\w\bp = \phi(\bp)$. Let
\[
A = A(\bp) = \bigcup_{h\in G_\bp} \geo{h(\zero)}{\bp},
\]
and define a bijection $\psi = \psi_\bp:A\to \w A := A(\w\bp)$ by
\[
\psi(\geo{h(\zero)}{\bp}_t) = \geo{\Phi(h)(\zero)}{\w\bp}_{\alpha_\bp t}.
\]
Then $\psi$ is a quasi-isometry.
\end{lemma}
\begin{proof}
Fix two points $x_i = \geo{h_i(\zero)}{\bp}_{t_i} \in A$, $i = 1,2$. Write $y_i = h_i(\zero)$, $i = 1,2$. Then
\begin{equation}
\label{distinhoroball}
\dist(x_1,x_2) \asymp_\plus |t_2 - t_1| \vee (\dist(y_1,y_2) - t_1 - t_2).
\end{equation}
(This can be seen e.g. by repeated application of Proposition \ref{propositionrips}(ii).) On the other hand, if we write $\w y_i = \Phi(h_i)(\zero)$, $\w t_i = \alpha_\bp t_i$, and $\w x_i = \geo{\w y_i}{\w\bp}_{\w t_i}$, then by \eqref{tukia2} we have $\dist(\w y_1,\w y_2) \asymp_\plus \alpha_\bp \dist(y_1,y_2)$; applying \eqref{distinhoroball} along with its tilded version, we see that $\dist(\w x_1,\w x_2) \asymp_\plus \alpha_\bp \dist(x_1,x_2)$.
\end{proof}

For $g(\bp)\in G(P) = \Lbp(G)$, write $A_{g(\bp)} = g(A_\bp)$ and $\psi_{g(\bp)} = \Phi(g)\circ\psi_\bp\circ g^{-1}$; then $\psi_{g(\bp)}: A_{g(\bp)}\to A_{\phi(g(\bp))}$ is a quasi-isometry, and the implied constants are independent of $g(\bp)$. Let
\[
S = S(G) = \bigcup_{\xi\in \Lbp(G)} A_\xi \supset G(\zero),
\]
and define $\psi:S\to \w S := S(\w G)$ by letting
\[
\psi(x) = \psi_\xi(x) \all \xi\in \Lbp(G) \all x\in A_\xi.
\]
Note that for $g\in G$, $\psi(g(\zero)) = \Phi(g)(\zero)$.
\begin{lemma}
\label{lemmapsiquasi}
$\psi$ is a quasi-isometry.
\end{lemma}
\begin{proof}
Fix two points $x_1,x_2\in S$. For each $i = 1,2$, write $x_i \in A_{g_i(\bp_i)}$ for some $g_i(\bp_i)\in \Lbp(G_i)$. If $g_1(\bp_1) = g_2(\bp_2)$, then $\dist(\psi(x_1),\psi(x_2)) \asymp_\plus \dist(x_1,x_2)$ by Lemma \ref{lemmatukia2}. Otherwise, let $t > 0$ be large enough so that the collection $\scrH = \{H_{g(\bp)} := g(H_{\bp,t}) : g\in G, \; \bp\in P\}$ is disjoint. Then $y_i := \geo{x_i}{g_i(\bp_i)}_t \in H_{g_i(\bp_i)}$. It follows that the geodesic $\geo{y_1}{y_2}$ intersects both $\del H_{g_1(\bp_1)}$ and $\del H_{g_2(\bp_2)}$ (cf. Figure \ref{figurepsiquasi}), say in the points $z_1,z_2$. By Lemma \ref{lemmaboundedparabolic}, there exist points $w_i\in g_i G_{\bp_i}(\zero)$ such that $\dist(z_i,w_i) \asymp_\plus 0$. To summarize, we have
\[
\dist(x_1,x_2) \asymp_{\plus,t} \dist(y_1,y_2) = \dist(z_1,z_2) + \sum_{i = 1}^2 \dist(y_i,z_i) \asymp_{\plus,t} \dist(w_1,w_2) + \sum_{i = 1}^2 \dist(x_i,w_i).
\]
As $x_i,w_i\in A_{g_i(\bp_i)}$, we have $\dist(\psi(x_i),\psi(w_i)) \asymp_\plus \dist(x_i,w_i)$ by Lemma \ref{lemmatukia2}. On the other hand, since $w_1,w_2\in G(\zero)$, we have $\dist(\w w_1,\w w_2) \asymp_{\plus,\times} \dist(w_1,w_2)$ by Lemma \ref{lemmatukia1} and Theorem \ref{theoremGF}(ii). (Here $\w x = \psi(x)$.) Thus,
\[
\dist(x_1,x_2) \asymp_{\plus,\times} \dist(\w w_1,\w w_2) + \sum_{i = 1}^2 \dist(\w x_i,\w w_i) \geq \dist(\w x_1,\w x_2).
\]
Since the situation is symmetric, the reverse inequality holds as well.
\end{proof}

\begin{figure}
\begin{center} 
\begin{tikzpicture}[line cap=round,line join=round,>=triangle 45,x=1.0cm,y=1.0cm]
\clip(-3.674,-3.239) rectangle (3.698,3.586);
\draw(0.0,0.0) circle (3.0cm);
\draw(0.0,2.0) circle (1.0cm);
\draw(-1.6765383941122631,-0.4770413063924982) circle (1.2569145430600153cm);
\draw (-2.884984519391117,-0.8227176446835247)-- (-0.8537077377087111,-0.24806336779043003);
\draw (0.0,3.0)-- (0.0,0.5410285663267452);
\draw [shift={(-4.432995402589801,3.34299557285843)}] plot[domain=5.372236452220289:5.900827991552244,variable=\t]({1.0*4.778027883453108*cos(\t r)+-0.0*4.778027883453108*sin(\t r)},{0.0*4.778027883453108*cos(\t r)+1.0*4.778027883453108*sin(\t r)});
\begin{scriptsize}
\draw[color=black] (0.51,2.093) node {$H_{g_2(p_2)}$};
\draw [fill=black] (-2.884984519391117,-0.8227176446835247) circle (1pt);
\draw[color=black] (-3.3,-0.887) node {$g_1(p_1)$};
\draw[color=black] (-2.0,0.305) node {$H_{g_1(p_1)}$};
\draw [fill=black] (0.0,3.0) circle (1pt);
\draw[color=black] (-0.013254332234670846,3.2) node {$g_2(p_2)$};
\draw [fill=black] (-0.8537077377087111,-0.24806336779043003) circle (1pt);
\draw[color=black] (-0.8250955547283336,-0.45) node {$x_1$};
\draw [fill=black] (0.0,0.5410285663267452) circle (1pt);
\draw[color=black] (0.30213211454936214,0.5204766261395189) node {$x_2$};
\draw [fill=black] (-1.5040817556787958,-0.43205612268083204) circle (1pt);
\draw[color=black] (-1.5209594597229017,-0.65) node {$y_1$};
\draw [fill=black] (0.0,1.56027231456143) circle (1pt);
\draw[color=black] (0.2855639263352058,1.5974088600596753) node {$y_2$};
\end{scriptsize}
\end{tikzpicture}
\caption{The proof of Lemma \ref{lemmapsiquasi}. The distance between $y_1$ and $y_2$ is broken up into three segments, each of which is asymptotically preserved upon applying $\psi$.}
\label{figurepsiquasi}
\end{center}
\end{figure}

\begin{figure}
\begin{center} 
\begin{tikzpicture}[line cap=round,line join=round,>=triangle 45,x=1.0cm,y=1.0cm]
\clip(-3.331733658928572,0.17471897426577723) rectangle (3.5450834715133723,4.0257365673132615);
\draw (-3.0,1.5)-- (3.0,1.5);
\draw [shift={(0.26067583481434903,0.08708448064792945)}] plot[domain=0.5343314243396922:2.8507321995474353,variable=\t]({1.0*2.0675167655199447*cos(\t r)+-0.0*2.0675167655199447*sin(\t r)},{0.0*2.0675167655199447*cos(\t r)+1.0*2.0675167655199447*sin(\t r)});
\draw (-1.06,1.36)-- (-1.06,3.6);
\draw (1.64,1.34)-- (1.64,3.6);
\draw (-1.2487264677147425,1.5)-- (-1.0637012613578571,1.3605243994346345);
\draw (1.641128303962467,1.3391423475348692)-- (1.77007813734344,1.5);
\begin{scriptsize}
\draw [fill=black] (2.04,1.14) circle (1pt);
\draw[color=black] (2.2055990913229415,1.0) node {$y_2$};
\draw [fill=black] (-1.72,0.68) circle (1pt);
\draw[color=black] (-1.8726524590782814,0.5933078430752864) node {$y_1$};
\draw [fill=black] (-1.2487264677147425,1.5) circle (1pt);
\draw[color=black] (-1.3823054984728558,1.7055582659119821) node {$z_1$};
\draw [fill=black] (1.77007813734344,1.5) circle (1pt);
\draw[color=black] (1.8826876782413198,1.7055582659119821) node {$z_2$};
\draw [fill=black] (-1.0637012613578571,1.3605243994346345) circle (1pt);
\draw[color=black] (-0.891,1.17) node {$w_1$};
\draw [fill=black] (1.641128303962467,1.3391423475348692) circle (1pt);
\draw[color=black] (1.511,1.16) node {$w_2$};
\draw[color=black] (0.03,3.9) node {$g(p) = \infty$};
\draw [fill=black] (-0.3562789079803063,2.0604046323337206) circle (1pt);
\draw[color=black] (-0.27005507563615844,2.2317842724153647) node {$x$};
\end{scriptsize}
\end{tikzpicture}
\caption{The proof of Lemma \ref{lemmaSCB}, in the upper half-space model. The thin triangles condition guarantees that $x$ is close to one of the geodesics $\geo{w_1}{g(p)}$, $\geo{w_2}{g(p)}$, both of which are contained in $S$.}
\label{figureSCB}
\end{center}
\end{figure}

\begin{lemma}
\label{lemmaSCB}
$S$ is cobounded in $C = \CC_\zero$.
\end{lemma}
\begin{proof}
Fix $x\in\CC_\zero$. If $x\notin \bigcup\scrH$, then $\dist(x,S)\leq \dist(x,G(\zero))\asymp_\plus 0$. So suppose $x\in H = H_{g(\bp)}$ for some $g\in G$, $\bp\in P$. Write $x\in \geo{y_1}{y_2}$ for some $y_1,y_2\in G(\zero)$. Then there exist $z_1,z_2\in \geo{y_1}{y_2}\cap \del H$ such that $x\in \geo{z_1}{z_2}$. By Lemma \ref{lemmaboundedparabolic}, there exist $w_1,w_2\in g G_\bp(\zero)$ such that $\dist(z_i,w_i)\asymp_\plus 0$. It follows that $\lb w_1|w_2\rb_x\asymp_\plus 0$. By Proposition \ref{propositionrips}, we have
\[
\dist(x,S) \leq \dist(x,S_{g(\bp)}) \leq \dist(x,\geo{w_1}{g(\bp)}\cup\geo{w_2}{g(\bp)}) \asymp_\plus 0
\]
(cf. Figure \ref{figureSCB}). This completes the proof.
\end{proof}

Thus, the embedding map from $S$ to $C$ is an equivariant quasi-isometry. Thus $S$, $C$, $\w S$, and $\w C$ are all equivariantly quasi-isometric. By \cite[Theorem 6.5]{BonkSchramm}, the quasi-isometry between $C$ and $\w C$ extends to a quasisymmetric homeomorphism between $\del C = \Lambda$ and $\del\w C = \w\Lambda$. This completes the proof of Theorem \ref{theoremtukia}(ii).

%
%

\subsection{Completion of the proof assuming only \eqref{tukia}}
We begin by recalling the Morse lemma:

\begin{definition}
\label{definitionquasigeodesic}
A path $\gamma:[a,b]\to X$ is a \emph{$K$-quasigeodesic} if for all $a\leq t_1 < t_2 \leq b$,
\[
\frac1K (t_2 - t_1) - K \leq \dist(\gamma(t_1),\gamma(t_2)) \leq K(t_2 - t_1) + K.
\]
(In other words, $\gamma$ is a $K$-quasigeodesic if $\dist(\gamma(t_1),\gamma(t_2)) \asymp_{\plus,\times} t_2 - t_1$, and the implied constants are both equal to $K$.)
\end{definition}

\begin{lemma}[Morse Lemma, {\cite[Theorem 9.38]{DrutuKapovich}}]
\label{lemmamorse}
For every $K > 0$, there exists $K_2 > 0$ such that the Hausdorff distance between any $K$-quasigeodesic $\gamma$ and the geodesic $\geo{\gamma(a)}{\gamma(b)}$ is at most $K_2$.
\end{lemma}

\begin{lemma}
\label{lemmaquasi}
Fix $h_1,\ldots,h_n\in E\cup F$, let $g_k = h_1\cdots h_k$ and $x_k = g_k(\zero)$ for all $k = 0,\ldots,n$, and suppose that
\begin{equation}
\label{quasi}
\dist(x_k,x_\ell) \asymp_\times \sum_{i = k + 1}^\ell \ell_0(h_i) \all 0\leq k < \ell \leq n.
\end{equation}
Then the path $\gamma = \bigcup_{k = 0}^{n - 1} \geo{x_k}{x_{k + 1}}$ is a $K$-quasigeodesic, where $K > 0$ is independent of $h_1,\ldots,h_n$.
\end{lemma}
\begin{proof}
Fix $0 <  k \leq \ell < n$ and points $z\in\geo{x_{k - 1}}{x_k}$, $w\in\geo{x_\ell}{x_{\ell + 1}}$. To show that $\gamma$ is a quasigeodesic, it suffices to show that
\begin{equation}
\label{ETSquasi}
\dist(z,w) \gtrsim_{\plus,\times} \dist(z,x_k) + \dist(x_k,x_\ell) + \dist(x_\ell,w).
\end{equation}
\begin{claim}
\label{claimquasi}
$\dist(z,w)\gtrsim_\plus \min(\dist(z,x_{k - 1}),\dist(z,x_k))$.
\end{claim}
\begin{subproof}
If $h_k\in F$, then $\dist(z,x_k) \leq \dist(x_{k - 1},x_k) \asymp_\plus 0$, so $\dist(z,w)\gtrsim_\plus \dist(z,x_k)$. Thus, suppose that $h_k\in E$; then $h_k\in G_\bp$ for some $\bp\in P$. Let $g = g_{k - 1}$; since $g^{-1}(z) \in \geo\zero{h_k(\zero)}$, by Proposition \ref{propositionrips}(i) we have $\dist(g^{-1}(z),\geo y\bp) \asymp_\plus 0$, where either $y = \zero$ or $y = h_k(\zero)$.
\begin{subclaim}
\label{subclaimquasi}
There exists $t > 0$ independent of $h_1,\ldots,h_n$ such that $g^{-1}(w)\notin H_{\bp,t}$.
\end{subclaim}
(Cf. Figure \ref{figurequasi}.)

\begin{figure}
\begin{center} 
\begin{tikzpicture}[line cap=round,line join=round,>=triangle 45,x=1.0cm,y=1.0cm]
\clip(-3.8491074720597833,-2.1266876752711386) rectangle (3.5797459638630347,3.7907291690819807);
\draw [shift={(0.0,-0.0)}] plot[domain=0.9045225827260511:3.1035675129286187,variable=\t]({1.0*3.6889564920177627*cos(\t r)+-0.0*3.6889564920177627*sin(\t r)},{0.0*3.6889564920177627*cos(\t r)+1.0*3.6889564920177627*sin(\t r)});
\draw(-1.3068803163851197,1.1823748155886806) circle (1.9266288587793043cm);
\draw [shift={(2.685346326431731,-2.9662333458959456)}] plot[domain=2.049120890160455:2.6001311772486844,variable=\t]({1.0*3.9319392309623806*cos(\t r)+-0.0*3.9319392309623806*sin(\t r)},{0.0*3.9319392309623806*cos(\t r)+1.0*3.9319392309623806*sin(\t r)});
\draw [shift={(4.767483130446365,-0.6701022273279595)}] plot[domain=2.648334556608196:2.8438017124837423,variable=\t]({1.0*4.071162818980422*cos(\t r)+-0.0*4.071162818980422*sin(\t r)},{0.0*4.071162818980422*cos(\t r)+1.0*4.071162818980422*sin(\t r)});
\draw [shift={(5.978973518661347,0.8697896686860462)}] plot[domain=2.9284471092129145:3.0609325657134683,variable=\t]({1.0*4.812997880472604*cos(\t r)+-0.0*4.812997880472604*sin(\t r)},{0.0*4.812997880472604*cos(\t r)+1.0*4.812997880472604*sin(\t r)});
\draw [shift={(3.7794856280301645,2.426664585692047)}] plot[domain=3.0967788621557717:3.353471214637198,variable=\t]({1.0*2.5618834704300895*cos(\t r)+-0.0*2.5618834704300895*sin(\t r)},{0.0*2.5618834704300895*cos(\t r)+1.0*2.5618834704300895*sin(\t r)});
\begin{scriptsize}
\draw [fill=black] (-2.723168609785217,2.4885242057654264) circle (1pt);
\draw[color=black] (-3.036353833196826,2.5787281286723056) node {$g(p)$};
\draw[color=black] (-1.53,2.75) node {$g(H_{p,t})$};
\draw [fill=black] (-0.6841540363190708,-0.9397552627666829) circle (1pt);
\draw[color=black] (-0.45369,-1.18728) node {$x_{k-1}$};
\draw [fill=black] (0.8755044855278509,0.524413961824313) circle (1pt);
\draw[color=black] (1.127,0.397) node {$x_k$};
\draw [fill=black] (1.181623953260003,1.257585677919404) circle (1pt);
\draw [fill=black] (1.2748920600343736,1.887908634247394) circle (1pt);
\draw[color=black] (1.4994382699732513,1.81) node {$x_{l}$};
\draw [fill=black] (1.2242241363497024,2.2425841000400935) circle (1pt);
\draw[color=black] (1.398,2.27) node {$w$};
\draw [fill=black] (1.2201742115269398,2.541433873479229) circle (1pt);
\draw[color=black] (1.38,2.7) node {$x_{l+1}$};
\draw [fill=black] (-0.17745538908965663,-0.2709573946084274) circle (1pt);
\draw[color=black] (-0.3556927085260205,-0.10193299599850501) node {$z$};
\end{scriptsize}
\end{tikzpicture}
\caption{In Subclaim \ref{subclaimquasi}, the geodesics $\geo{x_{k - 1}}{x_k}$ and $\geo{x_\ell}{x_{\ell + 1}}$ cannot penetrate the same cusp, thus guaranteeing some distance between $z$ and $w$.}
\label{figurequasi}
\end{center}
\end{figure}
\begin{subproof}
If $h_{\ell + 1}\in F$, then $\dist(g^{-1}(w),g^{-1}(x_{\ell + 1})) \asymp_\plus 0$, in which case the subclaim follows from the fact that $\bp$ is a bounded parabolic point. Thus suppose $h_{\ell + 1}\in E$; then $h_{\ell + 1}\in G_\eta$ for some $\eta\in P$. Let $k = g_\ell$; since $k^{-1}(w)\in \geo\zero{h_{\ell + 1}(\zero)}$, by Proposition \ref{propositionrips}(i) we have $\dist(k^{-1}(w),\geo p\eta) \asymp_\plus 0$, where either $p = \zero$ or $p = h_k(\zero)$. In particular $\busemann_\eta(\zero,k^{-1}(w))\gtrsim_\plus 0$, so by the disjointness of the family $\scrH$, there exists $t > 0$ such that $k^{-1}(w)\notin j(H_{\bq,t})$ for all $\bq\in P$ and $j\in G$ such that $j(\bq)\neq \eta$. In particular, letting $j = k^{-1} g = (h_k\cdots h_\ell)^{-1}$ and $\bq = \bp$, we have $g^{-1}(w)\notin H_{\bp,t}$ unless $j(\bp) = \eta$. But if $j(\bp) = \eta$, then $j = \id$ due to the minimality $P$, and this contradicts \eqref{quasi}.
\end{subproof}
It follows that
\begin{align*}
\dist(z,w) \geq \busemann_\bp(g^{-1}(w),g^{-1}(z))
&= \busemann_\bp(\zero,g^{-1}(z)) - \busemann_\bp(\zero,g^{-1}(w))\\
&= \dist(y,g^{-1}(z)) - \busemann_\bp(\zero,g^{-1}(w)) \gtrsim_\plus \dist(y,g^{-1}(z)) - t.
\end{align*}
Applying $g$ to both sides finishes the proof of Claim \ref{claimquasi}.
\end{subproof}

A similar argument shows that $\dist(z,w)\gtrsim_\plus \min(\dist(w,x_\ell),\dist(w,x_{\ell + 1}))$. Now let $y_1\in \{x_{k - 1},x_k\}$ and $y_2\in \{x_\ell,x_{\ell + 1}\}$ be such that
\begin{equation}
\label{zwzyzwwy}
\dist(z,w) \gtrsim_\plus \dist(z,y_1) \text{ and } \dist(z,w) \gtrsim_\plus \dist(w,y_2).
\end{equation}
Then the triangle inequality gives $\dist(z,w)\gtrsim_{\plus,\times} \dist(y_1,y_2)$. On the other hand, \eqref{quasi} implies that $\dist(y_1,y_2) \gtrsim_{\plus,\times} \dist(y_1,x_k) + \dist(x_k,x_\ell) + \dist(x_\ell,y_2)$. Combining with \eqref{zwzyzwwy} and using the triangle inequality gives \eqref{ETSquasi}.
\end{proof}

\begin{lemma}
\label{lemmagromovquasi}
For all $x,y,z\in G(\zero)$,
\[
\lb \w x | \w y\rb_{\w z} \asymp_{\plus,\times} \lb x | y \rb_z.
\]
\end{lemma}
\begin{proof}
Fix $g_1,g_2\in G$, and we will show that
\begin{equation}
\label{ETSgromovquasi}
\lb \w g_1(\zero) | \w g_2(\zero)\rb_\zero \lesssim_{\plus,\times} \lb g_1(\zero) | g_2(\zero)\rb_\zero.
\end{equation}
The reverse inequality will then follow by symmetry. By Theorem \ref{theoremGF}(ii), there exists a sequence $h_1,\ldots,h_n\in E\cup F$ such that $g_2 = g_1 h_1\cdots h_n$ and satisfying \eqref{quasi}. By Lemma \ref{lemmatukia1}, the sequence $\w h_1,\ldots,\w h_n\in \w E\cup \w F$ also satisfies \eqref{quasi}. Let $x_k = g_1 h_1 \cdots h_k(\zero)$. By Lemma \ref{lemmaquasi}, the paths
\begin{align*}
\gamma &= \bigcup_{k = 0}^{n - 1} \geo{x_k}{x_{k + 1}}\\
\w\gamma &= \bigcup_{k = 0}^{n - 1} \geo{\w x_k}{\w x_{k + 1}}
\end{align*}
are quasigeodesics. So by Lemma \ref{lemmamorse}, $\gamma$ and $\w\gamma$ lie within a bounded Hausdorff distance of the geodesics they represent, namely $\geo{x_0}{x_n}$ and $\geo{\w x_0}{\w x_n}$. Combining with Proposition \ref{propositionrips}(i), we have
\[
\lb g_1(\zero) | g_2(\zero)\rb_\zero = \lb x_0 | x_n \rb_\zero \asymp_\plus \dist(\zero,\geo{x_0}{x_n}) \asymp_\plus \dist(\zero,\gamma),
\]
and similarly for $\w\gamma$. So to prove \eqref{ETSgromovquasi}, we need to show that $\dist(\zero,\w\gamma) \lesssim_{\plus,\times} \dist(\zero,\gamma)$.

Fix $z\in\gamma$, and we will show that $\dox z \gtrsim_{\plus,\times} \dist(\zero,\w\gamma)$. Write $z\in \geo{x_{k - 1}}{x_k}$ for some $k = 1,\ldots,n$. By Proposition \ref{propositionrips}(i), we have
\begin{align*}
\dist(\zero,\w\gamma) &\leq \dist(\zero,\geo{\w x_{k - 1}}{\w x_k}) \asymp_\plus \lb \w x_{k - 1} | \w x_k \rb_\zero\\
\dox z &\geq \dist(\zero,\geo{x_{k - 1}}{x_k}) \asymp_\plus \lb x_{k - 1} | x_k \rb_\zero,
\end{align*}
so to complete the proof of Lemma \ref{lemmagromovquasi} it suffices to show that
\begin{equation}
\label{gromovquasi}
\lb \w x_{k - 1} | \w x_k \rb_\zero \lesssim_{\plus,\times} \lb x_{k - 1} | x_k \rb_\zero.
\end{equation}
Now, if $h_k\in F$, then $\dist(x_{k - 1},x_k) \asymp_\plus \dist(\w x_{k - 1},\w x_k) \asymp_\plus 0$, so \eqref{gromovquasi} follows from Theorem \ref{theoremGF}(ii). Thus, suppose that $h_k\in E$, and write $h_k\in G_\bp$ for some $\bp\in P$. Use the notations $g = g_1 h_1\cdots h_{k - 1}$ and $h = h_k$, so that $x_{k - 1} = g(\zero)$ and $x_k = gh(\zero)$. Then for $y = \zero,h(\zero)$, (h) of Proposition \ref{propositionbasicidentities} gives
\[
\lb y|\bp\rb_{g^{-1}(\zero)} \asymp_\plus \frac12[\dist(g^{-1}(\zero),y) - \busemann_\bp(g^{-1}(\zero),y)] \gtrsim_\plus \frac12\dist(g^{-1}(\zero),y) \asymp_\times \dox{g(y)},
\]
so by Gromov's inequality,
\[
\lb x_{k - 1} | x_k \rb_\zero = \lb g(\zero)|gh(\zero)\rb_\zero = \lb \zero|h(\zero)\rb_{g^{-1}(\zero)} \gtrsim_{\plus,\times} \dox{g(y)} \asymp_{\plus,\times} \dox{\w g(\w y)} \geq \lb \w x_{k - 1} | \w x_k \rb_\zero. 
\]
This demonstrates \eqref{gromovquasi} and completes the proof of Lemma \ref{lemmagromovquasi}.
\end{proof}

It follows that the map $\Phi$ sends Gromov sequences to Gromov sequences, so it induces an equivariant homeomorphism $\del\Phi:\Lambda\to \w\Lambda$. This completes the proof of Theorem \ref{theoremtukia}(i).

\subsection{Applications to finite-dimensional \CROSSs}
A particularly interesting case of Theorem \ref{theoremtukia} is when $X$ and $\w X$ are both finite-dimensional \CROSSs. In this case, \eqref{tukia} always holds, but \eqref{tukia2} does not; nevertheless, there is a reasonable sufficient condition for \eqref{tukia2} to hold. Specifically, we have the following:

\begin{proposition}
\label{propositionCROSSONCTquasi}
Let $X$ and $\w X$ be finite-dimensional \CROSSs, let $G\leq\Isom(X)$ and $\w G\leq \Isom(\w X)$ be geometrically finite groups, and let $\Phi:G\to\w G$ be a type-preserving isomorphism. Fix $\bp\in P$, and let $\w\bp = \phi(\bp)\in \Lbp(\w G)$ be the unique point such that $\Phi(G_\bp) = \w G_{\w\bp}$. Then
\begin{itemize}
\item[(i)] \eqref{tukia} holds.
\item[(ii)] Let $H\leq G_\bp$ be a nilpotent subgroup of finite index. If the underlying base fields of $X$ and $\w X$ are the same, say $\F$, and if $\rank([H,H]) = \dim_\R(\F) - 1$, then \eqref{tukia2} holds.
\end{itemize}
\end{proposition}

Before we begin the proof of Proposition \ref{propositionCROSSONCTquasi}, it will be necessary to understand the structure of a parabolic subgroup of $\Isom(X)$.

Let $X = \H = \H_\F^d$ be a finite-dimensional \CROSS, let $\bp = [(1,1,\0)]$, and let $J_\bp = \Stab(\Isom(X);\bp)$. Note that $J_\bp$ is a parabolic group in the sense of Lie theory, while it is a focal group according to the classification of Section \ref{sectionclassification}. To study the group $J_\bp$, we use the coordinate system generated by the basis
\[
\ff_0 = (\ee_0 + \ee_1)/2, \;\; \ff_1 = \ee_1 - \ee_0, \;\; \ff_i = \ee_i \;\;\; (i = 2,\ldots,d).
\]
In this coordinate system, the sesquilinear form $B_\QQ$ takes the form
\[
B_\QQ(\xx,\yy) = \wbar x_0 y_1 + \wbar x_1 y_0 + \sum_{i = 2}^d \wbar x_i y_i,
\]
the point $\bp$ takes the form $\bp = [\ff_0]$, and the group $J_\bp$ can be written (cf. Theorem \ref{theoremisometries}) as
\[
J_\bp = \left\{ h_{\lambda,a,\vv,\ww,m,\sigma} := \left[\begin{array}{lll}
\lambda & a & \ww^\dag\\
& \lambda^{-1} &\\
& \vv & m
\end{array}\right] \sigma^{d + 1} :
\begin{split}
\lambda > 0,\; a\in \F,\;  \vv,\ww\in \F^{d - 1},\\
m\in \SO(\F^{d - 1};\EE),\; \sigma\in\Aut(\F)
\end{split}\right\} \cap \Isom(X).
\]
Given $\lambda,a,\vv,\ww,m$, it is readily verified that $h_{\lambda,a,\vv,\ww,m}\in \Isom(X)$ if and only if
\[
2\lambda^{-1}\Re(a) + \|\vv\|^2 = 0 \text{ and } \lambda^{-1}\ww^\dag + \vv^\dag m = \0.
\]
Consequently, it makes sense to rewrite $J_\bp$ as
\[
J_\bp = \left\{ h_{\lambda,a,\vv,m,\sigma} := \left[\begin{array}{lll}
\lambda & a - \lambda\|\vv\|^2/2 & -\lambda\vv^\dag m\\
& \lambda^{-1} &\\
& \vv & m
\end{array}\right] \sigma^{d + 1} :
\begin{split}
\lambda > 0,\; a\in \Im(\F),\;  \vv\in \F^{d - 1},\\
m\in \SO(\F^{d - 1};\EE),\; \sigma\in\Aut(\F)
\end{split}\right\}.
\]
We can now define the \emph{Langlands decomposition} of $J_\bp$:
\begin{align*}
M_\bp &= \{h_{1,0,\0,m,\sigma} : m\in\SO(\F^{d - 1};\EE),\; \sigma\in \Aut(\F)\}\\
A_\bp &= \{h_{\lambda,0,\0,I_{d - 1},e} : \lambda > 0\}\\
N_\bp &= \{n(a,\vv) := h_{1,a,\vv,I_{d - 1},e} : a\in\Im(\F) ,\; \vv\in \F^{d - 1}\}\\
J_\bp &= M_\bp A_\bp N_\bp.
\end{align*}
We observe the following facts about the Langlands decomposition: the groups $M_\bp$ and $A_\bp$ commute with each other and normalize $N_\bp$, which is nilpotent of order at most $2$. Moreover, the subgroup $M_\bp N_\bp$ is exactly the kernel of the homomorphism $J_\bp\ni h\mapsto h'(\bp)$, where $h'$ denotes the metric derivative. Equivalently, $M_\bp N_\bp$ is the largest parabolic subgroup of $J_\bp$, where ``parabolic'' is interpreted in the sense of Section \ref{sectionclassification}.

Let's look a bit more closely at the internal structure of $N_\bp$. The composition law is given by
\begin{equation}
\label{compositionlaw}
n(a_1,\vv_1) n(a_2,\vv_2) = n(a_1 + a_2 + \Im B_\EE(\vv_2,\vv_1),\vv_1 + \vv_2),
\end{equation}
confirming that $N_\bp$ is nilpotent of order at most two, and that its commutator is given by
\[
Z_\bp = \{n(a,\0) : a\in\Im(\F)\}.
\]
Moreover, the map $\pi:n(a,\vv)\mapsto \vv\in \F^{d - 1}$ is a homomorphism whose kernel is $Z_\bp$.

Now let $H\leq M_\bp N_\bp$ be a discrete parabolic subgroup. By Margulis's lemma, $H$ is almost nilpotent, and so by \cite[Lemma 3.4]{CorletteIozzi}, there exist a finite index subgroup $H_2\subset H$ and a homomorphism $\psi:H_2\to N_\bp$ such that $\psi(h)(\zero) = h(\zero)$ for all $h\in H_2$. (Here $\zero = [\ee_0] = [2\ff_0 - \ff_1]$ as usual.) We then let $H_3 = \psi(H_2) \leq N_\bp$.
\begin{definition}
\label{definitionHregular}
The group $H$ is \emph{regular} if $\pi(H_3)$ is a discrete subgroup of $\F^{d - 1}$. If $H$ is regular, we define its \emph{quasi-commutator} to be the subgroup
\[
Z = Z(H) = \psi^{-1}(Z_\bp) = \Ker(\pi\circ\psi) \leq H.
\]
Note that in general, the quasi-commutator of $H$ cannot be determined from its algebraic structure; cf. Example \ref{exampletukia}. Nevertheless, since $\F^{d - 1}$ is abelian, the quasi-commutator of $H$ always contains the commutator of $H_2$.
\end{definition}

In general, if $H\leq\Isom(X)$ is a discrete parabolic subgroup, we can conjugate the fixed point of $H$ to $[(1,1,\0)]$, apply the above construction, and then conjugate back to get a subgroup $Z(H)\leq H$.

If $H$ is regular, then the quasi-commutator $Z\leq H$ can be used to give an algebraic description of the function $h\mapsto \dogo h$. Specifically, we have the following:

\begin{lemma}
\label{lemmaCROSSONCTquasi}
Let $\dist_H$ and $\dist_Z$ be any Cayley metrics on $H$ and $Z$, respectively.
\begin{itemize}
\item[(i)]
\begin{equation}
\label{CROSSONCTquasi1}
\dogo h \asymp_{\plus,\times} 0\vee \log\dist_H(e,h).
\end{equation}
\item[(ii)] If $H$ is regular, then
\begin{equation}
\label{CROSSONCTquasi2}
\dogo h \asymp_\plus \min_{z\in Z} \big(0\vee 2\log\dist_H(z,h) \vee \log\dist_Z(e,z)\big)  \all h\in H.
\end{equation}
\end{itemize}
\end{lemma}
\begin{proof}
Let $F\subset H$ be a finite set so that $H_2 F = H$, and let $H_3 = \psi(H_2)$. Then for all $h\in H$, we can write $h = h_2 f$ for some $h_2\in H_2$ and $f\in F$, and then
\begin{align*}
\dogo h &\asymp_\plus \dogo{h_2} = \dogo{\psi(h_2)}\\
\dist_H(z,h) &\asymp_\plus \dist_H(z,h_2) \asymp_\times \dist_{H_2}(z,h_2) = \dist_{H_3}(\psi(z),\psi(h_2))\\
\min_{z\in Z} \big(0\vee 2\log\dist_H(z,h) \vee \log\dist_Z(e,z)\big) & \asymp_\plus \min_{z\in \psi(Z)} \big(0\vee 2\log\dist_{H_3}(z,\psi(h_2)) \vee \log\dist_{\psi(Z)}(e,z)\big).
\end{align*}
Thus, we may without loss of generality assume that $H = H_3$, i.e. that $H\leq N_\bp$ and $Z_H = H\cap Z_\bp$. We can also without loss of generality assume that $\bp = [(1,1,\0)]$.

The following formula regarding the function $n(a,\vv)$ can be verified by direct computation (cf. \cite[(3.5)]{CorletteIozzi}):
\begin{equation}
\label{corletteiozzi}
\|n(a,\vv)\| \asymp_\plus 0\vee 2\log\|\vv\| \vee \log|a|
\end{equation}
On the other hand, iterating \eqref{compositionlaw} gives
\begin{equation}
\label{compositionlawbounds}
\begin{split}
\|\vv\| &\lesssim_\times \dist_H(e,n(a,\vv))\\
|a| &\lesssim_\times \dist_H(e,n(a,\vv))^2\\
|a| &\lesssim_\times \dist_Z(e,n(a,\0)).
\end{split}
\end{equation}
These formulas make it easy to verify the $\lesssim$ direction of \eqref{CROSSONCTquasi2}: given $h = n(a,\vv)\in H$ and $z = n(b,\0)\in Z$, we have
\begin{align*}
0\vee 2\log\dist_H(z,h) \vee \log\dist_Z(e,z)
&=_\pt 0\vee 2\log\dist_H(e,n(a - b,\vv)) \vee \log\dist_Z(e,n(b,\0))\\
&\geq_\pt 0\vee 2\log\big(\|\vv\|\vee \sqrt{|a - b|}\big)\vee \log|b|\\
&=_\pt 0\vee 2\log\|\vv\| \vee \log|a - b| \vee\log|b|\\
&\gtrsim_\plus 0\vee 2\log\|\vv\| \vee \log|a| \asymp_\plus \dogo{n(a,\vv)} = \dogo h.
\end{align*}
Setting $z = e$ yields the $\lesssim$ direction of \eqref{CROSSONCTquasi1}.

To prove the $\gtrsim$ directions, we will need the following easily verified fact:

\begin{fact}
\label{factFDVS}
If $V$ is a finite-dimensional vector space, $\Lambda\leq V$ is a discrete subgroup, and $\dist_\Lambda$ is a Cayley metric on $\Lambda$, then $\dist_\Lambda(\0,\vv) \asymp_\times \|\vv\|$ for all $\vv\in\Lambda$. Here $\|\cdot\|$ denotes any norm on $V$.
\end{fact}

To prove the $\gtrsim$ direction of \eqref{CROSSONCTquasi2}, assume that $H$ is regular, fix $h = n(a,\vv)\in H$, and let $F$ be a finite generating set for $H$. Since $H$ is regular, the group $\Lambda = \pi(H) \leq \F^{d - 1}$ is discrete. Since $\F^{d - 1}$ is a finite-dimensional vector space, Fact \ref{factFDVS} guarantees the existence of a sequence $f_1,\ldots,f_n\in F$ such that $\pi(f_1\cdots f_n) = \pi(h)$ and $n \lesssim_\times \|\vv\|$. Let $f = f_1\cdots f_n$ and let $z = hf^{-1} \in \pi^{-1}(0) = Z$, say $z = n(b,\0)$. Applying \eqref{compositionlaw} and the second equation of \eqref{compositionlawbounds}, we see that $|b| \lesssim_\times |a| \vee \|\vv\|^2 \vee n^2 \lesssim_\times |a| \vee \|\vv\|^2$. On the other hand, applying Fact \ref{factFDVS} to $Z_\bp$ gives $\dist_Z(e,z) \lesssim_\times |b|$. Thus
\begin{align*}
0\vee 2\log\dist_H(z,h) \vee \log\dist_Z(e,z) &=_\pt 0\vee 2\log\dist_H(e,f) \vee \log\dist_Z(e,z)\\
&\lesssim_\plus 0\vee 2\log(n)\vee \log|b|\\
&\lesssim_\plus 0\vee 2\log\|\vv\| \vee \log(|a| \vee \|\vv\|^2)\\
&=_\pt 0\vee 2\log\|\vv\| \vee \log|a| = \dogo h.
\end{align*}
This completes the proof of \eqref{CROSSONCTquasi2}.

To prove the $\gtrsim$ direction of \eqref{CROSSONCTquasi1}, let $\cl H$ and $\cl Z$ be the Zariski closures of $H$ and $Z$ in $N_\bp$, respectively. Then $\cl H/\cl Z$ and $\cl Z$ are abelian Lie groups, and therefore isomorphic to finite-dimensional vector spaces. Let $\w\pi:\cl H\to\cl H/\cl Z$ be the projection map. Note that $\|\w\pi(n(a,\vv))\| \lesssim_\times |a| \vee \|\vv\|$ for all $n(a,\vv)\in H$. Here $\|\cdot\|$ denotes any norm on $\cl H/\cl Z$.

Since $\cl Z$ is a vector space, the fact that $Z$ is Zariski dense in $\cl Z$ simply means that $Z$ is a lattice in $\cl Z$. In particular, $Z$ is cocompact in $\cl Z$, which implies that $\w\pi(H)$ is discrete. Fix $h = n(a,\vv)\in H$, and let $F$ be a finite generating set for $H$. Then by Fact \ref{factFDVS}, there exists a sequence $f_1,\ldots,f_n\in F$ such that $\w\pi(f_1)\cdots\w\pi(f_n) = \w\pi(h)$ and $n\lesssim_\times \|\w\pi(h)\|\lesssim_\times |a| \vee \|\vv\|$. Let $f = f_1\cdots f_n$ and let $z = h f^{-1} \in H\cap \w\pi^{-1}(0) = H\cap \cl Z = Z$, say $z = n(b,\0)$. Applying \eqref{compositionlaw} and the second equation of \eqref{compositionlawbounds}, we see that $|b| \lesssim_\times |a| \vee \|\vv\|^2 \vee n^2 \lesssim_\times |a|^2 \vee \|\vv\|^2$. On the other hand, applying Fact \ref{factFDVS} to $\cl Z$ gives $\dist_Z(e,z) \lesssim_\times |b|$. Thus
\begin{align*}
0\vee\log\dist_H(e,h) &\leq_\pt 0\vee \log\dist_H(e,f) \vee \log\dist_Z(e,z)\\
&\lesssim_\plus 0\vee \log(n)\vee \log|b|\\
&\lesssim_\plus 0\vee \log(|a|\vee \|\vv\|) \vee \log(|a|^2\vee \|\vv\|^2)\\
&\asymp_\times 0\vee 2\log\|\vv\| \vee \log|a| = \dogo h.
\end{align*}
This completes the proof of \eqref{CROSSONCTquasi1}.
\end{proof}

\begin{corollary}
\label{corollaryCROSSONCTquasi}
Let $X$ and $\w X$ be finite-dimensional \CROSSs, let $H\leq\Isom(X)$ and $\w H\leq \Isom(\w X)$ be parabolic groups with fixed points $\bp$ and $\w\bp$, respectively, and let $\Phi:H\to\w H$ be an isomorphism. Then
\begin{itemize}
\item[(i)] \eqref{tukia} holds.
\item[(ii)] If $H$ and $\w H$ are regular, then \eqref{tukia2} holds if and only if $\Phi(Z)$ is commensurable to $\w Z$. Here $Z = Z(H)$ and $\w Z = Z(\w H)$.
\end{itemize}
\end{corollary}
\begin{proof}
\eqref{tukia} follows immediately from \eqref{CROSSONCTquasi1}. 
Suppose that $H$ and $\w H$ are regular and that $\Phi(Z)$ is commensurable to $\w Z$. Since the right hand side of \eqref{CROSSONCTquasi2} depends on both $h$ and $Z$, let us write it as a function $R(h,Z)$. We then have
\[
\dogo h \asymp_\plus R(h,Z) = R(\w h, \Phi(Z)) \asymp_\plus R(\w h,\w Z) \asymp_\plus \|\w h\|.
\]
On the other hand, suppose that $\Phi(Z)$ and $\w Z$ are not commensurable. Without loss of generality, suppose that the index of $\Phi(Z)\cap \w Z$ in $\Phi(Z)$ is infinite. Since $\Phi(Z)$ is a finitely generated abelian group, it follows that there exists $\w h = \Phi(h)\in \Phi(Z)$ such that $\w h^n\notin \w Z$ for all $n\in\Z\butnot\{0\}$. Without loss of generality, suppose that $\w h\in \w H_2$; otherwise replace $h$ by an appropriate power. Then \eqref{corletteiozzi} implies that
\[
\|h^n\| \asymp_{\plus,h} \log(n) \text{ but } \|\w h^n\| \asymp_{\plus,h} 2\log(n).
\]
Thus \eqref{tukia2} fails along the sequence $(h_n)_1^\infty$.
\end{proof}

\begin{corollary}
\label{corollaryrealROSSONCTquasi}
In the context of Corollary \ref{corollaryCROSSONCTquasi}, if $X$ and $\w X$ are both real \ROSSs, then \eqref{tukia2} holds.
\end{corollary}
\begin{proof}
Since $\Im(\R) = \{0\}$, the group $Z_\bp$ is trivial and thus $Z$ and $\w Z$ are trivial as well; moreover, every discrete parabolic group is regular.
\end{proof}

\begin{corollary}
\label{corollarylatticequasi}
In the context of Corollary \ref{corollaryCROSSONCTquasi}, if we assume both that
\begin{itemize}
\item[(I)] $H$ is a lattice in $M_\bp N_\bp$, and that
\item[(II)] the underlying base fields of $X$ and $\w X$ are the same, or at least satisfy $\dim_\R(\F)\geq\dim_\R(\w\F)$,
\end{itemize}
then \eqref{tukia2} holds.
\end{corollary}
\begin{proof}
Let $H_2$, $\psi$, $H_3$, and $Z = Z(H)$ be as on page \pageref{lemmaCROSSONCTquasi}. Without loss of generality, we may assume that $H = H_3$ and $\w H = \w H_3$. Then $H$ is a lattice in $N_\bp$ and $\w H\leq \w N_{\w\bp}$.

Since $H$ is a lattice in $N_\bp$, $H$ is Zariski dense in $N_\bp$; this implies that $[H,H]$ is Zariski dense in $Z_\bp = [N_\bp,N_\bp]$. Thus, the rank of $[H,H]$ (and also of $\Phi([H,H]) = [\w H,\w H]$) is equal to $\dim_\R(\Im(\F)) = \dim_\R(\F) - 1$. Thus $\dim_\R(\w\F) - 1 = \rank([H,H]) \leq \dim(Z_\bp) = \dim_\R(\w\F) - 1$. Since by assumption $\dim_\R(\F)\geq\dim_\R(\w\F)$, equality holds. Thus $Z$ is a lattice in $Z_\bp$ and is commensurable to $[H,H]$. Similarly, $\w Z$ is a lattice in $\w Z_{\w\bp}$ and is commensurable to $[\w H,\w H]$. Thus, the groups $H$ and $\w H$ are regular. Finally, $\w Z$ is commensurable to $[\w H,\w H] = \Phi([H,H])$ which is commensurable to $\Phi(Z)$, so Corollary \ref{corollaryCROSSONCTquasi} finishes the proof.
\end{proof}

\begin{remark}
If $G\leq\Isom(X)$ is a lattice, then every parabolic subgroup $G_\bp$ satisfies (I).
\end{remark}

As an application, we generalize a rigidity result due to X. Xie \cite[Theorem 3.1]{Xie}:

\begin{corollary}
\label{corollaryxie}
Let $X$, $\w X$ be finite-dimensional \ROSSs\ over the same base field, with $X\neq\H_\R^2$. Let $G\leq\Isom(X)$ be a noncompact lattice, and let $\w G\leq\Isom(\w X)$ be a geometrically finite group, both torsion-free. Let $\Phi:G\to \w G$ be a type-preserving isomorphism. Then $\HD(\Lambda(\w G)) \geq \HD(\Lambda(G)) = \dim(\del X)$. Furthermore, equality holds if and only if $\w G$ stabilizes an isometric copy of $X$ in $\w X$.
\end{corollary}
\begin{proof}
Xie has observed that the main result of his paper generalizes to \ROSSs\ once one verifies that Tukia's isomorphism theorem and the Global Measure Formula both generalize to that setting (cf. \cite[p.1]{Xie}). We have just shown that Tukia's isomorphism theorem generalizes (to the present setting at least), and the Global Measure Formula has been shown to generalize by B. Schapira \cite[Th\'eor\`eme 3.2]{Schapira}.

Actually, we should mention a minor change that needs to be made to Xie's proof in the setting of \ROSSs: Since the Hausdorff and topological dimensions of the boundary of a nonreal \ROSS\ are not equal, \commariusz{Quite often both are equal to infinity. Does it cause a problem?} \comdavid{In the setting of this corollary, everything is finite-dimensional.} at the top of \cite[p.252]{Xie} one should use Pansu's lemma \cite[Proposition 6.5]{Pansu1}, \cite[Lemma 2.3(a)]{Xie} to deduce the lower bound on the Hausdorff dimension of $\Lambda(G_2)$ (i.e. \cite[p.252, line 4]{Xie}) rather than using Szpilrajn's inequality between Hausdorff and topological dimensions (cf. \cite[p.252, lines 2-3]{Xie}).
\end{proof}

Note that in Xie's proof, quasisymmetry is used in an essential way due to his use of Pansu's lemma \cite[Corollary 7.2]{Pansu1}, \cite[Lemma 2.3]{Xie}. Thus, the fact that the stronger asymptotic \eqref{tukia2} holds in the context of Corollary \ref{corollarylatticequasi} is essential to the proof of Corollary \ref{corollaryxie}. It remains to be answered whether Corollary \ref{corollaryxie} holds if we drop the assumption of identical base fields.

We end this section by giving an example of groups for which \eqref{tukia2} fails.

\begin{example}
\label{exampletukia}
Let $\H = \H_\C^3$, let $\bp = [(1,1,\0)]$, and define a homomorphism $\theta:\R^3\to N_\bp$ by $\theta(x,y,z) = n(xi ,(y,z))$, where $i = \sqrt{-1}$. Consider the parabolic groups $H,H',H''\leq N_\bp$ defined by
\begin{align*}
H &= \theta(\Z\times\Z\times\{0\})\\
H' &= \theta(\Lambda\times\{0\})\\
H'' &= \theta(\{0\}\times\Z\times\Z).
\end{align*}
In the middle equation, $\Lambda$ denotes a lattice in $\R^2$ which does not intersect the axes. Then the groups $H,H',H''$ are all isomorphic, but we will show below that \eqref{tukia2} cannot hold for any isomorphisms between them. This is accounted for in Corollary \ref{corollaryCROSSONCTquasi} as follows: The group $H'$ is irregular, so Corollary \ref{corollaryCROSSONCTquasi} does not apply; The groups $Z(H)$ and $Z(H'')$ are not almost isomorphic (the former is isomorphic to $\Z$ while the latter is isomorphic to $\{0\}$), so Corollary \ref{corollaryCROSSONCTquasi} does not apply.
\end{example}
\begin{proof}
Note that the function $\dogo\cdot$ is described on $\theta(\R^3)$ by
\[
\dogo{\theta(x,y,z)} \asymp_\plus 0\vee \log|x|\vee 2\log(|y|\vee|z|)
\]
(cf. \eqref{corletteiozzi}). Now let $h_1 = \theta((1,0,0)) \in H$, $h_2 = \theta((0,1,0))\in H$. Then
\[
\|h_i^n\| \asymp_\plus i \log(n);
\]
but if $\Phi$ is an isomorphism from $H$ to either $H'$ or $H''$, then
\[
\|\Phi(h_i)^n\| \asymp_\plus 2 \log(n).
\]
This demonstrates the failure of \eqref{tukia2}, as setting $h = h_1^n$ gives $\alpha_\bp = 2$ while setting $h = h_2^n$ gives $\alpha_\bp = 1$.

Next, let $\dist_{H'}$ and $\dist_{H''}$ be Cayley metrics on $H'$ and $H''$, respectively. Then for all $R\geq 1$,
\[
\sup_{\dist_{H'}(e,h') \leq R} \|h'\| \asymp_\plus 2\log(R) > \log(R) \asymp_\plus \inf_{\dist_{H'}(e,h') > R} \|h'\|.
\]
but
\[
\sup_{\dist_{H''}(e,h'') \leq R} \|h''\| \asymp_\plus \inf_{\dist_{H''}(e,h'') > R} \|h''\| \asymp_\plus 2\log(R)
\]
This demonstrates the failure of \eqref{tukia2} for any isomorphism between $H'$ and $H''$, as taking the supremum over a ball in the Cayley metric gives $\alpha_\bp = 1$, while taking the infimum over the complement of a ball in the Cayley metric gives $\alpha_\bp = 2$.
%
%
\end{proof}

%

%

\begin{remark}
\label{remarktukia}
The above proof actually shows more; namely, it shows that if $\Phi:G\to \w G$ is a type-preserving isomorphism so that for some $\bp\in\Lbp$, $G_\bp$ and $\w G_{\w\bp}$ are distinct elements of $\{H,H',H''\}$, then the equivariant boundary extension of $\Phi$ is not quasisymmetric.
\end{remark}
\begin{proof}
By contradiction suppose that the equivariant boundary extension $\phi:\Lambda\to\w\Lambda$ is quasisymmetric. Fix $\zeta\in \Lambda\butnot\{\bp\}$, and let $\w\zeta = \phi(\zeta)$. Then by equivariance, for each $h\in G_\bp$ we have
\[
\phi(h(\zeta)) = \w h(\w\zeta).
\]
Let $f:(0,\infty)\to(0,\infty)$ be as in Definition \ref{definitionquasisymmetric}, so that for all $\xi,\eta_1,\eta_2\in \Lambda$,
\[
\frac{\w\Dist(\w\xi,\w\eta_2)}{\w\Dist(\w\xi,\w\eta_1)} \leq f\left(\frac{\Dist(\xi,\eta_2)}{\Dist(\xi,\eta_1)}\right).
\]
Letting $\xi = \bp$ and $\eta_i = h_i(\zeta)$ gives
\[
\frac{\w\Dist(\w\bp,\w h_2(\w\zeta))}{\w\Dist(\w\bp,\w h_1(\w\zeta))} \leq f\left(\frac{\Dist(\bp,h_2(\zeta))}{\Dist(\bp,h_1(\zeta))}\right).
\]
But $\Dist(\bp,h_i(\zeta)) \asymp_{\times,\zeta} \Dist(\bp,h_i(\zero)) = e^{(1/2)\dogo{h_i}}$; thus
\[
\exp\left(\frac12 \left[\|\w h_2\| - \|\w h_1\|\right]\right) \leq f_2\exp\left(\frac12 \Big[\dogo{h_2} - \dogo{h_1}\Big]\right),
\]
where $f_2(t) = Cf(Ct)$ for some constant $C > 0$. Letting $f_3(t) = 2\log f_2(e^{(1/2)t})$ gives
\[
\|\w h_2\| - \|\w h_1\| \leq f_3(\dogo{h_2} - \dogo{h_1}).
\]
But this is readily seen to contradict the proof of Example \ref{exampletukia}.
\end{proof}

}
%
%

\chapter{Counterexamples}\label{sectionexamples}

In Chapter \ref{sectiondiscreteness} we defined various notions of discreteness and demonstrated some relations between them, and in Section \ref{subsectionpropositionpoincareregular} we related some of these notions to the modified Poincar\'e exponent $\w\delta$. In this chapter we give counterexamples to show that the relations which we did not prove are in fact false. Specifically, we prove that no more arrows can be added to Table \ref{figurediscreteness} (reproduced below as Table \ref{figurediscretenesscopy}), and that the discreteness hypotheses of Proposition \ref{propositionpoincareregular} cannot be weakened.

\begin{table}[h!]
\begin{tabular}{ | c | c c c c c c c | }
 \hline
 Finite dimensional			& SD & $\leftrightarrow$ & MD & $\leftrightarrow$ & WD && \\
 Riemannian manifold 	& $\uparrow$ & & & & $\updownarrow$ && \\
	 			& PrD & & & & COTD & $\leftrightarrow$ & UOTD\\
 \hline
	 		& SD & $\rightarrow$ & MD & $\rightarrow$ & WD && \\
 General metric space		& & $\nearrow$ & & $\searrow$ & && \\
				& PrD & & & & COTD && \\
 \hline
 Infinite dimensional				& SD & $\rightarrow$ & MD & $\rightarrow$ & WD && \\
 algebraic hyperbolic space 	& & $\nearrow$ & & & $\downarrow$ && \\
			& PrD & & & & COTD & $\rightarrow$ & UOTD \\
 \hline
	 		& SD & $\leftrightarrow$ & MD & $\leftrightarrow$ & COTD && \\
 Proper metric  space 		& $\uparrow$ & & & & $\downarrow$ &&\\
			& PrD & & & & WD && \\
 \hline
\end{tabular}
\vskip0.12cm
\caption{The relations between different notions of discreteness. COTD and UOTD stand for discrete with respect to the compact-open and uniform operator topologies respectively. All implications not listed have counterexamples, which are described below.}
\label{figurediscretenesscopy}
\end{table}


The examples are arranged roughly in order of discreteness level; the most discrete examples are listed first.

We note that many of the examples below are examples of elementary groups. In most cases, a nonelementary example can be achieved by taking the Schottky product with an approprate group; cf. Proposition \ref{propositionschottkyproduct}.

The notations $\BB = \del\E^\infty\butnot\{\infty\}\equiv\ell^2(\N)$ and $\what\cdot:\Isom(\BB)\to \Isom(\H^\infty)$ will be used without comment; cf. Section \ref{subsectionparabolic}.

%
%

\section{Embedding $\R$-trees into real hyperbolic spaces}
Many of the examples in this chapter are groups acting on $\R$-trees, but it turns out that there is a natural way to convert such an action into an action on a real hyperbolic space. Specifically, we have the following:

\begin{theorem}[Generalization of {\cite[Theorem 1.1]{BIM}}]
\label{theoremBIM}
Let $X$ be a separable $\R$-tree. Then for every $\lambda > 1$ there is an embedding $\Psi_\lambda:X\to\H^\infty$ and a homomorphism $\pi_\lambda : \Isom(X)\to \Isom(\H^\infty)$ such that:
\begin{itemize}
\item[(i)] The map $\Psi_\lambda$ is $\pi_\lambda$-equivariant and extends equivariantly to a boundary map $\Psi_\lambda:\del X\to \del\H^\infty$ which is a homeomorphism onto its image.
\item[(ii)] For all $x,y\in X$ we have
\begin{equation}
\label{BIM1}
\lambda^{\dist(x,y)} = \cosh \dist(\Psi_\lambda(x),\Psi_\lambda(y)).
\end{equation}
\item[(iii)]
\begin{equation}
\label{BIM2}
\Hull_1(\Psi_\lambda(\del X)) \subset B(\Psi_\lambda(X),\cosh^{-1}(\sqrt 2)).
\end{equation}
\item[(iv)] For any set $S\subset X$, the dimension of the smallest totally geodesic subspace $[V_S]\subset\H^\infty$ containing $\Psi_\lambda(S)$ is $\#(S) - 1$. Here cardinalities are interpreted in the weak sense: if $\#(S) = \infty$, then $\dim([V_S]) = \infty$ but $S$ may be uncountable even though $[V_S]$ is separable.
\end{itemize}
\end{theorem}
\begin{proof}
Let $\VV = \{\xx\in\R^X : x_v = 0 \text{ for all but finitely many $v\in X$}\}$, and define the bilinear form $B_\QQ$ on $\VV$ via the formula
\begin{equation}
\label{BIMdef}
B_\QQ(\xx,\yy) = -\sum_{v,w\in X} \lambda^{\dist(v,w)} x_v y_w.
\end{equation}
\begin{claim}
The associated quadratic form $\QQ(\xx) = B_\QQ(\xx,\xx)$ has signature $(\omega,1)$.
\end{claim}
\begin{subproof}
It suffices to show that $\QQ\given \ee_{v_0}^\perp$ is positive definite, where $v_0\in X$ is fixed. Indeed, fix $\xx\in\ee_{v_0}^\perp\butnot\{\0\}$, and we will show that $\QQ(\xx) > 0$. Now, the set $X_0 = \{v\in X : x_v\neq 0\}\cup\{v_0\}$ is finite. It follows that the convex hull of $X_0$ can be written in the form $X(V,E,\ell)$ for some finite acyclic weighted undirected graph $(V,E,\ell)$. Consider the subspace
\[
\VV_0 = \{\xx\in\ee_{v_0}^\perp : x_v = 0 \text{ for all $v\in X\butnot V$}\} \subset \VV,
\]
which contains $\xx$. We will construct an orthogonal basis for $\VV_0$ as follows. For each edge $(v,w)\in E$, let
\[
\ff_{v,w} = \ee_v - \lambda^{\dist(v,w)} \ee_w
\]
if $w\in \geo{v_0}v$; otherwise let $\ff_{v,w} = \ff_{w,v}$. This vector has the following key property:
\begin{equation}
\label{fvw}
\begin{split}
\text{For all $v'\in X$, if $\geo vw$ intersects $\geo{v_0}{v'}$ in}\\ 
\text{at most one point, then $B_\QQ(\ff_{v,w},\ee_{v'}) = 0$.}
\end{split}
\end{equation}
(The hypothesis implies that $w\in \geo{v'}v$ and thus $\dist(v,v') = \dist(v,w) + \dist(w,v')$.) In particular, letting $v' = v_0$ we see that $\ff_{v,w}\in \ee_{v_0}^\perp$. Moreover, the tree structure of $(V,E)$ implies that for any two edges $(v_1,w_1)\not\sim (v_2,w_2)$, we have either $\#(\geo{v_1}{w_1}\cap\geo{v_0}{v_2}) \leq 1$ or $\#(\geo{v_2}{w_2}\cap\geo{v_0}{v_1}) \leq 1$; either way, \eqref{fvw} implies that $B_\QQ(\ff_{v_1,w_1},\ff_{v_2,w_2}) = 0$. Finally, $\QQ(\ff_{v,w}) = \lambda^{2\dist(v,w)} - 1 > 0$ for all $(v,w)\in E$, so $\QQ\given \VV_0$ is positive definite. Thus $\QQ(\xx) > 0$; since $\xx\in\ee_{v_0}^\perp$ was arbitrary, $\QQ\given \ee_{v_0}^\perp$ is positive definite. This concludes the proof of the claim.
\end{subproof}
It follows that for any $v\in X$, the quadratic form
\[
B_{\QQ_v}(\xx,\yy) = B_\QQ(\xx,\yy) + 2B_\QQ(\xx,\ee_v) B_\QQ(\ee_v,\yy)
\]
is positive definite. We leave it as an exercise to show that for any $v_1,v_2\in X$, the norms induced by $\QQ_{v_1}$ and $\QQ_{v_2}$ are comparable. Let $\LL$ be the completion of $\VV$ with respect to any of these norms, and (abusing notation) let $B_\QQ$ denote the unique continuous extension of $B_\QQ$ to $\LL$. Since the map $X\ni v\mapsto \ee_v\in \LL$ is continuous with respect to the norms in question, $\LL$ is separable. On the other hand, since these norms are nondegenerate, we have $\dim(\lb \ee_v : v\in S\rb) = \#(S)$ for all $S\subset X$, and in particular $\dim(\LL) = \infty$. Thus $\LL$ is isomorphic to $\LL^\infty$, so $\H := \{[\xx]\in \proj(\LL) : \QQ(\xx) < 0\}$ is isomorphic to $\H^\infty$.

We define the embedding $\Psi_\lambda:X\to\H$ via the formula $\Psi_\lambda(v) = [\ee_v]$. \eqref{BIM1} now follows immediately from \eqref{BIMdef} and \eqref{distanceinL}. In particular, we have 
\[
\dist(\Psi_\lambda(v),\Psi_\lambda(w)) \asymp_\plus \log(\lambda) \dist(v,w), 
\]
which implies that $\Psi_\lambda$ extends naturally to a boundary map $\Psi_\lambda:\del X\to \del\H^\alpha$ which is a homeomorphism onto its image. Given any $g\in\Isom(X)$, we let $\pi_\lambda(g) = [T_g]\in\Isom(\H)$, where $T_g\in \O_\R(\LL;\QQ)$ is given by the formula $T_g(\ee_v) = \ee_{g(v)}$. Then $\Psi_\lambda$ and its extension are both $\pi_\lambda$-equivariant, demonstrating condition (i).

For $S\subset X$, we have $\dim(V_S) = \dim(\lb \ee_v : v\in S\rb) = \#(S)$ as noted above, and thus $\dim([V_S]) = \dim(V_S) - 1 = \#(S) - 1$. This demonstrates (iv).

It remains to show (iii). Fix $\xi,\eta\in\del X$ and $[\zz]\in\geo{\Psi_\lambda(\xi)}{\Psi_\lambda(\eta)}$. Write $\Psi_\lambda(\xi) = [\xx]$ and $\Psi_\lambda(\eta) = [\yy]$. Since $[\xx],[\yy]\in\del\H$ and $[\zz]\in\H$, we have $\QQ(\xx) = \QQ(\yy) = 0$, and we may choose $\xx$, $\yy$, and $\zz$ to satisfy $B_\QQ(\xx,\yy) = B_\QQ(\xx,\zz) = B_\QQ(\yy,\zz) = -1$. Since $[\zz]\in\geo{[\xx]}{[\yy]}$, we have $\zz = a\xx + b\yy$ for some $a,b \geq 0$; we must have $a = b = 1$ and thus $\QQ(\zz) = -2$.

Now, since $\Psi_\lambda(w) = [\ee_w]\to \Psi_\lambda(\xi) = [\xx]$ as $w\to\xi$, there exists a function $f:X\to\R$ such that $f(w)\ee_w\to \xx$ as $w\to\xi$. Fixing $v\in\geo\xi\eta$, we have
\[
B_\QQ(\xx,\ee_v) = \lim_{w\to\xi} f(w)B_\QQ(\ee_w,\ee_v) = -\lim_{w\to\xi} f(w) \lambda^{\dist(v,w)}.
\]
In particular $B_\QQ(\xx,\ee_{v_2}) = B_\QQ(\xx,\ee_{v_1}) \lambda^{\busemann_\xi(v_2,v_1)}$, which implies that there exists $v\in\geo\xi\eta$ such that $B_\QQ(\xx,\ee_v) = 0$. Similarly, there exists a function $g:X\to\R$ such that $g(w')\ee_{w'}\to\yy$ as $w'\to\eta$; we have
\begin{align*}
B_\QQ(\yy,\ee_v) &= -\lim_{w'\to\eta} g(w') \lambda^{\dist(v,w')}\\
-1 = B_\QQ(\xx,\yy) &= \lim_{\substack{w\to\xi \\ w'\to\eta}} f(w) g(w') B_\QQ(\ee_w,\ee_{w'})\\
&= -\lim_{\substack{w\to\xi \\ w'\to\eta}} f(w) g(w') \lambda^{\dist(w,w')}
= -B_\QQ(\xx,\ee_v) B_\QQ(\yy,\ee_v)\\
B_\QQ(\xx,\ee_v) &= B_\QQ(\yy,\ee_v) = -1,
\end{align*}
so $\ee_v = \zz + \ww$ for some $\ww\in \xx^\perp\cap \yy^\perp$. Since $\QQ(\zz) = -2$ and $\QQ(\ee_v) = -1$, we have $\QQ(\ww) = 1$ and thus
\[
\cosh\dist([\ee_v],[\zz]) = \frac{|B_\QQ(\ee_v,\zz)|}{\sqrt{|\QQ(\ee_v)|\cdot|\QQ(\zz)|}} = \frac{2}{\sqrt{1\cdot 2}} = \sqrt 2.
\]
In particular $\dist([\zz],\Psi_\lambda(X)) \leq \cosh^{-1}(\sqrt 2)$. 
\end{proof}

\begin{definition}
\label{definitionBIM}
Given an $\R$-tree $X$ and a parameter $\lambda > 1$, the maps $\Psi_\lambda$ and $\pi_\lambda$ will be called the \emph{BIM embedding} and the \emph{BIM representation} with parameter $\lambda$, respectively. (Here BIM stands for M. Burger, A. Iozzi, and N. Monod, who proved the special case of Theorem \ref{theoremBIM} where $X$ is an unweighted simplicial tree.)
\end{definition}

\begin{remark}
\label{remarkBIM}
Let $X$, $\lambda$, $\Psi_\lambda$, and $\pi_\lambda$ be as in Theorem \ref{theoremBIM}. Fix $\Gamma\leq\Isom(X)$, and suppose that $\Lambda_\Gamma = \del X$. Let $G = \pi_\lambda(\Gamma)\leq\Isom(\H^\infty)$.
\begin{itemize}
\item[(i)] \eqref{BIM2} implies that if $\Gamma$ is convex-cobounded in the sense of Definition \ref{definitionCCB} below, then $G$ is convex-cobounded as well. Moreover, we have 
\begin{equation*}
\Lr(G) = \del\Psi_\lambda(\Lr(\Gamma)) ~\text{and}~ \Lur(G) = \del\Psi_\lambda(\Lur(\Gamma)).
\end{equation*}
\item[(ii)] Since $\cosh(t) \asymp_\times e^t$ for all $t\geq 0$, \eqref{BIM1} implies that
\begin{align*}
\Sigma_s(G) = \sum_{\gamma\in\Gamma} e^{-s\dogo{\pi_\lambda(\gamma)}} &\asymp_\times \sum_{\gamma\in\Gamma} \cosh^{-s}(\dogo{\pi_\lambda(\gamma)})\\ 
&= \sum_{\gamma\in\Gamma} \lambda^{-s\dogo{\gamma}} = \Sigma_{s\log(\lambda)}(\Gamma)
\end{align*}
for all $s\geq 0$. In particular $\delta_G = \delta_\Gamma/\log(\lambda)$. A similar argument shows that $\w\delta_G = \w\delta_\Gamma/\log(\lambda)$, which implies that $G$ is Poincar\'e regular if and only if $\Gamma$ is.
\item[(iii)] $G$ is strongly discrete (resp. COT-discrete) if and only if $\Gamma$ is strongly discrete (resp. COT-discrete). However, this fails for weak discreteness; cf. Example \ref{exampleAutTBIM} below.
\end{itemize}
\end{remark}
\begin{proof}[Proof of (iii)]
The difficult part is showing that if $G$ is COT-discrete, then $\Gamma$ is as well. Suppose that $\Gamma$ is not COT-discrete. Then there exists a sequence $\Gamma\ni\gamma_n\to\id$ in the compact-open topology. Let $g_n = \pi_\lambda(\gamma_n)\in G \leq \Isom(\H^\infty) \equiv \O(\LL)$. Then the set
\[
\{\xx\in\LL : g_n(\xx)\to \xx\}
\]
contains $\Psi_\lambda(X)$. On the other hand, since the sequence $(g_n)_1^\infty$ is equicontinuous (Lemma \ref{lemmaoperatornorm}), this set is a closed linear subspace of $\LL$. Clearly, the only such subspace which contains $\Psi_\lambda(X)$ is $\LL$. Thus $g_n(x)\to x$ for all $x\in\H^\infty$, and so $g_n\to\id$ in the compact-open topology. Thus $G$ is not COT-discrete.
\end{proof}

We begin our list of examples with the following counterexample to an infinite-dimensional analogue of Margulis's lemma suggested in Remark \ref{remarkmargulislemma}:

\begin{example}
\label{examplemargulis}
Let $\Gamma = \F_2(\Z) = \lb \gamma_1,\gamma_2\rb$, and let $X$ be the Cayley graph of $\Gamma$. Let $\Phi:\Gamma\to \Isom(X)$ be the natural action. Then $H := \Phi(\Gamma)$ is nonelementary and strongly discrete. For each $\lambda > 1$, the image of $H$ under the BIM representation $\pi_\lambda$ is a nonelementary strongly discrete subgroup $G = \pi_\lambda(H)\leq \Isom(\H^\infty)$ generated by the elements $g_1 = \pi_\lambda\Phi(\gamma_1)$, $g_2 = \pi_\lambda\Phi(\gamma_2)$. But
\[
\cosh\dogo{g_i} = \lambda^{\dist(e,\gamma_i)} = \lambda,
\]
so by an appropriate choice of $\lambda$, $\dogo{g_i}$ can be made arbitrarily small. So for arbitrarily small $\epsilon$, we can find a free group $G\leq\Isom(\H^\infty)$ such that $G_\epsilon(\zero) = G$ is nonelementary. This provides a counterexample to a hypothetical infinite-dimensional analogue of Margulis's lemma, namely, the claim that there exists $\epsilon > 0$ such that for every strongly discrete $G\leq\Isom(\H^\infty)$, $G_\epsilon(\zero)$ is elementary.
\end{example}

\begin{remark}
\label{remarklengthspectrum}
If $\H$ is a finite-dimensional algebraic hyperbolic space and $G\leq\Isom(\H)$ is nonelementary, then a theorem of I. Kim \cite{Kim} states that the \emph{length spectrum} of $G$
\[
\L = \{\log g'(g_-) : g\in G \text{ is loxodromic}\}
\]
is not contained in any discrete subgroup of $\R$. Example \ref{examplemargulis} shows that this result does not generalize to infinite-dimensional algebraic hyperbolic spaces. Indeed, if $G\leq\Isom(\H^\infty)$ is as in Example \ref{examplemargulis} and if $g = \pi_\lambda(\gamma)\in G$, then \eqref{BIM1} implies that
\begin{align*}
\log g'(g_-) = \lim_{n\to\infty} \frac1n \dogo{g^n} &= \lim_{n\to\infty}\cosh^{-1}\lambda^{\dogo{\gamma^n}}\\
&= \log(\lambda)\lim_{n\to\infty}\frac1n \dogo{\gamma^n}\\ 
&= \log(\lambda)\log \gamma'(\gamma_-),
\end{align*}
demonstrating that $\L$ is contained in the discrete subgroup $\log(\lambda)\Z \leq\R$.
\end{remark}

\section{Strongly discrete groups with infinite Poincar\'e exponent}
We have already seen two examples of strongly discrete groups with infinite Poincar\'e exponent, namely the Edelstein-type Example \ref{exampleparabolicinfinite}, and the parabolic torsion Example \ref{exampleparabolictorsion}. We give three more examples here.

\begin{example}[A nonelementary strongly discrete group $G$ acting on a proper $\R$-tree $X$ and satisfying $\delta_G = \infty$]
\label{exampleinfinitepoincare}
Let $Y = \Rplus$, let $P = \N$, and for each $p = n\in P$ let
\[
\Gamma_p = \Z/n!\Z
\]
(or more generally, let $\Gamma_p$ be any sufficiently large finite group). Let $(X,G)$ be the geometric product of $Y$ with $(\Gamma_p)_{p\in P}$, as defined below in Example \ref{examplegeometricproducts}. By Proposition \ref{propositionschottkyRtree}, $X$ is proper, and $G = \lb G_p\rb_{p\in P}$ is a global weakly separated Schottky product. So by Corollary \ref{corollaryschottkySD}, $G$ is strongly discrete. Clearly, $G$ is nonelementary. Finally, $\delta_G = \infty$ because for all $s\geq 0$,
\[
\Sigma_s(G) \geq \sum_{p\in P} \sum_{g\in \Gamma_p\butnot\{e\}} e^{-s\dogo g} = \sum_{p\in P} \#(\Gamma_p\butnot\{e\}) e^{-2s\dox p} = \sum_{n\in\N} (n! - 1) e^{-2ns} = \infty.
\]
\end{example}

Applying a BIM representation gives:

\begin{example}[A nonelementary strongly discrete convex-cobounded group acting on $\H^\infty$ and satisfying $\delta = \infty$]
\label{exampleinfinitepoincareBIM}
Cf. Remark \ref{remarkBIM} and the example above.
\end{example}

\begin{example}[A parabolic strongly discrete group $G$ acting on $\H^\infty$ and satisfying $\delta_G = \infty$]
\label{examplehaagerup}
Since $\F_2(\Z)$ has the Haagerup property (Remark \ref{remarkhaagerup}), there is a homomorphism $\Phi:\F_2(\Z)\to\Isom(\BB)$ whose image $G = \Phi(\F_2(\Z))$ is strongly discrete. However, $\what G$ must have infinite Poincar\'e exponent by Corollary \ref{corollaryfiniteimpliesnilpotent}.
\end{example}

\section{Moderately discrete groups which are not strongly discrete}
We have already seen one example of a moderately discrete group which is not strongly discrete, namely the Edelstein-type Example \ref{examplevalette} (parabolic acting on $\H^\infty$). We give three more examples here, and we will give one more example in Section \ref{subsectionpoincareirregular}, namely Example \ref{exampleAutTparttwo}. All five examples are are also examples of properly discontinuous actions, so they also demonstrate that proper discontinuity does not imply strong discreteness. (The fact that moderate discreteness (or even strong discreteness) does not imply proper discontinuity can be seen e.g. from Examples \ref{exampleparabolictorsion}, \ref{exampleinfinitepoincare}, and \ref{exampleinfinitepoincareBIM}, all of which are generated by torsion elements.)


\begin{example}[A parabolic group which acts properly discontinuously on $\H^\infty$ but is not strongly discrete]
\label{exampletranslations}
Let $\Z^\infty\subset\BB = \ell^2(\N)$ denote the set of all infinite sequences in $\Z$ with only finitely many nonzero entries. Let
\[
G := \{ \xx \mapsto \xx + \nn : \nn\in\Z^\infty \} \subset \Isom(\BB).
\]
Then $G$ acts properly discontinuously, since $\|(\xx + \nn) - \xx\| \geq 1$ for all $\xx\in\BB$ and $\nn\in\Z^\infty\butnot\{\0\}$. On the other hand, $G$ is not strongly discrete since $\|\nn\| = 1$ for infinitely many $\nn\in\Z^\infty$. By Observation \ref{observationisomE}, these properties also hold for the Poincar\'e extension $\what G\leq\Isom(\H^\infty)$.
\end{example}

\begin{example}[A nonelementary group $G$ which acts properly discontinuously on a separable $\R$-tree $X$ but is not strongly discrete]
\label{examplenotstronglydiscrete}
Let $X$ be the Cayley graph of $\Gamma = \F_\infty(\Z)$ with respect to its standard generators, and let $\Phi:\Gamma\to\Isom(X)$ be the natural action. Then $G = \Phi(\Gamma)$ acts properly discontinuously on $X$. On the other hand, since by definition each generator $g\in G$ satisfies $\dogo g = 1$, $G$ is not strongly discrete.
\end{example}

Applying a BIM representation gives:

\begin{example}[A nonelementary group which acts properly discontinuously on $\H^\infty$ but is not strongly discrete]
\label{examplenotstronglydiscreteBIM}
Let $X$ and $G$ be as in Example \ref{examplenotstronglydiscrete}. Fix $\lambda > 1$ large to be determined, and let $\pi_\lambda:\Isom(X)\to \Isom(\H^\infty)$ be the corresponding BIM representation. By Remark \ref{remarkBIM}, the group $\pi_\lambda(G)$ is a nonelementary group which acts isometrically on $\H^\infty$ but is not strongly discrete. To complete the proof, we must show that $\pi_\lambda(G)$ acts properly discontinuously. By Proposition \ref{propositionSproductMDWDPD}, it suffices to show that $G = \prod_1^\infty \pi_\lambda(\gamma_i)^\Z$ is a global strongly separated Schottky group. And indeed, if we denote the generators of $\Gamma = \F_\infty(\Z)$ by $\gamma_i$ ($i\in\N$), and if we consider the balls $U_i^\pm = B(\Psi_\lambda((\gamma_i)_\pm),1/2)$ (taken with respect to the Euclidean metric), and if $\lambda$ is sufficiently large, then the sets $U_i = U_i^+\cup U_i^-$ form a global strongly separated Schottky system for $G$.
\end{example}

\begin{remark}
\label{remarkMDuncountable}
The groups of Examples \ref{examplenotstronglydiscrete}-\ref{examplenotstronglydiscreteBIM} can be easily modified to make the group $G$ uncountable at the cost of separability; let $X$ be the Cayley graph of $\F_{\#(\R)}(\Z)$ in Example \ref{examplenotstronglydiscrete}, and applying (a modification of) Theorem \ref{theoremBIM} gives an action on $\H^{\#(\R)}$.
\end{remark}

\begin{remark*}
By Proposition \ref{propositionpoincareregular}, the groups of Examples \ref{examplenotstronglydiscrete}-\ref{examplenotstronglydiscreteBIM} are all Poincar\'e regular and therefore satisfy $\HD(\Lur) = \infty$.
\end{remark*}

\section{Poincar\'e irregular groups}
\label{subsectionpoincareirregular}

We give six examples of Poincar\'e irregular groups, providing counterexamples to many conceivable generalizations of Proposition \ref{propositionpoincareregular}.

\begin{example}[A Poincar\'e irregular nonelementary group $G$ acting on a proper $\R$-tree $X$ which is weakly discrete but not COT-discrete]
\label{exampleAutT}
Let $X$ be the Cayley graph of $V = \F_2(\Z)$ (equivalently, let $X$ be the unique $3$-regular unweighted simplicial tree), and let $G = \Isom(X)$. Since $\#(\Stab(G;e)) = \infty$, $G$ is not strongly discrete, so by Proposition \ref{propositionparametricdiscreteness}, $G$ is also not COT-discrete. (The fact that $G$ is not COTD can also be deduced from Proposition \ref{propositionpoincareregular}, since we will soon show that $G$ is Poincar\'e irregular.)

On the other hand, suppose $x\in X$. Then either $x\in V$, or $x = ((v_x,w_x),t_x)$ for some $(v_x,w_x)\in E$ and $t_x\in (0,1)$. In the first case, we observe that $G(x) = V$, while in the second we observe that
\[
G(x) = \{((v,w),t_x) : (v,w)\in E\}.
\]
In either case $x$ is not an accumulation point of $G(x)$. Thus $G$ is weakly discrete.

To show that $G$ is Poincar\'e irregular, we first observe that $\delta = \infty$ since $G$ is not strongly discrete. On the other hand, Proposition \ref{propositionbasicmodified}(iv) can be used to compute that $\w\delta = \log_b(2)$. (Alternatively, one may use Theorem \ref{theorembishopjonesmodified} together with the fact that $\HD(\del X) = \log_b(2)$.)
\end{example}


\begin{remark*}
The group $G$ in Example \ref{exampleAutT} is uncountable. However, if $G$ is replaced by a countable dense subgroup (cf. Remark \ref{remarkpolish}) then the conclusions stated above will not be affected. This remark applies also to Examples \ref{exampleAutTBIM} and \ref{exampleAutTparttwo} below.
\end{remark*}

Applying a BIM representation to the group of Example \ref{exampleAutT} yields:

\begin{example}[A Poincar\'e irregular nonelementary group acting irreducibly on $\H^\infty$ which is UOT-discrete but not COT-discrete]
\label{exampleAutTBIM}
Let $G\leq\Isom(X)$ be as in Example \ref{exampleAutT} and let $\pi_\lambda:\Isom(X)\to \Isom(\H^\infty) \equiv \O(\LL)$ be a BIM representation. Remark \ref{remarkBIM} shows that the group $\pi_\lambda(G)$ is Poincar\'e irregular and is not COT-discrete. Note that it follows from either Proposition \ref{propositionparametricdiscreteness}(ii) or Proposition \ref{propositionpoincareregular} that $\pi_\lambda(G)$ is not weakly discrete, despite $G$ being weakly discrete.

To complete the proof, we must show that $\pi_\lambda(G)$ is UOT-discrete. Let $\Psi_\lambda:X\to\H^\infty\subset\LL$ be the BIM embedding corresponding to the BIM representation $\pi_\lambda$, and write $\zz = \Psi_\lambda(\zero)$; without loss of generality we may assume $\zz = (1,\0)$, so that $\QQ(\xx) = \|\xx\|^2$ for all $\xx\in\zz^\perp$.

Now fix $T = \pi_\lambda(g)\in \pi_\lambda(G)\butnot\{\id\}$, and we will show that $\|T - I\|\geq \min(\sqrt 2,\lambda - 1) > 0$. We consider two cases. If $g(\zero)\neq\zero$, then $\dogo g\geq 1$, which implies that $|B_\QQ(\zz,T\zz)| \geq \lambda$ and thus that $\|T\zz - \zz\| \geq |B_\QQ(\zz,T\zz - \zz)| \geq \lambda - 1$. So suppose $g(\zero) = \zero$. Since $g\neq\id$, we have $g(x) \neq x$ for some $x\in V$; choose such an $x$ so as to minimize $\dox x$. Letting $y = \geo\zero x_{\dox x - 1}$, the minimality of $\dox x$ implies that $g(y) = y$ (cf. Figure \ref{figureexampleAutTBIM}).

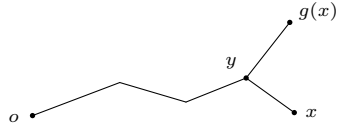
\begin{figure}
\begin{center}
\begin{tikzpicture}[line cap=round,line join=round,>=triangle 45,x=1.0cm,y=1.0cm]
\clip(-2.35,-0.03) rectangle (3.55,1.8);
\draw (-1.16,0.24)-- (0,0.68);
\draw (0,0.68)-- (0.88,0.42);
\draw (0.88,0.42)-- (1.68,0.74);
\draw (1.68,0.74)-- (2.26,1.48);
\draw (1.68,0.74)-- (2.32,0.28);
\begin{scriptsize}

\fill [color=black] (1.68,0.74) circle (1pt);
\draw[color=black] (1.48,0.94) node {$y$};

\fill [color=black] (-1.16,0.24) circle (1pt);
\draw[color=black] (-1.41,0.21) node {$\zero$};
\fill [color=black] (2.26,1.48) circle (1pt);
\draw[color=black] (2.64,1.63) node {$g(x)$};
\fill [color=black] (2.32,0.28) circle (1pt);
\draw[color=black] (2.55,0.27) node {$x$};
\end{scriptsize}
\end{tikzpicture}
\caption[Geometry of automorphisms of a simplicial tree]{The point $y$ is the center of the triangle $\Delta(\zero,x,g(x))$. Both $\zero$ and $y$ are fixed by $g$. Intuitively, this means that $g$ (really, $\pi_\lambda(g)$) must have a significant rotational component in order to ``swing up'' the point $x$ to the point $g(x)$.}
\label{figureexampleAutTBIM}
\end{center}
\end{figure}

Let $\xx = \Phi_\lambda(x)$ and $\yy = \Phi_\lambda(y)$, so that $T\yy = \yy$ but $T\xx\neq \xx$. Let $\ww_1 = \xx - \lambda \yy$ and $\ww_2 = T\ww_1 = T\xx - \lambda\yy$. An easy computation based on \eqref{BIM1} and \eqref{distanceinL} gives $B_\QQ(\zz,\ww_1) = B_\QQ(\zz,\ww_2) = B_\QQ(\ww_1,\ww_2) = 0$ (cf. \eqref{fvw}). It follows that
\begin{align*}
\|(T - I)\ww_1\| = \|\ww_2 - \ww_1\| &= \sqrt{\QQ(\ww_2 - \ww_1)}\\ 
&= \sqrt{\QQ(\ww_2) + \QQ(\ww_1)}\\ 
&= \sqrt{2\QQ(\ww_1)} = \sqrt 2 \|\ww_1\|,
\end{align*}
and thus $\|T - I\| \geq \sqrt 2$.
\end{example}

\begin{remark}
\label{remarkfocalfractal}
Let $G,\pi_\lambda$ be as above and fix $\xi\in \del X$. Then $\pi_\lambda(G_\xi)$ is a focal group acting irreducibly on $\H^\infty$ whose limit set is totally disconnected. This contrasts with the finite-dimensional situation, where any nondiscrete group (and thus any focal group) acting irreducibly on $\H^d$ is of the first kind \cite[Theorem 2]{Greenberg}.
\end{remark}


\begin{example}[A Poincar\'e irregular nonelementary group $G'$ acting properly discontinuously on a hyperbolic metric space $X'$]
\label{exampleAutTparttwo}
Let $G$ be the group described in Example \ref{exampleAutT}. Let $X' = G$ and let
\[
\dist'(g,h) := \begin{cases}
1\vee\dist(g(\zero),h(\zero)) & g\neq h\\
0 & g = h
\end{cases}.
\]
Since the orbit map $X'\ni g\to g(\zero)\in X$ is a quasi-isometric embedding, $(X',\dist')$ is a hyperbolic metric space. The left action of $G$ on $X'$ is isometric and properly discontinuous. Denote its image in $\Isom(X')$ by $G'$. Clearly $\delta_{G'} = \delta_G$ and $\w\delta_{G'} = \w\delta_G$ (the Poincar\'e exponent and modified Poincar\'e exponent do not depend on whether $G$ is acting on $X$ or on $X'$), so $G'$ is Poincar\'e irregular.
\end{example}

The next set of examples have a somewhat different flavor.

\begin{example}[A Poincar\'e irregular group $G$ acting on $\H^d$]
\label{examplenondiscrete}
Fix $2\leq d < \infty$, and let $G$ be any nondiscrete subgroup of $\Isom(\H^d)$. Then 
\[
\w\delta_G = \HD(\Lr) \leq \HD(\del\H^d) = d - 1.
\] 
On the other hand, since $G$ is not strongly discrete we have $\delta_G = \infty$. Thus $G$ is Poincar\'e irregular.
\end{example}

In Example \ref{examplenondiscrete}, $G$ could be a Lie subgroup with nontrivial connected component (e.g. $G = \Isom(\H^d)$, but this is not the only possibility - $G$ can even be finitely generated, as we now show:

\begin{lemma}
\label{lemmafreenotschottky}
Let $H$ be a connected algebraic group which contains a copy of the free group $\F_2(\Z)$. Then there exist $g_1,g_2\in H$ such that $G := \lb g_1,g_2\rb$ is a nondiscrete group isomorphic to $\F_2(\Z)$.
\end{lemma}
By Lemma \ref{lemmapingpong}, the group $G$ cannot be a Schottky product - thus this lemma provides an example of a free product which is not a Schottky product.
\begin{proof}
An orders-of-magnitude argument shows that there exists $\epsilon > 0$ such that for any $h_1,h_2\in H$ with $\dist(\id,h_i)\leq \epsilon$, we have $$\dist(\id,[h_1,h_2]) \leq \frac12\max_i \dist(\id,h_i),$$ where $[h_1,h_2]$ denotes the commutator of $h_1$ and $h_2$. Thus for any $g_1,g_2\in H$ such that $\dist(\id,g_i)\leq \epsilon$, letting
\[
h_1 = g_1, \;\; h_2 = g_2, \;\; h_{n + 2} = [h_n,h_{n + 1}]
\]
gives $h_n\to\id$. But the elements $h_n$ are the images of nontrivial words in the free group $\F_2(\Z)$ under the natural homomorphism, so if this homomorphism is injective then $G$ is not discrete. For each element $\gg\in \F_2(\Z)$, the set of homomorphisms $\pi:\F_2(\Z)\to H$ such that $\pi(\gg) = \id$ is a proper algebraic subset of the set of all homomorphisms, and therefore has measure zero. Thus for typical $g_1,g_2$ satisfying $\dist(\id,g_i)\leq \epsilon$, $G$ is a nondiscrete free group.
\end{proof}

Instead of a Lie subgroup of $\Isom(\H^d)$, we could also take a locally compact subgroup of $\Isom(\H^\infty)$; there are many interesting examples of such subgroups. In particular, one such example is given by the following theorem:

\begin{theorem}[Monod--Py representation theorem, {\cite[Theorems B and C]{MonodPy}}]
\label{theoremmonodpy}
For any $d\in\N$ and $0 < t < 1$, there exist an irreducible representation $\rho_t:\Isom(\H^d)\to\Isom(\H^\infty)$ and a $\rho_t$-equivariant embedding $f_t:\bord\H^d\to\bord\H^\infty$ such that
\begin{equation}
\label{monodpy}
\dist(f_t(x),f_t(y)) \asymp_\plus t \dist(x,y) \text{ for all $x,y\in\H^d$.}
\end{equation}
The pair $(\rho_t,f_t)$ is unique up to conjugacy.
\end{theorem}


\begin{example}[A Poincar\'e irregular nonelementary group $G$ acting irreducibly on $\H^\infty$]
\label{examplemonodpy}
Fix $d\in\N$ and $0 < t < 1$, and let $\rho_t$, $f_t$ be as in Theorem \ref{theoremmonodpy}. Let $\Gamma = \Isom(\H^d)$, and let $G = \rho_t(\Gamma)$. As $G$ is locally compact, the modified Poincar\'e exponent of $G$ can be computed using Definition \ref{definitionmodified1}:
\begin{align*}
\w\delta_G &= \inf\left\{s\geq 0 : \int_G e^{-s\dogo g}\;\dee g < \infty \right\}\\
&= \inf\left\{s\geq 0 : \int_\Gamma e^{-s\dogo{\rho_t(\gamma)}}\;\dee \gamma < \infty \right\}\\
&= \inf\left\{s\geq 0 : \int_\Gamma e^{-st\dogo{\gamma}}\;\dee \gamma < \infty \right\}\\
&= \frac{\w\delta_\Gamma}{t} = \frac{\HD(\Lambda_\Gamma)}{t} = \frac{d - 1}{t}\cdot
\end{align*}
On the other hand, since $G$ is convex-cobounded by \cite[Theorem D]{MonodPy}, Theorem \ref{theoremCCB} shows that $\Lambda_G = \Lr(G) = \Lur(G)$. (It may be verified that the strong discreteness assumption is not needed for those directions.) Combining with Theorem \ref{theorembishopjonesmodified}, we have
\[
\HD(\Lambda_G) = \HD(\Lr(G)) = \HD(\Lur(G)) = \frac{d - 1}{t} > d - 1 = \HD(\Lambda_\Gamma).
\]
In particular, it follows that the map $f_t:\Lambda_\Gamma\to\Lambda_G$ cannot be smooth or even Lipschitz. This contrasts with the smoothness of $f_t$ in the interior (see \cite[Theorem C(2)]{MonodPy}).
\end{example}

\begin{remark*}
The Hausdorff dimension of $\Lambda_G$ may also be computed directly from the formulas \eqref{monodpy} and \eqref{distanceasymptotic}, which imply that the map $f_t\given\Lambda_\Gamma$ and its inverse are H\"older continuous of exponents $t$ and $1/t$, respectively. However, the computation above gives a nice application of the Poincar\'e irregular case of Theorem \ref{theorembishopjonesmodified}.
\end{remark*}

In Examples \ref{examplenondiscrete} and \ref{examplemonodpy}, the group $G$ does not satisfy any of the discreteness conditions discussed in Chapter \ref{sectiondiscreteness}. Our next example satisfies a weak discreteness condition:

\begin{example}[A Poincar\'e irregular nonelementary COT-discrete group $G$ acting reducibly on $\H^\infty$ which is not weakly discrete]
\label{examplepoincareextension}
Let $\Gamma = \F_2(\Z)$ and let $\iota_1:\Gamma\to\Isom(\H^d) \equiv \O(\LL^{d + 1})$ be an injective homomorphism whose image is a nondiscrete group; this is possible by Lemma \ref{lemmafreenotschottky}. Define $\iota_2:\Gamma\to\O(\HH^\Gamma)$ by letting
\[
\iota_2(\gamma)[\ee_\delta] = \ee_{\gamma\delta}.
\]
Note that $\iota_2(\Gamma)$ is COT-discrete, since $\|\iota_2(\gamma)\ee_e - \ee_e\| = \sqrt 2$ for all $\gamma\in\Gamma\butnot\{e\}$.

The direct sum $\iota := \iota_1\oplus\iota_2:\Gamma\to\O(\LL^{d + 1}\times\HH^\Gamma)$ is an isometric action of $\Gamma$ on $\H^{\Gamma\dot\cup\{1,\ldots,d\}} \equiv \H^\infty$. Let $G = \iota(\Gamma)$. Since $\iota_1(\Gamma)$ is the restriction of $G$ to the invariant totally geodesic subspace $\H^d$, we have $\delta_G = \delta_{\iota_1(\Gamma)} = \infty$ and $\w\delta_G = \w\delta_{\iota_1(\Gamma)} < \infty$, so $G$ is Poincar\'e irregular. On the other hand, $G$ is COT-discrete because $\iota_2(\Gamma)$ is. Finally, the fact that $G$ is not weakly discrete can be seen from either Observation \ref{observationrestrictions} or Proposition \ref{propositionpoincareregular}.
\end{example}

\section{Miscellaneous counterexamples}
Our remaining examples include a COTD group which is not WD and a WD group which is not MD.

\begin{example}[A nonelementary COT-discrete group $G$ which acts irreducibly on $\H^\infty$ and satisfies $\delta_G = \w\delta_G = \infty$ but which is not weakly discrete]
\label{examplepoincareextensiontwo}
Let $G_1\leq\Isom(\H^\infty)$ be as in Example \ref{examplepoincareextension}, and let $g$ be a loxodromic isometry whose fixed points are $g_\pm = [\ee_0 \pm \ee_e]\in\del\H^{\Gamma\dot\cup\{1,\ldots,d\}}\subset \proj\LL^{\Gamma\dot\cup\{0,\ldots,d\}}$. Then for $n$ sufficiently large, the product $G = \lb G_1,(g^n)^\Z\rb$ is a global strongly separated Schottky product. By Lemma \ref{lemmapingpong}, $G$ is COT-discrete. Since $G$ contains $G_1$, $G$ is not weakly discrete.

The fact that $\w\delta_G = \infty$ follows from either Proposition \ref{propositionSproductexponent}(iii) or Proposition \ref{propositionpoincareregular}. So the only thing left to show is that $G$ acts irreducibly. We assume that the original group $\iota_1(\Gamma)$ acts irreducibly. Then if $[V]\subset\H^\infty$ is a $G$-invariant totally geodesic subspace containing the limit set of $G$, then $\LL^{d + 1}\subset V$ and so $V = \LL^{d + 1}\oplus V_2$ for some $V_2\subset \HH^\Gamma$. But $[\ee_0 + \ee_e]\in\Lambda_G$, so $\ee_e\in V_2$. The $G$-invariance of $[V]$ implies that $V_2$ is $\iota_2(\Gamma)$-invariant, and thus that $V_2 = \HH^\Gamma$ and so $[V] = \H^\infty$.
\end{example}
\begin{remark*}
Example \ref{examplepoincareextensiontwo} gives a good example of how Theorem \ref{theorembishopjonesmodified} gives interesting information even when $\w\delta_G = \infty$. Namely, in this example Theorem \ref{theorembishopjonesmodified} tells us that $\HD(\Lr) = \HD(\Lur) = \infty$, which is not at all obvious simply from looking at the group.
\end{remark*}

\begin{example}[An elliptic group $G$ acting on $\H^\infty$ which is weakly discrete but not moderately discrete]
\label{exampleMDWD}
Let $\HH = \ell^2(\Z)$, and let $T\in\O(\HH)$ be the shift map $T(\xx) = (x_{n + 1})_{n = 1}^\infty$. Let $G$ be the cyclic group $G = T^\Z\leq\O(\HH) \leq \Isom(\B^\infty)$. Since $g(\0) = \0$ for all $g\in G$, $G$ is not moderately discrete. On the other hand, fix $\xx\in\HH\butnot\{\0\}$. Then $T^n(\xx) \to \0$ weakly as $n\to \pm \infty$, so $\#\{n\in\Z : \|T^n(\xx) - \xx\| \leq \|\xx\|/2\} < \infty$. Thus $G$ is weakly discrete.
\end{example}

%
%

\chapter{$\R$-trees and their isometry groups} \label{sectionRtrees}

In this chapter we describe various ways to construct $\R$-trees which admit isometric actions. Section \ref{subsectionconeconstruction} describes the cone construction, in which one starts with an ultrametric space $(Z,\Dist)$ and builds an $\R$-tree $X$ whose Gromov boundary contains a point $\infty$ such that $(Z,\Dist) = (\del X\butnot\{\infty\},\Dist_{\infty,\zero})$. Sections \ref{subsectiongraphtheory} and \ref{subsectionRtreegeometric} are preliminaries for Section \ref{subsectionstapling}, which describes the ``stapling method'' in which one starts with a collection of $\R$-trees $(X_v)_{v\in V}$ and staples them together to get another $\R$-tree. We give three very general examples of the stapling method in which the resulting $\R$-tree admits a natural isometric action.

We recall that whenever we have an example of an $\R$-tree $X$ with an isometric action $\Gamma\leq\Isom(X)$, then we can get a corresponding example of a group of isometries of $\H^\infty$ by applying a BIM representation (Theorem \ref{theoremBIM}). Thus, the examples of this chapter contribute to our goal of understanding the behavior of isometry groups acting on $\H^\infty$.

\section{Construction of $\R$-trees by the cone method}
\label{subsectionconeconstruction}
The construction of hyperbolic metric spaces by cone methods has a long history; see e.g. \cite[1.8.A.(b)]{Gromov3}, \cite{TrotsenkoVaisala}, \cite[\67]{BonkSchramm}. The construction below does not appear to be equivalent to any of those existing in the literature, although our formula \eqref{coneconstruction} is similar to \cite[7.1]{BonkSchramm} (with the difference that their $+$ sign is replaced by a $\vee$; this change only works because we assume that $Z$ is ultrametric).

Let $(Z,\Dist)$ be a complete ultrametric space. Define an equivalence relation on $Z\times (0,\infty)$ by letting $(z_1,r_1) \sim (z_2,r_2)$ if $\dist(z_1,z_2)\leq r_1 = r_2$, and denote the equivalence class of $(z,r)$ by $\lb z,r\rb$. Let $X = Z\times (0,\infty)/\sim$, and define a distance function on $X$:
\begin{equation}
\label{coneconstruction}
\dist\big(\lb z_1,r_1\rb,\lb z_2,r_2\rb\big) = \log\left(\frac{r_1^2\vee r_2^2\vee \Dist^2(z_1,z_2)}{r_1 r_2}\right)
\end{equation}
(cf. Corollary \ref{corollarysullivansformula}). We call $(X,\dist)$ the \emph{cone} of $(Z,\Dist)$. Note that
\begin{equation}
\label{coneconstructiongromov}
\big\lb \lb z_1,r_1\rb | \lb z_2,r_2\rb\big\rb_{\lb z_0,r_0\rb} = \log\left(\frac{(r_0\vee r_1\vee \Dist(z_0,z_1))(r_0\vee r_2\vee \Dist(z_0,z_2))}{r_0(r_1\vee r_2\vee \Dist(z_1,z_2))}\right).
\end{equation}
\begin{theorem}
\label{theoremconeconstruction}
The cone $(X,\dist)$ is an $\R$-tree. Moreover, there exists a map $\iota:Z\to\del X$ such that $\del X\butnot \iota(Z)$ consists of one point, $\infty$, and such that $\Dist = \Dist_{\infty,\zero}\circ\iota$, where $\zero = \lb z_0,1\rb$ for any $z_0\in Z$.
\end{theorem}
\begin{proof}
Fix $x_i = \lb z_i,r_i\rb \in X$, $i = 1,2$, let $R = r_1\vee r_2\vee \Dist(z_1,z_2)$, and let $\gamma_i:[\log(r_i),\log(R)]\to X$ be defined by $\gamma_i(t) = \lb z_i,e^t\rb$. Then $\gamma_i$ parameterizes a geodesic connecting $x_i$ and $\lb z_i,R\rb$. Since $(z_1,R) \sim (z_2,R)$, the geodesics $\gamma_i$ can be concatenated, and their concatenation is a geodesic connecting $x_1$ and $x_2$. It can be verified that the collection of such geodesics satisfies the conditions of Lemma \ref{lemmaRtreeequivalent}. Thus $(X,\dist)$ is an $\R$-tree. (For an alternative proof that $(X,\dist)$ is an $\R$-tree, see Example \ref{exampleconeconstruction} below.)

Fix $z_0\in Z$. For all $z_1,z_2\in Z$ and $R > 0$, \eqref{coneconstructiongromov} gives
\[
\lim_{r_1,r_2\to 0} \big\lb \lb z_1,r_1\rb \big| \lb z_2,r_2\rb \big\rb_{\lb z_0,R\rb} = \sum_{i = 1}^2 \log(\sqrt R)\vee\log\left(\frac{\Dist(z_0,z_i)}{\sqrt R}\right) - \log\Dist(z_1,z_2).
\]
In particular, if $z_1 = z_2 = z$, then this shows that the sequence $\big(\lb z,1/n\rb\big)_1^\infty$ is a Gromov sequence. Let $\iota(z) = \left[\big(\lb z,1/n\rb\big)_1^\infty\right]$. Similarly, the sequence $\big(\lb z_0,n\rb\big)_1^\infty$ is a Gromov sequence; let $\infty = \left[\big(\lb z_0,1/n\rb\big)_1^\infty\right]$. Then Lemma \ref{lemmanearcontinuity} gives
\[
\lb \iota(z_1) | \iota(z_2)\rb_{(z_0,R)} = \sum_{i = 1}^2 \log(R)\vee\log\left(\frac{\Dist(z_0,z_i)}{R}\right) - \log\Dist(z_1,z_2)
\]
and thus
\[
-\log\Dist_{\infty,\zero}(\iota(z_1),\iota(z_2)) = \lim_{R\to\infty} \Big[ \lb \iota(z_1) | \iota(z_2)\rb_{(z_0,R)} - \log(R) \Big] = -\log\Dist(z_1,z_2),
\]
i.e. $\Dist_{\infty,\zero}\equiv \Dist$.

To complete the proof we need to show that $\del X = \iota(Z)\cup\{\infty\}$. Indeed, fix $\xi = \left[\big(\lb z_n,r_n\rb\big)_1^\infty\right]\in\del X$. Without loss of generality suppose that $r_n\to r\in [0,\infty]$ and $\Dist(z_0,z_n) \to R\in[0,\infty]$. If $r = \infty$ or $R = \infty$, then it follows from \eqref{coneconstructiongromov} that $\lb \lb z_n,r_n\rb | \infty\rb_{\lb z_0,1\rb} \to \infty$, i.e. $\xi = \infty$. Otherwise, it follows from \eqref{coneconstructiongromov} that
\begin{align*}
\infty &= \lim_{n,m \to\infty} \big\lb \lb z_n,r_n\rb | \lb z_m,r_m\rb \big\rb_{\lb z_0,1\rb}\\ 
&= 2\log(1\vee r\vee R) - \log \left( \lim_{n,m\to\infty} r_n\vee r_m\vee \Dist(z_n,z_m) \right),
\end{align*}
which implies that $r_n\vee r_m\vee \Dist(z_n,z_m) \tendsto{n,m} 0$, i.e. $r_n\to 0$ and $(z_n)_1^\infty$ is a Cauchy sequence. Since $Z$ is complete we can find a limit point $z_n\to z\in Z$. Then \eqref{coneconstructiongromov} shows that $\xi = \iota(z)$.
\end{proof}

\begin{corollary}
\label{corollaryultrametricRtree}
Every ultrametric space can be isometrically embedded into an $\R$-tree.
\end{corollary}
\begin{proof}
Let $(Y,\dist)$ be an ultrametric space, and without loss of generality suppose that $Y$ is complete. Let $Z = Y$, and let $\Dist(z_1,z_2) = e^{(1/2)\dist(z_1,z_2)}$. Then $(Z,\Dist)$ is a complete ultrametric space. Let $(X,\dist)$ be the cone of $(Z,\Dist)$; by Theorem \ref{theoremconeconstruction}, $X$ is an $\R$-tree. Now define an embedding $\iota:Y\to X$ via $\iota(y) = \lb y,1\rb$. Then
\[
\dist(\iota(y_1),\iota(y_2)) = 0\vee \log\Dist^2(y_1,y_2) = \dist(y_1,y_2),
\]
i.e. $\iota$ is an isometric embedding.
\end{proof}

\begin{remark}
Corollary \ref{corollaryultrametricRtree} can also be proven from \cite[Theorem 4.1]{BonkSchramm} by verifying directly that an ultrametric space satisfies Gromov's inequality with an implied constant of zero, and then proving that every geodesic metric space satisfying Gromov's inequality with an implied constant of zero is an $\R$-tree.

However, the proof of Corollary \ref{corollaryultrametricRtree} yields the additional information that the isometric image of $(Y,\dist)$ is contained in a horosphere, i.e.
\begin{equation}
\label{ultrametrichorosphere}
\busemann_\infty(\iota(y_1),\iota(y_2)) = 0 \all y_1,y_2\in Y,
\end{equation}
where $\infty$ is as in Theorem \ref{theoremconeconstruction}.
\end{remark}

\begin{remark}
\label{remarkconverseultrametric}
The converse of the cone construction also holds: if $(X,\dist)$ is an $\R$-tree and $\zero\in X$, $\xi\in\del X$, then $(\del X\butnot\{\xi\},\Dist_{\xi,\zero})$ and $(\{x\in X: \busemann_\xi(\zero,x) = 0\},\dist)$ are both ultrametric spaces.
\end{remark}
\begin{proof}
For all $x,y\in\EE_\xi$, we have $\Dist_\xi(x,y) = \exp\busemann_\xi(\zero,C(x,y,\xi)))$, where $C(x,y,\xi)$ denotes the center of the geodesic triangle $\Delta(x,y,\xi)$ (cf. Definition \ref{definitioncenterRtree}). It can be verified by drawing appropriate diagrams (cf. Figure \ref{figurequadruple}) that for all $x_1,x_2,x_3\in\EE_\xi$, there exists $i$ such that $C(x_i,x_j,\xi) = C(x_i,x_k,\xi)$ and $C(x_j,x_k,\xi) \in \geo{\xi}{C(x_i,x_j,\xi)}$ (where $j,k$ are chosen so that $\{i,j,k\} = \{1,2,3\}$), from which follows the ultrametric inequality for $\Dist_\xi$. Since $\Dist_\xi = e^{(1/2)\dist}$ on $\{x\in X: \busemann_\xi(\zero,x) = 0\}$, the space $(\{x\in X: \busemann_\xi(\zero,x) = 0\},\dist)$ is also ultrametric.
\end{proof}

\begin{theorem}
\label{theoremRtreeorbitalcounting}
Given an unbounded function $f:\CO 0\infty\to\N$, the following are equivalent:
\begin{itemize}
\item[(A)] $f$ is right-continuous and satisfies
\begin{equation}
\label{Rtreeorbitalcounting}
\forall R_1,R_2 \geq 0 \text{ such that } R_1\leq R_2, \;\; f(R_1)\text{ divides } f(R_2).
\end{equation}
\item[(B)] There exist an $\R$-tree $X$ (with a distinguished point $\zero$) and a parabolic group $G\leq\Isom(X)$ such that $\NN_{X,G} = f$.
\item[(C)] There exist an $\R$-tree $X$ (with a distinguished point $\zero$) and a parabolic group $G\leq\Isom(X)$ such that $\NN_{\EE_\bp,G} = f$, where $\bp$ is the global fixed point of $G$.
\end{itemize}
Moreover, in \text{(B)} and \text{(C)} the $\R$-tree $X$ may be chosen to be proper.
\end{theorem}
\begin{proof}[Proof of \text{(A) \implies (B)}]
Let $(\lambda_n)_1^\infty$ and $(N_n)_1^\infty$ be sequences such that
\[
f(\rho) = \prod_{\substack{n\in\N \\ \lambda_n \leq \rho}} N_n.
\]
The hypotheses on $f$ guarantee that $(N_n)_1^\infty$ can be chosen to be integers. Then for each $n\in\N$, let $\Gamma_n$ be a finite group of cardinality $N_n$, and let
\[
\Gamma
= \left\{(\gamma_n)_1^\infty \in \prod_{n\in\N} \Gamma_n : \gamma_n = e \text{ for all but finitely many $n$}\right\}.
\]
For each $(\gamma_n)_1^\infty \in \Gamma$ let
\begin{equation}
\label{infinitecoproduct}
\|(\gamma_n)_1^\infty\| = \max_{\substack{n\in\N \\ \gamma_n\neq e}} \lambda_n,
\end{equation}
with the understanding that $\|e\| = 0$. For each $\alpha,\beta\in \Gamma$ let $\dist(\alpha,\beta) = \|\alpha^{-1} \beta\|$. It is readily verified that $\dist$ is an ultrametric on $\Gamma$. Thus by Corollary \ref{corollaryultrametricRtree}, $(\Gamma,\dist)$ can be isometrically embedded into an $\R$-tree $(X,\dist)$. Since $\Gamma$ is proper, $X$ is proper. Moreover, the natural isometric action of $\Gamma$ on itself extends naturally to an isometric action on $X$. Denote this isometric action by $\phi$, and let $G = \phi(\Gamma)$. Then by \eqref{ultrametrichorosphere}, $G$ is a parabolic group with global fixed point $\infty$. If we let $\zero$ be the image of $e$ under the isometric embedding of $\Gamma$ into $X$, then $G$ satisfies
\[
\NN_{X,G}(\rho) = \#\{\gamma\in \Gamma : \|\gamma\| \leq \rho\} = \prod_{\substack{n\in\N \\ \lambda_n\leq \rho}} \#(\Gamma_n) = f(\rho).
\]
This completes the proof.
\end{proof}
\begin{proof}[Proof of \text{(B) \implies (A)}]
For each $\rho > 0$ let
\[
G_\rho = \{g\in G : \dist(\zero,g(\zero)) \leq \rho\}.
\]
Since $G(\zero)$ is an ultrametric space by Remark \ref{remarkconverseultrametric}, $G_\rho$ is a subgroup of $G$. Thus by Lagrange's theorem, the function $f(\rho) = \NN_{X,G}(\rho) = \#(G_\rho)$ satisfies \eqref{Rtreeorbitalcounting}. Since orbital counting functions are always right-continuous, this completes the proof.
\end{proof}
\begin{proof}[Proof of \text{(A) \iff (C)}]
Since the equation
\[
\NN_{\EE_\bp,G}(R) = \NN_{X,G}(2\log(R))
\]
holds for strongly hyperbolic spaces, including $\R$-trees (Observation \ref{observationeuclideanparabolicasymp}), and since condition (A) is invariant under the transformation $f\mapsto (R\mapsto f(2\log(R)))$, the equivalence (A) \iff (B) directly implies the equivalence (A) \iff (C).
\end{proof}

\begin{remark}
Applying a BIM representation (Theorem \ref{theoremBIM}) shows that if $f:\Rplus\to\N$ is an unbounded function satisfying (A) of Theorem \ref{theoremRtreeorbitalcounting}, then there exists a parabolic group $G\leq\Isom(\H^\infty)$ such that $\NN_{X,G} = f$. This improves a previous result of two of the authors \cite[Proposition A.2]{FSU4}.
\end{remark}

\section{Graphs with contractible cycles}
\label{subsectiongraphtheory}

In Section \ref{subsectionstapling}, we will describe a method of stapling together a collection of $\R$-trees $(X_v)_{v\in V}$ based on some data. This data will include a collection of edge pairings $E\subset V\times V\butnot\{(v,v) : v\in V\}$ that indicates which trees are to be stapled to each other. In this section, we describe the criterion which this collection of edge pairings needs to satisfy in order for the construction to work (Definition \ref{definitioncontractiblecycles}), and we analyze that criterion.

Let $(V,E)$ be an unweighted undirected graph, and let $\dist_E$ denote the path metric of $(V,E)$ (cf. Definition \ref{definitiongraphmetrization}). A sequence $(v_i)_0^n$ in $V$ will be called a \emph{path} if $(v_i,v_{i + 1})\in E\all i < n$. The path $(v_i)_0^n$ is said to \emph{connect} the vertices $v_0$ and $v_n$. The path $(v_i)_0^n$ is called a \emph{geodesic} if $n = \dist_E(v_0,v_n)$, in which case it is denoted $\geo{v_0}{v_n}$. Note that a sequence is a geodesic if and only if $\geo{v_0}{v_1}\ast\cdots\ast\geo{v_{n - 1}}{v_n}$ is a geodesic in the metrization $X(V,E)$ (cf. Definition \ref{definitiongraphmetrization}). Also, recall that a \emph{cycle} in $(V,E)$ is a finite sequence of distinct vertices $v_1,\ldots,v_n\in V$, with $n\geq 3$, such that $(v_1,v_2), (v_2,v_3), \ldots, (v_{n - 1},v_n), (v_n,v_1)\in E$ (cf. \eqref{cycle}).

\begin{definition}
\label{definitioncontractiblecycles}
The graph $(V,E)$ is said to have \emph{contractible cycles} if every cycle forms a complete graph, i.e. if for every cycle $(v_i)_0^n$ we have $(v_i,v_j)\in E \all i,j \text{ such that $v_i\neq v_j$}$.
\end{definition}

\begin{standingassumption}
In the remainder of this section, $(V,E)$ denotes a connected graph with contractible cycles.
\end{standingassumption}

\begin{lemma}
\label{lemmacontractiblecycles1}
For every $v,w\in V$ there exists a unique geodesic $\geo vw = (v_i)_0^n$ connecting $v$ and $w$; moreover, if $(w_j)_0^m$ is any path connecting $v$ and $w$, then the vertices $(v_i)_0^n$ appear in order (but not necessarily consecutively) in the sequence $(w_j)_0^m$.
\end{lemma}
\begin{proof}
\begin{claim}
\label{claimonevertex}
Let $(v_i)_0^n$ be a geodesic, and let $(w_j)_0^m$ be a path connecting $v_0$ and $v_n$. Suppose $n\geq 2$. Then there exist $i = 1,\ldots,n - 1$ and $j = 1,\ldots,m - 1$ such that $v_i = w_j$.
\end{claim}
\begin{subproof}
By contradiction suppose not, and without loss of generality suppose that $(w_j)_0^m$ is minimal with this property. Then the vertices $(w_j)_0^m$ are distinct, since if we had $w_{j_1} = w_{j_2}$ for some $j_1 < j_2$, we could replace $(w_j)_0^m$ by $(w_0,\ldots,w_{j_1 - 1},w_{j_1} = w_{j_2},w_{j_2 + 1},\ldots,w_m)$. Since $n\geq 2$, it follows that the path $(v_0,v_1,\ldots,v_n = w_m,w_{m - 1},\ldots,w_1,w_0 = v_0)$ is a cycle. But then $(v,w)\in E$, contradicting that $(v_i)_0^n$ is a geodesic of length $n\geq 2$.
\end{subproof}
\begin{claim}
\label{claimcontractiblecycles1}
Let $(v_i)_0^n$ be a geodesic, and let $(w_j)_0^m$ be a path connecting $v_0$ and $v_n$. Then the vertices $(v_i)_0^n$ appear in order in the sequence $(w_j)_0^m$.
\end{claim}
\begin{subproof}
We proceed by induction on $n$. The cases $n = 0$, $n = 1$ are trivial. Suppose the claim is true for all geodesics of length less than $n$. By Claim \ref{claimonevertex}, there exist $i_0 = 1,\ldots,n - 1$ and $j_0 = 1,\ldots,m - 1$ such that $v_i = w_j$. By the induction hypothesis, the vertices $(v_i)_0^{i_0}$ appear in order in the sequence $(w_j)_0^{j_0}$, and the vertices $(v_i)_{i_0}^n$ appear in order in the sequence $(w_j)_{j_0}^m$. Combining these facts yields the conclusion.
\end{subproof}
To finish the proof of Lemma \ref{lemmacontractiblecycles1}, it suffices to observe that if $(v_i)_0^n$ and $(w_j)_0^m$ are two geodesics connecting the same vertices $v$ and $w$, then by Claim \ref{claimcontractiblecycles1} the vertices $(v_i)_0^n$ appear in order in the sequence $(w_j)_0^m$, and the vertices $(w_j)_0^m$ appear in order in the sequence $(v_i)_0^n$. It follows that $(v_i)_0^n = (w_j)_0^m$, so geodesics are unique.
\end{proof}

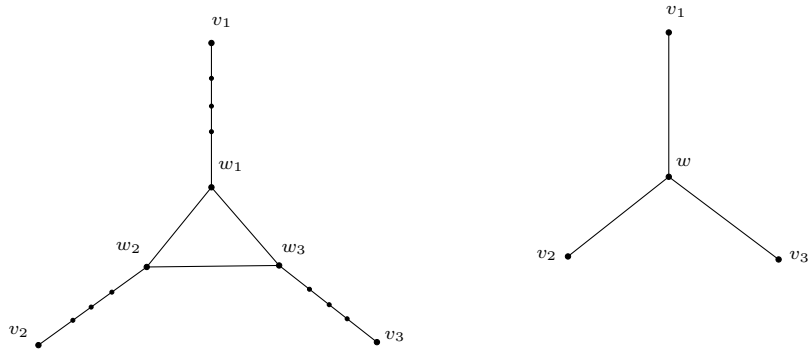
\begin{figure}
\begin{tikzpicture}[line cap=round,line join=round,>=triangle 45,x=1.0cm,y=1.0cm]
\clip(-5.799999999999999,-0.6200000000000001) rectangle (5.639999999999999,5.120000000000001);
\draw (-2.76,4.24)-- (-2.76,2.32);
\draw (-2.76,2.32)-- (-3.62,1.26);
\draw (-3.62,1.26)-- (-1.86,1.28);
\draw (-1.86,1.28)-- (-2.76,2.32);
\draw (-3.62,1.26)-- (-5.06,0.22);
\draw (-1.86,1.28)-- (-0.56,0.26);
\draw (3.32,4.38)-- (3.32,2.46);
\draw (3.32,2.46)-- (1.98,1.4);
\draw (3.32,2.46)-- (4.78,1.36);
\begin{scriptsize}
\draw [fill=black] (-2.76,4.24) circle (1pt);
\draw[color=black] (-2.6199999999999997,4.52) node {$v_1$};
\draw [fill=black] (-2.76,2.32) circle (1pt);
\draw[color=black] (-2.5,2.60) node {$w_1$};
\draw [fill=black] (-3.62,1.26) circle (1pt);
\draw[color=black] (-3.8599999999999994,1.52) node {$w_2$};
\draw [fill=black] (-1.86,1.28) circle (1pt);
\draw[color=black] (-1.6599999999999997,1.50) node {$w_3$};
\draw [fill=black] (-5.06,0.22) circle (1pt);
\draw[color=black] (-5.319999999999999,0.40) node {$v_2$};
\draw [fill=black] (-0.56,0.26) circle (1pt);
\draw[color=black] (-0.31999999999999995,0.34) node {$v_3$};
\draw [fill=black] (3.32,4.38) circle (1pt);
\draw[color=black] (3.42,4.66) node {$v_1$};
\draw [fill=black] (3.32,2.46) circle (1pt);
\draw[color=black] (3.5199999999999996,2.66) node {$w$};
\draw [fill=black] (1.98,1.4) circle (1pt);
\draw[color=black] (1.7,1.42) node {$v_2$};
\draw [fill=black] (4.78,1.36) circle (1pt);
\draw[color=black] (5.059999999999999,1.4) node {$v_3$};
\draw [fill=black] (-2.7599999999999993,3.77) circle (.75pt);
\draw [fill=black] (-2.7599999999999993,3.4) circle (.75pt);
\draw [fill=black] (-2.7599999999999993,3.06) circle (.75pt);
\draw [fill=black] (-4.08369168356998,0.9251115618661258) circle (.75pt);
\draw [fill=black] (-4.358985801217038,0.7262880324543608) circle (.75pt);
\draw [fill=black] (-4.602150101419878,0.5506693711967542) circle (.75pt);
\draw [fill=black] (-1.457011426897158,0.9638089657193085) circle (.75pt);
\draw [fill=black] (-1.1945736888368006,0.7578962789334895) circle (.75pt);
\draw [fill=black] (-0.9559419865221214,0.5706621740404336) circle (.75pt);
\end{scriptsize}
\end{tikzpicture}
\caption[Triangles in graphs with contractible cycles]{The two possibilities for a geodesic triangle in a graph with contractible cycles. Lemma \ref{lemmacontractiblecycles2} states that either the geodesic triangle looks like a triangle in an $\R$-tree (right figure), or there is 3-cycle in the ``center'' of the triangle (left figure).}
\label{figurecontractiblecycles}
\end{figure}

\begin{lemma}[Cf. Figure \ref{figurecontractiblecycles}]
\label{lemmacontractiblecycles2}
Fix $v_1,v_2,v_3\in V$ distinct. Then either
\begin{itemize}
\item[(1)] there exists $w\in V$ such that for all $i\neq j$, $\geo{v_i}{v_j} = \geo{v_i}w\ast \geo w{v_j}$, or
\item[(2)] there exists a cycle $w_1,w_2,w_3\in V$ such that for all $i\neq j$, $\geo{v_i}{v_j} = \geo{v_i}{w_i}\ast \geo{w_i}{w_j}\ast \geo{w_j}{v_j}$.
\end{itemize}
\end{lemma}
\begin{proof}
For each $i = 1,2,3$, let $n_i$ be the number of initial vertices on which the geodesics $\geo{v_i}{v_j}$ and $\geo{v_i}{v_k}$ agree, i.e.
\[
n_i = \max\{n : \geo{v_i}{v_j}_\ell = \geo{v_i}{v_k}_\ell \all \ell = 0,\ldots,n\},
\]
and let $w_i = \geo{v_i}{v_j}_{n_i}$. Here $j,k$ are chosen such that $\{i,j,k\} = \{1,2,3\}$. Then uniqueness of geodesics implies that the geodesics $\geo{w_i}{w_j}$, $i\neq j$ are disjoint except for their common endpoints. If $(w_i)_1^3$ are distinct, then the path $\geo{w_1}{w_2}\ast\geo{w_2}{w_3}\ast\geo{w_3}{w_1}$ is a cycle, and since $(V,E)$ has contractible cycles, this implies $(w_1,w_2),(w_2,w_3),(w_3,w_1)\in E$, completing the proof. Otherwise, we have $w_i = w_j$ for some $i\neq j$; letting $w = w_i = w_j$ completes the proof.
\end{proof}


\begin{figure}
\begin{tikzpicture}[line cap=round,line join=round,>=triangle 45,x=1.0cm,y=1.0cm]
\clip(0.2,2.0) rectangle (10.84,4.68);
\draw (0.48,3.18)-- (2.7,3.2);
\draw (4.04,3.18)-- (6.38,3.2);
\draw (7.26,3.14)-- (8.82,3.12);
\draw (8.82,3.12)-- (10.0,4.0);
\draw (10.0,4.0)-- (10.06,2.64);
\draw (10.06,2.64)-- (8.82,3.12);
\begin{scriptsize}
\draw [fill=black] (0.48,3.18) circle (1.5pt);
\draw[color=black] (0.6200000000000001,3.46) node {$u$};
\draw [fill=black] (1.9195228047394903,3.1929686739165724) circle (1.5pt);
\draw[color=black] (2.06,3.4800000000000004) node {$v_1$};
\draw [fill=black] (2.7,3.2) circle (1.5pt);
\draw[color=black] (2.8400000000000003,3.48) node {$v_2$};

\draw [fill=black] (4.04,3.18) circle (1.5pt);
\draw[color=black] (4.18,3.4600000000000004) node {$u$};
\draw [fill=black] (5.5397195032870705,3.1928181154127095) circle (1.5pt);
\draw[color=black] (5.68,3.4800000000000004) node {$v_2$};
\draw [fill=black] (6.38,3.2) circle (1.5pt);
\draw[color=black] (6.52,3.4800000000000004) node {$v_1$};

\draw [fill=black] (7.26,3.14) circle (1.5pt);
\draw[color=black] (7.3999999999999995,3.42) node {$u$};
\draw [fill=black] (8.82,3.12) circle (1.5pt);
\draw[color=black] (8.95,3.42) node {$w$};
\draw [fill=black] (10.0,4.0) circle (1.5pt);
\draw[color=black] (10.14,4.28) node {$v_1$};
\draw [fill=black] (10.06,2.64) circle (1.5pt);
\draw[color=black] (10.33,2.6) node {$v_2$};
\end{scriptsize}
\end{tikzpicture}
\caption[Triangles in graphs with contractible cycles: general case]{When the vertices $v_1$ and $v_2$ are adjacent, Corollary \ref{corollarycontractiblecycles} describes three possible pictures for the geodesic triangle $\Delta(u,v_1,v_2)$. In the rightmost figure, $w$ is the vertex adjacent to both $v_1$ and $v_2$ from which the paths $\geo u{v_1}$ and $\geo u{v_2}$ diverge.}
\label{figurecontractiblecycles2}
\end{figure}
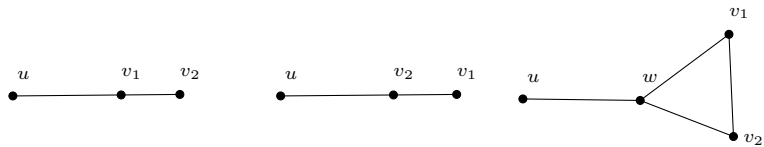

\begin{corollary}[Cf. Figure \ref{figurecontractiblecycles2}]
\label{corollarycontractiblecycles}
Fix $v_1,v_2,u\in V$ distinct such that $(v_1,v_2)\in E$. Then either $v_1\in\geo u{v_2}$, $v_2\in\geo u{v_1}$, or there exists $w\in V$ such that for each $i = 1,2$, $(w,v_i)\in E$ and $w\in \geo u{v_i}$.
\end{corollary}
\begin{proof}
Write $v_3 = u$, so that we can use the same notation as Lemma \ref{lemmacontractiblecycles2}. If we are in case (1), then the equation $\geo{v_1}{v_2} = \geo{v_1}w\ast \geo w{v_2}$ implies that $w\in\{v_1,v_2\}$, and so either $v_1 = w\in \geo u{v_2}$ or $v_2 = w\in \geo u{v_1}$. If we are in case 2, then the equation $\geo{v_1}{v_2} = \geo{v_1}{w_1}\ast\geo{w_1}{w_2}\ast\geo{w_2}{v_2}$ implies that $w_1 = v_1$ and $w_2 = v_2$. Letting $w = w_3$ completes the proof.
\end{proof}

\section{The nearest-neighbor projection onto a convex set}
\label{subsectionRtreegeometric}

Let $X$ be an $\R$-tree, and let $A\subset X$ be a nonempty closed convex set. Since $X$ is a CAT(-1) space, for each $z\in X$ there is a unique point $\pi(z)\in A$ such that $\dist(z,\pi(z)) = \dist(z,A)$, and the map $z\to\pi(z)$ is semicontracting (see e.g. \cite{BridsonHaefliger}). Since $X$ is an $\R$-tree, we can say more about this nearest-neighbor projection map $\pi$, as well as providing a simpler proof of its existence. In the following theorems, $X$ denotes an $\R$-tree.

\begin{lemma}
\label{lemmaconvexprojection1}
Let $A\subset X$ be a nonempty closed convex set. Then for each $z\in X$ there exists a unique point $\pi(z)\in A$ such that for all $x\in A$, $\pi(z)\in \geo zx$. Moreover, for all $z_1,z_2\in X$, we have
\begin{equation}
\label{convexprojection1}
\dist(\pi(z_1),\pi(z_2)) = 0\vee(\dist(z_1,z_2) - \dist(z_1,A) - \dist(z_2,A)).
\end{equation}
\end{lemma}
\begin{proof}
Since $A$ is nonempty and closed, there exists a point $\pi(z)\in A$ such that $\geo z{\pi(z)}\cap A = \{\pi(z)\}$. Fix $z\in A$. Since $C(x,z,\pi(z))\in \geo z{\pi(z)}\cap \geo x{\pi(x)}\subset \geo z{\pi(z)}\cap A$, we get $C(x,z,\pi(z)) = \pi(z)$, i.e. $\lb x|z\rb_{\pi(z)} = 0$, i.e. $\pi(z)\in \geo zx$. This completes the proof of existence; uniqueness is trivial.

To demonstrate the equation \eqref{convexprojection1}, we consider two cases:
\begin{itemize}
\item[Case 1:] If $\geo{z_1}{z_2}\cap A\neq\emptyset$, then $\pi(z_1)$ and $\pi(z_2)$ both lie on the geodesic $\geo{z_1}{z_2}$, so $\dist(\pi(z_1),\pi(z_2)) = \dist(z_1,z_2) - \dist(z_1,A) - \dist(z_2,A) \geq 0$.
\item[Case 2:]  Suppose that $\geo{z_1}{z_2}\cap A = \emptyset$; we claim that $\pi(z_1) = \pi(z_2)$. Indeed, by the definition of $\pi(z_2)$ we have $\pi(z_2)\in \geo{z_2}{\pi(z_1)}$, and by assumption we have $\pi(z_2)\notin \geo{z_1}{z_2}$, so we must have $\pi(z_2)\in \geo{z_1}{\pi(z_1)}$. But from the definition of $\pi(z_1)$, this can only happen if $\pi(z_1) = \pi(z_2)$. The proof is completed by noting that the triangle inequality gives $\dist(z_1,z_2) - \dist(z_1,A) - \dist(z_2,A) = \dist(z_1,z_2) - \dist(z_1,\pi(z_1)) - \dist(z_2,\pi(z_1)) \leq 0$.
\end{itemize}
\end{proof}

\begin{lemma}
\label{lemmaconvexprojection2}
Let $A_1, A_2\subset X$ be closed convex sets such that $A_1\cap A_2\neq\emptyset$. For each $i$ let $\pi_i:X\to A_i$ denote the nearest-neighbor projection map. Then for all $z\in X$, either $\pi_1(z)\in A_2$ or $\pi_2(z)\in A_1$. In particular, $\pi_1(A_2) \subset A_1\cap A_2$.
\end{lemma}
\begin{proof}
Let $x_1 = \pi_1(z)$ and $x_2 = \pi_2(z)$, and fix $y\in A_1\cap A_2$. By Lemma \ref{lemmaconvexprojection1}, $x_1,x_2\in \geo zy$. Without loss of generality assume $\dist(z,x_1)\leq \dist(z,x_2)$, so that $x_2\in \geo{x_1}y$. Since $A_1$ is convex, $x_2\in A_1$.
\end{proof}

\begin{lemma}
\label{lemmaconvexunion}
Let $A_1, A_2\subset X$ be closed convex sets such that $A_1\cap A_2\neq\emptyset$. Then $A_1\cup A_2$ is convex.
\end{lemma}
\begin{proof}
It suffices to show that if $x_1\in A_1$ and $x_2\in A_2$, then $\geo{x_1}{x_2}\subset A_1\cup A_2$. Since $x_2\in A_2$, Lemma \ref{lemmaconvexprojection1} shows that $\geo{x_1}{x_2}$ intersects the point $\pi_2(x_1)$. By Lemma \ref{lemmaconvexprojection2}, $\pi_2(x_1)\in A_1\cap A_2$. But then the two subsegments $\geo{x_1}{\pi_2(x_1)}$ and $\geo{\pi_2(x_1)}{x_2}$ are contained in $A_1\cup A_2$, so the entire geodesic $\geo{x_1}{x_2}$ is contained in $A_1\cup A_2$.
\end{proof}


\section{Constructing $\R$-trees by the stapling method}
\label{subsectionstapling}

We now describe the ``stapling method'' for constructing $\R$-trees. The following definition is phrased for arbitrary metric spaces.

\begin{definition}
\label{definitionstapledunion}
Let $(V,E)$ be an unweighted undirected graph, let $(X_v)_{v\in V}$ be a collection of metric spaces, and for each $(v,w)\in E$ fix a set $\set vw\subset X_v$ and an isometry $\bij vw:\set vw\to \set wv$ such that $\bij wv = \bij vw^{-1}$. Let $\sim$ be the equivalence relation on $\coprod_{v\in V} X_v$ defined by the relations
\[
x\sim \bij vw(x) \all (v,w)\in E \all x\in \set vw.
\]
Then the \emph{stapled union} of of the collection $(X_v)_{v\in V}$ with respect to the sets $(\set vw)_{(v,w)\in E}$ and the bijections $(\bij vw)_{(v,w)\in E}$ is the set 
\[
X =  \coprod_{v\in V}^{\text{st}} X_v := \coprod_{v\in V} X_v / \sim,
\] 
equipped with the \emph{path metric}
\begin{equation}
\label{pathmetricv3}
\dist\big(\lb v,x\rb,\lb w,y\rb\big) = \inf\left\{\sum_{i = 0}^n \dist_{v_i}(x_i,y_i) \left\vert
\begin{split}
v_0,\ldots,v_n\in V\\
(v_i,v_{i + 1}) \in E \all i < n\\
v_0 = v, \; v_n = w\\
y_i\in \set{v_i}{v_{i + 1}} \all i < n\\
x_{i + 1} = \bij{v_i}{v_{i + 1}}(y_i) \all i < n\\
x_0 = x, y_n = y
\end{split}
\right.\right\}.
\end{equation}
Note that $\dist$ is finite as long as the graph $(V,E)$ is connected. We leave it to the reader to verify that in this case, $\dist$ is a metric on $X$.
\end{definition}

\begin{example}
If for each $(v,w)\in E$ we fix a point $\point vw\in X_v$, then we can let $\set vw = \{\point vw\}$ and let $\bij vw$ be the unique bijection between $\{\point vw\}$ and $\{\point wv\}$.
\end{example}

Intuitively, the stapled union $\coprod_{v\in V}^{\text{st}} X_v$ is the metric space that results from starting with the spaces $(X_v)_{v\in V}$ and for each $(v,w)\in E$, stapling the set $\set vw\subset X_v$ with the set $\set wv\subset X_w$ along the bijection $\bij vw$.

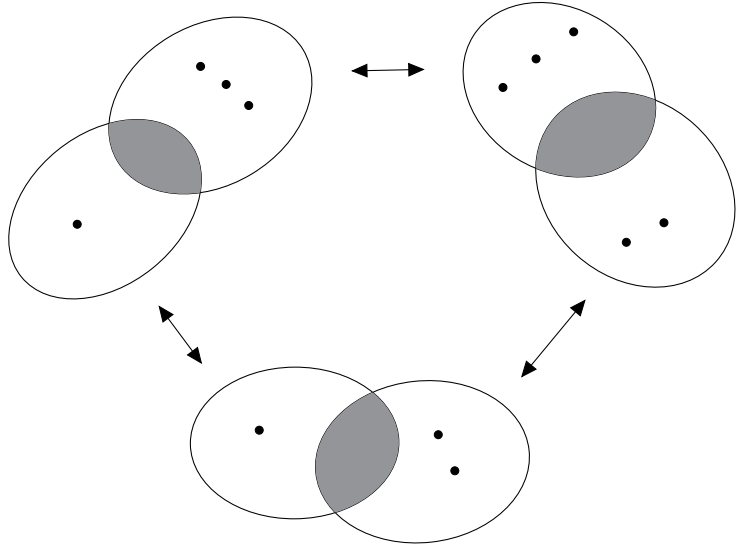
\begin{figure}
\begin{tikzpicture}[line cap=round,line join=round,>=triangle 45,x=1.0cm,y=1.0cm]
\clip(-6.280000000000002,-2.5000000000000018) rectangle (5.080000000000002,5.38);
\draw\ellipseA;
\draw\ellipseB;
\draw\ellipseC;
\draw\ellipseD;
\draw\ellipseE;
\draw\ellipseF;
\begin{scope}
\clip\ellipseA;
\fill[Gray]\ellipseB;
\end{scope}
\begin{scope}
\clip\ellipseC;
\fill[Gray]\ellipseD;
\end{scope}
\begin{scope}
\clip\ellipseE;
\fill[Gray]\ellipseF;
\end{scope}
\draw[<->] (-0.88,4.1)-- (0.1,4.12);
\draw[<->] (-2.86,0.2)-- (-3.44,0.98);
\draw[<->] (2.24,1.06)-- (1.38,0.02);
\begin{scriptsize}
\draw [fill=black] (-4.52,2.06) circle (1.5pt);
\draw [fill=black] (-2.88,4.16) circle (1.5pt);
\draw [fill=black] (-2.54,3.92) circle (1.5pt);
\draw [fill=black] (-2.24,3.64) circle (1.5pt);
\draw [fill=black] (1.14,3.88) circle (1.5pt);
\draw [fill=black] (1.58,4.26) circle (1.5pt);
\draw [fill=black] (2.08,4.62) circle (1.5pt);
\draw [fill=black] (3.28,2.08) circle (1.5pt);
\draw [fill=black] (2.78,1.82) circle (1.5pt);
\draw [fill=black] (0.28,-0.74) circle (1.5pt);
\draw [fill=black] (0.5,-1.22) circle (1.5pt);
\draw [fill=black] (-2.1,-0.68) circle (1.5pt);
\end{scriptsize}
\end{tikzpicture}
\caption[The consistency condition for stapling metric spaces]{In this diagram, the arrows represent the bijections $\bij{v_i}{v_j}$, while the ovals represent the sets $\set{v_i}{v_j}$. The consistency condition (Definition \ref{definitionconsistency}) states that (I) each of the shaded regions is nonempty, (IIa) shaded regions go to shaded regions, and (IIb) if you start in a shaded region and traverse the diagram, then you will get back to where you started.}
\label{figureconsistencycondition}
\end{figure}

\begin{definition}[Cf. Figure \ref{figureconsistencycondition}]
\label{definitionconsistency}
We say that the \emph{consistency condition} is satisfied if for every $3$-cycle $u,v,w\in V$, we have
\begin{itemize}
\item[(I)] $\set uv\cap \set uw \neq \emptyset$, and
\item[(II)] for all $z\in \set uv\cap \set uw$, we have
\begin{itemize}
\item[(a)] $\bij uw(z)\in \set wv$ and
\item[(b)] $\bij wv\bij uw(z) = \bij uv(z)$.
\end{itemize}
\end{itemize}
\end{definition}

Obviously, the consistency condition is satisfied whenever $(V,E)$ has no cycles. Theorem \ref{theoremschottkyproductRtree} and Examples \ref{exampleconeconstruction}-\ref{examplegeometricproducts} below show how it can be satisfied in many reasonable circumstances. Now we prove the main theorem of this chapter: for a connected graph with contractible cycles, the consistency condition implies that the stapled union of $\R$-trees is an $\R$-tree, if the staples are taken along convex sets. More precisely:

\begin{theorem}
\label{theoremstapledunion}
Let $(V,E)$ be a connected graph with contractible cycles, let $(X_v)_{v\in V}$ be a collection of $\R$-trees, and for each $(v,w)\in E$ let $\set vw\subset X_v$ be a nonempty closed convex set and let $\bij vw:\set vw\to \set wv$ be an isometry such that $\bij wv = \bij vw^{-1}$. Assume that the consistency condition is satisfied. Then
\begin{itemize}
\item[(i)] The stapled union $X = \coprod_{v\in V}^{\text{st}} X_v$ is an $\R$-tree.
\item[(ii)] The infimum in \eqref{pathmetricv3} is achieved when
\begin{itemize}
\item[(a)] $(v_i)_0^n = \geo vw$, and
\item[(b)] for each $i < n$, $y_i$ is the image of $x_i$ under the nearest-neighbor projection to $\set{v_i}{v_{i + 1}}$.
\end{itemize}
\end{itemize}
\end{theorem}
\begin{proof}
We prove part (ii) first. For each $(v,w)\in E$, let $\pi_{v,w}:X_v\to \set vw$ be the nearest-neighbor projection; then $\pi_{v,w}$ is $1$-Lipschitz. Now fix $v\in V$ arbitrary. We define a map $\pi_v:X\to X_v$ as follows. Fix $\wbar x = \lb w,x\rb\in X$, so that $x\in X_w$. Let $(v_i)_0^n = \geo vw$, and let
\[
\pi_v(\wbar x) = \pi_v(w,x) = \bij{v_1}{v_0}\pi_{v_1,v_0}\cdots\bij{v_n}{v_{n - 1}}\pi_{v_n,v_{n - 1}}(x).
\]
\begin{claim}
\label{claimpiv}
The map $\pi_v$ is well-defined.
\end{claim}
\begin{subproof}
Fix $(u,w)\in E$ and $x\in \set uw$ and let $y = \bij uw(x)$; we need to show that $\pi_v(u,x) = \pi_v(w,y)$. If $w\in \geo vu$ or $u\in \geo vw$ then the equality is trivial, so by Corollary \ref{corollarycontractiblecycles} we are reduced to proving the case where there exists $v'\in V$ such that $(v',w),(v',u)\in E$ and $v'\in \geo vw,\geo vu$. We have $\pi_v(u,x) = \pi_v(v',\bij u{v'}\pi_{u,v'}(x))$ and $\pi_v(w,y) = \pi_v(v',\bij w{v'}\pi_{w,v'}(y))$, so to complete the proof it suffices to show that
\begin{equation}
\label{NTSpiv}
\bij u{v'}\pi_{u,v'}(x) = \bij w{v'}\pi_{w,v'}(y).
\end{equation}
Since $u,v',w$ form a $3$-cycle, part (I) of the consistency condition gives $\set u{v'}\cap \set uw\neq\emptyset$. By Lemma \ref{lemmaconvexprojection2}, we have $x' := \pi_{u,v'}(x)\in \set u{v'}\cap \set uw$. Applying part (IIa) of the consistency condition gives $y'' := \bij uw(x') \in \set w{v'}$ and thus $\dist(x,\set u{v'}) = \dist(x,x') = \dist(y,y'') \leq \dist(y,\set w{v'})$. A symmetric argument gives $\dist(y,\set w{v'}) \leq \dist(x,\set u{v'})$, so we have equality and thus $y'' = y' := \pi_{w,v'}(y)$. Applying part (IIb) of the consistency condition gives $\bij u{v'}(x') = \bij w{v'}(y')$, i.e. \eqref{NTSpiv} holds.
\end{subproof}

Since for each $w\in V$ the map $X_w\ni x\mapsto \pi_v(w,x)\in X_v$ is $1$-Lipschitz, the map $\pi_v:X\to X_v$ is also $1$-Lipschitz.

Fix $\wbar x = \lb v,x\rb,\wbar y = \lb w,y\rb\in X$. Let $(v_i)_0^n$, $(x_i)_0^n$, and $(y_i)_0^n$ be as in (ii), i.e.
\begin{align*}
(v_i)_0^n = \geo vw, \text{ where } &\; x_0 = x,\\ 
&\; y_i = \pi_{v_i,v_{i + 1}}(x_i)  \all i < n,\\ 
&\; x_{i + 1} = \bij{v_i}{v_{i + 1}}(y_i) \all i < n, \text{ and }\\ 
&\; y_n = y.
\end{align*}
We define a function $f:X\to \R^{n + 1}$ as follows: for each $\wbar z\in X$, we let
\[
f(\wbar z) = \big(\dist_{v_i}(x_i,\pi_{v_i}(\wbar z))\big)_{i = 0}^n.
\]
Then $f$ is $1$-Lipschitz, when $\R^{n + 1}$ is interpreted as having the max norm.

\begin{claim}
\label{claimfS}
Fix $\wbar z\in X$ and $i = 0,\ldots,n - 1$. If $f_{i + 1}(\wbar z) > 0$, then $f_i(\wbar z) \geq r_i := \dist_{v_i}(x_i,y_i)$.
\end{claim}
\begin{subproof}
By contradiction, suppose that $f_{i + 1}(\wbar z) > 0$ but $f_i(\wbar z) < \dist_{v_i}(x_i,y_i)$. Then $z_{i + 1} := \pi_{v_{i + 1}}(\wbar z) \neq x_{i + 1}$, but $z_i := \pi_{v_i}(\wbar z) \in B(x_i,r_i)\butnot\{y_i\}$. In particular, $\pi_{v_{i + 1},v_i}(z_i) = y_i$, so
\begin{equation}
\label{fS1}
z_{i + 1} \neq \bij{v_i}{v_{i + 1}}\pi_{v_i,v_{i + 1}}(z_i).
\end{equation}
On the other hand, since $z_i \notin \set{v_i}{v_{i + 1}}$, we have 
\begin{equation}
\label{fS2}
z_i \neq \bij{v_{i + 1}}{v_i}\pi_{v_{i + 1},v_i}(z_{i + 1}).
\end{equation}
Write $\wbar z = \lb w,z\rb$. Then the definition of the maps $(\pi_v)_{v\in V}$ together with \eqref{fS1}, \eqref{fS2} implies that $v_i\notin\geo w{v_{i + 1}}$ and $v_{i + 1}\notin\geo w{v_i}$. Thus by Corollary \ref{corollarycontractiblecycles}, there exists $w'\in V$ such that $(w',v_i),(w',v_{i + 1})\in E$ and $w'\in\geo w{v_i},\geo w{v_{i + 1}}$. Let $z' = \pi_{w'}(\wbar z)$, so that $\bij{w'}{v_i}\pi_{w',v_i}(z') = z_i$ and $\bij{w'}{v_{i + 1}}\pi_{w',v_{i + 1}}(z') = z_{i + 1}$. Let $F = \bij{v_i}{w'}(B(x_i,r_i)\cap \set{v_i}{w'})$, and let $\pi_F:X_{w'}\to F$ be the nearest-neighbor projection map. By Lemma \ref{lemmaconvexprojection2}, either $\pi_F(z')\in \set{w'}{v_{i + 1}}$ or $\pi_{w',v_{i + 1}}(z')\in F$.
\begin{itemize}
\item[Case 1:] $\pi_F(z')\in \set{w'}{v_{i + 1}}$. Since $F\subset \set{w'}{v_i}$ and $\pi_{w',v_i}(z')\in F$, we have $\pi_{w',v_i}(z') = \pi_F(z')\in \set{w'}{v_{i + 1}}$ and then part (IIa) of the consistency condition gives $z_i = \bij{w'}{v_i}\pi_{w',v_i}(z') \in \set{v_i}{v_{i + 1}}$, a contradiction.
\item[Case 2:] $\pi_{w',v_{i + 1}}(z')\in F$. Since $F\subset \set{w'}{v_i}$, part (IIa) of the consistency condition gives $z_{i + 1} = \bij{w'}{v_{i + 1}}\pi_{w',v_{i + 1}}(z') \in \set{v_{i + 1}}{v_i}$ and $\bij{v_{i + 1}}{v_i}(z_{i + 1})\in \bij{w'}{v_i}(F) \subset B(x_i,r_i)$. But then $\bij{v_{i + 1}}{v_i}(z_{i + 1}) = y_i$ and thus $z_{i + 1} = x_{i + 1}$, a contradiction.
\end{itemize}
\end{subproof}

Thus $f(X)$ is contained in the set
\[
S = \{(t_i)_0^n : \forall i = 0,\ldots,n - 1 \; t_{i + 1} > 0 \; \Rightarrow t_i \geq r_i\} \subset \R^{n + 1}.
\]
Now the function $h:S\to\R$ defined by
\[
h\big((t_i)_0^n\big) = \max_{\substack{i \in\{0,\ldots,n\} \\ t_i > 0 \text{ if }i > 0}} [r_0 + \ldots + r_{i - 1} + t_i]
\]
is Lipschitz $1$-continuous with respect to the path metric of the max norm. Thus since $X$ is a path-metric space, $h\circ f:X\to\R$ is Lipschitz $1$-continuous. Thus
\[
\dist(\wbar x,\wbar y) \geq h\circ f(\wbar y) - h\circ f(\wbar x) \geq r_0 + \ldots + r_n = \sum_{i = 0}^n \dist_{v_i}(x_i,y_i),
\]
completing the proof of (ii).

For each $\wbar x = \lb v,x\rb,\wbar y = \lb w,y\rb\in X$, let
\[
\geo{\wbar x}{\wbar y}= \geo{x_0}{y_0}_{v_0}\ast \cdots \ast \geo{x_n}{y_n}_{v_n},
\]
where $\ast$ denotes the concatenation of geodesics, and $(v_i)_0^n$, $(x_i)_0^n$, and $(y_i)_0^n$ are as in (ii). Here $\geo xy_v$ denotes the image of the geodesic $\geo xy$ under the map $X_v\ni z\to \lb v,z\rb\in X$. Then by (ii), $\geo{\wbar x}{\wbar y}$ is a geodesic connecting $\wbar x$ and $\wbar y$. Thus we have a family of geodesics $(\geo{\wbar x}{\wbar y})_{\wbar x,\wbar y\in X}$.

We now prove that $X$ is an $\R$-tree, using the criteria of Lemma \ref{lemmaRtreeequivalent}. Condition (BII) is readily verified. So to complete the proof, we must demonstrate (BIII). Fix $\wbar x_1,\wbar x_2,\wbar x_3\in X$ distinct, and we show that two of the geodesics $\geo{\wbar x_i}{\wbar x_j}$ have a nontrivial intersection. Write $\wbar x_i = \lb v_i,x_i\rb$. If there is more than one possible choice, choose $(v_i)_1^3$ so as to minimize $\sum_{i\neq j} \dist_E(v_i,v_j)$.

Let $w_1,w_2,w_3\in V$ be as in Lemma \ref{lemmacontractiblecycles2}, with the convention that $w_1 = w_2 = w_3 = w$ if we are in Case 1 of Lemma \ref{lemmacontractiblecycles2}.
\begin{itemize}
\item[Case A:] For some $i$, $v_i\neq w_i$. Choose $j,k$ such that $i,j,k$ are distinct. Then there exists a vertex $w\in V$ adjacent to $v_i$ such that $w\in\geo{v_i}{v_j}\cap \geo{v_i}{v_k}$. The choice of $(v_i)_1^3$ guarantees that $x_i\notin \set{v_i}w$, so that $\geo{x_i}{\pi_{v_i,w}(x_i)}_{v_i}$ forms a common initial segment of the geodesics $\geo{\wbar x_i}{\wbar x_j}$ and $\geo{\wbar x_i}{\wbar x_k}$.
\item[Case B:] For all $i$, $v_i = w_i$. Then either $v_1 = v_2 = v_3$, or $v_1,v_2,v_3$ form a cycle.
\begin{itemize}
\item[Case B1:] Suppose that $v_1 = v_2 = v_3 = v$. Then since $X_v$ is an $\R$-tree, there exist distinct $i,j,k\in\{1,2,3\}$ such that the geodesics $\geo{x_i}{x_j}_v$ and $\geo{x_i}{x_k}_v$ have a common initial segment.
\item[Case B2:] Suppose that $v_1,v_2,v_3$ form a cycle. Then by part (I) of the consistency condition $\set{v_1}{v_2}\cap\set{v_1}{v_3}\neq\emptyset$, so by Lemma \ref{lemmaconvexunion} the set $F = \set{v_1}{v_2}\cup\set{v_1}{v_3}$ is convex. But the choice of $(v_i)_1^3$ guarantees that $x_1\notin F$, so that $\geo{x_1}{\pi_F(x_1)}_{v_1}$ forms a common initial segment of the geodesics $\geo{\wbar x_1}{\wbar x_2}$ and $\geo{\wbar x_1}{\wbar x_3}$.
\end{itemize}
\end{itemize}
\end{proof}

\section{Examples of $\R$-trees constructed using the stapling method}
\label{subsectionstaplingexamples}
We give three examples of ways to construct $\R$-trees using the stapling method so that the resulting $\R$-tree admits a natural isometric action.

\begin{example}[Cone construction again]
\label{exampleconeconstruction}
Let $(Z,\Dist)$ be a complete ultrametric space, let $V = Z$ and $E = V\times V\butnot\{(v,v) : v\in V\}$, and for each $v\in V$ let $X_v = \R$. For each $v,w\in V$ let $\set vw = \CO{\log\Dist(v,w)}\infty$, and let $\bij vw$ be the identity map. Since $(V,E)$ is a complete graph, it is connected and has contractible cycles. Part (IIa) of the consistency condition is equivalent to the ultrametric inequality for $\Dist$, while parts (I) and (IIb) are obvious. Thus we can consider the stapled union $X = \coprod_{v\in V}^{\text{st}} X_v$. One can verify that the stapled union is isometric to the $\R$-tree $X$ considered in the proof of Theorem \ref{theoremconeconstruction}. Indeed, the map $\lb z,t\rb \mapsto \lb z,e^t\rb$ provides the desired isometry. Note that the map $\iota$ constructed in Theorem \ref{theoremconeconstruction} can be described in terms of the stapled union as follows: For each $z\in Z$, $\iota(z)$ is the image of $-\infty$ under the isometric embedding of $X_z \equiv\R$ into $X$. (The image of $+\infty$ is $\infty$).
\end{example}

Our next example is a type of Schottky product which we call a ``pure Schottky product''. To describe it, it will be convenient to introduce the following terminology:

\begin{definition}
\label{definitiontreegeometric}
If $\Gamma$ is a group, a function $\|\cdot\|:\Gamma\to\Rplus$ is called \emph{tree-geometric} if there exist an $\R$-tree $X$, a distinguished point $\zero\in X$, and an isometric action $\phi:\Gamma\to\Isom(X)$ such that
\[
\dogo{\phi(\gamma)} = \|\gamma\| \all \gamma\in\Gamma.
\]
\end{definition}

\begin{example}
Theorem \ref{theoremRtreeorbitalcounting} gives a sufficient but not necessary condition for a function to be tree-geometric.
\end{example}

\begin{remark}
If the group $\Gamma$ is countable, then whenever $\Gamma$ is a tree-geometric function, the $\R$-tree $X$ can be chosen to be separable.
\end{remark}
\begin{proof}
Without loss of generality, we may replace $X$ by the convex hull of $\Gamma(\zero)$.
\end{proof}

\begin{figure}
\begin{center}
\begin{tabular}{@{}ll@{}}

\begin{tikzpicture}[line cap=round,line join=round,>=triangle 45,x=0.5cm,y=0.5cm]
\clip(-6.705221014292816,-5.620235521243468) rectangle (6.556270148687806,5.50828853160739);
\draw (0.0,3.0)-- (0.0,-3.0);
\draw (-3.0,0.0)-- (3.0,0.0);
\draw (-1.5,3.0)-- (0.0,3.0);
\draw (0.0,3.0)-- (1.5,3.0);
\draw (0.0,3.0)-- (0.0,4.5);
\draw (-3.0,0.0)-- (-3.0,1.5);
\draw (-3.0,0.0)-- (-4.5,0.0);
\draw (-3.0,0.0)-- (-3.0,-1.5);
\draw (0.0,-3.0)-- (-1.5,-3.0);
\draw (0.0,-3.0)-- (1.5,-3.0);
\draw (0.0,-3.0)-- (0.0,-4.5);
\draw (3.0,0.0)-- (3.0,1.5);
\draw (3.0,0.0)-- (3.0,-1.5);
\draw (0.0,4.5)-- (0.0,5.5);
\draw (0.0,4.5)-- (-1.0,4.5);
\draw (0.0,4.5)-- (1.0,4.5);
\draw (-1.5,3.0)-- (-1.5,4.0);
\draw (-1.5,3.0)-- (-2.5,3.0);
\draw (-1.5,3.0)-- (-1.5,2.0);
\draw (1.5,3.0)-- (1.5,4.0);
\draw (1.5,3.0)-- (2.5,3.0);
\draw (1.5,3.0)-- (1.5,2.0);
\draw (3.0,1.5)-- (3.0,2.5);
\draw (3.0,1.5)-- (2.0,1.5);
\draw (3.0,1.5)-- (4.0,1.5);
\draw (3.0,-1.5)-- (2.0,-1.5);
\draw (3.0,-1.5)-- (4.0,-1.5);
\draw (3.0,-1.5)-- (3.0,-2.5);
\draw (1.5,-3.0)-- (1.5,-2.0);
\draw (1.5,-3.0)-- (2.5,-3.0);
\draw (1.5,-3.0)-- (1.5,-4.0);
\draw (0.0,-4.5)-- (1.0,-4.5);
\draw (0.0,-4.5)-- (0.0,-5.5);
\draw (0.0,-4.5)-- (-1.0,-4.5);
\draw (-1.5,-3.0)-- (-2.5,-3.0);
\draw (-1.5,-3.0)-- (-1.5,-2.0);
\draw (-1.5,-3.0)-- (-1.5,-4.0);
\draw (-3.0,-1.5)-- (-2.0,-1.5);
\draw (-3.0,-1.5)-- (-3.0,-2.5);
\draw (-3.0,-1.5)-- (-4.0,-1.5);
\draw (-4.5,0.0)-- (-4.5,1.0);
\draw (-4.5,0.0)-- (-4.5,-1.0);
\draw (-4.5,0.0)-- (-5.5,0.0);
\draw (-3.0,1.5)-- (-4.0,1.5);
\draw (-3.0,1.5)-- (-3.0,2.5);
\draw (-3.0,1.5)-- (-2.0,1.5);
\draw (4.5,0.0)-- (4.5,1.0);
\draw (4.5,0.0)-- (4.5,-1.0);
\draw (4.5,0.0)-- (5.5,0.0);
\draw (3.0,0.0)-- (4.5,0.0);
\begin{scriptsize}
\draw [fill=black] (0.0,3.0) circle (1.0pt);
\draw [fill=black] (0.0,-3.0) circle (1.0pt);
\draw [fill=black] (-3.0,0.0) circle (1.0pt);
\draw [fill=black] (3.0,0.0) circle (1.0pt);
\draw [fill=black] (-1.5,3.0) circle (1.0pt);
\draw [fill=black] (1.5,3.0) circle (1.0pt);
\draw [fill=black] (0.0,4.5) circle (1.0pt);
\draw [fill=black] (-3.0,1.5) circle (1.0pt);
\draw [fill=black] (-4.5,0.0) circle (1.0pt);
\draw [fill=black] (-3.0,-1.5) circle (1.0pt);
\draw [fill=black] (-1.5,-3.0) circle (1.0pt);
\draw [fill=black] (1.5,-3.0) circle (1.0pt);
\draw [fill=black] (0.0,-4.5) circle (1.0pt);
\draw [fill=black] (3.0,1.5) circle (1.0pt);
\draw [fill=black] (3.0,-1.5) circle (1.0pt);
\draw [fill=black] (4.5,0.0) circle (1.0pt);
\draw [fill=black] (0.0,5.4) circle (1.0pt);
\draw [fill=black] (-1.0,4.5) circle (1.0pt);
\draw [fill=black] (1.0,4.5) circle (1.0pt);
\draw [fill=black] (-1.5,4.0) circle (1.0pt);
\draw [fill=black] (-2.5,3.0) circle (1.0pt);
\draw [fill=black] (-1.5,2.0) circle (1.0pt);
\draw [fill=black] (1.5,4.0) circle (1.0pt);
\draw [fill=black] (2.5,3.0) circle (1.0pt);
\draw [fill=black] (1.5,2.0) circle (1.0pt);
\draw [fill=black] (3.0,2.5) circle (1.0pt);
\draw [fill=black] (2.0,1.5) circle (1.0pt);
\draw [fill=black] (4.0,1.5) circle (1.0pt);
\draw [fill=black] (2.0,-1.5) circle (1.0pt);
\draw [fill=black] (4.0,-1.5) circle (1.0pt);
\draw [fill=black] (3.0,-2.5) circle (1.0pt);
\draw [fill=black] (1.5,-2.0) circle (1.0pt);
\draw [fill=black] (2.5,-3.0) circle (1.0pt);
\draw [fill=black] (1.5,-4.0) circle (1.0pt);
\draw [fill=black] (1.0,-4.5) circle (1.0pt);
\draw [fill=black] (0.0,-5.5) circle (1.0pt);
\draw [fill=black] (-1.0,-4.5) circle (1.0pt);
\draw [fill=black] (-2.5,-3.0) circle (1.0pt);
\draw [fill=black] (-1.5,-2.0) circle (1.0pt);
\draw [fill=black] (-1.5,-4.0) circle (1.0pt);
\draw [fill=black] (-2.0,-1.5) circle (1.0pt);
\draw [fill=black] (-3.0,-2.5) circle (1.0pt);
\draw [fill=black] (-4.0,-1.5) circle (1.0pt);
\draw [fill=black] (-4.5,1.0) circle (1.0pt);
\draw [fill=black] (-4.5,-1.0) circle (1.0pt);
\draw [fill=black] (-5.5,0.0) circle (1.0pt);
\draw [fill=black] (-4.0,1.5) circle (1.0pt);
\draw [fill=black] (-3.0,2.5) circle (1.0pt);
\draw [fill=black] (-2.0,1.5) circle (1.0pt);
\draw [fill=black] (4.5,1.0) circle (1.0pt);
\draw [fill=black] (4.5,-1.0) circle (1.0pt);
\draw [fill=black] (5.5,0.0) circle (1.0pt);
\draw [fill=black] (0.0,0.0) circle (1.0pt);
\end{scriptsize}
\end{tikzpicture}
\end{tabular}

\caption[The Cayley graph of $\F_2(\Z)$ as a pure Schottky product]{The Cayley graph of $\F_2(\Z)$, interpreted as the pure Schottky product $H_1\ast H_2$, where $H_1 = H_2 = \Z$ is interpreted as acting on $X_1 = X_2 = \R$ by translation. The horizontal lines correspond to copies of $\R$ which correspond to vertices of the form $(1,g)$, while the vertical lines correspond to copies of $\R$ which correspond to vertices of the form $(2,g)$. The intersection points between horizontal and vertical lines are the staples which hold the tree together.}
\label{figureschottkyproductRtree}
\end{center}
\end{figure}
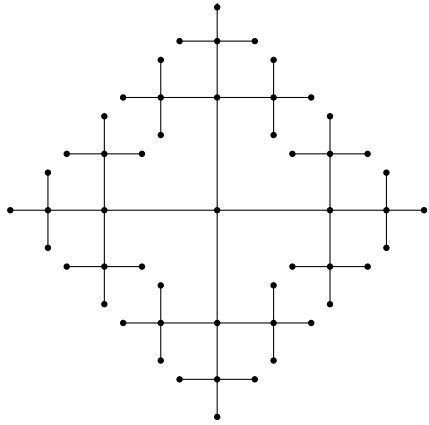

\begin{theorem}[Cf. Figure \ref{figureschottkyproductRtree}]
\label{theoremschottkyproductRtree}
Let $(H_j)_{j\in J}$ be a (possibly infinite) collection of groups and for each $j\in J$ let $\|\cdot\|:H_j\to\Rplus$ be a tree-geometric function. Then the function $\|\cdot\|:G = \ast_{j\in J} H_j \to\Rplus$ defined by
\begin{equation}
\label{schottkyproductRtree}
\|h_1 \cdots h_n\| := \|h_1\| + \cdots + \|h_n\|
\end{equation}
(assuming $h_1\ldots h_n$ is given in reduced form) is a tree-geometric function.
\end{theorem}
\begin{proof}
For each $j\in J$ write $H_j\leq\Isom(X_j)$ and $\|h\| = \dist(\zero_j,h(\zero_j))\all h\in H_j$ for some $\R7$-tree $X_j$ and for some distinguished point $\zero_j\in X_j$. Let $V = J\times G$, and for each $(j,g)\in V$ let $X_v = X_j$. Let
\begin{align*}
E_1 &= \{((j,g),(k,g)) : j\neq k, \; g\in G\}\\
E_2 &= \{((j,g),(j,gh)) : j\in J, \; g\in G,\; h\in H_j\butnot\{e\}\}\\
E &= E_1\cup E_2.
\end{align*}
\begin{claim}
\label{claimcycleE2}
Any cycle in $(V,E)$ is contained in a complete graph of one of the following forms:
\begin{align} \label{form1}
&\{(j,gh): h\in H_j\} \; (j\in J, g\in G \text{ fixed}),\\ \label{form2}
&\{(j,g) : j\in J\} \; (g\in G \text{ fixed}).
\end{align}
In particular, $(V,E)$ is a graph with contractible cycles.
\end{claim}
\begin{subproof}
Let $(v_i)_0^n$ be a cycle in $V$, and for each $i = 0,\ldots,n - 1$ let $e_i = (v_i,v_{i + 1})$. By contradiction suppose that $(v_i)_0^n$ is not contained in a complete graph of one of the forms \eqref{form1},\eqref{form2}. Without loss of generality suppose that $(v_i)_0^n$ is minimal with this property. Then no two consecutive edges $e_i,e_{i + 1}$ can lie in the same set $E_k$. After reindexing if necessary, we find ourselves in the position that $e_i\in E_2$ for $i$ even and $e_i\in E_1$ for $i$ odd. Write $v_0 = (j_1,g)$; then
\[
v_0 = (j_1,g), \; v_1 = (j_1,g h_1), \; v_2 = (j_2,g h_1), \; v_3 = (j_2,g h_1 h_2), \; \text{[etc.]}
\]
with $h_i\in H_{j_i}$, $j_i\neq j_{i + 1}$. Since $G$ is a free product, this contradicts that $v_n = v_0$.
\end{subproof}

For each $(v,w) = ((j,g),(k,g))\in E_1$, we let $\set vw = \{\zero_j\}$ and we let $\bij vw(\zero_j) = \zero_j$. For each $(v,w) = ((j,g),(j,gh))\in E_2$, we let $\set vw = X_j$ and we let $\bij vw = h^{-1}$. Claim \ref{claimcycleE2} then implies the consistency condition. Consider the stapled union $X = \coprod_{(j,g)\in V}^{\text{st}} X_j = \coprod_{(j,g)\in V} X_j / \sim$. Elements of $\coprod_{(j,g)\in V} X_j$ consist of pairs $((j,g),x)$, where $g\in G$ and $x\in X_j$. We will abuse notation by writing $((j,g),x) = (j,g,x)$ and $\lb (j,g),x\rb = \lb j,g,x\rb$. Then the ``staples'' are given by the relations
\begin{align*}
(j,g,\zero_j) &\sim  (k,g,\zero_k) \; [g\in G, \; j,k\in J],\\
\;\; (j,gh,x) & \sim (j,g,h(x)) \; [g\in G, \; j\in J, \; h\in H_j, \; x\in X_j].
\end{align*}
Now consider the following action of $G$ on $\coprod_{(j,g)\in V} X_j$:
\[
g_1\big((j,g_2,x)\big) = (j,g_1g_2,x).
\]
Since the ``staples'' are preserved by this action, it descends to an action on the stapled union $X$. To finish the proof, we need to show that $\dogovar g = \|g\|\all g\in G$, where $\zero = \lb j,e,\zero_j\rb \all j\in J$, and $\|\cdot\|$ is given by \eqref{schottkyproductRtree}. Indeed, fix $g\in G$ and write $g = h_1 \cdots h_n$, where for each $i = 1,\ldots,n$, $h_i\in H_{j_i}\butnot\{e\}$ for some $j\in J$, and $j_i\neq j_{i + 1}\all i$. For each $i = 0,\ldots,n$ let $g_i = h_1 \cdots h_i$, and for each $i = 1,\ldots,n$ let
\[
v_i^{(1)} = (j_i,g_{i - 1}), \; v_i^{(2)} = (j_i,g_i).
\]
Then the sequence $(v_1^{(1)},v_1^{(2)},v_2^{(1)},\ldots,v_n^{(1)},v_n^{(2)})$ is a geodesic whose endpoints are $(j_1,e)$ and $(j_n,g)$. We compute the sequences $(x_i^{(k)})$, $(y_i^{(k)})$ as in Theorem \ref{theoremstapledunion}(ii):
\[
x_i^{(1)} = \zero_{j_i}, \; y_i^{(1)} = \zero_{j_i}, \; x_i^{(2)} = h_i^{-1}(\zero_{j_i}), \; y_i^{(2)} = \zero_{j_i},
\]
It follows that
\[
\dogovar g = \sum_{i = 1}^n \sum_{j = 1}^2 \dist(x_i^{(j)},y_i^{(j)}) = \sum_{i = 1}^n \|h_i\| = \|g\|,
\]
which completes the proof.
\end{proof}

\begin{definition}
\label{definitionschottkyproductRtree}
Let $(H_j)_{j\in J}$ and $G$ be as in Theorem \ref{theoremschottkyproductRtree}. If we write $G\leq\Isom(X)$ and $\|g\| = \dist(\zero,g(\zero))\all g\in G$ for some $\R$-tree $X$ and some distinguished point $\zero\in X$, then we call $(X,G)$ the \emph{pure Schottky product} of $(H_j)_{j\in J}$. (It is readily verified that every pure Schottky product is a Schottky product.)
\end{definition}

\begin{proposition}
\label{propositionschottkyproductRtree}
The Poincar\'e set of a pure Schottky product $H_1\ast H_2$ can be computed by the formula
\[
s\in \Delta(H_1\ast H_2) \;\;\Leftrightarrow\;\; (\Sigma_s(H_1) - 1)(\Sigma_s(H_2) - 1) \geq 1.
\]
\end{proposition}
\begin{proof}
Let
\[
E = (H_1\butnot\{\id\})(H_2\butnot\{\id\}),
\]
so that
\[
G = \bigcup_{n\geq 0} H_2 E^n H_1.
\]
Then by \eqref{schottkyproductRtree}, we have for all $s\geq 0$
\begin{align*}
\Sigma_s(G)
= \sum_{g\in G} e^{-s\dogo g}
&= \sum_{n = 0}^\infty \; \sum_{h_0\in H_2} \; \sum_{g_1,\ldots,g_n\in E} \; \sum_{h_{n + 1}\in H_1} e^{-s[\|h_0\| + \sum_1^n \|g_i\| + \|h_{n + 1}\|]}\\
&= \Sigma_s(H_2) \Sigma_s(H_1) \sum_{n = 0}^\infty \left(\sum_{g\in E} e^{-s\|g\|}\right)^n\\
&= \Sigma_s(H_2) \Sigma_s(H_1) \sum_{n = 0}^\infty \big((\Sigma_s(H_1) - 1)(\Sigma_s(H_2) - 1)\big)^n.
\end{align*}
This completes the proof.
\end{proof}


Proposition \ref{propositionschottkyproductRtree} generalizes to the case of more than two groups as follows:

\begin{proposition}
\label{propositionschottkyproductRtree2}
The Poincar\'e set of a finite pure Schottky product 
\[
G = \ast_{j = 1}^k H_j
\] 
can be computed by the formula
\[
s\in \Delta(H_1\ast H_2) \;\;\Leftrightarrow\;\; \rho(A_s) \geq 1,
\]
where $\rho$ denotes spectral radius, and $A_s$ denotes the matrix whose $(j,j')$th entry is
\[
(A_s)_{j,j'} = \begin{cases}
\Sigma_s(H_j) & j'\neq j\\
0 & j' = j
\end{cases}.
\]
\end{proposition}
\begin{proof}
Let $J = \{1,\ldots,k\}$. Then
\[
G = \bigcup_{n = 0}^\infty \; \bigcup_{\substack{j_1,\ldots,j_n\in J \\ j_1\neq \cdots \neq j_n}} \{h_1\cdots h_n : h_1\in H_{j_1}, \cdots, h_n\in H_{j_n}\}.
\]
So by \eqref{schottkyproductRtree}, we have for all $s\geq 0$
\begin{align*}
\Sigma_s(G)
= \sum_{g\in G} e^{-s\dogo g}
&= \sum_{n = 0}^\infty \; \sum_{\substack{j_1,\ldots,j_n\in J \\ j_1\neq \cdots \neq j_n}} \; \sum_{h_1\in H_{j_1}} \cdots \sum_{h_n\in H_{j_n}} ~~~~ e^{-s\sum_1^n \dogo{h_i}}\\
&= \sum_{n = 0}^\infty \; \sum_{\substack{j_1,\ldots,j_n\in J \\ j_1\neq \cdots \neq j_n}} \; \prod_{i = 1}^n (\Sigma_s(H_{j_i}) - 1)\\
&= 1 + \sum_{n = 1}^\infty [1 \cdots 1] A_s^{n - 1} \left[\begin{array}{l}
\Sigma_s(H_1) - 1\\
\phantom{MM}\vdots\phantom{MM}\\
\Sigma_s(H_n) - 1
\end{array}\right]\\
&\begin{cases}
= \infty & \rho(A_s) \geq 1\\
< \infty & \rho(A_s) < 1
\end{cases}.
\end{align*}
This completes the proof.
\end{proof}
Note that only the last step (the series converges or diverges according to whether or not the spectral radius is at least one) uses the hypothesis that $J$ is finite.

Our last example of an $\R$-tree constructed using the stapling method is similar to the method of pure Schottky products, but differs in important ways:

\begin{example}[Geometric products]
\label{examplegeometricproducts}
Let $Y$ be an $\R$-tree, let $P\subset Y$ be a set, and let $(\Gamma_p)_{p\in P}$ be a collection of abstract groups. Let $\Gamma = \ast_{p\in P} \Gamma_p$. Let $V = \Gamma$, and let
\[
E = \{(\gamma,\gamma\alpha) : \gamma\in \Gamma , \; \alpha\in \Gamma_p\butnot\{e\}\}.
\]
For each $v\in V$, let $X_v = Y$. For each $(v,w) = (\gamma,\gamma\alpha)\in E$, where $\gamma\in\Gamma$ and $\alpha\in \Gamma_p\butnot\{e\}$, we let $\set vw = \{p\}$, and we let $\bij vw(p) = p$. In a manner similar to the proof of Claim \ref{claimcycleE2}, one can check that every cycle in $(V,E)$ is contained in one of the complete graphs $\gamma\Gamma_p \subset V$ ($\gamma\in\Gamma$, $p\in P$), so $(V,E)$ has contractible cycles. The consistency condition is trivial. Thus we can consider the stapled union $X = \coprod_{v\in V}^{\text{st}} X_v$, which admits a natural left action $\phi:\Gamma\to\Isom(X)$:
\[
\iota(\gamma)(\lb v,x\rb) = \lb\gamma v,x\rb.
\]
We let $G = \phi(\Gamma)$, and we call the pair $(X,G)$ the \emph{geometric product} of $Y$ with $(\Gamma_p)_{p\in P}$.
\end{example}

Note that if $(X,G)$ is the geometric product of $Y$ with $(\Gamma_y)_{y\in A}$, then for all $g = (p_1,\gamma_1)\cdots (p_n,\gamma_n)\in G$, we have
\begin{equation}
\label{geometricproducts}
\dogo g = \dist(\zero,p_1) + \sum_{i = 1}^{n - 1} \dist(p_i,p_{i + 1}) + \dist(p_n,\zero).
\end{equation}
To compare this formula with \eqref{schottkyproductRtree}, we observe that if $n = 1$, then we get $\|(a,\gamma)\| = 2\dist(\zero,a)$, so that
\begin{align*}
\dogo{(p_1,\gamma_1)} + \cdots + \dogo{(p_n,\gamma_n)} &= \sum_{i = 1}^n 2\dist(\zero,p_i)\\
& = \dist(\zero,p_1) + \sum_{i = 1}^{n - 1} [\dist(\zero,p_i) + \dist(\zero,p_{i + 1})] + \dist(\zero,p_n).
\end{align*}
So if $(X,G)$ is a geometric product, then the right hand side of \eqref{schottkyproductRtree} exceeds the left hand side by $\sum_{i = 1}^{n - 1} 2\lb p_i|p_{i + 1}\rb_\zero$. The formula \eqref{geometricproducts} is more complicated to deal with because its terms depend on the relation between the neighborhing points $p_i$ and $p_{i + 1}$, rather than just on the individual terms $p_i$. In particular, it is more difficult to compute the Poincar\'e exponent of a geometric product than it is to compute the Poincar\'e exponent of a group coming from Theorem \ref{theoremschottkyproductRtree}. We will investigate the issue of computing Poincar\'e exponents of geometric products in \cite{DSU2}, as well as other topics related to the geometry of these groups.

\begin{figure}
\begin{center}
\begin{tabular}{ll}

\begin{tikzpicture}[line cap=round,line join=round,>=triangle 45,scale=0.85]
\clip(-3.5,-1.0) rectangle (4.41,3.56);
\draw [dashed] (-3.0,3.0)-- (-3.0,3.5);
\draw (-3.0,3.0)-- (3.0,3.0);
\draw (3.0,3.0)-- (3.0,2.5);
\draw (-3.0,2.5)-- (3.0,2.5);
\draw (-3.0,2.5)-- (-3.0,2.0);
\draw (-3.0,2.0)-- (3.0,2.0);
\draw (3.0,2.0)-- (3.0,1.5);
\draw (3.0,1.5)-- (-3.0,1.5);
\draw (-3.0,1.5)-- (-3.0,1.0);
\draw (-3.0,1.0)-- (3.0,1.0);
\draw (3.0,1.0)-- (3.0,0.5);
\draw (3.0,0.5)-- (-3.0,0.5);
\draw (-3.0,0.5)-- (-3.0,0.0);
\draw (-3.0,0.0)-- (3.0,0.0);
\draw [dashed] (3.0,0.0)-- (3.0,-0.5);
\begin{scriptsize}
\draw [fill=black] (-3.0,3.0) circle (1pt);
\draw [fill=black] (3.0,3.0) circle (1pt);
\draw [fill=black] (3.0,2.5) circle (1pt);
\draw [fill=black] (-3.0,2.5) circle (1pt);
\draw [fill=black] (-3.0,2.0) circle (1pt);
\draw [fill=black] (3.0,2.0) circle (1pt);
\draw[color=black] (0.04288653053857866,2.8333178519676443) node {$(Y,bab)$};
\draw[color=black] (0.04288653053857866,2.342033980896144) node {$(Y,ba)$};
\draw[color=black] (0.04288653053857866,1.8358627197921742) node {$(Y,b)$};
\draw [fill=black] (3.0,1.5) circle (1pt);
\draw[color=black] (3.2,1.5) node {$1$};
\draw [fill=black] (-3.0,1.5) circle (1pt);
\draw[color=black] (-3.2,1.5) node {$0$};
\draw [fill=black] (-3.0,1.0) circle (1pt);
\draw [fill=black] (3.0,1.0) circle (1pt);
\draw [fill=black] (3.0,0.5) circle (1pt);
\draw [fill=black] (-3.0,0.5) circle (1pt);
\draw[color=black] (0.027999140506108963,1.339691458688204) node {$(Y,e)$};
\draw[color=black] (0.04288653053857866,0.8384075876167036) node {$(Y,a)$};
\draw[color=black] (0.013111750473639265,0.30246154644779404) node {$(Y,ab)$};
\draw [fill=black] (-3.0,0.0) circle (1pt);
\draw [fill=black] (3.0,0.0) circle (1pt);
\draw[color=black] (0.04288653053857866,-0.15904754455876693) node {$(Y,aba)$};
\end{scriptsize}
\end{tikzpicture}

\begin{tikzpicture}[line cap=round,line join=round,>=triangle 45,scale=0.76]
\clip(-3.5,-1.0) rectangle (4.41,3.56);
\draw (-3.0,3.0)-- (-3.0,3.0);
\draw (-3.0,3.0)-- (-3.0,3.0);
\draw (-3.0,3.0)-- (3.0,2.5);
\draw (3.0,2.5)-- (3.0,2.5);
\draw (-3.0,2.0)-- (3.0,2.5);
\draw (-3.0,2.0)-- (-3.0,2.0);
\draw (-3.0,2.0)-- (3.0,1.5);
\draw (3.0,1.5)-- (3.0,1.5);
\draw (3.0,1.5)-- (-3.0,1.0);
\draw (-3.0,1.0)-- (-3.0,1.0);
\draw (-3.0,1.0)-- (3.0,0.5);
\draw (3.0,0.5)-- (3.0,0.5);
\draw (3.0,0.5)-- (-3.0,0.0);
\draw (-3.0,0.0)-- (-3.0,0.0);
\draw (-3.0,0.0)-- (3.0,-0.5);
\draw (3.0,-0.5)-- (3.0,-0.5);
\draw [dashed] (-3.0,3.0)-- (-1.9817985138773002,3.478813453266429);
\draw [dashed] (3.0,-0.5)-- (1.6209498739803665,-1.0022909465069538);
\begin{scriptsize}
\draw [fill=black] (-3.0,3.0) circle (1pt);
\draw [fill=black] (-3.0,0.0) circle (1pt);
\draw [fill=black] (3.0,0.5) circle (1pt);
\draw [fill=black] (-3.0,1.0) circle (1pt);
\draw [fill=black] (3.0,1.5) circle (1pt);
\draw [fill=black] (-3.0,2.0) circle (1pt);
\draw [fill=black] (3.0,2.5) circle (1pt);
\draw [fill=black] (3.0,2.5) circle (1pt);
\draw [fill=black] (-3.0,2.0) circle (1pt);
\draw [fill=black] (3.0,1.5) circle (1pt);
\draw [fill=black] (-3.0,1.0) circle (1pt);
\draw [fill=black] (3.0,0.5) circle (1pt);
\draw [fill=black] (-3.0,0.0) circle (1pt);
\draw [fill=black] (3.0,-0.5) circle (1pt);
\end{scriptsize}
\end{tikzpicture}

\end{tabular}
\caption[An example of a geometric product]{The geometric product of $Y$ with $(\Gamma_p)_{p\in P}$, where $Y = [0,1]$, $P = \{0,1\}$, $\Gamma_0 = \{e,\gamma_0\} \equiv \Z_2$, and $\Gamma_1 = \{e,\gamma_1\} \equiv \Z_2$. In the left hand picture, copies of $Y$ are drawn as horizontal lines and identifications between points in different copies are drawn as vertical lines. The right hand picture is the result of stapling together certain pairs of points in the left hand picture.}
\end{center}
\end{figure}
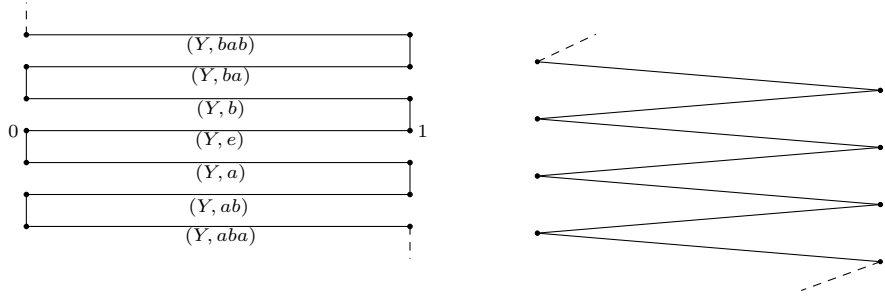

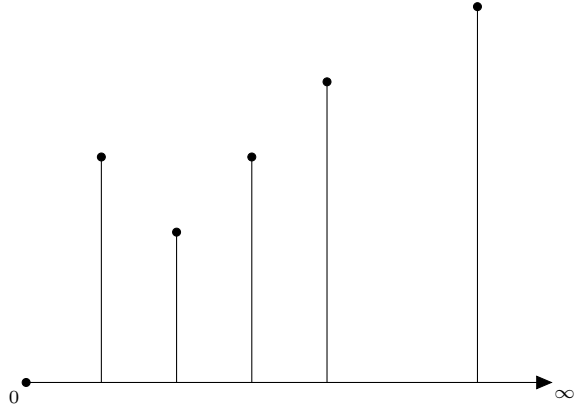
\begin{figure}
\begin{tikzpicture}[line cap=round,line join=round,>=triangle 45,x=1.0cm,y=1.0cm]
\clip(0.48,0.4000000000000001) rectangle (8.56,6.7200000000000015);
\draw[<-] (8.0,1.0)-- (1.0,1.0);
\draw (2.0,4.0)-- (2.0,1.0);
\draw [fill=black] (2.0,4.0) circle (1.5pt);
\draw (3.0,1.0)-- (3.0,3.0);
\draw [fill=black] (3.0,3.0) circle (1.5pt);
\draw (4.0,4.0)-- (4.0,1.0);
\draw [fill=black] (4.0,4.0) circle (1.5pt);
\draw (5.0,1.0)-- (5.0,5.0);
\draw [fill=black] (5.0,5.0) circle (1.5pt);
\draw (7.0,6.0)-- (7.0,1.0);
\draw [fill=black] (7.0,6.0) circle (1.5pt);
\begin{scriptsize}
\draw [fill=black] (1.0,1.0) circle (1.5pt);
\draw[color=black] (0.8400000000000001,0.8000000000000002) node {$0$};
\draw[color=black] (8.16,0.8400000000000002) node {$\infty$};
\end{scriptsize}
\end{tikzpicture}
\caption[Another example of a geometric product]{The set $Y$ of Example \ref{examplegeometricproduct1}. The points at the tops of the vertical lines are ``branch points'' which correspond to fixed points in the geometric product $(X,G)$. If a geodesic in the geometric product is projected down to $Y$, the result will be a sequence of geodesics, each of which starts and ends at one of the indicated points (either $\zero$, an element of $P$, or $\infty$).}
\label{figuregeometricproduct}
\end{figure}

\begin{example}[Cf. Figure \ref{figuregeometricproduct}]
\label{examplegeometricproduct1}
Let $(a_n)_1^\infty$ be an increasing sequence of nonnegative real numbers, and let $(b_n)_1^\infty$ be a sequence of nonnegative real numbers. Let
\[
Y = (\Rplus\times\{0\}) \cup \bigcup_{n = 1}^\infty (\{a_n\}\times [0,b_n])
\]
with the path metric induced from $\R^2$. Let $P = \{p_n : n\in\N\}$, where $p_n = (a_n,b_n)$. Then
\begin{equation}
\label{geometricproduct1}
\dist(p_n,p_m) = b_n + b_m + |a_n - a_m| \all m\neq n,
\end{equation}
so \eqref{geometricproducts} would become
\begin{align*}
\dogo g &= b_1 + a_1 + \sum_{i = 1}^{n - 1} [b_i + b_{i + 1} + |a_{i + 1} - a_i|] + b_n + a_n\\ 
&= \sum_{i = 1}^n 2b_i + a_1 + \sum_{i = 1}^{n - 1} |a_{i + 1} - a_i| + a_n.
\end{align*}
This formula exhibits clearly the fact that the relation between neighborhing points $p_i$ and $p_{i + 1}$ is involved, via the appearance of the term $|a_{i + 1} - a_i|$.
\end{example}

\begin{proposition}
\label{propositionschottkyRtree}
Let $(X,G)$ be the geometric product of $Y$ with $(\Gamma_p)_{p\in P}$, where $P\subset Y$.
\begin{itemize}
\item[(i)] If
\begin{equation}
\label{schottkyRtree1}
\inf\{\dist(y,z) : y,z\in E, y\neq z \} > 0,
\end{equation}
then $G = \lb G_a\rb_{a\in E}$ is a global weakly separated Schottky product. If furthermore
\begin{equation}
\label{schottkyRtree2}
\inf\{\wbar\Dist(y,z) : y,z\in E, y\neq z \} > 0,
\end{equation}
then $G$ is strongly separated.
\item[(ii)] $X$ is proper if and only if all three of the following hold: $Y$ is proper, $\#(\Gamma_a) < \infty$ for all $a\in E$, and $\#(E\cap B(\zero,\rho)) < \infty$ for all $\rho > 0$.
\end{itemize}
\end{proposition}

\begin{proof}[Proof of \text{(i)}]
Suppose that \eqref{schottkyRtree1} holds, and for each $p\in P$, let
\[
U_p = \{\lb g_1\cdots g_n,y\rb \in X : g_1\in G_p\} \cup \{\lb \id,y\rb : y\in B(p,\epsilon)\},
\]
where $\epsilon \leq \inf\{\dist(y,z) : y,z\in P, y\neq z\}/2$. Then $(U_p)_{p\in P}$ a global Schottky system for $G$. If \eqref{schottkyRtree2} also holds, then it is strongly separated, because 
\[
\inf\{\wbar\Dist(U_p,U_q) : p\neq q\} \geq \inf\{\wbar\Dist(y,z) : y,z\in P, y\neq z\} - 2\epsilon
\] 
can be made positive if $\epsilon$ is sufficiently small. Finally, if we go back to assuming only that \eqref{schottkyRtree1} holds, then $(U_p)_{p\in P}$ is still weakly separated, because \eqref{schottkyRtree2} holds for finite subsets.
\end{proof}

\begin{proof}[Proof of \text{(ii)}]
The necessity of these conditions is obvious; conversely, suppose they hold. Fix $\rho > 0$ and $x = \lb g,y\rb\in B_X(\zero,\rho)$; by \eqref{pathmetricv3}, we have
\[
\dist(\zero,p_1) + \dist(p_1,p_2) + \ldots + \dist(p_{n - 1},p_n) + \dist(p_n,y) \leq \rho,
\]
where $g = h_1\cdots h_n$, $h_i\in G_{p_i}\butnot\{\id\}$, $p_i\in P$, and $p_i\neq p_{i + 1}$ for all $i$. It follows that $\dox{p_i} \leq \rho$ for all $i = 1,\ldots,n$, i.e. $p_i\in P\cap B(\zero,\rho)$. In particular, letting $\epsilon = \min_{p,q\in P\cap B(\zero,\rho)} \dist(a,b)$, we have $(n - 1)\epsilon\leq \rho$, or equivalently $n\leq 1 + \rho/\epsilon$. It follows that
\[
g\in\bigcup_{n\leq 1 + \rho/\epsilon} \; \bigcup_{p_1,\ldots,p_n\in P\cap B(\zero,\rho)} (G_{p_1}\butnot\{\id\})\cdots (G_{p_n}\butnot\{\id\}),
\]
a finite set. Thus, $B_X(\zero,\rho)$ is contained in the union of finitely many compact sets of the form $B_Y(\zero,\rho)\times\{g\}\subset X$, and is therefore compact.
\end{proof}

\part{Patterson--Sullivan theory}
\label{partahlforsthurston}
This part will be divided as follows: In Chapter \ref{sectionconformal} we recall the definition of quasiconformal measures, and we prove basic existence and non-existence results. In Chapter \ref{sectionahlforsthurston}, we prove Theorem \ref{theoremahlforsthurstongeneral} (Patterson--Sullivan theorem for groups of divergence type). In Chapter \ref{sectionGFmeasures}, we investigate the geometry of quasiconformal measures of geometrically finite groups, and we prove a generalization of the Global Measure Formula (Theorem \ref{theoremglobalmeasure}) as well as giving various necessary and/or sufficient conditions for the Patterson--Sullivan measure of a geometrically finite group to be doubling (\6\ref{subsectiondoubling}) or exact dimensional (\6\ref{subsectionexactdimensional}).


%

\chapter{Conformal and quasiconformal measures}\label{sectionconformal}

\begin{standingassumption*}
\label{}
Throughout the final part of the monograph, i.e. in Chapters \ref{sectionconformal}-\ref{sectionGFmeasures},  we fix $(X,\dist,\zero,b)$ as in \6\ref{standingassumptions2}, and a \emph{group} $G\leq\Isom(X)$.
\end{standingassumption*}

\section{The definition}
Conformal measures, introduced by S. G. Patterson \cite{Patterson2} and D. P. Sullivan \cite{Sullivan_density_at_infinity}, are an important tool in studying the geometry of the limit set of a Kleinian group. Their definition can be generalized directly to the case of a group acting on a strongly hyperbolic metric space, but for a hyperbolic metric space which is not strongly hyperbolic, a multiplicative error term is required. Thus we make the following definition (cf. \cite[Definition 4.1]{Coornaert}):

\begin{definition}
\label{definitionquasiconformal}
For each $s\geq 0$, a nonzero measure\Footnote{In this monograph, ``measure'' always means ``nonnegative finite Borel measure''.} $\mu$ on $\del X$ is called \emph{$s$-quasiconformal}\Footnote{Not to be confused with the concept of a quasiconformal \emph{map}, cf. \cite{Heinonen2}.} if
\begin{equation}
\label{quasiconformal}
\mu(g(A)) \asymp_\times \int_A [\overline g'(\xi)]^s\;\dee\mu(\xi)
\end{equation}
for every $g\in G$ and for every Borel set $A\subset\del X$. If $X$ is strongly hyperbolic and if equality holds in \eqref{quasiconformal}, then $\mu$ is called \emph{$s$-conformal}.
\end{definition}
\begin{remark}
\label{remarkmuasymp}
For two measures $\mu_1,\mu_2$, write $\mu_1\asymp_\times \mu_2$ if $\mu_1$ and $\mu_2$ are in the same measure class and if the Radon--Nikodym derivative $\dee \mu_1/\dee \mu_2$ is bounded from above and below. Then a measure $\mu$ is $s$-quasiconformal if and only if
\[
\mu\circ g \asymp_\times [\overline g'(\xi)]^s\mu,
\]
and is $s$-conformal if $X$ is strongly hyperbolic and if equality holds.
\end{remark}

\begin{remark}
One might ask whether it is possible to generalize the notions of conformal and quasiconformal measures to semigroups. However, this appears to be difficult. The issue is that the condition \eqref{quasiconformal} is sometimes impossible to satisfy for measures supported on $\Lambda$ -- for example, it may happen that there exist $g_1,g_2\in G$ such that $g_1(\Lambda)\cap g_2(\Lambda) = \emptyset$, in which case letting $A = \del X\butnot\Lambda$ in \eqref{quasiconformal} shows both that $\Supp(\mu)\subset g_1(\Lambda)$ and that $\Supp(\mu)\subset g_2(\Lambda)$, and thus that $\mu = 0$. One may try to fix this by changing the formula \eqref{quasiconformal} somehow, but it is not clear what the details of this should be.
\end{remark}

\section{Conformal measures}
Before discussing quasiconformal measures, let us consider the relation between conformal measures and quasiconformal measures. Obviously, every conformal measure is quasiconformal. In the converse direction we have:

\begin{proposition}
\label{propositionQCtoC}
Suppose that $G$ is countable and that $X$ is strongly hyperbolic. Then for every $s\geq 0$, if $\mu$ is an $s$-quasiconformal measure, then there exists an $s$-conformal measure $\nu$ satisfying $\nu\asymp_\times \mu$.
\end{proposition}
\begin{proof}
For each $g\in G$, let $f_g:\del X\to(0,\infty)$ be a Radon--Nikodym derivative of $\mu\circ g$ with respect to $\mu$. Since $\mu$ is $s$-quasiconformal, we have for $\mu$-a.e. $\xi\in\del X$
\begin{equation}
\label{fgxiasymp}
f_g(\xi) \asymp_\times [g'(\xi)]^s.
\end{equation}
Since $G$ is countable, the set of $\xi\in\del X$ for which \eqref{fgxiasymp} holds for all $g\in G$ is of full $\mu$-measure. In particular, if
\[
f(\xi) = \sup_{g\in G} \frac{f_g(\xi)}{[g'(\xi)]^s},
\]
then $f(\xi)\asymp_\times 1$ for $\mu$-a.e. $\xi\in X$. Now for each $g,h\in G$, the equality $\mu\circ(gh) = (\mu\circ g)\circ h$ implies that
\[
f_{gh}(\xi) = f_g(h(\xi)) f_h(\xi) \text{ for $\mu$-a.e. $\xi\in\del X$.}
\]
Combining with the chain rule for metric derivatives, we have
\[
\frac{f_{gh}(\xi)}{[(gh)'(\xi)]^s} = \frac{f_g(h(\xi))}{[g'(h(\xi))]^s}\frac{f_h(\xi)}{[h'(\xi)]^s} \text{ for $\mu$-a.e. $\xi\in\del X$.}
\]
Note that we are using the strong hyperbolicity assumption here to get equality rather than a coarse asymptotic. Taking the supremum over all $g$ gives
\[
f(\xi) = f(h(\xi))\frac{f_h(\xi)}{[h'(\xi)]^s} \text{ for $\mu$-a.e. $\xi\in\del X$.}
\]
We now claim that $\nu := f\mu$ is an $s$-conformal measure. Indeed,
\[
\frac{\dee\nu\circ g}{\dee\nu}(\xi)
= \frac{f(g(\xi))}{f(\xi)}\frac{\dee\mu\circ g}{\dee\mu}(\xi)
= \frac{f(g(\xi))}{f(\xi)}f_g(\xi)
= [g'(\xi)]^s \text{ for $\mu$-a.e. $\xi\in\del X$.}
\]
\end{proof}

\section{Ergodic decomposition}

Let $\MM(\del X)$ denote the set of all measures on $\del X$, and let $\MM_1(\del X)$ denote the set of all probability measures on $\del X$.

\begin{definition}
\label{definitionergodic}
A measure $\mu\in\MM(\del X)$ is \emph{ergodic} if for every $G$-invariant Borel set $A\subset\del X$, we have $\mu(A) = 0$ or $\mu(\del X\butnot A) = 0$.
\end{definition}

It is often useful to be able to write a non-ergodic measure as the convex combination of ergodic measures. To make this rigorous, suppose that $X$ is complete and separable, so that $\bord X$ and $\del X$ are Polish spaces. Then $\del X$ together with its Borel $\sigma$-algebra forms a standard Borel space. Let $\BB$ denote the smallest $\sigma$-algebra on $\MM(\del X)$ with the following property:
\begin{property}
\label{propertymeasurablestructure}
For every bounded Borel-measurable function $f:\del X\to\R$, the function
\[
\mu \mapsto \int f\;\dee\mu
\]
is a $\BB$-measurable map from $\MM(\del X)$ to $\R$.
\end{property}
Then $(\MM(\del X),\BB)$ is a standard Borel space. We may now state the following theorem:

\begin{proposition}[Ergodic decomposition of quasiconformal measures]
\label{propositionergodicdecomposition}
We suppose that $G$ is countable and that $X$ is separable. 
Fix $s\geq 0$.
\begin{itemize}
\item[(i)] For every $s$-quasiconformal measure $\mu$, there is a measure $\what\mu$ on $\MM_1(\del X)$ which satisfies
\begin{equation}
\label{ergodicdecomposition}
\mu(A) = \int \nu(A)\;\dee\what\mu(\nu) \text{ for every Borel set $A\subset\del X$}
\end{equation}
and gives full measure to the set of ergodic $s$-quasiconformal measures.\Footnote{If $A$ is a non-measurable set, then a measure $\mu$ gives full measure to $A$ if and only if $A$ contains a measurable set of full $\mu$-measure. Thus we do not need to check whether or not the set of ergodic $s$-quasiconformal measures is a measurable set in $\MM_1(\del X)$.}
\item[(ii)] If $X$ is strongly hyperbolic, then for every $s$-conformal measure $\mu$, there is a unique measure $\what\mu$ on $\MM(\del X)$ which satisfies \eqref{ergodicdecomposition} and which gives full measure to the set of ergodic $s$-conformal measures.
\end{itemize}
\end{proposition}
\begin{remark}
Note that we have uniqueness in (ii) but not in (i).
\end{remark}
\begin{proof}[Proof of Proposition \ref{propositionergodicdecomposition}]
Both cases of the proposition are essentially special cases of \cite[Theorem 1.4]{GreschonigSchmidt}, as we now demonstrate:
\begin{itemize}
\item[(i)]
Let $\mu$ be an $s$-quasiconformal measure. Let $\varrho:G\times\del X\to\R$ satisfy \cite[(1.1)-(1.3)]{GreschonigSchmidt}. Then by \cite[Theorem 1.4]{GreschonigSchmidt}, there is a measure $\what\mu$ satisfying \eqref{ergodicdecomposition} supported on the set of ergodic probability measures which are ``$\varrho$-admissible'' (in the terminology of \cite{GreschonigSchmidt}). But by \cite[(1.1)]{GreschonigSchmidt}, we have $b^{\varrho(g,\xi)}\asymp_\times \overline g'(\xi)^s$ for $\mu$-a.e. $\xi\in\del X$, say for all $\xi\in\del X\butnot S$, where $\mu(S) = 0$. Then every $\varrho$-admissible measure $\nu$ satisfying $\nu(S) = 0$ is $s$-quasiconformal. But by \eqref{ergodicdecomposition}, $\nu(S) = 0$ for $\what\mu$-a.e. $\nu$, so $\what\mu$-a.e. $\nu$ is $s$-quasiconformal.
\item[(ii)] Let $\mu$ be an $s$-conformal measure. Let $\varrho:G\times\del X\to\R$ satisfy \cite[(1.1)-(1.3)]{GreschonigSchmidt}. Then by \cite[(1.1)]{GreschonigSchmidt}, we have $b^{\varrho(g,\xi)} = g'(\xi)^s$ for $\mu$-a.e. $\xi\in\del X$, say for all $\xi\in\del X\butnot S$, where $\mu(S) = 0$. Then for every measure $\nu$ satisfying $\nu(S) = 0$, $\nu$ is $\varrho$-admissible \emph{if and only if} $\nu$ is $s$-conformal. By \cite[Theorem 1.4]{GreschonigSchmidt}, there is a unique measure $\what\mu$ satisfying \eqref{ergodicdecomposition} supported on the set of $\varrho$-admissible ergodic probability measures; such a measure is also unique with respect to satisfying \eqref{ergodicdecomposition} being supported on the set of $s$-conformal ergodic measures.
\end{itemize}
\end{proof}

\begin{corollary}
Suppose that $G$ is countable and that $X$ is separable, and fix $s\geq 0$. If there is an $s$-(quasi)conformal measure, then there is an ergodic $s$-(quasi)conformal measure.
\end{corollary}

In the sequel, we will be concerned with when an $s$-quasiconformal measure is unique up to coarse asymptotic. This is closely connected with ergodicity:
\begin{proposition}
\label{propositionquasiconformaluniqueness}
Suppose that $G$ is countable and that $X$ is separable, and fix $s\geq 0$. Suppose that there is an $s$-quasiconformal measure $\mu$. The following are equivalent:
\begin{itemize}
\item[(A)] $\mu$ is unique up to coarse asymptotic i.e. $\mu\asymp_\times \w\mu$ for any $s$-quasiconformal measure $\w\mu$.
\item[(B)] Every $s$-quasiconformal measure is ergodic.
\end{itemize}
If in addition $X$ is strongly hyperbolic, then (A)-(B) are equivalent to
\begin{itemize}
\item[(C)] There is exactly one $s$-conformal probability measure.
\end{itemize}
\end{proposition}
\begin{proof}[Proof of \text{(A) \implies (B)}]
If $\mu$ is a non-ergodic $s$-quasiconformal measure, then there exists a $G$-invariant set $A\subset\del X$ such that $\mu(A),\mu(\del X\butnot A) > 0$. But then $\nu_1 = \mu\given A$ and $\nu_2 = \mu\given\del X\butnot A$ are non-asymptotic $s$-quasiconformal measures, a contradiction.
\end{proof}
\begin{proof}[Proof of \text{(B) \implies (A)}]
Suppose that $\mu_1,\mu_2$ are two $s$-quasiconformal measures. Then the measure $\mu = \mu_1 + \mu_2$ is also $s$-quasiconformal, and therefore ergodic. Let $f_i$ be a Radon--Nikodym derivative of $\mu_i$ with respect to $\mu$. Then for all $g\in G$,
\begin{equation}
\label{fiinvariant}
f_i \circ g(\xi) = \frac{\dee \mu_i\circ g}{\dee \mu\circ g}(\xi) \asymp_\times \frac{[\overline g'(\xi)]^s}{[\overline g'(\xi)]^s} \frac{\dee\mu_i}{\dee\mu}(\xi) = f_i(\xi) \text{ for $\mu$-a.e. $\xi\in\del X$}.
\end{equation}
It follows that
\[
h_i(\xi) := \sup_{g\in G} f_i\circ g(\xi) \asymp_\times f_i(\xi) \text{ for $\mu$-a.e. $\xi\in\del X$}.
\]
But the functions $h_i$ are $G$-invariant, so since $\mu$ is ergodic, they are constant $\mu$-a.e., say $h_i = c_i$. It follows that $\mu_i \asymp_\times c_i \mu$; since $\mu_i \neq 0$, we have $c_i > 0$ and thus $\mu_1\asymp_\times \mu_2$.
\end{proof}
\begin{proof}[Proof of \text{(B) \implies (C)}]
The existence of an $s$-conformal measure is guaranteed by Proposition \ref{propositionQCtoC}. If $\mu_1,\mu_2$ are two $s$-conformal measures, then the Radon--Nikodym derivatives $f_i = \dee\mu_i/\dee(\mu_1 + \mu_2)$ satisfy \eqref{fiinvariant} with equality, so $f_i = c_i$ for some constants $c_i$. It follows that $\mu_1 = (c_1/c_2) \mu_2$, and so if $\mu_1,\mu_2$ are probability measures then $\mu_1 = \mu_2$.
\end{proof}
\begin{proof}[Proof of \text{(C) \implies (A)}]
Follows immediately from Proposition \ref{propositionQCtoC}.
\end{proof}

\section{Quasiconformal measures}

We now turn to the deeper question of when a quasiconformal measure exists in the first place. To approach this question we begin with a fundamental geometrical lemma about quasiconformal measures:

\begin{lemma}[Sullivan's Shadow Lemma, cf. {\cite[Proposition 3]{Sullivan_density_at_infinity}}, {\cite[\61.1]{Roblin2}}]
\label{lemmasullivanshadow}
Fix $s\geq 0$, and let $\mu$ be a $s$-quasiconformal measure on $\del X$ which is not a pointmass. Then for all $\sigma > 0$ sufficiently large and for all $g\in G$,
\[
\mu(\Shad(g(\zero),\sigma)) \asymp_{\times,\sigma,\mu} b^{-s \dogo g}.
\]
\end{lemma}
\begin{proof}
We have
\begin{align*}
\mu(\Shad(g(\zero),\sigma))
&\asymp_{\times,\mu} \int_{g^{-1}(\Shad(g(\zero),\sigma))}\big(\overline g'\big)^s \;\dee\mu\\
\by{the definition of $s$-quasiconformality}\\
&=_{\phantom{\times,\mu}} \int_{\Shad_{g^{-1}(\zero)}(\zero,\sigma)}\big(\overline g'\big)^s \;\dee\mu\\
&\asymp_{\times,\sigma} \int_{\Shad_{g^{-1}(\zero)}(\zero,\sigma)}b^{-s\dogo g} \;\dee\mu\\
\by{the Bounded Distortion Lemma \ref{lemmaboundeddistortion}}\\
&=_{\phantom{\times,\mu}} b^{-s\dogo g}\mu\big(\Shad_{g^{-1}(\zero)}(\zero,\sigma)\big).
\end{align*}
Thus, to complete the proof, it is enough to show that
\[
\mu\big(\Shad_{g^{-1}(\zero)}(\zero,\sigma)\big) \asymp_{\times,\mu,\sigma} 1,
\]
assuming $\sigma$ is sufficiently large (depending on $\mu$). The upper bound is automatic since $\mu$ is finite. Now, since by assumption $\mu$ is not a pointmass, we have $\#(\Supp(\mu))\geq 2$. Choose distinct $\xi_1,\xi_2\in\Supp(\mu)$, and let $\epsilon = \Dist(\xi_1,\xi_2)/3$. By the Big Shadows Lemma \ref{lemmabigshadow}, we have
\[
\Diam(\del X\butnot \Shad_{g^{-1}(\zero)}(\zero,\sigma)) \leq \epsilon
\]
for all $\sigma > 0$ sufficiently large (independent of $g$). Now since 
$$\Dist(B(\xi_1,\epsilon),B(\xi_2,\epsilon))\geq \epsilon,$$ it follows that
\[
\exists i = 1,2 \;\; B(\xi_i,\epsilon)\subset \Shad_{g^{-1}(\zero)}(\zero,\sigma)
\]
and thus
\[
\mu\big(\Shad_{g^{-1}(\zero)}(\zero,\sigma)\big) \geq \min_{i = 1}^2 \mu\big(B(\xi_i,\epsilon)\big) > 0.
\]
The right hand side is independent of $g$, which completes the proof.
\end{proof}

Sullivan's Shadow Lemma suggests that in the theory of quasiconformal measures, there is a division between those measures which are pointmasses and those which are not. Let us first consider the easier case of a pointmass quasiconformal measure, and then move on to the more interesting theory of non-pointmass quasiconformal measures.

\subsection{Pointmass quasiconformal measures}
\begin{proposition}
\label{propositionquasiconformalpointmass}
A pointmass $\delta_\xi$ is $s$-quasiconformal if and only if
\begin{itemize}
\item[(I)] $\xi\in\del X$ is a global fixed point of $G$, and
\item[(II)] either
\begin{itemize}
\item[(IIA)] $\xi$ is neutral with respect to every $g\in G$, or
\item[(IIB)] $s = 0$.
\end{itemize}
\end{itemize}
\end{proposition}
\begin{proof} To begin we recall that $g'(\xi)$ denotes the dynamical derivative, cf. Proposition \ref{propositiondynamicalderivative}.
For each $\xi\in\del X$,
\begin{align*} 
\text{$\delta_\xi$ is $s$-quasiconformal}
&\Leftrightarrow \delta_\xi\circ g \asymp_\times (\overline g')^s\delta_\xi \all g\in G\\
&\Leftrightarrow \text{$g(\xi) = \xi$ and $[\overline g'(\xi)]^s \asymp_\times 1$} \all g\in G\\
&\Leftrightarrow \text{$g(\xi) = \xi$ and $[g'(\xi)]^s = 1$} \all g\in G\\
&\Leftrightarrow \text{$g(\xi) = \xi$ and ($g'(\xi) = 1$ or $s = 0$)} \all g\in G.\noreason
\end{align*}
\end{proof}

\begin{corollary}
\label{corollaryquasiconformalpointmass}
~
\begin{itemize}
\item[(i)] If $G$ is of general type, then no pointmass is $s$-quasiconformal for any $s\geq 0$.
\item[(ii)] If $G$ is loxodromic, then no pointmass is $s$-quasiconformal for any $s > 0$.
\end{itemize}
\end{corollary}

\subsection{Non-pointmass quasiconformal measures}
Next we will ask the following question: Given a group $G$, for what values of $s$ does a non-pointmass quasiconformal measure exist, and when is it unique up to coarse asymptotic? We first recall the situation in the Standard Case, where the answers are well-known. The first result is the Patterson--Sullivan theorem \cite[Theorem 1]{Sullivan_density_at_infinity}, which states that any discrete subgroup $G\leq\Isom(\H^d)$ admits a $\delta_G$-conformal measure supported on $\LambdaG$. It is unique up to a multiplicative constant if $G$ is of divergence type (\cite[Theorem 8.3.5]{Nicholls} together with Proposition \ref{propositionquasiconformaluniqueness}). The next result is negative, stating that if $s < \delta_G$, then $G$ admits no non-pointmass $s$-conformal measure. From these results and from Corollary \ref{corollaryquasiconformalpointmass}, it follows that if $G$ is of general type, then $\delta_G$ is the infimum of $s$ for which there exists an $s$-conformal measure \cite[Corollary 4]{Sullivan_density_at_infinity}. Finally, for $s > \delta_G$, an $s$-conformal measure on $\LambdaG$ exists if and only if $G$ is not convex-cocompact (\cite[Theorem 4.1]{AFT} for $\Leftarrow$, \cite[Theorem 4.4.1]{Nicholls} for $\Rightarrow$); no nontrivial conditions are known which guarantee uniqueness in this case.

We now generalize the above results to the setting of hyperbolic metric spaces, replacing the Poincar\'e exponent $\delta_G$ with the modified Poincar\'e exponent $\w\delta_G$, and the notion of divergence type with the notion of generalized divergence type. By Proposition \ref{propositionbasicmodified}(ii), our theorems will reduce to the known results in the case of a strongly discrete group.

We begin with the negative result, as its proof is the easiest:
\begin{proposition}[cf. {\cite[p.178]{Sullivan_density_at_infinity}}]
\label{propositionlessthandelta}
For any $s < \w\delta_G$, there does not exist a non-pointmass $s$-quasiconformal measure.
\end{proposition}
\begin{proof}
By contradiction, suppose that $\mu$ is a non-pointmass $s$-quasiconformal measure. Let $\sigma > 0$ be large enough so that Sullivan's Shadow Lemma \ref{lemmasullivanshadow} holds, and let $\tau > 0$ be the implied constant of \eqref{distbusemann} from the Intersecting Shadows Lemma \ref{lemmatau}. Let $S_{\tau + 1}$ be a maximal $(\tau + 1)$-separated subset of $G(\zero)$. Fix $n\in\N$, and let $A_n$ be the $n$th annulus $A_n = B(\zero,n)\butnot B(\zero,n - 1)$. Now by the Intersecting Shadows Lemma \ref{lemmatau}, the shadows $\big(\Shad(x,\sigma)\big)_{x\in S_{\tau + 1}\cap A_n}$ are disjoint, and so by Sullivan's Shadow Lemma \ref{lemmasullivanshadow}
\begin{align*}
1 \asymp_{\times,\mu} \mu(\del X)
&\geq
\sum_{x\in S_{\tau + 1}\cap A_n}\mu(\Shad(x,\sigma))\\
&\asymp_{\times,\sigma,\mu} \sum_{x\in S_{\tau + 1}\cap A_n}b^{-s\dox x}\\
& \asymp_{\times}
b^{-sn}\#(S_{\tau + 1}\cap A_n).
\end{align*}
Thus for all $t > s$,
\[
\Sigma_t(S_{\tau + 1})
\asymp_\times \sum_{n\in\N} b^{-tn}\#(S_{\tau + 1}\cap A_n)
\lesssim_{\times,\sigma,\mu} \sum_{n\in\N}b^{(s - t)n} < \infty.
\]
But this implies that $\w\delta_G\leq t$ (cf. \eqref{modifiedpoincaredef}); letting $t\searrow s$ gives $\w\delta_G\leq s$, contradicting our hypothesis.
\end{proof}

\begin{remark}
\label{remarkupperorbitalbound}
The above proof shows that if there exists a non-pointmass $\w\delta$-conformal measure, then
\[
\#(S_{\tau + 1}\cap A_n) \lesssim_\times b^{\w\delta n} \all n\geq 1.
\]
In particular, if $\w\delta > 0$ then summing over $n = 1,\ldots,N$ gives
\[
\#(S_{\tau + 1}\cap B(\zero,N)) \lesssim_\times b^{\w\delta N} \all n\geq 1.
\]
If $G$ is strongly discrete, then for all $\rho > 0$,
\begin{align*}
\NN_{X,G}(\rho) = \#\{g\in G : \dogo g \leq \rho\} 
&\lesssim_\times \#(S_{\tau + 1}\cap B(\zero,\rho + \tau + 1))\\
&\lesssim_\times b^{\delta\lceil\rho + \tau + 1\rceil}\\
&\asymp_\times b^{\delta\rho}.
\end{align*}
The bound $\NN_{X,G}(\rho) \lesssim_\times b^{\delta\rho}$ in fact holds without assuming the existence of a $\delta$-conformal measure; see Corollary \ref{corollaryupperorbitalbound}.
\end{remark}

Next we study hypotheses which guarantee the existence of a $\w\delta_G$-quasiconformal measure. In particular, we will show that if $\w\delta_G < \infty$ and if $G$ is of compact type or of generalized divergence type, then there exists a $\w\delta_G$-quasiconformal measure. The first case we consider now, while the case of a group of generalized divergence type will be considered in Chapter \ref{sectionahlforsthurston}.

\draftnewpage
\begin{theorem}[cf. {\cite[Th\'eor\`eme 5.4]{Coornaert}}]
\label{theorempattersonsullivangeneral}
Assume that $G$ is of compact type and that $\w\delta < \infty$. Then there exists a $\w\delta$-quasiconformal measure supported on $\Lambda$. If $X$ is strongly hyperbolic, then there exists a $\w\delta$-conformal measure supported on $\Lambda$.
\end{theorem}
\begin{remark}
Any group acting on a proper geodesic hyperbolic metric space is of compact type, so Theorem \ref{theorempattersonsullivangeneral} includes the case of proper geodesic hyperbolic metric spaces.
\end{remark}

\begin{remark}
Combining Theorem \ref{theorempattersonsullivangeneral} with Proposition \ref{propositionlessthandelta} and Corollary \ref{corollaryquasiconformalpointmass} shows that for $G$ nonelementary of compact type,
\[
\w\delta = \inf\{s > 0 : \text{there exists an $s$-quasiconformal measure supported on $\Lambda$}\},
\]
thus giving another geometric characterization of $\w\delta$ (the first being Theorem \ref{theorembishopjonesmodified}).
\end{remark}

Before proving Theorem \ref{theorempattersonsullivangeneral}, we recall the following lemma due to Patterson:

\begin{lemma}[{\cite[Lemma 3.1]{Patterson2}}]
\label{lemmapatterson}
Let $\AA = (a_n)_1^\infty$ be a sequence of positive real numbers, and let
\[
\delta = \delta(\AA) = \inf\left\{s \geq 0 : \sum_{n = 1}^\infty a_n^{-s} < \infty\right\}.
\]
Then there exists an increasing continuous function $k:(0,\infty)\to(0,\infty)$ such that:
\begin{itemize}
\item[(i)] The series
\[
\Sigma_{s,k}(\AA) = \sum_{n = 1}^\infty k(a_n) a_n^{-s}
\]
converges for $s > \delta$ and diverges for $s\leq \delta$.
\item[(ii)] There exists a decreasing function $\epsilon:(0,\infty)\to(0,\infty)$ such that for all $y > 0$ and $x > 1$,
\begin{equation}
\label{patterson}
k(xy) \leq x^{\epsilon(y)} k(y),
\end{equation}
and such that $\lim_{y\to\infty}\epsilon(y) = 0$.
\end{itemize}
\end{lemma}

\begin{proof}[Proof of Theorem \ref{theorempattersonsullivangeneral}]
By Proposition \ref{propositionbasicmodified}, there exist $\rho > 0$ and a maximal $\rho$-separated set $S_\rho\subset G(\zero)$ such that $\w\delta(G) = \delta(S_\rho)$; moreover, this $\rho$ may be chosen large enough so that $S_{\rho/2}$ does not contain a bounded infinite set, where $S_{\rho/2}$ is a $\rho/2$-separated set. Let $\AA = (a_n)_1^\infty$ be any indexing of the sequence $(b^{\dox x})_{x\in S_\rho}$, and let $k:(0,\infty)\to(0,\infty)$ be the function given by Lemma \ref{lemmapatterson}. For shorthand let
\begin{align*}
k(x) &= k(b^{\dox x})\\
\epsilon(x) &= \epsilon(b^{\dox x})\\
\Sigma_{s,k} &= \Sigma_{s,k}(\AA) = \sum_{x\in S_\rho} k(x) b^{-s\dox x}.
\end{align*}
Then $\Sigma_{s,k} < \infty$ if and only if $s > \w\delta$; moreover, the function $s\mapsto \Sigma_{s,k}$ is continuous. For each $s > \w\delta_G$, let
\begin{equation}
\label{musdef}
\mu_s = \frac{1}{\Sigma_{s,k}}\sum_{x\in S_\rho}k(x) b^{-s\dox x}\delta_x\in\MM_1(S_\rho\cup\LambdaG).
\end{equation}
Now since $G$ is of compact type, the set $S_\rho\cup\LambdaG$ is compact (cf. (B) of Proposition \ref{propositioncompacttype}). Thus by the Banach--Alaoglu theorem, the set $\MM_1(S_\rho\cup\LambdaG)$ is compact in the weak-* topology. So there exists a sequence $s_n\searrow\w\delta$ so that if we let $\mu_n = \mu_{s_n}$, then $\mu_n\to\mu\in\MM_1(S_\rho\cup\LambdaG)$. We will show that $\mu$ is $\w\delta_G$-quasiconformal and that $\Supp(\mu) = \LambdaG$.

\begin{claim}
\label{claimboundarysupport}
$\Supp(\mu) \subset\Lambda$.
\end{claim}
\begin{subproof}
Fix $R > 0$. Since $\delta(S_\rho) < \infty$, we have $\#(S_\rho\cap B(\zero,R)) < \infty$. Thus,
\begin{align*}
\mu(B(\zero,R)) \leq \limsup_{s\searrow\w\delta} \mu_s(B(\zero,R))
&\leq \limsup_{s\searrow\w\delta} \frac{\#(S_\rho\cap B(\zero,R)) k(b^R) b^{-\w\delta R}}{\Sigma_{s,k}}\\
&= \frac{\#(S_\rho\cap B(\zero,R)) k(b^R) b^{-\w\delta R}}{\infty} = 0.
\end{align*}
Letting $R\to\infty$ shows that $\mu(X) = 0$; thus $\Supp(\mu)\subset \cl{S_\rho\cup\Lambda\butnot X} = \Lambda$.
\end{subproof}

To complete the proof, we must show that $\mu$ is $\w\delta$-quasiconformal. Fix $g\in G$, and let
\[
\nu_g = [(\overline g')^{\w\delta} \mu] \circ g^{-1}.
\]
We want to show that $\nu_g \asymp_\times\mu$.

\begin{claim}
\label{claimconformalf}
For every continuous function $f:\bord X\to(0,\infty)$, we have
\begin{equation}
\label{conformalf}
\int f\;\dee \nu_g \asymp_\times \int f\;\dee\mu.
\end{equation}
\end{claim}
\begin{subproof}
Since $S_\rho\cup \LambdaG$ is compact, $\log_b(f)$ is uniformly continuous on $S_\rho\cup\LambdaG$ with respect to the metric $\wbar\Dist$. Let $\phi_f$ denote the modulus of continuity of $\log_b(f)$, so that
\begin{equation}
\label{uniformcontinuity}
\wbar\Dist(x,y) \leq r \;\;\Rightarrow\;\; \frac{f(x)}{f(y)} \leq b^{\phi_f(r)} \all x,y\in S_\rho\cup\LambdaG.
\end{equation}
For each $n\in\N$ let
\[
\nu_{g,n} = [(\overline g')^{s_n} \mu_n] \circ g^{-1},
\]
so that $\nu_{g,n} \tendsto{n,\times} \nu$.
Then
\begin{align*}
\nu_{g,n} &=_\pt \frac{1}{\Sigma_{s_n,k}} \sum_{x\in S_\rho} k(x) b^{-s_n \dox x} [(\overline g')^{s_n} \delta_x] \circ g^{-1}\\
&\asymp_\times \frac{1}{\Sigma_{s_n,k}} \sum_{x\in S_\rho} k(x) b^{-s_n \dox x} b^{s_n[\dox x - \dox{g(x)}]} [\delta_x \circ g^{-1}]\\
&=_\pt \frac{1}{\Sigma_{s_n,k}} \sum_{x\in S_\rho} b^{-s_n \dox{g(x)}} k(x) \delta_{g(x)}\\
&=_\pt \frac{1}{\Sigma_{s_n,k}} \sum_{x\in g(S_\rho)} b^{-s_n \dox x} k(g^{-1}(x)) \delta_x,
\end{align*}
and so
\begin{equation}
\label{nugnmun}
\frac{\int f\;\dee\nu_{g,n}}{\int f\;\dee\mu_n} \asymp_\times \frac{\sum_{x\in g(S_\rho)}b^{-s_n\dox x} k(g^{-1}(x)) f(x)}{\sum_{y\in S_\rho}b^{-s_n\dox y} k(y) f(y)}\cdot
\end{equation}
For each $x\in g(S_\rho)\subset G(\zero)$, there exists $y_x\in S_\rho$ such that $\dist(x,y_x) \leq \rho$.
\begin{observation}
$\#\{x : y_x = y\}$ is bounded independent of $y$ and $g$.
\end{observation}
\begin{subproof}
Write $y = h(\zero)$; then
\[
\#\{x : y_x = y\} \leq \#(g(S_\rho)\cap B(y,\rho)) = \#(h^{-1}g(S_\rho)\cap B(\zero,\rho)).
\]
But $S_\rho' := h^{-1}g(S_\rho)$ is a $\rho$-separated set. For each $x\in S_\rho'$, choose $z_x\in S_{\rho/2}$ such that $\dist(x,z_x) < \rho/2$; then the map $x\mapsto z_x$ is injective, so
\[
\#(S_\rho') \leq \#(S_{\rho/2}\cap B(\zero,2\rho)),
\]
which is bounded independent of $y$ and $g$.
\end{subproof}
Now
\[
\wbar\Dist(x,y_x) \leq b^{-\lb x|y_x\rb_\zero} \leq b^{\rho - \dox{y_x}};
\]
applying \eqref{uniformcontinuity} gives
\[
f(x) \leq b^{\phi_f(b^{\rho - \dox{y_x}})} f(y_x).
\]
On the other hand, by \eqref{patterson} we have
\[
k(g^{-1}(x)) \leq b^{\epsilon(y_x) [\rho + \dogo g]} k(y_x),
\]
and we also have
\[
b^{-s_n \dox x} \leq b^{s_n \rho} b^{-s_n \dox{y_x}}.
\]
Combining everything gives
\begin{align*}
&\sum_{x\in g(S_\rho)}b^{-s_n\dox x} k(g^{-1}(x)) f(x)\\
&\leq_\pt \sum_{x\in g(S_\rho)} \exp_b\left(s_n \rho + \epsilon(y_x)[\rho + \dogo g] + \phi_f(b^{\rho - \dox{y_x}})\right) b^{-s_n\dox{y_x}} k(y_x) f(y_x)\\
&\lesssim_\times \sum_{y\in S_\rho} \exp_b\left(\epsilon(y)[\rho + \dogo g] + \phi_f(b^{\rho - \dox y})\right) b^{-s_n\dox y} k(y) f(y),
\end{align*}
and taking the limit as $n\to\infty$ we have
\[
\int f(x) \;\dee \nu(x) \lesssim_\times \int \exp_b\left(\epsilon(y)[\rho + \dogo g] + \phi_f(b^{\rho - \dox y})\right) f(y)\;\dee\mu(y) = \int f(y)\;\dee\mu(y)
\]
since $\phi_f(b^{\rho - \dox y}) = \epsilon(y) = 0$ for all $y\in\del X$. A symmetric argument gives the converse direction.
\end{subproof}
Now let $C$ be the implied constant of \eqref{conformalf}. Then for every continuous function $f:X\to (0,\infty)$,
\[
C\int f\;\dee\nu - \int f\;\dee\mu \geq 0 \text{ and }C\int f \;\dee \mu - \int f\;\dee\nu \geq 0,
\]
i.e. the linear functionals $I_1[f] = C\int f\;\dee\nu - \int f\;\dee\mu$ and $I_2[f] = C\int f \;\dee \mu - \int f\;\dee\nu$ are positive. Thus by the Riesz representation theorem, there exist measures $\gamma_1,\gamma_2$ such that $I_{\gamma_i} = I_i$ ($i = 1,2$). The uniqueness assertion of the Riesz representation theorem then guarantees that
\begin{equation}
\label{riesz}
\gamma_1 + \mu = C\nu \text{ and } \gamma_2 + \nu = C\mu.
\end{equation}
In particular, $C\nu \geq \mu$, and $C\mu \geq \nu$. This completes the proof.
\end{proof}
%

\chapter{Patterson--Sullivan theorem for groups of divergence type}
\label{sectionahlforsthurston}

In this chapter, we prove Theorem \ref{theoremahlforsthurstongeneral}, which states that a nonelementary group of generalized divergence type possesses a $\w\delta$-quasiconformal measure.

\section{Samuel--Smirnov compactifications}
We begin by summarizing the theory of Samuel--Smirnov compactifications, which will be used in the proof of Theorem \ref{theoremahlforsthurstongeneral}.

\begin{proposition}
\label{propositionsamuelsmirnov}
Let $(Z,\Dist)$ be a complete metric space. Then there exists a compact Hausdorff space $\what Z$ together with a homeomorphic embedding $\iota:Z\to\what Z$ with the following property:
\begin{property}
\label{propertysamuelsmirnov}
\narrower
If $A,B\subset Z$, then $\cl A\cap \cl B \neq \emptyset$ if and only if $\Dist(A,B) = 0$. Here $\cl A$ and $\cl B$ denote the closures of $A$ and $B$ relative to $\what Z$.
\end{property}
\noindent The pair $(\what Z,\iota)$ is unique up to homeomorphism. Moreover, if $Z_1,Z_2$ are two complete metric spaces and if $f:Z_1\to Z_2$ is uniformly continuous, then there exists a unique continuous map $\fhat:\what Z_1\to \what Z_2$ such that $\iota\circ f = \fhat\circ \iota$. The reverse is also true: if $f$ admits such an extension, then $f$ is uniformly continuous.
\end{proposition}
The space $\what Z$ will be called the \emph{Samuel--Smirnov compactification} of $Z$.
\begin{proof}[Proof of Proposition \ref{propositionsamuelsmirnov}]
The metric $\Dist$ induces a proximity on $Z$ in the sense of \cite[Definition 1.7]{NaimpallyWarrack}. Then the existence and uniqueness of a pair $(\what Z,\iota)$ for which Property \ref{propertysamuelsmirnov} holds is guaranteed by \cite[Theorem 7.7]{NaimpallyWarrack}. The assertions concerning uniformly continuous maps follow from \cite[Theorem 7.10]{NaimpallyWarrack} and \cite[Theorem 4.4]{NaimpallyWarrack}, respectively (cf. \cite[Remark 4.8]{NaimpallyWarrack} and \cite[Definition 4.10]{NaimpallyWarrack}).
\end{proof}
\begin{remark}
\label{remarkstonecech}
The Samuel--Smirnov compactification may be compared with the Stone--\v Cech compactification, which is usually larger. The difference is that instead of Property \ref{propertysamuelsmirnov}, the Stone--\v Cech compactification has the property that for all $A,B\subset Z$, $\cl A\cap \cl B \neq \emptyset$ if and only if $\cl A\cap \cl B \cap Z \neq \emptyset$. Moreover, in the remarks following Property \ref{propertysamuelsmirnov}, ``uniformly continuous'' should be replaced with just ``continuous''.

We remark that if $\delta_G < \infty$ (i.e. if $G$ is of divergence type rather than of generalized divergence type), then the proof below works equally well if the Samuel--Smirnov compactification is replaced by the Stone--\v Cech compactification. This is not the case for the general proof; cf. Remark \ref{remarksamuelsmirnov}.
\end{remark}

To prove Theorem \ref{theoremahlforsthurstongeneral}, we will consider the Samuel--Smirnov compactification of the complete metric space $(\bord X,\wbar\Dist)$ (cf. Proposition \ref{propositionwbarDist}), which we will denote by $\what X$. For convenience of notation we will assume that $\bord X$ is a subset of $\what X$ and that $\iota:\bord X\to\what X$ is the inclusion map. As a point of terminology we will call points in $\bord X$ ``standard points'' and points in $\what X\butnot\bord X$ ``nonstandard points''.

\begin{remark}
Since $\wbar\Dist \asymp_\times \wbar\Dist_x$ for all $x\in X$, the Samuel--Smirnov compactification $\what X$ is independent of the basepoint $\zero$.
\end{remark}

At this point we can give a basic outline of the proof of Theorem \ref{theoremahlforsthurstongeneral}: First we will construct a measure $\what\mu$ on $\what X$ which satisfies the transformation equation \eqref{quasiconformal}. We will call such a measure $\what\mu$ a quasiconformal measure, although it is not \emph{a priori} a quasiconformal measure in the sense of Definition \ref{definitionquasiconformal}, as it is not necessarily supported on the set of standard points. Then we will use Thurston's proof of the Hopf--Tsuji--Sullivan theorem \cite[Theorem 4 of Section VII]{Ahlfors} (see also \cite[Theorem 2.4.6]{Nicholls}) to show that $\what\mu$ is supported on the nonstandard analogue of radial limit set. Finally, we will show that the nonstandard analogue of the radial limit set is actually a subset of $\bord X$, i.e. we will show that radial limit points are automatically standard. This demonstrates that $\what\mu$ is a measure on $\bord X$, and is therefore a \emph{bona fide} quasiconformal measure.

We now begin the preliminaries to the proof of Theorem \ref{theoremahlforsthurstongeneral}. As always $(X,\zero,b)$ denotes a Gromov triple. Let $\what X$ be the Samuel--Smirnov compactification of $\bord X$.

\begin{remark}
Throughout this chapter, $\cl S$ denotes the closure of a set $S$ taken with respect to $\what X$, not $\bord X$.
\end{remark}

\section{Extending the geometric functions to $\what X$}
We begin by extending the geometric functions $\dist(\cdot,\cdot)$, $\lb\cdot|\cdot\rb$, and $\busemann(\cdot,\cdot)$ to the Samuel--Smirnov compactification $\what X$. Extending $\dist(\cdot,\cdot)$ is the easiest:

\begin{observation}
If $x\in X$ is fixed, then the function $f_x:\bord X\to [0,1]$ defined by $f_x(y) = b^{-\dist(x,y)}$ is uniformly continuous by Remark \ref{remarkDistunifcont}. Thus by Proposition \ref{propositionsamuelsmirnov}, there exists a unique continuous extension $\fhat_x:\what X\to [0,1]$. We write 
\[
\dhat(x,\what y) = -\log_b\fhat_x(\what y).
\]
We define the \emph{extended boundary} of $X$ to be the set
\[
\what{\del X} := \{\xihat\in\what X:\dhat(\zero,\xihat) = \infty\}.
\]
Note that $\dhat(x,y) = \dist(x,y)$ if $x,y\in X$, and $\what{\del X} \cap \bord X = \del X$.
\end{observation}

\begin{warning*}
It is possible that $\what{\del X}\neq\cl{\del X}$.
\end{warning*}

On the other hand, extending the Gromov product to $\what X$ presents some difficulty, since the Gromov product is not necessarily continuous (cf. Example \ref{examplegromovnotcont}). Our solution is as follows: Fix $x\in X$ and $y\in\bord X$. Then by Remark \ref{remarkDistunifcont}, the map $\bord X\ni z\mapsto \Dist_x(y,z)$ is uniformly continuous, so by Proposition \ref{propositionsamuelsmirnov} it extends to a continuous map $\what X\ni \what z\mapsto \what\Dist_x(y,\what z)$. We define the \emph{Gromov product in $\what X$} via the formula
\[
\lb y|\what z\rb_x = -\log_b \what\Dist_x(y,\what z).
\]
Note that if $\what z\in \bord X$, then this notation conflicts with the previous definition of the Gromov product, but by Proposition \ref{propositionDist} the harm is only an additive asymptotic. We will ignore this issue in what follows.

\begin{observation}
Using (j) of Proposition \ref{propositionbasicidentities} we may define for each $x,y\in X$ the Busemann function
\[
\what\busemann_{\what z}(x,y) = \lb x|\what z\rb_y - \lb y|\what z\rb_x.
\]
Again, if $\what z\in\bord X$, then this definition conflicts with the previous one, but again the harm is only an additive asymptotic.
\end{observation}

\begin{remark}
We note that an appropriate analogue of Proposition \ref{propositionbasicidentities} (cf. also Corollary \ref{corollaryboundaryasymptotic}) holds on $\what X$. Specifically, each formula of Proposition \ref{propositionbasicidentities} holds with an additive asymptotic, as long as all expressions are defined. Note in particular that we have not defined the value of expressions which contain more than one nonstandard point. Such a definition would present additional difficulties (namely, noncommutativity of limits) which we choose to avoid.
\end{remark}

%
%
%
%

We are now ready to define the nonstandard analogue of the radial limit set:

\begin{definition}[cf. Definitions \ref{definitionshadow} and \ref{definitionradialconvergence}]
Given $x\in X$ and $\sigma > 0$, let
\[
\what\Shad(x,\sigma) = \{\xihat\in\what X : \lb\zero|\xihat\rb_x\leq\sigma\},
\]
so that $\what\Shad(x,\sigma)\cap\bord X = \Shad(x,\sigma)$. A sequence $(x_n)_1^\infty$ in $X$ will be said to converge to a point $\xihat\in \what{\del X}$ \emph{$\sigma$-radially} if $\dox{x_n}\to\infty$ and if $\xihat\in \what\Shad(x_n,\sigma)$ for all $n\in\N$. Note that in the definition of $\sigma$-radial convergence, we do not require that $x_n\to\xihat$ in the topology on $\what X$, although this can be seen from the proof of Lemma \ref{lemmaradialstandard} below.
\end{definition}

We conclude this section with the following lemma:

\begin{lemma}[Every radial limit point is a standard point]
\label{lemmaradialstandard}
Suppose that a sequence $(x_n)_1^\infty$ converges to a point $\xihat\in\what{\del X}$ $\sigma$-radially for some $\sigma > 0$. Then $\xihat\in\del X$.
\end{lemma}
\begin{proof}
We observe first that
\[
\lb x_n|\xihat\rb_\zero \asymp_\plus \dox{x_n} - \lb \zero|\xihat\rb_{x_n} \asymp_{\plus,\sigma} \dox{x_n} \tendsto n \dhat(\zero,\xihat) = \infty.
\]
Together with Gromov's inequality $\lb x_n|x_m\rb_\zero \gtrsim_\plus \min(\lb x_n|\xihat\rb_\zero,\lb x_m|\xihat\rb_\zero)$, this implies that $(x_n)_1^\infty$ is a Gromov sequence.

By the definition of the Gromov boundary, it follows that there exists a (standard) point $\eta\in\del X$ such that the sequence $(x_n)_1^\infty$ converges to $\eta$. Gromov's inequality now implies that $\lb \eta|\xihat\rb_\zero = \infty$. We claim now that $\xihat = \eta$, so that $\xihat$ is standard. By contradiction, suppose $\xihat\neq\eta$. Since $\what X$ is a Hausdorff space, it follows that there exist disjoint open sets $U,V\subset\what X$ containing $\xihat$ and $\eta$, respectively. Since $V$ contains a neighborhood of $\eta$, the function $f_{\zero,\eta}(z) = \lb \eta|z\rb_\zero$ is bounded from above on $\bord X\butnot V$. By continuity, $\fhat_{\zero,\eta}$ is bounded from above on $\cl{\bord X\butnot V}$. In particular, $\xihat\notin \cl{\bord X\butnot V}$. On the other hand $\xihat\notin \cl V$, since $\xihat$ is in the open set $U$ which is disjoint from $V$. It follows that $\xihat\notin \cl{\bord X} = \what X$, a contradiction.
\end{proof}
\begin{remark}
\label{remarkradialstandard}
In fact, the above proof shows that if
\begin{equation}
\label{strongconvergence}
\lb x_n|\xihat\rb_\zero \to \infty
\end{equation}
for some sequence $(x_n)_1^\infty$ in $X$ and some $\xihat\in\what{\del X}$, then $\xihat\in\del X$. However, there may be a sequence $(x_n)_1^\infty$ such that $x_n\to \xihat$ in the topology on $\what X$ but for which \eqref{strongconvergence} does not hold. In this case, we could have $\xihat\notin\del X$.
\end{remark}

\section{Quasiconformal measures on $\what X$}
We define the notion of a quasiconformal measure on $\what X$ as follows:

\begin{definition}[cf. Definition \ref{definitionquasiconformal}, Proposition \ref{propositionderivativebusemann}]
\label{definitionquasiconformalnonstandard}
For each $s\geq 0$, a Radon probability measure $\what\mu$ on $\what{\del X}$ is called \emph{$s$-quasiconformal} if
\[
\what\mu(\what g(A)) \asymp_\times \int_A b^{s\what\busemann_{\what\eta}(\zero,g^{-1}(\zero))}\;\dee\what\mu(\what\eta).
\]
for every $g\in G$ and for every Borel set $A\subset\what{\del X}$. Here $\what g$ denotes the unique continuous extension of $g$ to $\what X$ (cf. Proposition \ref{propositionsamuelsmirnov}).
\end{definition}

\begin{remark}
\label{remarkradon}
Note that we have added here the assumption that the measure $\what\mu$ is \emph{Radon}. Since the phrase ``Radon measure'' seems to have no generally accepted meaning in the literature, we should make clear that for us a (finite, nonnegative, Borel) measure $\mu$ on a compact Hausdorff space $Z$ is Radon if the following two conditions hold (cf. \cite[\67]{Folland}):
\begin{align*}
\mu(A) &= \inf\{\mu(U) : U\supset A, \;\; U\text{ open}\} \all A\subset Z \text{ Borel}\\
\mu(U) &= \sup\{\mu(K) : K\subset U, \;\; K\text{ compact}\} \all U\subset Z \text{ open.}
\end{align*}
The assumption of Radonness was not needed in Definition \ref{definitionquasiconformal}, since every measure on a compact metric space is Radon \cite[Theorem 7.8]{Folland}. However, the assumption is important in the present proof, since $\what X$ is not necessarily metrizable, and so it may have non-Radon measures.

On the other hand, the Radon condition itself is of no importance to us, except for the following facts:
\begin{itemize}
\item[(i)] The image of a Radon measure under a homeomorphism is Radon.
\item[(ii)] Every measure absolutely continuous to a Radon measure is Radon.
\item[(iii)] The sum of two Radon measures is Radon.
\item[(iv)] (Riesz representation theorem, \cite[Theorem 7.2]{Folland}) Let $Z$ be a compact Hausdorff space. For each measure $\mu$ on $Z$, let $I_\mu$ denote the nonnegative linear function
\[
I_\mu[f] := \int f\;\dee\mu.
\]
Then for every nonnegative linear functional $I:C(Z)\to \R$, there exists a unique Radon measure $\mu$ on $Z$ such that $I_\mu = I$. (If $\mu_1$ and $\mu_2$ are not both Radon, it is possible that $I_{\mu_1} = I_{\mu_2}$ while $\mu_1\neq\mu_2$.)
\end{itemize}
\end{remark}

We now state two lemmas which are nonstandard analogues of lemmas proven in Chapter \ref{sectionconformal}. We omit the parts of the proofs which are the same as in the standard case, reminding the reader that the important point is that no function is ever used which takes two nonstandard points as inputs. We begin by proving an analogue of Sullivan's shadow lemma:

\begin{lemma}[Sullivan's Shadow Lemma on $\what X$; cf. Lemma \ref{lemmasullivanshadow}]
\label{lemmasullivanshadownonstandard}
Fix $s\geq 0$, and let $\what\mu\in\MM(\what{\del X})$ be an $s$-quasiconformal measure which is not a pointmass supported on a standard point. Then for all $\sigma > 0$ sufficiently large and for all $g\in G$, we have
\[
\what\mu(\what\Shad(g(\zero),\sigma)) \asymp_\times b^{-s\dogo g}.
\]
\end{lemma}
\begin{proof}
Obvious modifications\Footnote{We remark that the expression $\overline g'(\xi)$ occuring in the proof of Lemma \ref{lemmasullivanshadow} should be replaced by $b^{-\what\busemann_{\what\xi}(\zero,g^{-1}(\zero))}$ as per Proposition \ref{propositionderivativebusemann}; of course, the expression $\overline g'(\what\xi)$ makes no sense, since $\what X$ is not a metric space. \label{footnotenonstandardderivative}} to the proof of Lemma \ref{lemmasullivanshadow} yield
\[
\what\mu(\what\Shad(g(\zero),\sigma)) \asymp_{\times,\mu,\sigma} b^{-s\dogo g}\what\mu\big(\what\Shad_{g^{-1}(\zero)}(\zero,\sigma)\big).
\]
So to complete the proof, we need to show that
\[
\what\mu\big(\what\Shad_{g^{-1}(\zero)}(\zero,\sigma)\big) \asymp_{\times,\mu,\sigma} 1,
\]
assuming $\sigma > 0$ is sufficiently large (depending on $\what\mu$). By contradiction, suppose that for each $n\in\N$ there exists $g_n\in G$ such that
\[
\what\mu\big(\what\Shad_{g_n^{-1}(\zero)}(\zero,n)\big) \leq \frac{1}{2^n}\cdot
\]
Then for $\what\mu$-a.e. $\xihat\in\what X$,
\[
\xihat\in \what\Shad_{g_n^{-1}(\zero)}(\zero,n) \text{ for all but finitely many $n$},
\]
which implies
\[
\lb g_n^{-1}(\zero)|\xihat\rb_\zero \gtrsim_\plus n \tendsto n \infty.
\]
By Remark \ref{remarkradialstandard}, it follows that $\xihat\in\del X$ and $g_n^{-1}(\zero) \to \xihat$. This implies that $\what\mu$ is a pointmass supported on the standard point $\lim_{n\to\infty} g_n^{-1}(\zero)$, contradicting our hypothesis.
\end{proof}

\begin{lemma}[cf. Theorem \ref{theorempattersonsullivangeneral}]
\label{lemmapattersonsullivannonstandard}
Assume that $\w\delta = \w\delta_G < \infty$. Then there exists a $\w\delta$-quasiconformal measure supported on $\what{\del X}$.
\end{lemma}
\begin{proof}
Let the measures $\mu_s$ be as in \eqref{musdef}. The compactness of $\what X$ replaces the assumption that $G$ is of compact type which occurs in Theorem \ref{theorempattersonsullivangeneral}, so there exists a sequence $s_n\searrow\w\delta$ such that $\mu_n := \mu_{s_n}\to \what\mu$ for some Radon measure $\what\mu\in \MM(\what X)$. Claim \ref{claimboundarysupport} shows that $\what\mu$ is supported on $\what{\del X}$.

To complete the proof, we must show that $\what\mu$ is $\w\delta$-quasiconformal. Fix $g\in G$ and a continuous function $f:\what X\to (0,\infty)$. The final assertion of Proposition \ref{propositionsamuelsmirnov} guarantees that $\log(f) \given \bord X$ is uniformly continuous, so the proof of Claim \ref{claimconformalf} shows that \eqref{conformalf} holds.

The equation \eqref{riesz} deserves some comment; it depends on the uniqueness assertion of the Riesz representation theorem, which, now that we are no longer in a metric space, holds only for Radon measures. But by Remark \ref{remarkradon}, all measures involved in \eqref{riesz} are Radon, so \eqref{riesz} still holds.
\end{proof}

\begin{remark}
\label{remarksamuelsmirnov}
In this lemma we used the final assertion of Proposition \ref{propositionsamuelsmirnov} in a nontrivial way. The proof of this lemma would not work for the Stone--\v Cech compactification, except in the case $\delta < \infty$, in which case the uniform continuity of $f$ is not necessary in the proof of Theorem \ref{theorempattersonsullivangeneral}.
\end{remark}

\begin{lemma}[Intersecting Shadows Lemma on $\what X$; cf. Lemma \ref{lemmatau}]
\label{lemmataunonstandard}
For each $\sigma > 0$, there exists $\tau = \tau_\sigma > 0$ such that for all $x,y,z\in X$ satisfying $\dist(z,y)\geq \dist(z,x)$ and $\what\Shad_z(x,\sigma)\cap \what\Shad_z(y,\sigma)\neq\emptyset$, we have
\begin{equation}
\label{shadcontainmentnonstandard}
\what\Shad_z(y,\sigma) \subset \what\Shad_z(x,\tau)
\end{equation}
and
\begin{equation}
\label{distbusemannnonstandard}
\dist(x,y) \asymp_{\plus,\sigma} \dist(z,y) - \dist(z,x).
\end{equation}
\end{lemma}
\begin{proof}
The proof of Lemma \ref{lemmatau} goes through with no modifications needed.
\end{proof}

\section{The main argument}

\begin{proposition}[Generalization/nonstandard version of Theorem \ref{theoremahlforsthurston}(A) \implies (B)]
\label{propositionahlforsthurston}
Let $\what\mu$ be a $\w\delta$-quasiconformal measure on $\what{\del X}$ which is not a pointmass supported on a standard point. If $G$ is of generalized divergence type, then $\what\mu(\Lr(G)) > 0$.
\end{proposition}
\begin{proof}
Fix $\sigma > 0$ large enough so that Sullivan's Shadow Lemma \ref{lemmasullivanshadownonstandard} holds. Let $\rho > 0$ be large enough so that there exists a maximal $\rho$-separated set $S_\rho\subset G(\zero)$ which has finite intersection with bounded sets (cf. Proposition \ref{propositionbasicmodified}(iii)). Let $(x_n)_1^\infty$ be an indexing of $S_\rho$. By Lemma \ref{lemmaradialstandard}, we have 
\[
\bigcap_{N\in\N}\bigcup_{n\geq N}\what\Shad(x_n,\sigma + \rho) \subset \Lr(G).
\] 
By contradiction suppose that $\what\mu(\Lr(G)) = 0$. Fix $\epsilon > 0$ small to be determined. Then there exists $N\in\N$ such that
\[
\what\mu\left(\bigcup_{n\geq N}\what\Shad(x_n,\sigma + \rho)\right) \leq \epsilon.
\]
Let $R = \rho + \max_{n < N}\dox{x_n}$. Then
\[
\what\mu\left(\bigcup_{\substack{g\in G \\ \dogo g >
R}}\what\Shad(g(\zero),\sigma)\right) \leq \epsilon. 
\]
We shall prove the following.

\begin{observation}
\label{observationhypothesis}
If $A\subset G(\zero)$ is any subcollection satisfying
\begin{itemize}
\item [(I)] $\dox x > R$ for all $x\in A$, and
\item [(II)] $(\what\Shad(x,\sigma))_{x\in A}$ are disjoint,
\end{itemize}
then
\begin{equation}
\label{hypothesis}
\sum_{x\in A}b^{-s\dox x} \lesssim_\times \epsilon.
\end{equation}
\end{observation}
\begin{subproof}
The disjointness condition guarantees that
\[
\sum_{x\in A}\what\mu(\what\Shad(x,\sigma)) \leq \what\mu\left(\bigcup_{\substack{g\in G \\ \dogo g > R}}\what\Shad(g_n(\zero),\sigma)\right) \leq \epsilon.
\]
Combining with Sullivan's Shadow Lemma \ref{lemmasullivanshadownonstandard} yields \eqref{hypothesis}.
\end{subproof}

Now choose $R' > R$ and $\sigma' > \sigma$ large to be determined. Let $S_{R'}$ be a maximal $R'$-separated subset of $G(\zero)$. For convenience we assume $\zero\in S_{R'}$. By Proposition \ref{propositionbasicmodified}(iv), if $R'$ is sufficiently large then $\Sigma_{\w\delta}(S_{R'}) = \infty$ if and only if $\w\delta$ is of generalized divergence type. So to complete the proof, it suffices to show that
\[
\Sigma_{\w\delta}(S_{R'}) < \infty.
\]
\begin{notation}
\label{notationpartialorder}
Let $(x_i)_1^\infty$ be an indexing of $S_{R'}$ such that $i < j$ implies $\dox{x_i} \leq \dox{x_j}$. For $x_i,x_j\in S_{R'}$ distinct, we write $x_i < x_j$ if
\begin{itemize}
\item [(I)] $i < j$ and
\item [(II)] $\what\Shad(x_i,\sigma')\cap\what\Shad(x_j,\sigma')\neq\emptyset$.
\end{itemize}
(This is just a notation, it does not mean that $<$ is a partial order on $S_{R'}$.)
\end{notation}
\begin{lemma}
\label{lemmaahlforsdisjoint}
If $R'$ and $\sigma'$ are sufficiently large (with $\sigma'$ chosen first), then
\[
x < y \;\;\Rightarrow\;\; \what\Shad_x(y,\sigma) \subset \what\Shad(y,\sigma').
\]
\end{lemma}
\begin{subproof}
Suppose $x < y$; then $\what\Shad(x,\sigma')\cap\what\Shad(y,\sigma')\neq\emptyset$. By the Intersecting Shadows Lemma \ref{lemmataunonstandard}, we have $\dist(x,y) \asymp_{\plus,\sigma'} \dox y - \dox x$. On the other hand, since $S_{R'}$ is $R'$-separated we have $\dist(x,y) \geq R'$. Thus
\[
\lb \zero|x\rb_y \gtrsim_{\plus,\sigma'} R'.
\]
Now for any $\xihat\in\what X$, we have
\[
\lb x|\xihat\rb_y \gtrsim_\plus \min(\lb \zero|\xihat\rb_y,\lb \zero|x\rb_y).
\]
Thus if $\xihat\in\what\Shad_x(y,\sigma)$, then
\[
\sigma \gtrsim_\plus \lb \zero|\xihat\rb_y \text{ or } \sigma \gtrsim_{\plus,\sigma'} R'.
\]
Let $\sigma'$ be $\sigma$ plus the implied constant of the first asymptotic, and then let $R'$ be $\sigma + 1$ plus the implied constant of the second asymptotic. Then the second asymptotic is automatically impossible, so
\[
\lb \zero|\xihat\rb_y \leq \sigma',
\]
i.e. $\xihat\in\what\Shad(y,\sigma')$.
\end{subproof}
If $x\in S_{R'}$ is fixed, let us call $y\in S_{R'}$ an \emph{immediate successor} of $x$ if $x < y$ but there is no $z$ such that $x < z < y$. We denote by $S_{R'}(x)$ the collection of all immediate successors of $x$.

\begin{lemma}
\label{lemmasuccessors}
For each $z\in S_{R'}$, we have
\begin{equation}
\label{successors}
\sum_{y\in S_{R'}(z)}b^{-s\dox y} \lesssim_\times \epsilon b^{-s\dox z}.
\end{equation}
\end{lemma}
\begin{subproof}
We claim first that the collection $(\what\Shad(y,\sigma'))_{y\in S_{R'}(z)}$ consists of mutually disjoint sets. Indeed, if $\what\Shad(y_1,\sigma')\cap \what\Shad(y_2,\sigma')\neq\emptyset$ for some distinct $y_1,y_2\in S_{R'}(z)$, then we would have either $z < y_1 < y_2$ or $z < y_2 < y_1$, contradicting the definition of immediate successor. Combining with Lemma \ref{lemmaahlforsdisjoint}, we see that the collection $(\what\Shad_z(y,\sigma))_{y\in S_{R'}(z)}$ also consists of mutually disjoint sets.

Fix $g\in G$ such that $g(\zero) = z$. We claim that the collection
\[
A = g^{-1}(S_{R'}(z))
\]
satisfies the hypotheses of Observation \ref{observationhypothesis}. Indeed, as $\zero\notin A$ (since $z\notin S_{R'}(z)$) and as $g$ is an isometry of $X$, (I) follows from the fact that $S_{R'}$ is $R'$-separated and $R' > R$. Since $\what\Shad(g^{-1}(y),\sigma) = g^{-1}(\what\Shad_z(y,\sigma))$ for all $y\in S_{R'}(z)$, the collection $(\what\Shad(x,\sigma))_{x\in A}$ consists of mutually disjoint sets, meaning that (II) holds. Thus, by Observation \ref{observationhypothesis}, we have
\[
\sum_{x\in A}b^{-s\dox x} \lesssim_\times \epsilon,
\]
or, since $g$ is an isometry of $X$ and $z=g(\zero)$,
\[
\sum_{y\in S_{R'}(z)}b^{-s\dist(z,y)} \lesssim_\times \epsilon.
\]
Inserting \eqref{distbusemannnonstandard} into the last inequality yields \eqref{successors}.
\end{subproof}

Using Lemma \ref{lemmasuccessors}, we complete the proof. Define the sequence $(S_n)_{n = 0}^\infty$ inductively as follows:
\begin{align*}
S_0 &= \{\zero\}, \\
S_{n + 1} &= \bigcup_{x\in S_n}S_{R'}(x). 
\end{align*}
Clearly, all immediate successors of all points of $\bigcup_{n\geq 0}S_n$ belong to $\bigcup_{n\geq 0}S_n$. We claim that
\[
S_{R'} = \bigcup_{n\geq 0}S_n.
\]
Indeed, let $(x_i)_1^\infty$ be the indexing of $S_{R'}$ considered in Notation \ref{notationpartialorder}, and by induction suppose that $x_i\in \bigcup_1^\infty S_n$ for all $i < j$. If $j = 0$, then $x_j = \zero\in S_0$. Otherwise, let $i < j$ be maximal satisfying $x_i < x_j$. Then $x_j$ is an immediate successor of $x_i\in \bigcup_1^\infty S_n$, so $x_j\in \bigcup_1^\infty S_n$.

Summing \eqref{successors} over all $x\in S_n$, we have
\[
\sum_{y\in S_{n + 1}}b^{-s\dox y} \lesssim_\times \epsilon
\sum_{x\in S_n}b^{-s\dox x}. 
\]
Set $\epsilon$ equal to $1/2$ divided by the implied constant, so that
\[
\sum_{y\in S_{n + 1}}b^{-s\dox y} \leq \frac12 \sum_{x\in S_n}b^{-s\dox x}.
\]
Applying the Ratio Test, we see that the series $\Sigma_{\w\delta}(S_{R'})$ converges, contradicting that $G$ was of generalized divergence type. 
\end{proof}

\begin{corollary}
\label{corollaryahlforsthurston}
Let $\what\mu$ be a $\w\delta$-quasiconformal measure on $\what{\del X}$. If $G$ is of generalized divergence type, then $\what\mu(\Lr(G)) = 1$.
\end{corollary}
\begin{proof}
By contradiction suppose not. Then $\what\nu := \what\mu\given \what{\del X}\butnot\Lr(G)$ is a $\w\delta$-quasiconformal measure on $\what{\del X}$ which gives zero measure to $\Lr(G)$, contradicting Proposition \ref{propositionahlforsthurston}.
\end{proof}

\section{End of the argument}
We now complete the proof of Theorem \ref{theoremahlforsthurstongeneral}:

\begin{proof}[Proof of Theorem \ref{theoremahlforsthurstongeneral}]
Let $\what\mu$ be the $\w\delta$-quasiconformal measure supported on $\what{\del X}$ guaranteed by Lemma \ref{lemmapattersonsullivannonstandard}. By Corollary \ref{corollaryquasiconformalpointmass}, $\what\mu$ is not a pointmass supported on a standard point. By Corollary \ref{corollaryahlforsthurston}, $\what\mu$ is supported on $\Lr(G) \subset \del X$. This completes the proof of the existence assertion.

Suppose that $\mu_1,\mu_2$ are two $\w\delta$-quasiconformal measures on $\del X$. By Corollary \ref{corollaryahlforsthurston}, $\mu_1$ and $\mu_2$ are both supported on $\Lr(G)$.

Suppose first that $\mu_1,\mu_2$ are supported on $\Lrsigma$ for some $\sigma > 0$. Fix an open set $U\subset\del X$. By the Vitali covering theorem, there exists a collection of disjoint shadows $(\Shad(g(\zero),\sigma))_{g\in A}$ contained in $U$ such that 
\[
\mu_1(U\butnot \bigcup_{g\in A}\Shad(g(\zero),\sigma)) = 0.
\]
Then we have 
\begin{align*}
\mu_1(U) = \sum_{g\in A} \mu_1(\Shad(g(\zero),\sigma)) 
&\asymp_{\times,\mu_1} \sum_{g\in A} b^{-s\dogo g}\\
&\asymp_{\times,\mu_2} \sum_{g\in A} \mu_2(\Shad(g(\zero),\sigma))\\
&\leq \mu_2(U).
\end{align*}
A similar argument shows that $\mu_2(U) \lesssim_\times \mu_1(U)$. Since $U$ was arbitrary, a standard approximation argument shows that $\mu_1\asymp_\times \mu_2$. It follows that any individual measure $\mu$ supported on $\Lrsigma$ is ergodic, because if $A$ is an invariant set with $0 < \mu(A) < 1$ then $\frac{1}{\mu(A)}\mu\given A$ and $\frac{1}{1 - \mu(A)}\mu\given (\Lr\butnot A)$ are two measures which are not asymptotic, a contradiction.

In the general case, define the function $f:\Lr\to \Rplus$ by
\[
f(\xi) = \sup\{\sigma > 0 : \exists g\in G \;\; g(\xi)\in\Lrsigma\}.
\]
By Proposition \ref{propositionnearinvariance}, $f(\xi) < \infty$ for all $\xi\in\Lr$. On the other hand, $f$ is $G$-invariant. Now let $\mu$ be a $\w\delta$-quasiconformal measure on $\Lr$. Then for each $\sigma_0 < \infty$ the measure $\mu\given f^{-1}([0,\sigma_0])$ is supported on $\Lambda_{\mathrm r, \sigma_0}$, and is therefore ergodic; thus $f$ is constant $\mu\given f^{-1}([0,\sigma_0])$-a.s. It is clear that this constant value is independent of $\sigma_0$ for large enough $\sigma_0$, so $f$ is constant $\mu$-a.s. Thus there exists $\sigma > 0$ such that $\mu$ is supported on $\Lrsigma$, and we can reduce to the previous case.
\end{proof}

\section{Necessity of the generalized divergence type assumption}
The proof of Theorem \ref{theoremahlforsthurstongeneral} makes crucial use of the generalized divergence type assumption, just as the proof of Theorem \ref{theorempattersonsullivangeneral} made crucial use of the compact type assumption. What happens if neither of these assumptions holds? Then there may not be a $\w\delta$-quasiconformal measures supported on the limit set, as we now show:

\begin{proposition}
\label{propositioncounterexample}
There exists a strongly discrete group of general type $G\leq\Isom(\H^\infty)$ satisfying $\delta < \infty$, such that there does not exist any quasiconformal measure supported on $\Lambda$.
\end{proposition}
\begin{proof}
The idea is to first construct such a group in an $\R$-tree, and then to use a BIM embedding (Theorem \ref{theoremBIM}) to get an example in $\H^\infty$. Fix a sequence of numbers $(a_k)_1^\infty$. For each $k$ let $\Gamma_k = \{e,\gamma_k\} \equiv \Z_2$, and let $\|\cdot\|:\Gamma_k\to\R$ be defined by $\|\gamma_k\| = a_k$, $\|e\| = 0$. Clearly, the function $\|\cdot\|$ is tree-geometric in the sense of Definition \ref{definitiontreegeometric}, so by Theorem \ref{theoremschottkyproductRtree}, the function $\|\cdot\|:\Gamma\to\Rplus$ defined by \eqref{schottkyproductRtree} is tree-geometric, where $\Gamma = \ast_{k\in\N} \Gamma_k$. So there exist an $\R$-tree $X$ and a homomorphism $\phi:\Gamma\to\Isom(X)$ such that $\dogo{\phi(\gamma)} = \|\gamma\| \all\gamma\in\Gamma$. Let $G = \phi(\Gamma)$.
\begin{claim}
If the sequences $(a_k)_1^\infty$ is chosen appropriately, then $G$ is of convergence type.
\end{claim}
\begin{subproof}
For $s\geq 0$ we have
\begin{align*}
\Sigma_s(G) - 1 &= \sum_{g\in G\butnot\{\id\}} e^{-s\dogo g}\\
&= \sum_{(k_1,\gamma_1)\cdots(k_n,\gamma_n)\in (\Gamma_E)^*\butnot\{\smallemptyset\}} \exp\big(-s\big[a_{k_1} +\ldots + a_{k_n}\big]\big)\\
&= \sum_{n\in\N} \sum_{k_1\neq k_2\neq \cdots \neq k_n} \sum_{\gamma_1\in \Gamma_{k_1}\butnot\{e\}}\cdots\sum_{\gamma_n\in \Gamma_{k_n}\butnot\{e\}} \exp\big(-s\big[a_{k_1} + \ldots + a_{k_n}\big]\big)\\
&= \sum_{n\in\N} \sum_{k_1\neq k_2\neq \cdots \neq k_n} \prod_{i = 1}^n e^{-s a_{k_i}}\\
\Sigma_s(G) &\leq 1 + \sum_{n\in\N} \left( \sum_{k\in\N} e^{-s a_k}\right)^n\\
\Sigma_s(G) &\geq 1 + \sum_{k\in\N} e^{-s a_k}.
\end{align*}
Thus, letting
\[
P_s = \sum_{k\in\N} e^{-s a_k},
\]
we have
\begin{equation}
\label{esak}
\begin{cases}
\Sigma_s(G) < \infty & \text{ if } P_s < 1 \\
\Sigma_s(G) = \infty & \text{ if } P_s = \infty
\end{cases}.
\end{equation}
Now clearly, there exists a sequence $(a_k)_1^\infty$ such that $P_{1/2} < 1$ but $P_s = \infty$ for all $s < 1/2$; for example, take $a_k = \log(k) + 2\log\log(k) + C$ for sufficiently large $C$.
\end{subproof}
\begin{claim}
$\Lambda(G) = \Lr(G)$.
\end{claim}
\begin{subproof}
For all $\xi\in\LambdaG$, the path traced by the geodesic ray $\geo\zero\xi$ in $X/G$ is the concatenation of infinitely many paths of the form $\geo\zero{g(\zero)}$, where $g\in \bigcup_{n\in\N}\phi(\Gamma_n)$. Each such path crosses $\zero$, so the path traced by the geodesic ray $\geo\zero\xi$ in $X/G$ crosses $\zero$ infinitely often. Equivalently, the geodesic ray $\geo\zero\xi$ crosses $G(\zero)$ infinitely often. By Proposition \ref{propositionradialconvergence}, this implies that $\xi\in\Lr(G)$.
\end{subproof}

Now let $\w G$ be the image of $G$ under a BIM representation (cf. Theorem \ref{theoremBIM}). By Remark \ref{remarkBIM}, $\w G$ is of convergence type and $\Lambda(\w G) = \Lr(\w G)$. The proof is completed by the following lemma:

\begin{lemma}
\label{lemmaconvergencetype}
If the group $G$ is of generalized convergence type and $\mu$ is a $\w\delta$-quasiconformal measure, then $\mu(\Lr) = 0$.
\end{lemma}
\begin{subproof}
Fix $\sigma > 0$ large enough so that Sullivan's Shadow Lemma \ref{lemmasullivanshadow} holds. Fix $\rho > 0$ and a maximal $\rho$-separated set $S_\rho\subset G(\zero)$ such that $\Sigma_{\w\delta}(S_\rho) < \infty$. Then
\[
\sum_{x\in S_\rho} \mu(\Shad(x,\rho + \sigma)) \asymp_{\times,\rho,\sigma} \sum_{x\in S_\rho} b^{-\w\delta\dox x} < \infty.
\]
On the other hand, $\Lrsigma\subset\limsup_{x\in S_\rho} \Shad(x,\rho + \sigma)$. So by the Borel--Cantelli lemma, $\mu(\Lrsigma) = 0$. Since $\sigma$ was arbitrary, $\mu(\Lr) = 0$.
\end{subproof}
\end{proof}

Combining Theorem \ref{theoremahlforsthurstongeneral} and Lemma \ref{lemmaconvergencetype} yields the following:

\begin{proposition}
\label{propositiondivtypedichotomy}
Let $G\leq\Isom(X)$ be a nonelementary group with $\w\delta < \infty$. Then the following are equivalent:
\begin{itemize}
\item[(A)] $G$ is of generalized divergence type.
\item[(B)] There exists a $\w\delta$-conformal measure $\mu$ on $\Lambda$ satisfying $\mu(\Lr) > 0$.
\item[(C)] Every $\w\delta$-conformal measure $\mu$ on $\Lambda$ satisfies $\mu(\Lr) = 1$.
\item[(D)] There exists a unique $\w\delta$-conformal measure $\mu$ on $\Lambda$, and it satisfies $\mu(\Lr) = 1$.
\end{itemize}
\end{proposition}

\section{Orbital counting functions of nonelementary groups}
Theorem \ref{theoremahlforsthurston} allows us to prove the following result which, on the face of it, does not involve quasiconformal measures at all:

\begin{corollary}
\label{corollaryupperorbitalbound}
Let $G\leq\Isom(X)$ be nonelementary and satisfy $\delta < \infty$. Then
\[
\NN_{X,G}(\rho) \lesssim_\times b^{\delta\rho} \all \rho\geq 0.
\]
\end{corollary}
\begin{proof}
If $G$ is of convergence type, then the bound is obvious, as
\[
b^{-\delta \rho} \NN_{X,G}(\rho) \leq \sum_{\substack{g\in G \\ \dogo g \leq \rho}} b^{-\delta\dogo g} \leq \Sigma_\delta(G) < \infty.
\]
On the other hand, if $G$ is of divergence type, then by Theorem \ref{theoremahlforsthurstongeneral}, there exists a $\delta$-conformal measure $\mu$ on $\Lambda$, which is not a pointmass by Corollary \ref{corollaryquasiconformalpointmass} and Proposition \ref{propositionnonelementaryequivalent}(C). Remark \ref{remarkupperorbitalbound} finishes the proof.
\end{proof}

We contrast this with a philosophically related result whose proof uses the Ahlfors regular measures constructed in the proof of Theorem \ref{theorembishopjonesregular}:

\begin{proposition}
Let $G\leq\Isom(X)$ be a nonelementary and strongly discrete. Then
\[
\frac{\log_b\NN_{X,G}(\rho)}{\rho} \tendsto{\rho\to\infty} \delta.
\]
\end{proposition}
\begin{proof}
Fix $s < \delta$. By Theorem \ref{theorembishopjonesregular}, there exist $\tau > 0$ and an Ahlfors $s$-regular measure $\mu$ supported on $\Lurtau$. Now fix $\rho > 0$. By definition, for all $\xi\in\Lurtau$ there exists $g\in G$ such that $\rho - \tau \leq \|g\| \leq \rho$ and $\lb \zero | \xi \rb_{g(\zero)} \leq \tau$. Equivalently,
\[
\Lurtau \subset \bigcup_{\substack{g\in G \\ \rho - \tau \leq \|g\| \leq \rho}} \Shad(g(\zero),\tau).
\]
Applying $\mu$ to both sides gives
\[
1 \leq \sum_{\substack{g\in G \\ \rho - \tau \leq \|g\| \leq \rho}} \mu\big(\Shad(g(\zero),\tau)\big).
\]
By the Diameter of Shadows Lemma, $\Diam(\Shad(g(\zero),\tau)) \lesssim_{\times,s} b^{-\|g\|}$ and thus since $\mu$ is Ahlfors $s$-regular,
\[
\mu\big(\Shad(g(\zero),\tau)\big) \lesssim_{\times,s} b^{-s\|g\|} \asymp_{\times,s} b^{-s\rho}.
\]
So
\[
1 \lesssim_{\times,s} b^{-s\rho}\#\{g\in G : \rho - \tau \leq \|g\| \leq \rho\}
\]
and thus
\[
\#\{g\in G : \|g\| \leq \rho\} \gtrsim_{\times,s} b^{s\rho}.
\]
Since $s < \delta$ was arbitrary, we get
\[
\liminf_{\rho\to\infty} \frac{\log_b\#\{g\in G : \|g\| \leq \rho\}}{\rho} \geq \delta.
\]
Combining with \eqref{poincarealternate} completes the proof.
\end{proof}

%

\chapter{Quasiconformal measures of geometrically finite groups}
\label{sectionGFmeasures}

In this chapter we investigate the $\delta$-quasiconformal measure or measures associated to a geometrically finite group. Note that since geometrically finite groups are of compact type (Theorem \ref{theoremGFcompact}), Theorem \ref{theorempattersonsullivangeneral} guarantees the existence of a $\delta$-quasiconformal measure $\mu$ on $\Lambda$. However, this measure is not necessarily unique (Corollary \ref{corollaryuniquenesscharacterization}); a sufficient condition for uniqueness is that $G$ is of divergence type (Theorem \ref{theoremahlforsthurstongeneral}). In Section \ref{subsectiondivergencetype}, we generalize a theorem of Dal'bo, Otal, and Peigne \cite[Th\'eor\`eme A]{DOP} which shows that ``most'' geometrically finite groups are of divergence type. In Sections \ref{subsectionglobalmeasure}-\ref{subsectionexactdimensional} we investigate the geometry of $\delta$-conformal measures; specifically, in Sections \ref{subsectionglobalmeasure}-\ref{subsectionglobalmeasureproof} we prove a generalization of the Global Measure Formula (Theorem \ref{theoremglobalmeasure}), in Sections \ref{subsectiondoubling} and \ref{subsectionexactdimensional} we investigate the questions of when the $\delta$-conformal measure of a geometrically finite group is doubling and exact dimensional, respectively.

\begin{standingassumptions}
\label{standingassumptionsGFmeasures}
In this chapter, we assume that
\begin{itemize}
\item[(I)] $X$ is regularly geodesic and strongly hyperbolic,
\item[(II)] $G\leq\Isom(X)$ is nonelementary and geometrically finite, and $\delta < \infty$.\Footnote{Note that by Corollary \ref{corollaryGF}(ii), we have $\delta < \infty$ if and only if $\delta_\bp < \infty$ for all $\bp\in P$.}
\end{itemize}
Moreover, we fix a complete set of inequivalent parabolic points $P\subset\Lbp$, and for each $\bp\in P$ we write $\delta_\bp = \delta(G_\bp)$, and let $S_\bp\subset\EE_\bp$ be a $\bp$-bounded set satisfying (A)-(C) of Lemma \ref{lemmaboundedparabolic}. Finally, we choose a number $t_0 > 0$ large enough so that if
\begin{align*}
H_\bp &= H_{\bp,t_0} = \{ x\in X : \busemann_\bp(\zero,x) > t_0 \}\\
\scrH &= \{ g(H_\bp) : \bp\in P, g\in G\},
\end{align*}
then the collection $\scrH$ is disjoint (cf. Proof of Theorem \ref{theoremGFcompact}\text{(B3) \implies (A)}).
\end{standingassumptions}


\section{Sufficient conditions for divergence type}
\label{subsectiondivergencetype}

In the Standard Case, all geometrically finite groups are of divergence type \cite[Proposition 2]{Sullivan_entropy}; however, once one moves to the more general setting of pinched Hadamard manifolds, one has examples of geometrically finite groups of convergence type \cite[Th\'eor\`eme C]{DOP}. On the other hand, Proposition \ref{propositiondivtypedichotomy} shows that for every $\delta$-conformal measure $\mu$, $G$ is of divergence type if and only if $\mu(\Lambda\butnot\Lr) = 0$. Now by Theorem \ref{theoremGFcompact}, $\Lambda\butnot\Lr = \Lbp = G(P)$, so the condition $\mu(\Lambda\butnot\Lr) = 0$ is equivalent to the condition $\mu(P) = 0$. To summarize:

\begin{observation}
\label{observationnonatomic}
The following are equivalent:
\begin{itemize}
\item[(A)] $G$ is of divergence type.
\item[(B)] There exists a $\delta$-conformal measure $\mu$ on $\Lambda$ satisfying $\mu(P) = 0$.
\item[(C)] Every $\delta$-conformal measure $\mu$ on $\Lambda$ satisfies $\mu(P) = 0$.
\item[(D)] There exists a unique $\delta$-conformal measure $\mu$ on $\Lambda$, and it satisfies $\mu(P) = 0$.
\end{itemize}
In particular, every convex-cobounded group is of divergence type.
\end{observation}

It is of interest to ask for sufficient conditions which are not phrased in terms of measures. We have the following:

\begin{theorem}[Cf. {\cite[Proposition 2]{Sullivan_entropy}}, {\cite[Th\'eor\`eme A]{DOP}}]
\label{theoremGFdivergencetype}
If $\delta > \delta_\bp$ for all $\bp\in P$, then $G$ is of divergence type.
\end{theorem}
\begin{proof}
We will demonstrate (B) of Observation \ref{observationnonatomic}. Let $\mu$ be the measure constructed in the proof of Theorem \ref{theorempattersonsullivangeneral}, fix $\bp\in P$, and we will show that $\mu(\bp) = 0$. In what follows, we use the same notation as in the proof of Theorem \ref{theorempattersonsullivangeneral}. Since $G$ is strongly discrete, we can let $\rho$ be small enough so that $S_\rho = G(\zero)$. For any neighborhood $U$ of $\bp$, we have
\begin{equation}
\label{pattersonsullivanbound}
\mu(\bp) \leq \liminf_{s\searrow\delta} \mu_s(U)
= \liminf_{s\searrow\delta} \frac{1}{\Sigma_{s,k}} \sum_{x\in G(\zero)\cap U} k(x) e^{-s\dox x}.
\end{equation}

\begin{lemma}
\label{lemmaSbpasymp}
\[
\lb h(\zero) | x\rb_\zero \asymp_\plus 0 \all x\in S_\bp.
\]
\end{lemma}
\begin{subproof}
Since $S_\bp$ is $\bp$-bounded, Gromov's inequality implies that
\[
\lb h(\zero) | x\rb_\zero \wedge \lb h(\zero)|\bp\rb_\zero \asymp_\plus 0
\]
for all $h\in G_\bp$ and $x\in S_\bp$. Denote the implied constant by $\sigma$. For all $h\in G_\bp$ such that $\lb h(\zero)|\bp\rb_\zero > \sigma$, we have $\lb h(\zero) | x\rb_\zero \leq \sigma \all x\in S_\bp$. Since this applies to all but finitely many $h\in G_\bp$, (c) of Proposition \ref{propositionbasicidentities} completes the proof.
\end{subproof}

Let $T$ be a transversal of $G_\bp\backslash G$ such that $T(\zero) \subset S_\bp$. Then by Lemma \ref{lemmaSbpasymp},
\[
\dox{h(x)} \asymp_\plus \dogo h + \dox x \all h\in G_\bp \all x\in T(\zero).
\]
Thus for all $s > \delta$ and $V\subset X$,
\begin{equation}
\label{SigmaVasymp}
\begin{split}
\sum_{x\in G(\zero)\cap U} k(x) e^{-s\dox x}
&=_\pt \sum_{h\in G_\bp} \sum_{x\in hT(\zero)\cap U} k(e^{\dox x}) e^{-s\dox x}\\
&\asymp_\times \sum_{h\in G_\bp} \sum_{x\in T(\zero)\cap h^{-1}(U)} k(e^{\dogo h + \dox x}) e^{-s[\dogo h + \dox x]}.
\end{split}
\end{equation}
Now fix $0 < \epsilon < \delta - \delta_\bp$, and note that by \eqref{patterson},
\[
k(R) \leq k(\lambda R) \lesssim_{\times,\epsilon} \lambda^\epsilon k(R) \all \lambda > 1 \all R \geq 1.
\]
Thus setting $V = U$ in \eqref{SigmaVasymp} gives
\[
\sum_{x\in G(\zero)\cap U} k(x) e^{-s\dox x} \lesssim_{\times,\epsilon}  \sum_{\substack{h\in G_\bp \\ h(S_\bp)\cap U\neq\smallemptyset}} e^{-(s - \epsilon)\dogo h} \sum_{x\in T(\zero)} k(x)e^{-s\dox x},
\]
while setting $V = X$ gives
\begin{align*}
\Sigma_{s,k} = \sum_{x\in G(\zero)} k(x) e^{-s\dox x} \gtrsim_\times \sum_{h\in G_\bp} e^{-s\dogo h} \sum_{x\in T(\zero)} k(x)e^{-s\dox x}.
\end{align*}
Dividing these inequalities and combining with \eqref{pattersonsullivanbound} gives
\[
\mu(\bp)
\lesssim_{\times,\epsilon} \liminf_{s\searrow\delta} \frac{1}{\Sigma_s(G_\bp)} \sum_{\substack{h\in G_\bp \\ h(S_\bp)\cap U \neq\smallemptyset}} e^{-(s - \epsilon)\dogo h}
= \frac{1}{\Sigma_\delta(G_\bp)} \sum_{\substack{h\in G_\bp \\ h(S_\bp)\cap U \neq\smallemptyset}} e^{-(\delta - \epsilon)\dogo h} ~.
\]
Note that the right hand series converges since $\delta - \epsilon > \delta_\bp$ by construction. As the neighborhood $U$ shrinks, the series converges to zero. This completes the proof.
\end{proof}

Combining Theorem \ref{theoremGFdivergencetype} with Proposition \ref{propositiondeltagtrdeltap} gives the following immediate corollary:

\begin{corollary}
If for all $\bp\in P$, $G_\bp$ is of divergence type, then $G$ is of divergence type.
\end{corollary}
Thus in some sense divergence type can be ``checked locally'' just like the properties of finite generation and finite Poincar\'e exponent (cf. Corollary \ref{corollaryGF}).
\begin{corollary}
Every convex-cobounded group is of divergence type.
\end{corollary}

\begin{remark}
It is somewhat awkward that it seems to be difficult or impossible to prove Theorem \ref{theoremGFdivergencetype} via any of the equivalent conditions of Observation \ref{observationnonatomic} other than (B). Specifically, the fact that the above argument works for the measure constructed in Theorem \ref{theorempattersonsullivangeneral} (the ``Patterson--Sullivan measure'') but not for other $\delta$-conformal measures seems rather asymmetric. However, after some thought one realizes that it would be impossible for a proof along similar lines to work for every $\delta$-conformal measure. This is because the above proof shows that the Patterson--Sullivan measure $\mu$ satisfies
\begin{equation}
\label{pattersonsullivanproperty}
\text{$\mu(\bp) = 0$ for all $\bp\in P$ satisfying $\delta > \delta_\bp$},
\end{equation}
but there are geometrically finite groups for which \eqref{pattersonsullivanproperty} does not hold for all $\delta$-conformal measures $\mu$. Specifically, one may construct geometrically finite groups of convergence type (cf. \cite[Th\'eor\`eme C]{DOP}) such that $\delta_\bp < \delta$ for some $\bp\in P$; the following proposition shows that there exists a $\delta$-conformal measure for which \eqref{pattersonsullivanproperty} fails:
\end{remark}

\begin{proposition}
\label{propositionGFconvergencetype}
If $G$ is of convergence type, then for each $\bp\in P$ there exists a $\delta$-conformal measure supported on $G(\bp)$.
\end{proposition}
\begin{proof}
Let
\[
\mu = \sum_{g(\bp)\in G(\bp)} [g'(\bp)]^\delta \delta_{g(\bp)};
\]
clearly $\mu$ is a $\delta$-conformal measure, but we may have $\mu(\del X) = \infty$. To prove that this is not the case, as before we let $T$ be a transversal of $G_\bp\backslash G$ such that $T(\zero)\subset S_\bp$. Then
\[
\mu(\del X) = \sum_{g(\bp)\in G(\bp)} [g'(\bp)]^\delta = \sum_{g\in T^{-1}} [g'(\bp)]^\delta \asymp_\times \sum_{g\in T^{-1}} e^{-\delta\dogo g} \leq \Sigma_\delta(G) < \infty.
\]
\end{proof}
Proposition \ref{propositionGFconvergencetype} yields the following characterization of when there exists a unique $\delta$-conformal measure:
\begin{corollary}
\label{corollaryuniquenesscharacterization}
The following are equivalent:
\begin{itemize}
\item[(A)] There exists a unique $\delta$-conformal measure on $\Lambda$.
\item[(B)] Either $G$ is of divergence type, or $\#(P) = 1$.
\end{itemize}
\end{corollary}

\section{The global measure formula}
\label{subsectionglobalmeasure}
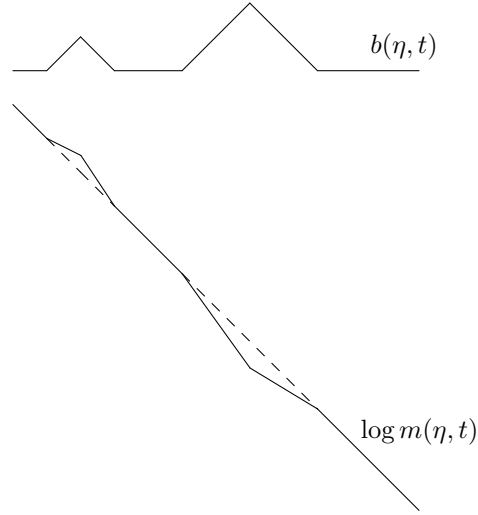
\begin{figure}
\begin{tikzpicture}[line cap=round,line join=round,>=triangle 45,scale=0.45]
\clip(-0.8063966880341881,-10.0638437660102) rectangle (14.110082748891738,6.841499595839061);
\draw (0.0,4.0)-- (1.0,4.0);
\draw (2.0,5.0)-- (3.0,4.0);
\draw (2.0,5.0)-- (1.0,4.0);
\draw (3.0,4.0)-- (5.0,4.0);
\draw (7.0,6.0)-- (5.0,4.0);
\draw (7.0,6.0)-- (9.0,4.0);
\draw (9.0,4.0)-- (12.0,4.0);
\draw (0.0,3.0)-- (1.0,2.0);
\draw [dash pattern=on 4pt off 4pt] (1.0,2.0)-- (2.0,1.0);
\draw [dash pattern=on 4pt off 4pt] (2.0,1.0)-- (3.0,0.0);
\draw (3.0,0.0)-- (5.0,-2.0);
\draw [dash pattern=on 4pt off 4pt] (5.0,-2.0)-- (9.0,-6.0);
\draw (9.0,-6.0)-- (12.0,-9.0);
\draw (1.0,2.0)-- (2.009795868802983,1.4948904390904918);
\draw (2.009795868802983,1.4948904390904918)-- (3.0,0.0);
\draw (5.0,-2.0)-- (7.006776519282634,-4.79037428534092);
\draw (7.006776519282634,-4.79037428534092)-- (9.0,-6.0);
\draw (10.285344431731245,5.388099035318084) node[anchor=north west] {$b(\eta,t)$};
\draw (10.0,-6.0) node[anchor=north west] {$\log m(\eta,t)$};
\end{tikzpicture}
\caption[Cusp excursion and ball measure functions]{A possible (approximate) graph of the functions $t\mapsto b(\eta,t)$ and $t\mapsto \log m(\eta,t)$ (cf. \eqref{metat} and \eqref{betat}). The graph indicates that there are at least two inequivalent parabolic points $\bp_1,\bp_2\in P$, which satisfy $\NN_{\bp_i}(R) \asymp_\times R^{2\delta} \II_{\bp_i}(R) \asymp_\times R^{k_i}$ for some $k_1 < 2\delta < k_2$. The dotted line in the second graph is just the line $y = -\delta t$.\newline
Note the relation between the two graphs, which may be either direct or inverted depending on the functions $\NN_\bp$. Specifically, the relation is direct for the first cusp but inverted for the second cusp.}
\label{figuremetat}
\end{figure}

In this section and the next, we fix a $\delta$-quasiconformal measure $\mu$, and ask the following geometrical question: Given $\eta\in\Lambda$ and $r > 0$, can we estimate $\mu(B(\eta,r))$? If $G$ is convex-cobounded, then we can show that $\mu$ is Ahlfors $\delta$-regular (Corollary \ref{corollaryCCBregular}), but in general the measure $\mu(B(\eta,r))$ will depend on the point $\eta$, in a manner described by the \emph{global measure formula}. To describe the global measure formula, we need to introduce some notation:

\begin{notation}
\label{notationmetat}
Given $\xi = g(\bp)\in\Lbp$, let $t_\xi > 0$ be the unique number such that
\[
H_\xi = H_{\xi,t_\xi} = g(H_\bp) = g(H_{\bp,t_0}),
\]
i.e. $t_\xi = t_0 + \busemann_\xi(\zero,g(\zero))$. (Note that $t_\bp = t_0$ for all $\bp\in P$.) Fix $\theta > 0$ large to be determined below (cf. Proposition \ref{propositionmetadecreasing}). For each $\eta\in\Lambda$ and $t > 0$, let $\eta_t = \geo\zero\eta_t$, and write
\begin{equation}
\label{metat}
m(\eta,t) = \begin{cases}
e^{-\delta t} & \eta_t\notin\bigcup(\scrH) \\
e^{-\delta t_\xi} [\II_\bp(e^{t - t_\xi - \theta}) + \mu(\bp)] & \eta_t\in H_\xi \text{ and }t \leq \lb \xi|\eta\rb_\zero\\
e^{-\delta (2\lb \xi|\eta\rb_\zero - t_\xi)} \NN_\bp(e^{2\lb \xi|\eta\rb_\zero - t - t_\xi - \theta}) & \eta_t\in H_\xi \text{ and } t > \lb \xi|\eta\rb_\zero
\end{cases}
\end{equation}
(cf. Figure \ref{figuremetat}.) Here we use the notation
\begin{align*}
\II_\bp(R) &= \sum_{\substack{h\in G_\bp \\ \|h\|_\bp \geq R}} \|h\|_\bp^{-2\delta}\\
\NN_\bp(R) &= \NN_{\EE_\bp,G_\bp}(R) = \#\{h\in G_\bp : \|h\|_\bp \leq R\}
\end{align*}
where
\[
\|h\|_\bp = \Dist_\bp(\zero,h(\zero)) = e^{(1/2)\dogo h}
\all h\in G_\bp.
\]
\end{notation}

\begin{theorem}[Global measure formula; cf. {\cite[Theorem 2]{StratmannVelani}} and {\cite[Th\'eor\`eme 3.2]{Schapira}}]
\label{theoremglobalmeasure}
For all $\eta\in\Lambda$ and $t > 0$,
\begin{equation}
\label{globalmeasureformula}
m(\eta,t + \sigma) \lesssim_\times \mu(B(\eta,e^{-t})) \lesssim_\times m(\eta,t - \sigma),
\end{equation}
where $\sigma > 0$ is independent of $\eta$ and $t$ (but may depend on $\theta$).
\end{theorem}

\begin{corollary}
\label{corollaryCCBregular}
If $G$ is convex-cobounded, then
\begin{equation}
\label{ahlforsregular}
\mu(B(\eta,r)) \asymp_\times r^\delta \all \eta\in\Lambda \all 0 < r \leq 1,
\end{equation}
i.e. $\mu$ is Ahlfors $\delta$-regular.
\end{corollary}
\begin{proof}
If $G$ is convex-cobounded then $\scrH = \emptyset$, so $m(\eta,t) = e^{-\delta t}\all\eta,t$, and thus \eqref{globalmeasureformula} reduces to \eqref{ahlforsregular}.
\end{proof}

\begin{remark}
Corollary \ref{corollaryCCBregular} can be deduced directly from Lemma \ref{lemmaCCBglobalmeasure} below.
\end{remark}

We will prove Theorem \ref{theoremglobalmeasure} in the next section. For now, we investigate more closely the function $t\mapsto m(\eta,t)$ defined by \eqref{metat}. The main result of this section is the following proposition, which will be used in the proof of Theorem \ref{theoremglobalmeasure}:

\begin{proposition}
\label{propositionmetadecreasing}
If $\theta$ is chosen sufficiently large, then for all $\eta\in\Lambda$ and $0 < t_1 < t_2$,
\begin{equation}
\label{metadecreasing}
m(\eta,t_2) \lesssim_{\times,\theta} m(\eta,t_1).
\end{equation}
\end{proposition}

The proof of Proposition \ref{propositionmetadecreasing} itself requires several lemmas.

\begin{lemma}
\label{lemmabetatinterpretation}
Fix $\xi,\eta\in \del X$ and $t > 0$, and let $x = \eta_t$. Then
\begin{equation}
\label{betatinterpretation}
\busemann_\xi(\zero,x) \asymp_\plus t \wedge (2\lb\xi|\eta\rb_\zero - t).
\end{equation}
\end{lemma}
\begin{proof}
Since $\lb \zero|\eta\rb_x = 0$, Gromov's inequality gives $\lb \zero|\xi\rb_x \wedge \lb \xi|\eta\rb_x \asymp_\plus 0$.

\begin{itemize}
\item[Case 1:] $\lb \zero|\xi\rb_x \asymp_\plus \zero$. In this case, by (h) of Proposition \ref{propositionbasicidentities},
\[
\busemann_\xi(\zero,x) = - \busemann_\xi(x,\zero) = -[2\lb \zero|\xi\rb_x - \dox x] \asymp_\plus \dox x = t,
\]
while (g) of Proposition \ref{propositionbasicidentities} gives
\[
\lb \xi|\eta\rb_\zero = \lb \xi|\eta\rb_x + \frac12[\busemann_\xi(\zero,x) + \busemann_\eta(\zero,x)] \gtrsim_\plus \frac12[t + t] = t;
\]
thus $\busemann_\xi(\zero,x) \asymp_\plus t \asymp_\plus t\wedge(2\lb \xi|\eta\rb_\zero - t)$.
\item[Case 2:] $\lb \xi|\eta\rb_x \asymp_\plus \zero$. In this case, (g) of Proposition \ref{propositionbasicidentities} gives
\[
\lb \xi|\eta\rb_\zero \asymp_\plus \frac12[\busemann_\xi(\zero,x) + \busemann_\eta(\zero,x)] = \frac12[\busemann_\xi(\zero,x) + t] \lesssim_\plus \frac12[t + t] = t;
\]
thus $\busemann_\xi(\zero,x) \asymp_\plus 2\lb \xi|\eta\rb_\zero - t \asymp_\plus t\wedge(2\lb \xi|\eta\rb_\zero - t)$.
\end{itemize}
\end{proof}

\begin{corollary}
\label{corollarybetat}
The function
\begin{equation}
\label{betat}
b(\eta,t) = \begin{cases}
0 & \eta_t\notin\bigcup(\scrH) \\
t\wedge(2\lb \xi|\eta\rb_\zero - t) - t_\xi & \eta_t\in H_\xi
\end{cases}
\end{equation}
satisfies
\begin{equation}
\label{betatdoubling}
b(\eta,t + \tau) \asymp_{\plus,\tau} b(\eta,t - \tau).
\end{equation}
\end{corollary}
\begin{proof}
Indeed, by Lemma \ref{lemmabetatinterpretation},
\begin{align*}
b(\eta,t) &\asymp_\plus \begin{cases}
0 & \eta_t\notin\bigcup(\scrH) \\
\busemann_\xi(\zero,\eta_t) - t_\xi & \eta_t\in H_\xi
\end{cases}\\
&=_\pt 0\vee \max_{\xi\in\Lbp} (\busemann_\xi(\zero,\eta_t) - t_\xi).
\end{align*}
The right hand side is $1$-Lipschitz continuous with respect to $t$, which demonstrates \eqref{betatdoubling}.
\end{proof}

\begin{lemma}
\label{lemmatop}
For all $\xi\in G(\bp) \subset \Lbp$, $\bp\in P$, there exists $g\in G$ such that
\begin{equation}
\label{top}
\text{$\xi = g(\bp)$, $\dogo g \asymp_\plus t_\xi$, and $\{\eta\in\del X  : \geo\zero\eta\cap H_\xi\neq\emptyset \} \subset \Shad(g(\zero),\sigma)$,}
\end{equation}
where $\sigma > 0$ is independent of $\xi$.
\end{lemma}
\begin{proof}
Write $\xi = g(\bp)$ for some $g\in G$. Since $x := \xi_{t_\xi} \in \del H_\xi$, Lemma \ref{lemmaboundedparabolic}(D) shows that
\[
\dist(g^{-1}(x),h(\zero)) \asymp_\plus 0
\]
for some $h\in G_\bp$. We claim that $gh$ is the desired isometry. Clearly $\dogo{gh} \asymp_\plus \dox x = t_\xi$. Fix $\eta\in\del X$ such that $\geo\zero\eta\cap H_\xi\neq\emptyset$, say $\eta_t\in H_\xi$. By Lemma \ref{lemmabetatinterpretation}, we have
\[
\dox x = t_\xi < \busemann_\xi(\zero,\eta_t) \asymp_\plus t\wedge(2\lb \xi|\eta\rb_\zero - t) \leq \lb \xi|\eta\rb_\zero \leq \lb x|\eta\rb_\zero,
\]
i.e. $\eta\in \Shad(x,\sigma) \subset \Shad(g(\zero),\sigma + \tau)$ for some $\sigma,\tau > 0$.
\end{proof}

\begin{proof}[Proof of Proposition \ref{propositionmetadecreasing}]
Fix $\eta\in\Lambda$ and $0 < t_1 < t_2$.
\begin{itemize}
\item[Case 1:] $\eta_{t_1},\eta_{t_2}\in H_\xi$ for some $\xi = g(\bp) \in \Lbp$, $g$ satisfying \eqref{top}. In this case, \eqref{metadecreasing} follows immediately from \eqref{metat} unless $t_1 \leq \lb \xi|\eta\rb_\zero < t_2$. If the latter holds, then
\begin{align*}
m(\eta,t_1) &\geq \lim_{t\nearrow\lb \xi|\eta\rb_\zero}m(\eta,t) = e^{-\delta t_\xi} [\II_\bp(e^{\lb \xi|\eta\rb_\zero - t_\xi - \theta}) + \mu(\bp)]\\
m(\eta,t_2) &\leq \lim_{t\searrow\lb \xi|\eta\rb_\zero}m(\eta,t) = e^{-\delta(2\lb \xi|\eta\rb_\zero - t_\xi)} \NN_\bp(e^{\lb \xi|\eta\rb_\zero - t_\xi - \theta}).
\end{align*}
Consequently, to demonstrate \eqref{metadecreasing} it suffices to show that
\begin{equation}
\label{ETSmetadecreasing}
\NN_\bp(e^t) \lesssim_{\times,\theta} e^{2\delta t}\II_\bp(e^t),
\end{equation}
where $t := \lb \xi|\eta\rb_\zero - t_\xi - \theta > 0$.

To demonstrate \eqref{ETSmetadecreasing}, let $\zeta = g^{-1}(\eta)\in\Lambda$. We have
\[
\lb \bp|\zeta\rb_\zero = \lb \xi|\eta\rb_{g(\zero)} \asymp_\plus \lb \xi|\eta\rb_\zero - \dogo g \asymp_\plus \lb \xi|\eta\rb_\zero - t_\xi = t + \theta
\]
and thus
\[
\Dist_\bp(\zero,\zeta) \asymp_\times e^{t + \theta}.
\]
Since $\bp$ is a bounded parabolic point, there exists $h_\zeta\in G_\bp$ such that $\Dist_\bp(h_\zeta(\zero),\zeta) \lesssim_\times 1$. Denoting all implied constants by $C$, we have
\begin{align*}
C^{-1} e^{t + \theta} - C \leq \Dist_\bp(\zero,\zeta) - \Dist_\bp(h_\zeta(\zero),\zeta) 
&\leq \|h_\zeta\|_\bp\\ 
&\leq \Dist_\bp(\zero,\zeta) + \Dist_\bp(h_\zeta(\zero),\zeta) \leq C e^{t + \theta} + C.
\end{align*}
Choosing $\theta \geq \log(4C)$, we have
\[
2 e^t \leq \|h_\zeta\|_\bp \leq 2C e^{t + \theta} \text{ unless $e^{t + \theta} \leq 2C^2$}.
\]
If $2 e^t \leq \|h_\zeta\|_\bp \leq 2C e^{t + \theta}$, then for all $h\in G_\bp$ satisfying $\|h\|_\bp\leq e^t$ we have $e^t \leq \|h_\zeta h\|_\bp \lesssim_{\times,\theta} e^t$; it follows that
\[
\II_\bp(e^t) \geq \sum_{h\in G_\bp} \|h_\zeta h\|_\bp^{-2\delta} \asymp_{\times,\theta} e^{-2\delta t} \NN_\bp(e^t),
\]
thus demonstrating \eqref{ETSmetadecreasing}. On the other hand, if $e^{t + \theta} \leq 2C^2$, then both sides of \eqref{ETSmetadecreasing} are bounded from above and below independent of $t$.
\item[Case 2:] No such $\xi$ exists. In this case, for each $i$ write $\eta_i\in H_{\xi_i}$ for some $\xi_i = g_i(\bp_i)\in\Lbp$ if such a $\xi_i$ exists. If $\xi_1$ exists, let $s_1 > t_1$ be the smallest number such that $\eta_{s_1}\in\del H_{\xi_1}$, and if $\xi_2$ exists, let $s_2 < t_2$ be the largest number such that $\eta_{s_2}\in \del H_{\xi_2}$. If $\xi_i$ does not exist, let $s_i = t_i$. Then $t_1 \leq s_1 \leq s_2 \leq t_2$. Since $m(\eta,s_i) = e^{-\delta s_i}$, we have $m(\eta,s_2) \leq m(\eta,s_1)$, so to complete the proof it suffices to show that
\begin{align*}
m(\eta,s_1) &\lesssim_{\times,\theta} m(\eta,t_1) \text{ and}\\
m(\eta,s_2) &\gtrsim_{\times,\theta} m(\eta,t_2).
\end{align*}
By Case 1, it suffices to show that
\begin{align*}
m(\eta,s_1) &\lesssim_\times \lim_{t\nearrow s_1}m(\eta,t) \text{ if $\xi_1$ exists, and}\\
m(\eta,s_2) &\gtrsim_\times \lim_{t\searrow s_2} m(\eta,t) \text{ if $\xi_2$ exists.}
\end{align*}
Comparing with \eqref{metat}, we see that the desired formulas are
\begin{align*}
e^{-\delta s_1} &\lesssim_\times e^{-\delta (2\lb \xi|\eta\rb_\zero - t_{\xi_1})} \NN_\bp(e^{2\lb \xi_1|\eta\rb_\zero - s_1 - t_{\xi_1}})\\
e^{-\delta s_2} &\gtrsim_\times  e^{-\delta t_{\xi_2}}[\II_\bp(e^{s_2 - t_{\xi_2}}) + \mu(\bp)],
\end{align*}
which follow upon observing that the definitions of $s_1$ and $s_2$ imply that $s_1 \asymp_\plus 2\lb \xi|\eta\rb_\zero - t_{\xi_1}$ and $s_2 \asymp_\plus t_{\xi_2}$ (cf. Lemma \ref{lemmabetatinterpretation}).
\end{itemize}
\end{proof}

\section{Proof of the global measure formula}
\label{subsectionglobalmeasureproof}
Although we have finished the proof of Proposition \ref{propositionmetadecreasing}, we still need a few lemmas before we can begin the proof of Theorem \ref{theoremglobalmeasure}. Throughout these lemmas, we fix $\bp\in P$, and let
\[
R_\bp = \sup_{x\in S_\bp} \Dist_\bp(\zero,x) < \infty.
\]
Here $S_\bp\subset\EE_\bp$ is a $\bp$-bounded set satisfying $\Lambda\butnot\{\bp\}\subset G_\bp(S_\bp)$, as in Standing Assumptions \ref{standingassumptionsGFmeasures}.

\begin{lemma}
\label{lemmaglobalmeasure1}
For all $A\subset G_\bp$,
\begin{equation}
\label{ASmeasure}
\mu\left(\bigcup_{h\in A}h(S_\bp)\right) \asymp_\times \sum_{h\in A} e^{-\delta\dogo h} = \sum_{h\in A} \|h\|_\bp^{-2\delta}.
\end{equation}
\end{lemma}
\begin{proof}
As the equality follows from Observation \ref{observationeuclideanparabolicasymp}, we proceed to demonstrate the asymptotic. By Lemma \ref{lemmaSbpasymp}, there exists $\sigma > 0$ such that $S_\bp \subset \Shad_{h^{-1}(\zero)}(\zero,\sigma)$ for all $h\in G_\bp$. Then by the Bounded Distortion Lemma \ref{lemmaboundeddistortion},
\[
\mu(h(S_\bp)) = \int_{S_\bp} (\overline h')^\delta \;\dee \mu \asymp_{\times,\sigma} e^{-\delta\dogo h} \mu(S_\bp) \asymp_\times e^{-\delta\dogo h}.
\]
(In the last asymptotic, we have used the fact that $\mu(S_\bp) > 0$, which follows from the fact that $\Lambda\butnot\{\bp\} \subset G_\bp(S_\bp)$ together with the fact that $\mu$ is not a pointmass (Corollary \ref{corollaryquasiconformalpointmass}).) Combining with the subadditivity of $\mu$ gives the $\lesssim$ direction of the first asymptotic of \eqref{ASmeasure}. To get the $\gtrsim$ direction, we observe that since $S_\bp$ is $\bp$-bounded, the strong discreteness of $G_\bp$ implies that $S_\bp\cap h(S_\bp)\neq\emptyset$ for only finitely many $h\in G_\bp$; it follows that the function $\eta\mapsto\#\{h\in G_\bp : \eta\in h(S_\bp)\}$ is bounded, and thus
\begin{align*}
\mu\left(\bigcup_{h\in A}h(S_\bp)\right) 
&\asymp_\times \int \#\{h\in G_\bp : \eta\in h(S_\bp)\} \;\dee\mu(\eta)\\
&= \sum_{h\in A}\mu(h(S_\bp))\\
&\asymp_\times \sum_{h\in A} e^{-\delta\dogo h} ~.
\end{align*}
\end{proof}

\begin{corollary}
\label{corollaryglobalmeasure11}
For all $r > 0$,
\begin{equation}
\label{muBbprbounds}
\II_\bp\left(\frac2r\right) \lesssim_\times \mu\big(B(\bp,r)\butnot\{\bp\}\big) \lesssim_\times \II_\bp\left(\frac1{2r}\right)
\end{equation}
\end{corollary}
\begin{proof}
Since
\[
\bigcup_{\substack{h\in G_\bp \\ \|h\|_\bp \geq R + R_\bp}} h(S_\bp)
\subset B(\bp,1/R)\butnot\{\bp\} = \EE_\bp \butnot B_\bp(\zero,R)
\subset \bigcup_{\substack{h\in G_\bp \\ \|h\|_\bp \geq R - R_\bp}} h(S_\bp),
\]
Lemma \ref{lemmaglobalmeasure1} gives
\[
\II_\bp\left(\frac1r + R_\bp\right) \lesssim_\times \mu(B(\bp,r)) \lesssim_\times \II_\bp\left(\frac1r - R_\bp\right),
\]
thus proving the lemma if $r \leq 1/(2R_\bp)$. But when $r > 1/(2 R_\bp)$, all terms of \eqref{muBbprbounds} are bounded from above and below independent of $r$.
\end{proof}

Adding $\mu(\bp)$ to all sides of \eqref{muBbprbounds} gives
\begin{equation}
\label{muBbprbounds2}
\II_\bp\left(\frac2r\right) + \mu(\bp) \lesssim_\times \mu(B(\bp,r)) \lesssim_\times \II_\bp\left(\frac1{2r}\right) + \mu(\bp).
\end{equation}

\begin{figure}
\begin{center}
\begin{tabular}{ll}
\begin{tikzpicture}[line cap=round,line join=round,>=triangle 45,scale=.81]
\clip(-2.7,-2.7) rectangle (4.3,4.3);
\draw(0.0,0.0) circle (4.0cm);
\draw(2.42888,0.1875) circle (1.551112790229163cm);
\draw (0.0,0.0)-- (4.0,0.0);
\draw(3.212205575454762,0.0) circle (0.3760152197599791cm);
\draw [dash pattern=on 3pt off 3pt] (0.0,0.0)-- (3.9725001940685334,0.4682330703030985);
\draw [dash pattern=on 3pt off 3pt] (0.0,0.0)-- (3.9725001940685334,-0.4682330703030985);
\draw [shift={(0.0,0.0)},line width=1.2000000000000002pt]  plot[domain=-0.11732726409329874:0.11732726409329915,variable=\t]({1.0*4.0*cos(\t r)+-0.0*4.0*sin(\t r)},{0.0*4.0*cos(\t r)+1.0*4.0*sin(\t r)});
\begin{scriptsize}
\draw [fill=black] (0.0,0.0) circle (1pt);
\draw[color=black] (-0.2491340591615706,-0.0014223831473892945) node {$\zero$};
\draw [fill=black] (4.0,0.0) circle (1pt);
\draw[color=black] (4.211681114182143,0.0) node {$\eta$};
\draw [fill=black] (3.98,0.3) circle (1pt);
\draw[color=black] (4.17,0.4) node {$\xi$};
\draw [fill=black] (3.384151714048233,0.0) circle (1pt);
\draw[color=black] (3.2344150960861695,-0.1380274618110051) node {$\eta_t$};
\end{scriptsize}
\end{tikzpicture}
\begin{tikzpicture}[line cap=round,line join=round,>=triangle 45,scale=.85]
\clip(-4.152066304477812,-0.497280356328431) rectangle (4.307443018392998,5.909099879057296);
\draw (-3.5,0.0)-- (3.5,0.0);
\draw (-3.5,1.5)-- (3.5,1.5);
\draw [shift={(0.40000000000000036,0.0)}] plot[domain=0.3996780691889681:3.141592653589793,variable=\t]({1.0*2.518491612056709*cos(\t r)+-0.0*2.518491612056709*sin(\t r)},{0.0*2.518491612056709*cos(\t r)+1.0*2.518491612056709*sin(\t r)});
\draw (2.4230669786242873,1.5)-- (2.1163456022508336,1.2434075389985795);
\draw (-2.1184916120567086,0.0)-- (-2.12,5.4);
\draw (2.1163456022508336,1.2434075389985795)-- (2.12,5.4);
\begin{scriptsize}
\draw [fill=black] (2.72,0.98) circle (1pt);
\draw[color=black] (2.9026704848505824,1.0) node {$\zero$};
\draw [fill=black] (-2.1184916120567086,0.0) circle (1pt);
\draw[color=black] (-2.284181946690644,-0.1267908969326299) node {$\eta$};
\draw [fill=black] (2.4230669786242873,1.5) circle (1pt);
\draw[color=black] (2.58,1.68) node {$\eta_{t_\xi}$};
\draw [fill=black] (2.1163456022508336,1.2434075389985795) circle (1pt);
\draw[color=black] (1.8073209579773947,1.092736906911882) node {$g(\zero)$};
\draw[color=black] (-0.1,5.662106906126763) node {$\xi=\infty$};
\draw [fill=black] (1.6640939374267694,2.1782714517159905) circle (1pt);
\draw[color=black] (1.420712647267375,2.1424570418666518) node {$\eta_t$};
\draw [dash pattern=on 1pt off 1pt] (1.6640939374267694,2.1782714517159905)-- (2.1173124637830534,2.3431371337242592);
\draw[color=black] (-3.084181946690644,1.867908969326299) node {$H_{\xi}$};
\end{scriptsize}
\end{tikzpicture}
\end{tabular}
\end{center}

\caption[Estimating measures of balls via information ``at infinity'']{Cusp excursion in the ball model (left) and upper half-space model (right). Since $\xi = g(\bp)\in B(\eta,e^{-t})$, our estimate of $\mu(B(\eta,e^{-t}))$ is based on the function $\II_\bp$, which captures information ``at infinity'' about the cusp $\bp$. In the right-hand picture, the measure of $B(\eta,e^{-t})$ can be estimated by considering the measure from the perspective of $g(\zero)$ of a small ball around $\xi$.}
\label{figurecuspright}
\end{figure}

\begin{corollary}[Cf. Figure \ref{figurecuspright}]
\label{corollaryglobalmeasure12}
Fix $\eta\in\Lambda$ and $t > 0$ such that $\eta_t\in H_\xi$ for some $\xi = g(\bp)\in\Lbp$ satisfying $t \leq \lb \xi|\eta\rb_\zero - \log(2)$. Then
\[
e^{-\delta t_\xi} [\II_\bp(e^{t - t_\xi + \sigma}) + \mu(\bp)] \lesssim_\times \mu\big(B(\eta,e^{-t})\big) \lesssim_\times e^{-\delta t_\xi} [\II_\bp(e^{t - t_\xi - \sigma}) + \mu(\bp)],
\]
where $\sigma > 0$ is independent of $\eta$ and $t$.
\end{corollary}
\begin{proof}
The inequality $\lb \xi|\eta\rb_\zero \geq t + \log(2)$ implies that
\[
B(\xi,e^{-t}/2) \subset B(\eta,e^{-t}) \subset B(\xi,2e^{-t}).
\]
Without loss of generality suppose that $g$ satisfies \eqref{top}. Since $t > t_\xi$, \eqref{ShadcontainsB} guarantees that $B(\xi,2e^{-t}) \subset \Shad(g(\zero),\sigma_0)$ for some $\sigma_0 > 0$ independent of $\eta$ and $t$. Then by the Bounded Distortion Lemma \ref{lemmaboundeddistortion}, we have
\begin{align*}
B\left(\bp, e^{-(t - t_\xi)}/(2C)\right) &\subset g^{-1}(B(\xi,e^{-t}/2))\\ 
&\subset g^{-1}(B(\eta,e^{-t}))\\
&\subset g^{-1}(B(\xi,2e^{-t}))\\
&\subset B\left(\bp,2C e^{-(t - t_\xi)}\right)
\end{align*}
for some $C > 0$, and thus
\[
e^{-\delta t_\xi} \mu\left(B\left(\bp, e^{-(t - t_\xi)}/(2C)\right)\right) \lesssim_\times  \mu(B(\eta,e^{-t})) \lesssim_\times e^{-\delta t_\xi} \mu\left(B\left(\bp,2C e^{-(t - t_\xi)}\right)\right).
\]
Combining with \eqref{muBbprbounds2} completes the proof.
\end{proof}

\begin{lemma}
\label{lemmaglobalmeasure2}
For all $\eta\in \Lambda\butnot\{\bp\}$ and $3R_\bp \leq R\leq \Dist_\bp(\zero,\eta)/2$,
\[
\Dist_\bp(\zero,\eta)^{-2\delta} \NN_\bp(R/2) \lesssim_\times \mu(B_\bp(\eta,R)) \lesssim_\times \Dist_\bp(\zero,\eta)^{-2\delta} \NN_\bp(2R).
\]
\end{lemma}
\begin{proof}
Since $\eta\in \Lambda\butnot\{\bp\} \subset G_\bp(S_\bp)$, there exists $h_\eta\in G_\bp$ such that $\eta\in h_\eta(S_\bp)$. Since
\[
\bigcup_{\substack{h\in G_\bp \\ \|h\|_\bp \leq R - R_\bp}} h_\eta h(S_\bp) \subset B_\bp(\eta,R) \subset \bigcup_{\substack{h\in G_\bp \\ \|h\|_\bp \leq R + R_\bp}} h_\eta h(S_\bp),
\]
Lemma \ref{lemmaglobalmeasure1} gives
\[
\sum_{\substack{h\in G_\bp \\ \|h\|_\bp \leq R - R_\bp}} \|h_\eta h\|_\bp^{-2\delta} \lesssim_\times \mu(B_\bp(\eta,R)) \lesssim_\times \sum_{\substack{h\in G_\bp \\ \|h\|_\bp \leq R + R_\bp}} \|h_\eta h\|_\bp^{-2\delta}.
\]
The proof will be complete if we can show that for each $h\in G_\bp$ such that $\|h\|_\bp \leq R + R_\bp$, we have
\begin{equation}
\label{heta}
\|h_\eta h\|_\bp \asymp_\times \Dist_\bp(\zero,\eta).
\end{equation}
And indeed,
\[
\Dist_\bp(\eta,h_\eta h(\zero)) \leq \Dist_\bp(\eta,h_\eta(\zero)) + \|h\|_\bp \leq R_\bp + (R + R_\bp) \leq \frac56 \Dist_\bp(\zero,\eta),
\]
demonstrating \eqref{heta} with an implied constant of $6$.
\end{proof}

\begin{corollary}
\label{corollaryglobalmeasure21}
For all $\eta\in\Lambda\butnot\{\bp\}$ and $6R_\bp \Dist(\bp,\eta)^2 \leq r \leq \Dist(\bp,\eta)/4$, we have
\[
\Dist(\bp,\eta)^{2\delta} \NN_\bp\left(\frac{r}{4\Dist(\bp,\eta)^2}\right)
\lesssim_\times \mu(B(\eta,r))
\lesssim_\times \Dist(\bp,\eta)^{2\delta} \NN_\bp\left(\frac{4r}{\Dist(\bp,\eta)^2}\right).
\]
\end{corollary}

\begin{proof}
By \eqref{GMVT2}, for every $\zeta\in B_\bp\left(\eta,\frac{r}{\Dist(\bp,\eta)(\Dist(\bp,\eta) + r)}\right)$ we have that
\begin{align*}
\Dist(\eta,\zeta) = \frac{\Dist_\bp(\eta,\zeta)}{\Dist_\bp(\zero,\eta)\Dist_\bp(\zero,\zeta)} 
&\leq \frac{\Dist_\bp(\eta,\zeta)}{\Dist_\bp(\zero,\eta)(\Dist_\bp(\zero,\eta) - \Dist_\bp(\eta,\zeta))}\\
&\leq \frac{\frac{r}{\Dist(\bp,\eta)(\Dist(\bp,\eta) + r)}}{\Dist_\bp(\zero,\eta)\left(\Dist_\bp(\zero,\eta) - \frac{r}{\Dist(\bp,\eta)(\Dist(\bp,\eta) + r)}\right)}\\ 
&= r. 
\end{align*}
Analogously, \eqref{GMVT2} also implies that for every $\zeta \in B(\eta,r)$ we have 
\begin{align*}
\Dist_\bp(\eta,\zeta) 
&= \frac{\Dist(\eta,\zeta)}{\Dist(\bp,\eta)\Dist(\bp,\zeta)}\\ 
&\leq \frac{r}{\Dist(\bp,\eta) \left( \Dist(\bp,\eta) - r \right)} \cdot
\end{align*}
Combining these inequalities gives us that 
\[
B_\bp\left(\eta,\frac{r}{\Dist(\bp,\eta)(\Dist(\bp,\eta) + r)}\right)
\subset B(\eta,r)
\subset B_\bp\left(\eta,\frac{r}{\Dist(\bp,\eta)(\Dist(\bp,\eta) - r)}\right)
\]
Now since $r\leq \Dist(\bp,\eta)/4$, we have
\[
B_\bp\left(\eta,\frac{r}{2\Dist(\bp,\eta)^2}\right)
\subset B(\eta,r)
\subset B_\bp\left(\eta,\frac{2r}{\Dist(\bp,\eta)^2}\right).
\]
On the other hand, since $6R_\bp \Dist(\bp,\eta)^2 \leq r \leq \Dist(\bp,\eta)/4$, we have
\[
3R_\bp \leq \frac{r}{2\Dist(\bp,\eta)^2} \leq \frac{2r}{\Dist(\bp,\eta)^2} \leq \frac{\Dist(\bp,\eta)}{2}
\]
whereupon Lemma \ref{lemmaglobalmeasure2} completes the proof.
\end{proof}

\begin{figure}
\begin{center}
\begin{tabular}{ll}
\begin{tikzpicture}[line cap=round,line join=round,>=triangle 45,scale=.81]
\clip(-2.7,-2.7) rectangle (4.3,4.3);
\draw(0.0,0.0) circle (4.0cm);
\draw(2.342806400391098,0.712964401564811) circle (1.551112790229163cm);
\draw (0.0,0.0)-- (4.0,0.0);
\draw(3.212205575454762,0.0) circle (0.3760152197599791cm);
\draw [dash pattern=on 3pt off 3pt] (0.0,0.0)-- (3.9725001940685334,0.4682330703030985);
\draw [dash pattern=on 3pt off 3pt] (0.0,0.0)-- (3.9725001940685334,-0.4682330703030985);
\draw [shift={(0.0,0.0)},line width=1.2000000000000002pt]  plot[domain=-0.11732726409329874:0.11732726409329915,variable=\t]({1.0*4.0*cos(\t r)+-0.0*4.0*sin(\t r)},{0.0*4.0*cos(\t r)+1.0*4.0*sin(\t r)});
\begin{scriptsize}
\draw [fill=black] (0.0,0.0) circle (1pt);
\draw[color=black] (-0.2491340591615706,-0.0014223831473892945) node {$\zero$};
\draw [fill=black] (4.0,0.0) circle (1pt);
\draw[color=black] (4.211681114182143,0.0) node {$\eta$};
\draw [fill=black] (3.825718764158789,1.1678509911642618) circle (1pt);
\draw[color=black] (4.057469971268139,1.2) node {$\xi$};
\draw [fill=black] (3.384151714048233,0.0) circle (1pt);
\draw[color=black] (3.2344150960861695,-0.1380274618110051) node {$\eta_t$};
\end{scriptsize}
\end{tikzpicture}

\begin{tikzpicture}[line cap=round,line join=round,>=triangle 45,scale=.85]
\clip(-4.152066304477812,-0.497280356328431) rectangle (4.307443018392998,5.909099879057296);
\draw (-3.5,0.0)-- (3.5,0.0);
\draw (-3.5,1.5)-- (3.5,1.5);
\draw [shift={(0.40000000000000036,0.0)}] plot[domain=0.3996780691889681:3.141592653589793,variable=\t]({1.0*2.518491612056709*cos(\t r)+-0.0*2.518491612056709*sin(\t r)},{0.0*2.518491612056709*cos(\t r)+1.0*2.518491612056709*sin(\t r)});
\draw (2.4230669786242873,1.5)-- (2.1163456022508336,1.2434075389985795);
\draw (-2.1184916120567086,0.0)-- (-2.12,5.4);
\draw (2.1163456022508336,1.2434075389985795)-- (2.12,5.4);
\draw [dash pattern=on 1pt off 1pt] (-1.5356104253500993,1.6112765998660834)-- (-2.1189952002490946,1.8028361012686287);
\begin{scriptsize}
\draw [fill=black] (2.72,0.98) circle (1pt);
\draw[color=black] (2.9026704848505824,1.0) node {$\zero$};
\draw [fill=black] (-2.1184916120567086,0.0) circle (1pt);
\draw[color=black] (-2.284181946690644,-0.1267908969326299) node {$\eta$};
\draw [fill=black] (2.4230669786242873,1.5) circle (1pt);
\draw[color=black] (2.58,1.68) node {$\eta_{t_\xi}$};
\draw [fill=black] (2.1163456022508336,1.2434075389985795) circle (1pt);
\draw[color=black] (1.8073209579773947,1.092736906911882) node {$g(\zero)$};
\draw[color=black] (-0.1,5.662106906126763) node {$\xi=\infty$};
\draw [fill=black] (-1.5356104253500993,1.6112765998660834) circle (1pt);
\draw[color=black] (-1.2,1.65) node {$\eta_t$};
\draw[color=black] (-3.084181946690644,1.867908969326299) node {$H_{\xi}$};
\end{scriptsize}
\end{tikzpicture}
\end{tabular}

\caption[Estimating measures of balls via ``local'' information]{Cusp excursions in the ball model (left) and upper half-space model (right). Since $\xi = g(\bp)\notin B(\eta,e^{-t})$, our estimate of $\mu(B(\eta,e^{-t}))$ is based on the function $\NN_\bp$, which captures ``local'' information about the cusp $\bp$. In the right-hand picture, the measure of $B(\eta,e^{-t})$ can be estimated by considering the measure from the perspective of $g(\zero)$ of a large ball around $\eta$ taken with respect to the $\Dist_\xi$-metametric.}
\label{figurecuspleft}
\end{center}
\end{figure}

\begin{corollary}[Cf. Figure \ref{figurecuspleft}]
\label{corollaryglobalmeasure22}
Fix $\eta\in\Lambda$ and $t > 0$ such that $\eta_t\in H_\xi$ for some $\xi = g(\bp)\in\Lbp$. If
\begin{equation}
\label{globalmeasure22hypothesis}
\lb \xi|\eta\rb_\zero + \tau \leq t \leq 2\lb \xi|\eta\rb_\zero - t_\xi - \tau,
\end{equation}
then
\begin{equation}
\label{globalmeasure22}
\begin{split}
e^{-\delta (2\lb \xi|\eta\rb_\zero - t_\xi)} \NN_\bp(e^{2\lb \xi|\eta\rb_\zero - t_\xi - t - \sigma})
&\lesssim_\times \mu\big(B(\eta,e^{-t})\big)\\ 
&\lesssim_\times e^{-\delta (2\lb \xi|\eta\rb_\zero - t_\xi)} \NN_\bp(e^{2\lb \xi|\eta\rb_\zero - t_\xi - t + \sigma})
\end{split}
\end{equation}
where $\sigma,\tau > 0$ are independent of $\eta$ and $t$.
\end{corollary}
\begin{proof}
Without loss of generality suppose that $g$ satisfies \eqref{top}, and write $\zeta = g^{-1}(\eta)$. Since $t > t_\xi$, \eqref{ShadcontainsB} guarantees that $B(\eta,e^{-t}) \subset \Shad(g(\zero),\sigma_0)$ for some $\sigma_0 > 0$ independent of $\eta$ and $t$. Then by the Bounded Distortion Lemma \ref{lemmaboundeddistortion}, we have
\[
B(\zeta, e^{-(t - t_\xi)}/C) \subset g^{-1}(B(\eta,e^{-t})) \subset B(\zeta,C e^{-(t - t_\xi)})
\]
for some $C > 0$, and thus
\[
e^{-\delta t_\xi} \mu(B(\zeta, e^{-(t - t_\xi)}/C)) \lesssim_\times  \mu(B(\eta,e^{-t})) \lesssim_\times e^{-\delta t_\xi} \mu(B(\zeta,C e^{-(t - t_\xi)})).
\]
If
\begin{equation}
\label{globalmeasure22hypothesisinternal}
6R_\bp \Dist(\bp,\eta)^2 \leq \frac{e^{-(t - t_\xi)}}{C} \leq C e^{-(t - t_\xi)} \leq \frac{\Dist(\bp,\zeta)}{4},
\end{equation}
then Corollary \ref{corollaryglobalmeasure21} guarantees that
\begin{align*}
e^{-\delta t_\xi} \Dist(\bp,\zeta)^{2\delta} \NN_\bp\left(\frac{e^{-(t - t_\xi)}}{4C\Dist(\bp,\zeta)^2}\right)
&\lesssim_\times \mu(B(\eta,e^{-t}))\\
&\lesssim_\times e^{-\delta t_\xi}\Dist(\bp,\zeta)^{2\delta} \NN_\bp\left(\frac{4C e^{-(t - t_\xi)}}{\Dist(\bp,\zeta)^2}\right).
\end{align*}
On the other hand, since $\xi,\eta\in\Shad(g(\zero),\sigma_0)$, the Bounded Distortion Lemma \ref{lemmaboundeddistortion} guarantees that
\[
\Dist(\bp,\zeta) \asymp_\times e^{t_\xi} \Dist(\xi,\eta) = e^{-(\lb \xi|\eta\rb_\zero - t_\xi)}.
\]
Denoting the implied constant by $K$, we deduce \eqref{globalmeasure22} with $\sigma = \log(4CK^2)$. The proof is completed upon observing that if $\tau = \log(4CK\vee 6R_\bp CK^2)$, then \eqref{globalmeasure22hypothesis} implies \eqref{globalmeasure22hypothesisinternal}.
\end{proof}

\begin{lemma}[Cf. Lemma \ref{lemmasullivanshadow}]
\label{lemmaCCBglobalmeasure}
Fix $\eta\in\Lambda$ and $t > 0$ such that $\eta_t\notin\bigcup(\scrH)$. Then
\[
\mu(B(\eta,e^{-t})) \asymp_\times e^{-\delta t}.
\]
\end{lemma}
\begin{proof}
By \eqref{sigmadef}, there exists $g\in G$ such that $\dist(g(\zero),\eta_t) \asymp_\plus 0$. By \eqref{ShadcontainsB}, we have $B(\eta,e^{-t}) \subset \Shad(g(\zero),\sigma)$ for some $\sigma > 0$ independent of $\eta,t$. It follows that
\[
\mu(B(\eta,e^{-t})) \asymp_\times e^{-\delta t} \mu(g^{-1}(B(\eta,e^{-t}))).
\]
To complete the proof it suffices to show that $\mu(g^{-1}(B(\eta,e^{-t})))$ is bounded from below. By the Bounded Distortion Lemma \ref{lemmaboundeddistortion},
\[
g^{-1}(B(\eta,e^{-t})) \supset B(g^{-1}(\eta),\epsilon)
\]
for some $\epsilon > 0$ independent of $\eta,t$. Now since $G$ is of compact type, we have
\[
\inf_{x\in\Lambda} \mu(B(x,\epsilon)) \geq \min_{x\in S_{\epsilon/2}} \mu(B(x,\epsilon/2)) > 0
\]
where $S_{\epsilon/2}$ is a maximal $\epsilon/2$-separated subset of $\Lambda$. This completes the proof.
\end{proof}

We are now ready to prove Theorem \ref{theoremglobalmeasure}:

\begin{proof}[Proof of Theorem \ref{theoremglobalmeasure}]
Let $\sigma_0$ denote the implied constant of \eqref{betatinterpretation}. Then by \eqref{metat}, for all $\eta\in\Lambda$, $t > 0$, and $\xi\in\Lbp$,
\begin{equation}
\label{metatsimple}
m(\eta,t) = \begin{cases}
e^{-\delta t_\xi} [\II_\bp(e^{t - t_\xi - \theta}) + \mu(\bp)] & t_\xi + \sigma_0 \leq t \leq \lb \xi|\eta\rb_\zero\\
e^{-\delta (2\lb \xi|\eta\rb_\zero - t_\xi)} \NN_\bp(e^{2\lb \xi|\eta\rb_\zero - t - t_\xi - \theta}) & \lb \xi|\eta\rb_\zero < t \leq 2\lb \xi|\eta\rb_\zero - t_\xi - \sigma_0\\
\text{unknown} & \text{otherwise}
\end{cases}
\end{equation}
Applying this formula to Corollaries \ref{corollaryglobalmeasure12} and \ref{corollaryglobalmeasure22} yields the following:

\begin{lemma}
\label{lemmaglobalmeasureinternal}
There exists $\tau\geq \sigma_0$ such that for all $\eta\in\Lambda$ and $t > 0$.
\begin{itemize}
\item[(i)] If for some $\xi$, $t_\xi + \tau \leq t \leq \lb \xi|\eta\rb_\zero - \tau$, then \eqref{globalmeasureformula} holds.
\item[(ii)] If for some $\xi$, $\lb \xi|\eta\rb_\zero + \tau \leq t \leq 2\lb \xi|\eta\rb_\zero - t_\xi - \tau$, then \eqref{globalmeasureformula} holds.
\end{itemize}
\end{lemma}
Now fix $\eta\in\Lambda$, and let
\[
A = \left\{t > 0 : \eta_t\notin\bigcup(\scrH)\right\} \cup \bigcup_{\xi\in\Lbp} [t_\xi + \tau,\lb \xi|\eta\rb_\zero - \tau]\cup \bigcup_{\xi\in\Lbp}  [\lb \xi|\eta\rb_\zero + \tau,2\lb \xi|\eta\rb_\zero - t_\xi - \tau].
\]
Then by Lemmas \ref{lemmaCCBglobalmeasure} and \ref{lemmaglobalmeasureinternal}, $\eqref{globalmeasureformula}_{\sigma = \tau}$ holds for all $t\in A$.
\begin{claim}
Every interval of length $2\tau$ intersects $A$.
\end{claim}
\begin{subproof}
If $[s - \tau,s + \tau]$ does not intersect $A$, then by connectedness, there exists $\xi\in\Lbp$ such that $\eta_t\in H_\xi$ for all $t\in [s - \tau, s + \tau]$. By Lemma \ref{lemmabetatinterpretation}, the fact that $\eta_{s\pm\tau}\in H_\xi$ implies that $t_\xi \leq s \leq 2\lb \xi|\eta\rb_\zero - t_\xi$ (since $\tau\geq \sigma_0$). If $s \leq \lb \xi|\eta\rb_\zero$, then $[s - \tau,s + \tau]\cap [t_\xi + \tau,\lb \xi|\eta\rb_\zero - \tau] \neq \emptyset$, while if $s \geq \lb \xi|\eta\rb_\zero$, then $[s - \tau,s + \tau]\cap [\lb \xi|\eta\rb_\zero + \tau,2\lb \xi|\eta\rb_\zero - t_\xi - \tau] \neq \emptyset$.
\end{subproof}
Thus for all $t > 0$, there exist $t_\pm\in A$ such that $t - 2\tau \leq t_- \leq t \leq t_+ \leq t - 2\tau$; then
\begin{align*}
m(\eta,t + 3\tau)
\lesssim_\times m(\eta,t_+ + \tau)
&\lesssim_\times \mu\big(B(\eta,e^{-t_+})\big)\\
&\leq_\pt \mu\big(B(\eta,e^{-t})\big)\\
&\leq_\pt \mu\big(B(\eta,e^{-t_-})\big)
\lesssim_\times m(\eta,t_- - \tau)
\lesssim_\times m(\eta,t - 3\tau),
\end{align*}
i.e. $\eqref{globalmeasureformula}_{\sigma = 3\tau}$ holds.
\end{proof}

\section{Groups for which $\mu$ is doubling}
\label{subsectiondoubling}

Recall that a measure $\mu$ is said to be \emph{doubling} if for all $\eta\in\Supp(\mu)$ and $r > 0$, $\mu(B(\eta,2r)) \asymp_\times \mu(B(\eta,r))$. In the Standard Case, the Global Measure Formula implies that the $\delta$-conformal measure of a geometrically finite group is always doubling (Example \ref{exampledoubling}). However, in general there are geometrically finite groups whose $\delta$-conformal measures are not doubling (Example \ref{exampledoublingnondoubling}). It is therefore of interest to determine necessary and sufficient conditions on a geometrically finite group for its $\delta$-conformal measure to be doubling. The Global Measure Formula immediately yields the following criterion:

\begin{lemma}
$\mu$ is doubling if and only if the function $m$ satisfies
\begin{equation}
\label{mdoubling}
m(\eta,t + \tau) \asymp_{\times,\tau} m(\eta,t - \tau) \all \eta\in\Lambda \all t,\tau > 0.
\end{equation}
\end{lemma}
\begin{proof}
If \eqref{mdoubling} holds, then \eqref{globalmeasureformula} reduces to
\begin{equation}
\label{globalmeasureformula2}
\mu(B(\eta,e^{-t})) \asymp_\times m(\eta,t),
\end{equation}
and then \eqref{mdoubling} shows that $\mu$ is doubling. On the other hand, if $\mu$ is doubling, then \eqref{globalmeasureformula} implies that
\[
m(\eta,t - \tau) \lesssim_\times \mu(B(\eta,e^{-(t - \tau - \sigma)})) \asymp_\times \mu(B(\eta,e^{-(t + \tau + \sigma)})) \lesssim_\times m(\eta,t + \tau);
\]
combining with Proposition \ref{propositionmetadecreasing} shows that \eqref{mdoubling} holds.
\end{proof}

Of course, the criterion \eqref{mdoubling} is not very useful by itself, since it refers to the complicated function $m$. In what follows we find more elementary necessary and sufficient conditions for doubling. First we must introduce some terminology.

\begin{definition}
\label{definitiondoubling}
A function $f:\CO 1\infty\to\CO 1\infty$ is called \emph{doubling} if there exists $\beta > 1$ such that
\begin{equation}
\label{dexp1}
f(\beta R) \lesssim_{\times,\beta} f(R) \all R\geq 1,
\end{equation}
and \emph{codoubling} if there exists $\beta > 1$ such that
\begin{equation}
\label{dexp2}
f(\beta R) - f(R) \gtrsim_{\times,\beta} f(R) \all R\geq 1.
\end{equation}
\end{definition}

\begin{observation}
\label{observationcodoubling}
If there exists $\beta > 1$ such that
\[
\NN_\bp(\beta R) > \NN_\bp(R) \all R\geq 1,
\]
then $\NN_\bp$ is codoubling.
\end{observation}
\begin{proof}
Fix $R\geq 1$; there exists $h\in G_\bp$ such that $2R < \|h\|_\bp \leq 2\beta R$. We have
\[
h\{j\in G_\bp : \|j\|_\bp \leq R\} \subset \{j\in G_\bp : R < \|j\|_\bp \leq (2\beta + 1) R\},
\]
and taking cardinalities gives
\[
\NN_\bp(R) \leq \NN_\bp\big((2\beta + 1) R\big) - \NN_\bp(R).
\]
\end{proof}

We are now ready to state a more elementary characterization of when $\mu$ is doubling:

%

\begin{proposition}
\label{propositiondoubling}
$\mu$ is doubling if and only if all of the following hold:
\begin{itemize}
\item[(I)] For all $\bp\in P$, $\NN_\bp$ is both doubling and codoubling.
\item[(II)] For all $\bp\in P$ and $R\geq 1$,
\begin{equation}
\label{INasymp2}
\II_\bp(R) \asymp_\times R^{-2\delta} \NN_\bp(R).
\end{equation}
\item[(III)] $G$ is of divergence type.
\end{itemize}
Moreover, \text{(II)} can be replaced by
\begin{itemize}
\item[(II$'$)] For all $\bp\in P$ and $R\geq 1$,
\begin{equation}
\label{INasymp3}
\w\II_\bp(R) := \sum_{k = 0}^\infty e^{-2\delta k} \NN_\bp(e^k R) \asymp_\times \NN_\bp(R).
\end{equation}
\end{itemize}
\end{proposition}

\begin{proof}[Proof that \text{(I)-(III)} imply $\mu$ doubling]
Fix $\eta\in\Lambda$ and $t,\tau > 0$, and we will demonstrate \eqref{mdoubling}. By (II), (III), and Observation \ref{observationnonatomic}, we have
\begin{equation}
\label{GMFdoubling1}
\begin{split}
m(\eta,t) &\asymp_\times \begin{cases}
e^{-\delta t} & \eta_t\notin\bigcup(\scrH) \\
e^{-\delta t_\xi} e^{-2\delta(t - t_\xi - \theta)}\NN_\bp(e^{t - t_\xi - \theta}) & \eta_t\in H_\xi \text{ and }t \leq \lb \xi|\eta\rb_\zero\\
e^{-\delta (2\lb \xi|\eta\rb_\zero - t_\xi)} \NN_\bp(e^{2\lb \xi|\eta\rb_\zero - t - t_\xi - \theta}) & \eta_t\in H_\xi \text{ and } t > \lb \xi|\eta\rb_\zero
\end{cases}\\
&\asymp_\times e^{-\delta t} \begin{cases}
1 & \eta_t\notin\bigcup(\scrH)\\
e^{-\delta b(\eta,t)} \NN_\bp(e^{b(\eta,t) - \theta}) & \eta_t \in H_{g(\bp)}
\end{cases}
\end{split}
\end{equation}
where $b(\eta,t)$ is as in \eqref{betat}. Let $t_\pm = t \pm \tau$. We split into two cases:
\begin{itemize}
\item[Case 1:] $\eta_{t_+},\eta_{t_-}\in H_{g(\bp)}$ for some $g(\bp)\in \Lbp$. In this case, \eqref{mdoubling} follows from Corollary \ref{corollarybetat} together with the fact that $\NN_\bp$ is doubling.
\item[Case 2:] $\eta_{t + s}\notin \bigcup(\scrH)$ for some $s\in [-\tau,\tau]$. In this case, Corollary \ref{corollarybetat} shows that $b(\eta,t_\pm) \asymp_{\plus,\tau} 0$ and thus
\[
m(\eta,t_+) \asymp_{\times,\tau} e^{-\delta t} \asymp_{\times,\tau} m(\eta,t_-).
\]
\end{itemize}
\end{proof}
\noindent Before continuing the proof of Proposition \ref{propositiondoubling}, we observe that
\begin{equation*}
\begin{split}
\II_\bp(R) + R^{-2\delta}\NN_\bp(R) \asymp_\times \sum_{h\in G_\bp} (R\vee \|h\|_\bp)^{-2\delta}
&\asymp_\times \sum_{h\in G_\bp} \sum_{k = 1}^\infty (e^k R)^{-2\delta} [e^k R \geq \|h\|_\bp]\\
&=_\pt \sum_{k = 1}^\infty (e^k R)^{-2\delta} \NN_\bp(e^k R)\\
&= R^{-2\delta} \w\II_\bp(R).
\end{split}
\end{equation*}
In particular, it follows that \eqref{INasymp3} is equivalent to
\begin{equation}
\label{INasymp4}
\II_\bp(R) \lesssim_\times R^{-2\delta} \NN_\bp(R).
\end{equation}
\begin{proof}[Proof that \text{(I) and (II$'$)} imply \text{(II)}]
Since $\NN_\bp$ is codoubling, let $\beta > 1$ be as in \eqref{dexp2}. Then
\[
\II_\bp(R) \geq \sum_{\substack{h\in G_\bp \\ R < \|h\|_\bp \leq \beta R}} (\beta R)^{-2\delta} = (\beta R)^{-2\delta} (\NN_\bp(\beta R) - \NN_\bp(R)) \gtrsim_{\times,\beta} R^{-2\delta} \NN_\bp(R).
\]
Combining with \eqref{INasymp4} completes the proof.
\end{proof}
\begin{proof}[Proof that $\mu$ doubling implies \text{(I)-(III)} and \text{(II$'$)}]
Since a doubling measure whose topological support is a perfect set cannot have an atomic part, we must have $\mu(P) = 0$ and thus by Observation \ref{observationnonatomic}, (III) holds. Since
\[
m(\bp,t) \asymp_{\times,\bp} \II_\bp(e^{t - t_0 - \theta}) + \mu(\bp) = \II_\bp(e^{t - t_0 - \theta})
\]
for all sufficiently large $t$, setting $\eta = \bp$ in \eqref{mdoubling} shows that the function $\II_\bp$ is doubling.

Fix $\eta\in\Lambda\butnot\{\bp\}$. Let $\sigma_0 > 0$ denote the implied constant of \eqref{betatinterpretation}. For $s\in [t_0 + \sigma_0 + \tau, \lb \bp|\eta\rb_\zero - \tau]$, plugging $t = 2\lb \bp|\eta\rb_\zero - s$ into \eqref{mdoubling} and simplifying using $\eqref{metatsimple}_{\xi = \bp}$ shows that
\begin{equation}
\label{Nbps}
\NN_\bp(e^{s - \tau - t_0 - \theta}) \asymp_{\times,\tau} \NN_\bp(e^{s + \tau - t_0 - \theta}).
\end{equation}
Since $\lb \bp|\eta\rb_\zero$ can be made arbitrarily large, \eqref{Nbps} holds for all $s\geq t_0 + \sigma_0 + \tau$. It follows that $\NN_\bp$ is doubling.

Next, we compare the values of $m(\eta,\lb \bp|\eta\rb_\zero \pm\tau)$. This gives (assuming $\lb \bp|\eta\rb_\zero > t_0 + \sigma_0 + \tau$)
\[
e^{-\delta t_0} \II_\bp(e^{\lb \bp|\eta\rb_\zero - \tau - t_0 - \theta}) \asymp_\times e^{-\delta (2\lb \bp|\eta\rb_\zero - t_0)} \NN_\bp(e^{\lb \bp|\eta\rb_\zero - \tau - t_0 - \theta}).
\]
Letting $R_\eta = \exp(\lb \bp|\eta\rb_\zero - \tau - t_0 - \theta)$, we have
\begin{equation}
\label{INasymp6}
\II_\bp(R_\eta) \asymp_\times R_\eta^{-2\delta} \NN_\bp(R_\eta).
\end{equation}
Now fix $\zeta\in\Lambda\butnot\{\bp\}$ and $h\in G_\bp$, and let $\eta = h(\zeta)$. Then $\Dist_\bp(h(\zero),\eta) \asymp_{\plus,\zeta} 0$, and thus the triangle inequality gives
\[
1 \leq \Dist_\bp(\zero,\eta) \asymp_{\plus,\zeta} \|h\|_\bp \geq 1,
\]
and so $R_\eta \asymp_\times \Dist_\bp(\zero,\eta) \asymp_{\times,\zeta} \|h\|_\bp$. Combining with \eqref{INasymp6} and the fact that the functions $\II_\bp$ and $\NN_\bp$ are doubling, we have
\begin{equation}
\label{INasymp7}
\II_\bp(\|h\|_\bp) \asymp_\times \|h\|_\bp^{-2\delta} \NN_\bp(\|h\|_\bp)
\end{equation}
for all $h\in G_\bp$.

Now fix $1\leq R_1 < R_2$ such that $\|h_i\|_\bp = R_i$ for some $h_1,h_2\in G_\bp$, but such that the formula $R_1 < \|h\|_\bp < R_2$ is not satisfied for any $h\in G_\bp$. Then
\[
\lim_{R\searrow R_1} \II_\bp(R) = \lim_{R\nearrow R_2} \II_\bp(R) \text{ and } \lim_{R\searrow R_1} \NN_\bp(R) = \lim_{R\nearrow R_2} \NN_\bp(R).
\]
On the other hand, applying \eqref{INasymp7} with $h = h_1,h_2$ gives
\[
\II_\bp(R_i) \asymp_\times R_i^{-2\delta} \NN_\bp(R_i).
\]
Since $\II_\bp$ and $\NN_\bp$ are doubling, we have
\begin{align*}
R_1^{-2\delta}
\asymp_\times \frac{\II_\bp(R_1)}{\NN_\bp(R_1)}
\asymp_\times \frac{\lim_{R\searrow R_1} \II_\bp(R)}{\lim_{R\searrow R_1} \NN_\bp(R)}
&= \frac{\lim_{R\nearrow R_2} \II_\bp(R)}{\lim_{R\nearrow R_2} \NN_\bp(R)}\\
&\asymp_\times \frac{\II_\bp(R_2)}{\NN_\bp(R_2)}
\asymp_\times R_2^{-2\delta}
\end{align*}
and thus $R_1 \asymp_\times R_2$. Since $R_1,R_2$ were arbitrary, Observation \ref{observationcodoubling} shows that $\NN_\bp$ is codoubling. This completes the proof of (I).

It remains to demonstrate (II) and (II$'$). Given any $R\geq 1$, since $\NN_\bp$ is codoubling, we may find $h\in G_\bp$ such that $\|h\|_\bp \asymp_\times R$; combining with \eqref{INasymp7} and the fact that $\II_\bp$ and $\NN_\bp$ are doubling gives \eqref{INasymp2} and \eqref{INasymp4}, demonstrating (II) and (II$'$).
\end{proof}

We note that the proof actually shows the following (cf. \eqref{GMFdoubling1}):

\begin{corollary}
\label{corollaryGMFdoubling}
If $\mu$ is doubling, then
\[
\mu(B(\eta,e^{-t})) \asymp_\times e^{-\delta t} \begin{cases}
1 & \eta_t\notin \bigcup(\scrH)\\
e^{-\delta b(\eta,t)} \NN_\bp(e^{b(\eta,t)}) & \eta_t\in H_{g(\bp)}
\end{cases}
\]
for all $\eta\in\Lambda$, $t > 0$. Here $b(\eta,t)$ is as in \eqref{betat}.
\end{corollary}

Although Proposition \ref{propositiondoubling} is the best necessary and sufficient condition we can give for doubling, in what follows we give necessary conditions and sufficient conditions which are more elementary (Proposition \ref{propositiondoubling2}), although the necessary conditions are not the same as the sufficient conditions. In practice these conditions are usually powerful enough to determine whether any given measure is doubling.

To state the result, we need the concept of the \emph{polynomial growth rate} of a function:

\begin{definition}[Cf. \eqref{growthrate}]
\label{definitiongrowthrate}
The \emph{(polynomial) growth rate} of a function $f:\CO 1\infty \to \CO 1\infty$ is the limit
\[
\dexp(f) := \lim_{\lambda,R\to\infty}\frac{\log f(\lambda R) - \log f(R)}{\log(\lambda)}
\]
if it exists. If the limit does not exist, then the numbers
\begin{align*}
\uexp(f) &:= \limsup_{\lambda,R\to\infty}\frac{\log f(\lambda R) - \log f(R)}{\log(\lambda)}\\
\lexp(f) &:= \liminf_{\lambda,R\to\infty}\frac{\log f(\lambda R) - \log f(R)}{\log(\lambda)}
\end{align*}
are the \emph{upper} and \emph{lower polynomial growth rates} of $f$, respectively.
\end{definition}

%

\begin{lemma}
\label{lemmadexp}
Let $f:\CO 1\infty\to\CO 1\infty$.
\begin{itemize}
\item[(i)] $f$ is doubling if and only if $\uexp(f) < \infty$.
\item[(ii)] $f$ is codoubling if and only if $\lexp(f) > 0$.
\item[(iii)]
\[
\lexp(f) \leq \liminf_{\lambda\to\infty} \frac{\log f(\lambda)}{\log(\lambda)} \leq \limsup_{\lambda\to\infty} \frac{\log f(\lambda)}{\log(\lambda)} \leq \uexp(f).
\]
In particular, $\lexp(\NN_\bp) \leq 2\delta_\bp \leq \uexp(\NN_\bp)$.
\end{itemize}
\end{lemma}
\begin{proof}[Proof of \text{(i)}]
Suppose that $f$ is doubling, and let $C > 1$ denote the implied constant of \eqref{dexp1}. Iterating gives
\[
f(\beta^n R) \leq C^n f(R) \all n\in\N \all R\geq 1
\]
and thus
\[
f(\lambda R) \lesssim_\times \lambda^{\log_\beta(C)} f(R) \all \lambda > 1 \all R\geq 1.
\]
It follows that $\uexp(f) \leq \log_\beta(C) < \infty$. The converse direction is trivial.
\end{proof}
\begin{proof}[Proof of \text{(ii)}]
Suppose that $f$ is codoubling, and let $C > 1$ denote the implied constant of \eqref{dexp2}. Then
\[
f(\beta R) \geq (1 + C^{-1}) f(R) \all R\geq 1.
\]
Iterating gives
\[
f(\beta^n R) \geq (1 + C^{-1})^n f(R) \all n\in\N \all R\geq 1
\]
and thus
\[
f(\lambda R) \gtrsim_\times \lambda^{\log_\beta(1 + C^{-1})} f(R) \all \lambda > 1 \all R\geq 1.
\]
It follows that $\lexp(f) \geq \log_\beta(1 + C^{-1}) > 0$. The converse direction is trivial.
\end{proof}
\begin{proof}[Proof of \text{(iii)}]
Let $R_n\to\infty$. For each $n\in\N$,
\[
\limsup_{\lambda\to\infty} \frac{\log f(\lambda R_n) - \log f(R_n)}{\log(\lambda)} = \overline s := \limsup_{\lambda\to\infty} \frac{\log f(\lambda)}{\log(\lambda)}\cdot
\]
Thus given $s < \overline s$, we may find a large number $\lambda_n > 1$ such that 
\[
\frac{\log f(\lambda_n R_n) - \log f(R_n)}{\log(\lambda_n)}  > s.
\] 
Since $\lambda_n,R_n\to\infty$ as $n\to\infty$, it follows that $\uexp(f) \geq s$; since $s$ was arbitrary, $\uexp(f) \geq \overline s$. A similar argument shows that $\lexp(f) \leq \underline s$.

Finally, when $f = \NN_\bp$, the equality $\overline s = \underline s = 2\delta_\bp$ is a consequence of \eqref{poincarealternate} and Observation \ref{observationeuclideanparabolicasymp}.
\end{proof}

We can now state our final result regarding criteria for doubling:

\begin{proposition}
\label{propositiondoubling2}
In the following list, \text{(A) \implies (B) \implies (C):}
\begin{itemize}
\item[(A)] For all $\bp\in P$, $0 < \lexp(\NN_\bp) \leq \uexp(\NN_\bp) < 2\delta$.
\item[(B)] $\mu$ is doubling.
\item[(C)] For all $\bp\in P$, $0 < \lexp(\NN_\bp) \leq \uexp(\NN_\bp) \leq 2\delta$.
\end{itemize}
\end{proposition}
\begin{proof}[Proof of \text{(A) \implies (B)}]
Suppose that (A) holds. Then by Lemma \ref{lemmadexp}, (I) of Proposition \ref{propositiondoubling} holds. Since $\delta_\bp \leq \uexp(\NN_\bp)/2 < \delta$ for all $\bp\in P$, Theorem \ref{theoremGFdivergencetype} implies that (III) of Proposition \ref{propositiondoubling} holds. To complete the proof, we need to show that (II$'$) of Proposition \ref{propositiondoubling} holds. Fix $s\in (\uexp(\NN_\bp),2\delta)$. Since $s > \uexp(\NN_\bp)$, we have
\[
\NN_\bp(\lambda R) \lesssim_{\times,s} \lambda^s \NN_\bp(R) \all \lambda > 1,\; R\geq 1
\]
and thus
\[
\NN_\bp(R) \leq \w\II_\bp(R) \lesssim_\times \sum_{k = 0}^\infty e^{-2\delta k} e^{sk} \NN_\bp(R) \asymp_\times \NN_\bp(R),
\]
demonstrating \eqref{INasymp3} and completing the proof.
\end{proof}
\begin{proof}[Proof of \text{(B) \implies (C)}]
Suppose $\mu$ is doubling. By (I) of Proposition \ref{propositiondoubling}, $\lexp(\NN_\bp) > 0$. On the other hand, by \eqref{INasymp3} we have
\[
\lambda^{-2\delta} \NN_\bp(\lambda R) \lesssim_\times \NN_\bp(R) \all \lambda > 1,\; R\geq 1
\]
and thus $\uexp(\NN_\bp) \leq 2\delta$.
\end{proof}

Proposition \ref{propositiondoubling} shows that if $G$ is a geometrically finite group with $\delta$-conformal measure $\mu$, then the question of whether $\mu$ is doubling is determined entirely by its parabolic subgroups $(G_\bp)_{\bp\in P}$ and its Poincar\'e set $\Delta_G$. A natural question is when the second input can be removed, that is: if we are told what the parabolic subgroups $(G_\bp)_{\bp\in P}$ are, can we sometimes determine whether $\mu$ is doubling without looking at $\Delta_G$? A trivial example is that if $\lexp(\NN_\bp) = 0$ or $\uexp(\NN_\bp) = \infty$ for some $\bp\in P$, then we automatically know that $\mu$ is not doubling. Conversely, the following definition and proposition describe when we can deduce that $\mu$ is doubling:

\begin{definition}
\label{definitionpredoubling}
A parabolic group $H\leq\Isom(X)$ with global fixed point $\bp\in\del X$ is \emph{pre-doubling} if
\begin{equation}
\label{predoubling}
0 < \lexp(\NN_{\EE_\bp,H}) \leq \uexp(\NN_{\EE_\bp,H}) = 2\delta_H < \infty
\end{equation}
and $H$ is of divergence type.
\end{definition}


\begin{proposition}
\label{propositionpredoubling}
~
\begin{itemize}
\item[(i)] If $G_\bp$ is pre-doubling for every $\bp\in P$, then $\mu$ is doubling.
\item[(ii)] Let $H\leq\Isom(X)$ be a parabolic subgroup, and let $g\in\Isom(X)$ be a loxodromic isometry such that $\lb g,H\rb$ is a strongly separated Schottky product. Then the following are equivalent:
\begin{itemize}
\item[(A)] $H$ is pre-doubling.
\item[(B)] For every $n\in\N$, the $\delta_n$-quasiconformal measure $\mu_n$ of $G_n = \lb g^n,H\rb$ is doubling. Here we assume that $\delta_n := \delta(G_n) < \infty$.
\end{itemize}
\end{itemize}
\end{proposition}
\begin{proof}[Proof of \text{(i)}]
For all $\bp\in P$, the fact that $G_\bp$ is of divergence type implies that $\delta > \delta_\bp$ (Proposition \ref{propositiondeltagtrdeltap}); combining with \eqref{predoubling} gives $0 < \lexp(\NN_\bp) \leq \uexp(\NN_\bp) < 2\delta$. Proposition \ref{propositiondoubling2} completes the proof.
\end{proof}

\begin{proof}[Proof of \text{(ii)}]
Since (up to equivalence) the only parabolic point of $G_n$ is the global fixed point of $H$ (Proposition \ref{propositionSproductGF}), the implication (A) \implies (B) follows from part (i). Conversely, suppose that (B) holds. Then by Proposition \ref{propositiondoubling2}, we have
\[
0 < \lexp(\NN_{\EE_\bp,H}) \leq \uexp(\NN_{\EE_\bp,H}) \leq 2\delta_n < \infty.
\]
Since $\delta_n\to \delta_H$ as $n\to\infty$ (Proposition \ref{propositionSproductexponent}(iv)), taking the limit and combining with the inequality $2\delta_H \leq \uexp(\NN_{\EE_\bp,H})$ yields \eqref{predoubling}. On the other hand, by Proposition \ref{propositiondoubling}, for each $n$, $G_n$ is of divergence type, so applying Proposition \ref{propositionSproductexponent}(iv) again, we see that $H$ is of divergence type.
%
\end{proof}

\begin{example}
\label{exampledoubling}
If
\begin{equation}
\label{Nbppowerlaw}
\NN_\bp(R) \asymp_\times R^{2\delta_\bp} \all \bp\in P,
\end{equation}
then the groups $(G_\bp)_{\bp\in P}$ are pre-doubling, and thus by Proposition \ref{propositionpredoubling}(i), $\mu$ is doubling. Combining with Corollary \ref{corollaryGMFdoubling} gives
\[
\mu(B(\eta,e^{-t})) \asymp_\times e^{-\delta t} \begin{cases}
1 & \eta_t\notin \bigcup(\scrH)\\
e^{(2\delta_\xi - \delta) b(\eta,t)} & \eta_t\in H_\xi
\end{cases}.
\]
This generalizes B. Schapira's global measure formula \cite[Th\'eor\`eme 3.2]{Schapira} to the setting of regularly geodesic strongly hyperbolic metric spaces.
\end{example}

We remark that the asymptotic \eqref{Nbppowerlaw} is satisfied whenever $X$ is a finite-dimensional algebraic hyperbolic space; see e.g. \cite[Lemma 3.5]{Newberger}. In particular, specializing Schapira's global measure formula to the settings of finite-dimensional algebraic hyperbolic spacess and finite-dimensional real hyperbolic spaces give the global measure formulas of Newberger \cite[Main Theorem]{Newberger} and Stratmann--Velani--Sullivan \cite[Theorem 2]{StratmannVelani}, \cite[Theorem on p.271]{Sullivan_entropy}, respectively.

By contrast, when $X = \H = \H^\infty$, the asymptotic \eqref{Nbppowerlaw} is usually not satisfied. Let us summarize the various behaviors that we have seen for the orbital counting functions of parabolic groups acting on $\H^\infty$, and their implications for doubling:

\begin{examples}[Examples of doubling and non-doubling Patterson--Sullivan measures of geometrically finite subgroups of $\Isom(\H^\infty)$]
\label{exampledoublingnondoubling}
~
\begin{itemize}
\item[1.] In the proof of Theorem \ref{theoremnilpotentembedding} (cf. Remark \ref{remarknilpotentembedding}), we saw that if $\Gamma$ is a finitely generated virtually nilpotent group and if $f:\CO 1\infty\to\CO 1\infty$ is a function satisfying
\[
\alpha_\Gamma < \lexp(f) \leq \uexp(f) < \infty,
\]
then there exists a parabolic group $H\leq\Isom(\H^\infty)$ isomorphic to $\Gamma$ whose orbital counting function is asymptotic to $f$. Now, a group $H$ constructed in this way may or may not be pre-doubling; it depends on the chosen function $f$. We note that by applying Proposition \ref{propositionpredoubling}(ii) to such a group, one can construct examples of geometrically finite subgroups of $\Isom(\H^\infty)$ whose Patterson--Sullivan measures are not doubling. On the other hand, for any parabolic group $H$ constructed in this way, if $H$ is embedded into a geometrically finite group $G$ with sufficiently large Poincar\'e exponent (namely $2\delta_G > \uexp(f)$), then the Patterson--Sullivan measure of $G$ may be doubling (assuming that no other parabolic subgroups of $G$ are causing problems).
\item[2.] In Theorem \ref{theoremRtreeorbitalcounting}, we showed that if $f:\Rplus\to\N$ satisfies the condition
\[
\forall 0\leq R_1 \leq R_2 \; \text{$f(R_1)$ divides $f(R_2)$},
\]
then there exists a parabolic subgroup of $\Isom(\H^\infty)$ whose orbital counting function is equal to $f$. This provides even more examples of parabolic groups which are not pre-doubling. In particular, it provides examples of parabolic groups $H$ which satisfy either $\lexp(\NN_H) = 0$ or $\uexp(\NN_H) = \infty$ (cf. Example \ref{exampleparabolictorsion}); such groups cannot be embedded into any geometrically finite group with a doubling Patterson--Sullivan measure.
\end{itemize}
\end{examples}

Note that example 2 can be used to construct a geometrically finite group acting isometrically on an $\R$-tree which does not have a doubling Patterson--Sullivan measure. On the other hand, example 1 has no analogue in $\R$-trees by Remark \ref{remarkparabolicRtree}.




\section{Exact dimensionality of $\mu$}
\label{subsectionexactdimensional}
We now turn to the question of the fractal dimensions of the measure $\mu$. We recall that the \emph{Hausdorff dimension} and \emph{packing dimension} of a measure $\mu$ on $\del X$ are defined by the formulas
\begin{align*}
\HD(\mu) &= \inf\left\{\HD(A): \mu(\del X\butnot A) = 0\right\}\\
\PD(\mu) &= \inf\left\{\PD(A): \mu(\del X\butnot A) = 0\right\}.
\end{align*}
If $G$ is of convergence type, then $\mu$ is atomic, so $\HD(\mu) = \PD(\mu) = 0$. Consequently, for the remainder of this chapter we make the
\begin{standingassumption}
\label{assumptiondivergencetype}
$G$ is of divergence type.
\end{standingassumption}
Given this assumption, it is natural to expect that $\HD(\mu) = \PD(\mu) = \delta$. Indeed, the inequality $\HD(\mu)\leq\delta$ follows immediately from Theorems \ref{theorembishopjonesregular} and \ref{theoremGFcompact}, and in the Standard Case equality holds \cite[Proposiiton 4.10]{StratmannVelani}. Even stronger than the equalities $\HD(\mu) = \PD(\mu) = \delta$, it is natural to expect that $\mu$ is \emph{exact dimensional}:

\begin{definition}
\label{definitionexactdimensional}
A measure $\mu$ on a metric space $(Z,\Dist)$ is called \emph{exact dimensional of dimension $s$} if the limit
\begin{equation}
\label{dmudef}
d_\mu(\eta) := \lim_{t\to\infty} \frac1t \log\frac{1}{\mu(B(\eta,e^{-t}))}
\end{equation}
exists and equals $s$ for $\mu$-a.e. $\eta\in Z$.
\end{definition}

For example, every Ahlfors $s$-regular measure is exact dimensional of dimension $s$.

If the limit in \eqref{dmudef} does not exist, then we denote the lim inf by $\underline d_\mu(\eta)$ and the lim sup by $\overline d_\mu(\eta)$.

\begin{proposition}[{\cite[\68]{MSU}}]
\label{propositionexactdimensional}
For any measure $\mu$ on a metric space $(Z,\Dist)$,
\begin{align*}
\HD(\mu) &= \esssup_{\eta\in Z} \underline d_\mu(\eta)\\
\PD(\mu) &= \esssup_{\eta\in Z} \overline d_\mu(\eta).
\end{align*}
In particular, if $\mu$ is exact dimensional of dimension $s$, then
\[
\HD(\mu) = \PD(\mu) = s.
\]
\end{proposition}

Combining Proposition \ref{propositionexactdimensional} with Lemma \ref{lemmaCCBglobalmeasure} and Observation \ref{observationnonatomic} immediately yields the following:

\begin{observation}
\label{observationHDdeltaPD}
If $\mu$ is the Patterson--Sullivan measure of a geometrically finite group of divergence type, then
\[
\HD(\mu) \leq \delta \leq \PD(\mu).
\]
In particular, if $\mu$ is exact dimensional, then $\mu$ is exact dimensional of dimension $\delta$.
\end{observation}

It turns out that $\mu$ is not necessarily exact dimensional (Example \ref{examplenotexactdim}), but counterexamples to exact dimensionality must fall within a very narrow window (Theorem \ref{theoremhlogseries}), and in particular if $\mu$ is doubling then $\mu$ is exact dimensional (Corollary \ref{corollarydoublingexactdim}). As a first step towards these results, we will show that exact dimensionality is equivalent to a certain Diophantine condition. For this, we need to recall some results from \cite{FSU4}.

\subsection{Diophantine approximation on $\Lambda$}
\label{subsubsectiondiophantineapproximation}
Classically, Diophantine approximation is concerned with the approximation of a point $x\in\R\butnot\Q$ by a rational number $p/q\in\Q$. The two important quantities are the \emph{error term} $|x - p/q|$ and the \emph{height} $q$. Given a function $\Psi:\N\to\Rplus$, the point $x\in\R\butnot\Q$ is said to be \emph{$\Psi$-approximable} if
\[
\left|x - \frac pq\right| \leq \Psi(q) \text{ for infinitely many $p/q\in\Q$}.
\]
In the setting of a group acting on a hyperbolic metric space, we can instead talk about \emph{dynamical} Diophantine approximation, which is concerned with the approximation of a point $\eta\in\Lambda$ by points $g(\xi)\in G(\xi)$, where $\xi\in\Lambda$ is a distinguished point. For this to make sense, one needs a new definition of error and height: the error term is defined to be $\Dist(g(\xi),\eta)$, and the height is defined to be $b^{\dogo g}$. (If there is more than one possibility for $g$, it may be chosen so as to minimize the height.) Some motivation for these definitions comes from considering classical Diophantine approximation as a special case of dynamical Diophantine approximation which occurs when $X = \H^2$ and $G = \SL_2(\Z)$; see e.g. \cite[Observation 1.15]{FSU4} for more details. Given a function $\Phi:\Rplus\to (0,\infty)$, the point $\eta\in\Lambda$ is said to be \emph{$\Phi,\xi$-well approximable} if for every $K > 0$ there exists $g\in G$ such that
\[
\Dist(g(\xi),\eta) \leq \Phi(K b^{\dogo g}) \text{ for infinitely many $g\in G$}
\]
(cf. \cite[Definition 1.36]{FSU4}). Moreover, $\eta$ is said to be \emph{$\xi$-very well approximable} if
\[
\omega_\xi(\eta) := \limsup_{\substack{g\in G \\ g(\xi)\to\eta}} \frac{-\log_b \Dist(g(\xi),\eta)}{\dogo g} > 1
\]
(cf. \cite[p.9]{FSU4}). The set of $\Phi,\xi$-well approximable points is denoted $\WA_{\Phi,\xi}$, while the set of $\xi$-very well approximable points is denoted $\VWA_\xi$. Finally, a point $\eta$ is said to be \emph{Liouville} if $\omega_\xi(\eta) = \infty$; the set of Liouville points is denoted $\Liou_\xi$.

In the following theorems, we return to the setting of Standing Assumptions \ref{standingassumptionsGFmeasures} and \ref{assumptiondivergencetype}.

\begin{theorem}[Corollary of {\cite[Theorem 8.1]{FSU4}}]
\label{theoremFSUkhinchin}
Fix $\bp\in P$, and let $\Phi:\Rplus\to(0,\infty)$ be a function such that the function $t\mapsto t\Phi(t)$ is nonincreasing. Then
\begin{itemize}
\item[(i)] $\mu(\WA_{\Phi,\bp}) = 0\text{ or }1$ according to whether the series
\begin{equation}
\label{WAseries}
\sum_{g\in G} e^{-\delta \dogo g} \II_\bp\left(\frac{1}{e^{\dogo g} \Phi(K e^{\dogo g})}\right)
\end{equation}
converges for some $K > 0$ or diverges for all $K > 0$, respectively.
\item[(ii)] $\mu(\VWA_\bp) = 0\text{ or }1$ according to whether the series
\begin{equation}
\label{VWAseries}
\SigmaVWA(\bp,\kappa) := \sum_{g\in G} e^{-\delta \dogo g} \II_\bp(e^{\kappa\dogo g})
\end{equation}
converges for all $\kappa > 0$ or diverges for some $\kappa > 0$, respectively.
\item[(iii)] $\mu(\Liou_\bp) = 0\text{ or }1$ according to whether the series $\SigmaVWA(\bp,\kappa)$ converges for some $\kappa > 0$ or diverges for all $\kappa > 0$, respectively.
\end{itemize}
\end{theorem}
\begin{proof}
Standing Assumption \ref{assumptiondivergencetype}, Theorem \ref{theoremahlforsthurstongeneral}, and Observation \ref{observationnonatomic} imply that $\mu$ is ergodic and that $\mu(\bp) = 0$, thus verifying the hypotheses of \cite[Theorem 8.1]{FSU4}. Theorem \ref{theoremglobalmeasure} shows that
\[
\II_\bp(C_1/r) \lesssim_{\times,\bp} \mu(B(\bp,r)) \lesssim_{\times,\bp} \II_\bp(C_2/r)
\]
for some constants $C_1 \geq 1 \geq C_2 > 0$. Thus for all $K > 0$,
\begin{align*}
&\sum_{g\in G} e^{-\delta\dogo g} \II_\bp\left(\frac{1}{e^{\dogo g} \Phi(K C_1 e^{\dogo g})}\right)\\
&\leq_\pt \sum_{g\in G} e^{-\delta\dogo g} \II_\bp\left(\frac{C_1}{e^{\dogo g} \Phi(K e^{\dogo g})}\right)\\
&\lesssim_\times \text{\cite[(8.1)]{FSU4}}\\
&\lesssim_\times \sum_{g\in G} e^{-\delta\dogo g} \II_\bp\left(\frac{C_2}{e^{\dogo g} \Phi(K e^{\dogo g})}\right)\\
&\leq_\pt \sum_{g\in G} e^{-\delta\dogo g} \II_\bp\left(\frac{1}{e^{\dogo g} \Phi((K/C_1) e^{\dogo g})}\right).
\end{align*}
Thus, \cite[(8.1)]{FSU4} diverges for all $K > 0$ if and only if \eqref{WAseries} diverges for all $K > 0$. This completes the proof of (i). To demonstrate (ii) and (iii), simply note that $\VWA_\bp = \bigcup_{c > 0} \WA_{\Phi_c,\bp}$ and $\Liou_\bp = \bigcap_{c > 0} \WA_{\Phi_c,\bp}$, where $\Phi_c(t) = t^{-(1 + c)}$, and apply (i). The constant $K$ may be absorbed by a slight change of $\kappa$.
\end{proof}

\begin{theorem}[Corollary of {\cite[Theorem 7.1]{FSU4}}]
\label{theoremFSUjarnik}
For all $\xi\in\Lambda$ and $c > 0$,
\[
\HD(\WA_{\Phi_c,\xi}) \leq \frac{\delta}{1 + c},
\]
where $\Phi_c(t) = t^{-(1 + c)}$ as above. In particular, $\HD(\Liou_\xi) = 0$, and $\VWA_\xi$ can be written as the countable union of sets of Hausdorff dimension strictly less than $\delta$.
\end{theorem}
\noindent (No proof is needed as this follows directly from \cite[Theorem 7.1]{FSU4}.)

There is a relation between dynamical Diophantine approximation by the orbits of parabolic points and the lengths of cusp excursions along geodesics. A well-known example is that a point $\eta\in\Lambda$ is dynamically badly approximable with respect to every parabolic point if and only if the geodesic $\geo\zero\eta$ has bounded cusp excursion lengths \cite[Proposition 1.21]{FSU4}. The following observation is in a similar vein:

\begin{observation}
\label{observationVWAbetat}
For $\eta\in\Lambda$, we have:
\begin{align*}
\eta\in\bigcup_{p\in P}\VWA_\bp
&\;\;\Leftrightarrow\;\; \limsup_{\substack{\xi\in\Lbp \\ t_\xi\to\infty}} \frac{\lb\xi|\eta\rb - t_\xi}{t_\xi} > 0
\;\;\Leftrightarrow\;\; \limsup_{t\to\infty} \frac{b(\eta,t)}{t} > 0\\
\eta\in\bigcup_{p\in P}\Liou_\bp
&\;\;\Leftrightarrow\;\; \limsup_{\substack{\xi\in\Lbp \\ t_\xi\to\infty}} \frac{\lb\xi|\eta\rb - t_\xi}{t_\xi} = \infty
\;\;\Leftrightarrow\;\; \limsup_{t\to\infty} \frac{b(\eta,t)}{t} = 1.
\end{align*}
\end{observation}
\begin{proof}
If $\xi = g(\bp)\in\Lbp$, then $\dogo g \gtrsim_\plus t_\xi$, with $\asymp_\plus$ for at least one value of $g$ (Lemma \ref{lemmatop}). Thus
\begin{align*}
\max_{\bp\in P}\omega_\bp(\eta) = \max_{\bp\in P}\limsup_{\substack{g\in G \\ g(\bp)\to\eta}} \frac{\log\Dist(g(\bp),\eta)}{\dogo g}
= \limsup_{\substack{\xi\in\Lbp \\ \xi\to\eta}} \frac{\lb \xi|\eta\rb}{t_\xi},
\end{align*}
so
\begin{equation}
\label{VWAbetat1}
\limsup_{\substack{\xi\in\Lbp \\ t_\xi\to\infty}} \frac{\lb\xi|\eta\rb - t_\xi}{t_\xi} = \max_{\bp\in P}\omega_\bp(\eta) - 1.
\end{equation}
On the other hand, it is readily verified that if $\geo\zero\eta$ intersects $H_\xi$, then the function $f(t) = b(\eta,t)/t$ attains its maximum at $t = \lb \xi|\eta\rb_\zero$, at which $f(t) = \lb \xi|\eta\rb_\zero - t_\xi$. Thus we have that
\begin{equation}
\label{VWAbetat2}
\begin{split}
\limsup_{t\to\infty} \frac{b(\eta,t)}{t}
= \limsup_{\substack{\xi\in\Lbp \\ t_\xi\to\infty}} \sup_{\substack{t > 0 \\ \eta_t\in H_\xi}} \frac{b(\eta,t)}{t}
&= \limsup_{\substack{\xi\in\Lbp \\ t_\xi\to\infty}} \frac{\lb \xi|\eta\rb_\zero - t_\xi}{\lb \xi|\eta\rb_\zero}\\
&= 1 - \frac{1}{\max_{\bp\in P}\omega_\bp(\eta)}
\end{split}
\end{equation}
Since
\[
\max_{\bp\in P}\omega_\bp(\eta) \begin{cases}
= \infty & \eta\in\bigcup_{p\in P}\Liou_\bp\\
\in (1,\infty) & \eta\in\bigcup_{p\in P}\VWA_\bp\butnot \bigcup_{p\in P}\Liou_\bp\\
= 1 & \eta\notin\bigcup_{p\in P}\VWA_\bp
\end{cases}
\]
applying \eqref{VWAbetat1} and \eqref{VWAbetat2} completes the proof.
\end{proof}

We are now ready to state our main theorem regarding the relation between exact dimensionality and dynamical Diophantine approximation:

\begin{theorem}
\label{theoremequivalentVWA}
The following are equivalent:
\begin{itemize}
\item[(A)] $\mu(\VWA_\bp) = 0\all \bp\in P$.
\item[(B)] $\mu$ is exact dimensional.
\item[(C)] $\HD(\mu) = \delta$.
\item[(D)] $\mu(\VWA_\xi) = 0\all \xi\in\Lambda$.
\end{itemize}
\end{theorem}
The implication (B) \implies (C) is part of Proposition \ref{propositionexactdimensional}, while (C) \implies (D) is an immediate consequence of Theorem \ref{theoremFSUjarnik}, and (D) \implies (A) is trivial. Thus we demonstrate (A) \implies (B):

\begin{proof}[Proof of \text{(A) \implies (B)}]
Fix $\eta\in\Lambda\butnot\bigcup_{\bp\in P}\VWA_\bp$ and $t > 0$. Suppose that $\eta_t\in H_\xi$ for some $\xi\in\Lbp$. Let $t_- < t < t_+$ satisfy
\[
t_- \asymp_\plus t_\xi, \; t_+ \asymp_\plus 2\lb \xi|\eta\rb_\zero - t_\xi, \text{ and } \eta_{t_\pm} \notin\bigcup(\scrH).
\]
Then by Lemma \ref{lemmaCCBglobalmeasure},
\[
\mu(B(\eta,e^{-t_\pm})) \asymp_\times e^{-\delta t_\pm}.
\]
In particular
\begin{equation}
\label{dt-dt+bounds}
\delta t_- \lesssim_\plus \log\frac{1}{\mu(B(\eta,e^{-t}))} \lesssim_\plus \delta t_+.
\end{equation}
Now, by Observation \ref{observationVWAbetat}, we have
\[
\frac{t_+ - t_-}{t} \leq \frac{2(\lb \xi|\eta\rb_\zero - t_\xi + (\text{constant}))}{t_\xi} \to 0 \text{ as } t\to\infty.
\]
Since $t_- < t < t_+$, it follows that $t_-/t,t_+/t\to 1$ as $t\to \infty$. Combining with \eqref{dt-dt+bounds} gives $d_\mu(\eta) = \delta$ (cf. \eqref{dmudef}). But by assumption (A), this is true for $\mu$-a.e. $\eta\in\Lambda$. Thus $\mu$ is exact dimensional.
\end{proof}

\subsection{Examples and non-examples of exact dimensional measures}
Combining Theorems \ref{theoremequivalentVWA} and \ref{theoremFSUkhinchin} gives a necessary and sufficient condition for $\mu$ to be exact dimensional in terms of the convergence or divergence of a family of series. We can ask how often this condition is satisfied. Our first result shows that it is almost always satisfied:

\begin{theorem}
\label{theoremhlogseries}
If for all $\bp\in P$, the series
\begin{equation}
\label{hlogseries}
\sum_{h\in G_\bp} e^{-\delta\dogo h} \dogo h
\asymp_\times \sum_{h\in G_\bp} \|h\|_\bp^{-2\delta} \log\|h\|_\bp
\asymp_\times \sum_{k = 0}^\infty e^{-2\delta k} k \NN_\bp(e^k)
\end{equation}
converges, then $\mu$ is exact dimensional.
\end{theorem}
\begin{proof}
Fix $\bp\in P$ and $\kappa > 0$. We have
\begin{align*}
\SigmaVWA(\bp,\kappa)
&=_\pt \sum_{g\in G} e^{-\delta\dogo g} \sum_{\substack{h\in G_\bp \\ \dogo h \geq \kappa\dogo g/2}} e^{-\delta\dogo h}\\
&=_\pt \sum_{h\in G_\bp} e^{-\delta\dogo h} \sum_{\substack{g\in G \\ \dogo g \leq 2\dogo h/\kappa}} e^{-\delta\dogo g}\\
&\asymp_\times \sum_{h\in G_\bp} e^{-\delta\dogo h} \sum_{k\leq 2\dogo h/\kappa + 1} e^{-\delta k} \#\{g\in G : k - 1 \leq \dogo g < k\} \noreason\\
&\leq_\pt \sum_{h\in G_\bp} e^{-\delta\dogo h} \sum_{k\leq 2\dogo h/\kappa + 1} e^{-\delta k} \NN_{X,G}(k) \noreason\\
&\lesssim_\times \sum_{h\in G_\bp} e^{-\delta\dogo h} \sum_{k\leq 2\dogo h/\kappa + 1} 1 \by{Corollary \ref{corollaryupperorbitalbound}}\\
&\asymp_\times \sum_{h\in G_\bp} e^{-\delta\dogo h} \dogo h.
\end{align*}
So if \eqref{hlogseries} converges, so does $\SigmaVWA(\bp,\kappa)$, and thus by Theorems \ref{theoremFSUkhinchin} and \ref{theoremequivalentVWA}, $\mu$ is exact dimensional.
\end{proof}

\begin{corollary}
\label{corollarydeltagtrdeltap}
If for all $\bp\in P$, $\delta_\bp < \delta$, then $\mu$ is exact dimensional.
\end{corollary}
\begin{proof}
In this case, the series \eqref{hlogseries} converges, as it is dominated by $\Sigma_s(G_\bp)$ for any $s\in (\delta_\bp,\delta)$.
\end{proof}

\begin{remark}
\label{remarkdeltagtrdeltap}
Combining with Proposition \ref{propositiondeltagtrdeltap} shows that if $\mu$ is not exact dimensional, then
\[
\sum_{h\in G_\bp} e^{-\delta\dogo h} < \infty = \sum_{h\in G_\bp} e^{-\delta\dogo h} \dogo h
\]
for some $\bp\in P$. Equivalently,
\[
\sum_{k = 0}^\infty e^{-2\delta k} \NN_\bp(e^k) < \infty = \sum_{k = 0}^\infty e^{-2\delta k} k \NN_\bp(e^k).
\]
This creates a very ``narrow window'' for the orbital counting function $\NN_\bp$.
\end{remark}

\begin{corollary}
\label{corollarydoublingexactdim}
If $\mu$ is doubling, then $\mu$ is exact dimensional.
\end{corollary}
\begin{proof}
If $\mu$ is doubling, then
\begin{align*}
\sum_{k = 0}^\infty e^{-2\delta k} k \NN_\bp(e^k)
&=_\pt \sum_{k = 1}^\infty \sum_{\ell = 0}^\infty e^{-2\delta (k + \ell)} \NN_\bp(e^{k + \ell})\\
&=_\pt \sum_{k = 1}^\infty e^{-2\delta k} \w\II_\bp(e^k)\\
&\asymp_\times \sum_{k = 1}^\infty e^{-2\delta k} \NN_\bp(e^k). \by{Proposition \ref{propositiondoubling}}
\end{align*}
Remark \ref{remarkdeltagtrdeltap} completes the proof.
\end{proof}

Our next theorem shows that in certain circumstances, the converse holds in Theorem \ref{theoremhlogseries}. Specifically:

\begin{theorem}
\label{theoremhlogseriesconverse}
Suppose that $X$ is an $\R$-tree and that $G$ is the pure Schottky product (cf. Definition \ref{definitionschottkyproductRtree}) of a parabolic group $H$ with a lineal group $J$. Let $\bp$ be the global fixed point of $H$, so that $P = \{\bp\}$ is a complete set of inequivalent parabolic points for $G$ (Proposition \ref{propositionSproductGF}). Suppose that the series \eqref{hlogseries} diverges. Then $\mu$ is not exact dimensional; moreover, $\mu(\Liou_\bp) = 1$ and $\dim_H(\mu) = 0$.
\end{theorem}
\begin{example}
\label{examplenotexactdim}
To see that the hypotheses of this theorem are not vacuous, fix $\delta > 0$ and let
\[
f(R) = \frac{R^{2\delta}}{\log^2(R)},
\]
or more generally, let $f$ be any increasing function such that $\sum_1^\infty e^{-2\delta k} k f(e^k)$ diverges but $\sum_1^\infty e^{-2\delta k} f(e^k)$ converges. By Theorem \ref{theoremRtreeorbitalcounting}, there exists an $\R$-tree $X$ and a parabolic group $H\leq\Isom(X)$ such that $\NN_{\EE_\bp,H} \asymp_\times f$. Then the series \eqref{hlogseries} diverges, but $\Sigma_\delta(H) < \infty$. Thus, there exists a unique $r > 0$ such that
\[
2\sum_{n = 1}^\infty e^{-\delta r} = \frac{1}{\Sigma_\delta(H) - 1}\cdot
\]
Let $J = r\Z$, interpreted as a group acting by translations on the $\R$-tree $\R$, and let $G$ be the pure Schottky product of $H$ and $J$. Then $\Sigma_\delta(J) - 1 = 2\sum_{n = 1}^\infty e^{-\delta r}$, so $(\Sigma_\delta(H) - 1)(\Sigma_\delta(J) - 1) = 1$. Since the map $s\mapsto (\Sigma_s(H) - 1)(\Sigma_s(J) - 1)$ is decreasing, it follows from Proposition \ref{propositionschottkyproductRtree} that $\Delta(G) = [0,\delta]$. In particular, $G$ is of divergence type, so Standing Assumption \ref{assumptiondivergencetype} is satisfied.
\end{example}
\begin{remark}
Applying a BIM embedding allows us to construct an example acting on $\H^\infty$.
\end{remark}
\begin{proof}[Proof of Theorem \ref{theoremhlogseriesconverse}]
As in the proof of Proposition \ref{propositionschottkyproductRtree} we let
\[
E = (H\butnot\{\id\})(J\butnot\{\id\}),
\]
so that
\[
G = \bigcup_{n\geq 0} J E^n H.
\]
Define a measure $\theta$ on $E$ via the formula
\[
\theta = \sum_{g\in E} e^{-\delta\dogo g} \delta_g.
\]
By Proposition \ref{propositionschottkyproductRtree}, the fact that $G$ is of divergence type (Standing Assumption \ref{assumptiondivergencetype}), and the fact that $\Sigma_\delta(J),\Sigma_\delta(H) < \infty$ (Proposition \ref{propositiondeltagtrdeltap}), $\theta$ is a probability measure. The Patterson--Sullivan measure of $G$ is related to $\theta$ by the formula
\[
\mu = \frac{1}{\Sigma_\delta(J) - 1}\sum_{j\in J} e^{-\delta\dogo j} j_* \pi_*[\theta^\N],
\]
where $\pi:E^\N \to \Lambda_G$ is the coding map.

Next, we use a theorem proven independently by H. Kesten and A. Raugi,\Footnote{We are grateful to ``cardinal'' of \url{http://mathoverflow.net} and J. P. Conze, respectively, for these references.} which we rephrase here in the language of measure theory:

\begin{theorem}[\cite{Kesten}; see also \cite{Raugi}]
\label{theoremkesten}
Let $\theta$ be a probability measure on a set $E$, and let $f:E\to\R$ be a function such that
\[
\int |f(x)| \;\dee\theta(x) = \infty.
\]
Then for $\theta^\N$-a.e. $(x_n)_1^\infty \in \R^\N$,
\[
\limsup_{n\to\infty} \frac{|f(x_{n + 1})|}{|\sum_1^n f(x_i)|} = \infty.
\]
\end{theorem}

Letting $f(g) = \dogo g$, the theorem applies to our measure $\theta$, because our assumption that \eqref{hlogseries} diverges is equivalent to the assertion that $\int f(x) \;\dee\theta(x) = \infty$. Now fix $j\in J$, and let $(g_n)_1^\infty \in E^\N$ be a $\theta^\N$-typical point. Then the limit point
\[
\eta = \lim_{n\to\infty} j g_1\cdots g_n(\zero)
\]
represents a typical point with respect to the Patterson--Sullivan measure $\mu$. By Theorem \ref{theoremkesten}, we have
\[
\limsup_{n\to\infty} \frac{\dogo{g_{n + 1}}}{\sum_1^n \dogo{g_i}} = \infty.
\]
Write $g_i = h_i j_i$ for each $i$. Then $\dogo{g_i} = \dogo{h_i} + \dogo{j_i}$. Since $\int \dogo{j}\dee\theta(hj) < \infty$, the law of large numbers implies that $\lim_{n \to \infty} \frac{\dogo{j_{n + 1}}}{\sum_1^n \dogo{j_i}} < \infty$, so
\[
\limsup_{n\to\infty} \frac{\dogo{h_{n + 1}}}{\sum_1^n \dogo{g_i}} = \infty.
\]
But $\dogo{h_{n + 1}}$ represents the length of the excursion of the geodesic $\geo\zero\eta$ into the cusp corresponding to the parabolic point $g_1\cdots g_n(\bp)$. Combining with Observation \ref{observationVWAbetat} shows that $\eta\in\Liou_\bp$. Since $\eta$ was a $\mu$-typical point, this shows that $\mu(\Liou_\bp) = 1$. By Theorem \ref{theoremFSUjarnik}, this implies that $\HD(\mu) = 0$. By Observation \ref{observationHDdeltaPD}, $\mu$ is not exact dimensional.
\end{proof}

\appendix

\chapter{Open problems}
\label{appendixopenproblems}

\begin{problem}[Cf. Chapter \ref{sectionmodified}, Section \ref{subsectionpoincareirregular}]
Do there exist a hyperbolic metric space $X$ and a group $G$ such that $\w\delta(G) < \delta(G) = \infty$, but $\w\delta(G)$ cannot be ``computed from'' the modified Poincar\'e exponent of a locally compact group via Definition \ref{definitionmodified1}? This question is vague because a more precise version might be contradicted by Example \ref{examplepoincareextension}, in which a group $G$ is constructed such that $\w\delta(G) < \delta(G) = \infty$ but the closure of $G$ (in the compact-open topology) is not locally compact. In this case, $\w\delta(G)$ cannot be computed from $\w\delta(\cl G)$, but there is still a locally compact group ``hidden'' in the argument, namely the closure of $G\given \H^d = \iota_1(\Gamma)$. Is there any Poincar\'e irregular group whose construction is not somehow ``based on'' a locally compact group?
\end{problem}

\begin{problem}[Cf. Theorem \ref{theorembishopjonesmodified}]
If $G$ is a Poincar\'e irregular parabolic group, does the modified Poincar\'e exponent $\w\delta(G)$ have a geometric significance? Theorem \ref{theorembishopjonesmodified} does not apply directly since $G$ is elementary. It is tempting to claim that
\begin{equation}
\label{wdeltaparabolic} \tag{A.1}
\w\delta(G) = \inf\{\HD(\Lr(H)) : H\geq G \text{ nonelementary}\}
\end{equation}
(under some reasonable hypotheses about the isometry group of the space in question), but it seems that the right hand side is equal to infinity in most cases due to Proposition \ref{propositionSproductexponent}(iii). Note that by contrast, \eqref{wdeltaparabolic} is usually true for Poincar\'e regular groups; for example, it holds in the Standard Case \cite{Beardon1}.
\end{problem}

\begin{problem}[Cf. Chapter \ref{sectionparabolic}, Remark \ref{remarknilpotentembedding}]
Given a virtually nilpotent group $\Gamma$ which is not virtually abelian, determine whether there exists a homomorphism $\Phi:\Gamma\to \Isom(\BB)$ such that $\delta(\Phi(\Gamma)) = \alpha(\Gamma)/2$, where both quantities are defined in Section \ref{subsectionparabolicexponent}. Intuitively, this corresponds to the existence an equivariant embedding of $\Gamma$ into $\BB$ which approaches infinity ``as fast as possible''. It is known \cite[Theorem 1.3]{CTV1} that such an embedding cannot be quasi-isometric, but this by itself does not imply the non-existence of a homomorphism with the desired property.
\end{problem}

\begin{problem}[Cf. Chapter \ref{sectionparabolic}, Remark \ref{remarkheisenbergdeltainfty}]
Does there exist a strongly discrete parabolic subgroup of $\Isom(\H^\infty)$ isomorphic to the Heisenberg group which has infinite Poincar\'e exponent?\Footnote{It has been pointed out to us by X. Xie that this question is answered affirmatively by \cite[Proposition 3.10]{CTV1}, letting the function $f$ in that proposition be any function whose growth is sublinear.}
\end{problem}

\begin{problem}[Cf. Chapter \ref{sectionGF}, Section \ref{subsectionCBCCB}]
Is there any form of discreteness for which there exists a cobounded subgroup of $\Isom(\H)$ (for example, UOT discreteness)? If so, what is the strongest such form of discreteness?
\end{problem}


\begin{problem}[Cf. Chapter \ref{sectionGFmeasures}]
Can Theorem \ref{theoremhlogseriesconverse} be improved as follows?
\begin{conjecture*}
Let $X$ be a hyperbolic metric space and let $G\leq\Isom(X)$ be a geometrically finite group such that for some $\bp\in\Lbp$, the series \eqref{hlogseries} diverges. Then the $\delta$-quasiconformal measure $\mu$ is not exact dimensional.
\end{conjecture*}
What if some of the hypotheses of this conjecture are strengthened, e.g. $X$ is strongly hyperbolic (e.g. $X = \H^\infty$), or $G$ is a Schottky product of a parabolic group with a lineal group?
\end{problem}

\chapter{Index of defined terms}
\label{appendixindex}

See also Conventions \ref{conventionimplied}-\ref{conventionstandard} on pages \pageref{conventionimplied}, \pageref{conventionH}, and \pageref{conventiontriples}.

\begin{itemize}
\item\emph{acts irreducibly}: Definition \ref{definitionreducibly}, p.\pageref{definitionreducibly}
\item\emph{acts properly discontinuously}: Definition \ref{definitionproperdiscontinuity}, p.\pageref{definitionproperdiscontinuity}
\item\emph{acts reducibly}: Definition \ref{definitionreducibly}, p.\pageref{definitionreducibly}
\item\emph{algebraic hyperbolic space}: Definition \ref{definitionROSSONCT}, p.\pageref{definitionROSSONCT}
\item\emph{attracting fixed point}: Definition \ref{definitionderivativefixedpoint}, p.\pageref{definitionderivativefixedpoint}
\item\emph{ball model}: \6\ref{subsubsectionballmodel}, p.\pageref{subsubsectionballmodel}
\item\emph{bi-infinite geodesic}: Definition \ref{definitiongeodesic}, p.\pageref{definitiongeodesic}
\item\emph{BIM embedding}: Definition \ref{definitionBIM}, p.\pageref{definitionBIM}
\item\emph{BIM representation}: Definition \ref{definitionBIM}, p.\pageref{definitionBIM}
\item\emph{bordification}: Definition \ref{definitiongromovboundary}, p.\pageref{definitiongromovboundary}
\item\emph{$\xi$-bounded}: Definition \ref{definitionxibounded}, p.\pageref{definitionxibounded}
\item\emph{bounded parabolic point}: Definition \ref{definitionboundedparabolic}, p.\pageref{definitionboundedparabolic}
\item\emph{Busemann function}: \eqref{busemanndef}, p.\pageref{busemanndef}
\item\emph{CAT(-1) inequality}: \eqref{CAT}, p.\pageref{CAT}
\item\emph{CAT(-1) space}: Definition \ref{definitionCAT}, p.\pageref{definitionCAT}
\item\emph{Cayley graph}: Example \ref{examplecayleygraph}, p.\pageref{examplecayleygraph}
\item\emph{Cayley hyperbolic plane}: Remark \ref{remarkH2O}, p.\pageref{remarkH2O}
\item\emph{Cayley metric}: Example \ref{examplecayleygraph}, p.\pageref{examplecayleygraph}
\item\emph{center} (of a triangle in an $\R$-tree): Definition \ref{definitioncenterRtree}, p.\pageref{definitioncenterRtree}
\item\emph{center} (of a horoball): Definition \ref{definitionhoroball}, p.\pageref{definitionhoroball}
\item\emph{cobounded}: Definition \ref{definitionCB}, p.\pageref{definitionCB}
\item\emph{codoubling} (function): Definition \ref{definitiondoubling}, p.\pageref{definitiondoubling}
\item\emph{convergence type}: Definition \ref{definitiondivergencetype}, p.\pageref{definitiondivergencetype}
\item\emph{compact-open topology (COT)}: p.\pageref{subsectiontopologies}
\item\emph{compact type, semigroup of}: Definition \ref{definitioncompacttype}, p.\pageref{definitioncompacttype}
\item\emph{comparison point}: p.\pageref{subsectionCAT}, Definition \ref{definitiontriangles}, p.\pageref{definitiontriangles}
\item\emph{comparison triangle}: Example \ref{exampletreetriangle}, p.\pageref{exampletreetriangle}; Definition \ref{definitiontriangles}, p.\pageref{definitiontriangles}
\item\emph{compatible} (regarding a metametric and a topology): Definition \ref{definitioncompatiblemetametric}, p.\pageref{definitioncompatiblemetametric}
\item\emph{complete set of inequivalent parabolic points}: Definition \ref{definitioncompleteset}, p.\pageref{definitioncompleteset}
\item\emph{cone}: \eqref{coneconstruction}, p.\pageref{coneconstruction}
\item\emph{conformal measure}: Definition \ref{definitionquasiconformal}, p.\pageref{definitionquasiconformal}
\item\emph{conical convergence}: p.\pageref{subsectionmodes}
\item\emph{connected graph}: Definition \ref{definitiongraphmetrization}, p.\pageref{definitiongraphmetrization}
\item\emph{contractible cycles} (property of a graph): Definition \ref{definitioncontractiblecycles}, p.\pageref{definitioncontractiblecycles}
\item\emph{convex-cobounded}: Definition \ref{definitionCCB}, p.\pageref{definitionCCB}
\item\emph{convex hull}: Definition \ref{definitionconvexhull}, p.\pageref{definitionconvexhull}
\item\emph{convex}: \eqref{convex}, p.\pageref{convex}
\item\emph{convex core}: Definition \ref{definitionconvexcore}, p.\pageref{definitionconvexcore}
\item\emph{cycle}: \eqref{cycle}, p.\pageref{cycle}
\item\emph{Dirichlet domain}: Definition \ref{definitiondirichletdomain}, p.\pageref{definitiondirichletdomain}
\item\emph{divergence type}: Definition \ref{definitiondivergencetype}, p.\pageref{definitiondivergencetype}
\item\emph{domain of reflexivity}: Definition \ref{definitionmetametric}, p.\pageref{definitionmetametric}
\item\emph{doubling} (metric space): Footnote \ref{footnotedoubling}, p.\pageref{footnotedoubling}
\item\emph{doubling} (function): Definition \ref{definitiondoubling}, p.\pageref{definitiondoubling}
\item\emph{doubling} (measure): Section \ref{subsectiondoubling}, p.\pageref{subsectiondoubling}
\item\emph{dynamical derivative}: Proposition \ref{propositiondynamicalderivative}, p.\pageref{propositiondynamicalderivative}
\item\emph{Edelstein-type isometry}: Definition \ref{definitionedelstein}, p.\pageref{definitionedelstein}
\item\emph{elementary}: Definition \ref{definitionelementary}, p.\pageref{definitionelementary}
\item\emph{elliptic isometry}: Definition \ref{definitionclassification1}, p.\pageref{definitionclassification1}
\item\emph{elliptic semigroup}: Definition \ref{definitionclassification2}, p.\pageref{definitionclassification2}
\item\emph{ergodic}: Definition \ref{definitionergodic}, p.\pageref{definitionergodic}
\item\emph{equivalent} (for Gromov sequences): Definition \ref{definitiongromovsequence}, p.\pageref{definitiongromovsequence}
\item\emph{extended visual metric}: Proposition \ref{propositionwbarDist}, p.\pageref{propositionwbarDist}
\item\emph{fixed point (neutral/attracting/repelling)}: Definition \ref{definitionclassification1}, p.\pageref{definitionclassification1}
\item\emph{fixed point (parabolic)}: Definition \ref{definitionparabolic}, p.\pageref{definitionparabolic}
\item\emph{focal semigroup}: Definition \ref{definitionCCMT}, p.\pageref{definitionCCMT}
\item\emph{free group}: Remark \ref{remarkfreegroup}, p.\pageref{remarkfreegroup}
\item\emph{free product}: Section \ref{subsectionfreeproducts}, p.\pageref{subsectionfreeproducts}
\item\emph{free semigroup}: Remark \ref{remarkfreegroup}, p.\pageref{remarkfreegroup}
\item\emph{general type, semigroup of}: Definition \ref{definitionCCMT}, p.\pageref{definitionCCMT}
\item\emph{generalized convergence type}: Definition \ref{definitionmodifiedexponent}, p.\pageref{definitionmodifiedexponent}
\item\emph{generalized divergence type}: Definition \ref{definitionmodifiedexponent}, p.\pageref{definitionmodifiedexponent}
\item\emph{generalized polar coordinate functions}: Definition \ref{definitionpolarcoordinates}, p.\pageref{definitionpolarcoordinates}
\item\emph{geodesic metric space}: Remark \ref{remarkgeodesicmetric}, p.\pageref{remarkgeodesicmetric}
\item\emph{geodesic segment}: Remark \ref{remarkgeodesicmetric}, p.\pageref{remarkgeodesicmetric}
\item\emph{geodesic triangle}: p.\pageref{subsectionCAT}, Definition \ref{definitiontriangles}, p.\pageref{definitiontriangles}
\item\emph{geometric product}: Example \ref{examplegeometricproducts}, p.\pageref{examplegeometricproducts}
\item\emph{geodesic path}: Section \ref{subsectiongraphtheory}, p.\pageref{subsectiongraphtheory}
\item\emph{geodesic ray/line}: Definition \ref{definitiongeodesic}, p.\pageref{definitiongeodesic}
\item\emph{geometric realization}: Definition \ref{definitiongraphmetrization}, p.\pageref{definitiongraphmetrization}
\item\emph{geometric graph}: Definition \ref{definitiongraphmetrization}, p.\pageref{definitiongraphmetrization}
\item\emph{geometrically finite}: Definition \ref{definitionGF}, p.\pageref{definitionGF}
\item\emph{Gromov boundary}: Definition \ref{definitiongromovboundary}, p.\pageref{definitiongromovboundary}
\item\emph{Gromov hyperbolic}: Definition \ref{definitiongromovhyperbolic}, p.\pageref{definitiongromovhyperbolic}
\item\emph{Gromov's inequality}: \eqref{gromov}, p.\pageref{gromov}
\item\emph{Gromov product}: \eqref{gromovproduct}, p.\pageref{gromovproduct}
\item\emph{Gromov sequence}: Definition \ref{definitiongromovsequence}, p.\pageref{definitiongromovsequence}
\item\emph{Gromov triple}: Definition \ref{definitiongromovtriple}, p.\pageref{definitiongromovtriple}
\item\emph{global fixed points}: Notation \ref{notationglobalfixedpoints}, p.\pageref{notationglobalfixedpoints}
\item\emph{growth rate}: \eqref{bassguivarch}, p.\pageref{bassguivarch}; Definition \ref{definitiongrowthrate}, p.\pageref{definitiongrowthrate}
\item\emph{global Schottky product}: Definition \ref{definitionschottkyproduct}, p.\pageref{definitionschottkyproduct}
\item\emph{group of isometries}: p.\pageref{subsectionisometries}
\item\emph{Haagerup property}: \6\ref{subsubsectionhaagerup}, p.\pageref{subsubsectionhaagerup}
\item\emph{half-space}: Remark \ref{remarkhalfspace}, p.\pageref{remarkhalfspace}
\item\emph{half-space model}: \6\ref{subsubsectionE}, p.\pageref{subsubsectionE}
\item\emph{horoball}: Definition \ref{definitionhoroball}, p.\pageref{definitionhoroball}
\item\emph{horospherical convergence}: Definition \ref{definitionhorosphericalconvergence}, p.\pageref{definitionhorosphericalconvergence}
\item\emph{horospherical limit set}: Definition \ref{definitionlimitset}, p.\pageref{definitionlimitset}
\item\emph{hyperbolic}: Definition \ref{definitiongromovhyperbolic}, p.\pageref{definitiongromovhyperbolic}
\item\emph{hyperboloid model}: \6\ref{subsectionhyperboloidmodel}, p.\pageref{subsectionhyperboloidmodel}
\item\emph{implied constant}: Convention \ref{conventionimplied}, p.\pageref{conventionimplied}
\item\emph{inward focal}: Definition \ref{definitionfocal}, p.\pageref{definitionfocal}
\item\emph{irreducible action}: Definition \ref{definitionreducibly}, p.\pageref{definitionreducibly}
\item\emph{isomorphism} (between pairs $(X,\bord X)$ and $(Y,\bord Y)$): p.\pageref{subsectiontotallygeodesic}
\item\emph{length spectrum}: Remark \ref{remarklengthspectrum}, p.\pageref{remarklengthspectrum}
\item\emph{limit set} (of a semigroup): Definition \ref{definitionlimitset}, p.\pageref{definitionlimitset}
\item\emph{limit set} (of a partition structure): Definition \ref{definitionpartitionlimitset}, p.\pageref{definitionpartitionlimitset}
\item\emph{lineal semigroup}: Definition \ref{definitionCCMT}, p.\pageref{definitionCCMT}
\item\emph{Lorentz boosts}: \eqref{lorentzboost}, p.\pageref{lorentzboost}
\item\emph{lower central series}: \6\ref{subsubsectionnilpotent}, p.\pageref{subsubsectionnilpotent}
\item\emph{lower polynomial growth rate}: Definition \ref{definitiongrowthrate}, p.\pageref{definitiongrowthrate}
\item\emph{loxodromic isometry}: Definition \ref{definitionclassification1}, p.\pageref{definitionclassification1}
\item\emph{loxodromic semigroup}: Definition \ref{definitionclassification2}, p.\pageref{definitionclassification2}
\item\emph{Margulis's lemma}: Proposition \ref{propositionmargulislemma}, p.\pageref{propositionmargulislemma}
\item\emph{metametric}: Definition \ref{definitionmetametric}, p.\pageref{definitionmetametric}
\item\emph{metric derivative}: p.\pageref{subsectionderivatives}, p.\pageref{subsubsectionderivativemaps}
\item\emph{moderately discrete (MD)}: Definition \ref{definitiondiscreteness}, p.\pageref{definitiondiscreteness}
\item\emph{modified Poincar\'e exponent}: Definition \ref{definitionmodifiedexponent}, p.\pageref{definitionmodifiedexponent}
\item\emph{natural action}: (on a Cayley graph) Remark \ref{remarknaturalaction}, p.\pageref{remarknaturalaction}
\item\emph{natural map} (from a free product): Section \ref{subsectionfreeproducts}, p.\pageref{subsectionfreeproducts}
\item\emph{$\rho$-net}: Footnote \ref{footnoterhonet}, p.\pageref{footnoterhonet}
\item\emph{neutral fixed point}: Definition \ref{definitionderivativefixedpoint}, p.\pageref{definitionderivativefixedpoint}
\item\emph{nilpotent}: \6\ref{subsubsectionnilpotent}, p.\pageref{subsubsectionnilpotent}
\item\emph{nilpotency class}: \6\ref{subsubsectionnilpotent}, p.\pageref{subsubsectionnilpotent}
\item\emph{nonelementary}: Definition \ref{definitionelementary}, p.\pageref{definitionelementary}
\item\emph{orbital counting function}: Remark \ref{remarkorbitalcounting}, p.\pageref{remarkorbitalcounting}
\item\emph{outward focal}: Definition \ref{definitionfocal}, p.\pageref{definitionfocal}
\item\emph{parabolic isometry}: Definition \ref{definitionclassification1}, p.\pageref{definitionclassification1}
\item\emph{parabolic fixed point}: Definition \ref{definitionparabolic}, p.\pageref{definitionparabolic}
\item\emph{parabolic semigroup}: Definition \ref{definitionclassification2}, p.\pageref{definitionclassification2}
\item\emph{parameterization} (of a geodesic): Remark \ref{remarkgeodesicmetric}, p.\pageref{remarkgeodesicmetric}
\item\emph{partition structure}: Definition \ref{definitionpartitionstructure}, p.\pageref{definitionpartitionstructure}
\item\emph{path}: Section \ref{subsectiongraphtheory}, p.\pageref{subsectiongraphtheory}
\item\emph{path metric}: Definition \ref{definitiongraphmetrization}, p.\pageref{definitiongraphmetrization}, \eqref{pathmetricv3}, p.\pageref{pathmetricv3}
\item\emph{Poincar\'e exponent}: Definition \ref{definitionpoincareexponent}, p.\pageref{definitionpoincareexponent}
\item\emph{Poincar\'e extension}: Observation \ref{observationpoincareextension}, p.\pageref{observationpoincareextension}
\item\emph{Poincar\'e integral}: \eqref{IsG}, p.\pageref{IsG}
\item\emph{Poincar\'e regular/irregular}: p.\pageref{corollarygeneralizedtypes}
\item\emph{Poincar\'e set}: Notation \ref{notationpoincareset}, p.\pageref{notationpoincareset}
\item\emph{Poincar\'e series}: Definition \ref{definitionpoincareexponent}, p.\pageref{definitionpoincareexponent}
\item\emph{polynomial growth rate}: \eqref{bassguivarch}, p.\pageref{bassguivarch}; Definition \ref{definitiongrowthrate}, p.\pageref{definitiongrowthrate}
\item\emph{pre-doubling} (parabolic group): Definition \ref{definitionpredoubling}, p.\pageref{definitionpredoubling}
\item\emph{proper}: Remark \ref{remarkproper}, p.\pageref{remarkproper}
\item\emph{properly discontinuous (PrD)}: Definition \ref{definitionproperdiscontinuity}, p.\pageref{definitionproperdiscontinuity}
\item\emph{pure Schottky product}: Definition \ref{definitionschottkyproductRtree}, p.\pageref{definitionschottkyproductRtree}
\item\emph{quasiconformal measure}: Definition \ref{definitionquasiconformal}, p.\pageref{definitionquasiconformal}
\item\emph{quasiconvex core}: Definition \ref{definitionconvexcore}, p.\pageref{definitionconvexcore}
\item\emph{quasi-isometry/quasi-isometric}: Definition \ref{definitionquasiisometry}, p.\pageref{definitionquasiisometry}
\item\emph{radial convergence}: Definition \ref{definitionradialconvergence}, p.\pageref{definitionradialconvergence}
\item\emph{radial limit set}: Definition \ref{definitionlimitset}, p.\pageref{definitionlimitset}
\item\emph{Radon}: Remark \ref{remarkradon}, p.\pageref{remarkradon}
\item\emph{rank} (of an abelian group): \6\ref{subsubsectionnilpotent}, p.\pageref{subsubsectionnilpotent}
\item\emph{reducible action}: Definition \ref{definitionreducibly}, p.\pageref{definitionreducibly}
\item\emph{regularly geodesic}: Definition \ref{definitionregularlygeodesic}, p.\pageref{definitionregularlygeodesic}
\item\emph{repelling fixed point}: Definition \ref{definitionderivativefixedpoint}, p.\pageref{definitionderivativefixedpoint}
\item\emph{Samuel--Smirnov compactification}: Proposition \ref{propositionsamuelsmirnov}, p.\pageref{propositionsamuelsmirnov}
\item\emph{Schottky group}: Definition \ref{definitionschottkygroup}, p.\pageref{definitionschottkygroup}
\item\emph{Schottky position}: Definition \ref{definitionschottkyproduct}, p.\pageref{definitionschottkyproduct}
\item\emph{Schottky product}: Definition \ref{definitionschottkyproduct}, p.\pageref{definitionschottkyproduct}
\item\emph{Schottky semigroup}: Definition \ref{definitionschottkygroup}, p.\pageref{definitionschottkygroup}
\item\emph{Schottky system}: Definition \ref{definitionschottkyproduct}, p.\pageref{definitionschottkyproduct}
\item\emph{$\rho$-separated set}: Footnote \ref{footnoterhoseparated}, p.\pageref{footnoterhoseparated}
\item\emph{sesquilinear form}: p.\pageref{subsectionROSSONCT}
\item\emph{shadow}: Definition \ref{definitionshadow}, p.\pageref{definitionshadow}
\item\emph{similarity}: Observation \ref{observationpoincareextension}, p.\pageref{observationpoincareextension}
\item\emph{simplicial tree}: Definition \ref{definitionweightedtree}, p.\pageref{definitionweightedtree}
\item\emph{$\F$-skew linear}: \eqref{sigmaTdef}, p.\pageref{sigmaTdef}
\item\emph{skew-symmetric}: p.\pageref{subsectionROSSONCT}
\item\emph{Standard Case}: Convention \ref{conventionstandard}, p.\pageref{conventionstandard}
\item\emph{standard parameterization}: p.\pageref{lemmageodesicconvergence}
\item\emph{stapled union}: Definition \ref{definitionstapledunion}, p.\pageref{definitionstapledunion}
\item\emph{strong operator topology (SOT)}: p.\pageref{subsectiontopologies}
\item\emph{strongly discrete (SD)}: Definition \ref{definitiondiscreteness}, p.\pageref{definitiondiscreteness}, Remark \ref{remarkstronglydiscrete}, p.\pageref{remarkstronglydiscrete}
\item\emph{strongly (Gromov) hyperbolic}: Definition \ref{definitionstronglyhyperbolic}, p.\pageref{definitionstronglyhyperbolic}
\item\emph{strongly separated Schottky group/product/system}: Definition \ref{definitionstrongseparation}, p.\pageref{definitionstrongseparation}
\item\emph{substructure} (of a partition structure): Definition \ref{definitionpartitionsubstructure}, p.\pageref{definitionpartitionsubstructure}
\item\emph{$s$-thick}: Definition \ref{definitionpartitionstructure}, p.\pageref{definitionpartitionstructure}
\item\emph{topological discreteness}: Definition \ref{definitionparametricdiscreteness}, p.\pageref{definitionparametricdiscreteness}
\item\emph{totally geodesic subset}: Definition \ref{definitiontotallygeodesic}, p.\pageref{definitiontotallygeodesic}
\item\emph{tree, simplicial}: Definition \ref{definitionweightedtree}, p.\pageref{definitionweightedtree}
\item\emph{tree} (on $\N$): Definition \ref{definitiontree}, p.\pageref{definitiontree}
\item\emph{$\R$-tree}: Definition \ref{definitionRtree}, p.\pageref{definitionRtree}
\item\emph{$\Z$-tree}: Definition \ref{definitionweightedtree}, p.\pageref{definitionweightedtree}
\item\emph{tree-geometric}: Definition \ref{definitiontreegeometric}, p.\pageref{definitiontreegeometric}
\item\emph{tree triangle}: p.\pageref{subsectionCAT}
\item\emph{Tychonoff topology}: p.\pageref{propositionCOTSOT}
\item\emph{uniform operator topology (UOT)}: p.\pageref{subsectiontopologies}
\item\emph{uniformly radial convergence}: Definition \ref{definitionradialconvergence}, p.\pageref{definitionradialconvergence}
\item\emph{uniformly radial limit set}: Definition \ref{definitionlimitset}, p.\pageref{definitionlimitset}
\item\emph{uniquely geodesic metric space}: Remark \ref{remarkgeodesicmetric}, p.\pageref{remarkgeodesicmetric}
\item\emph{unweighted simplicial tree}: Definition \ref{definitionweightedtree}, p.\pageref{definitionweightedtree}
\item\emph{upper polynomial growth rate}: Definition \ref{definitiongrowthrate}, p.\pageref{definitiongrowthrate}
\item\emph{virtually nilpotent}: \6\ref{subsubsectionnilpotent}, p.\pageref{subsubsectionnilpotent}
\item\emph{visual metric}: p.\pageref{propositionDist}
\item\emph{weakly discrete (WD)}: Definition \ref{definitiondiscreteness}, p.\pageref{definitiondiscreteness}
\item\emph{weakly separated Schottky group/product/system}: Definition \ref{definitionstrongseparation}, p.\pageref{definitionstrongseparation}
\item\emph{weighted Cayley graph}: Example \ref{examplecayleygraph}, p.\pageref{examplecayleygraph}
\item\emph{weighted undirected graph}: Definition \ref{definitiongraphmetrization}, p.\pageref{definitiongraphmetrization}
\end{itemize}

\backmatter
\bibliographystyle{amsplain}
\bibliography{bibliography}

\printindex

\end{document}